\documentstyle[12pt,fleqn,makeidx,epsfig]{book}

\makeindex

\begin{document}

\pagestyle{plain}

\title{Complex Numbers in $n$ Dimensions}

\author{Silviu Olariu\\
Institute of Physics and Nuclear Engineering, Tandem Laboratory\\
76900 Magurele, P.O. Box MG-6, Bucharest, Romania\\
e-mail: olariu@ifin.nipne.ro}

\date{23 August 2000}

\maketitle

\tableofcontents

\newpage

\section*{Introduction}

A regular, two-dimensional complex number $x+iy$ can be represented
geometrically by the modulus $\rho=(x^2+y^2)^{1/2}$ and by the polar angle
$\theta=\arctan(y/x)$. The modulus $\rho$ is multiplicative and the polar angle
$\theta$ is additive upon the multiplication of ordinary complex numbers.

The quaternions of Hamilton are a system of hypercomplex numbers defined in
four dimensions, the multiplication being a noncommutative operation, \cite{1}
and many other hypercomplex systems are possible, \cite{2a}-\cite{2b} but these interesting
hypercomplex systems do not have all the required properties of regular,
two-dimensional complex numbers which rendered possible the development of the
theory of functions of a complex variable.

Two distinct systems of hypercomplex numbers in $n$ dimensions will be
described in this work, for which the multiplication is associative and
commutative, and which are rich enough in properties such that exponential 
and trigonometric forms
exist and the concepts of analytic n-complex function,  contour integration and
residue can be defined. \cite{2c} The n-complex numbers described in this work have the
form $u=x_0+h_1x_1+\cdots+h_{n-1}x_{n-1}$, where $h_1,...,h_{n-1}$ are the
hypercomplex bases and the variables $x_0,...,x_{n-1}$ are real numbers, unless
otherwise stated. If the n-complex number $u$ is represented by the point $A$
of coordinates $x_0,x_1,...,x_{n-1}$, the position of the point $A$ can be
described with the aid of the modulus $d=(x_0^2+x_1^2+\cdots+x_{n-1}^2)^{1/2}$
and of $n-1$ angular variables. 

The first type of hypercomplex numbers described in this work is
characterized by the presence, in an even number of dimensions
$n\geq 4$ , of two polar
axes, and by the presence, in an odd number of dimensions, of one polar axis.
Therefore, these numbers will be called polar hypercomplex numbers in $n$
dimensions. One polar axis is the normal through the origin $O$ to the
hyperplane $v_+=0$, where $v_+=x_0+x_1+\cdots+x_{n-1}$.  In an even number $n$
of dimensions, the second polar axis is the normal through the origin $O$ to
the hyperplane $v_-=0$, where $v_-=x_0-x_1+\cdots+x_{n-2}-x_{n-1}$.  Thus, in
addition to the distance $d$, the position of the point $A$ can be specified,
in an even number of dimensions, by 2 polar angles $\theta_+,\theta_-$, by
$n/2-2$ planar angles $\psi_k$, and by $n/2-1$ azimuthal angles $\phi_k$. In an
odd number of dimensions, the position of the point $A$ is specified by $d$, by
1 polar angle $\theta_+$, by $(n-3)/2$ planar angles $\psi_{k-1}$, and by
$(n-1)/2$ azimuthal angles $\phi_k$.  The multiplication rules for the polar
hypercomplex bases $h_1,...,h_{n-1}$ are $h_j h_k =h_{j+k}$ if $0\leq j+k\leq
n-1$, and $h_jh_k=h_{j+k-n}$ if $n\leq j+k\leq 2n-2$, where $h_0=1$.

The other type of hypercomplex numbers described in this work exists as a
distinct entity only when the number of dimensions $n$ of the space is even.
The position of the point $A$ is specified, in addition to the distance $d$, by
$n/2-1$ planar angles $\psi_k$ and by $n/2$ azimuthal angles $\phi_k$.  These
numbers will be called planar hypercomplex numbers.  The multiplication
rules for the planar hypercomplex bases $h_1,...,h_{n-1}$ are $h_j h_k
=h_{j+k}$ if $0\leq j+k\leq n-1$, and $h_jh_k=-h_{j+k-n}$ if $n\leq j+k\leq
2n-2$, where $h_0=1$. For $n=2$, the planar hypercomplex numbers become the
usual 2-dimensional complex numbers $x+iy$.

The development of analytic functions of hypercomplex variables was rendered
possible by the existence of an exponential form of the n-complex numbers.
The azimuthal angles $\phi_k$, which are cyclic variables,
appear in these forms at the exponent, and lead to the concept of n-dimensional
hypercomplex residue.  Expressions are given for the elementary functions of
n-complex variable.  In particular, the exponential function of an n-complex
number is expanded in terms of functions called in this work n-dimensional
cosexponential functions of the polar and respectively planar type.  The polar
cosexponential functions are a generalization to $n$ dimensions of the
hyperbolic functions 
$\cosh y, \sinh y$, and the planar cosexponential functions are a
generalization to $n$ dimensions of the trigonometric functions $\cos y, \sin
y$. Addition theorems and other relations are obtained for the n-dimensional
cosexponential functions.

Many of the properties of 2-dimensional complex functions can be extended to
hypercomplex numbers in $n$ dimensions.  Thus, the functions $f(u)$ of an
n-complex variable which are defined by power series have derivatives
independent of the direction of approach to the point under consideration. If
the n-complex function $f(u)$ of the n-complex variable $u$ is written in terms
of the real functions $P_k(x_0,...,x_{n-1}), k=0,...,n-1$, then relations of
equality exist between the partial derivatives of the functions $P_k$.  The
integral $\int_A^B f(u) du$ of an n-complex function between two points $A,B$
is independent of the path connecting $A,B$, in regions where $f$ is regular.
If $f(u)$ is an analytic n-complex function, then $\oint_\Gamma f(u)du/(u-u_0)$
is expressed in this work in terms of the n-dimensional hypercomplex residue
$f(u_0)$.

In the case of polar complex numbers, a polynomial can be written as a product
of linear or 
quadratic factors, although several factorization are in general possible.  In
the case of planar hypercomplex numbers, a polynomial can always be written as
a product of 
linear factors, although, again, several factorization are in general possible.

The work presents a detailed analysis of the hypercomplex numbers in 2, 3 and 4
dimensions, then presents the properties of hypercomplex numbers in 5 and
6 dimensions, and it continues with a detailed analysis of polar and planar
hypercomplex numbers in $n$ dimensions. 
The essence of this work is the interplay between the algebraic, the geometric 
and the analytic facets of the relations.

\chapter{Hyperbolic Complex Numbers in Two Dimensions}

A system of hypercomplex numbers in 2 dimensions is described in this chapter,
for which the multiplication is associative and commutative, and for which
an exponential form and the concepts of analytic twocomplex
function and contour integration can be defined.
The twocomplex numbers introduced in this chapter have 
the form $u=x+\delta y$, the variables $x, y$ being real 
numbers.  
The multiplication rules for the complex units $1, \delta$ are
$1\cdot \delta=\delta, \delta^2=1$. In a geometric representation, the
twocomplex number $u$ is represented by the point $A$ of coordinates $(x,y).$
The product of two twocomplex numbers is equal to zero if both numbers are
equal to zero, or if one of
the twocomplex numbers lies on the line $x=y$ and the other on the 
line $x=-y$.

The exponential form of a twocomplex number, defined for $x+y>0,
x-y>0$, is $u=\rho\exp(\delta\lambda/2)$,
where the amplitude is $\rho=(x^2-y^2)^{1/2}$ and the argument is
$\lambda=\ln\tan\theta$, $\tan\theta=(x+y)/(x-y)$,
$0<\theta<\pi/2$. 
The trigonometric form of a twocomplex number is
$u=d\sqrt{\sin 2\theta}\exp\{(1/2)\delta \ln\tan\theta\}$,
where $d^2=x^2+y^2$.
The amplitude $\rho$ is equal to
zero on the lines $x=\pm y$. The division 
$1/(x+\delta y) $is possible provided that $\rho\not= 0$.
If $u_1=x_1+\delta y_1, 
u_2=x_2+\delta y_2$ are twocomplex
numbers of amplitudes and arguments $\rho_1,\lambda_1$ and respectively
$\rho_2, \lambda_2$, then the amplitude and the argument $\rho,
\lambda$ of
the product twocomplex number
$u_1u_2=x_1x_2+y_1y_2+\delta (x_1y_2+y_1x_2)$
are $\rho=\rho_1\rho_2,
\lambda=\lambda_1+\lambda_2$. 
Thus, the amplitude $\rho$ is a
multiplicative quantity
and the argument $\lambda$ is an additive quantity upon the
multiplication of twocomplex numbers, which reminds the properties of 
ordinary, two-dimensional complex numbers.

Expressions are given for the elementary functions of twocomplex variable.
Moreover, it is shown that the region of convergence of series of powers of
twocomplex variables is a rectangle  having the sides parallel to the bisectors
$x=\pm y$ . 

A function $f(u)$ of the twocomplex variable 
$u=x+\delta y$ can be defined by a corresponding 
power series. It will be shown that the function $f$ has a
derivative $\lim_{u\rightarrow u_0} [f(u)-f(u_0)]/(u-u_0)$ independent of the
direction of approach of $u$ to $u_0$. If the twocomplex function $f(u)$
of the twocomplex variable $u$ is written in terms of 
the real functions $P(x,y),Q(x,y)$ of real variables $x,y$ as
$f(u)=P(x,y)+\delta Q(x,y)$, then relations of equality 
exist between partial derivatives of the functions $P,Q$, and the 
functions $P,Q$ are solutions of the two-dimensional wave equation.

It will also be shown that the integral $\int_A^B f(u) du$ of a twocomplex
function between two points $A,B$ is independent of the 
path connecting the points $A,B$.

A polynomial $u^n+a_1 u^{n-1}+\cdots+a_{n-1} u +a_n  $ can be
written as a 
product of linear or quadratic factors, although the factorization may not be
unique. 

The twocomplex numbers described in this chapter are a particular case for 
$n=2$ of the polar complex numbers in $n$ dimensions discussed in Sec. 6.1.

\section{Operations with hyperbolic twocomplex numbers}

A hyperbolic complex number in two dimensions is determined by its two
components $(x,y)$. The sum 
of the hyperbolic twocomplex numbers $(x,y)$ and
$(x^\prime,y^\prime)$ is the hyperbolic twocomplex
number $(x+x^\prime,y+y^\prime)$.\index{sum!twocomplex} 
The product of the hyperbolic twocomplex numbers
$(x,y)$ and $(x^\prime,y^\prime)$ 
is defined in this chapter to be the hyperbolic twocomplex
number
$(xx^\prime+yy^\prime,
xy^\prime+yx^\prime)$.\index{product!twocomplex}


Twocomplex numbers and their operations can be represented by  writing the
twocomplex number $(x,y)$ as  
$u=x+\delta y$, where $\delta$ 
is a basis for which the multiplication rules are 
\begin{equation}
1\cdot\delta=\delta,\: \delta^2=1 .
\label{2-1}
\end{equation}\index{complex units!twocomplex}
Two twocomplex numbers $u=x+\delta y, 
u^\prime=x^\prime+\delta y^\prime$ are equal, 
$u=u^\prime$, if and only if $x=x^\prime, y=y^\prime. $
If $u=x+\delta y, 
u^\prime=x^\prime+\delta y^\prime$
are twocomplex numbers, 
the sum $u+u^\prime$ and the 
product $uu^\prime$ defined above can be obtained by applying the usual
algebraic rules to the sum 
$(x+\delta y)+ 
(x^\prime+\delta y^\prime)$
and to the product 
$(x+\delta y)
(x^\prime+\delta y^\prime)$,
and grouping of the resulting terms,
\begin{equation}
u+u^\prime=x+x^\prime+\delta(y+y^\prime),
\label{2-1a}
\end{equation}\index{sum!twocomplex} 
\begin{equation}
uu^\prime=
xx^\prime+yy^\prime+
\delta(xy^\prime+yx^\prime)
\label{2-1b}
\end{equation}\index{product!twocomplex}

If $u,u^\prime,u^{\prime\prime}$ are twocomplex numbers, the multiplication is
associative 
\begin{equation}
(uu^\prime)u^{\prime\prime}=u(u^\prime u^{\prime\prime})
\label{2-2}
\end{equation}
and commutative
\begin{equation}
u u^\prime=u^\prime u ,
\label{2-3}
\end{equation}
as can be checked through direct calculation.
The twocomplex zero is $0+\delta\cdot 0,$ 
denoted simply 0, 
and the twocomplex unity is $1+\delta\cdot 0,$ 
denoted simply 1.

The inverse of the twocomplex number 
$u=x+\delta y$ is a twocomplex number $u^\prime=x^\prime+\delta y^\prime$
having the property that
\begin{equation}
uu^\prime=1 .
\label{2-4}
\end{equation}
Written on components, the condition, Eq. (\ref{2-4}), is
\begin{equation}
\begin{array}{c}
xx^\prime+yy^\prime=1,\\
yx^\prime+xy^\prime=0.\\
\end{array}
\label{2-5}
\end{equation}
The system (\ref{2-5}) has the solution
\begin{equation}
x^\prime=\frac{x}{\nu} ,
\label{2-6a}
\end{equation}
\begin{equation}
y^\prime=
-\frac{y}{\nu} ,
\label{2-6b}
\end{equation}\index{inverse!twocomplex}
provided that $\nu\not= 0$, where
\begin{equation}
\nu=x^2-y^2 .
\label{2-6e}
\end{equation}\index{inverse, determinant!twocomplex}

The quantity $\nu$ can be written as
\begin{equation}
\nu=v_+v_-  ,
\label{2-7}
\end{equation}
where
\begin{equation}
v_+=x+y, \: v_-= x-y.
\label{2-8}
\end{equation}\index{canonical variables!twocomplex}
The variables $v_+, v_-$ will be called canonical hyperbolic twocomplex
variables. 
Then a twocomplex number $u=x+\delta y$ has an inverse,
unless 
\begin{equation}
v_+=0 ,\:\:{\rm or}\:\: v_-=0.
\label{2-9}
\end{equation}

For arbitrary values of the variables $x,y$, the quantity $\nu$ can be
positive or negative. If $\nu\geq 0$, the quantity 
\begin{equation}
\rho=\nu^{1/2},\;\nu>0,
\label{2-9b}
\end{equation}\index{amplitude!twocomplex}
will be called amplitude of the twocomplex number $x+\delta y$.
The normals of the lines in Eq. (\ref{2-9}) are orthogonal to
each other. Because of conditions (\ref{2-9}) these lines will
be also called the nodal lines. \index{nodal lines!twocomplex}
It can be shown that 
if $uu^\prime=0$ then either $u=0$, or $u^\prime=0$, or 
one of the twocomplex numbers $u, u^\prime$ is of the form $x+\delta x$ and the
other is of the form $x-\delta x$.\index{divisors of zero!twocomplex}

\section{Geometric representation of hyperbolic twocomplex numbers}

The twocomplex number $x+\delta y$ can be represented by 
the point $A$ of coordinates $(x,y)$. 
If $O$ is the origin of the two-dimensional space $x,y$, the distance 
from $A$ to the origin $O$ can be taken as
\begin{equation}
d^2=x^2+y^2 .
\label{2-12}
\end{equation}\index{distance!twocomplex}
The distance $d$ will be called modulus of the twocomplex number $x+\delta y$. 

Since
\begin{equation}
(x+y)^2+(x-y)^2=2d^2 ,
\label{2-12b}
\end{equation}\index{modulus, canonical variables!twocomplex}
$x+y$ and $x-y$ can be written as
\begin{equation}
x+y=\sqrt{2}d\sin\theta, \;x-y=\sqrt{2}d\cos\theta,
\label{2-12c}
\end{equation}
so that
\begin{equation}
x=d\sin(\theta+\pi/4),\; y=-d\cos(\theta+\pi/4) .
\label{2-12d}
\end{equation}

If $u=x+\delta y, u_1=x_1+\delta y_1,
u_2=x_2+\delta y_2$, and $u=u_1u_2$, and if
\begin{equation}
v_{j+}=x_j+y_j, \; v_{j-}= x_j-y_j , \;2d_j^2=v_{j+}^2+v_{j-}^2,\;
x_j+y_j=\sqrt{2}d_j\sin\theta_j, \;x_j-y_j=d_j\sqrt{2}\cos\theta_j,
\label{2-13}
\end{equation}
for $j=1,2$,
it can be shown that
\begin{equation}
v_+=v_{1+}v_{2+} ,\:\:
v_-=v_{1-}v_{2-}, \:\:\tan\theta=\tan\theta_1\tan\theta_2.
\label{2-14}
\end{equation}\index{transformation of variables!twocomplex}
The relations (\ref{2-14}) are a consequence of the identities
\begin{eqnarray}
(x_1x_2+y_1y_2)+(x_1y_2+y_1x_2)=(x_1+y_1)(x_2+y_2),
\label{2-15}
\end{eqnarray}
\begin{eqnarray}
(x_1x_2+y_1y_2)-(x_1y_2+y_1x_2)=(x_1-y_1)(x_2-y_2).
\label{2-16}
\end{eqnarray}

A consequence of Eqs. (\ref{2-14}) is that if $u=u_1u_2$, then
\begin{equation}
\nu=\nu_1\nu_2 ,
\label{2-19}
\end{equation}
where
\begin{equation}
\nu_j=v_{j+} v_{j-}, 
\label{2-20}
\end{equation}
for $j=1,2$.
If $\nu>0,\nu_1>0,\nu_2>0$, then
\begin{equation}
\rho=\rho_1\rho_2 ,
\label{2-20b}
\end{equation}\index{transformation of variables!twocomplex}
where
\begin{equation}
\rho_j=\nu_j^{1/2} , 
\label{2-20c}
\end{equation}
for $j=1,2$.

The twocomplex numbers
\begin{equation}
e_+=\frac{1+\delta}{2},\:
e_-=\frac{1-\delta}{2},\:
\label{2-21}\index{canonical base!twocomplex}
\end{equation}
are orthogonal,
\begin{equation}
e_+e_-=0,
\label{2-22}
\end{equation}
and have also the property that
\begin{equation}
e_+^2=e_+, \;e_-^2=e_-.
\label{2-23}
\end{equation}
The ensemble $e_+, e_-$ will be called the canonical hyperbolic twocomplex
base. 
The twocomplex number $u=x+\delta y$ can be written as
\begin{equation}
x+\delta y
=(x+y)e_++(x-y)e_-,  
\label{2-24}\index{canonical form!twocomplex}
\end{equation}
or, by using Eq. (\ref{2-8}),
\begin{equation}
u=v_+e_++v_- e_-,
\label{2-25}
\end{equation}\index{canonical form!twocomplex}
which will be called the canonical form of the hyperbolic twocomplex number.
Thus, if $u_j=v_{j+}e_++v_{j-} e_-, \:j=1,2$, and $u=u_1u_2$, then
the multiplication of the hyperbolic twocomplex numbers is expressed by the
relations (\ref{2-14}).

The relation (\ref{2-19}) for the product of twocomplex numbers can 
be demonstrated also by using a representation of the multiplication of the 
twocomplex numbers by matrices, in which the twocomplex number $u=x+\delta
y$ is represented by the matrix
\begin{equation}
\left(\begin{array}{cc}
x&y\\
y&x\\
\end{array}\right) .
\label{2-26}
\end{equation}\index{matrix representation!twocomplex}
The product $u=x+\delta y$ of the twocomplex numbers
$u_1=x_1+\delta y_1, u_2=x_2+\delta y_2$, can be represented by the matrix
multiplication  
\begin{equation}
\left(\begin{array}{cc}
x&y\\
y&x\\
\end{array}\right) =
\left(\begin{array}{cc}
x_1&y_1\\
y_1&x_1\\
\end{array}\right) 
\left(\begin{array}{cc}
x_2&y_2\\
y_2&x_2\\
\end{array}\right) .
\label{2-27}
\end{equation}
It can be checked that
\begin{equation}
{\rm det}\left(\begin{array}{cccc}
x&y\\
y&x\\
\end{array}\right) =
\nu .
\label{2-28}
\end{equation}
The identity (\ref{2-19}) is then a consequence of the fact the determinant 
of the product of matrices is equal to the product of the determinants 
of the factor matrices. 

\section{Exponential and trigonometric forms of a twocomplex number}

The exponential function of the hypercomplex variable $u$ can be defined by the
series
\begin{equation}
\exp u = 1+u+u^2/2!+u^3/3!+\cdots . 
\label{2-29}
\end{equation}\index{exponential function, hypercomplex!definition}
It can be checked by direct multiplication of the series that
\begin{equation}
\exp(u+u^\prime)=\exp u \cdot \exp u^\prime . 
\label{2-30}
\end{equation}\index{exponential function, hypercomplex!addition theorem}
The series for the exponential function and the addition theorem have the same
form for all systems of hypercomplex numbers discussed in this work.
If $u=x+\delta y$, then  $\exp u$ can be calculated as
$\exp u=\exp x \cdot \exp (\delta y) $. According to Eq. (\ref{2-1}), 
\begin{equation}
\delta^{2m}=1, \delta^{2m+1}=\delta, 
\label{2-31}
\end{equation}\index{complex units, powers of!twocomplex}
where $m$ is a natural number,
so that $\exp (\delta y)$ can be
written as 
\begin{equation}
\exp (\delta y) = \cosh y +\delta \sinh y .
\label{2-32}
\end{equation}\index{exponential, expression!twocomplex}
From Eq. (\ref{2-32}) it can be inferred that
\begin{eqnarray}
(\cosh t +\delta \sinh t)^m=\cosh mt +\delta \sinh mt .
\label{2-33}
\end{eqnarray}

The twocomplex numbers $u=x+\delta y$ for which
$v_+=x+y>0, \: v_-= x-y>0$ can be written in the form 
\begin{equation}
x+\delta y=e^{x_1+\delta y_1} .
\label{2-34}
\end{equation}
The expressions of $x_1, y_1$ as functions of 
$x, y$ can be obtained by
developing $e^{\delta y_1}$ with the aid of
Eq. (\ref{2-32}) and separating the hypercomplex components, 
\begin{equation}
x=e^{x_1}\cosh y_1 ,
\label{2-35}
\end{equation}
\begin{equation}
y=e^{x_1}\sinh y_1 ,
\label{2-36}
\end{equation}
It can be shown from Eqs. (\ref{2-35})-(\ref{2-36}) that
\begin{equation}
x_1=\frac{1}{2} \ln(v_+ v_- ) , \:
y_1=\frac{1}{2}\ln\frac{v_+}{v_- }.\:
\label{2-39}
\end{equation}
The twocomplex number $u$ can thus be written as
\begin{equation}
u=\rho\exp(\delta \lambda),
\label{2-40}
\end{equation}
where the amplitude is $\rho=(x^2-y^2)^{1/2}$ and the argument is
$\lambda=(1/2)\ln\{(x+y)/(x-y)\} $, for $x+y>0, x-y>0$.
The expression (\ref{2-40}) can be written with the aid of the variables $d,
\theta$, Eq. (\ref{2-12c}), as
\begin{equation}
u=\rho\exp\left(\frac{1}{2}\delta \ln\tan\theta\right),
\label{2-40b}
\end{equation}
which is the exponential form of the twocomplex number $u$, where
$0<\theta<\pi/2$. \index{exponential form!twocomplex}

The relation between the amplitude $\rho$ and the distance $d$ is
\begin{equation}
\rho=d\sin^{1/2}2\theta.
\label{2-40c}
\end{equation}
Substituting this form of $\rho$ in Eq. (\ref{2-40b})
yields 
\begin{equation}
u=d\sin^{1/2}2\theta\exp\left(\frac{1}{2}\delta \ln\tan\theta\right),
\label{2-40d}
\end{equation}\index{trigonometric form!expression}
which is the trigonometric form of the twocomplex number $u$.

\section{Elementary functions of a twocomplex variable}

The logarithm $u_1$ of the twocomplex number $u$, $u_1=\ln u$, can be defined
for $v_+>0, v_->0$ as the solution of the equation
\begin{equation}
u=e^{u_1} ,
\label{2-41}
\end{equation}
for $u_1$ as a function of $u$. From Eq. (\ref{2-40b}) it results that 
\begin{equation}
\ln u=\ln\rho+\frac{1}{2}\delta \ln\tan\theta
\label{2-42}
\end{equation}\index{logarithm!twocomplex}
It can be inferred from Eqs. (\ref{2-42}) and (\ref{2-14}) that
\begin{equation}
\ln(u_1u_2)=\ln u_1+\ln u_2 .
\label{2-43}
\end{equation}
The explicit form of Eq. (\ref{2-42}) is 
\begin{eqnarray}
\ln (x+\delta y)=
\frac{1}{2}(1+\delta)\ln(x+y)
+\frac{1}{2}(1-\delta)\ln(x-y),
\label{2-45}
\end{eqnarray}
so that the relation (\ref{2-42}) can be written with the aid of 
Eq. (\ref{2-21}) as 
\begin{equation}
\ln u = e_+\ln v_+ + e_-\ln v_-.
\label{2-44}
\end{equation}

The power function $u^n$ can be defined for $v_+>0, v_->0$ and real values
of $n$ as 
\begin{equation}
u^n=e^{n\ln u} .
\label{2-46}
\end{equation}
It can be inferred from Eqs. (\ref{2-46}) and (\ref{2-43}) that
\begin{equation}
(u_1u_2)^n=u_1^n\:u_2^n .
\label{2-47}
\end{equation}
Using the expression (\ref{2-44}) for $\ln u$ and the relations (\ref{2-22}) and
(\ref{2-23}) it can be shown that
\begin{eqnarray}
(x+\delta y)^n=
\frac{1}{2}(1+\delta)(x+y)^n
+\frac{1}{2}(1-\delta)(x-y)^n.
\label{2-48}
\end{eqnarray}\index{power function!twocomplex}
For integer $n$, the relation (\ref{2-48}) is valid for any $x,y$. The
relation (\ref{2-48}) for $n=-1$ is
\begin{equation}
\frac{1}{x+\delta y}=
\frac{1}{2}\left(\frac{1+\delta}{x+y}
+\frac{1-\delta}{x-y}
\right) .
\label{2-49}
\end{equation}

The trigonometric functions $\cos u$ and $\sin u $ of the hypercomplex variable
$u$ are defined by the series
\begin{equation}
\cos u = 1 - u^2/2!+u^4/4!+\cdots, 
\label{2-50}
\end{equation}
\begin{equation}
\sin u=u-u^3/3!+u^5/5! +\cdots .
\label{2-51}
\end{equation}\index{trigonometric functions, hypercomplex!definitions}
It can be checked by series multiplication that the usual addition theorems
hold for the hypercomplex numbers $u_1, u_2$,
\begin{equation}
\cos(u_1+u_2)=\cos u_1\cos u_2 - \sin u_1 \sin u_2 ,
\label{2-52}
\end{equation}
\begin{equation}
\sin(u_1+u_2)=\sin u_1\cos u_2 + \cos u_1 \sin u_2 .
\label{2-53}
\end{equation}\index{trigonometric functions, hypercomplex!addition theorems}
The series for the trigonometric functions and the addition theorems have the
same form for all systems of hypercomplex numbers discussed in this work.
The cosine and sine functions of the hypercomplex variables $\delta y$ can be
expressed as 
\begin{equation}
\cos\delta y=\cos y, \: \sin\delta y=\delta\sin y.
\label{2-54}
\end{equation}\index{trigonometric functions, expressions!twocomplex}
The cosine and sine functions of a twocomplex number $x+\delta y$ can then be
expressed in terms of elementary functions with the aid of the addition
theorems Eqs. (\ref{2-52}), (\ref{2-53}) and of the expressions in  Eq. 
(\ref{2-54}). 

The hyperbolic functions $\cosh u$ and $\sinh u $ of the hypercomplex variable
$u$ are defined by the series
\begin{equation}
\cosh u = 1 + u^2/2!+u^4/4!+\cdots, 
\label{2-57}
\end{equation}
\begin{equation}
\sinh u=u+u^3/3!+u^5/5! +\cdots .
\label{2-58}
\end{equation}\index{hyperbolic functions, hypercomplex!definitions}
It can be checked by series multiplication that the usual addition theorems
hold for the hypercomplex numbers $u_1, u_2$,
\begin{equation}
\cosh(u_1+u_2)=\cosh u_1\cosh u_2 + \sinh u_1 \sinh u_2 ,
\label{2-59}
\end{equation}
\begin{equation}
\sinh(u_1+u_2)=\sinh u_1\cosh u_2 + \cosh u_1 \sinh u_2 .
\label{2-60}
\end{equation}\index{hyperbolic functions, hypercomplex!addition theorems}
The series for the hyperbolic functions and the addition theorems have the
same form for all systems of hypercomplex numbers discussed in this work.
The $\cosh$ and $\sinh$ functions of the hypercomplex variable $\delta y$ can
be expressed as 
\begin{equation}
\cosh\delta y=\cosh y, \: \sinh\delta y=\delta\sinh y.
\label{2-61}
\end{equation}\index{hyperbolic functions, expressions!twocomplex}
The hyperbolic cosine and sine functions of a twocomplex number $x+\delta y$
can then be 
expressed in terms of elementary functions with the aid of the addition
theorems Eqs. (\ref{2-59}), (\ref{2-60}) and of the expressions in  Eq. 
(\ref{2-61}). 

\section{Twocomplex power series}

A twocomplex series is an infinite sum of the form
\begin{equation}
a_0+a_1+a_2+\cdots+a_n+\cdots , 
\label{2-64}
\end{equation}\index{series!twocomplex}
where the coefficients $a_n$ are twocomplex numbers. The convergence of 
the series (\ref{2-64}) can be defined in terms of the convergence of its 2 real
components. The convergence of a twocomplex series can however be studied
using twocomplex variables. The main criterion for absolute convergence 
remains the comparison theorem, but this requires a number of inequalities
which will be discussed further.

The modulus of a twocomplex number $u=x+\delta y$ can be defined as 
\begin{equation}
|u|=(x^2+y^2)^{1/2} ,
\label{2-65}
\end{equation}\index{modulus, definition!twocomplex}
so that according to Eq. (\ref{2-12}) $d=|u|$. Since $|x|\leq |u|, |y|\leq |u|$,
a property of  
absolute convergence established via a comparison theorem based on the modulus
of the series (\ref{2-64}) will ensure the absolute convergence of each real
component of that series.

The modulus of the sum $u_1+u_2$ of the twocomplex numbers $u_1, u_2$ fulfils
the inequality
\begin{equation}
||u_1|-|u_2||\leq |u_1+u_2|\leq |u_1|+|u_2| .
\label{2-66}
\end{equation}\index{modulus, inequalities!twocomplex}
For the product the relation is 
\begin{equation}
|u_1u_2|\leq \sqrt{2}|u_1||u_2| ,
\label{2-67}
\end{equation}
which replaces the relation of equality extant for regular complex numbers.
The equality in Eq. (\ref{2-67}) takes place for $x_1=y_1,x_2=y_2$ or
$x_1=-y_1,x_2=-y_2$.
In particular
\begin{equation}
|u^2|\leq \sqrt{2}|u|^2 .
\label{2-68}
\end{equation}
The inequality in Eq. (\ref{2-67}) implies that
\begin{equation}
|u^m|\leq 2^{(m-1)/2}|u|^m .
\label{2-69}
\end{equation}
From Eqs. (\ref{2-67}) and (\ref{2-69}) it results that
\begin{equation}
|au^m|\leq 2^{m/2} |a| |u|^m .
\label{2-70}
\end{equation}

A power series of the twocomplex variable $u$ is a series of the form
\begin{equation}
a_0+a_1 u + a_2 u^2+\cdots +a_l u^l+\cdots .
\label{2-71}
\end{equation}\index{power series!twocomplex}
Since
\begin{equation}
\left|\sum_{l=0}^\infty a_l u^l\right| \leq  \sum_{l=0}^\infty
2^{l/2}|a_l| |u|^l ,
\label{2-72}
\end{equation}
a sufficient condition for the absolute convergence of this series is that
\begin{equation}
\lim_{l\rightarrow \infty}\frac{\sqrt{2}|a_{l+1}||u|}{|a_l|}<1 .
\label{2-73}
\end{equation}\index{convergence of power series!twocomplex}
Thus the series is absolutely convergent for 
\begin{equation}
|u|<c_0,
\label{2-74}
\end{equation}
where 
\begin{equation}
c_0=\lim_{l\rightarrow\infty} \frac{|a_l|}{\sqrt{2}|a_{l+1}|} ,
\label{2-75}
\end{equation}

The convergence of the series (\ref{2-71}) can be also studied with the aid of
the formula (\ref{2-48}) which, for integer values of $l$, is valid for any $x,
y, z, t$. If $a_l=a_{lx}+\delta a_{ly}$, and
\begin{eqnarray}
A_{l+}=a_{lx}+a_{ly}, A_{l-}= a_{lx}-a_{ly} , 
\label{2-76}
\end{eqnarray}
it can be shown with the aid of relations (\ref{2-22}) and (\ref{2-23}) that
\begin{equation}
a_l e_+=A_{l+} e_+, \: a_l e_-=A_{l-} e_-, \: 
\label{2-77}
\end{equation}
so that the expression of the series (\ref{2-71}) becomes
\begin{equation}
\sum_{l=0}^\infty \left(A_{l+}  v_+^l e_++
A_{l-}  v_-^le_-\right) ,
\label{2-78}
\end{equation}
where the quantities $v_+, v_-$
have been defined in Eq. (\ref{2-8}).
The sufficient conditions for the absolute convergence of the series 
in Eq. (\ref{2-78}) are that
\begin{equation}
\lim_{l\rightarrow \infty}\frac{|A_{l+1,+}||v_+|}{|A_{l+}|}<1,
\lim_{l\rightarrow \infty}\frac{|A_{l+1,-}||v_-|}{|A_{l-}|}<1.
\label{2-79}
\end{equation}
Thus the series in Eq. (\ref{2-78}) is absolutely convergent for 
\begin{equation}
|x+y|<c_+,\:
|x-y|<c_-,
\label{2-80}
\end{equation}\index{convergence, region of!twocomplex}
where 
\begin{equation}
c_+=\lim_{l\rightarrow\infty} \frac{|A_{l+}|}{|A_{l+1,+}|} ,\:
c_-=\lim_{l\rightarrow\infty} \frac{|A_{l-}|}{|A_{l+1,-}|} .
\label{2-81}
\end{equation}
The relations (\ref{2-80}) show that the region of convergence of the series
(\ref{2-78}) is a rectangle having the sides parallel to the bisectors $x=\pm y$.
It can be shown that $c_0=(1/\sqrt{2}){\rm min}(c,c^\prime)$, where
${\rm min}(c,c^\prime)$ designates the smallest of the numbers $c, c^\prime$.
Since $|u|^2=(v_+^2+v_-^2)/2$, it can be seen that the circular region of
convergence defined in Eqs. (\ref{2-74}), (\ref{2-75})
is included in the parallelogram defined in Eqs. (\ref{2-80}) and (\ref{2-81}).

\section{Analytic functions of twocomplex variables}

The derivative  
of a function $f(u)$ of the hypercomplex variables $u$ is
defined as a function $f^\prime (u)$ having the property that
\begin{equation}
|f(u)-f(u_0)-f^\prime (u_0)(u-u_0)|\rightarrow 0 \:\:{\rm as} 
\:\:|u-u_0|\rightarrow 0 . 
\label{2-gs88}
\end{equation}\index{derivative, hypercomplex!definition}
If the difference $u-u_0$ is not parallel to one of the nodal hypersurfaces,
the definition in Eq. (\ref{2-gs88}) can also 
be written as
\begin{equation}
f^\prime (u_0)=\lim_{u\rightarrow u_0}\frac{f(u)-f(u_0)}{u-u_0} .
\label{2-gs89}
\end{equation}
The derivative of the function $f(u)=u^m $, with $m$ an integer, 
is $f^\prime (u)=mu^{m-1}$, as can be seen by developing $u^m=[u_0+(u-u_0)]^m$
as
\begin{equation}
u^m=\sum_{p=0}^{m}\frac{m!}{p!(m-p)!}u_0^{m-p}(u-u_0)^p,
\label{2-gs90}
\end{equation}\index{derivative, hypercomplex!of power function}
and using the definition (\ref{2-gs88}).

If the function $f^\prime (u)$ defined in Eq. (\ref{2-gs88}) is independent of the
direction in space along which $u$ is approaching $u_0$, the function $f(u)$ 
is said to be analytic, analogously to the case of functions of regular complex
variables. \cite{3} \index{analytic function, hypercomplex!definition}
The function $u^m$, with $m$ an integer, 
of the hypercomplex variable $u$ is analytic, because the
difference $u^m-u_0^m$ is always proportional to $u-u_0$, as can be seen from
Eq. (\ref{2-gs90}). Then series of
integer powers of $u$ will also be analytic functions of the hypercomplex
variable $u$, and this result holds in fact for any commutative algebra. 
\index{analytic function, hypercomplex!power series}

If an analytic function is defined by a series around a certain point, for
example $u=0$, as
\begin{equation}
f(u)=\sum_{k=0}^\infty a_k u^k ,
\label{2-gs91a}
\end{equation}\index{analytic function, hypercomplex!expansion in series}
an expansion of $f(u)$ around a different point $u_0$,
\begin{equation}
f(u)=\sum_{k=0}^\infty c_k (u-u_0)^k ,
\label{2-gs91aa}
\end{equation}
can be obtained by
substituting in Eq. (\ref{2-gs91a}) the expression of $u^k$ according to Eq.
(\ref{2-gs90}). Assuming that the series are absolutely convergent so that the
order of the terms can be modified and ordering the terms in the resulting
expression according to the increasing powers of $u-u_0$ yields
\begin{equation}
f(u)=\sum_{k,l=0}^\infty \frac{(k+l)!}{k!l!}a_{k+l} u_0^l (u-u_0)^k .
\label{2-gs91b}
\end{equation}
Since the derivative of order $k$ at $u=u_0$ of the function $f(u)$ , Eq.
(\ref{2-gs91a}), is 
\begin{equation}
f^{(k)}(u_0)=\sum_{l=0}^\infty \frac{(k+l)!}{l!}a_{k+l} u_0^l ,
\label{2-gs91c}
\end{equation}
the expansion of $f(u)$ around $u=u_0$, Eq. (\ref{2-gs91b}), becomes
\begin{equation}
f(u)=\sum_{k=0}^\infty \frac{1}{k!} f^{(k)}(u_0)(u-u_0)^k ,
\label{2-gs91d}
\end{equation}
which has the same form as the series expansion of 2-dimensional complex
functions. 
The relation (\ref{2-gs91d}) shows that the coefficients in the series expansion,
Eq. (\ref{2-gs91aa}), are
\begin{equation}
c_k=\frac{1}{k!}f^{(k)}(u_0) .
\label{2-gs92}
\end{equation}

The rules for obtaining the derivatives and the integrals of the basic
functions can 
be obtained from the series of definitions and, as long as these series
expansions have the same form as the corresponding series for the
2-dimensional complex functions, the rules of derivation and integration remain
unchanged. \index{rules for derivation and integration}
The relations (\ref{2-gs88})-(\ref{2-gs92}) have the same form for all systems of
hypercomplex numbers discussed in this work.

If the twocomplex function $f(u)$
of the twocomplex variable $u$ is written in terms of 
the real functions $P(x,y),Q(x,y)$ of real
variables $x,y$ as 
\begin{equation}
f(u)=P(x,y)+\delta Q(x,y), 
\label{2-87}
\end{equation}\index{functions, real components!twocomplex}
then relations of equality 
exist between partial derivatives of the functions $P,Q$. These relations
can be obtained by writing the derivative of the function $f$ as
\begin{eqnarray}
\lim_{\Delta x,\Delta y\rightarrow 0}\frac{1}{\Delta x+\delta \Delta y } 
\left[\frac{\partial P}{\partial x}\Delta x+
\frac{\partial P}{\partial y}\Delta y
+\delta\left(\frac{\partial Q}{\partial x}\Delta x+
\frac{\partial Q}{\partial y}\Delta y\right) 
\right] ,
\label{2-88}
\end{eqnarray}\index{derivative, independence of direction!twocomplex}
where the difference $u-u_0$ in Eq. (\ref{2-gs89}) is 
$u-u_0=\Delta x+\delta\Delta y$. 
The relations between the partials derivatives of the functions $P, Q$ are
obtained by setting successively in   
Eq. (\ref{2-88}) $\Delta x\rightarrow 0, \Delta y=0$;
then $\Delta x= 0, \Delta y\rightarrow 0$. The
relations are 
\begin{equation}
\frac{\partial P}{\partial x} = \frac{\partial Q}{\partial y} ,
\label{2-89}
\end{equation}
\begin{equation}
\frac{\partial Q}{\partial x} = \frac{\partial P}{\partial y} .
\label{2-90}
\end{equation}{relations between partial derivatives!twocomplex}

The relations (\ref{2-89})-(\ref{2-90}) are analogous to the Riemann relations
for the real and imaginary components of a complex function. It can be shown
from Eqs. (\ref{2-89})-(\ref{2-90}) that the components $P, Q$ are solutions
of the equations 
\begin{equation}
\frac{\partial^2 P}{\partial x^2}-\frac{\partial^2 P}{\partial y^2}=0,
\label{2-93}
\end{equation}
\begin{equation}
\frac{\partial^2 Q}{\partial x^2}-\frac{\partial^2 Q}{\partial y^2}=0,
\label{2-94}
\end{equation}\index{relations between second-order derivatives!twocomplex}
As can be seen from Eqs. (\ref{2-93})-(\ref{2-94}), the components $P, Q$ of
an analytic function of twocomplex variable are solutions of the wave 
equation with respect to the variables $x,y$.

\section{Integrals of twocomplex functions}

The singularities of twocomplex functions arise from terms of the form
$1/(u-u_0)^m$, with $m>0$. Functions containing such terms are singular not
only at $u=u_0$, but also at all points of the lines
passing through $u_0$ and which are parallel to the nodal lines. 

The integral of a twocomplex function between two points $A, B$ along a path
situated in a region free of singularities is independent of path, which means
that the integral of an analytic function along a loop situated in a region
free from singularities is zero,
\begin{equation}
\oint_\Gamma f(u) du = 0.
\label{2-105}
\end{equation}
Using the expression, Eq. (\ref{2-87})
for $f(u)$ and the fact that $du=dx+\delta  dy$, the
explicit form of the integral in Eq. (\ref{2-105}) is
\begin{eqnarray}
\oint _\Gamma f(u) du = \oint_\Gamma
[(Pdx+Qdy)+\delta(Qdx+Pdy)]
\label{2-106}
\end{eqnarray}\index{integrals, path!twocomplex}
If the functions $P, Q$ are regular on the surface $\Sigma$
enclosed by the loop $\Gamma$,
the integral along the loop $\Gamma$ can be transformed with the aid of the
theorem of Stokes in an integral over the surface $\Sigma$ of terms of the form
$\partial P/\partial y -  \partial Q/\partial x$ and
$\partial P/\partial x -  \partial Q/\partial y$ 
which are equal to zero by Eqs. (\ref{2-89})-(\ref{2-90}), and this proves Eq.
(\ref{2-105}). 

The exponential form of the twocomplex numbers, Eq. (\ref{2-40}), contains no
cyclic variable, and therefore the concept of residue is not applicable to the
twocomplex numbers defined in Eqs. (\ref{2-1}).

\section{Factorization of twocomplex polynomials}

A polynomial of degree $m$ of the twocomplex variable 
$u=x+\delta y$ has the form
\begin{equation}
P_m(u)=u^m+a_1 u^{m-1}+\cdots+a_{m-1} u +a_m ,
\label{2-106b}
\end{equation}
where the constants are in general twocomplex numbers.
If $a_m=a_{mx}+\delta a_{my}$, and with the
notations of Eqs. (\ref{2-8}) and (\ref{2-76}) applied for $0, 1, \cdots, m$ , the
polynomial $P_m(u)$ can be written as 
\begin{eqnarray}
\lefteqn{P_m= \left[v_+^m 
+A_{1+} v_+^{m-1}+\cdots+A_{m-1,+} v_++ A_{m+} \right] e_+\nonumber}\\
&&+\left[v_-^m 
+A_{1-} v_-^{m-1} +\cdots+A_{m-1,-} v_-+ A_{m-}
\right]e_-. 
\label{2-107}
\end{eqnarray}\index{polynomial, canonical variables!twocomplex}
Each of the polynomials of degree $m$ with real coefficients in Eq. (\ref{2-107})
can be written as a product
of linear or quadratic factors with real coefficients, or as a product of
linear factors which, if imaginary, appear always in complex conjugate pairs.
Using the latter form for the simplicity of notations, the polynomial $P_m$
can be written as
\begin{equation}
P_m=\prod_{l=1}^m (v_+-v_{l+})e_+
+\prod_{l=1}^m (v_--v_{l-})e_-,
\label{2-108}
\end{equation}
where the quantities $v_{l+}$ appear always in complex conjugate pairs, and the
same is true for the quantities $v_{l-}$.
Due to the properties in Eqs. (\ref{2-22}) and (\ref{2-23}),
the polynomial $P_m(u)$ can be written as a product of factors of
the form  
\begin{equation}
P_m(u)=\prod_{l=1}^m \left[(v_+-v_{l+})e_+
+(v_--v_{l-})e_-
\right].
\label{2-109}
\end{equation}\index{polynomial, factorization!twocomplex}
This relation can be written with the aid of Eqs. (\ref{2-25}) as
\begin{eqnarray}
P_m(u)=\prod_{l=1}^m (u-u_l),
\label{2-110}
\end{eqnarray}
where
\begin{eqnarray}
u_l=e_+ v_{l+} +e_- v_{l-}, 
\label{2-111}
\end{eqnarray}
for $l=1,...,m$.
The roots $v_{l+}$ and the roots $v_{l-}$ 
defined in Eq. (\ref{2-108}) may be ordered arbitrarily, which
means that Eq. (\ref{2-111}) gives sets of $m$ roots
$u_1,...,u_m$ of the polynomial $P_m(u)$, 
corresponding to the various ways in which the roots $v_{l+}, v_{l-}$
are ordered according to $l$ in each group. Thus, while the hypercomplex
components in Eq. (\ref{2-108}) taken 
separately have unique factorizations, the polynomial $P_m(u)$ can be written
in many different ways as a product of linear factors. 

If $P(u)=u^2-1$, the degree is $m=2$, the coefficients of the polynomial are
$a_1=0, a_2=-1$, the twocomplex components of $a_2$ are $a_{2x}=-1, a_{2y}=0$,
the components $A_{2+}, A_{2-}$ are $A_{2+}=-1, A_{2-}=-1$.
The expression, Eq. (\ref{2-107}), of $P(u)$ 
is $P(u)=e_+(v_+^2-1)+e_- (v_-^2-1)$, and
the factorization in Eq. (\ref{2-110}) is $u^2-1=(u-u_1)(u-u_2)$, where 
$u_1=\pm e_+\pm e_-, u_2=-u_1$. The factorizations are thus
$u^2-1=(u+1)(u-1)$ and $u^2-1=(u+\delta)(u-\delta)$. 
It can be checked that 
$(\pm e_+\pm e_-)^2=e_++e_-=1$.

\section{Representation of hyperbolic twocomplex complex numbers by irreducible
matrices} 

If the matrix in Eq. (\ref{2-26}) representing the twocomplex number $u$ is
called $U$, and 
\begin{equation}
T=\left(
\begin{array}{cc}
\frac{1}{\sqrt{2}}   &  \frac{1}{\sqrt{2}}   \\
-\frac{1}{\sqrt{2}}  &  \frac{1}{\sqrt{2}}   \\
\end{array}\right),
\label{2-112}
\end{equation}
it can be checked that 
\begin{equation}
TUT^{-1}=\left(
\begin{array}{cc}
x+y   & 0     \\
0     & x-y   \\
\end{array}\right).
\label{2-113}
\end{equation}\index{representation by irreducible matrices!twocomplex}
The relations for the variables $v_+=x+y, v_-=x-y$ for the
multiplication of twocomplex numbers have been written in Eq. (\ref{2-14}). The
matrix
$T U T^{-1}$  provides an irreducible representation
\cite{4} of the twocomplex numbers $u=x+\delta y$, in terms of matrices with
real coefficients.

\chapter{Complex Numbers in Three Dimensions}

A system of hypercomplex numbers in three dimensions is described in this chapter,
for which the multiplication is associative and commutative, which have
exponential and trigonometric forms, and for which the concepts
of analytic tricomplex 
function,  contour integration and residue can be defined.
The tricomplex numbers introduced in this chapter have 
the form $u=x+hy+kz$, the variables $x, y$ and $z$ being real 
numbers.  
The multiplication rules for the complex units $h, k$ are
$h^2=k, \:  k^2=h,\:  hk=1$. In a geometric representation, the
tricomplex number $u$ is represented by the point $P$ of coordinates $(x,y,z)$.
If $O$ is the origin of the $x,y,z$
axes, $(t)$ the trisector line $x=y=z$ of the positive octant and 
$\Pi$ the plane $x+y+z=0$ passing through the origin $(O)$ and
perpendicular to $(t)$, then the tricomplex number $u$ can be described by
the projection $s$ of the segment $OP$ along the line $(t)$, by the distance
$D$ from $P$ to the line $(t)$, and by the azimuthal angle $\phi$ of the
projection of $P$ on the plane $\Pi$, measured from an angular origin defined
by the intersection of the plane determined by the line $(t)$ and the x axis,
with the plane $\Pi$. 
The amplitude $\rho$ of a tricomplex number is defined as 
$\rho=(x^3+y^3+z^3-3xyz)^{1/3}$, the polar 
angle $\theta$ of $OP$ with respect to the trisector line $(t)$ is given by
$\tan\theta=D/s$, and $d^2=x^2+y^2+z^2$.  
The amplitude $\rho$ is equal to
zero on the trisector line $(t)$ and on the plane $\Pi$. The division 
$1/(x+hy+kz) $is possible provided that $\rho\not= 0$.
The product of two tricomplex numbers is equal to zero if both numbers are
equal to zero, or if one of
the tricomplex numbers lies in the $\Pi$ plane and the other on the $(t)$
line. 

If $u_1=x_1+hy_1+kz_1, 
u_2=x_2+hy_2+kz_2$ are tricomplex
numbers of amplitudes and angles $\rho_1,\theta_1,\phi_1$ and respectively
$\rho_2, \theta_2, \phi_2$, then the amplitude and the angles $\rho,
\theta, \phi$ for
the product tricomplex number
$u_1u_2=x_1x_2+y_1z_2+y_2z_1+h(z_1z_2+x_1y_2+y_1x_2)+k(y_1y_2+x_1z_2+z_1x_2)$
are $\rho=\rho_1\rho_2,
\tan\theta=\tan\theta_1\tan\theta_2/\sqrt{2}, 
\phi=\phi_1+\phi_2$. 
Thus, the amplitude $\rho$ and  $(\tan\theta)/\sqrt{2}$ are
multiplicative quantities  
and the angle $\phi$ is an additive quantity upon the
multiplication of tricomplex numbers, which reminds the properties of 
ordinary, two-dimensional complex numbers.

For the description of the exponential function of a tricomplex variable, it
is useful to define the cosexponential functions 
${\rm cx}(\xi)=1+\xi^3/3!+\xi^6/6!\cdots, 
{\rm mx}(\xi)=\xi+\xi^4/4!+\xi^7/7!\cdots, 
{\rm px}(\xi)=\xi^2/2+\xi^5/5!+\xi^8/8!\cdots $, where p and m stand for plus
and respectively minus, as a reference to the sign of a phase shift in the
expressions of these functions. These functions fulfil the relation ${\rm cx}^3
\xi +{\rm px}^3 \xi+{\rm mx}^3 \xi -3 {\rm cx} \xi \:
{\rm px} \xi\: {\rm mx} \xi =1$. 

The exponential form of a tricomplex number is
$u=\rho$
$\exp\left[(1/3)(h+k)\ln(\sqrt{2}/\tan\theta)\right.$ 
$\left.+(1/3)(h-k)\phi\right]$,
and the trigonometric form of the tricomplex number is
$u=d\sqrt{3/2}$
$\left\{(1/3)(2-h-k)\sin\theta+\right.$
$\left.(1/3)(1+h+k)\sqrt{2}\cos\theta\right\}$
$\exp\left\{(h-k)\phi/\sqrt{3}\right\}$. 

Expressions are given for the elementary functions of tricomplex variable.
Moreover, it is shown that the region of convergence of series of powers of
tricomplex variables are cylinders with the axis parallel to the trisector
line. 
A function $f(u)$ of the tricomplex variable 
$u=x+hy+kz$ can be defined by a corresponding 
power series. It will be shown that the function $f(u)$ has a
derivative at $u_0$ independent of the
direction of approach of $u$ to $u_0$. If the tricomplex function $f(u)$
of the tricomplex variable $u$ is written in terms of 
the real functions $F(x,y,z),G(x,y,z),H(x,y,z)$ of real variables $x,y,z$ as
$f(u)=F(x,y,z)+hG(x,y,z)+kH(x,y,z)$, then relations of equality 
exist between partial derivatives of the functions $F,G,H$, and the differences
$F-G, F-H, G-H$ are solutions of the equation of Laplace.

It will be shown that the integral $\int_A^B f(u) du$ of a regular tricomplex
function between two points $A,B$ is independent of the three-dimensional
path connecting the points $A,B$.
If $f(u)$ is an analytic tricomplex function, then $\oint_\Gamma f(u)
du/(u-u_0) = 
2\pi ( h-k) f(u_0)$ if the integration loop is threaded by the parallel through
$u_0$ to the line $(t)$.

A tricomplex polynomial  $u^m+a_1 u^{m-1}+\cdots+a_{m-1} u +a_m $ can be
written as a 
product of linear or quadratic factors, although the factorization may not be
unique. 

The tricomplex numbers described in this chapter are a particular case for 
$n=3$ of the polar complex numbers in $n$ dimensions discussed in Sec. 6.1.

\section{Operations with tricomplex numbers}

A tricomplex number is determined by its three components $(x,y,z)$. The sum
of the tricomplex numbers $(x,y,z)$ and $(x^\prime,y^\prime,z^\prime)$ is the tricomplex
number $(x+x^\prime,y+y^\prime,z+z^\prime)$. \index{sum!tricomplex}
The product of the tricomplex numbers
$(x,y,z)$ and $(x^\prime,y^\prime,z^\prime)$ is defined in this chapter to be the tricomplex
number $(xx^\prime+yz^\prime+zy^\prime,zz^\prime+xy^\prime+yx^\prime,yy^\prime+xz^\prime+zx^\prime)$.
\index{product!tricomplex}

Tricomplex numbers and their operations can be represented by writing the
tricomplex number $(x,y,z)$ as  
$u=x+hy+kz$, where $h$ and $k$ are bases for which the
multiplication rules are 
\begin{equation}
 h^2=k, \:  k^2=h, \: 1\cdot h=h,\:  1\cdot k =k, \:  hk=1. 
\label{3-1}
\end{equation}\index{complex units!tricomplex}
Two tricomplex numbers $u=x+hy+kz, u^\prime=x^\prime+hy^\prime+
kz^\prime$ are equal, $u=u^\prime$, if and only if $x=x^\prime, y=y^\prime,
z=z^\prime$. 
If $u=x+hy+kz, 
u^\prime=x^\prime+hy^\prime+kz^\prime$ are tricomplex
numbers, the sum $u+u^\prime$ and the 
product $uu^\prime$ defined above can be obtained by applying the usual
algebraic rules to the sum 
$(x+hy+kz)+(x^\prime+hy^\prime+kz^\prime)$
and to the product 
$(x+hy+kz)(x^\prime+hy^\prime+kz^\prime)$,
and grouping of the resulting terms,
\begin{equation}
u+u^\prime=x+x^\prime+h(y+y^\prime)+k(z+z^\prime),
\label{3-1a}
\end{equation}
\begin{equation}\index{sum!tricomplex}
uu^\prime=xx^\prime+yz^\prime+zy^\prime+h(zz^\prime+xy^\prime+yx^\prime)+k(yy^\prime+xz^\prime+zx^\prime) .
\label{3-1b}
\end{equation}\index{product!tricomplex}
If $u,u^\prime,u^{\prime\prime}$ are tricomplex numbers, the multiplication is associative
\begin{equation}
(uu^\prime)u^{\prime\prime}=u(u^\prime u^{\prime\prime})
\label{3-2}
\end{equation}
and commutative
\begin{equation}
u u^\prime=u^\prime u ,
\label{3-3}
\end{equation}
as can be checked through direct calculation.
The tricomplex zero is $0+h\cdot 0+k\cdot 0,$ denoted simply 0, 
and the tricomplex unity is $1+h\cdot 0+k\cdot 0,$ 
denoted simply 1.

The inverse of the tricomplex number $u=x+hy+kz$ is a
tricomplex number $u^\prime=x^\prime+y^\prime+z^\prime$ having the property that
\begin{equation}
uu^\prime=1 .
\label{3-4}
\end{equation}
Written on components, the condition, Eq. (\ref{3-4}), is
\begin{equation}
\begin{array}{c}
xx^\prime+zy^\prime+yz^\prime=1 ,\\
yx^\prime+xy^\prime+zz^\prime=0 ,\\
zx^\prime+yy^\prime+xz^\prime=0 .
\end{array}
\label{3-5}
\end{equation}
The system (\ref{3-5}) has the solution
\begin{equation}
x^\prime=\frac{x^2-yz}{x^3+y^3+z^3-3xyz} ,
\label{3-6a}
\end{equation}
\begin{equation}
y^\prime=\frac{z^2-xy}{x^3+y^3+z^3-3xyz} ,
\label{3-6b}
\end{equation}
\begin{equation}
z^\prime=\frac{y^2-xz}{x^3+y^3+z^3-3xyz} ,
\label{3-6c}
\end{equation}\index{inverse!tricomplex}
provided that $x^3+y^3+z^3-3xyz\not=0$. Since
\begin{equation}
x^3+y^3+z^3-3xyz=(x+y+z)(x^2+y^2+z^2-xy-xz-yz) ,
\label{3-7}
\end{equation}\index{inverse, determinant!tricomplex}
a tricomplex number $x+hy+kz$ has an inverse, unless
\begin{equation}
x+y+z=0  
\label{3-8}
\end{equation}
or
\begin{equation}
x^2+y^2+z^2-xy-xz-yz=0 . 
\label{3-9}
\end{equation}

\begin{figure}
\begin{center}
\epsfig{file=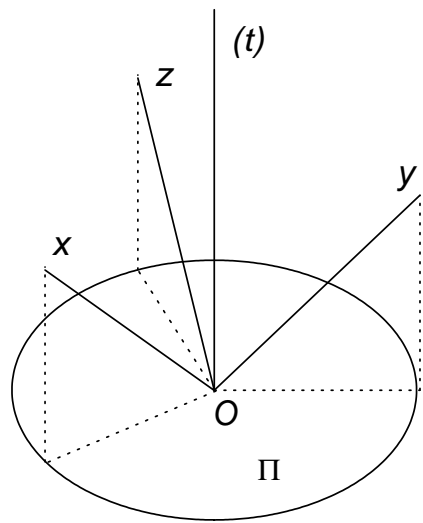,width=12cm}
\caption{Nodal plane $\Pi$, of equation $x+y+z=0$, and trisector line $(t)$, of
equation $x=y=z$, both passing through the origin $O$ of the rectangular axes
$x, y, z$.}
\label{fig1}
\end{center}
\end{figure}

The relation in Eq. (\ref{3-8}) represents the plane $\Pi$ perpendicular
to the trisector line $(t)$ of the $x,y,z$ axes, and passing through 
the origin $O$ of the axes. The plane $\Pi$, shown in Fig. \ref{fig1}, intersects
the $xOy$ plane along the line $z=0, x+y=0$, it intersect the 
$yOz$ plane along the line $x=0, y+z=0$, and it intersects the
$xOz$ plane along the line $y=0, x+z=0$.
The condition (\ref{3-9}) is equivalent to
$(x-y)^2+(x-z)^2+(y-z)^2=0$, which for real $x,y,z$ means that 
$x=y=z$, which represents the trisector line $(t)$ of the axes $x,y,z$. 
The trisector line $(t)$ is perpendicular to the plane $\Pi$. Because of
conditions (\ref{3-8}) and (\ref{3-9}), the trisector line $(t)$ and the plane
$\Pi$ will be also called nodal line and respectively nodal plane.
\index{nodal line, tricomplex}\index{nodal plane, tricomplex}

It can be shown that 
if $uu^\prime=0$ then either $u=0$, or $u^\prime=0$, or one of the tricomplex
numbers $u, u^\prime$ belongs to the trisector line $(t)$ and the other belongs to
the nodal plane $\Pi$.\index{divisors of zero!tricomplex}

\section{Geometric representation of tricomplex numbers}

The tricomplex number $x+hy+kz$ can be represented by the point
$P$ of coordinates $(x,y,z)$. If $O$ is the origin of the axes, 
then the projection $s=OQ$ of the line $OP$ on the trisector line $x=y=z$,
\index{trisector line}
which has the unit tangent $(1/\sqrt{3},1/\sqrt{3},1/\sqrt{3})$, is
\begin{equation}
s=\frac{1}{\sqrt{3}}(x+y+z) . 
\label{3-10}
\end{equation}
The distance $D=PQ$ from $P$ to the trisector line $x=y=z$, calculated as the 
distance from the point $(x,y,z)$ to the point $Q$ of coordinates 
$[(x+y+z)/3,(x+y+z)/3,(x+y+z)/3]$, is
\begin{equation}
D^2=\frac{2}{3}(x^2+y^2+z^2-xy-xz-yz) .
\label{3-11}
\end{equation}
The quantities $s$ and $D$ are shown in Fig. \ref{fig2}, where the plane through the
point $P$ and perpendicular to the trisector line $(t)$ intersects the $x$ axis
at point $A$ of coordinates $(x+y+z,0,0)$, the $y$ axis at point $B$ of
coordinates $(0,x+y+z,0)$, and the $z$ axis at point $C$ of coordinates
$(0,0,x+y+z)$.  The azimuthal angle $\phi$ \index{azimuthal angle, tricomplex}
of the tricomplex number $x+hy+kz$
is defined as the angle in the plane $\Pi$ of the projection of $P$ on this
plane, measured from the line of intersection of the plane determined by the
line $(t)$ and the x axis with the plane $\Pi$, $0\leq\phi<2\pi$. 
The expression of $\phi$ in terms of $x,y,z$
can be obtained in a system of coordinates defined by the unit vectors
\begin{equation}
\xi_1: \frac{1}{\sqrt{6}}(2,-1,-1) ; \: \xi_2: \frac{1}{\sqrt{2}}(0,1,-1) ;
\: \xi_3: \frac{1}{\sqrt{3}}(1,1,1) ,
\label{3-12}
\end{equation}
and having the point $O$ as origin.
The relation between the
coordinates of $P$ in the systems $\xi_1,\xi_2,\xi_3$ and $x,y,z$ can be
written in the form
\begin{equation}
\left(\begin{array}{c}
\xi_1\\ \xi_2\\ \xi_3\end{array}\right)= 
\left(\begin{array}{ccc}
\frac{2}{\sqrt{6}}&-\frac{1}{\sqrt{6}}&-\frac{1}{\sqrt{6}}\\
0&\frac{1}{\sqrt{2}}&-\frac{1}{\sqrt{2}}\\
\frac{1}{\sqrt{3}}&\frac{1}{\sqrt{3}}&\frac{1}{\sqrt{3}}
\end{array}\right)
\left(\begin{array}{c}
x\\
y\\
z
\end{array}\right) .
\label{3-13}
\end{equation}
The components of the vector $OP$ in the system $\xi_1,\xi_2,\xi_3$ 
can be obtained with the aid of Eq. (\ref{3-13}) as  
\begin{equation}
(\xi_1,\xi_2,\xi_3)=\left(\frac{1}{\sqrt{6}}(2x-y-z),
\frac{1}{\sqrt{2}}(y-z),\frac{1}{\sqrt{3}}(x+y+z)\right) .
\label{3-14}
\end{equation}
The expression of the angle $\phi$ as a function of $x,y,z$ is then
\begin{equation}
\cos\phi=\frac{2x-y-z}{2(x^2+y^2+z^2-xy-xz-yz)^{1/2}} ,
\label{3-15}
\end{equation}
\begin{equation}
\sin\phi=\frac{\sqrt{3}(y-z)}{2(x^2+y^2+z^2-xy-xz-yz)^{1/2}} .
\label{3-16}
\end{equation}
It can be seen from Eqs. (\ref{3-15}),(\ref{3-16}) that the angle of points on the
$x$ axis is $\phi=0$, the angle of points on the $y$ axis is $\phi=2\pi/3$, and
the angle of points on the $z$ axis is $\phi=4\pi/3$.
The angle
$\phi$ is shown in 
Fig. \ref{fig2} in the plane parallel to $\Pi$, passing through $P$. The axis
$Q\xi_1^{\parallel}$ is parallel to the axis $O\xi_1$, 
the axis $Q\xi_2^{\parallel}$ is parallel to the axis $O\xi_2$, 
and the axis $Q\xi_3^{\parallel}$ is parallel to the axis $O\xi_3$, 
so that, in the plane $ABC$, the angle $\phi$ is measured from the line $QA$.
The angle $\theta$ \index{polar angle, tricomplex}
between the line $OP$ and the trisector line $(t)$ is given by
\begin{equation}
\tan\theta=\frac{D}{s},
\label{3-16a}
\end{equation}
where $0\leq\theta\leq\pi$.
It can be checked that
\begin{equation}
d^2=D^2+s^2 ,
\label{3-16b}
\end{equation}
where 
\begin{equation}
d^2=x^2+y^2+z^2,
\label{3-16c}
\end{equation}{\index{distance!tricomplex}
so that
\begin{equation}
D=d\sin\theta,\;s=d\cos\theta .
\label{3-16d}
\end{equation}

\begin{figure}
\begin{center}
\epsfig{file=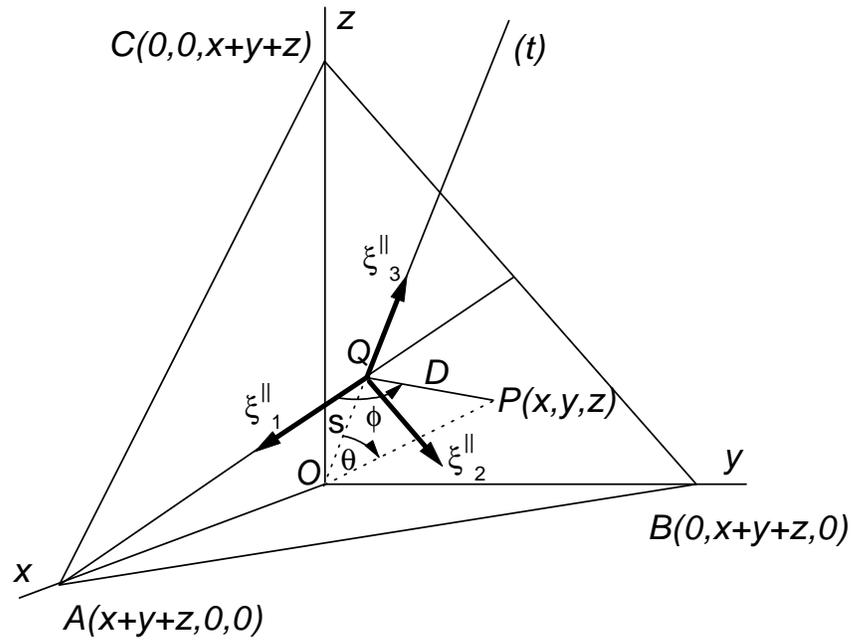,width=12cm}
\caption{Tricomplex variables $s, d, \theta, \phi$ for the tricomplex number
$x+hy+kz$, represented by the point $P(x,y,z)$. 
The azimuthal angle $\phi$ is shown in
in the plane parallel to $\Pi$, passing through $P$, which intersects 
the trisector line $(t)$ at $Q$ and the axes
of coordinates $x,y,z$ at the points $A, B, C$. The orthogonal axes
$\xi^\parallel _1,
\xi^\parallel_2, \xi^\parallel_3$ have the origin at $Q$.}
\label{fig2}
\end{center}
\end{figure}

The relations (\ref{3-10}), (\ref{3-11}), (\ref{3-15})-(\ref{3-16a}) can be used to
determine the associated 
projection $s$, the distance $D$, the polar angle $\theta$ with the trisector line
$(t)$  and the angle $\phi$ in the
$\Pi$ plane for the tricomplex number $x+hy+kz$.
It can be shown that if $u_1=x_1+hy_1+kz_1, 
u_2=x_2+hy_2+kz_2$ are tricomplex
numbers of projections, distances and angles $s_1,D_1,\theta_1, \phi_1$ and
respectively 
$s_2, D_2, \theta_2, \phi_2$, then the projection $s$, distance $D$ and the
angle $\theta, \phi$ for 
the product tricomplex number $u_1u_2
=x_1x_2+y_1z_2+y_2z_1+h(z_1z_2+x_1y_2+y_1x_2)+k(y_1y_2+x_1z_2+z_1x_2)$
are 
\begin{equation}
s=\sqrt{3}s_1s_2, \: D=\sqrt{\frac{3}{2}}D_1D_2, 
\: \tan\theta=\frac{1}{\sqrt{2}}\tan\theta_1\tan\theta_2\:, \phi=\phi_1+\phi_2 .
\label{3-17}
\end{equation}\index{transformation of variables!tricomplex}
The relations (\ref{3-17}) are consequences of the identities
\begin{eqnarray}
\lefteqn{(x_1x_2+y_1z_2+y_2z_1)+(z_1z_2+x_1y_2+y_1x_2)+(y_1y_2+x_1z_2+z_1x_2)}  
\nonumber\\
 & & =(x_1+y_1+z_1)(x_2+y_2+z_2) , 
\label{3-18}
\end{eqnarray}
\begin{eqnarray}
\lefteqn{(x_1x_2+y_1z_2+y_2z_1)^2+(z_1z_2+x_1y_2+y_1x_2)^2+(y_1y_2+x_1z_2+z_1x_2)^2} 
\nonumber\\
 & & -(x_1x_2+y_1z_2+y_2z_1)(z_1z_2+x_1y_2+y_1x_2)
-(x_1x_2+y_1z_2+y_2z_1)(y_1y_2+x_1z_2+z_1x_2) \nonumber\\
 & & -(z_1z_2+x_1y_2+y_1x_2)+(y_1y_2+x_1z_2+z_1x_2)  \nonumber\\
 & & =(x_1^2+y_1^2+z_1^2-x_1y_1-x_1z_1-y_1z_1) 
 (x_2^2+y_2^2+z_2^2-x_2y_2-x_2z_2-y_2z_2) ,
\label{3-19}
\end{eqnarray}
\begin{eqnarray}
\lefteqn{\frac{2x_1-y_1-z_1}{2}\frac{2x_2-y_2-z_2}{2}-\frac{\sqrt{3}}{2}(y_1-z_1) 
\frac{\sqrt{3}}{2}(y_2-z_2)}\nonumber\\
& & =\frac{1}{2}[2(x_1x_2+y_1z_2+z_1y_2)
-(z_1z_2+x_1y_2+y_1x_2)-(y_1y_2+x_1z_2+z_1x_2)] ,
\label{3-20}
\end{eqnarray}
\begin{eqnarray}
\lefteqn{\frac{\sqrt{3}}{2}(y_1-z_1)\frac{2x_2-y_2-z_2}{2}
+\frac{\sqrt{3}}{2}(y_2-z_2)\frac{2x_1-y_1-z_1}{2}}\nonumber\\
& & =\frac{\sqrt{3}}{2} [(z_1z_2+x_1y_2+y_1x_2)-(y_1y_2+x_1z_2+z_1x_2)] .
\label{3-21}
\end{eqnarray}
The relation (\ref{3-18}) shows that if $u$ is in the plane
$\Pi$, such that $x+y+z=0$, 
then the product $uu^\prime$ is also in the plane $\Pi$ for any $u^\prime$.
The relation (\ref{3-19}) shows that if $u$
is on the trisector line $(t)$, such that $x^2+y^2+z^2-xy-xz-yz=0$,
then $uu^\prime$ is also on the trisector line $(t)$ for any $u^\prime$.
If $u, u^\prime$ are points in the plane $x+y+z=1$, then the product
$uu^\prime$ is also in that plane, and if $u, u^\prime$ are points of the
cylindrical surface $x^2+y^2+z^2-xy-xz-yz=1$, then $uu^\prime$ is also in that
cylindrical surface. This means that if $u, u^\prime$ are points on the circle
$x+y+z=1, x^2+y^2+z^2-xy-xz-yz=1,$  
which is perpendicular to the trisector line, is situated at
a distance $1/\sqrt{3}$ from the origin and has the radius $\sqrt{2/3}$,
then the tricomplex product $uu^\prime$ is also on the same circle. This
invariant circle for the multiplication of tricomplex numbers is described by
the equations  
\begin{equation}
x=\frac{1}{3}+\frac{2}{3}\cos\phi,\:
y=\frac{1}{3}-\frac{1}{3}\cos\phi+\frac{1}{\sqrt{3}}\sin\phi,\:
z=\frac{1}{3}-\frac{1}{3}\cos\phi-\frac{1}{\sqrt{3}}\sin\phi.
\label{3-21b}
\end{equation}
It has the center at the point (1/3,1/3,1/3) and passes through the points
(1,0,0), (0,1,0) and (0,0,1), as shown in Fig. \ref{fig3}. 

\begin{figure}
\begin{center}
\epsfig{file=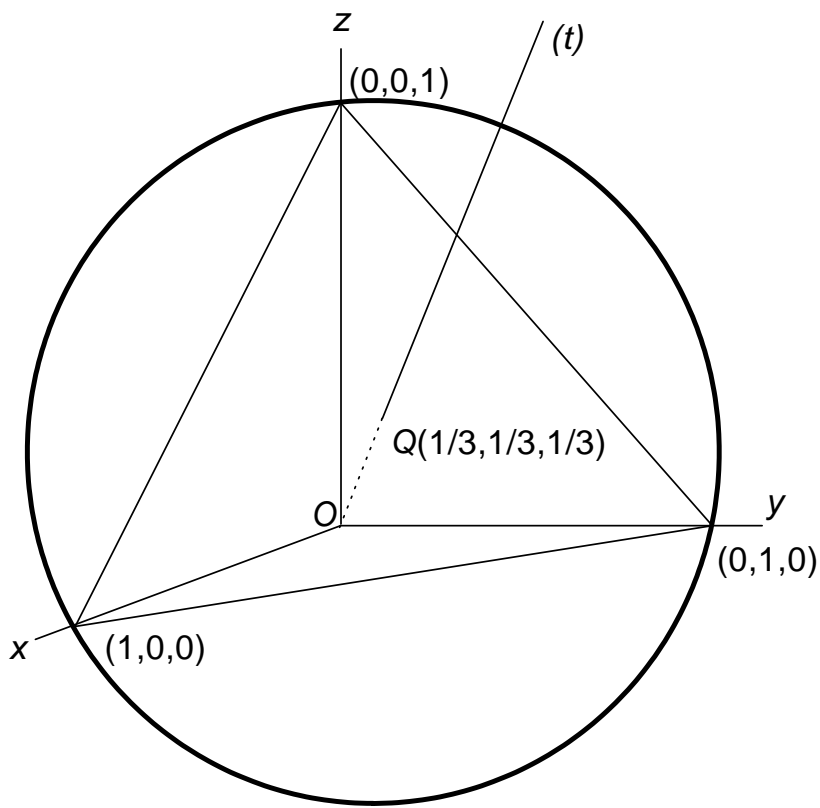,width=12cm}
\caption{Invariant circle for the multiplication of tricomplex numbers, lying
in a plane perpendicular to the trisector line and
passing through the points (1,0,0), (0,1,0) and (0,0,1). The center of the
circle is at the point (1/3,1/3,1/3), and its radius is $(2/3)^{1/2}$ .} 
\label{fig3}
\end{center}
\end{figure}

An important quantity is the amplitude $\rho$ defined as
$\rho=\nu^{1/3}$, so that
\begin{equation}
\rho^3=x^3+y^3+z^3-3xyz .
\label{3-22}
\end{equation}\index{amplitude!tricomplex}
The amplitude $\rho$ of the product $u_1u_2$ of the tricomplex numbers
$u_1,u_2$ of amplitudes $\rho_1,\rho_2$ is
\begin{equation}
\rho=\rho_1\rho_2 ,
\label{3-23}
\end{equation}\index{transformation of variables!tricomplex}
as can be seen from the identity
\begin{eqnarray}
\lefteqn{(x_1x_2+y_1z_2+y_2z_1)^3+(z_1z_2+x_1y_2+y_1x_2)^3+(y_1y_2+x_1z_2+z_1x_2)^3} 
\nonumber\\
 &&-3(x_1x_2+y_1z_2+y_2z_1)(z_1z_2+x_1y_2+y_1x_2)(y_1y_2+x_1z_2+z_1x_2)\nonumber\\ 
 &&=(x_1^3+y_1^3+z_1^3-3x_1y_1z_1)(x_2^3+y_2^3+z_2^3-3x_2y_2z_2) .
\label{3-24}
\end{eqnarray}
The identity in Eq. (\ref{3-24}) can be demonstrated with the aid of 
Eqs. (\ref{3-7}), (\ref{3-18}) and (\ref{3-19}). Another method
would be to use the representation of the multiplication of the tricomplex
numbers by matrices, in which the tricomplex number $u=x+hy+kz$ is
represented by the matrix
\begin{equation}
\left(\begin{array}{ccc}
x&y&z\\
z&x&y\\
y&z&x 
\end{array}\right) .
\label{3-25}
\end{equation}\index{matrix representation!tricomplex}
The product $u=x+hy+kz$ of the tricomplex numbers $u_1=x_1+hy_1+kz_1,
u_2=x_2+hy_2+kz_2$, is represented by the matrix multiplication
\begin{equation}
\left(\begin{array}{ccc}
x&y&z\\
z&x&y\\
y&z&x 
\end{array}\right) =
\left(\begin{array}{ccc}
x_1&y_1&z_1\\
z_1&x_1&y_1\\
y_1&z_1&x_1 
\end{array}\right) 
\left(\begin{array}{ccc}
x_2&y_2&z_2\\
z_2&x_2&y_2\\
y_2&z_2&x_2 
\end{array}\right) .
\label{3-26}
\end{equation}
If 
\begin{equation}
\nu={\rm det}\left(\begin{array}{ccc}
x&y&z\\
z&x&y\\
y&z&x 
\end{array}\right),
\label{3-27}
\end{equation}
it can be checked that
\begin{equation}
\nu=x^3+y^3+z^3-3xyz .
\label{3-27b}
\end{equation}
The identity (\ref{3-24}) is then a consequence of the fact the determinant 
of the product of matrices is equal to the product of the determinants 
of the factor matrices. 

It can be seen from Eqs. (\ref{3-10}) and (\ref{3-11}) that
\begin{equation}
x^3+y^3+z^3-3xyz=\frac{3\sqrt{3}}{2}sD^2 ,
\label{3-28}
\end{equation}\index{inverse, determinant!tricomplex}
which can be written with the aid of relations (\ref{3-16d}) and (\ref{3-22}) as
\begin{equation}
\rho=\frac{3^{1/2}}{2^{1/3}}d\sin^{2/3}\theta\cos^{1/3}\theta.
\label{3-28b}
\end{equation}
This means that the surfaces of constant $\rho$ are surfaces of rotation
having the trisector line $(t)$ as axis, as shown in Fig. \ref{fig4}.

\begin{figure}
\begin{center}
\epsfig{file=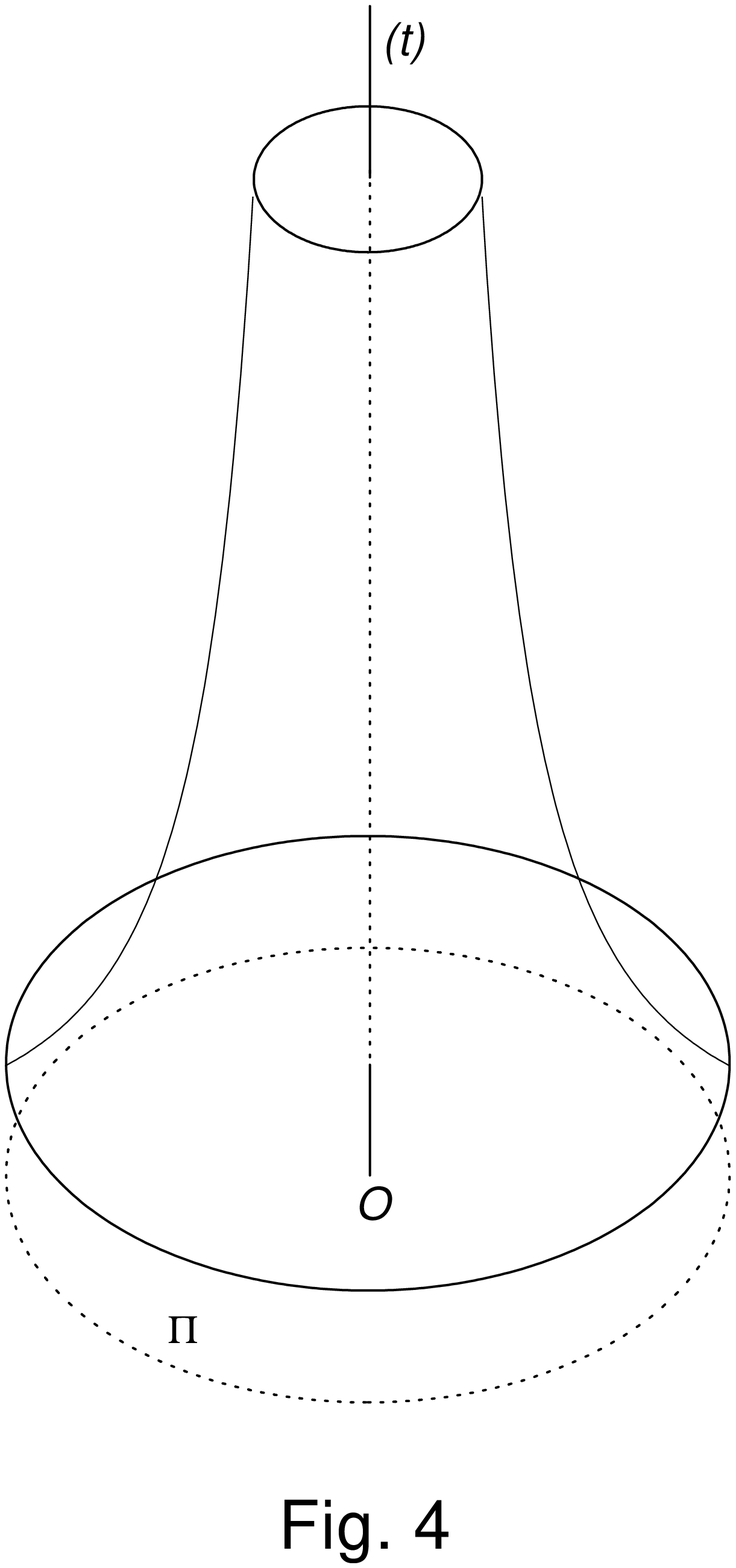,width=12cm}
\caption{Surfaces of constant $\rho$, which are surfaces of rotation
having the trisector line $(t)$ as axis.}
\label{fig4}
\end{center}
\end{figure}

\section{The tricomplex cosexponential functions}

The exponential function of a hypercomplex variable $u$ and the addition
theorem for the exponential function have been written in Eqs. 
~(1.35)-(1.36).
If $u=x+hy+kz$, then $\exp u$ can be calculated as $\exp u=\exp x \cdot
\exp (hy) \cdot \exp (kz)$. According to Eq. (\ref{3-1}), $h^2=k, h^3=1, k^2=h,
k^3=1$, and in general 
\begin{equation}
h^{3m}=1, h^{3m+1}=h, h^{3m+2}=k, k^{3m}=1, k^{3m+1}=k, k^{3m+2}=h ,
\label{3-31}
\end{equation}\index{complex units, powers of!tricomplex}
where $n$ is a natural number,
so that $\exp (hy)$ and $\exp(kz)$ can be written as
\begin{equation}
\exp (hy) = {\rm cx} \:y + h\: {\rm mx} \:y + k\: {\rm px} \:y , 
\label{3-32}
\end{equation}
\begin{equation}
\exp (kz) = {\rm cx} \:z + h \:{\rm px} \:z + k \:{\rm mx} \:z , 
\label{3-33}
\end{equation}\index{exponential, expressions!tricomplex}
where the functions cx, mx, px, which will be called in this chapter polar
cosexponential functions, are defined by the series 
\begin{equation}
{\rm cx \:y} = 1 + y^3/3! +y^6/6!+\cdots
\label{3-34}
\end{equation}
\begin{equation}
{\rm mx} \:y = y + y^4/4! + y^7/7!+\cdots
\label{3-35}
\end{equation}
\begin{equation}
{\rm px} \:y = y^2/2! + y^5/5! +y^8/8!+\cdots .
\label{3-36}
\end{equation}\index{cosexponential functions, tricomplex!definitions}
From the series definitions it can be seen that ${\rm cx} \: 0 =1, 
{\rm mx} \: 0=0, {\rm px} \: 0=0.$
The tridimensional polar cosexponential functions belong to the class
of the polar n-dimensional cosexponential functions
$g_{nk}$,  and ${\rm cx}=g_{30}, {\rm mx}=g_{31}, {\rm px}=g_{32}$.
It can be checked that
\begin{equation}
{\rm cx} \: y +{\rm px} \: y+{\rm mx} \: y =\exp y .
\label{3-36a}
\end{equation}
By expressing the fact that $\exp(hy+hz)=\exp(hy)\cdot\exp(hz)$ with the aid of
the cosexponential functions (\ref{3-34})-(\ref{3-36}) the following addition
theorems can be obtained
\begin{equation}
{\rm cx}\:(y+z)={\rm cx}\:y \:{\rm cx}\:z+{\rm mx}\:y 
\:{\rm px}\:z + {\rm px}\:y \:{\rm mx}\:z  ,
\label{3-38}
\end{equation}
\begin{equation}
{\rm mx}\:(y+z)={\rm px}\:y \:{\rm px}\:z+{\rm cx}\:y 
\:{\rm mx}\:z + {\rm mx}\:y \:{\rm cx}\:z  ,
\label{3-39}
\end{equation}
\begin{equation}
{\rm px}\:(y+z)={\rm mx}\:y \:{\rm mx}\:z+{\rm cx}\:y 
\:{\rm px}\:z + {\rm px}\:y \:{\rm cx}\:z  .
\label{3-40}
\end{equation}\index{cosexponential functions, tricomplex!addition theorems}
For $y=z$, Eqs. (\ref{3-38})-(\ref{3-40}) yield
\begin{equation}
{\rm cx}\: 2y={\rm cx}^2\:y +2 \:{\rm mx}\:y 
\:{\rm px}\:z ,
\label{3-38a}
\end{equation}
\begin{equation}
{\rm mx}\:2y={\rm px}^2 \:y +2\:{\rm cx}\:y 
\:{\rm mx}\:z ,
\label{3-39a}
\end{equation}
\begin{equation}
{\rm px}\:2y={\rm mx}^2\:y +2\:{\rm cx}\:y 
\:{\rm px}\:z .
\label{3-40a}
\end{equation}
The cosexponential functions are neither even nor odd functions. 
For $z=-y$, Eqs. (\ref{3-38})-(\ref{3-40}) yield
\begin{equation}
{\rm cx}\:y \:{\rm cx}\:(-y)+{\rm mx}\:y 
\:{\rm px}\:(-y) + {\rm px}\:y \:{\rm mx}\:(-y) =1  ,
\label{3-38b}
\end{equation}
\begin{equation}
{\rm px}\:y \:{\rm px}\:(-y)+{\rm cx}\:y 
\:{\rm mx}\:(-y) + {\rm mx}\:y \:{\rm cx}\:(-y)=0  ,
\label{3-39b}
\end{equation}
\begin{equation}
{\rm mx}\:y \:{\rm mx}\:(-y)+{\rm cx}\:y 
\:{\rm px}\:(-y) + {\rm px}\:y \:{\rm cx}\:(-y)=0  .
\label{3-40b}
\end{equation}

Expressions of the cosexponential functions in terms of regular exponential and
cosine functions can be obtained by considering the series expansions for
$e^{(h+k)y}$ and $e^{(h-k)y}$. These expressions can be obtained by calculating
first $(h+k)^n$ and $(h-k)^n$. It can be shown that
\begin{equation}
(h+k)^m=\frac{1}{3}\left[(-1)^{m-1}+2^m\right](h+k)+
\frac{2}{3}\left[(-1)^m+2^{m-1}\right] ,
\label{3-41}
\end{equation}
\begin{equation}
(h-k)^{2m}=(-1)^{m-1}3^{m-1}(k+k-2), \:(h-k)^{2m+1}=(-1)^m 3^m (h-k) ,
\label{3-42}
\end{equation}
where $n$ is a natural number. Then
\begin{equation}
e^{(h+k)y}=(h+k)\left(-\frac{1}{3}e^{-y}+\frac{1}{3}e^{2y}\right)
+\frac{2}{3}e^{-y}+\frac{1}{3}e^{2y} .
\label{3-43}
\end{equation}
As a corollary, the following identities can be obtained from Eq. (\ref{3-43}) by
writing $e^{(h+k)y}=e^{hy}e^{ky}$ and expressing $e^{hy}$ and $e^{ky}$ in terms
of cosexponential functions via Eqs. (\ref{3-32}) and (\ref{3-33}),
\begin{equation}
{\rm cx}^2 \: y +{\rm mx}^2 \: y+{\rm px}^2 \: y=\frac{2}{3} e^{-y}
+\frac{1}{3}e^{2y} , 
\label{3-43a}
\end{equation}
\begin{equation}
 {\rm cx} \: y \:\: {\rm mx} \: y
+ {\rm cx} \: y \:\: {\rm px} \: y
+ {\rm mx} \: y \:\: {\rm px} \: y  = -\frac{1}{3}e^{-y}+\frac{1}{3}
e^{2y} . 
\label{3-43b}
\end{equation}
From Eqs. (\ref{3-43a}) and (\ref{3-43b}) it results that
\begin{eqnarray}
\lefteqn{{\rm cx}^2 \: y +{\rm mx}^2 \: y+{\rm px}^2 \: y
}\nonumber\\ 
& & - {\rm cx} \: y \:\: {\rm mx} \: y
- {\rm cx} \: y \:\: {\rm px} \: y
- {\rm mx} \: y \:\: {\rm px} \: y  = \exp(-y) . 
\label{3-36b}
\end{eqnarray}
Then from Eqs. (\ref{3-7}), (\ref{3-36a}) and (\ref{3-36b}) it follows that
\begin{equation}
{\rm cx}^3 \: y +{\rm mx}^3 \: y+{\rm px}^3 \: y
 - 3 {\rm cx} \: y \:\: {\rm mx} \: y \:\:{\rm px} \: y =1 .
\label{3-37}
\end{equation}

Similarly,
\begin{equation}
e^{(h-k)y}=\frac{1}{3}(1+h+k)+\frac{1}{3}(2-h-k)\cos(\sqrt{3}y)
+\frac{1}{\sqrt{3}}(h-k)\sin(\sqrt{3}y) .
\label{3-44}
\end{equation}
The last relation can also be written as
\begin{eqnarray}
e^{(h-k)y}&=&\frac{1}{3}+\frac{2}{3}\cos(\sqrt{3}y)+h\left[\frac{1}{3}
+\frac{2}{3}\cos\left(\sqrt{3}y-\frac{2\pi}{3}\right)\right]\nonumber\\
&&+k\left[\frac{1}{3}+\frac{2}{3}\cos\left(\sqrt{3}y+\frac{2\pi}{3}\right)\right].
\label{3-44a}
\end{eqnarray}
As a corollary, the following identities can be obtained from Eq. (\ref{3-44}) by
writing $e^{(h-k)y}=e^{hy}e^{-ky}$ and expressing $e^{hy}$ and $e^{-ky}$ in terms
of cosexponential functions via Eqs. (\ref{3-32}) and (\ref{3-33}),
\begin{equation}
 {\rm cx} \: y \:\: {\rm cx} \: (-y)
+ {\rm mx} \: y \:\: {\rm mx} \: (-y)
+ {\rm px} \: y \:\: {\rm px} \: (-y)  =
\frac{1}{3}+\frac{2}{3}\cos(\sqrt{3}y) , 
\label{3-44b}
\end{equation}
\begin{eqnarray}
 {\rm cx} \: y \:\: {\rm px} \: (-y)
+ {\rm mx} \: y \:\: {\rm cx} \: (-y)
+ {\rm px} \: y \:\: {\rm mx} \: (-y) 
=\frac{1}{3}+\frac{2}{3}\cos\left(\sqrt{3}y-\frac{2\pi}{3}\right)
\label{3-44d}
\end{eqnarray}
\begin{eqnarray}
 {\rm cx} \: y \:\: {\rm mx} \: (-y)
+ {\rm mx} \: y \:\: {\rm px} \: (-y)
+ {\rm px} \: y \:\: {\rm cx} \: (-y) 
=\frac{1}{3}+\frac{2}{3}\cos\left(\sqrt{3}y+\frac{2\pi}{3}\right)
\label{3-44c}
\end{eqnarray}

Expressions of $e^{2hy}$ in terms of the regular exponential and
cosine functions can be obtained by the multiplication of the expressions of
$e^{(h+k)y}$ and $e^{(h-k)y}$ from Eqs. (\ref{3-43}) and (\ref{3-44}).
At the same time, Eq. (\ref{3-32}) gives an expression of $e^{2hy}$ in terms 
of cosexponential functions. By equating the real and hypercomplex parts of 
these two forms of $e^{2y}$ and then
replacing $2y$ by $y$ gives the expressions of the cosexponential functions as
\begin{equation}
{\rm cx} \:y =\frac{1}{3}\: e^y+\frac{2}{3}\cos\left(\frac{\sqrt{3}}{2}y\right)
\:e^{-y/2},
\label{3-45}
\end{equation}
\begin{equation}
{\rm mx} \:y =\frac{1}{3}\: e^y+\frac{2}{3}\cos\left(\frac{\sqrt{3}}{2}y-\frac{2\pi}{3}\right)
\:e^{-y/2} ,
\label{3-46}
\end{equation}
\begin{equation}
{\rm px} \:y =\frac{1}{3}\: e^y+\frac{2}{3}\cos\left(\frac{\sqrt{3}}{2}y+\frac{2\pi}{3}\right)
\:e^{-y/2} .
 \label{3-47}
\end{equation}\index{cosexponential functions, tricomplex!expressions}
It is remarkable that the series in Eqs. (\ref{3-34})-(\ref{3-36}), in which the 
terms are either of the form $y^{3m}$, or $y^{3m+1}$, or $y^{3m+2}$, 
can be expressed
in terms of elementary functions whose power series are not subject to such 
restrictions. 
The cosexponential functions differ by the
phase of the cosine function in their expression, and the designation of the
functions in Eqs. (\ref{3-46}) and (\ref{3-47}) as mx and px
refers respectively to the minus or plus sign of the phase term $2\pi/3$. 
The graphs of the cosexponential functions are shown in Fig. \ref{fig5}.

\begin{figure}
\begin{center}
\epsfig{file=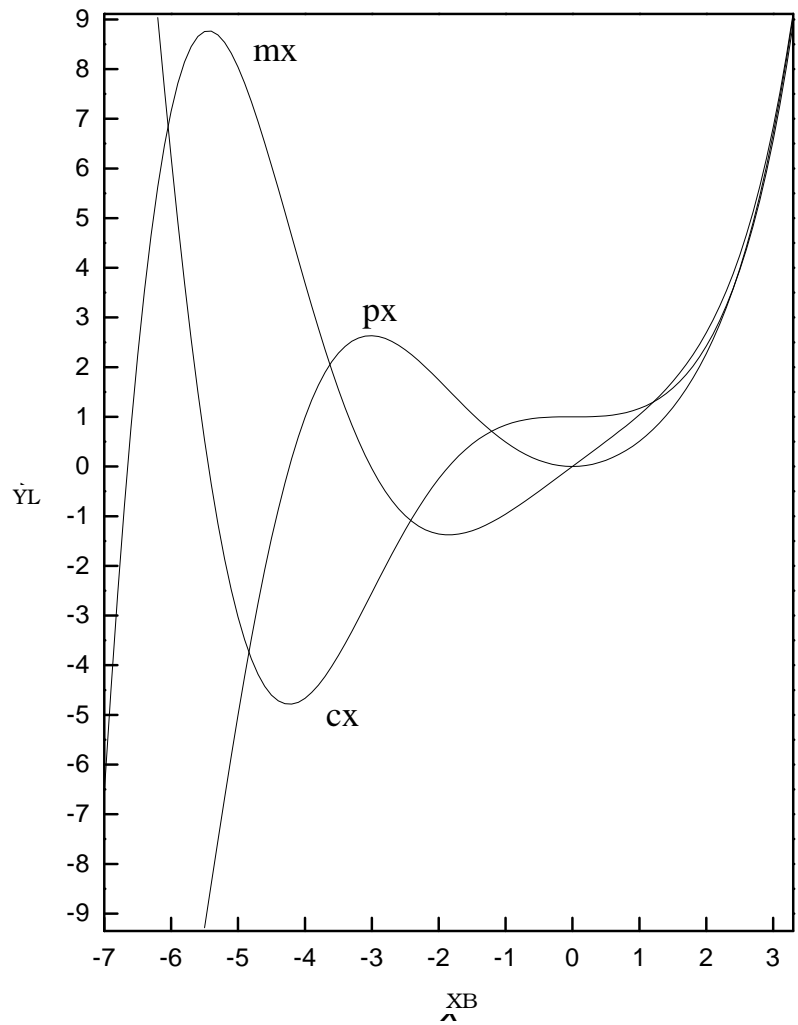,width=12cm}
\caption{Graphs of the cosexponential functions ${\rm cx}, {\rm mx},{\rm px}$.}
\label{fig5}
\end{center}
\end{figure}

It can be checked that the cosexponential functions are solutions of the
third-order differential equation
\begin{equation}
\frac{d^3\zeta}{du^3}=\zeta ,
\label{3-47b}
\end{equation}\index{cosexponential functions, tricomplex!differential equations}
whose solutions are of the form $\zeta(u)=A\:{\rm cx}\: u+B\:{\rm mx}\: u
+C\:{\rm px}\: u.$
It can also be checked that the derivatives of the cosexponential functions are
related by
\begin{equation}
\frac{d{\rm px}}{du}={\rm mx}, \:
\frac{d{\rm mx}}{du}={\rm cx}, \:
\frac{d{\rm cx}}{du}={\rm px} .
\label{3-47c}
\end{equation}

\section{Exponential and trigonometric forms of tricomplex numbers}

If for a tricomplex number $u=x+ky+kz$ another tricomplex number 
$u_1=x_1+hy_1+kz_1$ exists 
such that
\begin{equation}
x+hy+kz=e^{x_1+hy_1+kz_1} ,
\label{3-48}
\end{equation}
then $u_1$ is said to be the logarithm of $u$,
\begin{equation}
u_1=\ln u .
\label{3-49}
\end{equation}
The expressions of $x_1, y_1, z_1$ as functions of $x, y, z$ can be obtained by
developing $e^{hy_1}$ and $e^{kz_1}$ with the aid of Eqs. (\ref{3-32}) and
(\ref{3-33}), by multiplying these expressions and separating the hypercomplex
components, 
\begin{equation}
x=e^{x_1}[ {\rm cx} \: y_1 \:\: {\rm cx} \: z_1
+ {\rm mx} \: y_1 \:\: {\rm mx} \: z_1
+ {\rm px} \: y_1 \:\: {\rm px} \: z_1 ] ,
\label{3-50}
\end{equation}
\begin{equation}
y=e^{x_1}[ {\rm cx} \: y_1 \:\: {\rm px} \: z_1
+ {\rm mx} \: y_1 \:\: {\rm cx} \: z_1
+ {\rm px} \: y_1 \:\: {\rm mx} \: z_1 ] ,
\label{3-51}
\end{equation}
\begin{equation}
z=e^{x_1}[ {\rm cx} \: y_1 \:\: {\rm mx} \: z_1
+ {\rm px} \: y_1 \:\: {\rm cx} \: z_1
+ {\rm mx} \: y_1 \:\: {\rm px} \: z_1 ] ,
\label{3-52}
\end{equation}
Using Eq. (\ref{3-24}) with the substitutions $x_1\rightarrow{\rm cx}\: y_1,
y_1\rightarrow{\rm mx}\: y_1, z_1\rightarrow{\rm px}\: y_1, 
x_2\rightarrow{\rm cx}\: z_1, y_2\rightarrow{\rm px}\: z_1,
z_2\rightarrow{\rm mx}\: z_1$ and then the identity (\ref{3-37}) yields
\begin{equation}
x^3+y^3+z^3-3xyz=e^{3x_1} , 
\label{3-53}
\end{equation}
whence
\begin{equation}
x_1=\frac{1}{3}\ln(x^3+y^3+z^3-3xyz) .
\label{3-54}
\end{equation}
The logarithm in Eq. (\ref{3-54}) exists as a real function for $x+y+z>0$.
A further relation can be obtained by summing Eqs. (\ref{3-50})-(\ref{3-52}) and
then using the addition theorems (\ref{3-38})-(\ref{3-40}) 
\begin{equation}
\frac{x+y+z}{(x^3+y^3+z^3-3xyz)^{1/3}}={\rm cx}\:(y_1+z_1)+
{\rm mx}\:(y_1+z_1)+{\rm px}\:(y_1+z_1) .
\label{3-55}
\end{equation}
The sum in Eq. (\ref{3-55}) is according to Eq. (\ref{3-36a}) $e^{y_1+z_1}$, so
that 
\begin{equation}
y_1+z_1=\ln\frac{x+y+z}{(x^3+y^3+z^3-3xyz)^{1/3}} .
\label{3-56}
\end{equation}
The logarithm in Eq. (\ref{3-56}) is defined for points which are not on the
trisector line $(t)$, so that $x^2+y^2+z^2-xy-xz-yz\not= 0$.
Substituting in Eq. (\ref{3-48}) the expression of $x_1$, Eq. (\ref{3-54}), and of
$z_1$ as a function of $x, y, z, y_1,$ Eq. (\ref{3-56}), yields
\begin{equation}
\frac{u}{\rho}
\exp\left[-k\ln\left(\frac{\sqrt{2}s}{D}\right)^{2/3}\right]=e^{(h-k)y_1} ,
\label{3-57}
\end{equation}
where the quantities $\rho, s$ and $D$ have been defined in Eqs.
(\ref{3-22}),(\ref{3-10}) and (\ref{3-11}).
Developing the exponential functions in the left-hand side of Eq. (\ref{3-57}) 
with the aid of Eq. (\ref{3-33}) and using the expressions of the cosexponential
functions, Eqs. (\ref{3-45})-(\ref{3-47}), and using the relation (\ref{3-44}) for
the right-hand side of Eq. (\ref{3-57}) yields for the real part
\begin{equation}
\frac{\left(x-\frac{y+z}{2}\right)\cos\left[\frac{1}{\sqrt{3}}\ln
\left(\frac{\sqrt{2}s}{D}\right)\right]
-\frac{\sqrt{3}}{2}(y-z)\sin
\left[\frac{1}{\sqrt{3}}\ln
\left(\frac{\sqrt{2}s}{D}\right)\right]}
{(x^2+y^2+z^2-xy-xz-yz)^{1/2}}=\cos(\sqrt{3}y_1) ,
\label{3-58}
\end{equation}
which can also be written as
\begin{equation}
\cos\left[\frac{1}{\sqrt{3}}\ln\left(\frac{\sqrt{2}s}{D}\right)+\phi\right]
=\cos(\sqrt{3}y_1)
\label{3-59}
\end{equation}
where $\phi$ is the angle defined in Eqs. (\ref{3-15}) and (\ref{3-16}).
Thus
\begin{equation}
y_1=\frac{1}{3}\ln\left(\frac{\sqrt{2}s}{D}\right)+\frac{1}{\sqrt{3}}\phi .
\label{3-60}
\end{equation}
The exponential form of the tricomplex number $u$ is then
\begin{equation}
u=\rho\:
\exp\left[\frac{1}{3}(h+k)\ln\frac{\sqrt{2}}{\tan\theta}+
\frac{1}{\sqrt{3}}(h-k)\phi\right] ,
\label{3-61}
\end{equation}\index{exponential form!tricomplex}
where $\theta$ is the angle between the line $OP$ 
connecting the origin to the point $P$ of coordinates $(x,y,z)$ and the
trisector line $(t)$, defined in Eq. (\ref{3-16a}) and shown in Fig. \ref{fig2}.
The exponential in Eq. (\ref{3-61}) can be expanded with the aid of Eq.
(\ref{3-43}) and (\ref{3-44a}) as
\begin{equation}
\exp\left[\frac{1}{3}(h+k)\ln\frac{\sqrt{2}}{\tan\theta}\right] 
=\frac{2-h-k}{3}\left(\frac{\tan\theta}{\sqrt{2}}\right)^{1/3}
+\frac{1+h+k}{3}\left(\frac{\sqrt{2}}{\tan\theta}\right)^{2/3} ,
\label{3-62}
\end{equation}
so that 
\begin{equation}
x+hy+kz=\rho
\left[\frac{2-h-k}{3}\left(\frac{\tan\theta}{\sqrt{2}}\right)^{1/3}
+\frac{1+h+k}{3}\left(\frac{\sqrt{2}}{\tan\theta}\right)^{2/3}\right]
\exp\left(\frac{h-k}{\sqrt{3}}\phi\right) .
\label{3-63}
\end{equation}
Substituting in Eq. (\ref{3-63}) the expression of the amplitude $\rho$, Eq.
(\ref{3-28b}), yields
\begin{equation}
u=d\sqrt{\frac{3}{2}}
\left(\frac{2-h-k}{3}\sin\theta
+\frac{1+h+k}{3}\sqrt{2}\cos\theta\right)
\exp\left(\frac{h-k}{\sqrt{3}}\phi\right) ,
\label{3-63b}
\end{equation}\index{trigonometric form!tricomplex}
which is the trigonometric form of the tricomplex number $u$.
As can be seen from Eq. (\ref{3-63b}), the tricomplex number $x+hy+kz$ is
written as the product of the modulus $d$, of a factor depending on the
polar angle $\theta$ with respect to the trisector line, and of a factor depending 
of the azimuthal angle $\phi$ in the plane $\Pi$ perpendicular to the
trisector line. 
The exponential in Eq. (\ref{3-63b}) can be expanded further with the aid of
Eq. (\ref{3-44a}) as 
\begin{equation}
\exp\left(\frac{1}{\sqrt{3}}(h-k)\phi\right) 
=\frac{1+h+k}{3}+\frac{2-h-k}{3}\cos\phi +\frac{h-k}{\sqrt{3}}\sin\phi ,
\label{3-64}
\end{equation}
so that the tricomplex number $x+hy+kz$ can also be written, after
multiplication of the factors, in the form
\begin{eqnarray}
\lefteqn{x+hy+kz=
\frac{2-h-k}{3}(x^2+y^2+z^2-xy-xz-yz)^{1/2}\cos\phi}\nonumber\\
&&+\frac{h-k}{\sqrt{3}}(x^2+y^2+z^2-xy-xz-yz)^{1/2}\sin\phi
+\frac{1+h+k}{3}(x+y+z)
\label{3-65}
\end{eqnarray}\index{trigonometric form!tricomplex}
The validity of Eq. (\ref{3-65}) can be checked by substituting the expressions
of $\cos\phi$ and $\sin\phi$ from Eqs. (\ref{3-15}) and (\ref{3-16}).
 
\section{Elementary functions of a tricomplex variable}

It can be shown with the aid of Eq. (\ref{3-61}) that 
\begin{eqnarray}
\lefteqn{(x+hy+kz)^m=
\rho^m
\left[\frac{2-h-k}{3}\left(\frac{\tan\theta}{\sqrt{2}}\right)^{m/3}
+\frac{1+h+k}{3}\left(\frac{\sqrt{2}}{\tan\theta}\right)^{2m/3}\right]
\exp\left(\frac{h-k}{\sqrt{3}}m\phi\right) ,\nonumber}\\
&&
\label{3-66}
\end{eqnarray}\index{power function!tricomplex}
or equivalently
\begin{eqnarray}
\lefteqn{(x+hy+kz)^m=
\frac{2-h-k}{3}(x^2+y^2+z^2-xy-xz-yz)^{m/2}\cos(m\phi)}
\nonumber\\
&&+\frac{h-k}{\sqrt{3}}(x^2+y^2+z^2-xy-xz-yz)^{m/2}\sin(m\phi)
+\frac{1+h+k}{3}(x+y+z)^m,\nonumber\\
&&
\label{3-67}
\end{eqnarray}
which are valid for real values of $m$. Thus Eqs. (\ref{3-66}) or (\ref{3-67})
define the power function $u^m$ of the tricomplex variable
$u$. 

The power function is multivalued unless $m$ is an integer. 
It can be inferred from Eq. (\ref{3-61}) that, for integer values of $m$,
\begin{equation}
(uu^\prime)^m=u^m\:u^{\prime m} .
\label{3-68}
\end{equation}
For natural
$m$, Eq. (\ref{3-67}) can be checked with the aid of relations (\ref{3-15}) and
(\ref{3-16}). For example if $m=2$, it can be checked that the right-hand side of
Eq. (\ref{3-67}) is equal to $(x+hy+kz)^2=x^2+2yz+h(z^2+2xz)+k(y^2+2xz)$.

The logarithm $u_1$ of the tricomplex number $u$, $u_1=\ln u$, can be defined
as the solution of Eq. (\ref{3-48}) for $u_1$ as a function of $u$. 
For $x+y+z>0$, from Eq. (\ref{3-61}) it results that 
\begin{equation}
\ln u=\ln \rho+
\frac{1}{3}(h+k)\ln\left(\frac{\tan\theta}{\sqrt{2}}\right)+
\frac{1}{\sqrt{3}}(h-k)\phi .
\label{3-69}
\end{equation}\index{logarithm!tricomplex}
It can be checked with the aid of Eqs. (\ref{3-17}) and (\ref{3-23}) that
\begin{equation}
\ln(uu^\prime)=\ln u+\ln u^\prime ,
\label{3-69b}
\end{equation}
which is valid up to integer multiples of $2\pi(h-k)/\sqrt{3}$.

The trigonometric functions of the hypercomplex variable
$u$ and the addition theorems for these functions have been written in Eqs.
~(1.57)-(1.60). 
The trigonometric functions of the hypercomplex variables $hy, ky$ can be
expressed in terms of the cosexponential functions as
\begin{equation}
\cos(hy)=\frac{1}{2}[{\rm cx}\:(iy)+{\rm cx}\: (-iy)]+ 
\frac{1}{2}h[{\rm mx}\:(iy)+{\rm mx}\: (-iy)]+
\frac{1}{2}k[{\rm px}\:(iy)+{\rm px}\: (-iy)] ,
\label{3-74}
\end{equation}\index{trigonometric functions, expressions!tricomplex}
\begin{equation}
\cos(ky)=\frac{1}{2}[{\rm cx}\:(iy)+{\rm cx}\: (-iy)]+ 
\frac{1}{2}h[{\rm px}\:(iy)+{\rm px}\: (-iy)]+
\frac{1}{2}k[{\rm mx}\:(iy)+{\rm mx}\: (-iy)] ,
\label{3-75}
\end{equation}
\begin{equation}
\sin(hy)=\frac{1}{2i}[{\rm cx}\:(iy)-{\rm cx}\: (-iy)]+ 
\frac{1}{2i}h[{\rm mx}\:(iy)-{\rm mx}\: (-iy)]+
\frac{1}{2i}k[{\rm px}\:(iy)-{\rm px}\: (-iy)] ,
\label{3-76}
\end{equation}
\begin{equation}
\sin(ky)=\frac{1}{2i}[{\rm cx}\:(iy)-{\rm cx}\: (-iy)]+ 
\frac{1}{2i}h[{\rm px}\:(iy)-{\rm px}\: (-iy)]+
\frac{1}{2i}k[{\rm mx}\:(iy)-{\rm mx}\: (-iy)] ,
\label{3-77}
\end{equation}
where $i$ is the imaginary unit. Using the expressions of the cosexponential
functions in Eqs. (\ref{3-45})-(\ref{3-47}) gives expressions of the 
trigonometric functions of $hy, hz$ as
\begin{eqnarray}
\lefteqn{\cos(hy) =\frac{1}{3}\cos y+
\frac{2}{3}\cosh\left(\frac{\sqrt{3}}{2}y\right)\cos\frac{y}{2}+}\nonumber\\ 
&&+h\left[\frac{1}{3}\cos y
-\frac{1}{3}\cosh\left(\frac{\sqrt{3}}{2}y\right)\cos\frac{y}{2} 
+\frac{1}{\sqrt{3}}\sinh\left(\frac{\sqrt{3}}{2}y\right)\sin\frac{y}{2}
\nonumber\right] \\
&&+k\left[\frac{1}{3}\cos y
-\frac{1}{3}\cosh\left(\frac{\sqrt{3}}{2}y\right)\cos\frac{y}{2} 
-\frac{1}{\sqrt{3}}\sinh\left(\frac{\sqrt{3}}{2}y\right)\sin\frac{y}{2}
\right],
\label{3-78}
\end{eqnarray}
\begin{eqnarray}
\lefteqn{\cos(ky) =\frac{1}{3}\cos y+
\frac{2}{3}\cosh\left(\frac{\sqrt{3}}{2}y\right)\cos\frac{y}{2}+
}\nonumber\\ 
&&+h\left[\frac{1}{3}\cos y
-\frac{1}{3}\cosh\left(\frac{\sqrt{3}}{2}y\right)\cos\frac{y}{2} 
-\frac{1}{\sqrt{3}}\sinh\left(\frac{\sqrt{3}}{2}y\right)\sin\frac{y}{2}
\nonumber\right] \\
&&+k\left[\frac{1}{3}\cos y
-\frac{1}{3}\cosh\left(\frac{\sqrt{3}}{2}y\right)\cos\frac{y}{2} 
+\frac{1}{\sqrt{3}}\sinh\left(\frac{\sqrt{3}}{2}y\right)\sin\frac{y}{2}
\right],
\label{3-79}
\end{eqnarray}
\begin{eqnarray}
\lefteqn{\sin(hy) =\frac{1}{3}\sin y
-\frac{2}{3}\cosh\left(\frac{\sqrt{3}}{2}y\right)\sin\frac{y}{2}}
\nonumber\\ 
&&
+h\left[\frac{1}{3}\sin y
+\frac{1}{3}\cosh\left(\frac{\sqrt{3}}{2}y\right)\sin\frac{y}{2} 
+\frac{1}{\sqrt{3}}\sinh\left(\frac{\sqrt{3}}{2}y\right)\cos\frac{y}{2}
\nonumber\right] \\
&&
+k\left[\frac{1}{3}\sin y
+\frac{1}{3}\cosh\left(\frac{\sqrt{3}}{2}y\right)\sin\frac{y}{2} 
-\frac{1}{\sqrt{3}}\sinh\left(\frac{\sqrt{3}}{2}y\right)\cos\frac{y}{2}
\right],
\label{3-80}
\end{eqnarray}
\begin{eqnarray}
\lefteqn{\sin(ky) =\frac{1}{3}\sin y
-\frac{2}{3}\cosh\left(\frac{\sqrt{3}}{2}y\right)\sin\frac{y}{2}}\nonumber\\ 
&&
+h\left[\frac{1}{3}\sin y
+\frac{1}{3}\cosh\left(\frac{\sqrt{3}}{2}y\right)\sin\frac{y}{2} 
-\frac{1}{\sqrt{3}}\sinh\left(\frac{\sqrt{3}}{2}y\right)\cos\frac{y}{2}
\nonumber\right] \\
&&
+k\left[\frac{1}{3}\sin y
+\frac{1}{3}\cosh\left(\frac{\sqrt{3}}{2}y\right)\sin\frac{y}{2} 
+\frac{1}{\sqrt{3}}\sinh\left(\frac{\sqrt{3}}{2}y\right)\cos\frac{y}{2}
\right].
\label{3-81}
\end{eqnarray}
The trigonometric functions of a tricomplex number $x+hy+kz$ can then be
expressed in terms of elementary functions with the aid of the addition
theorems Eqs. (1.59), (1.60) and of the expressions in  Eqs. 
(\ref{3-78})-(\ref{3-81}).

The hyperbolic functions of the hypercomplex variable
$u$ and the addition theorems for these functions have been written in Eqs.
~(1.62)-(1.65). 
The hyperbolic functions of the hypercomplex variables $hy, ky$ can be
expressed in terms of the elementary functions as
\begin{eqnarray}
\lefteqn{\cosh(hy) =\frac{1}{3}\cosh y+
\frac{2}{3}\cos\left(\frac{\sqrt{3}}{2}y\right)\cosh\frac{y}{2}+}\nonumber\\ 
&&+h\left[\frac{1}{3}\cosh y
-\frac{1}{3}\cos\left(\frac{\sqrt{3}}{2}y\right)\cosh\frac{y}{2} 
-\frac{1}{\sqrt{3}}\sin\left(\frac{\sqrt{3}}{2}y\right)\sinh\frac{y}{2}
\nonumber\right] \\
&&+k\left[\frac{1}{3}\cosh y
-\frac{1}{3}\cos\left(\frac{\sqrt{3}}{2}y\right)\cosh\frac{y}{2} 
+\frac{1}{\sqrt{3}}\sin\left(\frac{\sqrt{3}}{2}y\right)\sinh\frac{y}{2}\right],
\label{3-81e}
\end{eqnarray}\index{hyperbolic functions, expressions!tricomplex}
\begin{eqnarray}
\lefteqn{\cosh(ky) =\frac{1}{3}\cosh y+
\frac{2}{3}\cos\left(\frac{\sqrt{3}}{2}y\right)\cosh\frac{y}{2}+}\nonumber\\ 
&&+h\left[\frac{1}{3}\cosh y
-\frac{1}{3}\cos\left(\frac{\sqrt{3}}{2}y\right)\cosh\frac{y}{2} 
+\frac{1}{\sqrt{3}}\sin\left(\frac{\sqrt{3}}{2}y\right)\sinh\frac{y}{2}
\nonumber\right] \\
&&+k\left[\frac{1}{3}\cosh y
-\frac{1}{3}\cos\left(\frac{\sqrt{3}}{2}y\right)\cosh\frac{y}{2} 
-\frac{1}{\sqrt{3}}\sin\left(\frac{\sqrt{3}}{2}y\right)\sinh\frac{y}{2}
\right],
\label{3-81f}
\end{eqnarray}
\begin{eqnarray}
\lefteqn{\sinh(hy) =\frac{1}{3}\sinh y
-\frac{2}{3}\cos\left(\frac{\sqrt{3}}{2}y\right)\sinh\frac{y}{2}
}\nonumber\\ 
&&+h\left[\frac{1}{3}\sinh y
+\frac{1}{3}\cos\left(\frac{\sqrt{3}}{2}y\right)\sinh\frac{y}{2} 
+\frac{1}{\sqrt{3}}\sin\left(\frac{\sqrt{3}}{2}y\right)\cosh\frac{y}{2}
\nonumber\right] \\
&&+k\left[\frac{1}{3}\sinh y
+\frac{1}{3}\cos\left(\frac{\sqrt{3}}{2}y\right)\sinh\frac{y}{2} 
-\frac{1}{\sqrt{3}}\sin\left(\frac{\sqrt{3}}{2}y\right)\cosh\frac{y}{2}
\right],
\label{3-81g}
\end{eqnarray}
\begin{eqnarray}
\lefteqn{\sinh(ky) =\frac{1}{3}\sinh y
-\frac{2}{3}\cos\left(\frac{\sqrt{3}}{2}y\right)\sinh\frac{y}{2}
}\nonumber\\ 
&&+h\left[\frac{1}{3}\sinh y
+\frac{1}{3}\cos\left(\frac{\sqrt{3}}{2}y\right)\sinh\frac{y}{2} 
-\frac{1}{\sqrt{3}}\sin\left(\frac{\sqrt{3}}{2}y\right)\cosh\frac{y}{2}
\nonumber\right] \\
&&+k\left[\frac{1}{3}\sinh y
+\frac{1}{3}\cos\left(\frac{\sqrt{3}}{2}y\right)\sinh\frac{y}{2} 
+\frac{1}{\sqrt{3}}\sin\left(\frac{\sqrt{3}}{2}y\right)\cosh\frac{y}{2}
\right].
\label{3-81h}
\end{eqnarray}
The hyperbolic functions of a tricomplex number $x+h y+k z$ can then
be expressed in terms of the elementary functions with the aid of the addition
theorems Eqs. (1.64), (1.65) and of the expressions in  Eqs. 
(\ref{3-81e})-(\ref{3-81h}).


\section{Tricomplex power series}

A tricomplex series is an infinite sum of the form
\begin{equation}
a_0+a_1+a_2+\cdots+a_l+\cdots , 
\label{3-82}
\end{equation}\index{series!tricomplex}
where the coefficients $a_n$ are tricomplex numbers. The convergence of 
the series (\ref{3-82}) can be defined in terms of the convergence of its 3 real
components. The convergence of a tricomplex series can however be studied
using tricomplex variables. The main criterion for absolute convergence 
remains the comparison theorem, but this requires a number of inequalities
which will be discussed further.

The modulus of a tricomplex number $u=x+hy+kz$ can be defined as
\begin{equation}
|u|=(x^2+y^2+z^2)^{1/2} .
\label{3-83}
\end{equation}\index{modulus, definition!tricomplex}
Since $|x|\leq |u|, |y|\leq |u|, |z|\leq |u|$, a property of absolute 
convergence established via a comparison theorem based on the modulus of the
series (\ref{3-82}) will ensure the absolute convergence of each real component
of that series.

The modulus of the sum $u_1+u_2$ of the tricomplex numbers $u_1, u_2$ fulfils
the inequality
\begin{equation}
||u_1|-|u_2||\leq |u_1+u_2|\leq |u_1|+|u_2| .
\label{3-84}
\end{equation}\index{modulus, inequalities!tricomplex}
For the product the relation is 
\begin{equation}
|u_1u_2|\leq \sqrt{3}|u_1||u_2| ,
\label{3-85}
\end{equation}
which replaces the relation of equality extant for regular complex numbers.
The equality in Eq. (\ref{3-85}) takes place for $x_1=y_1=z_1$ and 
$x_2=y_2=z_2$, i.e
when both tricomplex numbers lie on the trisector line $(t)$.
Using Eq. (\ref{3-65}), the relation (\ref{3-85}) can be written equivalently as 
\begin{equation}
\frac{2}{3}\delta_1^2\delta_2^2+\frac{1}{3}\sigma_1^2\sigma_2^2\leq 
3\left(\frac{2}{3}\delta_1^2+\frac{1}{3}\sigma_1^2\right)
\left(\frac{2}{3}\delta_2^2+\frac{1}{3}\sigma_2^2\right) ,
\label{3-85a}
\end{equation}
where $\delta_j^2=x_j^2+y_j^2+z_j^2-x_jy_j-x_jz_j-y_jz_j, \sigma_j=x_j+y_j+z_j,
j=1,2$, the equality taking place for $\delta_1=0,\delta_2=0$.
A particular form of Eq. (\ref{3-85}) is
\begin{equation}
|u^2|\leq \sqrt{3} |u|^2 ,
\label{3-86}
\end{equation}
and it can be shown that
\begin{equation}
|u^l|\leq 3^{(l-1)/2}|u|^l ,
\label{3-87}
\end{equation}
the equality in Eqs. (\ref{3-86}) and (\ref{3-87}) taking place for $x=y=z$.
It can be shown from Eq. (\ref{3-67}) that
\begin{equation}
|u^l|^2=\frac{2}{3}\delta^{2l}+\frac{1}{3}\sigma^{2l} ,
\label{3-87a}
\end{equation}
where $\delta^2=x^2+y^2+z^2-xy-xz-yz, \sigma=x+y+z$.
Then Eq. (\ref{3-87}) can also be written as
\begin{equation}
\frac{2}{3}\delta^{2l}+\frac{1}{3}\sigma^{2l}\leq
3^{l-1}\left(\frac{2}{3}\delta^2+\frac{1}{3}\sigma^2\right)^l ,
\label{3-87b}
\end{equation}
the equality taking place for $\delta=0$.
From Eqs. (\ref{3-85}) and (\ref{3-87}) it results that
\begin{equation}
|au^l|\leq 3^{l/2} |a| |u|^l .
\label{3-88}
\end{equation}
It can also be shown that
\begin{equation}
\left|\frac{1}{u}\right|\geq\frac{1}{|u|} ,
\label{3-89}
\end{equation}
the equality taking place for $\sigma^2=\delta^2$, or $xy+xz+yz=0$.

A power series of the tricomplex variable $u$ is a series of the form
\begin{equation}
a_0+a_1 u + a_2 u^2+\cdots +a_l u^l+\cdots .
\label{3-90}
\end{equation}\index{power series!tricomplex}
Since
\begin{equation}
\left|\sum_{l=0}^\infty a_l u^l\right| \leq \sum_{l=0}^\infty
3^{l/2}|a_l||u|^l ,
\label{3-91}
\end{equation}
a sufficient condition for the absolute convergence of this series is that
\begin{equation}
\lim_{n\rightarrow \infty}
\frac{\sqrt{3}|a_{l+1}||u|}{|a_l|}<1 .
\label{3-92}
\end{equation}
Thus the series is absolutely convergent for 
\begin{equation}
|u|<c_0,
\label{3-new86}
\end{equation}\index{convergence of power series!tricomplex}
where 
\begin{equation}
c_0=\lim_{l\rightarrow\infty} \frac{|a_l|}{\sqrt{3}|a_{l+1}|} .
\label{3-new87}
\end{equation}

The convergence of the series (\ref{3-90}) can be also studied with the aid of
a transformation which explicits the transverse and longitudinal parts
of the variable $u$ and of the constants $a_l$, 
\begin{equation}
x+hy+kz=v_1 e_1+\tilde v_1 \tilde e_1+v_+ e_+ , 
\label{3-116}
\end{equation}
where
\begin{equation}
v_1=\frac{2x-y-z}{2}, \:\:\tilde v_1=\frac{\sqrt{3}}{2}(y-z), \:\: v_+=x+y+z,
\label{3-117a}
\end{equation}\index{canonical variables!tricomplex}
and
\begin{equation}
e_1=\frac{2-h-k}{3},\:\:\tilde e_1=\frac{h-k}{\sqrt{3}},\:\:e_+=\frac{1+h+k}{3} .
\label{3-117b}
\end{equation}\index{canonical base!tricomplex}
The variables $v_1, \tilde v_1, v_+$ will be called 
canonical tricomplex variables,
$e_1, \tilde e_1, e_+$ will be called the  canonical tricomplex base,
and Eq. (\ref{3-116}) gives the canonical form of a tricomplex number.
In the geometric representation of Fig. \ref{fig6}, 
$e_1, \tilde e_1$ are situated in the plane $\Pi$, and $e_+$ is lying on the trisector
line $(t)$.  It can be checked that
\begin{equation}
e_1^2=e_1, \:\:\tilde e_1^2=-e_1,\:\: e_1\tilde e_1=\tilde e_1,\:\: e_1e_+=0, \:\:\tilde e_1e_+=0, 
\:\:e_+^2=e_+. 
\label{3-118}
\end{equation}
The moduli of the bases in Eq. (\ref{3-118}) are
\begin{equation}
|e_1|=\sqrt{\frac{2}{3}},\;|\tilde e_1|=\sqrt{\frac{2}{3}},\;|e_+|=\sqrt{\frac{1}{3}},
\label{3-118b}
\end{equation}
and it can be checked that
\begin{equation}
|x+hy+kz|^2=\frac{2}{3}(v_1^2+\tilde v_1^2)+\frac{1}{3}v_+^2 .
\label{3-118c}
\end{equation}\index{modulus, canonical variables!tricomplex}

\begin{figure}
\begin{center}
\epsfig{file=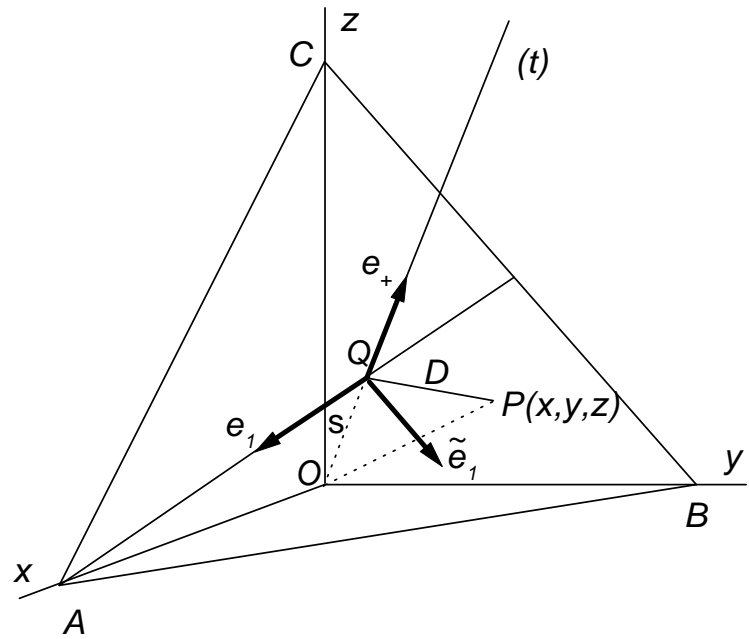,width=12cm}
\caption{Unit vectors $e_1, \tilde e_1, e_+$ of the orthogonal system of coordinates
with origin at $Q$.
The plane parallel to $\Pi$ passing through $P$ intersects 
the trisector line $(t)$ at $Q$ and the axes
of coordinates $x,y,z$ at the points $A, B, C$. }
\label{fig6}
\end{center}
\end{figure}

If $u=u^\prime u^{\prime\prime}$, the transverse and longitudinal
components are related by
\begin{equation}
v_1=v_1^\prime v_1^{\prime\prime}-\tilde v_1^\prime \tilde v_1^{\prime\prime}, \:\:
\tilde v_1=v_1^\prime \tilde v_1^{\prime\prime}+\tilde v_1^\prime v_1^{\prime\prime}, \:\:
v_+=v_+^\prime v_+^{\prime\prime} ,
\label{3-119}
\end{equation}\index{transformation of variables!tricomplex}
which show that, upon multiplication, the transverse components obey the same
rules as the real and imaginary components of usual, two-dimensional complex
numbers, and the rule for the longitudinal component is that of the regular
multiplication of numbers.

If the constants in Eq. (\ref{3-90}) are $a_l=p_l+hq_l+kr_l$, and
\begin{equation}
a_{l1}=\frac{2p_l-q_l-r_l}{2}, \:\:\tilde a_{l1}=\frac{\sqrt{3}}{2}(q_l-r_l), 
\:\: a_{l+}=p_l+q_l+r_l,
\label{3-121}
\end{equation}
where $p_0=1, q_0=0, r_0=0$,
the series (\ref{3-90}) can be written as
\begin{equation}
\sum_{l=0}^\infty \left[a_{l1}e_1+\tilde a_{l1}\tilde e_1)(v_1e_1+\tilde v_1\tilde e_1)^l
+e_+ a_{l+} v_+^l\right].
\label{3-121b}
\end{equation}
The series in Eq. (\ref{3-121b}) is absolutely convergent for 
\begin{equation}
|v_+|<c_+,\:
(v_1^2+\tilde v_1^2)^{1/2}<c_1,
\label{3-n90}
\end{equation}\index{convergence, region of!tricomplex}
where 
\begin{equation}
c_+=\lim_{l\rightarrow\infty} \frac{|a_{l+}|}{|a_{l+1,u}|} ,\:
c_1=\lim_{l\rightarrow\infty} \frac
{\left(a_{l1}^2+ \tilde a_{l1}^2\right)^{1/2}}
{\left(a_{l+1,1}^2+ a_{l+1,2}^2\right)^{1/2}} .
\label{3-n91}
\end{equation}
The relations (\ref{3-n90}) and (\ref{3-118c}) show that the region of convergence
of the series (\ref{3-90}) is a cylinder
of radius $c_1\sqrt{2/3}$ and height $2c_+/\sqrt{3}$, having the
trisector line $(t)$ as axis and the origin as center, which can be called
cylinder of convergence, as shown in Fig. \ref{fig7}.

\begin{figure}
\begin{center}
\epsfig{file=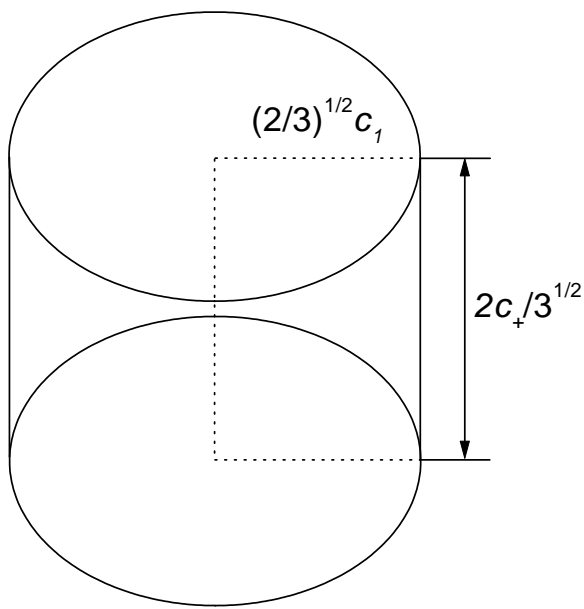,width=12cm}
\caption{Cylinder of convergence of tricomplex series, of radius
$ (2/3)^{1/2} c_1$ 
and height $2c_+/3^{1/2}$, having the axis parallel to the trisector line, and
the regions 
of space delimited by the plane, cylindrical and conical surfaces.}
\label{fig7}
\end{center}
\end{figure}

It can be shown that $c_1=(1/\sqrt{3})\;{\rm
min}(c_+,c_1)$, where ${\rm min}$ designates the smallest of
the numbers $c_+,c_1$. Using the expression of $|u|$ in
Eq. (\ref{3-118}), it can be seen that the spherical region of
convergence defined in Eqs. (\ref{3-new86}), (\ref{3-new87}) is a subset of the
cylindrical region of convergence defined in Eqs. (\ref{3-n90}) and (\ref{3-n91}).


\section{Analytic functions of tricomplex variables}

The analytic functions of the hypercomplex variable $u$ and the series 
expansion of functions have been discussed in Eqs. ~(1.85)-(1.93).
If the tricomplex function $f(u)$
of the tricomplex variable $u$ is written in terms of 
the real functions $F(x,y,z),G(x,y,z),H(x,y,z)$ of real variables $x,y,z$ as
\begin{equation}
f(u)=F(x,y,z)+hG(x,y,z)+kH(x,y,z), 
\label{3-98d}
\end{equation}\index{functions, real components!tricomplex}
then relations of equality 
exist between partial derivatives of the functions $F,G,H$. These relations can
be obtained by writing the derivative of the function $f$ as
\begin{eqnarray}
\lefteqn{\frac{1}{\Delta x+h\Delta y +k\Delta z} 
\left[\frac{\partial F}{\partial x}\Delta x+
\frac{\partial F}{\partial y}\Delta y+
\frac{\partial F}{\partial z}\Delta z
+h\left(\frac{\partial G}{\partial x}\Delta x+
\frac{\partial G}{\partial y}\Delta y+
\frac{\partial G}{\partial z}\Delta z\right) \right.}\nonumber\\
&&\left. 
+k\left(\frac{\partial H}{\partial x}\Delta x+
\frac{\partial H}{\partial y}\Delta y+
\frac{\partial H}{\partial z}\Delta z\right)\right] ,
\label{3-99}
\end{eqnarray}\index{derivative, independence of direction!tricomplex}
where the difference appearing in Eq. {(98)} is
$u-u_0=\Delta x+h\Delta y +k\Delta z$. The relations between 
the partials 
derivatives of the functions $F, G, H$ are obtained by setting successively in 
Eq. (\ref{3-99}) $\Delta x\rightarrow 0, \Delta y=0,  \Delta z=0$; then
$\Delta x= 0, \Delta y\rightarrow 0,  \Delta z=0$; and 
$\Delta x=0, \Delta y=0,  \Delta z\rightarrow 0$. The relations are
\begin{equation}
\frac{\partial F}{\partial x} = \frac{\partial G}{\partial y},\:\:
\frac{\partial G}{\partial x} = \frac{\partial H}{\partial y},\:\:
\frac{\partial H}{\partial x} = \frac{\partial F}{\partial y},
\label{3-100}
\end{equation}
\begin{equation}
\frac{\partial F}{\partial x} = \frac{\partial H}{\partial z},\:\:
\frac{\partial G}{\partial x} = \frac{\partial F}{\partial z},\:\:
\frac{\partial H}{\partial x} = \frac{\partial G}{\partial z},
\label{3-101}
\end{equation}
\begin{equation}
\frac{\partial G}{\partial y} = \frac{\partial H}{\partial z},\:\:
\frac{\partial H}{\partial y} = \frac{\partial F}{\partial z},\:\:
\frac{\partial F}{\partial y} = \frac{\partial G}{\partial z} .
\label{3-102}
\end{equation}\index{relations between partial derivatives!tricomplex}
The relations (\ref{3-100})-(\ref{3-102}) are analogous to the Riemann relations
for the real and imaginary components of a complex function. It can be shown
from Eqs. (\ref{3-100})-(\ref{3-102}) that the components $F$ solutions of the
equations 
\begin{equation}
\frac{\partial^2 F}{\partial x^2}-\frac{\partial^2 F}{\partial y\partial z}=0,
\:\: 
\frac{\partial^2 F}{\partial y^2}-\frac{\partial^2 F}{\partial x\partial z}=0,
\:\: 
\frac{\partial^2 F}{\partial z^2}-\frac{\partial^2 F}{\partial x\partial y}=0,
\label{3-103}
\end{equation}\index{relations between second-order derivatives!tricomplex}
\begin{equation}
\frac{\partial^2 G}{\partial x^2}-\frac{\partial^2 G}{\partial y\partial z}=0,
\:\: 
\frac{\partial^2 G}{\partial y^2}-\frac{\partial^2 G}{\partial x\partial z}=0,
\:\: 
\frac{\partial^2 G}{\partial z^2}-\frac{\partial^2 G}{\partial x\partial y}=0,
\label{3-104}
\end{equation}
\begin{equation}
\frac{\partial^2 H}{\partial x^2}-\frac{\partial^2 H}{\partial y\partial z}=0,
\:\: 
\frac{\partial^2 H}{\partial y^2}-\frac{\partial^2 H}{\partial x\partial z}=0,
\:\: 
\frac{\partial^2 H}{\partial z^2}-\frac{\partial^2 H}{\partial x\partial y}=0.
\label{3-105}
\end{equation}
It can also be shown that the
differences $F-G, F-H, G-H$ are solutions of the equation of Laplace,
\begin{equation}
\Delta (F-G)=0,\:\:\Delta (F-H)=0,\:\:\Delta
(G-H)=0,\:\:\Delta=\frac{\partial^2}{\partial x^2}+\frac{\partial^2}{\partial
y^2}+\frac{\partial^2}{\partial z^2} 
\label{3-106}
\end{equation}

If a geometric transformation is considered in which to a point $u$ is
associated the point $f(u)$, 
it can be shown that the tricomplex function $f(u)$ transforms a straight line
parallel to the trisector line $(t)$ in a straight line parallel to $(t)$, and
transforms a plane parallel to the nodal plane $\Pi$ in a plane parallel to
$\Pi$. 
A transformation generated by a tricomplex function $f(u)$ does not conserve in
general the angle of intersecting lines.

\section{Integrals of tricomplex functions}

The singularities of tricomplex functions arise from terms of the form
$1/(u-a)^m$, with $m>0$. Functions containing such terms are singular not
only at $u=a$, but also at all points of a plane $(\Pi_a)$
through the point $a$ and
parallel to the nodal plane $\Pi$ and at all points of a
straight line $(t_a)$ 
passing through $a$ and parallel to the trisector line $(t)$.

The integral of a tricomplex function between two points $A, B$ along a path
situated in a region free of singularities is independent of path, which means
that the integral of an analytic function along a loop situated in a region
free from singularities is zero,
\begin{equation}
\oint_\Gamma f(u) du = 0,
\label{3-107}
\end{equation}
where it is supposed that a surface $S$ spanning the 
closed loop $\Gamma$ is not intersected by any of
the planes and is not threaded by any of the lines associated with the
singularities of the function $f(u)$. Using the expression, Eq. (\ref{3-98d})
for $f(u)$ and the fact that $du=dx+hdy+kdz$, the explicit form of the 
integral in Eq. (\ref{3-107}) is
\begin{equation}
\oint _\Gamma f(u) du = \oint_\Gamma
[Fdx+Hdy+Gdz+h(Gdx+Fdy+Hdz)+k(Hdx+Gdy+Fdz)] .
\label{3-108}
\end{equation}\index{integrals, path!tricomplex}
If the functions $F, G, H$ are regular on a surface $S$
spanning the loop $\Gamma$,
the integral along the loop $\Gamma$ can be transformed with the aid of the
theorem of Stokes in an integral over the surface $S$ of terms of the form
$\partial H/\partial x -  \partial F/\partial y, \:\:
\partial G/\partial x - \partial F/\partial z,\:\:
\partial G/\partial y - \partial H/\partial z, \dots$
which are equal to zero by Eqs. (\ref{3-100})-(\ref{3-102}), and this proves Eq.
(\ref{3-107}). 

If there are singularities on the surface $S$, the integral $\oint f(u) du$
is not necessarily equal to zero. If $f(u)=1/(u-a)$ and the loop $\Gamma_a$ is
situated in the half-space above the plane $(\Pi_a)$ and encircles once
the line $(t_a)$, then 
\begin{equation}
\oint_{\Gamma_a} \frac{du}{u-a}=\frac{2\pi}{\sqrt{3}}(h-k) .
\label{3-109}
\end{equation}\index{poles and residues!tricomplex}
This is due to the fact that the integral of $1/(u-a)$ 
along the loop $\Gamma_a$ is equal to the integral of $1/(u-a)$ along a circle
$(C_a)$ with the center on the line $(t_a)$ and perpendicular to this line,
as shown in Fig. \ref{fig8}. 
\begin{equation}
\oint_{\Gamma_a} \frac{du}{u-a}=\oint_{C_a} \frac{du}{u-a} ,
\label{3-110}
\end{equation}
this being a corrolary of Eq. (\ref{3-107}). The integral on the right-hand side
of Eq. (\ref{3-110}) can be evaluated with the aid of the trigonometric form 
Eq. (\ref{3-63}) of the tricomplex quantity $u-a$, so that
\begin{equation}
\frac{du}{u-a}=\frac{h-k}{\sqrt{3}}d\phi ,
\label{3-111}
\end{equation}
which by integration over $d\phi$ from 0 to $2\pi$ yields Eq. (\ref{3-109}).

\begin{figure}
\begin{center}
\epsfig{file=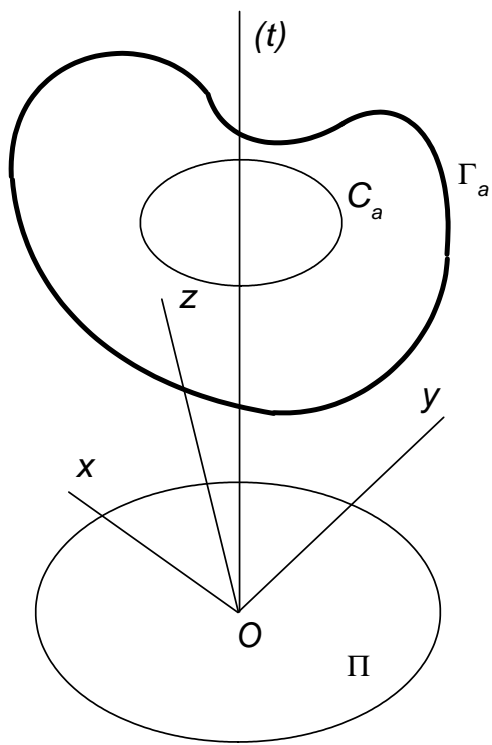,width=12cm}
\caption{The integral of $1/(u-a)$ 
along the loop $\Gamma_a$ is equal to the integral of $1/(u-a)$ along a circle
$(C_a)$ with the center on the line $(t_a)$ and perpendicular to this line.}
\label{fig8}
\end{center}
\end{figure}

The integral  $\oint_{\Gamma_a} du(u-a)^m$, 
with $m$ an integer number not equal to -1, is 
equal to zero, because $\int du (u-a)^m=(u-a)^{m+1}/(m+1)$, and 
$(u-a)^{m+1}/(m+1)$ is singlevalued,
\begin{equation}
\oint_{\Gamma_a}du(u-a)^m=0, \:\:{\rm for}\: m\: {\rm integer},
\: m\not=-1. 
\label{3-112}
\end{equation}
If $f(u)$ is an analytic tricomplex function which can be expanded in a
series as written in Eq. (1.89), with $u_0=a$, 
and the expansion holds on the curve
$\Gamma$ and on a surface spanning $\Gamma$, then from Eqs. (\ref{3-111}) and
(\ref{3-112}) it follows that
\begin{equation}
\oint_\Gamma \frac{f(u)du}{u-a}=\frac{2\pi}{\sqrt{3}}(h-k)f(a) .
\label{3-113}
\end{equation}
Substituting in the right-hand side of 
Eq. (\ref{3-113}) the expression of $f(u)$ in terms of the real 
components $F, G, H$,  Eq. (\ref{3-98d}), at $u=a$, yields
\begin{equation}
\oint_\Gamma \frac{f(u)du}{u-a}=\frac{2\pi}{\sqrt{3}}[H-G+h(F-H)+k(G-F)] .
\label{3-114}
\end{equation}
Since the sum of the real components in the paranthesis from the right-hand
side of Eq. (\ref{3-114}) is equal to zero, this equation
determines only the differences between the components $F, G, H$. 
If $f(u)$ can be expanded as written in Eq. (1.89), with $u_0=a$, on 
$\Gamma$ and on a surface spanning $\Gamma$, then from Eqs. (\ref{3-109}) and
(\ref{3-112}) it also results that
\begin{equation}
\oint_\Gamma \frac{f(u)du}{(u-a)^{m+1}}
=\frac{2\pi}{\sqrt{3}m!}(h-k)f^{(m)}(a) ,
\label{3-115}
\end{equation}
where the fact that has been used that the derivative $f^{(m)}(a)$ of order $m$
of $f(u)$ at $u=a$ is related to the expansion coefficient in 
Eq. (1.89), with $u_0=a$,
according to Eq. (1.93), with $u_0=a$. The relation (\ref{3-115}) can also be obtained by
successive derivations of Eq. (\ref{3-113}).

If a function $f(u)$ is expanded in positive and negative powers of $u-u_j$,
where $u_j$ are fourcomplex constants, $j$ being an index, the integral of $f$
on a closed loop $\Gamma$ is determined by the terms in the expansion of $f$
which are of the form $a_j/(u-u_j)$,
\begin{equation}
f(u)=\cdots+\sum_j\frac{a_j}{u-u_j}+\cdots .
\label{3-115b}
\end{equation}
In Eq. (\ref{3-115b}), $u_j$ is the pole and $a_j$ is the residue relative to the
pole $u_j$.
Then the integral of $f$ on a closed loop $\Gamma$ is
\begin{equation}
\oint_\Gamma f(u) du = 
\frac{2\pi}{\sqrt{3}}(h-k) \sum_j \:\;
{\rm int}(u_{j\Pi},\Gamma_{\Pi}) a_j ,
\label{3-115c}
\end{equation}
where the functional int($M,C$), defined for a point $M$ and a
closed curve $C$ in a two-dimensional plane, is given by 
\begin{equation}
{\rm int}(M,C)=\left\{
\begin{array}{l}
1 \;\:{\rm if} \;\:M \;\:{\rm is \;\:an \;\:interior \;\:point \;\:of} \;\:C ,\\ 
0 \;\:{\rm if} \;\:M \;\:{\rm is \;\:exterior \;\:to}\:\; C,\\
\end{array}\right. 
\label{3-115d}
\end{equation}
and $u_{j\Pi},\Gamma_{\Pi}$ are the projections of the point $u_j$ and
of the curve $\Gamma$ on the nodal plane $\Pi$, as shown in Fig. \ref{fig9}.

\begin{figure}
\begin{center}
\epsfig{file=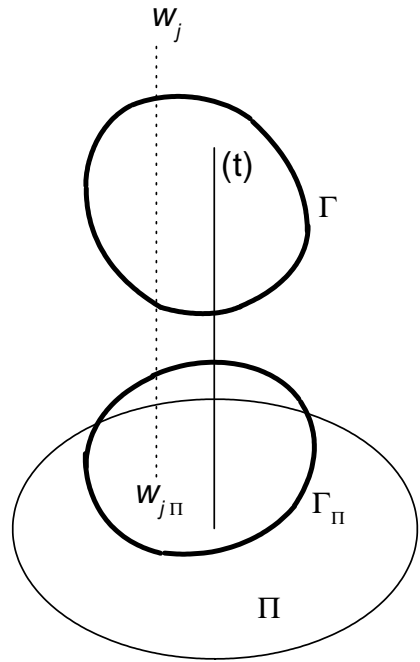,width=12cm}
\caption{Integration path $\Gamma$, pole $u_j$ and their projections
$\Gamma_{\Pi}, u_{j\Pi}$ on the nodal plane $\Pi$.}
\label{fig9}
\end{center}
\end{figure}


\section{Factorization of tricomplex polynomials}

A polynomial of degree $m$ of the tricomplex variable $u=x+hy+kz$ has the form
\begin{equation}
P_m(u)=u^m+a_1 u^{m-1}+\cdots+a_{m-1} u +a_m ,
\label{3-115e}
\end{equation}
where the constants are in general tricomplex numbers, $a_l=p_l+hq_l+kr_l$,
$l=1,\cdots,m$.  
In order to write the polynomial $P_m(u)$ as a product of factors, the variable
$u$ and the constants $a_l$ will be written in the form which explicits the
transverse and longitudinal components, 
\begin{equation}
P_m(u)=\sum_{l=0}^m (a_{l1}e_1+\tilde a_{l1}\tilde e_1)(v_1e_1+\tilde v_1\tilde e_1)^{m-l}
+e_+\sum_{l=0}^m a_{l+} v_+^{m-l} ,
\label{3-120}
\end{equation}\index{polynomial, canonical variables!tricomplex}
where the constants have been defined previously in Eq. (\ref{3-121}).
Due to the properties in Eq. (\ref{3-118}), the transverse part of the
polynomial $P_m(u)$ can be written as a product of linear factors of the form
\begin{equation}
\sum_{l=0}^m (a_{l1}e_1+\tilde a_{l1}\tilde e_1)(v_1e_1+\tilde v_1\tilde e_1)^{m-l}=
\prod_{l=1}^m [(v_1-v_{l1})e_1+(\tilde v_1-\tilde v_{l1})\tilde e_1] ,
\label{3-122}
\end{equation}
where the quantities $v_{l1}, \tilde v_{l1}$ are real numbers.
The longitudinal part of $P_m(u)$, Eq. (\ref{3-120}), can be written as a product
of linear or quadratic factors with real coefficients, or as a product of
linear factors which, if imaginary, appear always in complex conjugate pairs.
Using the latter form for the simplicity of notations, the longitudinal part
can be written as
\begin{equation}
\sum_{l=0}^m a_{l+} v_+^{m-l}=\prod_{l=1}^m (v_+-v_{l+}) ,
\label{3-123}
\end{equation}
where the quantities $v_{l+}$ appear always in complex conjugate pairs.
Due to the orthogonality of the transverse and longitudinal components, Eq.
(\ref{3-118}), the polynomial $P_m(u)$ can be written as a product of factors of
the form  
\begin{equation}
P_m(u)=\prod_{l=1}^m [(v_1-v_{l1})e_1+(\tilde v_1-\tilde v_{l1})\tilde e_1
+ (v_+-v_{l+})e_+] .
\label{3-124}
\end{equation}\index{polynomial, factorization!tricomplex}
These relations can be written with the aid of Eqs. (\ref{3-116}) as
\begin{eqnarray}
P_m(u)=\prod_{l=1}^m (u-u_l),
\label{3-125}
\end{eqnarray}
where
\begin{eqnarray}
u_l=v_{l1}e_1+\tilde v_{l1}\tilde e_1+ v_{l+}e_+ .
\label{3-126}
\end{eqnarray}
The roots $v_{l+}$ and the roots $v_{l1}e_1+\tilde v_{l1}\tilde e_1$
defined in Eq. (\ref{3-122}) may be ordered arbitrarily.
This means that Eq. (\ref{3-126}) gives sets of $m$ roots
$u_1,...,u_m$ of the polynomial $P_m(u)$, 
corresponding to the various ways in which the roots 
$v_{l+}, v_{l1}e_1+\tilde v_{l1}\tilde e_1$ are ordered according to $l$ in
each group. Thus, while the tricomplex components in Eq. (\ref{3-120}) taken
separately have unique factorizations, the polynomial $P_m(q)$ can be written
in many different ways as a product of linear factors. 

If $P(u)=u^2-1$, the degree is $m=2$, the coefficients of the polynomial are
$a_1=0, a_2=-1$, the coefficients defined in Eq. (\ref{3-121}) are 
$a_{21}=-1, a_{22}=0, a_{2u}=-1$. 
The expression of $P(u)$, Eq. (\ref{3-120}), is  
$(e_1v_1+\tilde e_1\tilde v_1)^2-e_1+e_+ (v_+^2-1)$. The factorizations in Eqs. (\ref{3-122})
and (\ref{3-123}) are
$(e_1v_1+\tilde e_1\tilde v_1)^2-e_1=[e_1(v_1+1)+\tilde e_1\tilde v_1][e_1(v_1-1)+\tilde e_1\tilde v_1]$ and 
$v_+^2-1=(v_++1)(v_+-1)$. The factorization of $P(u)$, Eq. (\ref{3-125}), is
$P(u)=(u-u_1)(u-u_2)$, where according to Eq. (\ref{3-126}) the roots are
$u_1=\pm e_1\pm  e_+, u_2=-u_1$. If $e_1$ and $e_+$ are expressed with the
aid of Eq. (\ref{3-117b}) in terms of $h$ and $k$, the factorizations of $P(u)$
are obtained as
\begin{equation}
u^2-1=(u+1)(u-1) ,
\label{3-127}
\end{equation}
or as
\begin{equation}
u^2-1=\left(u+\frac{1-2h-2k}{3}\right) 
\left(u-\frac{1-2h-2k}{3}\right) .
\label{3-128}
\end{equation}
It can be checked that 
$(\pm e_1\pm  e_+)^2=e_1+e_+=1$.

\section{Representation of tricomplex complex numbers by irreducible
matrices} 

If the matrix in Eq. (\ref{3-25}) representing the tricomplex number $u$ is
called $U$, and 
\begin{equation}
T=\left(
\begin{array}{ccc}
\sqrt{\frac{2}{3}}     &  -\frac{1}{\sqrt{6}} & -\frac{1}{\sqrt{6}} \\
0                      &  \frac{1}{\sqrt{2}}   & -\frac{1}{\sqrt{2}} \\
\frac{1}{\sqrt{3}}     & \frac{1}{\sqrt{3}}   & \frac{1}{\sqrt{3}}   \\
\end{array}\right),
\label{3-129}
\end{equation}
which is the matrix appearing in Eq. (\ref{3-13}), it can be checked that 
\begin{equation}
TUT^{-1}=\left(
\begin{array}{ccc}
x-\frac{y+z}{2}          & \frac{\sqrt{3}}{2}(y-z) & 0  \\
-\frac{\sqrt{3}}{2}(y-z) &x-\frac{y+z}{2}          & 0  \\
  0                      & 0                       &x+y+z          
\end{array}\right).
\label{3-130}
\end{equation}\index{representation by irreducible matrices}
The relations for the variables $ x-(y+z)/2, (\sqrt{3}/2)(y-z)$ and 
$x+y+z$ for the multiplication of tricomplex numbers have been written in
Eqs. (\ref{3-18}), (\ref{3-20}) and (\ref{3-21}). The matrices
$T U T^{-1}$  provide an irreducible representation
\cite{4} of the tricomplex numbers $u=x+hy+kz$, in terms of matrices with real
coefficients.

\chapter{Commutative Complex Numbers in Four Dimensions}

Systems of hypercomplex numbers in 4 dimensions of the form 
$u=x+\alpha y+\beta z+\gamma t$ are described in this chapter,
where the variables $x, y,z$ and $t$ are real numbers,
for which the multiplication is both associative and commutative.
The product of two
fourcomplex numbers is equal to zero if both numbers are 
equal to zero, or if the numbers belong to certain four-dimensional hyperplanes
as discussed further in this chapter. 
The fourcomplex numbers have exponential and trigonometric forms,
and the concepts of analytic fourcomplex 
function,  contour integration and residue can be defined. 
Expressions are given for the elementary functions of fourcomplex variable.
The functions $f(u)$ of fourcomplex variable defined by power series, 
have derivatives $\lim_{u\rightarrow u_0} [f(u)-f(u_0)]/(u-u_0)$ independent of
the direction of approach of $u$ to $u_0$. If the fourcomplex function $f(u)$
of the fourcomplex variable $u$ is written in terms of 
the real functions $P(x,y,z,t),Q(x,y,z,t),R(x,y,z,t),S(x,y,z,t)$, then
relations of equality  
exist between partial derivatives of the functions $P,Q,R,S$. 
The integral $\int_A^B f(u) du$ of a fourcomplex
function between two points $A,B$ is independent of the path connecting $A,B$.

Four distinct types
of hypercomplex numbers are studied, as discussed further.
In Sec. 3.1, the multiplication rules for the complex units $\alpha, \beta$ and
$\gamma$ are 
$\alpha^2=-1, \:\beta^2=-1, \:\gamma^2=1, \alpha\beta=\beta\alpha=-\gamma,\:
\alpha\gamma=\gamma\alpha=\beta, \:\beta\gamma=\gamma\beta=\alpha$.
The exponential form of a fourcomplex number is
$u=\rho\exp\left[\gamma \ln\tan\psi+
(1/2)\alpha(\phi+\chi)+(1/2)\beta(\phi-\chi) 
\right]$ , where the amplitude is
$\rho^4=\left[(x+t)^2+(y+z)^2\right]\left[(x-t)^2+(y-z)^2\right]$ ,
$\phi, \chi$ are azimuthal angles, $0\leq\phi<2\pi, 0\leq\chi<2\pi$, 
and $\psi$ is a planar angle, $0<\psi\leq\pi/2$.
The trigonometric form of a fourcomplex number is
$u=d[\cos(\psi-\pi/4)+\gamma\sin(\psi-\pi/4)]$
$\exp\left[(1/2)\alpha(\phi+\chi)+(1/2)\beta(\phi-\chi)\right]$, 
where $d^2=x^2+y^2+z^2+t^2$.
The amplitude $\rho$ and $\tan\psi$ are multiplicative 
and the angles $\phi, \chi$ 
are additive upon the multiplication of fourcomplex numbers.
Since there are two cyclic variables, $\phi$ and $\chi$, these fourcomplex
numbers are called circular fourcomplex numbers.
If $f(u)$ is an analytic fourcomplex function, then 
$\oint_\Gamma f(u)du/(u-u_0)=
\pi[(\alpha+\beta) \:\;{\rm int}(u_{0\xi\upsilon},\Gamma_{\xi\upsilon})
+ (\alpha-\beta)\:\;{\rm int}(u_{0\tau\zeta},\Gamma_{\tau\zeta})]\;f(u_0)$,
where the functional ${\rm int}$ takes the values 0 or 1 depending on the
relation between the projections of the point $u_0$ and of the curve $\Gamma$
on certain planes.
A polynomial  can be written as a 
product of linear or quadratic factors, although the factorization may not be
unique.

In Sec. 3.2, the multiplication rules for the complex units $\alpha, \beta$ and
$\gamma$ are
$\alpha^2=1, \:\beta^2=1, \:\gamma^2=1, \alpha\beta=\beta\alpha=\gamma,\:
\alpha\gamma=\gamma\alpha=\beta, \:\beta\gamma=\gamma\beta=\alpha$. 
The exponential form of a fourcomplex number, which can be defined for
$s=x+y+z+t>0, \: s^\prime= x-y+z-t>0 , \: s^{\prime\prime}=x+y-z-t>0,  \:
s^{\prime\prime\prime}=x-y-z+t>0$, is
$u=
\mu\exp(\alpha y_1+\beta z_1+\gamma t_1)$, where the amplitude is
$\mu=(s s^\prime s^{\prime\prime}s^{\prime\prime\prime})^{1/4}$ and the
arguments are
$y_1=(1/4)\ln(ss^{\prime\prime}/s^\prime s^{\prime\prime\prime}),\:
z_1=(1/4)\ln(ss^{\prime\prime\prime}/s^\prime s^{\prime\prime}),\:
t_1=(1/4)\ln(ss^\prime/s^{\prime\prime}s^{\prime\prime\prime})$.
Since there is no cyclic variable, these fourcomplex numbers are called
hyperbolic fourcomplex numbers. 
The amplitude $\mu$ is multiplicative and the arguments $y_1, z_1, t_1$
are additive upon the multiplication of fourcomplex numbers.
A polynomial can be
written as a 
product of linear or quadratic factors, although the factorization may not be
unique. 

In Sec. 3.3, the multiplication rules for the complex units $\alpha, \beta$ and
$\gamma$ are
$\alpha^2=\beta, \:\beta^2=-1, \:\gamma^2=-\beta,
\alpha\beta=\beta\alpha=\gamma,\: 
\alpha\gamma=\gamma\alpha=-1, \:\beta\gamma=\gamma\beta=-\alpha$ .
The exponential function of a fourcomplex number can be expanded in terms of
the four-dimensional cosexponential functions
$f_{40}(x)=1-x^4/4!+x^8/8!-\cdots ,\:\:
f_{41}(x)=x-x^5/5!+x^9/9!-\cdots,\:\:
f_{42}(x)=x^2/2!-x^6/6!+x^{10}/10!-\cdots ,\:\:
f_{43}(x)=x^3/3!-x^7/7!+x^{11}/11!-\cdots$. Expressions are obtained for the
four-dimensional cosexponential functions in terms of elementary functions.
The exponential form of a fourcomplex number is
$u=\rho\exp\left\{(1/2)(\alpha-\gamma)\ln\tan\psi\right.$
$\left.+(1/2)[\beta+(\alpha+\gamma)/\sqrt{2}]\phi
-(1/2)[\beta-(\alpha+\gamma)/\sqrt{2}]\chi
\right\}$, where the amplitude is \\
$\rho^4=\left\{\left[x+(y-t)/\sqrt{2}\right]^2+\left[z+(y+t)/\sqrt{2}\right]^2\right\}$
$\left\{ \left[x-(y-t)/\sqrt{2}\right]^2+\left[z-(y+t)/\sqrt{2}\right]^2\right\}$, \\
$\phi, \chi$ are azimuthal angles, $0\leq\phi<2\pi, 0\leq\chi<2\pi$, 
and $\psi$ is a planar angle, $0\leq\psi\leq\pi/2$.
The trigonometric form of a fourcomplex number is
$u=d\left[\cos\left(\psi-\pi/4\right)\right.$
$\left.+(1/\sqrt{2})(\alpha-\gamma)\sin\left(\psi-\pi/4\right)\right]$
$\exp\left\{(1/2)[\beta+(\alpha+\gamma)/\sqrt{2}]\phi\right.$
$\left.-(1/2)[\beta-(\alpha+\gamma)/\sqrt{2}]\chi\right\}$,
where $d^2=x^2+y^2+z^2+t^2$.
The amplitude $\rho$ and $\tan\psi$ are multiplicative 
and the angles $\phi, \chi$ 
are additive upon the multiplication of fourcomplex numbers.
There are two cyclic variables, $\phi$ and $\chi$, so that these fourcomplex
numbers are also of a circular type. In order to distinguish
them from the circular hypercomplex numbers, these are called
planar fourcomplex numbers.
If $f(u)$ is an analytic fourcomplex function, then 
$\oint_\Gamma f(u)du/(u-u_0)=
\pi\left[\left(\beta+(\alpha+\gamma)/\sqrt{2}\right) \;{\rm
int}(u_{0\xi\upsilon},\Gamma_{\xi\upsilon}) \right.$
$\left.+ \left(\beta-(\alpha+\gamma)/\sqrt{2}\right)\;{\rm
int}(u_{0\tau\zeta},\Gamma_{\tau\zeta})\right]\;f(u_0)$ , 
where the functional ${\rm int}$ takes the values 0 or 1 depending on the
relation between the projections of the point $u_0$ and of the curve $\Gamma$
on certain planes.
A polynomial can be
written as a 
product of linear or quadratic factors, although the factorization may not be
unique. 
The fourcomplex numbers described in this chapter are a particular case for 
$n=4$ of the planar hypercomplex numbers in $n$ dimensions discussed in Sec. 6.2.

In Sec. 3.4, the multiplication rules for the complex units $\alpha, \beta$ and
$\gamma$ are
$\alpha^2=\beta, \:\beta^2=1, \:\gamma^2=\beta,
\alpha\beta=\beta\alpha=\gamma,\: 
\alpha\gamma=\gamma\alpha=1, \:\beta\gamma=\gamma\beta=\alpha$.
The product of two fourcomplex numbers is equal to zero if both numbers are
equal to zero, or if the numbers belong to certain four-dimensional hyperplanes
described further in this section. 
The exponential function of a fourcomplex number can be expanded in terms of
the four-dimensional cosexponential functions
$g_{40}(x)=1+x^4/4!+x^8/8!+\cdots ,\:\:
g_{41}(x)=x+x^5/5!+x^9/9!+\cdots,\:\:
g_{42}(x)=x^2/2!+x^6/6!+x^{10}/10!+\cdots ,\:\:
g_{43}(x)=x^3/3!+x^7/7!+x^{11}/11!+\cdots$. Addition theorems and other
relations are obtained for these four-dimensional cosexponential functions.
The exponential form of a fourcomplex number, which can be defined for
$x+y+z+t>0, x-y+z-t>0$, is
$u=\rho\exp\left[(1/4)(\alpha+\beta+\gamma) \ln(\sqrt{2}/\tan\theta_+)
-(1/4)(\alpha-\beta+\gamma) \ln(\sqrt{2}/\tan\theta_-)
+(\alpha-\gamma)\phi/2\right]$,\\
where $\rho=(\mu_+\mu_-)^{1/2},\:\: \mu_+^2=(x-z)^2+(y-t)^2,\:\:
\mu_-^2=(x+z)^2-(y+t)^2$, $e_+=(1+\alpha+\beta+\gamma)/4,
e_-=(1-\alpha+\beta-\gamma)/4$,
$e_1=(1-\beta)/2, \tilde e_1=(\alpha-\gamma)/2$,
the polar angles are 
$\tan\theta_+=\sqrt{2}\mu_+/v_+,  
\tan\theta_-=\sqrt{2}\mu_+/v_-$, $0\leq\theta_+\leq\pi, 
0\leq\theta_-\leq\pi$,
and the azimuthal angle $\phi$ is defined by the relations 
$x-y=\mu_+\cos\phi,\:\:z-t=\mu_+\sin\phi$, $0\leq\phi<2\pi$. 
The trigonometric form of the fourcomplex number $u$ is
$u=d\sqrt{2}$
$\left(1+1/\tan^2\theta_++1/\tan^2\theta_-\right)^{-1/2}
\left\{e_1+e_+\sqrt{2}/\tan\theta_+
+e_-\sqrt{2}/\tan\theta_-\right\}$
$\exp\left[\tilde e_1\phi\right]$.
The amplitude $\rho$ and $\tan\theta_+/\sqrt{2}, \tan\theta_-/\sqrt{2}$,
are multiplicative, and the azimuthal angle $\phi$ is additive
upon the multiplication of fourcomplex numbers.
There is only one cyclic variable, $\phi$, 
and there are two axes $v_+, v_-$ which play an important role in the
description of these numbers, so that
these hypercomplex numbers are called polar fourcomplex numbers.
If $f(u)$ is an analytic fourcomplex function, then 
$\oint_\Gamma f(u)du/(u-u_0)=
\pi(\beta-\gamma) \;{\rm int}(u_{0\xi\upsilon},\Gamma_{\xi\upsilon})f(u_0)$ ,
where the functional ${\rm int}$ takes the values 0 or 1 depending on the
relation between the projections of the point $u_0$ and of the curve $\Gamma$
on certain planes.
A polynomial can be
written as a 
product of linear or quadratic factors, although the factorization may not be
unique.
The fourcomplex numbers described in this section are a particular case for 
$n=4$ of the polar hypercomplex numbers in $n$ dimensions discussed in Sec. 6.1.

\section{Circular complex numbers in four dimensions}

\subsection{Operations with circular fourcomplex numbers}

A circular fourcomplex number is determined by its four components $(x,y,z,t)$.
The sum of the circular fourcomplex numbers $(x,y,z,t)$ and
$(x^\prime,y^\prime,z^\prime,t^\prime)$ is the circular fourcomplex
number $(x+x^\prime,y+y^\prime,z+z^\prime,t+t^\prime)$. \index{sum!circular
fourcomplex} 
The product of the circular fourcomplex numbers
$(x,y,z,t)$ and $(x^\prime,y^\prime,z^\prime,t^\prime)$ 
is defined in this section to be the circular fourcomplex
number
$(xx^\prime-yy^\prime-zz^\prime+tt^\prime,
xy^\prime+yx^\prime+zt^\prime+tz^\prime,
xz^\prime+zx^\prime+yt^\prime+ty^\prime,
xt^\prime+tx^\prime-yz^\prime-zy^\prime)$.\index{product!circular fourcomplex}

Circular fourcomplex numbers and their operations can be represented by  writing the
circular fourcomplex number $(x,y,z,t)$ as  
$u=x+\alpha y+\beta z+\gamma t$, where $\alpha, \beta$ and $\gamma$ 
are bases for which the multiplication rules are 
\begin{equation}
\alpha^2=-1, \:\beta^2=-1, \:\gamma^2=1, \alpha\beta=\beta\alpha=-\gamma,\:
\alpha\gamma=\gamma\alpha=\beta, \:\beta\gamma=\gamma\beta=\alpha .
\label{1}
\end{equation}\index{complex units!circular fourcomplex}
Two circular fourcomplex numbers $u=x+\alpha y+\beta z+\gamma t, 
u^\prime=x^\prime+\alpha y^\prime+\beta z^\prime+\gamma t^\prime$ are equal, 
$u=u^\prime$, if and only if $x=x^\prime, y=y^\prime,
z=z^\prime, t=t^\prime$. 
If 
$u=x+\alpha y+\beta z+\gamma t, 
u^\prime=x^\prime+\alpha y^\prime+\beta z^\prime+\gamma t^\prime$
are circular fourcomplex numbers, 
the sum $u+u^\prime$ and the 
product $uu^\prime$ defined above can be obtained by applying the usual
algebraic rules to the sum 
$(x+\alpha y+\beta z+\gamma t)+ 
(x^\prime+\alpha y^\prime+\beta z^\prime+\gamma t^\prime)$
and to the product 
$(x+\alpha y+\beta z+\gamma t)
(x^\prime+\alpha y^\prime+\beta z^\prime+\gamma t^\prime)$,
and grouping of the resulting terms,
\begin{equation}
u+u^\prime=x+x^\prime+\alpha(y+y^\prime)+\beta(z+z^\prime)+\gamma(t+t^\prime),
\label{1a}
\end{equation}\index{sum!circular fourcomplex}
\begin{eqnarray}
\lefteqn{uu^\prime=
xx^\prime-yy^\prime-zz^\prime+tt^\prime+
\alpha(xy^\prime+yx^\prime+zt^\prime+tz^\prime)+
\beta(xz^\prime+zx^\prime+yt^\prime+ty^\prime)}\nonumber\\
&&+\gamma(xt^\prime+tx^\prime-yz^\prime-zy^\prime).
\label{1b}
\end{eqnarray}\index{product!circular fourcomplex}

If $u,u^\prime,u^{\prime\prime}$ are circular fourcomplex numbers, the multiplication 
is associative
\begin{equation}
(uu^\prime)u^{\prime\prime}=u(u^\prime u^{\prime\prime})
\label{2}
\end{equation}
and commutative
\begin{equation}
u u^\prime=u^\prime u ,
\label{3}
\end{equation}
as can be checked through direct calculation.
The circular fourcomplex zero is $0+\alpha\cdot 0+\beta\cdot 0+\gamma\cdot 0,$ 
denoted simply 0, 
and the circular fourcomplex unity is $1+\alpha\cdot 0+\beta\cdot 0+\gamma\cdot 0,$ 
denoted simply 1.

The inverse of the circular fourcomplex number 
$u=x+\alpha y+\beta z+\gamma t$ is a circular fourcomplex number
$u^\prime=x^\prime+\alpha y^\prime+\beta z^\prime+\gamma t^\prime$
having the property that
\begin{equation}
uu^\prime=1 .
\label{4}
\end{equation}
Written on components, the condition, Eq. (\ref{4}), is
\begin{equation}
\begin{array}{c}
xx^\prime-yy^\prime-zz^\prime+tt^\prime=1,\\
yx^\prime+xy^\prime+tz^\prime+zt^\prime=0,\\
zx^\prime+ty^\prime+xz^\prime+yt^\prime=0,\\
tx^\prime-zy^\prime-yz^\prime+xt^\prime=0 .
\end{array}
\label{5}
\end{equation}
The system (\ref{5}) has the solution
\begin{equation}
x^\prime=\frac{x(x^2+y^2+z^2-t^2)-2yzt}
{\rho^4} ,
\label{6a}
\end{equation}\index{inverse!circular fourcomplex}
\begin{equation}
y^\prime=
\frac{y(-x^2-y^2+z^2-t^2)+2xzt}
{\rho^4} ,
\label{6b}
\end{equation}
\begin{equation}
z^\prime=\frac{z(-x^2+y^2-z^2-t^2)+2xyt}
{\rho^4} ,
\label{6c}
\end{equation}
\begin{equation}
t^\prime=\frac{t(-x^2+y^2+z^2+t^2)-2xyz}
{\rho^4} ,
\label{6d}
\end{equation}
provided that $\rho\not=0, $ where
\begin{equation}
\rho^4=x^4+y^4+z^4+t^4+2(x^2y^2+x^2z^2-x^2t^2-y^2z^2+y^2t^2+z^2t^2)-8xyzt .
\label{6e}
\end{equation}\index{amplitude!circular fourcomplex}
The quantity $\rho$ will be called amplitude of the circular fourcomplex number
$x+\alpha y+\beta z +\gamma t$.
Since
\begin{equation}
\rho^4=\rho_+^2\rho_-^2 ,
\label{7a}
\end{equation}
where
\begin{equation}
\rho_+^2=(x+t)^2+(y+z)^2, \: \rho_-^2=(x-t)^2+(y-z)^2 ,
\label{7b}
\end{equation}
a circular fourcomplex number $u=x+\alpha y+\beta z+\gamma t$ has an inverse, unless
\begin{equation}
x+t=0,\: y+z=0 ,  
\label{8}
\end{equation}
or
\begin{equation}
x-t=0,\: y-z=0 .  
\label{9}
\end{equation}

Because of conditions (\ref{8})-(\ref{9}) these 2-dimensional surfaces
will be called nodal hyperplanes.\index{nodal hyperplanes!circular fourcomplex} 
It can be shown that 
if $uu^\prime=0$ then either $u=0$, or $u^\prime=0$, 
or one of the circular fourcomplex numbers is of the form $x+\alpha y+\beta y+\gamma x$
and the other of the form $x^\prime+\alpha y^\prime-\beta y^\prime -\gamma
x^\prime$.

\subsection{Geometric representation of circular fourcomplex numbers}

The circular fourcomplex number $x+\alpha y+\beta z+\gamma t$ can be represented by 
the point $A$ of coordinates $(x,y,z,t)$. 
If $O$ is the origin of the four-dimensional space $x,y,z,t,$ the distance 
from $A$ to the origin $O$ can be taken as
\begin{equation}
d^2=x^2+y^2+z^2+t^2 .
\label{10}
\end{equation}\index{distance!circular fourcomplex}
The distance $d$ will be called modulus of the circular fourcomplex number $x+\alpha
y+\beta z +\gamma t$, $d=|u|$. 
The orientation in the four-dimensional space of the line $OA$ can be specified
with the aid of three angles $\phi, \chi, \psi$
defined with respect to the rotated system of axes
\begin{equation}
\xi=\frac{x+t}{\sqrt{2}}, \: \tau=\frac{x-t}{\sqrt{2}}, \:
\upsilon=\frac{y+z}{\sqrt{2}}, \: \zeta=\frac{y-z}{\sqrt{2}} .
\label{11}
\end{equation}
The variables $\xi, \upsilon, \tau, \zeta$ will be called canonical 
circular fourcomplex variables.\index{canonical variables!circular fourcomplex}
The use of the rotated axes $\xi, \upsilon, \tau, \zeta$ 
for the definition of the angles $\phi, \chi, \psi$ 
is convenient for the expression of the circular fourcomplex numbers
in exponential and trigonometric forms, as it will be discussed further.
The angle $\phi$ is the angle between the projection of $A$ in the plane
$\xi,\upsilon$ and the $O\xi$ axis, $0\leq\phi<2\pi$,  
$\chi$ is the angle between the projection of $A$ in the plane $\tau,\zeta$ and
the  $O\tau$ axis, $0\leq\chi<2\pi$,
and $\psi$ is the angle between the line $OA$ and the plane $\tau O \zeta
$, $0\leq \psi\leq\pi/2$, 
as shown in Fig. \ref{fig10}. The angles $\phi$ and $\chi$ will be called azimuthal
angles, the angle $\psi$ will be called planar angle .
\index{azimuthal angles!circular fourcomplex}
\index{planar angle!circular fourcomplex}
The fact that $0\leq \psi\leq\pi/2$ means that $\psi$ has
the same sign on both faces of the two-dimensional hyperplane $\upsilon O
\zeta$. The components of the point $A$
in terms of the distance $d$ and the angles $\phi, \chi, \psi$ are thus
\begin{equation}
\frac{x+t}{\sqrt{2}}=d\cos\phi \sin\psi , 
\label{12a}
\end{equation}
\begin{equation}
\frac{x-t}{\sqrt{2}}=d\cos\chi \cos\psi , 
\label{12b}
\end{equation}
\begin{equation}
\frac{y+z}{\sqrt{2}}=d\sin\phi \sin\psi , 
\label{12c}
\end{equation}
\begin{equation}
\frac{y-z}{\sqrt{2}}=d\sin\chi \cos\psi .
\label{12d}
\end{equation}
It can be checked that $\rho_+=\sqrt{2}d\sin\psi, \rho_-=\sqrt{2}d\cos\psi$.
The coordinates $x,y,z,t$ in terms of the variables $d, \phi, \chi,
\psi$ are
\begin{equation}
x=\frac{d}{\sqrt{2}}(\cos\psi\cos\chi+\sin\psi\cos\phi),
\label{12e}
\end{equation}
\begin{equation}
y=\frac{d}{\sqrt{2}}(\cos\psi\sin\chi+\sin\psi\sin\phi),
\label{12g}
\end{equation}
\begin{equation}
z=\frac{d}{\sqrt{2}}(-\cos\psi\sin\chi+\sin\psi\sin\phi),
\label{12h}
\end{equation}
\begin{equation}
t=\frac{d}{\sqrt{2}}(-\cos\psi\cos\chi+\sin\psi\cos\phi).
\label{12f}
\end{equation}
The angles $\phi, \chi, \psi$ can be expressed in terms of the coordinates
$x,y,z,t$ as 
\begin{equation}
\sin\phi = (y+z)/\rho_+ ,\: \cos\phi = (x+t)/\rho_+ ,
\label{13a}
\end{equation}
\begin{equation}
\sin\chi = (y-z)/\rho_- ,\: \cos\chi = (x-t)/\rho_- ,
\label{13b}
\end{equation}
\begin{equation}
\tan\psi=\rho_+/\rho_- .
\label{13c}
\end{equation}
The nodal hyperplanes are $\xi O\upsilon$, for which $\tau=0, \zeta=0$, and
$\tau O\zeta$, for which $\xi=0, \upsilon=0$.
For points in 
the nodal hyperplane $\xi O\upsilon$ the planar angle is $\psi=\pi/2$, 
for points in the nodal hyperplane $\tau O\zeta$ the planar angle is $\psi=0$.

\begin{figure}
\begin{center}
\epsfig{file=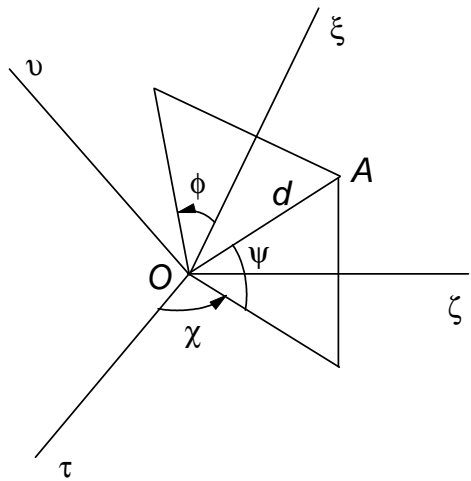,width=12cm}
\caption{Azimuthal angles $\phi, \chi$ and planar angle $\psi$ of the
fourcomplex number $x+\alpha y +\beta z+\gamma t$, represented by the point
$A$, situated at a distance $d$ from the origin $O$.}
\label{fig10}
\end{center}
\end{figure}

It can be shown that if $u_1=x_1+\alpha y_1+\beta z_1+\gamma t_1, 
u_2=x_2+\alpha y_2+\beta z_2+\gamma t_2$ are circular fourcomplex
numbers of amplitudes and angles $\rho_1, \phi_1, \chi_1, \psi_1$ and
respectively $\rho_2, \phi_2, \chi_2, \psi_2$, then the amplitude $\rho$ and
the angles $\phi, \chi, \psi$ of the product circular fourcomplex number $u_1u_2$
are 
\begin{equation}
\rho=\rho_1\rho_2, 
\label{14a}
\end{equation}\index{transformation of variables!circular fourcomplex}
\begin{equation}
 \phi=\phi_1+\phi_2, \: \chi=\chi_1+\chi_2, \: \tan\psi=\tan\psi_1\tan\psi_2 . 
\label{14b}
\end{equation}
The relations (\ref{14a})-(\ref{14b}) are consequences of the definitions
(\ref{6e})-(\ref{7b}), (\ref{13a})-(\ref{13c}) and of the identities
\begin{eqnarray}
\lefteqn{[(x_1x_2-y_1y_2-z_1z_2+t_1t_2)+(x_1t_2+t_1x_2-y_1z_2-z_1y_2)]^2
\nonumber}\\
&& +[(x_1y_2+y_1x_2+z_1t_2+t_1z_2)+(x_1z_2+z_1x_2+y_1t_2+t_1y_2)]^2\nonumber\\
&&=[(x_1+t_1)^2+(y_1+z_1)^2][(x_2+t_2)^2+(y_2+z_2)^2] ,
\label{15}
\end{eqnarray}
\begin{eqnarray}
\lefteqn{[(x_1x_2-y_1y_2-z_1z_2+t_1t_2)-(x_1t_2+t_1x_2-y_1z_2-z_1y_2)]^2
\nonumber}\\
&& +[(x_1y_2+y_1x_2+z_1t_2+t_1z_2)-(x_1z_2+z_1x_2+y_1t_2+t_1y_2)]^2\nonumber\\
&&=[(x_1-t_1)^2+(y_1-z_1)^2][(x_2-t_2)^2+(y_2-z_2)^2] ,
\label{16}
\end{eqnarray}
\begin{eqnarray}
\lefteqn{(x_1x_2-y_1y_2-z_1z_2+t_1t_2)+(x_1t_2+t_1x_2-y_1z_2-z_1y_2)
\nonumber}\\
&& =(x_1+t_1)(x_2+t_2)-(y_1+z_1)(y_2+z_2) ,
\label{17}
\end{eqnarray}
\begin{eqnarray}
\lefteqn{(x_1x_2-y_1y_2-z_1z_2+t_1t_2)-(x_1t_2+t_1x_2-y_1z_2-z_1y_2)
\nonumber}\\
&& =(x_1-t_1)(x_2-t_2)-(y_1-z_1)(y_2-z_2) ,
\label{18}
\end{eqnarray}
\begin{eqnarray}
\lefteqn{(x_1y_2+y_1x_2+z_1t_2+t_1z_2)+(x_1z_2+z_1x_2+y_1t_2+t_1y_2)
\nonumber}\\
&&=(y_1+z_1)(x_2+t_2)+(y_2+z_2)(x_2+t_2) ,
\label{19}
\end{eqnarray}
\begin{eqnarray}
\lefteqn{(x_1y_2+y_1x_2+z_1t_2+t_1z_2)-(x_1z_2+z_1x_2+y_1t_2+t_1y_2)
\nonumber}\\
&&=(y_1-z_1)(x_2-t_2)+(y_2-z_2)(x_2-t_2) .
\label{20}
\end{eqnarray}
The identities (\ref{15}) and (\ref{16}) can also be written as
\begin{equation}
\rho_+^2=\rho_{1+}\rho_{2+} ,
\label{21a}
\end{equation}
\begin{equation}
\rho_-^2=\rho_{1-}\rho_{2-} ,
\label{21b}
\end{equation}
where
\begin{equation}
\rho_{j+}^2=(x_j+t_j)^2+(y_j+z_j)^2, \: \rho_{j-}^2=(x_j-t_j)^2+(y_j-z_j)^2 ,\:
j=1,2. 
\label{22}
\end{equation}

The fact that the amplitude of the product is equal to the product of the 
amplitudes, as written in Eq. (\ref{14a}), can 
be demonstrated also by using a representation of the multiplication of the 
circular fourcomplex numbers by matrices, in which the circular fourcomplex number $u=x+\alpha
y+\beta z+\gamma t$ is represented by the matrix
\begin{equation}
A=\left(\begin{array}{cccc}
x&y&z&t\\
-y&x&t&-z\\
-z&t&x&-y\\
t&z&y&x 
\end{array}\right) .
\label{23}
\end{equation}\index{matrix representation!circular fourcomplex}
The product $u=x+\alpha y+\beta z+\gamma t$ of the circular fourcomplex numbers
$u_1=x_1+\alpha y_1+\beta z_1+\gamma t_1, u_2=x_2+\alpha y_2+\beta z_2+\gamma
t_2$, can be represented by the matrix multiplication 
\begin{equation}
A=A_1A_2.
\label{24}
\end{equation}
It can be checked that the determinant ${\rm det}(A)$ of the matrix $A$ is
\begin{equation}
{\rm det}A = \rho^4 .
\label{25}
\end{equation}
The identity (\ref{14a}) is then a consequence of the fact the determinant 
of the product of matrices is equal to the product of the determinants 
of the factor matrices.

\subsection{The exponential and trigonometric forms of circular 
fourcomplex numbers}

The exponential function of a hypercomplex variable $u$ and the addition
theorem for the exponential function have been written in Eqs. 
~(1.35)-(1.36).
If $u=x+\alpha y+\beta z+\gamma t$, then  $\exp u$ can be calculated as
$\exp u=\exp x \cdot \exp (\alpha y) \cdot \exp (\beta z) \cdot \exp (\gamma
t)$. According to Eq. (\ref{1}), 
\begin{equation}
\alpha^{2m}=(-1)^m, \alpha^{2m+1}=(-1)^m\alpha, 
\beta^{2m}=(-1)^m, \beta^{2m+1}=(-1)^m\beta, 
\gamma^m=1,  
\label{28}
\end{equation}\index{complex units, powers of!circular fourcomplex}
where $m$ is a natural number,
so that $\exp (\alpha y), \: \exp(\beta z)$ and $\exp(\gamma t)$ can be written as
\begin{equation}
\exp (\alpha y) = \cos y +\alpha \sin y , \:
\exp (\beta z) = \cos z +\beta \sin z , 
\label{29}
\end{equation}\index{exponential, expressions!circular fourcomplex}
and
\begin{equation}
\exp (\gamma t) = \cosh t +\gamma \sinh t . 
\label{30}
\end{equation}
From Eqs. (\ref{29})-(\ref{30}) it can be inferred that
\begin{equation}
\begin{array}{l}
(\cos y +\alpha \sin y)^m=\cos my +\alpha \sin my , \\
(\cos z +\beta \sin z)^m=\cos mz +\beta \sin mz , \\
(\cosh t +\gamma \sinh t)^m=\cosh mt +\gamma \sinh mt . 
\end{array}
\label{30b}
\end{equation}

Any circular fourcomplex number $u=x+\alpha y+\beta z+\gamma t$ can be writen in the
form 
\begin{equation}
x+\alpha y+\beta z+\gamma t=e^{x_1+\alpha y_1+\beta z_1+\gamma t_1} .
\label{31}
\end{equation}
The expressions of $x_1, y_1, z_1, t_1$ as functions of 
$x, y, z, t$ can be obtained by
developing $e^{\alpha y_1}, e^{\beta z_1}$ and $e^{\gamma t_1}$ with the aid of
Eqs. (\ref{29}) and (\ref{30}), by multiplying these expressions and separating
the hypercomplex components, 
\begin{equation}
x=e^{x_1}(\cos y_1\cos z_1\cosh t_1-\sin y_1\sin z_1\sinh t_1) ,
\label{32}
\end{equation}
\begin{equation}
y=e^{x_1}(\sin y_1\cos z_1\cosh t_1+\cos y_1\sin z_1\sinh t_1) ,
\label{33}
\end{equation}
\begin{equation}
z=e^{x_1}(\cos y_1\sin z_1\cosh t_1+\sin y_1\cos z_1\sinh t_1) ,
\label{34}
\end{equation}
\begin{equation}
t=e^{x_1}(-\sin y_1\sin z_1\cosh t_1+\cos y_1\cos z_1\sinh t_1) ,
\label{35}
\end{equation}
From Eqs. (\ref{32})-(\ref{35}) it can be shown by direct calculation that
\begin{equation}
x^2+y^2+z^2+t^2=e^{2x_1}\cosh 2t_1 ,
\label{36}
\end{equation}
\begin{equation}
2(xt+yz)=e^{2x_1}\sinh 2t_1 ,
\label{37}
\end{equation}
so that
\begin{equation}
e^{4x_1}=(x^2+y^2+z^2+t^2)^2-4(xt+yz)^2 .
\label{38}
\end{equation}
By comparing the expression in the right-hand side of Eq. (\ref{38}) with the
expression of $\rho$, Eq. ({\ref{7a}), it can be seen that
\begin{equation}
e^{x_1}=\rho . 
\label{39}
\end{equation}
The variable $t_1$ is then given by
\begin{equation}
\cosh 2t_1=\frac{d^2}{\rho^2}, \: \sinh 2t_1=\frac{2(xt+yz)}{\rho^2} .
\label{40}
\end{equation}
From the fact that  $\rho^4=d^4-4(xt+yz)^2$ it follows that $d^2/\rho^2\geq 1$,
so that Eq. (\ref{40}) has always a real solution, and $t_1=0$ for $xt+yz=0$.
It can be shown similarly that
\begin{equation}
\cos 2y_1 = \frac{x^2-y^2+z^2-t^2}{\rho^2} , \: 
\sin 2y_1 = \frac{2(xy-zt)}{\rho^2} ,
\label{41}
\end{equation}
\begin{equation}
\cos 2z_1 = \frac{x^2+y^2-z^2-t^2}{\rho^2} , \: 
\sin 2z_1 = \frac{2(xz-yt)}{\rho^2} .
\label{42}
\end{equation}
It can be shown that $(x^2-y^2+z^2-t^2)^2\leq\rho^4$, the equality taking place
for $xy=zt$, and $(x^2+y^2-z^2-t^2)^2\leq\rho^4$, the equality taking place
for $xz=yt$, so that Eqs. (\ref{41}) and Eqs. (\ref{42}) have always real
solutions. 

The variables
\begin{equation}
y_1 =\frac{1}{2}\arcsin \frac{2(xy-zt)}{\rho^2} ,\;
z_1 =\frac{1}{2}\arcsin \frac{2(xz-yt)}{\rho^2} ,\;
t_1 =\frac{1}{2}{\rm argsinh}\frac{2(xt+yz)}{\rho^2} 
\label{40b}
\end{equation}
are additive upon the multiplication of circular fourcomplex numbers, as can be seen
from the identities
\begin{eqnarray}
\lefteqn{(x x^{\prime}-y y^{\prime}-z z^{\prime}+t t^{\prime})(x y^{\prime}+y x^{\prime}+z t^{\prime}+t z^{\prime})\nonumber}\\
&&-(x z^{\prime}+z x^{\prime}+y t^{\prime}+t y^{\prime})(x t^{\prime}-y z^{\prime}-z y^{\prime}+t x^{\prime})\nonumber\\
&&=(x y -z t )(x^{\prime 2}-y^{\prime 2}+z^{\prime 2}-t^{\prime 2})+(x^{ 2}-y^{ 2}+z^{ 2}-t^{ 2})(x^{\prime}y^{\prime}-z^{\prime}t^{\prime}),
\label{51ca}
\end{eqnarray}
\begin{eqnarray}
\lefteqn{(x x^{\prime}-y y^{\prime}-z z^{\prime}+t t^{\prime})(x z^{\prime}+z x^{\prime}+y t^{\prime}+t y^{\prime})
\nonumber}\\
&&-(x y^{\prime}+y x^{\prime}+z t^{\prime}+t z^{\prime})(x t^{\prime}-y z^{\prime}-z y^{\prime}+t x^{\prime})\nonumber\\
&&=(x z -y t )(x^{\prime 2}+y^{\prime 2}-z^{\prime 2}-t^{\prime 2})+(x^{ 2}+y^{ 2}-z^{ 2}-t^{ 2})(x^{\prime}z^{\prime}-y^{\prime}t^{\prime}),
\label{51cb}
\end{eqnarray}
\begin{eqnarray}
\lefteqn{(x x^{\prime}-y y^{\prime}-z z^{\prime}+t t^{\prime})(x t^{\prime}-y z^{\prime}-z y^{\prime}+t x^{\prime})
\nonumber}\\
&&+(x y^{\prime}+y x^{\prime}+z t^{\prime}+t z^{\prime})(x z^{\prime}+z x^{\prime}+y t^{\prime}+t y^{\prime})\nonumber\\
&&=(x t +y z )(x^{\prime 2}+y^{\prime 2}+z^{\prime 2}+t^{\prime 2})+(x^{ 2}+y^{ 2}+z^{ 2}+t^{ 2})(x^{\prime}t^{\prime}+y^{\prime}z^{\prime}).
\label{51cc}
\end{eqnarray}

The expressions appearing in Eqs. (\ref{40})-(\ref{42}) can be calculated in
terms of the angles $\phi, \chi, \psi$ with the aid of Eqs.
(\ref{12e})-(\ref{12f}) as
\begin{equation}
\frac{d^2}{\rho^2}=\frac{1}{\sin
2\psi},\;\frac{2(xt+yz)}{\rho^2}=-\frac{1}{\tan 2\psi},
\label{51cd}
\end{equation}
\begin{equation}
\frac{x^2-y^2+z^2-t^2}{\rho^2}=\cos(\phi+\chi),\;\frac{2(xy-zt)}{\rho^2}=\sin(\phi+\chi),
\label{51ce}
\end{equation}
\begin{equation}
\frac{x^2+y^2-z^2-t^2}{\rho^2}=\cos(\phi-\chi),\;\frac{2(xz-yt)}{\rho^2}=\sin(\phi-\chi).
\label{51cf}
\end{equation}
Then from Eqs. (\ref{40})-(\ref{42}) and (\ref{51cd})-(\ref{51cf}) it results that
\begin{equation}
y_1 =\frac{\phi+\chi}{2},\; z_1=\frac{\phi-\chi}{2},\;t_1=\frac{1}{2}\ln\tan\psi,
\label{51cg}
\end{equation}
so that the circular fourcomplex number $u$, Eq. (\ref{31}), can be written as
\begin{equation}
u=\rho\exp\left[\alpha\frac{\phi+\chi}{2}+\beta\frac{\phi-\chi}{2}
+\gamma \frac{1}{2}\ln\tan\psi\right] .
\label{51}
\end{equation}\index{exponential form!circular fourcomplex}
In Eq. (\ref{51}) the circular fourcomplex number $u=x+\alpha y+\beta z+\gamma t$ is
written as the product of the amplitude $\rho$ and of an exponential function,
and therefore this form of $u$ will be called the exponential form of the
circular fourcomplex number.
It can be checked that
\begin{equation}
\exp\left(\frac{\alpha+\beta}{2}\phi\right)=\frac{1-\gamma}{2}
+\frac{1+\gamma}{2}\cos\phi+\frac{\alpha+\beta}{2}\sin\phi ,
\label{51b}
\end{equation}
\begin{equation}
\exp\left(\frac{\alpha-\beta}{2}\chi\right)=\frac{1+\gamma}{2}
+\frac{1-\gamma}{2}\cos\chi+\frac{\alpha-\beta}{2}\sin\phi ,
\label{51c}
\end{equation}
which shows that $e^{(\alpha+\beta)\phi/2}$ and $e^{(\alpha-\beta)\chi/2}$ are periodic
functions of $\phi$ and respectively $\chi$, with period $2\pi$.

The relations between the variables $y_1, z_1, t_1$ and the angles $\phi, \chi,
\psi$ can be obtained alternatively by
substituting in Eqs. (\ref{32})-(\ref{35}) the expression $e^{x_1}=d/(\cosh
2t_1)^{1/2} $, Eq. (\ref{36}), and summing and subtracting of the relations,
\begin{equation}
\frac{x+t}{\sqrt{2}}=d\cos(y_1+z_1)\sin(\eta+\pi/4) ,
\label{43}
\end{equation}
\begin{equation}
\frac{x-t}{\sqrt{2}}=d\cos(y_1-z_1)\cos(\eta+\pi/4) ,
\label{44}
\end{equation}
\begin{equation}
\frac{y+z}{\sqrt{2}}=d\sin(y_1+z_1)\sin(\eta+\pi/4) ,
\label{45}
\end{equation}
\begin{equation}
\frac{y-z}{\sqrt{2}}=d\sin(y_1-z_1)\cos(\eta+\pi/4) ,
\label{46}
\end{equation}
where the variable $\eta$ is defined by the relations
\begin{equation}
\frac{\cosh t_1}{(\cosh 2t_1)^{1/2}}=\cos\eta, \: 
\frac{\sinh t_1}{(\cosh 2t_1)^{1/2}}=\sin\eta ,
\label{47}
\end{equation}
and when $-\infty<t_1<\infty$, the range of the variable $\eta$ is
$-\pi/4\leq\eta\leq \pi/4$.
The comparison of Eqs. (\ref{12a})-(\ref{12d}) and (\ref{43})-(\ref{46}) shows
that 
\begin{equation}
\phi=y_1+z_1, \: \chi=y_1-z_1, \: \psi=\eta+\pi/4 .
\label{48}
\end{equation}

It can be shown with the aid of Eq. (\ref{30}) that
\begin{equation}
\exp \left(\frac{1}{2}\gamma\ln\tan\psi \right) 
=\frac{1}{(\sin 2\psi)^{1/2}}\left[
\cos(\psi-\pi/4)+\gamma\sin(\psi-\pi/4)\right] . 
\label{30bb}
\end{equation}
The circular fourcomplex number $u$, Eq. (\ref{51}), can then be written equivalently as
\begin{equation}
u=d\{\cos(\psi-\pi/4)+\gamma\sin(\psi-\pi/4)\}
\exp\left[\alpha\frac{\phi+\chi}{2}+\beta\frac{\phi-\chi}{2}\right].
\label{52}
\end{equation}\index{trigonometric form!circular fourcomplex}
In Eq. (\ref{52}), the circular fourcomplex number 
$u=x+\alpha y+\beta z+\gamma t$ is
written as the product of the modulus $d$ and of factors depending on the
geometric angles $\phi, \chi$ and $\psi$, and this form will be called the
trigonometric form of the circular fourcomplex number.

If $u_1, u_2$ are circular fourcomplex numbers of moduli and angles $d_1, \phi_1,
\chi_1, \psi_1$ and respectively $d_2, \phi_2, \chi_2, \psi_2$, the product of
the planar factors can be calculated to be
\begin{eqnarray}
\lefteqn{[\cos(\psi_1-\pi/4)+\gamma\sin(\psi_1-\pi/4)] 
[\cos(\psi_2-\pi/4)+\gamma\sin(\psi_2-\pi/4)]}\nonumber\\
&&=[\cos(\psi_1-\psi_2)-\gamma\cos(\psi_1+\psi_2)] .
\label{53}
\end{eqnarray}
The right-hand side of Eq. (\ref{53}) can be written as
\begin{eqnarray}
\lefteqn{\cos(\psi_1-\psi_2)-\gamma\cos(\psi_1+\psi_2)}\nonumber\\
&&=[2(\cos^2\psi_1\cos^2\psi_2+\sin^2\psi_1\sin^2\psi_2)]^{1/2}
[\cos(\psi-\pi/4)+\gamma\sin(\psi-\pi/4)] ,
\label{54}
\end{eqnarray}
where the angle $\psi$, determined by the condition that
\begin{equation}
\tan(\psi-\pi/4)=-\cos(\psi_1+\psi_2)/\cos(\psi_1-\psi_2)
\label{55}
\end{equation}
is given by $\tan\psi=\tan\psi_1\tan\psi_2$ ,
which is consistent with Eq. (\ref{14b}).
It can be checked that the modulus $d$ of the product $u_1u_2$ is 
\begin{equation}
d=\sqrt{2}d_1d_2\left(\cos^2\psi_1\cos^2\psi_2+\sin^2\psi_1\sin^2\psi_2\right)^{1/2} .
\label{56}
\end{equation}

\subsection{Elementary functions of a circular fourcomplex variable}

The logarithm $u_1$ of the circular fourcomplex number $u$, $u_1=\ln u$, can be defined
as the solution of the equation
\begin{equation}
u=e^{u_1} ,
\label{57}
\end{equation}
written explicitly previously in Eq. (\ref{31}), for $u_1$ as a function of
$u$. From Eq. (\ref{51}) it results that 
\begin{equation}
\ln u=\ln \rho+\frac{1}{2}\gamma \ln\tan\psi
+\alpha\frac{\phi+\chi}{2}+\beta\frac{\phi-\chi}{2} .
\label{58}
\end{equation}\index{logarithm!circular fourcomplex}
It can be inferred from Eqs. (\ref{14a}) and (\ref{14b}) that
\begin{equation}
\ln(uu^\prime)=\ln u+\ln u^\prime ,
\label{59}
\end{equation}
up to multiples of $\pi(\alpha+\beta)$ and $\pi(\alpha-\beta)$.

The power function $u^n$ can be defined for real values of $n$ as
\begin{equation}
u^m=e^{m\ln u} .
\label{60}
\end{equation}\index{power function!circular fourcomplex}
The power function is multivalued unless $n$ is an integer. 
For integer $n$, it can be inferred from Eq. (\ref{59}) that
\begin{equation}
(uu^\prime)^m=u^n\:u^{\prime m} .
\label{61}
\end{equation}
If, for example, $m=2$, it can be checked with the aid of Eq. (\ref{52})
that Eq. (\ref{60}) gives indeed $(x+\alpha y+\beta z+\gamma t)^2=
x^2-y^2-z^2+t^2+2\alpha(xy+zt)+2\beta(xz+yt)+2\gamma(2xt-yz)$.

The trigonometric functions of the hypercomplex variable
$u$ and the addition theorems for these functions have been written in Eqs.
~(1.57)-(1.60). 
The cosine and sine functions of the hypercomplex variables $\alpha y, 
\beta z$ and $ \gamma t$ can be expressed as
\begin{equation}
\cos\alpha y=\cosh y, \: \sin\alpha y=\alpha\sinh y, 
\label{66}
\end{equation}\index{trigonometric functions, expressions!circular fourcomplex}
\begin{equation}
\cos\beta y=\cosh y, \: \sin\beta y=\beta\sinh y, 
\label{67}
\end{equation}
\begin{equation}
\cos\gamma y=\cos y, \: \sin\gamma y=\gamma\sin y .
\label{68}
\end{equation}
The cosine and sine functions of a circular fourcomplex number $x+\alpha y+\beta
z+\gamma t$ can then be
expressed in terms of elementary functions with the aid of the addition
theorems Eqs. (1.59), (1.60) and of the expressions in  Eqs. 
(\ref{66})-(\ref{68}).

The hyperbolic functions of the hypercomplex variable
$u$ and the addition theorems for these functions have been written in Eqs.
~(1.62)-(1.65). 
The hyperbolic cosine and sine functions of the hypercomplex variables 
$\alpha y,  \beta z$ and $ \gamma t$ can be expressed as
\begin{equation}
\cosh\alpha y=\cos y, \: \sinh\alpha y=\alpha\sin y, 
\label{73}
\end{equation}\index{hyperbolic functions, expressions!circular fourcomplex}
\begin{equation}
\cosh\beta y=\cos y, \: \sinh\beta y=\beta\sin y, 
\label{74}
\end{equation}
\begin{equation}
\cosh\gamma y=\cosh y, \: \sinh\gamma y=\gamma\sinh y .
\label{75}
\end{equation}
The hyperbolic cosine and sine functions of a circular fourcomplex number 
$x+\alpha y+\beta z+\gamma t$ can then be
expressed in terms of elementary functions with the aid of the addition
theorems Eqs. (1.64), (1.65) and of the expressions in  Eqs. 
(\ref{73})-(\ref{75}).

\subsection{Power series of circular fourcomplex variables}

A circular fourcomplex series is an infinite sum of the form
\begin{equation}
a_0+a_1+a_2+\cdots+a_n+\cdots , 
\label{76}
\end{equation}\index{series!circular fourcomplex}
where the coefficients $a_n$ are circular fourcomplex numbers. The convergence of 
the series (\ref{76}) can be defined in terms of the convergence of its 4 real
components. The convergence of a circular fourcomplex series can however be studied
using circular fourcomplex variables. The main criterion for absolute convergence 
remains the comparison theorem, but this requires a number of inequalities
which will be discussed further.

The modulus of a circular fourcomplex number $u=x+\alpha y+\beta z+\gamma t$ can be
defined as 
\begin{equation}
|u|=(x^2+y^2+z^2+t^2)^{1/2} ,
\label{77}
\end{equation}\index{modulus, definition!circular fourcomplex}
so that, according to Eq. (\ref{10}), $d=|u|$. Since $|x|\leq |u|, |y|\leq |u|,
|z|\leq |u|, |t|\leq |u|$, a property of 
absolute convergence established via a comparison theorem based on the modulus
of the series (\ref{76}) will ensure the absolute convergence of each real
component of that series.

The modulus of the sum $u_1+u_2$ of the circular fourcomplex numbers $u_1, u_2$ fulfils
the inequality
\begin{equation}
||u_1|-|u_2||\leq |u_1+u_2|\leq |u_1|+|u_2| .
\label{78}
\end{equation}\index{modulus, inequalities!circular fourcomplex}
For the product the relation is 
\begin{equation}
|u_1u_2|\leq \sqrt{2}|u_1||u_2| ,
\label{79}
\end{equation}
which replaces the relation of equality extant for regular complex numbers.
The equality in Eq. (\ref{79}) takes place for $x_1=t_1, y_1=z_1,
x_2=t_2, y_2=z_2$ or $x_1=-t_1, y_1=-z_1,
x_2=-t_2, y_2=-z_2$. In Eq. (\ref{56}), this corresponds to $\psi_1=0,\psi_2=0$
or $\psi_1=\pi/2, \psi_2=\pi/2$.
The modulus of a product, which has the property that
$0\leq|u_1u_2|$, becomes equal to zero for
$x_1=t_1, y_1=z_1, x_2=-t_2, y_2=-z_2$ or $x_1=-t_1, y_1=-z_1,
x_2=t_2, y_2=z_2$, as discussed after Eq. (\ref{9}).
In Eq. (\ref{56}), the latter situation corresponds to $\psi_1=0,\psi_2=\pi/2$
or $\psi_1=0, \psi_2=\pi/2$.

It can be shown that
\begin{equation}
x^2+y^2+z^2+t^2\leq|u^2|\leq \sqrt{2}(x^2+y^2+z^2+t^2) .
\label{80}
\end{equation}
The left relation in Eq. (\ref{80}) becomes an equality, 
$x^2+y^2+z^2+t^2=|u^2|$, for $xt+yz=0$. This condition corresponds to
$\psi_1=\psi_2=\pi/4$ in Eq. (\ref{56}).
The inequality in Eq. (\ref{79}) implies that
\begin{equation}
|u^l|\leq 2^{(l-1)/2}|u|^l .
\label{81}
\end{equation}
From Eqs. (\ref{79}) and (\ref{81}) it results that
\begin{equation}
|au^l|\leq 2^{l/2} |a| |u|^l .
\label{82}
\end{equation}

A power series of the circular fourcomplex variable $u$ is a series of the form
\begin{equation}
a_0+a_1 u + a_2 u^2+\cdots +a_l u^l+\cdots .
\label{83}
\end{equation}\index{power series!circular fourcomplex}
Since
\begin{equation}
\left|\sum_{l=0}^\infty a_l u^l\right| \leq  \sum_{l=0}^\infty
2^{l/2}|a_l| |u|^l ,
\label{84}
\end{equation}
a sufficient condition for the absolute convergence of this series is that
\begin{equation}
\lim_{l\rightarrow \infty}\frac{\sqrt{2}|a_{l+1}||u|}{|a_l|}<1 .
\label{85}
\end{equation}
Thus the series is absolutely convergent for 
\begin{equation}
|u|<c,
\label{86}
\end{equation}\index{convergence of power series!circular fourcomplex}
where 
\begin{equation}
c=\lim_{l\rightarrow\infty} \frac{|a_l|}{\sqrt{2}|a_{l+1}|} .
\label{87}
\end{equation}

The convergence of the series (\ref{83}) can be also studied with the aid of
the transformation 
\begin{equation}
x+\alpha y+\beta z+\gamma t=\sqrt{2}(e_1\xi+\tilde e_1 \upsilon+e_2\tau
+\tilde e_2\zeta) , 
\label{87b}
\end{equation}\index{canonical form!circular fourcomplex}
where $\xi,\upsilon, \tau, \zeta$ have been defined in Eq. (\ref{11}),
and
\begin{equation}
e_1=\frac{1+\gamma}{2},\:\:\tilde e_1=\frac{\alpha+\beta}{2},\:\:
e_2=\frac{1-\gamma}{2},\:\:\tilde e_2=\frac{\alpha-\beta}{2}.
\label{87c}
\end{equation}\index{canonical base!circular fourcomplex}
The ensemble $e_1, \tilde e_1, e_2, \tilde e_2$ will be called the canonical
circular fourcomplex base, and Eq. (\ref{87b}) gives the canonical form of the
circular fourcomplex number.
It can be checked that
\begin{eqnarray}
\lefteqn{e_1^2=e_1, \:\:\tilde e_1^2=-e_1,\:\: e_1\tilde e_1=\tilde e_1,\:\:
e_2^2=e_2, \:\:\tilde e_2^2=-e_2,\:\: e_2\tilde e_2=\tilde e_2,\:\:\nonumber}\\
&&e_1e_2=0,\:\: \tilde e_1\tilde e_2=0, \:\:e_1\tilde e_2=0, \:\:
e_2\tilde e_1=0.  
\label{87d}
\end{eqnarray}
The moduli of the bases in Eq. (\ref{87c}) are
\begin{equation}
|e_1|=\frac{1}{\sqrt{2}},\;|\tilde e_1|=\frac{1}{\sqrt{2}},\;
|e_2|=\frac{1}{\sqrt{2}},\;|\tilde e_2|=\frac{1}{\sqrt{2}},
\label{87e}
\end{equation}
and it can be checked that
\begin{equation}
|x+\alpha y+\beta z+\gamma t|^2=\xi^2+\upsilon^2+\tau^2+\zeta^2.
\label{87f}
\end{equation}\index{modulus, canonical variables!circular fourcomplex}

If $u=u^\prime u^{\prime\prime}$, the components $\xi,\upsilon, \tau, \zeta$
are related, according to Eqs. (\ref{17})-(\ref{20}) by 
\begin{equation}
\xi=\sqrt{2}(\xi^\prime \xi^{\prime\prime}-\upsilon^\prime \upsilon^{\prime\prime}), \:\:
\upsilon=\sqrt{2}(\xi^\prime \upsilon^{\prime\prime}+\upsilon^\prime \xi^{\prime\prime}), \:\:
\tau=\sqrt{2}(\tau^\prime \tau^{\prime\prime}-\zeta^\prime \zeta^{\prime\prime}), \:\:
\zeta=\sqrt{2}(\tau^\prime \zeta^{\prime\prime}+\zeta^\prime \tau^{\prime\prime}), \:\:
\label{87g}
\end{equation}\index{transformation of variables!circular fourcomplex}
which show that, upon multiplication, the components $\xi,\upsilon$ and $\tau,
\zeta$ obey, up to a normalization constant, the same
rules as the real and imaginary components of usual, two-dimensional complex
numbers.

If the coefficients in Eq. (\ref{83}) are 
\begin{equation}
a_l= a_{l0}+\alpha a_{l1}+\beta a_{l2}+\gamma a_{l3}, 
\label{n88a}
\end{equation}
and
\begin{equation}
A_{l1}=a_{l0}+a_{l3},\;
\tilde A_{l1}=a_{l1}+a_{l2},\;
A_{l2}=a_{l0}-a_{l3},\;
\tilde A_{l2}=a_{l1}-a_{l2},
\label{n88b}
\end{equation}
the series (\ref{83}) can be written as
\begin{equation}
\sum_{l=0}^\infty 2^{l/2}\left[
(e_1 A_{l1}+\tilde e_1\tilde A_{l1})(e_1 \xi+\tilde e_1 \upsilon)^l 
+(e_2 A_{l2}+\tilde e_2\tilde A_{l2})(e_2 \tau+\tilde e_2 \zeta)^l 
\right].
\label{n89a}
\end{equation}
Thus, the series in Eqs. (\ref{83}) and (\ref{n89a}) are
absolutely convergent for   
\begin{equation}
\rho_+<c_1, \;\rho_-<c_2,
\label{n90}
\end{equation}\index{convergence, region of!circular fourcomplex}
where 
\begin{equation}
c_1=\lim_{l\rightarrow\infty} \frac
{\left[A_{l1}^2+\tilde A_{l1}^2\right]^{1/2}}
{\sqrt{2}\left[A_{l+1,1}^2+\tilde A_{l+1,1}^2\right]^{1/2}},\;\;
c_2=\lim_{l\rightarrow\infty} \frac
{\left[A_{l2}^2+\tilde A_{l2}^2\right]^{1/2}}
{\sqrt{2}\left[A_{l+1,2}^2+\tilde A_{l+1,2}^2\right]^{1/2}}.
\label{n91}
\end{equation}

It can be shown that $c=(1/\sqrt{2}){\rm
min}(c_1,c_2)$, where ${\rm min}$ designates the smallest of
the numbers $c_1,c_2$. Using the expression of $|u|$ in
Eq. (\ref{87f}),  it can be seen that the spherical region of
convergence defined in Eqs. (\ref{86}), (\ref{87}) is included in the
cylindrical region of convergence defined in Eqs. (\ref{n90}) and (\ref{n91}).

\subsection{Analytic functions of circular fourcomplex variables}

The analytic functions of the hypercomplex variable $u$ and the series 
expansion of functions have been discussed in Eqs. ~(1.85)-(1.93).
If the fourcomplex function $f(u)$
of the fourcomplex variable $u$ can be expressed in terms of 
the real functions $P(x,y,z,t)$, $Q(x,y,z,t)$, $R(x,y,z,t)$, $S(x,y,z,t)$ of real
variables $x,y,z,t$ as 
\begin{equation}
f(u)=P(x,y,z,t)+\alpha Q(x,y,z,t)+\beta R(x,y,z,t)+\gamma S(x,y,z,t),
\label{g16}
\end{equation}\index{functions, real components!circular fourcomplex}
then relations of equality 
exist between partial derivatives of the functions $P,Q,R,S$. These relations
can be obtained by writing the derivative of the function $f$ as
\begin{eqnarray}
\lim_{u\rightarrow u_0}
\lefteqn{\frac{1}{\Delta x+\alpha \Delta y +\beta\Delta z+\gamma\Delta t} 
\left[\frac{\partial P}{\partial x}\Delta x+
\frac{\partial P}{\partial y}\Delta y+
\frac{\partial P}{\partial z}\Delta z+
\frac{\partial P}{\partial t}\Delta t\right.\nonumber}\\
&&+\alpha\left(\frac{\partial Q}{\partial x}\Delta x+
\frac{\partial Q}{\partial y}\Delta y+
\frac{\partial Q}{\partial z}\Delta z 
+\frac{\partial Q}{\partial t}\Delta t \right)
+\beta\left(\frac{\partial R}{\partial x}\Delta x+
\frac{\partial R}{\partial y}\Delta y+
\frac{\partial R}{\partial z}\Delta z
+\frac{\partial R}{\partial t}\Delta t
\right)\nonumber\\
&&\left.+\gamma\left(\frac{\partial S}{\partial x}\Delta x+
\frac{\partial S}{\partial y}\Delta y+
\frac{\partial S}{\partial z}\Delta z
+\frac{\partial S}{\partial t}\Delta t\right)\right] ,
\label{g17}
\end{eqnarray}\index{derivative, independence of direction!fourcomplex}
where the difference appearing in Eq. (1.86) is
$u-u_0=\Delta x+\alpha\Delta y +\beta\Delta z+\gamma\Delta t$. 
These relations have the same form for all systems of hypercomplex numbers
discussed in this work.

For the present system of hypercomplex numbers, the
relations between the partial derivatives of the functions $P, Q, R, S$ are
obtained by setting succesively in   
Eq. (\ref{g17}) $\Delta x\rightarrow 0, \Delta y=\Delta z=\Delta t=0$;
then $ \Delta y\rightarrow 0, \Delta x=\Delta z=\Delta t=0;$  
then $  \Delta z\rightarrow 0,\Delta x=\Delta y=\Delta t=0$; and finally
$ \Delta t\rightarrow 0,\Delta x=\Delta y=\Delta z=0 $. 
The relations are \index{relations between partial derivatives!circular fourcomplex}
\begin{equation}
\frac{\partial P}{\partial x} = \frac{\partial Q}{\partial y} =
\frac{\partial R}{\partial z} = \frac{\partial S}{\partial t},
\label{95}
\end{equation}
\begin{equation}
\frac{\partial Q}{\partial x} = -\frac{\partial P}{\partial y} =
-\frac{\partial S}{\partial z} = \frac{\partial R}{\partial t},
\label{96}
\end{equation}
\begin{equation}
\frac{\partial R}{\partial x} = -\frac{\partial S}{\partial y} =
-\frac{\partial P}{\partial z} = \frac{\partial Q}{\partial t},
\label{97}
\end{equation}
\begin{equation}
\frac{\partial S}{\partial x} = \frac{\partial R}{\partial y} =
\frac{\partial Q}{\partial z} = \frac{\partial P}{\partial t}.
\label{98}
\end{equation}
The relations (\ref{95})-(\ref{98}) are analogous to the Riemann relations
for the real and imaginary components of a complex function. It can be shown
from Eqs. (\ref{95})-(\ref{98}) that the component $P$ is a solution
of the equations 
\begin{equation}
\frac{\partial^2 P}{\partial x^2}+\frac{\partial^2 P}{\partial y^2}=0,
\:\: 
\frac{\partial^2 P}{\partial x^2}+\frac{\partial^2 P}{\partial z^2}=0,
\:\:
\frac{\partial^2 P}{\partial y^2}+\frac{\partial^2 P}{\partial t^2}=0,
\:\:
\frac{\partial^2 P}{\partial z^2}+\frac{\partial^2 P}{\partial t^2}=0,
\:\:
\label{99}
\end{equation}\index{relations between second-order derivatives!circular fourcomplex}
\begin{equation}
\frac{\partial^2 P}{\partial x^2}-\frac{\partial^2 P}{\partial t^2}=0,
\:\:
\frac{\partial^2 P}{\partial y^2}-\frac{\partial^2 P}{\partial z^2}=0,
\label{100}
\end{equation}
and the components $Q, R, S$ are solutions of similar equations.

As can be seen from Eqs. (\ref{99})-(\ref{100}), the components $P, Q, R, S$ of
an analytic function of circular fourcomplex variable are harmonic 
with respect to the pairs of variables $x,y; x,z; y,t$ and $ z,t$, and are
solutions of the wave 
equation with respect to the pairs of variables $x,t$ and $y,z$.
The components $P, Q, R, S$ are also solutions of the mixed-derivative
equations 
\begin{equation}
\frac{\partial^2 P}{\partial x\partial y}=\frac{\partial^2 P}{\partial
z\partial t} ,
\:\: 
\frac{\partial^2 P}{\partial x\partial z}=\frac{\partial^2 P}{\partial
y\partial t} ,
\:\: 
\frac{\partial^2 P}{\partial x\partial t}=-\frac{\partial^2 P}{\partial
y\partial z} ,
\label{107}
\end{equation}
and the components $Q, R, S$ are solutions of similar equations.


\subsection{Integrals of functions of circular fourcomplex variables}

The singularities of circular fourcomplex functions arise from terms of the form
$1/(u-u_0)^m$, with $m>0$. Functions containing such terms are singular not
only at $u=u_0$, but also at all points of the two-dimensional hyperplanes
passing through $u_0$ and which are parallel to the nodal hyperplanes. 

The integral of a circular fourcomplex function between two points $A, B$ along a path
situated in a region free of singularities is independent of path, which means
that the integral of an analytic function along a loop situated in a region
free from singularities is zero,
\begin{equation}
\oint_\Gamma f(u) du = 0,
\label{111}
\end{equation}
where it is supposed that a surface $\Sigma$ spanning 
the closed loop $\Gamma$ is not intersected by any of
the two-dimensional hyperplanes associated with the
singularities of the function $f(u)$. Using the expression, Eq. (\ref{g16}),
for $f(u)$ and the fact that $du=dx+\alpha  dy+\beta dz+\gamma dt$, the
explicit form of the integral in Eq. (\ref{111}) is
\begin{eqnarray}
\lefteqn{\oint _\Gamma f(u) du = \oint_\Gamma
[(Pdx-Qdy-Rdz+Sdt)+\alpha(Qdx+Pdy+Sdz+Rdt)\nonumber}\\
&&+\beta(Rdx+Sdy+Pdz+Qdt)+\gamma(Sdx-Rdy-Qdz+Pdt)] .
\label{112}
\end{eqnarray}\index{integrals, path!circular fourcomplex}
If the functions $P, Q, R, S$ are regular on a surface $\Sigma$
spanning the loop $\Gamma$,
the integral along the loop $\Gamma$ can be transformed with the aid of the
theorem of Stokes in an integral over the surface $\Sigma$ of terms of the form
$\partial P/\partial y +  \partial Q/\partial x, \:\:
\partial P/\partial z + \partial R/\partial x,\:\:
\partial P/\partial t - \partial S/\partial x, \:\:
\partial Q/\partial z -  \partial R/\partial y, \:\:
\partial Q/\partial t + \partial S/\partial y,\:\:
\partial R/\partial t + \partial S/\partial z$ and of similar terms arising
from the $\alpha, \beta$ and $\gamma$ components, 
which are equal to zero by Eqs. (\ref{95})-(\ref{98}), and this proves Eq.
(\ref{111}). 

The integral of the function $(u-u_0)^m$ on a closed loop $\Gamma$ is equal to
zero for $m$ a positive or negative integer not equal to -1,
\begin{equation}
\oint_\Gamma (u-u_0)^m du = 0, \:\: m \:\:{\rm integer},\: m\not=-1 .
\label{112b}
\end{equation}
This is due to the fact that $\int (u-u_0)^m du=(u-u_0)^{m+1}/(m+1), $ and to the
fact that the function $(u-u_0)^{m+1}$ is singlevalued for $m$ an integer.

The integral $\oint du/(u-u_0)$ can be calculated using the exponential form 
(\ref{51}),
\begin{equation}
u-u_0=\rho\exp\left(\alpha\frac{\phi+\chi}{2}+\beta\frac{\phi-\chi}{2}
+\gamma \ln\tan\psi\right) ,
\label{113}
\end{equation}
so that 
\begin{equation}
\frac{du}{u-u_0}=\frac{d\rho}{\rho}+\frac{\alpha+\beta}{2}d\phi
+\frac{\alpha-\beta}{2}d\chi+\gamma d\ln\tan\psi .
\label{114}
\end{equation}
Since $\rho$ and $\ln\tan\psi$ are singlevalued variables, it follows that
$\oint_\Gamma d\rho/\rho =0, \oint_\Gamma d\ln\tan\psi=0$. On the other hand, $\phi$
and $\chi$ are cyclic variables, so that they may give a contribution to the
integral around the closed loop $\Gamma$.
Thus, if $C_+$ is a circle of radius $r$
parallel to the $\xi O\upsilon$ plane, and the
projection of the center of this circle on the $\xi O\upsilon$ plane
coincides with the projection of the point $u_0$ on this plane, the points
of the circle $C_+$ are described according to Eqs.
(\ref{11})-(\ref{12d}) by the equations
\begin{eqnarray}
\lefteqn{\xi=\xi_0+r \sin\psi\cos\phi , \:
\upsilon=\upsilon_0+r \sin\psi\sin\phi , \:
\tau=\tau_0+r\cos\psi \cos\chi , \nonumber}\\
&&\zeta=\zeta_0+r \cos\psi\sin\chi , 
\label{115}
\end{eqnarray}
for constant values of $\chi$ and $\psi, \:\psi\not=0, \pi/2$, where
$u_0=x_0+\alpha y_0+\beta 
z_0+\gamma t_0$,  and $\xi_0, \upsilon_0, \tau_0, \zeta_0$ are calculated from
$x_0, y_0, z_0, t_0$ according to Eqs. (\ref{11}).
Then
\begin{equation}
\oint_{C_+}\frac{du}{u-u_0}=\pi(\alpha+\beta) .
\label{116}
\end{equation}
If $C_-$ is a circle of radius $r$
parallel to the $\tau O\zeta$ plane,
and the projection of the center of this circle on the $\tau O\zeta$ plane
coincides with the projection of the point $u_0$ on this plane, the points
of the circle $C_-$ are described by the same Eqs. (\ref{115}) 
but for constant values of $\phi$ and $\psi, \:\psi\not=0, \pi/2$
Then
\begin{equation}
\oint_{C_-}\frac{du}{u-u_0}=\pi(\alpha-\beta) .
\label{117}
\end{equation}\index{poles and residues!circular fourcomplex}
The expression of $\oint_\Gamma du/(u-u_0)$ can be written as a single
equation with the aid of a functional int($M,C$) defined for a point $M$ and a
closed curve $C$ in a two-dimensional plane, such that 
\begin{equation}
{\rm int}(M,C)=\left\{
\begin{array}{l}
1 \;\:{\rm if} \;\:M \;\:{\rm is \;\:an \;\:interior \;\:point \;\:of} \;\:C ,\\ 
0 \;\:{\rm if} \;\:M \;\:{\rm is \;\:exterior \;\:to}\:\; C .\\
\end{array}\right.
\label{118}
\end{equation}
With this notation the result of the integration along a closed path
$\Gamma$ can be written as 
\begin{equation}
\oint_\Gamma\frac{du}{u-u_0}=
\pi(\alpha+\beta) \:\;{\rm int}(u_{0\xi\upsilon},\Gamma_{\xi\upsilon})
+\pi (\alpha-\beta)\:\;{\rm int}(u_{0\tau\zeta},\Gamma_{\tau\zeta}),
\label{119}
\end{equation}
where $u_{0\xi\upsilon}, u_{0\tau\zeta}$ and $\Gamma_{\xi\upsilon},
\Gamma_{\tau\zeta}$ are respectively the projections of the point $u_0$ and of
the loop $\Gamma$ on the planes $\xi \upsilon$ and $\tau \zeta$.

If $f(u)$ is an analytic circular fourcomplex function which can be expanded in a
series as written in Eq. (1.89), and the expansion holds on the curve
$\Gamma$ and on a surface spanning $\Gamma$, then from Eqs. (\ref{112b}) and
(\ref{119}) it follows that
\begin{equation}
\oint_\Gamma \frac{f(u)du}{u-u_0}=
\pi[(\alpha+\beta) \:\;{\rm int}(u_{0\xi\upsilon},\Gamma_{\xi\upsilon})
+ (\alpha-\beta)\:\;{\rm int}(u_{0\tau\zeta},\Gamma_{\tau\zeta})]\;f(u_0) ,
\label{120}
\end{equation}
where $\Gamma_{\xi\upsilon}, \Gamma_{\tau\zeta}$ are the projections of 
the curve $\Gamma$ on the planes $\xi \upsilon$ and respectively $\tau \zeta$,
as shown in Fig. \ref{fig11}.
Substituting in the right-hand side of 
Eq. (\ref{120}) the expression of $f(u)$ in terms of the real 
components $P, Q, R, S$, Eq. (\ref{g16}), yields
\begin{eqnarray}
\lefteqn{\oint_\Gamma \frac{f(u)du}{u-u_0}=
\pi [-(1+\gamma)(Q+R)+(\alpha+\beta)(P+S)] 
\:\;{\rm int}(u_{0\xi\upsilon},\Gamma_{\xi\upsilon})\nonumber}\\
&&+\pi [-(1-\gamma)(Q-R)+(\alpha-\beta)(P-S)] 
\:\;{\rm int}(u_{0\tau\zeta},\Gamma_{\tau\zeta}) ,
\label{121}
\end{eqnarray}
where $P, Q, R, S$ are the values of the components of $f$ at $u=u_0$.

\begin{figure}
\begin{center}
\epsfig{file=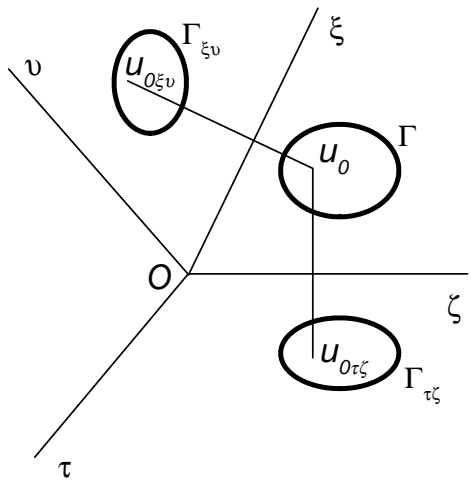,width=12cm}
\caption{Integration path $\Gamma$ and the pole $u_0$, and their projections
$\Gamma_{\xi\upsilon}, \Gamma_{\tau\zeta}$ and $u_{0\xi\upsilon},
u_{0\tau\zeta}$ on the planes $\xi \upsilon$ and respectively $\tau \zeta$.}
\label{fig11}
\end{center}
\end{figure}

If $f(u)$ can be expanded as written in Eq. (1.89) on 
$\Gamma$ and on a surface spanning $\Gamma$, then from Eqs. (\ref{112b}) and
(\ref{119}) it also results that
\begin{equation}
\oint_\Gamma \frac{f(u)du}{(u-u_0)^{m+1}}=
\frac{\pi}{m!}[(\alpha+\beta) \:\;{\rm int}(u_{0\xi\upsilon},\Gamma_{\xi\upsilon})
+ (\alpha-\beta)\:\;{\rm int}(u_{0\tau\zeta},\Gamma_{\tau\zeta})]\;
f^{(m)}(u_0) ,
\label{122}
\end{equation}
where it has been used the fact that the derivative $f^{(m)}(u_0)$ of order $m$
of $f(u)$ at $u=u_0$ is related to the expansion coefficient in Eq. (1.89)
according to Eq. (1.93).

If a function $f(u)$ is expanded in positive and negative powers of $u-u_j$,
where $u_j$ are circular fourcomplex constants, $j$ being an index, the integral of $f$
on a closed loop $\Gamma$ is determined by the terms in the expansion of $f$
which are of the form $a_j/(u-u_j)$,
\begin{equation}
f(u)=\cdots+\sum_j\frac{a_j}{u-u_j}+\cdots
\label{123}
\end{equation}
Then the integral of $f$ on a closed loop $\Gamma$ is
\begin{equation}
\oint_\Gamma f(u) du = 
\pi(\alpha+\beta) \sum_j{\rm int}(u_{j\xi\upsilon},\Gamma_{\xi\upsilon})a_j
+ \pi(\alpha-\beta)\sum_j{\rm int}(u_{j\tau\zeta},\Gamma_{\tau\zeta})a_j.
\label{124}
\end{equation}


\subsection{Factorization of circular fourcomplex polynomials}

A polynomial of degree $m$ of the circular fourcomplex variable 
$u=x+\alpha y+\beta z+\gamma t$ has the form
\begin{equation}
P_m(u)=u^m+a_1 u^{m-1}+\cdots+a_{m-1} u +a_m ,
\label{125}
\end{equation}
where the constants are in general circular fourcomplex numbers.

It can be shown that any circular fourcomplex polynomial has a circular fourcomplex root, whence
it follows that a polynomial of degree $m$ can be written as a product of
$m$ linear factors of the form $u-u_j$, where the circular fourcomplex numbers $u_j$ are
the roots of the polynomials, although the factorization may not be unique, 
\begin{equation}
P_m(u)=\prod_{j=1}^m (u-u_j) .
\label{126}
\end{equation}\index{roots!circular fourcomplex}

The fact that any circular fourcomplex polynomial has a root can be shown by considering
the transformation of a fourdimensional sphere with the center at the origin by
the function $u^m$. The points of the hypersphere of radius $d$ are of the form
written in Eq. (\ref{52}), with $d$ constant and $0\leq\phi<2\pi,
0\leq\chi<2\pi, 0\leq\psi\leq \pi/2$. The point $u^m$ is
\begin{equation}
u^m=d^m\exp\left(\alpha m\frac{\phi+\chi}{2}+\beta m\frac{\phi-\chi}{2}\right)
[\cos(\psi-\pi/4)+\gamma\sin(\psi-\pi/4)]^m .
\label{127}
\end{equation}
It can be shown with the aid of Eq. (\ref{56}) that
\begin{equation}
\left|u\exp
\left(\alpha \frac{\phi+\chi}{2}+\beta \frac{\phi-\chi}{2}\right)\right|
=|u|,
\label{127bb}
\end{equation}
so that
\begin{eqnarray}
\lefteqn{\left|[\cos(\psi-\pi/4)+\gamma\sin(\psi-\pi/4)]^m
\exp\left(\alpha m\frac{\phi+\chi}{2}+\beta m\frac{\phi-\chi}{2}\right)\right|
\nonumber}\\
&&=\left|[\cos(\psi-\pi/4)+\gamma\sin(\psi-\pi/4)]^m\right| .
\label{127b}
\end{eqnarray}
The right-hand side of Eq. (\ref{127b}) is
\begin{equation}
|(\cos\epsilon+\gamma\sin\epsilon)^m|^2
=\sum_{k=0}^m C_{2m}^{2k}\cos^{2m-2k}\epsilon\sin^{2k}\epsilon ,
\label{128}
\end{equation}
where $\epsilon=\psi-\pi/4$, 
and since $C_{2m}^{2k}\geq C_m^k$, it can be concluded that
\begin{equation}
|(\cos\epsilon+\gamma\sin\epsilon)^m|^2\geq 1 .
\label{129}
\end{equation}
Then
\begin{equation}
d^m\leq |u^m|\leq 2^{(m-1)/2} d^m ,
\label{129b}
\end{equation}
which shows that the image of a four-dimensional sphere via the transformation
operated by the function $u^m$ is a finite hypersurface.

If $u^\prime=u^m$, and
\begin{equation}
u^\prime=d^\prime
[\cos(\psi^\prime-\pi/4)+\gamma\sin(\psi^\prime-\pi/4)]
\exp\left(\alpha \frac{\phi^\prime+\chi^\prime}{2}+\beta
\frac{\phi^\prime-\chi^\prime}{2}\right),
\label{130}
\end{equation}
then 
\begin{equation}
\phi^\prime=m\phi, \: \chi^\prime=m\chi, \: \tan\psi^\prime=\tan^m\psi .
\label{131}
\end{equation}
Since for any values of the angles $\phi^\prime, \chi^\prime, \psi^\prime$
there is a set of solutions $\phi, \chi, \psi$ of Eqs. (\ref{131}), and since
the image of the hypersphere is a finite hypersurface, it follows that the
image of the four-dimensional sphere via the function $u^m$ is also a closed
hypersurface. A continuous hypersurface is called closed when any ray issued
from the 
origin intersects that surface at least once in the finite part of the space.

A transformation of the four-dimensional space by the polynomial $P_m(u)$
will be considered further. By this transformation, a hypersphere of radius $d$
having the center at the origin is changed into a certain finite closed
surface, as discussed previously. 
The transformation of the four-dimensional space by the polynomial $P_m(u)$
associates to the point $u=0$ the point $f(0)=a_m$, and the image of a
hypersphere of very large radius $d$ can be represented with good approximation
by the image of that hypersphere by the function $u^m$. 
The origin of the axes is an inner
point of the latter image. If the radius of the hypersphere is now reduced
continuously from the initial very large values to zero, the image hypersphere
encloses initially the origin, but the image shrinks to $a_m$ when the radius
approaches the value zero.  Thus, the
origin is initially inside the image hypersurface, and it lies outside the
image hypersurface when the radius of the hypersphere tends to zero. Then since
the image hypersurface is closed, the image surface must intersect at some
stage the origin of the axes, which means that there is a point $u_1$ such that
$f(u_1)=0$. The factorization in Eq. (\ref{126}) can then be obtained by
iterations.

The roots of the polynomial $P_m$ can be obtained by the following method.
If the constants in Eq. (\ref{125}) are $a_l=a_{l0}+\alpha a_{l1}
+\beta a_{l2}+\gamma a_{l3}$, and
with the 
notations of Eq. (\ref{n88b}), the polynomial $P_m(u)$ can be written as
\index{polynomial, canonical variables!circular fourcomplex}
\begin{eqnarray}
\lefteqn{P_m=\sum_{l=0}^{m} 2^{(m-l)/2}
(e_1 A_{l1}+\tilde e_1\tilde A_{l1})(e_1 \xi+\tilde e_1 \upsilon)^{m-l}\nonumber}\\
&&+\sum_{l=0}^{m} 2^{(m-l)/2}
(e_2 A_{l2}+\tilde e_2\tilde A_{l2})(e_2 \tau+\tilde e_2 \zeta)^{m-l} ,
\label{126a}
\end{eqnarray}
where the constants $A_{lk}, \tilde A_{lk}, k=1,2,$ are real numbers.
Each of the polynomials of degree $m$ in $e_1 \xi+\tilde e_1\upsilon, 
e_2 \tau+\tilde e_2\zeta$
in Eq. (\ref{126a}) 
can always be written as a product of linear factors of the form
$e_1 (\xi-\xi_p)+\tilde e_1(\upsilon- \upsilon_p)$ and respectively
$e_2 (\tau-\tau_p)+\tilde e_2(\zeta- \zeta_p)$, where the
constants $\xi_p, \upsilon_p, \tau_p, \zeta_p$ are real,
\begin{eqnarray}
\lefteqn{\sum_{l=0}^{m} 2^{(m-l)/2}
(e_1 A_{l1}+\tilde e_1\tilde A_{l1})(e_1 \xi+\tilde e_1 \upsilon)^{m-l}
=\prod_{p=1}^{m}2^{m/2}\left\{e_1 (\xi-\xi_p)+\tilde e_1(\upsilon- \upsilon_p)
\right\},\nonumber}\\
&&
\label{126bb}
\end{eqnarray}
\begin{eqnarray}
\lefteqn{\sum_{l=0}^{m} 2^{(m-l)/2}
(e_2 A_{l2}+\tilde e_2\tilde A_{l2})(e_2 \tau+\tilde e_2 \zeta)^{m-l}
=\prod_{p=1}^{m}2^{m/2}\left\{e_2 (\tau-\tau_p)+\tilde e_2(\zeta- \zeta_p)
\right\}.\nonumber}\\
&&
\label{126bc}
\end{eqnarray}

Due to the relations  (\ref{87d}),
the polynomial $P_m(u)$ can be written as a product of factors of
the form 
\begin{eqnarray}
P_m(u)=\prod_{p=1}^m 2^{m/2}\left\{e_1 (\xi-\xi_p)+\tilde e_1(\upsilon- \upsilon_p)
+e_2 (\tau-\tau_p)+\tilde e_2(\zeta- \zeta_p)\right\}.
\label{128b}
\end{eqnarray}\index{polynomial, factorization!circular fourcomplex}
This relation can be written with the aid of Eq. (\ref{87b}) in the form 
(\ref{126}), where
\begin{eqnarray}
u_p=\sqrt{2}(e_1 \xi_p+\tilde e_1 \upsilon_p
+e_2 \tau_p+\tilde e_2 \zeta_p) .
\label{129bx}
\end{eqnarray}
The roots  
$e_1 \xi_p+\tilde e_1 \upsilon_p$ and $e_2 \tau_p+\tilde e_2 \zeta_p$
defined in Eqs. (\ref{126bb}) and respectively (\ref{126bc}) may be ordered
arbitrarily. This means that Eq. 
(\ref{129bx}) gives sets of $m$ roots 
$u_1,...,u_m$ of the polynomial $P_m(u)$, 
corresponding to the various ways in which the roots 
$e_1 \xi_p+\tilde e_1 \upsilon_p$ and $e_2 \tau_p+\tilde e_2 \zeta_p$
are ordered according to $p$ for each polynomial. 
Thus, while the hypercomplex components in Eqs. (\ref{126bb}), Eqs.
(\ref{126bc}) taken separately have unique factorizations, the polynomial
$P_m(u)$ can be written in many different ways as a product of linear factors. 
The result of the circular fourcomplex integration, Eq. (\ref{124}), is however unique. 

If, for example, $P(u)=u^2+1$, the possible factorizations are
$P=(u-\tilde e_1-\tilde e_2)(u+\tilde e_1+\tilde e_2)$ and
$P=(u-\tilde e_1+\tilde e_2)(u+\tilde e_1-\tilde e_2)$, which can also be
written as $u^2+1=(u-\alpha)(u+\alpha)$ or as
$u^2+1=(u-\beta)(u+\beta)$. The result of the circular fourcomplex integration, Eq.
(\ref{124}), is however unique. 
It can be checked
that $(\pm \tilde e_1\pm\tilde e_2)^2=
-e_1-e_2=-1$.

\subsection{Representation of circular 
fourcomplex numbers by irreducible matrices}

If $T$ is the unitary matrix,
\begin{equation}
T =\left(
\begin{array}{cccc}
\frac{1}{\sqrt{2}}&0                 &0                 &\frac{1}{\sqrt{2}}\\
0                 &\frac{1}{\sqrt{2}}&\frac{1}{\sqrt{2}}&   0              \\
\frac{1}{\sqrt{2}}& 0                & 0                &-\frac{1}{\sqrt{2}} \\
0                 &\frac{1}{\sqrt{2}}&-\frac{1}{\sqrt{2}}&   0               \\
\end{array}
\right),
\label{129x}
\end{equation}
it can be shown 
that the matrix $T U T^{-1}$ has the form 
\begin{equation}
T U T^{-1}=\left(
\begin{array}{cc}
V_1      &     0    \\
0        &     V_2  \\
\end{array}
\right),
\label{129y}
\end{equation}\index{representation by irreducible matrices!circular fourcomplex}
where $U$ is the matrix in Eq. (\ref{23}) used to represent the circular fourcomplex
number $u$. In Eq. (\ref{129y}), $V_1, V_2$ are
the matrices
\begin{equation}
V_1=\left(
\begin{array}{cc}
x+t    &   y+z   \\
-y-z   &   x+t   \\
\end{array}\right),\;\;
V_2=\left(
\begin{array}{cc}
x-t    &   y-z   \\
-y+z   &   x-t   \\
\end{array}\right).
\label{130x}
\end{equation}
In Eq. (\ref{129y}), the symbols 0 denote the matrix
\begin{equation}
\left(
\begin{array}{cc}
0   &  0   \\
0   &  0   \\
\end{array}\right).
\label{131x}
\end{equation}
The relations between the variables $x+t,y+z$,$x-t,y-z$ for the multiplication
of circular fourcomplex numbers have been written in Eqs. (\ref{17})-(\ref{20}). The
matrix 
$T U T^{-1}$ provides an irreducible representation
\cite{4} of the circular fourcomplex number $u$ in terms of matrices with real
coefficients.

\section{Hyperbolic complex numbers in four dimensions}

\subsection{Operations with hyperbolic fourcomplex numbers}


A hyperbolic fourcomplex number is determined by its four components
$(x,y,z,t)$. The sum 
of the hyperbolic fourcomplex numbers $(x,y,z,t)$ and
$(x^\prime,y^\prime,z^\prime,t^\prime)$ is the hyperbolic fourcomplex
number $(x+x^\prime,y+y^\prime,z+z^\prime,t+t^\prime)$. \index{sum!hyperbolic fourcomplex}
The product of the hyperbolic fourcomplex numbers
$(x,y,z,t)$ and $(x^\prime,y^\prime,z^\prime,t^\prime)$ 
is defined in this section to be the hyperbolic fourcomplex
number
$(xx^\prime+yy^\prime+zz^\prime+tt^\prime,
xy^\prime+yx^\prime+zt^\prime+tz^\prime,
xz^\prime+zx^\prime+yt^\prime+ty^\prime,
xt^\prime+tx^\prime+yz^\prime+zy^\prime)$.\index{product!hyperbolic fourcomplex}


Hyperbolic fourcomplex numbers and their operations can be represented by  writing the
hyperbolic fourcomplex number $(x,y,z,t)$ as  
$u=x+\alpha y+\beta z+\gamma t$, where $\alpha, \beta$ and $\gamma$ 
are bases for which the multiplication rules are 
\begin{equation}
\alpha^2=1, \:\beta^2=1, \:\gamma^2=1, \alpha\beta=\beta\alpha=\gamma,\:
\alpha\gamma=\gamma\alpha=\beta, \:\beta\gamma=\gamma\beta=\alpha .
\label{h1}
\end{equation}\index{complex units!hyperbolic fourcomplex}
Two hyperbolic fourcomplex numbers $u=x+\alpha y+\beta z+\gamma t, 
u^\prime=x^\prime+\alpha y^\prime+\beta z^\prime+\gamma t^\prime$ are equal, 
$u=u^\prime$, if and only if $x=x^\prime, y=y^\prime,
z=z^\prime, t=t^\prime$. 
If 
$u=x+\alpha y+\beta z+\gamma t, 
u^\prime=x^\prime+\alpha y^\prime+\beta z^\prime+\gamma t^\prime$
are hyperbolic fourcomplex numbers, 
the sum $u+u^\prime$ and the 
product $uu^\prime$ defined above can be obtained by applying the usual
algebraic rules to the sum 
$(x+\alpha y+\beta z+\gamma t)+ 
(x^\prime+\alpha y^\prime+\beta z^\prime+\gamma t^\prime)$
and to the product 
$(x+\alpha y+\beta z+\gamma t)
(x^\prime+\alpha y^\prime+\beta z^\prime+\gamma t^\prime)$,
and grouping of the resulting terms,
\begin{equation}
u+u^\prime=x+x^\prime+\alpha(y+y^\prime)+\beta(z+z^\prime)+\gamma(t+t^\prime),
\label{h1a}
\end{equation}
\begin{eqnarray}
\lefteqn{uu^\prime=
xx^\prime+yy^\prime+zz^\prime+tt^\prime+
\alpha(xy^\prime+yx^\prime+zt^\prime+tz^\prime)+
\beta(xz^\prime+zx^\prime+yt^\prime+ty^\prime)\nonumber}\\
&&+\gamma(xt^\prime+tx^\prime+yz^\prime+zy^\prime).
\label{h1b}
\end{eqnarray}

If $u,u^\prime,u^{\prime\prime}$ are hyperbolic fourcomplex numbers, the multiplication is associative
\begin{equation}
(uu^\prime)u^{\prime\prime}=u(u^\prime u^{\prime\prime})
\label{h2}
\end{equation}
and commutative
\begin{equation}
u u^\prime=u^\prime u ,
\label{h3}
\end{equation}
as can be checked through direct calculation.
The hyperbolic fourcomplex zero is $0+\alpha\cdot 0+\beta\cdot 0+\gamma\cdot 0,$ 
denoted simply 0, 
and the hyperbolic fourcomplex unity is $1+\alpha\cdot 0+\beta\cdot 0+\gamma\cdot 0,$ 
denoted simply 1.

The inverse of the hyperbolic fourcomplex number 
$u=x+\alpha y+\beta z+\gamma t$ is a hyperbolic fourcomplex number
$u^\prime=x^\prime+\alpha y^\prime+\beta z^\prime+\gamma t^\prime$
having the property that
\begin{equation}
uu^\prime=1 .
\label{h4}
\end{equation}
Written on components, the condition, Eq. (\ref{h4}), is
\begin{equation}
\begin{array}{c}
xx^\prime+yy^\prime+zz^\prime+tt^\prime=1,\\
yx^\prime+xy^\prime+tz^\prime+zt^\prime=0,\\
zx^\prime+ty^\prime+xz^\prime+yt^\prime=0,\\
tx^\prime+zy^\prime+yz^\prime+xt^\prime=0 .
\end{array}
\label{h5}
\end{equation}
The system (\ref{h5}) has the solution\index{inverse!hyperbolic fourcomplex}
\begin{equation}
x^\prime=\frac{x(x^2-y^2-z^2-t^2)+2yzt}{\nu} ,
\label{h6a}
\end{equation}
\begin{equation}
y^\prime=
\frac{y(-x^2+y^2-z^2-t^2)+2xzt}{\nu} ,
\label{h6b}
\end{equation}
\begin{equation}
z^\prime=\frac{z(-x^2-y^2+z^2-t^2)+2xyt}{\nu} ,
\label{h6c}
\end{equation}
\begin{equation}
t^\prime=\frac{t(-x^2-y^2-z^2+t^2)+2xyz}{\nu} ,
\label{h6d}
\end{equation}
provided that $\nu\not= 0$, where
\begin{equation}
\nu=x^4+y^4+z^4+t^4-2(x^2y^2+x^2z^2+x^2t^2+y^2z^2+y^2t^2+z^2t^2)+8xyzt .
\label{h6e}
\end{equation}\index{inverse, determinant!hyperbolic fourcomplex}

The quantity $\nu$ can be written as
\begin{equation}
\nu=ss^\prime s^{\prime\prime}s^{\prime\prime\prime} ,
\label{h7}
\end{equation}
where
\begin{equation}
s=x+y+z+t, \: s^\prime= x-y+z-t , \: s^{\prime\prime}=x+y-z-t,  \:
s^{\prime\prime\prime}=x-y-z+t .
\label{h8}
\end{equation}
The variables $s, s^\prime, s^{\prime\prime}, s^{\prime\prime\prime}$ will be
called canonical  
hyperbolic fourcomplex variables.\index{canonical variables!hyperbolic fourcomplex}

Then a hyperbolic fourcomplex number $u=x+\alpha y+\beta z+\gamma t$ has an inverse,
unless 
\begin{equation}
s=0 ,\:\:{\rm or}\:\: s^\prime=0, \:\:{\rm or}\:\: s^{\prime\prime}=0, \:\:
{\rm or}\:\: s^{\prime\prime\prime}=0 . 
\label{h9}
\end{equation}

For arbitrary values of the variables $x,y,z,t$, the quantity $\nu$ can be
positive or negative. If $\nu\geq 0$, the quantity $\mu=\nu^{1/4}$
will be called amplitude of the hyperbolic fourcomplex number
$x+\alpha y+\beta z +\gamma t$.\index{amplitude!hyperbolic fourcomplex}
The normals of the hyperplanes in Eq. (\ref{h9}) are orthogonal to
each other. Because of conditions (\ref{h9}) these hyperplanes will
be also called the nodal hyperplanes. 
It can be shown that 
if $uu^\prime=0$ then either $u=0$, or $u^\prime=0$, or $q, q^\prime$ belong to
pairs of orthogonal hypersurfaces as described further.
Thus, divisors of zero exist if
one of the hyperbolic fourcomplex
numbers $u, u^\prime$ belongs to one of the nodal hyperplanes 
and the other hyperbolic fourcomplex 
number belongs to the straight line through the origin
which is normal to that hyperplane,
\begin{equation}
x+y+z+t=0,\:\:{\rm and} \:\: x^\prime=y^\prime=z^\prime=t^\prime,  
\label{h10a}
\end{equation}\index{divisors of zero!hyperbolic fourcomplex}
or
\begin{equation}
x-y+z-t=0, \:\:{\rm and}\:\: x^\prime=-y^\prime=z^\prime=-t^\prime, 
\label{h10b}
\end{equation}
or
\begin{equation}
x+y-z-t=0, \:\:{\rm and}\:\: x^\prime=y^\prime=-z^\prime=-t^\prime, 
\label{h10c}
\end{equation}
or
\begin{equation}
x-y-z+t=0, \:\:{\rm and}\:\: x^\prime=-y^\prime=-z^\prime=t^\prime.
\label{h10d}
\end{equation}
Divisors of zero also exist if the hyperbolic fourcomplex numbers $u,u^\prime$ belong to
different members of the pairs 
of two-dimensional hypersurfaces listed further,
\begin{equation}
x+y=0,\:z+t=0\:\: {\rm and}\:\:  x^\prime-y^\prime=0 , \: z^\prime-t^\prime=0,
\label{h11a}
\end{equation}
or
\begin{equation}
x+z=0,\:y+t=0\:\: {\rm and}\:\:  x^\prime-z^\prime=0 , \: y^\prime-t^\prime=0,
\label{h11b}
\end{equation}
or
\begin{equation}
y+z=0,\:x+t=0\:\: {\rm and}\:\:  y^\prime-z^\prime=0 , \: x^\prime-t^\prime=0.
\label{h11c}
\end{equation}


\subsection{Geometric representation of hyperbolic fourcomplex numbers}

The hyperbolic fourcomplex number $x+\alpha y+\beta z+\gamma t$ can be
represented by the point $A$ of coordinates $(x,y,z,t)$. 
If $O$ is the origin of the four-dimensional space $x,y,z,t,$ the distance 
from $A$ to the origin $O$ can be taken as
\begin{equation}
d^2=x^2+y^2+z^2+t^2 .
\label{h12}
\end{equation}\index{distance!hyperbolic fourcomplex}
The distance $d$ will be called modulus of the hyperbolic fourcomplex number $x+\alpha
y+\beta z +\gamma t$, $d=|u|$.\index{modulus!hyperbolic fourcomplex} 

If $u=x+\alpha y+\beta z +\gamma t, u_1=x_1+\alpha y_1+\beta z_1 +\gamma t_1,
u_2=x_2+\alpha y_2+\beta z_2 +\gamma t_2$, and $u=u_1u_2$, and if
\begin{equation}
s_j=x_j+y_j+z_j+t_j, \: s_j^\prime= x_j-y_j+z_j-t_j , \:
s_j^{\prime\prime}=x_j+y_j-z_j-t_j,  \: 
s_j^{\prime\prime\prime}=x_j-y_j-z_j+t_j , 
\label{h13}
\end{equation}
for $j=1,2$, it can be shown that
\begin{equation}
s=s_1s_2 ,\:\:
s^\prime=s_1^\prime s_2^\prime, \:\:
s^{\prime\prime}=s_1^{\prime\prime}s_2^{\prime\prime}, \:\:
s^{\prime\prime\prime}=s_1^{\prime\prime\prime}s_2^{\prime\prime\prime} . 
\label{h14}
\end{equation}\index{transformation of variables!hyperbolic fourcomplex}
The relations (\ref{h14}) are a consequence of the identities
\begin{eqnarray}
\lefteqn{(x_1x_2+y_1y_2+z_1z_2+t_1t_2)+(x_1y_2+y_1x_2+z_1t_2+t_1z_2)
\nonumber}\\
&&+(x_1z_2+z_1x_2+y_1t_2+t_1y_2)+(x_1t_2+t_1x_2+y_1z_2+z_1y_2)\nonumber\\
&&=(x_1+y_1+z_1+t_1)(x_2+y_2+z_2+t_2)
\label{h15}
\end{eqnarray}
\begin{eqnarray}
\lefteqn{(x_1x_2+y_1y_2+z_1z_2+t_1t_2)-(x_1y_2+y_1x_2+z_1t_2+t_1z_2)
\nonumber}\\
&&+(x_1z_2+z_1x_2+y_1t_2+t_1y_2)-(x_1t_2+t_1x_2+y_1z_2+z_1y_2)\nonumber\\
&&=(x_1-y_1+z_1-t_1)(x_2-y_2+z_2-t_2)
\label{h16}
\end{eqnarray}
\begin{eqnarray}
\lefteqn{(x_1x_2+y_1y_2+z_1z_2+t_1t_2)+(x_1y_2+y_1x_2+z_1t_2+t_1z_2)
\nonumber}\\
&&-(x_1z_2+z_1x_2+y_1t_2+t_1y_2)-(x_1t_2+t_1x_2+y_1z_2+z_1y_2)\nonumber\\
&&=(x_1+y_1-z_1-t_1)(x_2+y_2-z_2-t_2)
\label{h17}
\end{eqnarray}
\begin{eqnarray}
\lefteqn{(x_1x_2+y_1y_2+z_1z_2+t_1t_2)-(x_1y_2+y_1x_2+z_1t_2+t_1z_2)
\nonumber}\\
&&-(x_1z_2+z_1x_2+y_1t_2+t_1y_2)+(x_1t_2+t_1x_2+y_1z_2+z_1y_2)\nonumber\\
&&=(x_1-y_1-z_1+t_1)(x_2-y_2-z_2+t_2)
\label{h18}
\end{eqnarray}
A consequence of the relations (\ref{h14}) is that if $u=u_1u_2$, then
\begin{equation}
\nu=\nu_1\nu_2 ,
\label{h19}
\end{equation}
where
\begin{equation}
\nu_j=s_j s_j^\prime s_j^{\prime\prime} s_j^{\prime\prime\prime} , j=1,2.
\label{h20}
\end{equation} 

The hyperbolic fourcomplex numbers
\begin{eqnarray}
\lefteqn{e=\frac{1+\alpha+\beta+\gamma}{4},\:
e^\prime=\frac{1-\alpha+\beta-\gamma}{4},\:
e^{\prime\prime}=\frac{1+\alpha-\beta-\gamma}{4},\:
e^{\prime\prime\prime}=\frac{1-\alpha-\beta+\gamma}{4}\nonumber}\\
&&
\label{h21}
\end{eqnarray}\index{canonical base!hyperbolic fourcomplex}
are orthogonal,
\begin{equation}
ee^\prime=0,\:ee^{\prime\prime}=0,\:ee^{\prime\prime\prime}=0,\:e^\prime
e^{\prime\prime}=0,\:e^\prime
e^{\prime\prime\prime}=0,\:e^{\prime\prime}e^{\prime\prime\prime}=0,  
\label{h22}
\end{equation}
and have also the property that
\begin{equation}
e^2=e, \: e^{\prime 2}=e^{\prime}, \:
e^{\prime\prime 2}=e^{\prime\prime}, \:e^{\prime\prime\prime
2}=e^{\prime\prime\prime} . 
\label{h23}
\end{equation}
The hyperbolic fourcomplex number $u=x+\alpha y+\beta z+\gamma t$ can be written as
\begin{equation}
x+\alpha y+\beta z+\gamma t
=(x+y+z+t)e+(x-y+z-t)e^\prime+(x+y-z-t)e^{\prime\prime}+(x-y-z+t)e^{\prime\prime\prime},  
\label{h24}
\end{equation}
or, by using Eq. (\ref{h8}),
\begin{equation}
u=se+s^\prime
e^\prime+s^{\prime\prime}e^{\prime\prime}+s^{\prime\prime\prime}e^{\prime\prime\prime}. 
\label{h25}
\end{equation}\index{canonical form!hyperbolic fourcomplex}
The ensemble $e, e^\prime, e^{\prime\prime}, e^{\prime\prime\prime}$ will be
called the canonical 
hyperbolic fourcomplex base, and Eq. (\ref{h25}) gives the canonical form of the
hyperbolic fourcomplex number.
Thus, if $u_j=s_je+s_j^\prime e^\prime+s_j^{\prime\prime}e^{\prime\prime}
+s_j^{\prime\prime\prime}e^{\prime\prime\prime}, \:j=1,2$, and $u=u_1u_2$, then
the multiplication of the hyperbolic fourcomplex numbers is expressed by the
relations (\ref{h14}).
The moduli of the bases $e, e^\prime, e^{\prime\prime}, e^{\prime\prime\prime}$
are
\begin{equation}
|e|=\frac{1}{2},\; |e^\prime|=\frac{1}{2},\; 
|e^{\prime\prime}|=\frac{1}{2},\; |e^{\prime\prime\prime}|=\frac{1}{2}.
\label{h25b}
\end{equation}
The distance $d$, Eq. (\ref{h12}), is given by
\begin{equation}
d^2=\frac{1}{4}\left(s^2+s^{\prime 2}
+s^{\prime\prime 2}+s^{\prime\prime\prime 2}\right). 
\label{h25c}
\end{equation}\index{modulus, canonical variables!hyperbolic fourcomplex}

The relation (\ref{h25c}) shows that the variables $s, s^{\prime},
s^{\prime\prime}, s^{\prime\prime\prime}$ can be written as 
\begin{equation}
s=2d\cos\psi\cos\phi,\; s^\prime=2d\cos\psi\sin\phi,\; 
s^{\prime\prime}=2d\sin\psi\cos\chi,\; s=2d\sin\psi\sin\chi,
\label{h25d}
\end{equation}
where $\phi$ is the azimuthal angle in the $s,s^\prime$ plane,
$0\leq\phi<2\pi$, 
$\chi$ is the azimuthal angle in the $s^{\prime\prime},
s^{\prime\prime\prime}$ plane, $0\leq\chi<2\pi$,
and $\psi$ is the angle between the line $OA$ and the plane $ss^\prime$,
$0\leq\psi\leq\pi/2$.
The variables $x,y,z,t$ can be expressed in terms of the distance $d$ and the
angles $\phi, \chi, \psi$ as
\begin{equation}
\begin{array}{c}
x=(d/2)(\cos\psi\cos\phi+\cos\psi\sin\phi 
+\sin\psi\cos\chi+\sin\psi\sin\chi),\\
y=(d/2)(\cos\psi\cos\phi-\cos\psi\sin\phi 
+\sin\psi\cos\chi-\sin\psi\sin\chi),\\
z=(d/2)(\cos\psi\cos\phi+\cos\psi\sin\phi 
-\sin\psi\cos\chi-\sin\psi\sin\chi),\\
t=(d/2)(\cos\psi\cos\phi-\cos\psi\sin\phi 
-\sin\psi\cos\chi+\sin\psi\sin\chi).
\end{array}
\label{h25e}
\end{equation}

If $u=u_1u_2$, and the hypercomplex numbers $u_1, u_2$ are described by the
variables $d_1, \phi_1, \chi_1, \psi_1$ and respectively $d_2, \phi_2, \chi_2,
\psi_2$, then from Eq. (\ref{h25d}) it results that
\begin{eqnarray}
\lefteqn{\tan\phi=\tan\phi_1\tan\phi_2,\;
\tan\chi=\tan\chi_1\tan\chi_2,\nonumber}\\
&&\frac{\tan^2\psi\sin 2\chi}{\sin 2\phi}=
\frac{\tan^2\psi_1\sin 2\chi_1}{\sin 2\phi_1}
\frac{\tan^2\psi_2\sin 2\chi_2}{\sin 2\phi_2}.
\label{h25f}
\end{eqnarray}

The relation (\ref{h19}) for the product of hyperbolic fourcomplex numbers can 
be demonstrated also by using a representation of the multiplication of the 
hyperbolic fourcomplex numbers by matrices, in which the hyperbolic fourcomplex number $u=x+\alpha
y+\beta z+\gamma t$ is represented by the matrix
\begin{equation}
A=\left(\begin{array}{cccc}
x&y&z&t\\
y&x&t&z\\
z&t&x&y\\
t&z&y&x 
\end{array}\right) .
\label{h26}
\end{equation}\index{matrix representation!hyperbolic fourcomplex}
The product $u=x+\alpha y+\beta z+\gamma t$ of the hyperbolic fourcomplex numbers
$u_1=x_1+\alpha y_1+\beta z_1+\gamma t_1, u_2=x_2+\alpha y_2+\beta z_2+\gamma
t_2$, can be represented by the matrix multiplication 
\begin{equation}
A=A_1A_2.
\label{h27}
\end{equation}
It can be checked that the determinant ${\rm det}(A)$ of the matrix $A$ is
\begin{equation}
{\rm det}A = \nu .
\label{h28}
\end{equation}
The identity (\ref{h19}) is then a consequence of the fact the determinant 
of the product of matrices is equal to the product of the determinants 
of the factor matrices.

\subsection{Exponential form of a hyperbolic fourcomplex number}

The exponential function of a hypercomplex variable $u$ and the addition
theorem for the exponential function have been written in Eqs. 
~(1.35)-(1.36).
If $u=x+\alpha y+\beta z+\gamma t$, then  $\exp u$ can be calculated as
$\exp u=\exp x \cdot \exp (\alpha y) \cdot \exp (\beta z) \cdot \exp (\gamma
t)$. According to Eqs. (\ref{h1}), 
\begin{equation}
\alpha^{2m}=1, \alpha^{2m+1}=\alpha, 
\beta^{2m}=1, \beta^{2m+1}=\beta, 
\gamma^{2m}=1, \gamma^{2m+1}=\gamma, 
\label{h31}
\end{equation}\index{complex units, powers of!hyperbolic fourcomplex}
where $m$ is a natural number,
so that $\exp (\alpha y), \: \exp(\beta z)$ and $\exp(\gamma t)$ can be
written as 
\begin{equation}
\exp (\alpha y) = \cosh y +\alpha \sinh y , \:
\exp (\beta z) = \cosh y +\beta \sinh z , \:
\exp (\gamma t) = \cosh t +\gamma \sinh t . 
\label{h32}
\end{equation}\index{exponential, expressions!hyperbolic fourcomplex}
From Eqs. (\ref{h32}) it can be inferred that
\begin{eqnarray}
\lefteqn{(\cosh t +\alpha \sinh t)^m=\cosh mt +\alpha \sinh mt ,\:
(\cosh t +\beta \sinh t)^m=\cosh mt +\beta \sinh mt ,\nonumber}\\
&&(\cosh t +\gamma \sinh t)^m=\cosh mt +\gamma \sinh mt . 
\label{h33}
\end{eqnarray}

The hyperbolic fourcomplex numbers $u=x+\alpha y+\beta z+\gamma t$ for which
$s=x+y+z+t>0, \: s^\prime= x-y+z-t>0 , \: s^{\prime\prime}=x+y-z-t>0,  \:
s^{\prime\prime\prime}=x-y-z+t>0$ can be written in the form 
\begin{equation}
x+\alpha y+\beta z+\gamma t=e^{x_1+\alpha y_1+\beta z_1+\gamma t_1} .
\label{h34}
\end{equation}
The conditions $s=x+y+z+t>0, \: s^\prime= x-y+z-t>0 , \:
s^{\prime\prime}=x+y-z-t>0,  \: s^{\prime\prime\prime}=x-y-z+t>0$ 
correspond in Eq. (\ref{h25d}) to a range of angles $0<\phi<\pi/2,
0<\chi<\pi/2, 0<\psi\leq\pi/2$. 
The expressions of $x_1, y_1, z_1, t_1$ as functions of 
$x, y, z, t$ can be obtained by
developing $e^{\alpha y_1}, e^{\beta z_1}$and $e^{\gamma t_1}$ with the aid of
Eqs. (\ref{h32}), by multiplying these expressions and separating
the hypercomplex components, 
\begin{equation}
x=e^{x_1}(\cosh y_1\cosh z_1\cosh t_1+\sinh y_1\sinh z_1\sinh t_1) ,
\label{h35}
\end{equation}
\begin{equation}
y=e^{x_1}(\sinh y_1\cosh z_1\cosh t_1+\cosh y_1\sinh z_1\sinh t_1) ,
\label{h36}
\end{equation}
\begin{equation}
z=e^{x_1}(\cosh y_1\sinh z_1\cosh t_1+\sinh y_1\cosh z_1\sinh t_1) ,
\label{h37}
\end{equation}
\begin{equation}
t=e^{x_1}(\sinh y_1\sinh z_1\cosh t_1+\cosh y_1\cosh z_1\sinh t_1) ,
\label{h38}
\end{equation}
It can be shown from Eqs. (\ref{h35})-(\ref{h38}) that
\begin{equation}
x_1=\frac{1}{4} \ln(s s^\prime s^{\prime\prime}s^{\prime\prime\prime}) , \:
y_1=\frac{1}{4}\ln\frac{ss^{\prime\prime}}{s^\prime s^{\prime\prime\prime}},\:
z_1=\frac{1}{4}\ln\frac{ss^\prime}{s^{\prime\prime}s^{\prime\prime\prime}},\:
t_1=\frac{1}{4}\ln\frac{ss^{\prime\prime\prime}}{s^\prime s^{\prime\prime}}.
\label{h39}
\end{equation}
The exponential form of the hyperbolic fourcomplex number $u$ can be written
as
\begin{equation}
u=\mu\exp\left(
\frac{1}{4}\alpha \ln\frac{ss^{\prime\prime}}{s^\prime s^{\prime\prime\prime}}
+\frac{1}{4}\beta \ln\frac{ss^\prime}{s^{\prime\prime}s^{\prime\prime\prime}}
+\frac{1}{4}\gamma \ln\frac{ss^{\prime\prime\prime}}{s^\prime s^{\prime\prime}}
\right),
\label{h40}
\end{equation}\index{exponential form!hyperbolic fourcomplex}
where
\begin{equation}
\mu=(s s^\prime s^{\prime\prime}s^{\prime\prime\prime})^{1/4}.
\label{h40b}
\end{equation}
The exponential form of the hyperbolic fourcomplex number $u$ can be written
with the aid of the relations (\ref{h25d}) as
\begin{equation}
u=\mu\exp\left(\frac{1}{4}\alpha \ln\frac{1}{\tan\phi\tan\chi}
+\frac{1}{4}\beta\ln\frac{\sin 2\phi}{\tan^2\psi \sin 2\chi}
+\frac{1}{4}\gamma \ln\frac{\tan \chi}{\tan\phi}\right).
\label{h40x}
\end{equation}
The amplitude $\mu$ can be expressed in terms of the distance $d$ with the aid
of Eqs. (\ref{h25d}) as
\begin{equation}
\mu=d\sin^{1/2}2\psi\sin^{1/4}2\phi\sin^{1/4}2\chi.
\label{h40c}
\end{equation}
The hypercomplex number can be written as
\begin{eqnarray}
\lefteqn{u=d\sin^{1/2}2\psi\sin^{1/4}2\phi\sin^{1/4}2\chi
\exp\left(\frac{1}{4}\alpha \ln\frac{1}{\tan\phi\tan\chi}
+\frac{1}{4}\beta\ln\frac{\sin 2\phi}{\tan^2\psi \sin 2\chi}
+\frac{1}{4}\gamma \ln\frac{\tan \chi}{\tan\phi}\right),\nonumber}\\
&&
\label{h40d}
\end{eqnarray}\index{trigonometric form!hyperbolic fourcomplex}
which is the trigonometric form of the hypercomplex number $u$.


\subsection{Elementary functions of a hyperbolic fourcomplex variable}

The logarithm $u_1$ of the hyperbolic fourcomplex number $u$, $u_1=\ln u$, can be defined
for $s>0, s^\prime>0, s^{\prime\prime}>0, s^{\prime\prime\prime}>0$ 
as the solution of the equation
\begin{equation}
u=e^{u_1} ,
\label{h41}
\end{equation}
for $u_1$ as a function of
$u$. From Eq. (\ref{h40}) it results that 
\begin{equation}
\ln u=\frac{1}{4}\ln\mu+
+\frac{1}{4}\alpha \ln\frac{ss^{\prime\prime}}{s^\prime s^{\prime\prime\prime}}
+\frac{1}{4}\beta \ln\frac{ss^\prime}{s^{\prime\prime}s^{\prime\prime\prime}}
+\frac{1}{4}\gamma \ln\frac{ss^{\prime\prime\prime}}{s^\prime s^{\prime\prime}}.
\label{h42}
\end{equation}\index{logarithm!hyperbolic fourcomplex}
Using the expression in Eq. (\ref{h40x}), the logarithm can be written as
\begin{equation}
\ln u=\frac{1}{4}\ln\mu
+\frac{1}{4}\alpha \ln\frac{1}{\tan\phi\tan\chi}
+\frac{1}{4}\beta\ln\frac{\sin 2\phi}{\tan^2\psi \sin 2\chi}
+\frac{1}{4}\gamma \ln\frac{\tan \chi}{\tan\phi}.
\label{h42b}
\end{equation}

It can be inferred from Eqs. (\ref{h42}) and (\ref{h14}) that
\begin{equation}
\ln(u_1u_2)=\ln u_1+\ln u_2 .
\label{h43}
\end{equation}
The explicit form of Eq. (\ref{h42}) is 
\begin{eqnarray}
\lefteqn{\ln (x+\alpha y+\beta z+\gamma t)=
\frac{1}{4}(1+\alpha+\beta+\gamma)\ln(x+y+z+t)\nonumber}\\
&&+\frac{1}{4}(1-\alpha+\beta-\gamma)\ln(x-y+z-t)
+\frac{1}{4}(1+\alpha-\beta-\gamma)\ln(x+y-z-t)\nonumber\\
&&+\frac{1}{4}(1-\alpha-\beta+\gamma)\ln(x-y-z+t) .
\label{h45}
\end{eqnarray}
The relation (\ref{h45}) can be written with the aid of Eq. (\ref{h21}) as
\begin{equation}
\ln u = e\ln s + e^\prime\ln s^\prime+e^{\prime\prime}\ln s^{\prime\prime}
+e^{\prime\prime\prime}\ln s^{\prime\prime\prime}.
\label{h44}
\end{equation}

The power function $u^n$ can be defined for $s>0, s^\prime>0,
s^{\prime\prime}>0, 
s^{\prime\prime\prime}>0$ and real values of $n$ as
\begin{equation}
u^n=e^{n\ln u} .
\label{h46}
\end{equation}\index{power function!hyperbolic fourcomplex}
It can be inferred from Eqs. (\ref{h46}) and (\ref{h43}) that
\begin{equation}
(u_1u_2)^n=u_1^n\:u_2^n .
\label{h47}
\end{equation}
Using the expression (\ref{h45}) for $\ln u$ and the relations (\ref{h22}) and
(\ref{h23}) it can be shown that
\begin{eqnarray}
\lefteqn{(x+\alpha y+\beta z+\gamma t)^n=
\frac{1}{4}(1+\alpha+\beta+\gamma)(x+y+z+t)^n \nonumber}\\
&&+\frac{1}{4}(1-\alpha+\beta-\gamma)(x-y+z-t)^n
+\frac{1}{4}(1+\alpha-\beta-\gamma)(x+y-z-t)^n \nonumber\\
&&+\frac{1}{4}(1-\alpha-\beta+\gamma)(x-y-z+t)^n .
\label{h48}
\end{eqnarray}\index{power function!hyperbolic fourcomplex}
For integer $n$, the relation (\ref{h48}) is valid for any $x,y,z,t$. The
relation (\ref{h48}) for $n=-1$ is
\begin{eqnarray}
\lefteqn{\frac{1}{x+\alpha y+\beta z+\gamma t}=
\frac{1}{4}\left(\frac{1+\alpha+\beta+\gamma}{x+y+z+t}
+\frac{1-\alpha+\beta-\gamma}{x-y+z-t}
+\frac{1+\alpha-\beta-\gamma}{x+y-z-t}
+\frac{1-\alpha-\beta+\gamma}{x-y-z+t}\right) .\nonumber}\\
&&
\label{h49}
\end{eqnarray}

The trigonometric functions of the hypercomplex variable
$u$ and the addition theorems for these functions have been written in Eqs.
~(1.57)-(1.60). 
The cosine and sine functions of the hypercomplex variables $\alpha y, 
\beta z$ and $ \gamma t$ can be expressed as
\begin{equation}
\cos\alpha y=\cos y, \: \sin\alpha y=\alpha\sin y, 
\label{h54}
\end{equation}\index{trigonometric functions, expressions!hyperbolic fourcomplex}
\begin{equation}
\cos\beta y=\cos y, \: \sin\beta y=\beta\sin y, 
\label{h55}
\end{equation}
\begin{equation}
\cos\gamma y=\cos y, \: \sin\gamma y=\gamma\sin y .
\label{h56}
\end{equation}
The cosine and sine functions of a hyperbolic fourcomplex number $x+\alpha y+\beta
z+\gamma t$ can then be
expressed in terms of elementary functions with the aid of the addition
theorems Eqs. (1.59), (1.60) and of the expressions in  Eqs. 
(\ref{h54})-(\ref{h56}). 

The hyperbolic functions of the hypercomplex variable
$u$ and the addition theorems for these functions have been written in Eqs.
~(1.62)-(1.65). 
The hyperbolic cosine and sine functions of the hypercomplex variables $\alpha y, 
\beta z$ and $ \gamma t$ can be expressed as
\begin{equation}
\cosh\alpha y=\cosh y, \: \sinh\alpha y=\alpha\sinh y, 
\label{h61}
\end{equation}
\begin{equation}\index{hyperbolic functions, expressions!hyperbolic fourcomplex}
\cosh\beta y=\cosh y, \: \sinh\beta y=\beta\sinh y, 
\label{h62}
\end{equation}
\begin{equation}
\cosh\gamma y=\cosh y, \: \sinh\gamma y=\gamma\sinh y .
\label{h63}
\end{equation}
The hyperbolic cosine and sine functions of a hyperbolic fourcomplex number $x+\alpha y+\beta
z+\gamma t$ can then be
expressed in terms of elementary functions with the aid of the addition
theorems Eqs. (1.64), (1.65) and of the expressions in  Eqs. 
(\ref{h61})-(\ref{h63}).

\subsection{Power series of hyperbolic fourcomplex variables}

A hyperbolic fourcomplex series is an infinite sum of the form
\begin{equation}
a_0+a_1+a_2+\cdots+a_l+\cdots , 
\label{h64}
\end{equation}\index{series!hyperbolic fourcomplex}
where the coefficients $a_l$ are hyperbolic fourcomplex numbers. The convergence of 
the series (\ref{h64}) can be defined in terms of the convergence of its 4 real
components. The convergence of a hyperbolic fourcomplex series can however be studied
using hyperbolic fourcomplex variables. The main criterion for absolute convergence 
remains the comparison theorem, but this requires a number of inequalities
which will be discussed further.

The modulus of a hyperbolic fourcomplex number $u=x+\alpha y+\beta z+\gamma t$
can be defined as 
\begin{equation}
|u|=(x^2+y^2+z^2+t^2)^{1/2} ,
\label{h65}
\end{equation}\index{modulus, definition!hyperbolic fourcomplex}
so that according to Eq. (\ref{h12}) $d=|u|$. Since $|x|\leq |u|, |y|\leq |u|,
|z|\leq |u|, |t|\leq |u|$, a property of 
absolute convergence established via a comparison theorem based on the modulus
of the series (\ref{h64}) will ensure the absolute convergence of each real
component of that series.

The modulus of the sum $u_1+u_2$ of the hyperbolic fourcomplex numbers $u_1, u_2$ fulfils
the inequality
\begin{equation}
||u_1|-|u_2||\leq |u_1+u_2|\leq |u_1|+|u_2| .
\label{h66}
\end{equation}\index{modulus, inequalities!hyperbolic fourcomplex}
For the product the relation is 
\begin{equation}
|u_1u_2|\leq 2|u_1||u_2| ,
\label{h67}
\end{equation}
which replaces the relation of equality extant for regular complex numbers.
The equality in Eq. (\ref{h67}) takes place for $x_1^2=y_1^2=z_1^2=t_1^2$ and
$x_2/x_1=y_2/y_1=z_2/z_1=t_2/t_1$.
In particular
\begin{equation}
|u^2|\leq 2(x^2+y^2+z^2+t^2) .
\label{h68}
\end{equation}
The inequality in Eq. (\ref{h67}) implies that
\begin{equation}
|u^l|\leq 2^{l-1}|u|^l .
\label{h69}
\end{equation}
From Eqs. (\ref{h67}) and (\ref{h69}) it results that
\begin{equation}
|au^l|\leq 2^l |a| |u|^l .
\label{h70}
\end{equation}

A power series of the hyperbolic fourcomplex variable $u$ is a series of the form
\begin{equation}
a_0+a_1 u + a_2 u^2+\cdots +a_l u^l+\cdots .
\label{h71}
\end{equation}\index{power series!hyperbolic fourcomplex}
Since
\begin{equation}
\left|\sum_{l=0}^\infty a_l u^l\right| \leq  \sum_{l=0}^\infty
2^l|a_l| |u|^l ,
\label{h72}
\end{equation}
a sufficient condition for the absolute convergence of this series is that
\begin{equation}
\lim_{l\rightarrow \infty}\frac{2|a_{l+1}||u|}{|a_l|}<1 .
\label{h73}
\end{equation}
Thus the series is absolutely convergent for 
\begin{equation}
|u|<c_0,
\label{h74}
\end{equation}\index{convergence of power series!hyperbolic fourcomplex}
where 
\begin{equation}
c_0=\lim_{l\rightarrow\infty} \frac{|a_l|}{2|a_{l+1}|} .
\label{h75}
\end{equation}

The convergence of the series (\ref{h71}) can be also studied with the aid of
the formula (\ref{h48}) which, for integer values of $l$, is valid for any $x,
y, z, t$. If $a_l=a_{lx}+\alpha a_{ly}+\beta a_{lz}+\gamma a_{lt}$, and
\begin{eqnarray}
\label{h76a}
A_l=a_{lx}+a_{ly}+a_{lz}+a_{lt}, \\
A_l^\prime= a_{lx}-a_{ly}+a_{lz}-a_{lt} , \\
A_l^{\prime\prime}=a_{lx}+a_{ly}-a_{lz}-a_{lt},\\
A_l^{\prime\prime\prime}=a_{lx}-a_{ly}-a_{lz}+a_{lt} ,
\label{h76}
\end{eqnarray}
it can be shown with the aid of relations (\ref{h22}) and (\ref{h23}) that
\begin{equation}
a_l e=A_l e, \: a_l e^\prime=A_l^\prime e^\prime, \: 
a_l e^{\prime\prime}=A_l^{\prime\prime}e^{\prime\prime}, \:
a_l e^{\prime\prime\prime}=A_l^{\prime\prime\prime}e^{\prime\prime\prime} ,
\label{h77}
\end{equation}
so that the expression of the series (\ref{h71}) becomes
\begin{equation}
\sum_{l=0}^\infty \left(A_l  s^l e+
A_l^\prime  s^{\prime l}e^\prime+
A_l^{\prime\prime}s^{\prime\prime l}e^{\prime\prime}+
A_l^{\prime\prime\prime}s^{\prime\prime\prime l}e^{\prime\prime\prime}\right) ,
\label{h78}
\end{equation}
where the quantities $s, s^\prime, s^{\prime\prime}, s^{\prime\prime\prime}$
have been defined in Eq. (\ref{h8}).
The sufficient conditions for the absolute convergence of the series 
in Eq. (\ref{h78}) are that
\begin{equation}
\lim_{l\rightarrow \infty}\frac{|A_{l+1}||s|}{|A_l|}<1,
\lim_{l\rightarrow \infty}\frac{|A_{l+1}^\prime||s^\prime|}{|A_l^\prime|}<1,
\lim_{l\rightarrow \infty}\frac{|A_{l+1}^{\prime\prime}||s^{\prime\prime}|}{|A_l^{\prime\prime}|}<1,
\lim_{l\rightarrow \infty}\frac{|A_{l+1}^{\prime\prime\prime}||s^{\prime\prime\prime}|}{|A_l^{\prime\prime\prime}|}<1,
\label{h79}
\end{equation}
Thus the series in Eq. (\ref{h78}) is absolutely convergent for 
\begin{equation}
|x+y+z+t|<c,\:
|x-y+z-t|<c^\prime,\:
|x+y-z-t|<c^{\prime\prime},\:
|x-y-z+t|<c^{\prime\prime\prime},
\label{h80}
\end{equation}\index{convergence, region of!hyperbolic fourcomplex}
where 
\begin{equation}
c=\lim_{l\rightarrow\infty} \frac{|A_l|}{|A_{l+1}|} ,\:
c^\prime=\lim_{l\rightarrow\infty} \frac{|A_l^\prime|}{|A_{l+1}^\prime|} ,\:
c^{\prime\prime}=\lim_{l\rightarrow\infty} \frac{|A_l^{\prime\prime}|}{|A_{l+1}^{\prime\prime}|} ,\:
c^{\prime\prime\prime}=\lim_{l\rightarrow\infty} \frac{|A_l^{\prime\prime\prime}|}{|A_{l+1}^{\prime\prime\prime}|} .
\label{h81}
\end{equation}
The relations (\ref{h80}) show that the region of convergence of the series
(\ref{h78}) is a four-dimensional parallelepiped.
It can be shown that $c_0=(1/2){\rm min}(c,c^\prime,c^{\prime\prime},
c^{\prime\prime\prime})$, where
${\rm min}$ designates the smallest of the numbers $c, c^\prime,
c^{\prime\prime},c^{\prime\prime\prime}$.
Using Eq. (\ref{h25c}), it can be seen that the circular region of
convergence defined in Eqs. (\ref{h74}), (\ref{h75})
is included in the parallelogram defined in Eqs. (\ref{h80}) and (\ref{h81}).

\subsection{Analytic functions of hyperbolic fourcomplex variables}

The fourcomplex function $f(u)$ of the fourcomplex variable $u$ has
been expressed in Eq. (\ref{g16}) in terms of 
the real functions $P(x,y,z,t),Q(x,y,z,t),R(x,y,z,t), S(x,y,z,t)$ of real
variables $x,y,z,t$. The
relations between the partial derivatives of the functions $P, Q, R, S$ are
obtained by setting succesively in   
Eq. (\ref{g17}) $\Delta x\rightarrow 0, \Delta y=\Delta z=\Delta t=0$;
then $ \Delta y\rightarrow 0, \Delta x=\Delta z=\Delta t=0;$  
then $  \Delta z\rightarrow 0,\Delta x=\Delta y=\Delta t=0$; and finally
$ \Delta t\rightarrow 0,\Delta x=\Delta y=\Delta z=0 $. 
The relations are \index{relations between partial derivatives!hyperbolic fourcomplex}
\begin{equation}
\frac{\partial P}{\partial x} = \frac{\partial Q}{\partial y} =
\frac{\partial R}{\partial z} = \frac{\partial S}{\partial t},
\label{h89}
\end{equation}
\begin{equation}
\frac{\partial Q}{\partial x} = \frac{\partial P}{\partial y} =
\frac{\partial S}{\partial z} = \frac{\partial R}{\partial t},
\label{h90}
\end{equation}
\begin{equation}
\frac{\partial R}{\partial x} = \frac{\partial S}{\partial y} =
\frac{\partial P}{\partial z} = \frac{\partial Q}{\partial t},
\label{h91}
\end{equation}
\begin{equation}
\frac{\partial S}{\partial x} = \frac{\partial R}{\partial y} =
\frac{\partial Q}{\partial z} = \frac{\partial P}{\partial t}.
\label{h92}
\end{equation}

The relations (\ref{h89})-(\ref{h92}) are analogous to the Riemann relations
for the real and imaginary components of a complex function. It can be shown
from Eqs. (\ref{h89})-(\ref{h92}) that the component $P$ is a solution
of the equations \index{relations between second-order derivatives!hyperbolic fourcomplex}
\begin{equation}
\frac{\partial^2 P}{\partial x^2}-\frac{\partial^2 P}{\partial y^2}=0,
\:\: 
\frac{\partial^2 P}{\partial x^2}-\frac{\partial^2 P}{\partial z^2}=0,
\:\:
\frac{\partial^2 P}{\partial y^2}-\frac{\partial^2 P}{\partial t^2}=0,
\:\:
\frac{\partial^2 P}{\partial z^2}-\frac{\partial^2 P}{\partial t^2}=0,
\:\:
\label{h93}
\end{equation}
\begin{equation}
\frac{\partial^2 P}{\partial x^2}-\frac{\partial^2 P}{\partial t^2}=0,
\:\:
\frac{\partial^2 P}{\partial y^2}-\frac{\partial^2 P}{\partial z^2}=0,
\label{h94}
\end{equation}
and the components $Q, R, S$ are solutions of similar equations.
As can be seen from Eqs. (\ref{h93})-(\ref{h94}), the components $P, Q, R, S$ of
an analytic function of hyperbolic fourcomplex variable are solutions of the wave 
equation with respect to pairs of the variables $x,y,z,t$.
The component $P$ is also a solution of the mixed-derivative
equations 
\begin{equation}
\frac{\partial^2 P}{\partial x\partial y}=\frac{\partial^2 P}{\partial
z\partial t} ,
\:\: 
\frac{\partial^2 P}{\partial x\partial z}=\frac{\partial^2 P}{\partial
y\partial t} ,
\:\: 
\frac{\partial^2 P}{\partial x\partial t}=\frac{\partial^2 P}{\partial
y\partial z} ,
\end{equation}
and the components $Q, R, S$ are solutions of similar equations.

\subsection{Integrals of functions of hyperbolic fourcomplex variables}

The singularities of hyperbolic fourcomplex functions arise from terms of the form
$1/(u-u_0)^m$, with $m>0$. Functions containing such terms are singular not
only at $u=u_0$, but also at all points of the two-dimensional hyperplanes
passing through $u_0$ and which are parallel to the nodal hyperplanes. 

The integral of a hyperbolic fourcomplex function between two points $A, B$ along a path
situated in a region free of singularities is independent of path, which means
that the integral of an analytic function along a loop situated in a region
free from singularities is zero,
\begin{equation}
\oint_\Gamma f(u) du = 0,
\label{h105}
\end{equation}
where it is supposed that a surface $\Sigma$ spanning 
the closed loop $\Gamma$ is not intersected by any of
the two-dimensional hyperplanes associated with the
singularities of the function $f(u)$. Using the expression, Eq. (\ref{g16}),
for $f(u)$ and the fact that $du=dx+\alpha  dy+\beta dz+\gamma dt$, the
explicit form of the integral in Eq. (\ref{h105}) is
\begin{eqnarray}
\lefteqn{\oint _\Gamma f(u) du = \oint_\Gamma
[(Pdx+Qdy+Rdz+Sdt)+\alpha(Qdx+Pdy+Sdz+Rdt)\nonumber}\\
&&+\beta(Rdx+Sdy+Pdz+Qdt)+\gamma(Sdx+Rdy+Qdz+Pdt)] .
\label{h106}
\end{eqnarray}\index{integrals, path!hyperbolic fourcomplex}
If the functions $P, Q, R, S$ are regular on a surface $\Sigma$
spanning the loop $\Gamma$,
the integral along the loop $\Gamma$ can be transformed with the aid of the
theorem of Stokes in an integral over the surface $\Sigma$ of terms of the form
$\partial P/\partial y -  \partial Q/\partial x, \:\:
\partial P/\partial z - \partial R/\partial x,\:\:
\partial P/\partial t - \partial S/\partial x, \:\:
\partial Q/\partial z -  \partial R/\partial y, \:\:
\partial Q/\partial t - \partial S/\partial y,\:\:
\partial R/\partial t - \partial S/\partial z$ and of similar terms arising
from the $\alpha, \beta$ and $\gamma$ components, 
which are equal to zero by Eqs. (\ref{h89})-(\ref{h92}), and this proves Eq.
(\ref{h105}). 

The exponential form of the hyperbolic fourcomplex numbers, Eq. (\ref{h40x}), contains no
cyclic variable, and therefore the concept of residue is not applicable to the
hyperbolic fourcomplex numbers defined in Eqs. (\ref{h1}).

\subsection{Factorization of hyperbolic fourcomplex polynomials}

A polynomial of degree $m$ of the hyperbolic fourcomplex variable 
$u=x+\alpha y+\beta z+\gamma t$ has the form
\begin{equation}
P_m(u)=u^m+a_1 u^{m-1}+\cdots+a_{m-1} u +a_m ,
\label{h106b}
\end{equation}
where the constants are in general hyperbolic fourcomplex numbers.
If $a_m=a_{mx}+\alpha a_{my}+\beta a_{mz}+\gamma a_{mt}$, and with the
notations of Eqs. (\ref{h8}) and (\ref{h76a})-(\ref{h76}) applied for $l=0, 1, \cdots, m$ , the
polynomial $P_m(u)$ can be written as 
\begin{eqnarray}
\lefteqn{P_m= \left[s^m 
+A_1 s^{m-1}+\cdots+A_{m-1} s+ A_m \right] e
+\left[s^{\prime m} 
+A_1^\prime s^{\prime m-1} +\cdots+A_{m-1}^\prime s^\prime+ A_m^\prime \right]e^\prime\nonumber}\\
 & &
+\left[s^{\prime\prime m} 
+A_1^{\prime\prime} s^{\prime\prime m-1} +\cdots+A_{m-1}^{\prime\prime} s^{\prime\prime}+ A_m^{\prime\prime} \right]e^{\prime\prime}
+\left[s^{\prime\prime\prime m} 
+A_1^{\prime\prime\prime} s^{\prime\prime\prime m-1} +\cdots+A_{m-1}^{\prime\prime\prime} s^{\prime\prime\prime}+ A_m^{\prime\prime\prime} \right]e^{\prime\prime\prime}.\nonumber \\
&&
\label{h107}
\end{eqnarray}\index{polynomial, canonical variables!hyperbolic fourcomplex}
Each of the polynomials of degree $m$ with real coefficients in Eq. (\ref{h107})
can be written as a product
of linear or quadratic factors with real coefficients, or as a product of
linear factors which, if imaginary, appear always in complex conjugate pairs.
Using the latter form for the simplicity of notations, the polynomial $P_m$
can be written as
\begin{equation}
P_m=\prod_{l=1}^m (s-s_l)e
+\prod_{l=1}^m (s^\prime-s_l^\prime)e^\prime
+\prod_{l=1}^m (s^{\prime\prime}-s_l^{\prime\prime})e^{\prime\prime}
+\prod_{l=1}^m
(s^{\prime\prime\prime}-s_l^{\prime\prime\prime})e^{\prime\prime\prime} ,
\label{h108}
\end{equation}
where the quantities $s_l$ appear always in complex conjugate pairs, and the
same is true for the quantities $s_l^\prime$, for the quantitites $
s_l^{\prime\prime}$, and for the quantities $s_l^{\prime\prime\prime}$. 
Due to the properties in Eqs. (\ref{h22}) and (\ref{h23}),
the polynomial $P_m(u)$ can be written as a product of factors of
the form  
\begin{equation}
P_m(u)=\prod_{l=1}^m \left[(s-s_l)e
+(s^\prime-s_l^\prime)e^\prime
+(s^{\prime\prime}-s_l^{\prime\prime})e^{\prime\prime}
+(s^{\prime\prime\prime}-s_l^{\prime\prime\prime})e^{\prime\prime\prime}
\right].
\label{h109}
\end{equation}\index{polynomial, factorization!hyperbolic fourcomplex}
These relations can be written with the aid of Eq. (\ref{h24}) as
\begin{eqnarray}
P_m(u)=\prod_{p=1}^m (u-u_p) ,
\label{h128c}
\end{eqnarray}
where
\begin{eqnarray}
u_p=s_p e+s_p^{\prime}e^{\prime}
+s_p^{\prime\prime}e^{\prime\prime}
+s_p^{\prime\prime\prime}e^{\prime\prime\prime}.
\label{h128d}
\end{eqnarray}
The roots $s_p, s_p^{\prime}, s_p^{\prime\prime}, s_p^{\prime\prime\prime}$
of the corresponding polynomials in Eq. (\ref{h108}) may be ordered arbitrarily.
This means that Eq. (\ref{h128d}) gives sets of $m$ roots
$u_1,...,u_m$ of the polynomial $P_m(u)$, 
corresponding to the various ways in which the roots
$s_p, s_p^{\prime}, s_p^{\prime\prime}, s_p^{\prime\prime\prime}$
are ordered according to $p$ in each
group. Thus, while the hypercomplex components in Eq. (\ref{h107}) taken
separately have unique factorizations, the polynomial $P_m(u)$ can be written
in many different ways as a product of linear factors. 

If $P(u)=u^2-1$, 
the factorization in Eq. (\ref{h128c}) is $u^2-1=(u-u_1)(u-u_2)$, where 
$u_1=\pm  e\pm  e^\prime\pm  e^{\prime\prime}\pm  e^{\prime\prime\prime}$, 
$u_2=-u_1$, 
so that there are 8 distinct factorizations of $u^2-1$,
\begin{eqnarray}
\begin{array}{l}
u^2-1=\left(u-e-e^{\prime}-e^{\prime\prime}-e^{\prime\prime\prime}\right)
\left(u+e+e^{\prime}+e^{\prime\prime}+e^{\prime\prime\prime}\right),\\
u^2-1=\left(u-e-e^{\prime}-e^{\prime\prime}+e^{\prime\prime\prime}\right)
\left(u+e+e^{\prime}+e^{\prime\prime}-e^{\prime\prime\prime}\right),\\
u^2-1=\left(u-e-e^{\prime}+e^{\prime\prime}-e^{\prime\prime\prime}\right)
\left(u+e+e^{\prime}-e^{\prime\prime}+e^{\prime\prime\prime}\right),\\
u^2-1=\left(u-e+e^{\prime}-e^{\prime\prime}-e^{\prime\prime\prime}\right)
\left(u+e-e^{\prime}+e^{\prime\prime}+e^{\prime\prime\prime}\right),\\
u^2-1=\left(u-e-e^{\prime}+e^{\prime\prime}+e^{\prime\prime\prime}\right)
\left(u+e+e^{\prime}-e^{\prime\prime}-e^{\prime\prime\prime}\right),\\
u^2-1=\left(u-e+e^{\prime}-e^{\prime\prime}+e^{\prime\prime\prime}\right)
\left(u+e-e^{\prime}+e^{\prime\prime}-e^{\prime\prime\prime}\right),\\
u^2-1=\left(u-e+e^{\prime}+e^{\prime\prime}-e^{\prime\prime\prime}\right)
\left(u+e-e^{\prime}-e^{\prime\prime}+e^{\prime\prime\prime}\right),\\
u^2-1=\left(u-e+e^{\prime}+e^{\prime\prime}+e^{\prime\prime\prime}\right)
\left(u+e-e^{\prime}-e^{\prime\prime}-e^{\prime\prime\prime}\right).
\end{array}
\label{h129}
\end{eqnarray}
It can be checked that 
$\left\{\pm  e\pm  e^\prime\pm  e^{\prime\prime}
\pm  e^{\prime\prime\prime}\right\}^2= 
e+e^{\prime}+e^{\prime\prime}+e^{\prime\prime\prime}=1$.

\subsection{Representation of hyperbolic 
fourcomplex numbers by irreducible matrices}

If $T$ is the unitary matrix,
\begin{equation}
T =\left(
\begin{array}{cccc}
\frac{1}{{2}}&\frac{1}{{2}}  &\frac{1}{{2}}    &\frac{1}{{2}}    \\
\frac{1}{{2}}  &-\frac{1}{{2}}  &\frac{1}{{2}} &-\frac{1}{{2}}   \\
\frac{1}{{2}}&\frac{1}{{2}}  &-\frac{1}{{2}}    &-\frac{1}{{2}}    \\
\frac{1}{{2}}&-\frac{1}{{2}}  &-\frac{1}{{2}}    &\frac{1}{{2}}    \\
\end{array}
\right),
\label{h129x}
\end{equation}
it can be shown 
that the matrix $T U T^{-1}$ has the form 
\begin{equation}
T U T^{-1}=\left(
\begin{array}{cccc}
x+y+z+t &    0    &     0     &  0       \\
0       & x-y+z-t &     0     &  0       \\
0       &    0    &  x+y-z-t  &  0       \\
0       &    0    &     0     &  x-y-z+t \\
\end{array}
\right),
\label{h129y}
\end{equation}\index{representation by irreducible matrices!hyperbolic fourcomplex}
where $U$ is the matrix in Eq. (\ref{h26}) used to represent the hyperbolic fourcomplex
number $u$. 
The relations between the variables $x+y+z+t, x-y+z-t, x+y-z-t, x-y-z+t$
for the multiplication
of hyperbolic fourcomplex numbers have been written in Eqs. (\ref{h15})-(\ref{h18}).
The matrix $T U T^{-1}$ provides an irreducible representation
\cite{4} of the hyperbolic fourcomplex number $u$ in terms of matrices with real
coefficients.

\section{Planar complex numbers in four dimensions}

\subsection{Operations with planar fourcomplex numbers}

A planar fourcomplex number is determined by its four components $(x,y,z,t)$.
The sum of the planar fourcomplex numbers $(x,y,z,t)$ and
$(x^\prime,y^\prime,z^\prime,t^\prime)$ is the planar fourcomplex
number $(x+x^\prime,y+y^\prime,z+z^\prime,t+t^\prime)$. \index{sum!planar fourcomplex}
The product of the planar fourcomplex numbers
$(x,y,z,t)$ and $(x^\prime,y^\prime,z^\prime,t^\prime)$ 
is defined in this section to be the planar fourcomplex
number
$(xx^\prime-yt^\prime-zz^\prime-ty^\prime,
xy^\prime+yx^\prime-zt^\prime-tz^\prime,
xz^\prime+yy^\prime+zx^\prime-tt^\prime,
xt^\prime+yz^\prime+zy^\prime+tx^\prime)$.\index{product!planar fourcomplex}

Planar fourcomplex numbers and their operations can be represented by  writing the
planar fourcomplex number $(x,y,z,t)$ as  
$u=x+\alpha y+\beta z+\gamma t$, where $\alpha, \beta$ and $\gamma$ 
are bases for which the multiplication rules are 
\begin{equation}
\alpha^2=\beta, \:\beta^2=-1, \:\gamma^2=-\beta,
\alpha\beta=\beta\alpha=\gamma,\: 
\alpha\gamma=\gamma\alpha=-1, \:\beta\gamma=\gamma\beta=-\alpha .
\label{c2-1}
\end{equation}\index{complex units!planar fourcomplex}
Two planar fourcomplex numbers $u=x+\alpha y+\beta z+\gamma t, 
u^\prime=x^\prime+\alpha y^\prime+\beta z^\prime+\gamma t^\prime$ are equal, 
$u=u^\prime$, if and only if $x=x^\prime, y=y^\prime,
z=z^\prime, t=t^\prime$. 
If 
$u=x+\alpha y+\beta z+\gamma t, 
u^\prime=x^\prime+\alpha y^\prime+\beta z^\prime+\gamma t^\prime$
are planar fourcomplex numbers, 
the sum $u+u^\prime$ and the 
product $uu^\prime$ defined above can be obtained by applying the usual
algebraic rules to the sum 
$(x+\alpha y+\beta z+\gamma t)+ 
(x^\prime+\alpha y^\prime+\beta z^\prime+\gamma t^\prime)$
and to the product 
$(x+\alpha y+\beta z+\gamma t)
(x^\prime+\alpha y^\prime+\beta z^\prime+\gamma t^\prime)$,
and grouping of the resulting terms,
\begin{equation}
u+u^\prime=x+x^\prime+\alpha(y+y^\prime)+\beta(z+z^\prime)+\gamma(t+t^\prime),
\label{c2-1a}
\end{equation}
\begin{eqnarray}\index{sum!planar fourcomplex}
\lefteqn{uu^\prime=
xx^\prime-yt^\prime-zz^\prime-ty^\prime+
\alpha(xy^\prime+yx^\prime-zt^\prime-tz^\prime)+
\beta(xz^\prime+yy^\prime+zx^\prime-tt^\prime)\nonumber}\\
&&+\gamma(xt^\prime+yz^\prime+zy^\prime+tx^\prime).
\label{c2-1b}
\end{eqnarray}\index{product!planar fourcomplex}

If $u,u^\prime,u^{\prime\prime}$ are planar fourcomplex numbers, the multiplication 
is associative
\begin{equation}
(uu^\prime)u^{\prime\prime}=u(u^\prime u^{\prime\prime})
\label{c2-2}
\end{equation}
and commutative
\begin{equation}
u u^\prime=u^\prime u ,
\label{c2-3}
\end{equation}
as can be checked through direct calculation.
The planar fourcomplex zero is $0+\alpha\cdot 0+\beta\cdot 0+\gamma\cdot 0,$ 
denoted simply 0, 
and the planar fourcomplex unity is $1+\alpha\cdot 0+\beta\cdot 0+\gamma\cdot 0,$ 
denoted simply 1.

The inverse of the planar fourcomplex number 
$u=x+\alpha y+\beta z+\gamma t$ is a planar fourcomplex number
$u^\prime=x^\prime+\alpha y^\prime+\beta z^\prime+\gamma t^\prime$
having the property that
\begin{equation}
uu^\prime=1 .
\label{c2-4}
\end{equation}
Written on components, the condition, Eq. (\ref{c2-4}), is
\begin{equation}
\begin{array}{c}
xx^\prime-ty^\prime-zz^\prime-yt^\prime=1,\\
yx^\prime+xy^\prime-tz^\prime-zt^\prime=0,\\
zx^\prime+yy^\prime+xz^\prime-tt^\prime=0,\\
tx^\prime+zy^\prime+yz^\prime+xt^\prime=0.
\end{array}
\label{c2-5}
\end{equation}
The system (\ref{c2-5}) has the solution
\begin{equation}
x^\prime=\frac{x(x^2+z^2)-z(y^2-t^2)+2xyt}
{\rho^4} ,
\label{c2-6a}
\end{equation}
\begin{equation}\index{inverse!planar fourcomplex}
y^\prime=-\frac{y(x^2-z^2)+t(y^2+t^2)+2xzt}
{\rho^4} ,
\label{c2-6c}
\end{equation}
\begin{equation}
z^\prime=
\frac{-z(x^2+z^2)+x(y^2-t^2)+2zyt}
{\rho^4} ,
\label{c2-6b}
\end{equation}
\begin{equation}
t^\prime=-\frac{t(x^2-z^2)+y(y^2+t^2)-2xyz}
{\rho^4} ,
\label{c2-6d}
\end{equation}
provided that $\rho\not=0, $ where
\begin{equation}
\rho^4=x^4+z^4+y^4+t^4+2x^2z^2+2y^2t^2+4x^2yt-4xy^2z+4xzt^2-4yz^2t .
\label{c2-6e}
\end{equation}\index{inverse, determinant!planar fourcomplex}
\index{amplitude!planar fourcomplex}
The quantity $\rho$ will be called amplitude of the planar fourcomplex number
$x+\alpha y+\beta z +\gamma t$.
Since
\begin{equation}
\rho^4=\rho_+^2\rho_-^2 ,
\label{c2-7a}
\end{equation}
where
\begin{equation}
\rho_+^2=\left(x+\frac{y-t}{\sqrt{2}}\right)^2+\left(z+\frac{y+t}{\sqrt{2}}\right)^2,
\: \rho_-^2= \left(x-\frac{y-t}{\sqrt{2}}\right)^2+\left(z-\frac{y+t}{\sqrt{2}}\right)^2, 
\label{c2-7b}
\end{equation}
a planar fourcomplex number $u=x+\alpha y+\beta z+\gamma t$ has an inverse, unless
\begin{equation}
x+\frac{y-t}{\sqrt{2}}=0,\: z+\frac{y+t}{\sqrt{2}}=0 ,  
\label{c2-8}
\end{equation}
or
\begin{equation}
x-\frac{y-t}{\sqrt{2}}=0,\: z-\frac{y+t}{\sqrt{2}}=0 .  
\label{c2-9}
\end{equation}

Because of conditions (\ref{c2-8})-(\ref{c2-9}) these 2-dimensional hypersurfaces
will be called nodal hyperplanes. \index{nodal hyperplanes!planar fourcomplex}
It can be shown that 
if $uu^\prime=0$ then either $u=0$, or $u^\prime=0$, 
or one of the planar fourcomplex numbers is of the form $x+\alpha
(x+z)/\sqrt{2} +\beta z-\gamma (x-z)/\sqrt{2}$ 
and the other of the form $x^\prime-\alpha (x^\prime+z^\prime)/\sqrt{2}
+\beta z^\prime+\gamma (x^\prime-z^\prime)/\sqrt{2}$.

\subsection{Geometric representation of planar fourcomplex numbers}

The planar fourcomplex number $x+\alpha y+\beta z+\gamma t$ can be represented by 
the point $A$ of coordinates $(x,y,z,t)$. 
If $O$ is the origin of the four-dimensional space $x,y,z,t,$ the distance 
from $A$ to the origin $O$ can be taken as
\begin{equation}
d^2=x^2+y^2+z^2+t^2 .
\label{c2-10}
\end{equation}\index{distance!planar fourcomplex}
The distance $d$ will be called modulus of the planar fourcomplex number $x+\alpha
y+\beta z +\gamma t$, $d=|u|$. 
The orientation  in the four-dimensional space of the line $OA$ can be specified
with the aid of three angles $\phi, \chi, \psi$
defined with respect to the rotated system of axes
\begin{equation}
\xi=\frac{x}{\sqrt{2}}+\frac{y-t}{2}, \: 
\tau=\frac{x}{\sqrt{2}}-\frac{y-t}{2}, \:
\upsilon=\frac{z}{\sqrt{2}}+\frac{y+t}{2}, \: 
\zeta=-\frac{z}{\sqrt{2}}+\frac{y+t}{2} .
\label{c2-11}
\end{equation}\index{canonical variables!planar fourcomplex}
The variables $\xi, \upsilon, \tau, \zeta$ will be called canonical 
planar fourcomplex variables.
The use of the rotated axes $\xi, \upsilon, \tau, \zeta$ 
for the definition of the angles $\phi, \chi, \psi$ 
is convenient for the expression of the planar fourcomplex numbers
in exponential and trigonometric forms, as it will be discussed further.
The angle $\phi$ is the angle between the projection of $A$ in the plane
$\xi,\upsilon$ and the $O\xi$ axis, $0\leq\phi<2\pi$,  
$\chi$ is the angle between the projection of $A$ in the plane $\tau,\zeta$ and
the  $O\tau$ axis, $0\leq\chi<2\pi$,
and $\psi$ is the angle between the line $OA$ and the plane $\tau O \zeta
$, $0\leq \psi\leq\pi/2$, 
as shown in Fig. \ref{fig10}.
The definition of the variables in this section is different from the
definition used for the circular fourcomplex numbers, because the definition of
the rotated axes in Eq. (\ref{c2-11}) is different from the definition of the rotated
circular axes, Eq. (\ref{11}). 
The angles $\phi$ and $\chi$ will be called azimuthal
angles, \index{azimuthal angles!planar fourcomplex}
the angle $\psi$ will be called planar angle.\index{planar angle!planar fourcomplex}
The fact that $0\leq \psi\leq\pi/2$ means that $\psi$ has
the same sign on both faces of the two-dimensional hyperplane $\upsilon O
\zeta$. The components of the point $A$
in terms of the distance $d$ and the angles $\phi, \chi, \psi$ are thus
\begin{equation}
\frac{x}{\sqrt{2}}+\frac{y-t}{2}=d\cos\phi \sin\psi , 
\label{c2-12a}
\end{equation}
\begin{equation}
\frac{x}{\sqrt{2}}-\frac{y-t}{2}=d\cos\chi \cos\psi , 
\label{c2-12b}
\end{equation}
\begin{equation}
\frac{z}{\sqrt{2}}+\frac{y+t}{2}=d\sin\phi \sin\psi , 
\label{c2-12c}
\end{equation}
\newpage
\setlength{\oddsidemargin}{-0.2cm}      
\begin{equation}
-\frac{z}{\sqrt{2}}+\frac{y+t}{2}=d\sin\chi \cos\psi .
\label{c2-12d}
\end{equation}
It can be checked that $\rho_+=\sqrt{2}d\sin\psi, \rho_-=\sqrt{2}d\cos\psi$.
The coordinates $x,y,z,t$ in terms of the variables $d, \phi, \chi,
\psi$ are
\begin{equation}
x=\frac{d}{\sqrt{2}}(\cos\phi\sin\psi+\cos\chi\cos\psi),
\label{c2-12e}
\end{equation}
\begin{equation}
y=\frac{d}{\sqrt{2}}[\sin(\phi+\pi/4)\sin\psi+\sin(\chi-\pi/4)\cos\psi],
\label{c2-12g}
\end{equation}
\begin{equation}
z=\frac{d}{\sqrt{2}}(\sin\phi\sin\psi-\sin\chi\cos\psi),
\label{c2-12f}
\end{equation}
\begin{equation}
t=\frac{d}{\sqrt{2}}[-\cos(\phi+\pi/4)\sin\psi+\cos(\chi-\pi/4)\cos\psi].
\label{c2-12h}
\end{equation}
The angles $\phi, \chi, \psi$ can be expressed in terms of the coordinates
$x,y,z,t$ as 
\begin{equation}
\sin\phi = \frac{z+(y+t)/\sqrt{2}}{\rho_+} ,\: 
\cos\phi = \frac{x+(y-t)/\sqrt{2}}{\rho_+} ,
\label{c2-13a}
\end{equation}
\begin{equation}
\sin\chi = \frac{-z+(y+t)/\sqrt{2}}{\rho_-} ,\: 
\cos\chi = \frac{x-(y-t)/\sqrt{2}}{\rho_-} ,
\label{c2-13b}
\end{equation}
\begin{equation}
\tan\psi=\rho_+/\rho_- .
\label{c2-13c}
\end{equation}
The nodal hyperplanes are $\xi O\upsilon$, for which $\tau=0, \zeta=0$, and
$\tau O\zeta$, for which $\xi=0, \upsilon=0$.
For points in 
the nodal hyperplane $\xi O\upsilon$ the planar angle is $\psi=\pi/2$, 
for points in the nodal hyperplane $\tau O\zeta$ the planar angle is $\psi=0$.

It can be shown that if $u_1=x_1+\alpha y_1+\beta z_1+\gamma t_1, 
u_2=x_2+\alpha y_2+\beta z_2+\gamma t_2$ are planar fourcomplex
numbers of amplitudes and angles $\rho_1, \phi_1, \chi_1, \psi_1$ and
respectively $\rho_2, \phi_2, \chi_2, \psi_2$, then the amplitude $\rho$ and
the angles $\phi, \chi, \psi$ of the product planar fourcomplex number $u_1u_2$
are 
\begin{equation}
\rho=\rho_1\rho_2, 
\label{c2-14a}
\end{equation}\index{transformation of variables!planar fourcomplex}
\begin{equation}
 \phi=\phi_1+\phi_2, \: \chi=\chi_1+\chi_2, \: \tan\psi=\tan\psi_1\tan\psi_2 . 
\label{c2-14b}
\end{equation}
The relations (\ref{c2-14a})-(\ref{c2-14b}) are consequences of the definitions
(\ref{c2-6e})-(\ref{c2-7b}), (\ref{c2-13a})-(\ref{c2-13c}) and of the identities
\begin{eqnarray}
\lefteqn{\left[(x_1x_2-z_1z_2-y_1t_2-t_1y_2)
+\frac{(x_1y_2+y_1x_2-z_1t_2-t_1z_2)-(x_1t_2+t_1x_2+z_1y_2+y_1z_2)}{\sqrt{2}}
\right]^2
\nonumber}\\
&&+\left[(x_1z_2+z_1x_2+y_1y_2-t_1t_2)
+\frac{(x_1y_2+y_1x_2-z_1t_2-t_1z_2)+(x_1t_2+t_1x_2+z_1y_2+y_1z_2)}{\sqrt{2}}
\right]^2
\nonumber\\
&&=\left[\left(x_1+\frac{y_1-t_1}{\sqrt{2}}\right)^2
+\left(z_1+\frac{y_1+t_1}{\sqrt{2}}\right)^2\right]
\left[\left(x_2+\frac{y_2-t_2}{\sqrt{2}}\right)^2
+\left(z_2+\frac{y_2+t_2}{\sqrt{2}}\right)^2\right],
\label{c2-15}
\end{eqnarray}
\newpage
\setlength{\oddsidemargin}{0.9cm}      
\begin{eqnarray}
\lefteqn{\left[(x_1x_2-z_1z_2-y_1t_2-t_1y_2)
-\frac{(x_1y_2+y_1x_2-z_1t_2-t_1z_2)-(x_1t_2+t_1x_2+z_1y_2+y_1z_2)}{\sqrt{2}}
\right]^2
\nonumber}\\
&&+\left[(x_1z_2+z_1x_2+y_1y_2-t_1t_2)
-\frac{(x_1y_2+y_1x_2-z_1t_2-t_1z_2)+(x_1t_2+t_1x_2+z_1y_2+y_1z_2)}{\sqrt{2}}
\right]^2
\nonumber\\
&&=\left[\left(x_1-\frac{y_1-t_1}{\sqrt{2}}\right)^2
+\left(z_1-\frac{y_1+t_1}{\sqrt{2}}\right)^2\right]
\left[\left(x_2-\frac{y_2-t_2}{\sqrt{2}}\right)^2
+\left(z_2-\frac{y_2+t_2}{\sqrt{2}}\right)^2\right],
\label{c2-16}
\end{eqnarray}
\begin{eqnarray}
\lefteqn{(x_1x_2-z_1z_2-y_1t_2-t_1y_2)
+\frac{(x_1y_2+y_1x_2-z_1t_2-t_1z_2)-(x_1t_2+t_1x_2+z_1y_2+y_1z_2)}{\sqrt{2}}
\nonumber}\\
&&=\left(x_1+\frac{y_1-t_1}{\sqrt{2}}\right)
\left(x_2+\frac{y_2-t_2}{\sqrt{2}}\right)
-\left(z_1+\frac{y_1+t_1}{\sqrt{2}}\right)
\left(z_2+\frac{y_2+t_2}{\sqrt{2}}\right) ,
\label{c2-17}
\end{eqnarray}
\begin{eqnarray}
\lefteqn{(x_1z_2+z_1x_2+y_1y_2-t_1t_2)
+\frac{(x_1y_2+y_1x_2-z_1t_2-t_1z_2)+(x_1t_2+t_1x_2+z_1y_2+y_1z_2)}{\sqrt{2}}
\nonumber}\\
&&=\left(z_1+\frac{y_1+t_1}{\sqrt{2}}\right)
\left(x_2+\frac{y_2-t_2}{\sqrt{2}}\right)
+\left(x_1+\frac{y_1-t_1}{\sqrt{2}}\right)
\left(z_2+\frac{y_2+t_2}{\sqrt{2}}\right) ,
\label{c2-18}
\end{eqnarray}
\begin{eqnarray}
\lefteqn{(x_1x_2-z_1z_2-y_1t_2-t_1y_2)
-\frac{(x_1y_2+y_1x_2-z_1t_2-t_1z_2)-(x_1t_2+t_1x_2+z_1y_2+y_1z_2)}{\sqrt{2}}
\nonumber}\\
&&=\left(x_1-\frac{y_1-t_1}{\sqrt{2}}\right)
\left(x_2-\frac{y_2-t_2}{\sqrt{2}}\right)
-\left(-z_1+\frac{y_1+t_1}{\sqrt{2}}\right)
\left(-z_2+\frac{y_2+t_2}{\sqrt{2}}\right) ,
\label{c2-19}
\end{eqnarray}
\begin{eqnarray}
\lefteqn{-(x_1z_2+z_1x_2+y_1y_2-t_1t_2)
+\frac{(x_1y_2+y_1x_2-z_1t_2-t_1z_2)+(x_1t_2+t_1x_2+z_1y_2+y_1z_2)}{\sqrt{2}}
\nonumber}\\
&&=\left(-z_1+\frac{y_1+t_1}{\sqrt{2}}\right)
\left(x_2-\frac{y_2-t_2}{\sqrt{2}}\right)
+\left(x_1-\frac{y_1-t_1}{\sqrt{2}}\right)
\left(-z_2+\frac{y_2+t_2}{\sqrt{2}}\right) .
\label{c2-20}
\end{eqnarray}
The identities (\ref{c2-15}) and (\ref{c2-16}) can also be written as
\begin{equation}
\rho_+^2=\rho_{1+}\rho_{2+} ,
\label{c2-21a}
\end{equation}
\begin{equation}
\rho_-^2=\rho_{1-}\rho_{2-} ,
\label{c2-21b}
\end{equation}
where
\begin{equation}
\rho_{j+}^2=\left(x_j+\frac{y_j-t_j}{\sqrt{2}}\right)^2
+\left(z_j+\frac{y_j+t_j}{\sqrt{2}}\right)^2,
\: \rho_{j-}^2=\left(x_j-\frac{y_j-t_j}{\sqrt{2}}\right)^2
+\left(z-\frac{y_j+t_j}{\sqrt{2}}\right)^2 ,
\label{c2-22}
\end{equation}
for $j=1,2$.

The fact that the amplitude of the product is equal to the product of the 
amplitudes, as written in Eq. (\ref{c2-14a}), can 
be demonstrated also by using a representation of the multiplication of the 
planar fourcomplex numbers by matrices, in which the planar fourcomplex number $u=x+\alpha
y+\beta z+\gamma t$ is represented by the matrix
\begin{equation}
A=\left(\begin{array}{cccc}
x &y &z &t\\
-t&x &y &z\\
-z&-t&x &y\\
-y&-z&-t&x 
\end{array}\right) .
\label{c2-23}
\end{equation}\index{matrix representation!planar fourcomplex}
The product $u=x+\alpha y+\beta z+\gamma t$ of the planar fourcomplex numbers
$u_1=x_1+\alpha y_1+\beta z_1+\gamma t_1, u_2=x_2+\alpha y_2+\beta z_2+\gamma
t_2$, can be represented by the matrix multiplication 
\begin{equation}
A=A_1A_2.
\label{c2-24}
\end{equation}
It can be checked that the determinant ${\rm det}(A)$ of the matrix $A$ is
\begin{equation}
{\rm det}A = \rho^4 .
\label{c2-25}
\end{equation}
The identity (\ref{c2-14a}) is then a consequence of the fact the determinant 
of the product of matrices is equal to the product of the determinants 
of the factor matrices. 

\subsection{The planar fourdimensional cosexponential functions}

The exponential function of a hypercomplex variable $u$ and the addition
theorem for the exponential function have been written in Eqs. 
~(1.35)-(1.36).
If $u=x+\alpha y+\beta z+\gamma t$, then $\exp u$ can be calculated as 
$\exp u=\exp x \cdot \exp (\alpha y) \cdot \exp (\beta z)\cdot \exp (\gamma t)$.
According to Eq. (\ref{c2-1}),  
\begin{eqnarray}
\begin{array}{l}
\alpha^{8m}=1, \alpha^{8m+1}=\alpha, \alpha^{8m+2}=\beta, \alpha^{8m+3}=\gamma, \\
\alpha^{8m+4}=-1, \alpha^{8m+5}=-\alpha, \alpha^{8m+6}=-\beta, \alpha^{8m+7}=-\gamma, \\
\beta^{4m}=1, \beta^{4m+1}=\beta, \beta^{4m+2}=-1, \beta^{4m+3}=-\beta, \\
\gamma^{8m}=1, \gamma^{8m+1}=\gamma, \gamma^{8m+2}=-\beta, \gamma^{8m+3}=\alpha, \\
\gamma^{8m+4}=-1, \gamma^{8m+5}=-\gamma, \gamma^{8m+6}=\beta, \gamma^{8m+7}=-\alpha, \\
\end{array}
\label{c2-28}
\end{eqnarray}\index{complex units, powers of!planar fourcomplex}
where $n$ is a natural number,
so that $\exp (\alpha y), \: \exp(\beta z)$ and $\exp(\gamma t)$ can be written
as 
\begin{equation}
\exp (\beta z) = \cos z +\beta \sin z ,
\label{c2-29}
\end{equation}\index{exponential, expressions!planar fourcomplex}
and
\begin{equation}
\exp (\alpha y) = f_{40}(y)+\alpha f_{41}(y)+\beta f_{42}(y) +\gamma f_{43}(y) ,  
\label{c2-30a}
\end{equation}
\begin{equation}
\exp (\gamma t) = f_{40}(t)+\gamma f_{41}(t)-\beta f_{42}(t) +\alpha f_{43}(t) ,
\label{c2-30b}
\end{equation}
where the four-dimensional cosexponential functions $f_{40}, f_{41}, f_{42}, f_{43}$ are
defined by the series 
\begin{equation}
f_{40}(x)=1-x^4/4!+x^8/8!-\cdots ,
\label{c2-30c}
\end{equation}\index{cosexponential functions, planar fourcomplex!definitions}
\begin{equation}
f_{41}(x)=x-x^5/5!+x^9/9!-\cdots,
\label{c2-30d}
\end{equation}
\begin{equation}
f_{42}(x)=x^2/2!-x^6/6!+x^{10}/10!-\cdots ,
\label{c2-30e}
\end{equation}
\begin{equation}
f_{43}(x)=x^3/3!-x^7/7!+x^{11}/11!-\cdots .
\label{c2-30f}
\end{equation}
The functions $f_{40}, f_{42}$ are even, the functions $f_{41}, f_{43}$ are odd,
\begin{equation}
f_{40}(-u)=f_{40}(u),\:f_{42}(-u)=f_{42}(u),\:f_{41}(-u)=-f_{41}(u),\:f_{43}(-u)=-f_{43}(u).
\label{c2-30feo}
\end{equation}\index{cosexponential functions, planar fourcomplex!parity}

Addition theorems for the four-dimensional cosexponential functions can be
obtained from the relation $\exp \alpha(x+y)=\exp\alpha x\cdot\exp\alpha y $, by
substituting the expression of the exponentials as given in Eq. (\ref{c2-30a}),
\begin{equation}
f_{40}(x+y)=f_{40}(x)f_{40}(y)-f_{41}(x)f_{43}(y)-f_{42}(x)f_{42}(y)-f_{43}(x)f_{41}(y) ,
\label{c2-30g}
\end{equation}
\begin{equation}\index{cosexponential functions, planar fourcomplex!addition theorems}
f_{41}(x+y)=f_{40}(x)f_{41}(y)+f_{41}(x)f_{40}(y)-f_{42}(x)f_{43}(y)-f_{43}(x)f_{42}(y) ,
\label{c2-30h}
\end{equation}
\begin{equation}
f_{42}(x+y)=f_{40}(x)f_{42}(y)+f_{41}(x)f_{41}(y)+f_{42}(x)f_{40}(y)-f_{43}(x)f_{43}(y) ,
\label{c2-30i}
\end{equation}
\begin{equation}
f_{43}(x+y)=f_{40}(x)f_{43}(y)+f_{41}(x)f_{42}(y)+f_{42}(x)f_{41}(y)+f_{43}(x)f_{40}(y) .
\label{c2-30j}
\end{equation}
For $x=y$ the relations (\ref{c2-30g})-(\ref{c2-30j}) take the form
\begin{equation}
f_{40}(2x)=f_{40}^2(x)-f_{42}^2(x)-2f_{41}(x)f_{43}(x) ,
\label{c2-30gg}
\end{equation}
\begin{equation}
f_{41}(2x)=2f_{40}(x)f_{41}(x)-2f_{42}(x)f_{43}(x) ,
\label{c2-30hh}
\end{equation}
\begin{equation}
f_{42}(2x)=f_{41}^2(x)-f_{43}^2(x)+2f_{40}(x)f_{42}(x) ,
\label{c2-30ii}
\end{equation}
\begin{equation}
f_{43}(2x)=2f_{40}(x)f_{43}(x)+2f_{41}(x)f_{42}(x) .
\label{c2-30jj}
\end{equation}
For $x=-y$ the relations (\ref{c2-30g})-(\ref{c2-30j}) and (\ref{c2-30feo}) yield
\begin{equation}
f_{40}^2(x)-f_{42}^2(x)+2f_{41}(x)f_{43}(x)=1 ,
\label{c2-30eog}
\end{equation}
\begin{equation}
f_{41}^2(x)-f_{43}^2(x)-2f_{40}(x)f_{42}(x)=0 .
\label{c2-30eoi}
\end{equation}
From Eqs. (\ref{c2-29})-(\ref{c2-30b}) it can be shown that, for $m$ integer,
\begin{equation}
(\cos z +\beta \sin z)^m=\cos mz +\beta \sin mz ,
\label{c2-2929}
\end{equation}
and
\begin{equation}
[f_{40}(y)+\alpha f_{41}(y)+\beta f_{42}(y) +\gamma f_{43}(y)]^m= 
f_{40}(my)+\alpha f_{41}(my)+\beta f_{42}(my) +\gamma f_{43}(my),  
\label{c2-30a30a}
\end{equation}
\begin{equation}
[f_{40}(t)+\gamma f_{41}(t)-\beta f_{42}(t) +\alpha f_{43}(t)]^m=
f_{40}(mt)+\gamma f_{41}(mt)-\beta f_{42}(mt) +\alpha f_{43}(mt) .
\label{c2-30b30b}
\end{equation}

\newpage
\setlength{\oddsidemargin}{0.1cm}       
Since
\begin{equation}
(\alpha-\gamma)^{2m}=2^m,\: (\alpha-\gamma)^{2m+1}=2^m(\alpha-\gamma) ,
\label{c2-30k}
\end{equation}
it can be shown from the definition of the exponential function, Eq.
(1.35) that
\begin{equation}
\exp(\alpha-\gamma)x=\cosh\sqrt{2}x+\frac{\alpha-\gamma}{\sqrt{2}}
\sinh\sqrt{2}x .
\label{c2-30l}
\end{equation}
Substituting in the relation $\exp(\alpha-\gamma)x=\exp\alpha x\exp(-\gamma x)$ the
expression of the exponentials from Eqs. (\ref{c2-30a}), (\ref{c2-30b}) and
(\ref{c2-30l}) yields  
\begin{equation}
f_{40}^2+f_{41}^2+f_{42}^2+f_{43}^2=\cosh\sqrt{2}x ,
\label{c2-30m}
\end{equation}
\begin{equation}
f_{40}f_{41}-f_{40}f_{43}+f_{41}f_{42}+f_{42}f_{43}=\frac{1}{\sqrt{2}}\sinh\sqrt{2}x ,
\label{c2-30n}
\end{equation}
where $f_{40}, f_{41}, f_{42}, f_{43}$ are functions of $x$.
From relations (\ref{c2-30m}) and (\ref{c2-30n}) it can be inferred that
\begin{equation}
\left(f_{40}+\frac{f_{41}-f_{43}}{\sqrt{2}}\right)^2+
\left(f_{42}+\frac{f_{41}+f_{43}}{\sqrt{2}}\right)^2=\exp\sqrt{2}x ,
\label{c2-30o}
\end{equation}
\begin{equation}
\left(f_{40}-\frac{f_{41}-f_{43}}{\sqrt{2}}\right)^2+
\left(f_{42}-\frac{f_{41}+f_{43}}{\sqrt{2}}\right)^2=\exp(-\sqrt{2}x) ,
\label{c2-30p}
\end{equation}
which means that
\begin{eqnarray}
\lefteqn{\left[\left(f_{40}+\frac{f_{41}-f_{43}}{\sqrt{2}}\right)^2+
\left(f_{42}+\frac{f_{41}+f_{43}}{\sqrt{2}}\right)^2\right]
\left[\left(f_{40}-\frac{f_{41}-f_{43}}{\sqrt{2}}\right)^2+
\left(f_{42}-\frac{f_{41}+f_{43}}{\sqrt{2}}\right)^2\right]=1 .\nonumber}\\
&&
\label{c2-30q}
\end{eqnarray}
An equivalent form of Eq. (\ref{c2-30q}) is
\begin{equation}
f_{40}^4+f_{41}^4+f_{42}^4+f_{43}^4+2(f_{40}^2f_{42}^2+f_{41}^2f_{43}^2)
+4(f_{40}^2f_{41}f_{43}+
f_{40}f_{42}f_{43}^2-f_{40}f_{41}^2f_{42}-f_{41}f_{42}^2f_{43})=1.
\label{c2-30r}
\end{equation}
The form of this relation is similar to the expression in Eq. (\ref{c2-6e}).
Similarly, since
\begin{equation}
(\alpha+\gamma)^{2m}=(-1)^m2^m,\: (\alpha+\gamma)^{2m+1}=(-1)^m2^m(\alpha+\gamma) ,
\label{c2-30k30k}
\end{equation}
it can be shown from the definition of the exponential function, Eq.
(1.35) 
that
\begin{equation}
\exp(\alpha+\gamma)x=\cos\sqrt{2}x+\frac{\alpha+\gamma}{\sqrt{2}}
\sin\sqrt{2}x .
\label{c2-30l30l}
\end{equation}
Substituting in the relation $\exp(\alpha+\gamma)x=\exp\alpha x\exp\gamma x$ the
expression of the exponentials from Eqs. (\ref{c2-30a}), (\ref{c2-30b}) and
(\ref{c2-30l30l}) yields  
\begin{equation}
f_{40}^2-f_{41}^2+f_{42}^2-f_{43}^2=\cos\sqrt{2}x ,
\label{c2-30m30m}
\end{equation}
\newpage
\setlength{\oddsidemargin}{0.9cm}       
\begin{equation}
f_{40}f_{41}+f_{40}f_{43}-f_{41}f_{42}+f_{42}f_{43}=\frac{1}{\sqrt{2}}\sin\sqrt{2}x ,
\label{c2-30n30n}
\end{equation}
where $f_{40}, f_{41}, f_{42}, f_{43}$ are functions of $x$.

Expressions of the four-dimensional cosexponential functions
(\ref{c2-30c})-(\ref{c2-30f}) can be obtained using the fact that
$[(1+i)/\sqrt{2}]^4=-1$, so that\index{cosexponential functions, planar fourcomplex!expressions} 
\begin{equation}
f_{40}(x)=\frac{1}{2}\left(\cosh\frac{1+i}{\sqrt{2}}x+
\cos\frac{1+i}{\sqrt{2}}x\right),
\label{c2-30s}
\end{equation}
\begin{equation}
f_{41}(x)=\frac{1}{\sqrt{2}(1+i)}\left(\sinh\frac{1+i}{\sqrt{2}}x+
\sin\frac{1+i}{\sqrt{2}}x\right),
\label{c2-30t}
\end{equation}
\begin{equation}
f_{42}(x)=\frac{1}{2i}\left(\cosh\frac{1+i}{\sqrt{2}}x-
\cos\frac{1+i}{\sqrt{2}}x\right),
\label{c2-30u}
\end{equation}
\begin{equation}
f_{43}(x)=\frac{1}{\sqrt{2}(-1+i)}\left(\sinh\frac{1+i}{\sqrt{2}}x-
\sin\frac{1+i}{\sqrt{2}}x\right).
\label{c2-30w}
\end{equation}
Using the addition theorems for the functions in the right-hand sides of Eqs.
(\ref{c2-30s})-(\ref{c2-30w}), the expressions of the four-dimensional cosexponential
functions become
\begin{equation}
f_{40}(x)=\cos\frac{x}{\sqrt{2}}\cosh\frac{x}{\sqrt{2}} ,
\label{c2-30c30c}
\end{equation}
\begin{equation}
f_{41}(x)=\frac{1}{\sqrt{2}}
\left(\sin\frac{x}{\sqrt{2}}\cosh\frac{x}{\sqrt{2}}+
\sinh\frac{x}{\sqrt{2}}\cos\frac{x}{\sqrt{2}}\right),
\label{c2-30d30d}
\end{equation}
\begin{equation}
f_{42}(x)=\sin\frac{x}{\sqrt{2}}\sinh\frac{x}{\sqrt{2}}, 
\label{c2-30e30e}
\end{equation}
\begin{equation}
f_{43}(x)=\frac{1}{\sqrt{2}}
\left(\sin\frac{x}{\sqrt{2}}\cosh\frac{x}{\sqrt{2}}-
\sinh\frac{x}{\sqrt{2}}\cos\frac{x}{\sqrt{2}}\right) .
\label{c2-30f30f}
\end{equation}
It is remarkable that the series in Eqs. (\ref{c2-30c30c})-(\ref{c2-30f30f}), in
which the  
terms are either of the form $x^{4m}$, or $x^{4m+1}$, or $x^{4m+2}$,
or $x^{4m+3}$ can be expressed
in terms of elementary functions whose power series are not subject to such 
restrictions. 
The graphs of the four-dimensional cosexponential functions are shown in 
Fig. \ref{fig12}.

\begin{figure}
\begin{center}
\epsfig{file=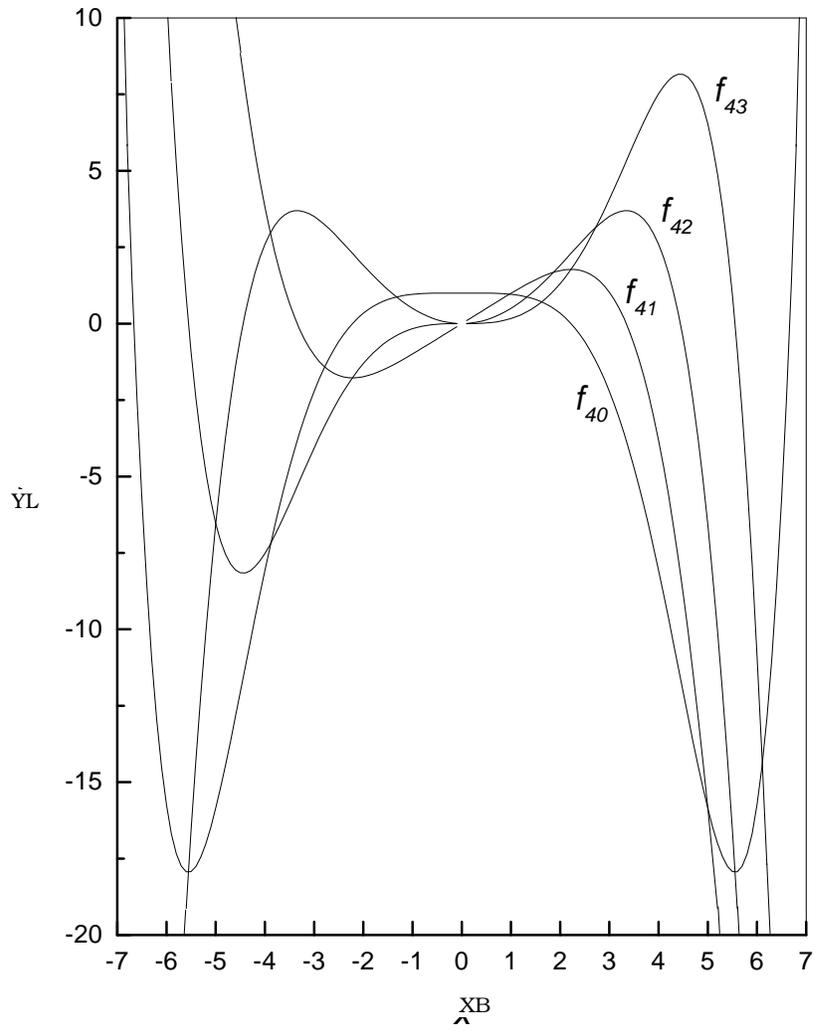,width=12cm}
\caption{The planar fourdimensional cosexponential functions $f_{40}, f_{41}, f_{42}, f_{43}$.}
\label{fig12}
\end{center}
\end{figure}

It can be checked that the cosexponential functions are solutions of the
fourth-order differential equation
\begin{equation}
\frac{{\rm d}^4\zeta}{{\rm d}u^4}=-\zeta ,
\label{c2-30x}
\end{equation}\index{cosexponential functions, planar fourcomplex!differential equations}
whose solutions are of the form $\zeta(u)=Af_{40}(u)+Bf_{41}(u)+Cf_{42}(u)+Df_{43}(u).$
It can also be checked that the derivatives of the cosexponential functions are
related by
\begin{equation}
\frac{df_{40}}{du}=-f_{43}, \:
\frac{df_{41}}{du}=f_{40}, \:
\frac{df_{42}}{du}=f_{41} ,
\frac{df_{43}}{du}=f_{42} .
\label{c2-30x30x}
\end{equation}


\subsection{The exponential and trigonometric forms of planar
fourcomplex numbers}

Any planar fourcomplex number $u=x+\alpha y+\beta z+\gamma t$ can be writen in the
form 
\begin{equation}
x+\alpha y+\beta z+\gamma t=e^{x_1+\alpha y_1+\beta z_1+\gamma t_1} .
\label{c2-31}
\end{equation}
The expressions of $x_1, y_1, z_1, t_1$ as functions of 
$x, y, z, t$ can be obtained by
developing $e^{\alpha y_1}, e^{\beta z_1}$ and $e^{\gamma t_1}$ with the aid of
Eqs. (\ref{c2-29})-(\ref{c2-30b}), by multiplying these expressions and separating
the hypercomplex components, and then substituting the expressions of the
four-dimensional cosexponential functions, Eqs. (\ref{c2-30c30c})-(\ref{c2-30f30f}), 
\begin{equation}
x=e^{x_1}\left(\cos z_1\cos \frac{y_1+t_1}{\sqrt{2}}\cosh \frac{y_1-t_1}{\sqrt{2}}
-\sin z_1\sin \frac{y_1+t_1}{\sqrt{2}}\sinh \frac{y_1-t_1}{\sqrt{2}}\right) ,
\label{c2-32}
\end{equation}
\begin{equation}
y=e^{x_1}\left[\sin \left(z_1+\frac{\pi}{4}\right)\cos \frac{y_1+t_1}{\sqrt{2}}\sinh \frac{y_1-t_1}{\sqrt{2}}
+\cos \left(z_1+\frac{\pi}{4}\right)\sin \frac{y_1+t_1}{\sqrt{2}}\cosh \frac{y_1-t_1}{\sqrt{2}}\right] ,
\label{c2-34}
\end{equation}
\begin{equation}
z=e^{x_1}\left(\cos z_1\sin \frac{y_1+t_1}{\sqrt{2}}\sinh \frac{y_1-t_1}{\sqrt{2}}
+\sin z_1\cos \frac{y_1+t_1}{\sqrt{2}}\cosh \frac{y_1-t_1}{\sqrt{2}}\right) ,
\label{c2-33}
\end{equation}
\begin{equation}
t=e^{x_1}\left[-\cos \left(z_1+\frac{\pi}{4}\right)\cos \frac{y_1+t_1}{\sqrt{2}}\sinh \frac{y_1-t_1}{\sqrt{2}}
+\sin \left(z_1+\frac{\pi}{4}\right)\sin \frac{y_1+t_1}{\sqrt{2}}\cosh \frac{y_1-t_1}{\sqrt{2}}\right] ,
\label{c2-35}
\end{equation}
The relations (\ref{c2-32})-(\ref{c2-35}) can be rewritten as
\begin{equation}
x+\frac{y-t}{\sqrt{2}}=e^{x_1}\cos\left(z_1+\frac{y_1+t_1}{\sqrt{2}}\right)e^{(y_1-t_1)/\sqrt{2}},
\label{c2-36}
\end{equation}
\begin{equation}
z+\frac{y+t}{\sqrt{2}}=e^{x_1}\sin\left(z_1+\frac{y_1+t_1}{\sqrt{2}}\right)e^{(y_1-t_1)/\sqrt{2}},
\label{c2-37}
\end{equation}
\begin{equation}
x-\frac{y-t}{\sqrt{2}}=e^{x_1}\cos\left(z_1-\frac{y_1+t_1}{\sqrt{2}}\right)e^{-(y_1-t_1)/\sqrt{2}},
\label{c2-38}
\end{equation}
\begin{equation}
z-\frac{y+t}{\sqrt{2}}=e^{x_1}\sin\left(z_1-\frac{y_1+t_1}{\sqrt{2}}\right)e^{-(y_1-t_1)/\sqrt{2}}.
\label{c2-39}
\end{equation}
By multiplying the sum of the squares of the first two and of the last
two relations (\ref{c2-36})-(\ref{c2-39}) it results that
\begin{equation}
e^{4x_1}=\rho_+^2\rho_-^2 ,
\label{c2-40}
\end{equation}
or
\begin{equation}
e^{x_1}=\rho .
\label{c2-41}
\end{equation}
By summing the squares of all relations (\ref{c2-36})-(\ref{c2-39}) it results that
\begin{equation}
d^2=\rho^2\cosh\left[\sqrt{2}(y_1-t_1)\right].
\label{c2-42}
\end{equation}
Then the quantities $y_1, z_1, t_1$ can be expressed in terms of the angles
$\phi, \chi, \psi$ defined in Eqs. (\ref{c2-12a})-(\ref{c2-12d}) as
\begin{equation}
z_1+\frac{y_1+t_1}{\sqrt{2}}=\phi ,
\label{c2-43}
\end{equation}
\begin{equation}
-z_1+\frac{y_1+t_1}{\sqrt{2}}=\chi ,
\label{c2-44}
\end{equation}
\begin{equation}
\frac{e^{(y_1-t_1)/\sqrt{2}}}{\sqrt{2}\left[\cosh\sqrt{2}(y_1-t_1)\right]^{1/2}}=\sin\psi,\:\:
\frac{e^{-(y_1-t_1)/\sqrt{2}}}{\sqrt{2}\left[\cosh\sqrt{2}(y_1-t_1)\right]^{1/2}}=\cos\psi.
\label{c2-45}
\end{equation}
From Eq. (\ref{c2-45}) it results that
\begin{equation}
y_1-t_1=\frac{1}{\sqrt{2}}\ln\tan\psi,
\label{c2-45b}
\end{equation}
so that 
\begin{equation}
y_1=\frac{\phi+\chi}{2\sqrt{2}}+\frac{1}{2\sqrt{2}}\ln\tan\psi,\;
z_1=\frac{\phi-\chi}{2},\;
t_1=\frac{\phi+\chi}{2\sqrt{2}}-\frac{1}{2\sqrt{2}}\ln\tan\psi.
\label{c2-45c}
\end{equation}
Substituting the expressions of the quantities $x_1, y_1, z_1, t_1$ in Eq.
(\ref{c2-31}) yields 
\begin{equation}
u=\rho\exp\left[\frac{1}{2\sqrt{2}}(\alpha-\gamma)\ln\tan\psi
+\frac{1}{2}\left(\beta+\frac{\alpha+\gamma}{\sqrt{2}}\right)\phi
-\frac{1}{2}\left(\beta-\frac{\alpha+\gamma}{\sqrt{2}}\right)\chi
\right],
\label{c2-46}
\end{equation}\index{exponential form!planar fourcomplex}
which will be called the exponential form of the planar fourcomplex number $u$.
It can be checked that
\begin{equation}
\exp\left[\frac{1}{2}\left(\beta+\frac{\alpha+\gamma}{\sqrt{2}}\right)\phi\right]
=\frac{1}{2}-\frac{\alpha-\gamma}{2\sqrt{2}}
+\left(\frac{1}{2}+\frac{\alpha-\gamma}{2\sqrt{2}}\right)\cos\phi+
\left(\frac{\beta}{2}+\frac{\alpha+\gamma}{2\sqrt{2}}\right)\sin\phi ,
\label{c2-51b}
\end{equation}
\begin{equation}
\exp\left[-\frac{1}{2}\left(\beta-\frac{\alpha+\gamma}{\sqrt{2}}\right)\chi\right]
=\frac{1}{2}+\frac{\alpha-\gamma}{2\sqrt{2}}
+\left(\frac{1}{2}-\frac{\alpha-\gamma}{2\sqrt{2}}\right)\cos\chi-
\left(\frac{\beta}{2}-\frac{\alpha+\gamma}{2\sqrt{2}}\right)\sin\chi ,
\label{c2-51c}
\end{equation}
which shows that $e^{[\beta+(\alpha+\gamma)/\sqrt{2}]\phi/2}$ and 
$e^{-[\beta-(\alpha+\gamma)/\sqrt{2}]\chi/2}$ are periodic
functions of $\phi$ and respectively $\chi$, with period $2\pi$.

The exponential of the logarithmic term in Eq. (\ref{c2-46}) can be expanded with
the aid of the relation (\ref{c2-30l}) as
\begin{equation}
\exp\left[\frac{1}{2\sqrt{2}}(\alpha-\gamma)\ln\tan\psi\right]=
\frac{1}{(\sin2\psi)^{1/2}}\left[\cos\left(\psi-\frac{\pi}{4}\right)
+\frac{\alpha-\gamma}{\sqrt{2}}\sin\left(\psi-\frac{\pi}{4}\right)\right].
\label{c2-47}
\end{equation}
Since according to Eq. (\ref{c2-13c}) $\tan\psi=\rho_+/\rho_-$, then 
\begin{equation}
\sin\psi\cos\psi=\frac{\rho_+\rho_-}{\rho_+^2+\rho_-^2} ,
\label{c2-48}
\end{equation}
and it can be checked that
\begin{equation}
\rho_+^2+\rho_-^2=2d^2,
\label{c2-49}
\end{equation}
where $d$ has been defined in Eq. (\ref{c2-10}). Thus
\begin{equation}
\rho^2=d^2\sin2\psi,
\label{c2-50}
\end{equation}
so that the planar fourcomplex number $u$ can be written as
\begin{equation}
u=d\left[\cos\left(\psi-\frac{\pi}{4}\right)
+\frac{\alpha-\gamma}{\sqrt{2}}\sin\left(\psi-\frac{\pi}{4}\right)\right]
\exp\left[\frac{1}{2}\left(\beta+\frac{\alpha+\gamma}{\sqrt{2}}\right)\phi
-\frac{1}{2}\left(\beta-\frac{\alpha+\gamma}{\sqrt{2}}\right)\chi\right],
\label{c2-51}
\end{equation}\index{trigonometric form!planar fourcomplex}
which will be called the trigonometric form of the planar fourcomplex number
$u$. 

If $u_1, u_2$ are planar fourcomplex numbers of moduli and angles $d_1, \phi_1,
\chi_1, \psi_1$ and respectively $d_2, \phi_2, \chi_2, \psi_2$, the product of
the factors depending on the planar angles can be calculated to be
\begin{eqnarray}
\lefteqn{[\cos(\psi_1-\pi/4)+\frac{\alpha-\gamma}{\sqrt{2}}\sin(\psi_1-\pi/4)] 
[\cos(\psi_2-\pi/4)+\frac{\alpha-\gamma}{\sqrt{2}}\sin(\psi_2-\pi/4)]}\nonumber\\
&&=[\cos(\psi_1-\psi_2)-\frac{\alpha-\gamma}{\sqrt{2}}\cos(\psi_1+\psi_2)] .
\label{c2-52-53}
\end{eqnarray}
The right-hand side of Eq. (\ref{c2-52-53}) can be written as
\begin{eqnarray}
\lefteqn{\cos(\psi_1-\psi_2)-\frac{\alpha-\gamma}{\sqrt{2}}\cos(\psi_1+\psi_2)
=[2(\cos^2\psi_1\cos^2\psi_2+\sin^2\psi_1\sin^2\psi_2)]^{1/2}
[\cos(\psi-\pi/4)\nonumber}\\
&&+\frac{\alpha-\gamma}{\sqrt{2}}\sin(\psi-\pi/4)] ,
\label{c2-54}
\end{eqnarray}
where the angle $\psi$, determined by the condition that
\begin{equation}
\tan(\psi-\pi/4)=-\cos(\psi_1+\psi_2)/\cos(\psi_1-\psi_2)
\label{c2-55}
\end{equation}
is given by $\tan\psi=\tan\psi_1\tan\psi_2$ ,
which is consistent with Eq. (\ref{c2-14b}).
The modulus $d$ of the product $u_1u_2$ is then
\begin{equation}
d=\sqrt{2}d_1d_2\left(\cos^2\psi_1\cos^2\psi_2+\sin^2\psi_1\sin^2\psi_2\right)^{1/2} .
\label{c2-56}
\end{equation}\index{transformation of variables!planar fourcomplex}

\subsection{Elementary functions of planar fourcomplex variables}

The logarithm $u_1$ of the planar fourcomplex number $u$, $u_1=\ln u$, can be defined
as the solution of the equation
\begin{equation}
u=e^{u_1} ,
\label{c2-57}
\end{equation}
written explicitly previously in Eq. (\ref{c2-31}), for $u_1$ as a function of
$u$. From Eq. (\ref{c2-46}) it results that 
\begin{equation}
\ln u=\ln \rho+\frac{1}{2\sqrt{2}}(\alpha-\gamma)\ln\tan\psi+
\frac{1}{2}\left(\beta+\frac{\alpha+\gamma}{\sqrt{2}}\right)\phi
-\frac{1}{2}\left(\beta-\frac{\alpha+\gamma}{\sqrt{2}}\right)\chi ,
\label{c2-58}
\end{equation}\index{logarithm!planar fourcomplex}
which is multivalued because of the presence of the terms proportional to
$\phi$ and $\chi$.
It can be inferred from Eqs. (\ref{c2-14a}) and (\ref{c2-14b}) that
\begin{equation}
\ln(uu^\prime)=\ln u+\ln u^\prime ,
\label{c2-59}
\end{equation}
up to multiples of $\pi[\beta+(\alpha+\gamma)/\sqrt{2}]$ and 
$\pi[\beta-(\alpha+\gamma)/\sqrt{2}]$.

The power function $u^m$ can be defined for real values of $n$ as
\begin{equation}
u^m=e^{m\ln u} .
\label{c2-60}
\end{equation}\index{power function!planar fourcomplex}
The power function is multivalued unless $n$ is an integer. 
For integer $n$, it can be inferred from Eq. (\ref{c2-59}) that
\begin{equation}
(uu^\prime)^m=u^m\:u^{\prime m} .
\label{c2-61}
\end{equation}
If, for example, $m=2$, it can be checked with the aid of Eq. (\ref{c2-51}) that
Eq. (\ref{c2-60}) gives indeed $(x+\alpha y+\beta z+\gamma t)^2=
x^2-z^2-2yt+2\alpha(xy-zt)+\beta(y^2-t^2+2xz)+2\gamma(xt+yz)$.

The trigonometric functions of the hypercomplex variable
$u$ and the addition theorems for these functions have been written in Eqs.
~(1.57)-(1.60). 
The cosine and sine functions of the hypercomplex variables $\alpha y, 
\beta z$ and $ \gamma t$ can be expressed as
\begin{equation}
\cos\alpha y=f_{40}(y)-\beta f_{42}(y), \: \sin\alpha y=\alpha f_{41}(y)-\gamma f_{43}(y), 
\label{c2-67}
\end{equation}\index{trigonometric functions, expressions!planar fourcomplex}
\begin{equation}
\cos\beta z=\cosh z, \: \sin\beta z=\beta\sinh z, 
\label{c2-66}
\end{equation}
\begin{equation}
\cos\gamma t=f_{40}(t)+\beta f_{42}(t), \: \sin\gamma t=\gamma f_{41}(t)-\alpha f_{43}(t) .
\label{c2-68}
\end{equation}
The cosine and sine functions of a planar fourcomplex number $x+\alpha y+\beta
z+\gamma t$ can then be
expressed in terms of elementary functions with the aid of the addition
theorems Eqs. (1.59), (1.60) and of the expressions in  Eqs. 
(\ref{c2-67})-(\ref{c2-68}). 

The hyperbolic functions of the hypercomplex variable
$u$ and the addition theorems for these functions have been written in Eqs.
~(1.62)-(1.65). 
The hyperbolic cosine and sine functions of the hypercomplex variables $\alpha y, 
\beta z$ and $ \gamma t$ can be expressed as
\begin{equation}
\cosh\alpha y=f_{40}(y)+\beta f_{42}(y), \: \sinh\alpha y=\alpha f_{41}(y)+\gamma f_{43}(y), 
\label{c2-74}
\end{equation}\index{hyperbolic functions, expressions!planar fourcomplex}
\begin{equation}
\cosh\beta z=\cos z, \: \sinh\beta z=\beta\sin z, 
\label{c2-73}
\end{equation}
\begin{equation}
\cosh\gamma t=f_{40}(t)-\beta f_{42}(t), \: \sinh\gamma t=\gamma f_{41}(t)+\alpha f_{43}(t) .
\label{c2-75}
\end{equation}
The hyperbolic cosine and sine functions of a planar fourcomplex number $x+\alpha y+\beta
z+\gamma t$ can then be
expressed in terms of elementary functions with the aid of the addition
theorems Eqs. (1.64), (1.65) and of the expressions in  Eqs. 
(\ref{c2-74})-(\ref{c2-75}).

\subsection{Power series of planar fourcomplex variables}

A planar fourcomplex series is an infinite sum of the form
\begin{equation}
a_0+a_1+a_2+\cdots+a_l+\cdots , 
\label{c2-76}
\end{equation}\index{series!planar fourcomplex}
where the coefficients $a_l$ are planar fourcomplex numbers. The convergence of 
the series (\ref{c2-76}) can be defined in terms of the convergence of its 4 real
components. The convergence of a planar fourcomplex series can however be studied
using planar fourcomplex variables. The main criterion for absolute convergence 
remains the comparison theorem, but this requires a number of inequalities
which will be discussed further.

The modulus of a planar fourcomplex number $u=x+\alpha y+\beta z+\gamma t$ can be
defined as 
\begin{equation}
|u|=(x^2+y^2+z^2+t^2)^{1/2} ,
\label{c2-77}
\end{equation}\index{modulus, definition!planar fourcomplex}
so that, according to Eq. (\ref{c2-10}), $d=|u|$. Since $|x|\leq |u|, |y|\leq |u|,
|z|\leq |u|, |t|\leq |u|$, a property of 
absolute convergence established via a comparison theorem based on the modulus
of the series (\ref{c2-76}) will ensure the absolute convergence of each real
component of that series.

The modulus of the sum $u_1+u_2$ of the planar fourcomplex numbers $u_1, u_2$ fulfils
the inequality
\begin{equation}
||u_1|-|u_2||\leq |u_1+u_2|\leq |u_1|+|u_2| .
\label{c2-78}
\end{equation}\index{modulus, inequalities!planar fourcomplex}
For the product the relation is 
\begin{equation}
|u_1u_2|\leq \sqrt{2}|u_1||u_2| ,
\label{c2-79}
\end{equation}
as can be shown from Eq. (\ref{c2-56}). The relation (\ref{c2-79}) replaces the
relation of equality extant for regular complex numbers. 
The equality in Eq. (\ref{c2-79}) takes place for 
$\cos^2(\psi_1-\psi_2)=1, \:\cos^2(\psi_1+\psi_2)=1,$ which means that
$x_1+(y_1-t_1)/\sqrt{2}=0, \:z_1+(y_1+t_1)/\sqrt{2}=0,
\:x_2+(y_2-t_2)/\sqrt{2}=0, \:z_2+(y_2+t_2)/\sqrt{2}=0,$ or
$x_1-(y_1-t_1)/\sqrt{2}=0, \:z_1-(y_1+t_1)/\sqrt{2}=0,
\:x_2-(y_2-t_2)/\sqrt{2}=0, \:z_2-(y_2+t_2)/\sqrt{2}=0$.
The modulus of the product, which has the property that $0\leq |u_1u_2|$,
becomes equal to zero for
$\cos^2(\psi_1-\psi_2)=0, \:\cos^2(\psi_1+\psi_2)=0,$ which means that
$x_1+(y_1-t_1)/\sqrt{2}=0, \:z_1+(y_1+t_1)/\sqrt{2}=0,
\:x_2-(y_2-t_2)/\sqrt{2}=0, \:z_2-(y_2+t_2)/\sqrt{2}=0$,
or
$x_1-(y_1-t_1)/\sqrt{2}=0, \:z_1-(y_1+t_1)/\sqrt{2}=0,
\:x_2+(y_2-t_2)/\sqrt{2}=0, \:z_2+(y_2+t_2)/\sqrt{2}=0$.
as discussed after Eq. (\ref{c2-9}).

It can be shown that
\begin{equation}
x^2+y^2+z^2+t^2\leq|u^2|\leq \sqrt{2}(x^2+y^2+z^2+t^2) .
\label{c2-80}
\end{equation}
The left relation in Eq. (\ref{c2-80}) becomes an equality for $\sin^2 2\psi=1$,
when $\rho_+=\rho_-$, which means that $x(y-t)+z(y+t)=0$.
The right relation in Eq. (\ref{c2-80}) becomes an equality for $\sin^2 2\psi=0$,
when $x+(y-t)/\sqrt{2}=0, \:z+(y+t)/\sqrt{2}=0,$
or $x-(y-t)/\sqrt{2}=0, \:z-(y+t)/\sqrt{2}=0.$
The inequality in Eq. (\ref{c2-79}) implies that
\begin{equation}
|u^l|\leq 2^{(l-1)/2}|u|^l .
\label{c2-81}
\end{equation}
From Eqs. (\ref{c2-79}) and (\ref{c2-81}) it results that
\begin{equation}
|au^l|\leq 2^{l/2} |a| |u|^l .
\label{c2-82}
\end{equation}

A power series of the planar fourcomplex variable $u$ is a series of the form
\begin{equation}
a_0+a_1 u + a_2 u^2+\cdots +a_l u^l+\cdots .
\label{c2-83}
\end{equation}\index{power series!planar fourcomplex}
Since
\begin{equation}
\left|\sum_{l=0}^\infty a_l u^l\right| \leq  \sum_{l=0}^\infty
2^{l/2}|a_l| |u|^l ,
\label{c2-84}
\end{equation}
a sufficient condition for the absolute convergence of this series is that
\begin{equation}
\lim_{l\rightarrow \infty}\frac{\sqrt{2}|a_{l+1}||u|}{|a_l|}<1 .
\label{c2-85}
\end{equation}
Thus the series is absolutely convergent for 
\begin{equation}
|u|<c,
\label{c2-86}
\end{equation}\index{convergence of power series!planar fourcomplex}
where 
\begin{equation}
c=\lim_{l\rightarrow\infty} \frac{|a_l|}{\sqrt{2}|a_{l+1}|} .
\label{c2-87}
\end{equation}

The convergence of the series (\ref{c2-83}) can be also studied with the aid of
the transformation 
\begin{equation}
x+\alpha y+\beta z+\gamma t=\sqrt{2}(e_1\xi+\tilde e_1 \upsilon+e_2\tau
+\tilde e_2\zeta) , 
\label{c2-87b}
\end{equation}\index{canonical form!planar fourcomplex}
where $\xi,\upsilon, \tau, \zeta$ have been defined in Eq. (\ref{c2-11}),
and
\begin{equation}
e_1=\frac{1}{2}+\frac{\alpha-\gamma}{2\sqrt{2}},\:\:
\tilde e_1=\frac{\beta}{2}+\frac{\alpha+\gamma}{2\sqrt{2}},\:\:
e_2=\frac{1}{2}-\frac{\alpha-\gamma}{2\sqrt{2}},\:\:
\tilde e_2=-\frac{\beta}{2}+\frac{\alpha+\gamma}{2\sqrt{2}}.
\label{c2-87c}
\end{equation}\index{canonical base!planar fourcomplex}
It can be checked that
\begin{eqnarray}
\lefteqn{e_1^2=e_1, \:\:\tilde e_1^2=-e_1,\:\: e_1\tilde e_1=\tilde e_1,\:\:
e_2^2=e_2, \:\:\tilde e_2^2=-e_2,\:\: e_2\tilde e_2=\tilde e_2,\:\:\nonumber}\\
&&e_1e_2=0,\:\: \tilde e_1\tilde e_2=0, \:\:e_1\tilde e_2=0, \:\:
e_2\tilde e_1=0.  
\label{c2-87d}
\end{eqnarray}
The moduli of the bases in Eq. (\ref{c2-87c}) are
\begin{equation}
|e_1|=\frac{1}{\sqrt{2}},\;|\tilde e_1|=\frac{1}{\sqrt{2}},\;
|e_2|=\frac{1}{\sqrt{2}},\;|\tilde e_2|=\frac{1}{\sqrt{2}},
\label{c2-87e}
\end{equation}
and it can be checked that
\begin{equation}
|x+\alpha y+\beta z+\gamma t|^2=\xi^2+\upsilon^2+\tau^2+\zeta^2.
\label{c2-87f}
\end{equation}\index{modulus, canonical variables!planar fourcomplex}
The ensemble $e_1, \tilde e_1, e_2, \tilde e_2$ will be called the canonical
planar fourcomplex base, and Eq. (\ref{c2-87b}) gives the canonical form of the
planar fourcomplex number.

If $u=u^\prime u^{\prime\prime}$, the components $\xi,\upsilon, \tau, \zeta$
are related, according to Eqs. (\ref{c2-17})-(\ref{c2-20}) by 
\begin{equation}
\xi=\sqrt{2}(\xi^\prime \xi^{\prime\prime}-\upsilon^\prime \upsilon^{\prime\prime}), \:\:
\upsilon=\sqrt{2}(\xi^\prime \upsilon^{\prime\prime}+\upsilon^\prime \xi^{\prime\prime}), \:\:
\tau=\sqrt{2}(\tau^\prime \tau^{\prime\prime}-\zeta^\prime \zeta^{\prime\prime}), \:\:
\zeta=\sqrt{2}(\tau^\prime \zeta^{\prime\prime}+\zeta^\prime \xi^{\prime\prime}), \:\:
\label{c2-87g}
\end{equation}\index{transformation of variables!planar fourcomplex}
which show that, upon multiplication, the components $\xi,\upsilon$ and $\tau,
\zeta$ obey, up to a normalization constant, the same
rules as the real and imaginary components of usual, two-dimensional complex
numbers.

If the coefficients in Eq. (\ref{c2-83}) are 
\begin{equation}
a_l= a_{l0}+\alpha a_{l1}+\beta a_{l2}+\gamma a_{l3}, 
\label{c2-n88a}
\end{equation}
and
\begin{equation}
A_{l1}=a_{l0}+\frac{a_{l1}-a_{l3}}{\sqrt{2}},\;
\tilde A_{l1}=a_{l2}+\frac{a_{l1}+a_{l3}}{\sqrt{2}},\;
A_{l2}=a_{l0}-\frac{a_{l1}-a_{l3}}{\sqrt{2}},\;
\tilde A_{l2}=-a_{l2}+\frac{a_{l1}+a_{l3}}{\sqrt{2}},
\label{c2-n88b}
\end{equation}
the series (\ref{c2-83}) can be written as
\begin{equation}
\sum_{l=0}^\infty 2^{l/2}\left[
(e_1 A_{l1}+\tilde e_1\tilde A_{l1})(e_1 \xi+\tilde e_1\upsilon)^l 
+(e_2 A_{l2}+\tilde e_2\tilde A_{l2})(e_2 \tau+\tilde e_2\zeta)^l 
\right].
\label{c2-n89a}
\end{equation}
Thus, the series in Eqs. (\ref{c2-83}) and (\ref{c2-n89a}) are
absolutely convergent for   
\begin{equation}
\rho_+<c_1, \;\rho_-<c_2,
\label{c2-n90}
\end{equation}\index{convergence, region of!planar fourcomplex}
where 
\begin{equation}
c_1=\lim_{l\rightarrow\infty} \frac
{\left[A_{l1}^2+\tilde A_{l1}^2\right]^{1/2}}
{\sqrt{2}\left[A_{l+1,1}^2+\tilde A_{l+1,1}^2\right]^{1/2}},\;\;
c_2=\lim_{l\rightarrow\infty} \frac
{\left[A_{l2}^2+\tilde A_{l2}^2\right]^{1/2}}
{\sqrt{2}\left[A_{l+1,2}^2+\tilde A_{l+1,2}^2\right]^{1/2}}.
\label{c2-n91}
\end{equation}

It can be shown that $c=(1/\sqrt{2}){\rm
min}(c_1,c_2)$, where ${\rm min}$ designates the smallest of
the numbers $c_1,c_2$. Using the expression of $|u|$ in
Eq. (\ref{c2-87f}),  it can be seen that the spherical region of
convergence defined in Eqs. (\ref{c2-86}), (\ref{c2-87}) is included in the
cylindrical region of convergence defined in Eqs. (\ref{c2-n90}) and (\ref{c2-n91}).

\subsection{Analytic functions of planar fourcomplex variables}

The fourcomplex function $f(u)$ of the fourcomplex variable $u$ has
been expressed in Eq. (\ref{g16}) in terms of 
the real functions $P(x,y,z,t),Q(x,y,z,t),R(x,y,z,t), S(x,y,z,t)$ of real
variables $x,y,z,t$.
The
relations between the partial derivatives of the functions $P, Q, R, S$ are
obtained by setting succesively in   
Eq. (\ref{g17}) $\Delta x\rightarrow 0, \Delta y=\Delta z=\Delta t=0$;
then $ \Delta y\rightarrow 0, \Delta x=\Delta z=\Delta t=0;$  
then $  \Delta z\rightarrow 0,\Delta x=\Delta y=\Delta t=0$; and finally
$ \Delta t\rightarrow 0,\Delta x=\Delta y=\Delta z=0 $. 
The relations are 
\begin{equation}
\frac{\partial P}{\partial x} = \frac{\partial Q}{\partial y} =
\frac{\partial R}{\partial z} = \frac{\partial S}{\partial t},
\label{c2-95}
\end{equation}\index{relations between partial derivatives!planar fourcomplex}
\begin{equation}
\frac{\partial Q}{\partial x} = \frac{\partial R}{\partial y} =
\frac{\partial S}{\partial z} = -\frac{\partial P}{\partial t},
\label{c2-97}
\end{equation}
\begin{equation}
\frac{\partial R}{\partial x} = \frac{\partial S}{\partial y} =
-\frac{\partial P}{\partial z} = -\frac{\partial Q}{\partial t},
\label{c2-96}
\end{equation}
\begin{equation}
\frac{\partial S}{\partial x} =-\frac{\partial P}{\partial y} =
-\frac{\partial Q}{\partial z} =-\frac{\partial R}{\partial t}.
\label{c2-98}
\end{equation}
The relations (\ref{c2-95})-(\ref{c2-98}) are analogous to the Riemann relations
for the real and imaginary components of a complex function. It can be shown
from Eqs. (\ref{c2-95})-(\ref{c2-98}) that the component $P$ is a solution
of the equations 
\begin{equation}
\frac{\partial^2 P}{\partial x^2}+\frac{\partial^2 P}{\partial z^2}=0,
\:\: 
\frac{\partial^2 P}{\partial y^2}+\frac{\partial^2 P}{\partial t^2}=0,
\:\:
\label{c2-99}
\end{equation}\index{relations between second-order derivatives!planar fourcomplex}
and the components $Q, R, S$ are solutions of similar equations.
As can be seen from Eqs. (\ref{c2-99}), the components $P, Q, R, S$ of
an analytic function of planar fourcomplex variable are harmonic 
with respect to the pairs of variables $x,y$ and $ z,t$.
The component $P$ is also a solution of the mixed-derivative
equations
\begin{equation}
\frac{\partial^2 P}{\partial x^2}=-\frac{\partial^2 P}{\partial y\partial t},
\:\: 
\frac{\partial^2 P}{\partial y^2}=\frac{\partial^2 P}{\partial x\partial z},
\:\:
\frac{\partial^2 P}{\partial z^2}=\frac{\partial^2 P}{\partial y\partial t},
\:\:
\frac{\partial^2 P}{\partial t^2}=-\frac{\partial^2 P}{\partial x\partial z},
\:\:
\label{c2-106b}
\end{equation}
and the components $Q, R, S$ are solutions of similar equations.
The component $P$ is also a solution of the mixed-derivative
equations 
\begin{equation}
\frac{\partial^2 P}{\partial x\partial y}=-\frac{\partial^2 P}{\partial
z\partial t} ,
\:\: 
\frac{\partial^2 P}{\partial x\partial t}=\frac{\partial^2 P}{\partial
y\partial z} ,
\label{c2-107}
\end{equation}
and the components $Q, R, S$ are solutions of similar equations.

\subsection{Integrals of functions of planar fourcomplex variables}

The singularities of planar fourcomplex functions arise from terms of the form
$1/(u-u_0)^m$, with $m>0$. Functions containing such terms are singular not
only at $u=u_0$, but also at all points of the two-dimensional hyperplanes
passing through $u_0$ and which are parallel to the nodal hyperplanes. 

The integral of a planar fourcomplex function between two points $A, B$ along a path
situated in a region free of singularities is independent of path, which means
that 
\begin{equation}
\oint_\Gamma f(u) du = 0,
\label{c2-111}
\end{equation}
where it is supposed that a surface $\Sigma$ spanning 
the closed loop $\Gamma$ is not intersected by any of
the two-dimensional hyperplanes associated with the
singularities of the function $f(u)$. Using the expression, Eq. (\ref{g16})
for $f(u)$ and the fact that $du=dx+\alpha  dy+\beta dz+\gamma dt$, the
explicit form of the integral in Eq. (\ref{c2-111}) is
\begin{eqnarray}
\lefteqn{\oint _\Gamma f(u) du = \oint_\Gamma
[(Pdx-Sdy-Rdz-Qdt)+\alpha(Qdx+Pdy-Sdz-Rdt)\nonumber}\\
&&+\beta(Rdx+Qdy+Pdz-Sdt)+\gamma(Sdx+Rdy+Qdz+Pdt)] .
\label{c2-112}
\end{eqnarray}\index{integrals, path!planar fourcomplex}
If the functions $P, Q, R, S$ are regular on a surface $\Sigma$
spanning the loop $\Gamma$,
the integral along the loop $\Gamma$ can be transformed with the aid of the
theorem of Stokes in an integral over the surface $\Sigma$ of terms of the form
$\partial P/\partial y + \partial S/\partial x,\:\:
\partial P/\partial z +  \partial R/\partial x, \:\:
\partial P/\partial t + \partial Q/\partial x, \:\:
\partial R/\partial y -  \partial S/\partial z, \:\:
\partial S/\partial t - \partial Q/\partial y, \:\:
\partial R/\partial t - \partial Q/\partial z$ 
and of similar terms arising
from the $\alpha, \beta$ and $\gamma$ components, 
which are equal to zero by Eqs. (\ref{c2-95})-(\ref{c2-98}), and this proves Eq.
(\ref{c2-111}). 

The integral of the function $(u-u_0)^m$ on a closed loop $\Gamma$ is equal to
zero for $m$ a positive or negative integer not equal to -1,
\begin{equation}
\oint_\Gamma (u-u_0)^m du = 0, \:\: m \:\:{\rm integer},\: m\not=-1 .
\label{c2-112b}
\end{equation}
This is due to the fact that $\int (u-u_0)^m du=(u-u_0)^{m+1}/(m+1), $ and to
the fact that the function $(u-u_0)^{m+1}$ is singlevalued for $m$ an integer.

The integral $\oint du/(u-u_0)$ can be calculated using the exponential form 
(\ref{c2-46}),
\begin{eqnarray}
\lefteqn{u-u_0=\rho\exp\left[\frac{1}{2\sqrt{2}}(\alpha-\gamma)\ln\tan\psi
+\frac{1}{2}\left(\beta+\frac{\alpha+\gamma}{\sqrt{2}}\right)\phi
-\frac{1}{2}\left(\beta-\frac{\alpha+\gamma}{\sqrt{2}}\right)\chi\right],\nonumber}\\
&&
\label{c2-113}
\end{eqnarray}
so that 
\begin{equation}
\frac{du}{u-u_0}=\frac{d\rho}{\rho}+
\frac{1}{2\sqrt{2}}(\alpha-\gamma)d\ln\tan\psi
+\frac{1}{2}\left(\beta+\frac{\alpha+\gamma}{\sqrt{2}}\right)d\phi
-\frac{1}{2}\left(\beta-\frac{\alpha+\gamma}{\sqrt{2}}\right)d\chi .
\label{c2-114}
\end{equation}
Since $\rho$ and $\psi$ are singlevalued variables, it follows that
$\oint_\Gamma d\rho/\rho =0, \oint_\Gamma d\ln\tan\psi=0$. On the other hand,
$\phi$ and $\chi$ are cyclic variables, so that they may give a contribution to
the integral around the closed loop $\Gamma$.
Thus, if $C_+$ is a circle of radius $r$
parallel to the $\xi O\upsilon$ plane, whose
projection of the center of this circle on the $\xi O\upsilon$ plane
coincides with the projection of the point $u_0$ on this plane, the points
of the circle $C_+$ are described according to Eqs.
(\ref{c2-11})-(\ref{c2-12d}) by the equations
\begin{eqnarray}
\lefteqn{\xi=\xi_0+r \sin\psi\cos\phi , \:
\upsilon=\upsilon_0+r \sin\psi\sin\phi , \:
\tau=\tau_0+r\cos\psi \cos\chi , \nonumber}\\
&&\zeta=\zeta_0+r \cos\psi\sin\chi , 
\label{c2-115}
\end{eqnarray}
for constant values of $\chi$ and $\psi, \:\psi\not=0, \pi/2$, where
$u_0=x_0+\alpha y_0+\beta 
z_0+\gamma t_0$,  and $\xi_0, \upsilon_0, \tau_0, \zeta_0$ are calculated from
$x_0, y_0, z_0, t_0$ according to Eqs. (\ref{c2-11}).
Then
\begin{equation}
\oint_{C_+}\frac{du}{u-u_0}=\pi\left(\beta+\frac{\alpha+\gamma}{\sqrt{2}}\right). 
\label{c2-116}
\end{equation}
If $C_-$ is a circle of radius $r$
parallel to the $\tau O\zeta$ plane,
whose projection of the center of this circle on the $\tau O\zeta$ plane
coincides with the projection of the point $u_0$ on this plane, the points
of the circle $C_-$ are described by the same Eqs. (\ref{c2-115}) 
but for constant values of $\phi$ and $\psi, \:\psi\not=0, \pi/2$.
Then
\begin{equation}
\oint_{C_-}\frac{du}{u-u_0}=-\pi\left(\beta-\frac{\alpha+\gamma}{\sqrt{2}}\right) .
\label{c2-117}
\end{equation}
\newpage
\setlength{\oddsidemargin}{-1cm}      
The expression of $\oint_\Gamma du/(u-u_0)$ can be written as a single
equation with the aid of the functional int($M,C$) defined in Eq. (\ref{118}) as
\begin{equation}
\oint_\Gamma\frac{du}{u-u_0}=
\pi\left(\beta+\frac{\alpha+\gamma}{\sqrt{2}}\right) \;{\rm
int}(u_{0\xi\upsilon},\Gamma_{\xi\upsilon}) 
-\pi \left(\beta-\frac{\alpha+\gamma}{\sqrt{2}}\right)\;{\rm
int}(u_{0\tau\zeta},\Gamma_{\tau\zeta}), 
\label{c2-119}
\end{equation}\index{poles and residues!planar fourcomplex}
where $u_{0\xi\upsilon}, u_{0\tau\zeta}$ and $\Gamma_{\xi\upsilon},
\Gamma_{\tau\zeta}$ are respectively the projections of the point $u_0$ and of
the loop $\Gamma$ on the planes $\xi \upsilon$ and $\tau \zeta$.

If $f(u)$ is an analytic planar fourcomplex function which can be expanded in a
series as written in Eq. (1.89), and the expansion holds on the curve
$\Gamma$ and on a surface spanning $\Gamma$, then from Eqs. (\ref{c2-112b}) and
(\ref{c2-119}) it follows that
\begin{equation}
\oint_\Gamma \frac{f(u)du}{u-u_0}=
\pi\left[\left(\beta+\frac{\alpha+\gamma}{\sqrt{2}}\right) \;{\rm
int}(u_{0\xi\upsilon},\Gamma_{\xi\upsilon}) 
- \left(\beta-\frac{\alpha+\gamma}{\sqrt{2}}\right)\;{\rm
int}(u_{0\tau\zeta},\Gamma_{\tau\zeta})\right]\;f(u_0) , 
\label{c2-120}
\end{equation}
where $\Gamma_{\xi\upsilon}, \Gamma_{\tau\zeta}$ are the projections of 
the curve $\Gamma$ on the planes $\xi \upsilon$ and respectively $\tau \zeta$,
as shown in Fig. \ref{fig11}. As remarked previously,
the definition of the variables in this section is different from the
former definition for the circular hypercomplex numbers. 

Substituting in the right-hand side of 
Eq. (\ref{c2-120}) the expression of $f(u)$ in terms of the real 
components $P, Q, R, S$, Eq. (\ref{g16}), yields
\begin{eqnarray}
\lefteqn{\oint_\Gamma \frac{f(u)du}{u-u_0}\nonumber}\\
&&=\pi \left[
\left(\beta+\frac{\alpha+\gamma}{\sqrt{2}}\right)P
-\left(1+\frac{\alpha-\gamma}{\sqrt{2}}\right)R
-\left(\gamma-\frac{1-\beta}{\sqrt{2}}\right)Q
-\left(\alpha-\frac{1+\beta}{\sqrt{2}}\right)S
\right] 
{\rm int}\left(u_{0\xi\upsilon},\Gamma_{\xi\upsilon}\right)\nonumber\\
&&-\pi \left[
\left(\beta-\frac{\alpha+\gamma}{\sqrt{2}}\right)P
-\left(1-\frac{\alpha-\gamma}{\sqrt{2}}\right)R
-\left(\gamma-\frac{1+\beta}{\sqrt{2}}\right)Q
-\left(\alpha-\frac{1-\beta}{\sqrt{2}}\right)S
\right] 
{\rm int}(u_{0\tau\zeta},\Gamma_{\tau\zeta}) ,\nonumber\\
&&
\label{c2-121}
\end{eqnarray}
where $P, Q, R, S$ are the values of the components of $f$ at $u=u_0$.

If $f(u)$ can be expanded as written in Eq. (1.89) on 
$\Gamma$ and on a surface spanning $\Gamma$, then from Eqs. (\ref{c2-112b}) and
(\ref{c2-119}) it also results that
\begin{equation}
\oint_\Gamma \frac{f(u)du}{(u-u_0)^{m+1}}=
\frac{\pi}{m!}\left[\left(\beta+\frac{\alpha+\gamma}{\sqrt{2}}\right) \;{\rm
int}(u_{0\xi\upsilon},\Gamma_{\xi\upsilon}) 
- \left(\beta-\frac{\alpha+\gamma}{\sqrt{2}}\right)\;{\rm
int}(u_{0\tau\zeta},\Gamma_{\tau\zeta})\right]\; 
f^{(m)}(u_0) ,
\label{c2-122}
\end{equation}
where it has been used the fact that the derivative $f^{(m)}(u_0)$ of order $n$
of $f(u)$ at $u=u_0$ is related to the expansion coefficient in Eq. (1.89)
according to Eq. (1.93).

If a function $f(u)$ is expanded in positive and negative powers of $u-u_j$,
where $u_j$ are planar fourcomplex constants, $j$ being an index, the integral of $f$
on a closed loop $\Gamma$ is determined by the terms in the expansion of $f$
which are of the form $a_j/(u-u_j)$,
\begin{equation}
f(u)=\cdots+\sum_j\frac{a_j}{u-u_j}+\cdots
\label{c2-123}
\end{equation}
Then the integral of $f$ on a closed loop $\Gamma$ is
\begin{equation}
\oint_\Gamma f(u) du = 
\pi\left(\beta+\frac{\alpha+\gamma}{\sqrt{2}}\right) \sum_j{\rm
int}(u_{j\xi\upsilon},\Gamma_{\xi\upsilon})a_j 
- \pi\left(\beta-\frac{\alpha+\gamma}{\sqrt{2}}\right)\sum_j{\rm
int}(u_{j\tau\zeta},\Gamma_{\tau\zeta})a_j .
\label{c2-124}
\end{equation}      
\newpage
\setlength{\oddsidemargin}{0.9cm}

\subsection{Factorization of planar fourcomplex polynomials}

A polynomial of degree $m$ of the planar fourcomplex variable 
$u=x+\alpha y+\beta z+\gamma t$ has the form
\begin{equation}
P_m(u)=u^m+a_1 u^{m-1}+\cdots+a_{m-1} u +a_m ,
\label{c2-125}
\end{equation}
where the constants are in general planar fourcomplex numbers.

It can be shown that any planar fourcomplex polynomial has a planar fourcomplex root, whence
it follows that a polynomial of degree $m$ can be written as a product of
$m$ linear factors of the form $u-u_j$, where the planar fourcomplex numbers $u_j$ are
the roots of the polynomials, although the factorization may not be unique, 
\begin{equation}
P_m(u)=\prod_{j=1}^m (u-u_j) .
\label{c2-126}
\end{equation}\index{roots!planar fourcomplex}

The fact that any planar fourcomplex polynomial has a root can be shown by considering
the transformation of a fourdimensional sphere with the center at the origin by
the function $u^m$. The points of the hypersphere of radius $d$ are of the form
written in Eq. (\ref{c2-51}), with $d$ constant and $\phi, \chi, \psi$
arbitrary. The point $u^m$ is
\begin{eqnarray}
\lefteqn{u^m=d^m\left[\cos\left(\psi-\frac{\pi}{4}\right)
+\frac{\alpha-\gamma}{\sqrt{2}}\sin\left(\psi-\frac{\pi}{4}\right)\right]^m\nonumber}\\
&&\exp\left[\frac{1}{2}\left(\beta+\frac{\alpha+\gamma}{\sqrt{2}}\right)m\phi
-\frac{1}{2}\left(\beta-\frac{\alpha+\gamma}{\sqrt{2}}\right)m\chi\right].
\label{c2-127}
\end{eqnarray}
It can be shown with the aid of Eq. (\ref{c2-56}) that
\begin{equation}
\left|u
\exp\left[\frac{1}{2}\left(\beta+\frac{\alpha+\gamma}{\sqrt{2}}\right)\phi
-\frac{1}{2}\left(\beta-\frac{\alpha+\gamma}{\sqrt{2}}\right)\chi\right]\right|
=|u|,
\label{c2-127bb}
\end{equation}
so that
\begin{eqnarray}
\lefteqn{\left|
\left[\cos(\psi-\pi/4)+\frac{\alpha-\gamma}\sin(\psi-\pi/4)\right]^m
\exp\left[\frac{1}{2}\left(\beta+\frac{\alpha+\gamma}{\sqrt{2}}\right)m\phi
-\frac{1}{2}\left(\beta-\frac{\alpha+\gamma}{\sqrt{2}}\right)m\chi\right]
\right|\nonumber}\\
&&=\left|\left(\cos(\psi-\pi/4)+\frac{\alpha-\gamma}{\sqrt{2}}\sin(\psi-\pi/4)\right)^m\right| .
\label{c2-127b}
\end{eqnarray}
The right-hand side of Eq. (\ref{c2-127b}) is
\begin{equation}
\left|\left(\cos\epsilon+\frac{\alpha-\gamma}{\sqrt{2}}\sin\epsilon\right)^m\right|^2
=\sum_{k=0}^m C_{2m}^{2k}\cos^{2m-2k}\epsilon\sin^{2k}\epsilon ,
\label{c2-128}
\end{equation}
where $\epsilon=\psi-\pi/4$, 
and since $C_{2m}^{2k}\geq C_m^k$, it can be concluded that
\begin{equation}
\left|\left(\cos\epsilon+\frac{\alpha-\gamma}{\sqrt{2}}\sin\epsilon\right)^m\right|^2\geq 1 .
\label{c2-129}
\end{equation}
Then
\begin{equation}
d^m\leq |u^m|\leq 2^{(m-1)/2} d^m ,
\label{c2-129b}
\end{equation}
which shows that the image of a four-dimensional sphere via the transformation
operated by the function $u^m$ is a finite hypersurface.

If $u^\prime=u^m$, and
\begin{eqnarray}
\lefteqn{u^\prime=d^\prime
\left[\cos(\psi^\prime-\pi/4)
+\frac{\alpha-\gamma}{\sqrt{2}}\sin(\psi^\prime-\pi/4)\right]\nonumber}\\
&&\exp\left[
\frac{1}{2}\left(\beta+\frac{\alpha+\gamma}{\sqrt{2}}\right)\phi^\prime
-\frac{1}{2}\left(\beta-\frac{\alpha+\gamma}{\sqrt{2}}\right)\chi^\prime
\right],
\label{c2-130}
\end{eqnarray}
then 
\begin{equation}
\phi^\prime=m\phi, \: \chi^\prime=m\chi, \: \tan\psi^\prime=\tan^m\psi .
\label{c2-131}
\end{equation}
Since for any values of the angles $\phi^\prime, \chi^\prime, \psi^\prime$
there is a set of solutions $\phi, \chi, \psi$ of Eqs. (\ref{c2-131}), and since
the image of the hypersphere is a finite hypersurface, it follows that the
image of the four-dimensional sphere via the function $u^m$ is also a closed
hypersurface. A continuous hypersurface is called closed when any ray issued from the
origin intersects that surface at least once in the finite part of the space.

A transformation of the four-dimensional space by the polynomial $P_m(u)$
will be considered further. By this transformation, a hypersphere of radius $d$
having the center at the origin is changed into a certain finite closed
surface, as discussed previously. 
The transformation of the four-dimensional space by the polynomial $P_m(u)$
associates to the point $u=0$ the point $f(0)=a_m$, and the image of a
hypersphere of very large radius $d$ can be represented with good approximation
by the image of that hypersphere by the function $u^m$. 
The origin of the axes is an inner
point of the latter image. If the radius of the hypersphere is now reduced
continuously from the initial very large values to zero, the image hypersphere
encloses initially the origin, but the image shrinks to $a_m$ when the radius
approaches the value zero.  Thus, the
origin is initially inside the image hypersurface, and it lies outside the
image hypersurface when the radius of the hypersphere tends to zero. Then since
the image hypersurface is closed, the image surface must intersect at some
stage the origin of the axes, which means that there is a point $u_1$ such that
$f(u_1)=0$. The factorization in Eq. (\ref{c2-126}) can then be obtained by
iterations.


The roots of the polynomial $P_m$ can be obtained by the following method.
If the constants in Eq. (\ref{c2-125}) are $a_l=a_{l0}+\alpha a_{l1}
+\beta a_{l2}+\gamma a_{l3}$, and with the 
notations of Eq. (\ref{c2-n88b}), the polynomial $P_m(u)$ can be written as
\begin{eqnarray}
\lefteqn{P_m=\sum_{l=0}^{m} 2^{(m-l)/2}
(e_1 A_{l1}+\tilde e_1\tilde A_{l1})(e_1 \xi+\tilde e_1 \upsilon)^{m-l}\nonumber}\\
&&+\sum_{l=0}^{m} 2^{(m-l)/2}
(e_2 A_{l2}+\tilde e_2\tilde A_{l2})(e_2 \tau+\tilde e_2 \zeta)^{m-l} ,
\label{c2-126a}
\end{eqnarray}\index{polynomial, canonical variables!planar fourcomplex}
where the constants $A_{lk}, \tilde A_{lk}, k=1,2$ are real numbers.
Each of the polynomials of degree $m$ in $e_1 \xi+\tilde e_1\upsilon, 
e_2 \tau+\tilde e_2\zeta$
in Eq. (\ref{c2-126a}) 
can always be written as a product of linear factors of the form
$e_1 (\xi-\xi_p)+\tilde e_1(\upsilon- \upsilon_p)$ and respectively
$e_2 (\tau-\tau_p)+\tilde e_2(\zeta- \zeta_p)$, where the
constants $\xi_p, \upsilon_p, \tau_p, \zeta_p$ are real,
\begin{eqnarray}
\sum_{l=0}^{m} 2^{(m-l)/2}
(e_1 A_{l1}+\tilde e_1\tilde A_{l1})(e_1 \xi+\tilde e_1 \upsilon)^{m-l}
=\prod_{p=1}^{m}2^{m/2}\left\{e_1 (\xi-\xi_p)+\tilde e_1(\upsilon- \upsilon_p)
\right\},
\label{c2-126bb}
\end{eqnarray}
\begin{eqnarray}
\sum_{l=0}^{m} 2^{(m-l)/2}
(e_2 A_{l2}+\tilde e_2\tilde A_{l2})(e_2 \tau+\tilde e_2 \zeta)^{m-l}
=\prod_{p=1}^{m}2^{m/2}\left\{e_2 (\tau-\tau_p)+\tilde e_2(\zeta- \zeta_p)
\right\}.
\label{c2-126bc}
\end{eqnarray}

Due to the relations  (\ref{c2-87d}),
the polynomial $P_m(u)$ can be written as a product of factors of
the form 
\begin{eqnarray}
P_m(u)=\prod_{p=1}^m 2^{m/2}\left\{e_1 (\xi-\xi_p)+\tilde e_1(\upsilon- \upsilon_p)
+e_2 (\tau-\tau_p)+\tilde e_2(\zeta- \zeta_p)\right\}.
\label{c2-128b}
\end{eqnarray}\index{polynomial, factorization!planar fourcomplex}
This relation can be written with the aid of Eq. (\ref{c2-87b}) in the form 
(\ref{c2-126}), where
\begin{eqnarray}
u_p=\sqrt{2}(e_1 \xi_p+\tilde e_1 \upsilon_p
+e_2 \tau_p+\tilde e_2 \zeta_p) .
\label{c2-129bx}
\end{eqnarray}
The roots  
$e_1 \xi_p+\tilde e_1 \upsilon_p$ and $e_2 \tau_p+\tilde e_2 \zeta_p$
defined in Eqs. (\ref{c2-126bb}) and respectively (\ref{c2-126bc}) may be ordered
arbitrarily. This means that Eq. 
(\ref{c2-129bx}) gives sets of $m$ roots 
$u_1,...,u_m$ of the polynomial $P_m(u)$, 
corresponding to the various ways in which the roots 
$e_1 \xi_p+\tilde e_1 \upsilon_p$ and $e_2 \tau_p+\tilde e_2 \zeta_p$
are ordered according to $p$ for each polynomial. 
Thus, while the hypercomplex components in Eqs. (\ref{c2-126bb}),
(\ref{c2-126bc}) taken separately have unique factorizations, the polynomial
$P_m(u)$ can be written in many different ways as a product of linear factors. 
The result of the planar fourcomplex integration, Eq. (\ref{c2-124}), is however unique.

If, for example, $P(u)=u^2+1$, the possible factorizations are
$P=(u-\tilde e_1-\tilde e_2)(u+\tilde e_1+\tilde e_2)$ and
$P=(u-\tilde e_1+\tilde e_2)(u+\tilde e_1-\tilde e_2)$ which can also be
written 
as $u^2+1=(u-\beta)(u+\beta)$ or as
$u^2+1=\left\{u-(\alpha+\gamma)/\sqrt{2}\right\}
\left\{(u+(\alpha+\gamma)/\sqrt{2}\right\}$. The result of the planar fourcomplex
integration, Eq. (\ref{c2-124}), is however unique. 
It can be checked
that $(\pm \tilde e_1\pm\tilde e_2)^2=
-e_1-e_2=-1$.

\subsection{Representation of planar fourcomplex numbers 
by irreducible matrices}

If $T$ is the unitary matrix,
\begin{equation}
T =\left(
\begin{array}{cccc}
\frac{1}{\sqrt{2}}&\frac{1}{2}       &0                  &-\frac{1}{2}  \\
0                 &\frac{1}{2}       &\frac{1}{\sqrt{2}} &\frac{1}{2} \\
\frac{1}{\sqrt{2}}&-\frac{1}{2}      &0                  &\frac{1}{2}  \\
0                 &\frac{1}{2}       &-\frac{1}{\sqrt{2}}&\frac{1}{2} \\
\end{array}
\right),
\label{c2-129x}
\end{equation}
it can be shown 
that the matrix $T U T^{-1}$ has the form 
\begin{equation}
T U T^{-1}=\left(
\begin{array}{cc}
V_1      &     0    \\
0        &     V_2  \\
\end{array}
\right),
\label{c2-129y}
\end{equation}\index{representation by irreducible matrices!planar fourcomplex}
where $U$ is the matrix in Eq. (\ref{c2-23}) used to represent the planar fourcomplex
number $u$. In Eq. (\ref{c2-129y}), $V_1, V_2$ are
the matrices
\begin{equation}
V_1=\left(
\begin{array}{cc}
x+\frac{y-t}{\sqrt{2}}    &  z+\frac{y+t}{\sqrt{2}}   \\
-z-\frac{y+t}{\sqrt{2}}     &  x+\frac{y-t}{\sqrt{2}}   \\
\end{array}\right),\;\;
V_2=\left(
\begin{array}{cc}
x-\frac{y-t}{\sqrt{2}}    &-z+\frac{y+t}{\sqrt{2}}   \\
z-\frac{y+t}{\sqrt{2}}  &   x-\frac{y-t}{\sqrt{2}}   \\
\end{array}\right).
\label{c2-130x}
\end{equation}
In Eq. (\ref{c2-129y}), the symbols 0 denote the matrix
\begin{equation}
\left(
\begin{array}{cc}
0   &  0   \\
0   &  0   \\
\end{array}\right).
\label{c2-131x}
\end{equation}
The relations between the variables 
$x+(y-t)/\sqrt{2}, z+(y+t)/\sqrt{2}$,
$x-(y-t)/\sqrt{2},-z+(y+t)/\sqrt{2}$ for the multiplication
of planar fourcomplex numbers have been written in Eqs. (\ref{c2-17})-(\ref{c2-20}). The
matrix 
$T U T^{-1}$ provides an irreducible representation
\cite{4} of the planar fourcomplex number $u$ in terms of matrices with real
coefficients.

\section{Polar complex numbers in four dimensions}

\subsection{Operations with polar fourcomplex numbers}

A polar fourcomplex number is determined by its four components $(x,y,z,t)$. The sum
of the polar fourcomplex numbers $(x,y,z,t)$ and
$(x^\prime,y^\prime,z^\prime,t^\prime)$ is the polar fourcomplex
number $(x+x^\prime,y+y^\prime,z+z^\prime,t+t^\prime)$. \index{sum!polar fourcomplex}
The product of the polar fourcomplex numbers
$(x,y,z,t)$ and $(x^\prime,y^\prime,z^\prime,t^\prime)$ 
is defined in this section to be the polar fourcomplex
number
$(xx^\prime+yt^\prime+zz^\prime+ty^\prime,
xy^\prime+yx^\prime+zt^\prime+tz^\prime,
xz^\prime+yy^\prime+zx^\prime+tt^\prime,
xt^\prime+yz^\prime+zy^\prime+tx^\prime)$.\index{product!polar fourcomplex}
Polar fourcomplex numbers and their operations can be represented by  writing the
polar fourcomplex number $(x,y,z,t)$ as  
$u=x+\alpha y+\beta z+\gamma t$, where $\alpha, \beta$ and $\gamma$ 
are bases for which the multiplication rules are 
\begin{equation}
\alpha^2=\beta, \:\beta^2=1, \:\gamma^2=\beta,
\alpha\beta=\beta\alpha=\gamma,\: 
\alpha\gamma=\gamma\alpha=-1, \:\beta\gamma=\gamma\beta=\alpha .
\label{ch1}
\end{equation}\index{complex units!polar fourcomplex}
Two polar fourcomplex numbers $u=x+\alpha y+\beta z+\gamma t, 
u^\prime=x^\prime+\alpha y^\prime+\beta z^\prime+\gamma t^\prime$ are equal, 
$u=u^\prime$, if and only if $x=x^\prime, y=y^\prime,
z=z^\prime, t=t^\prime$. 
If 
$u=x+\alpha y+\beta z+\gamma t, 
u^\prime=x^\prime+\alpha y^\prime+\beta z^\prime+\gamma t^\prime$
are polar fourcomplex numbers, 
the sum $u+u^\prime$ and the 
product $uu^\prime$ defined above can be obtained by applying the usual
algebraic rules to the sum 
$(x+\alpha y+\beta z+\gamma t)+ 
(x^\prime+\alpha y^\prime+\beta z^\prime+\gamma t^\prime)$
and to the product 
$(x+\alpha y+\beta z+\gamma t)
(x^\prime+\alpha y^\prime+\beta z^\prime+\gamma t^\prime)$,
and grouping of the resulting terms,
\begin{equation}
u+u^\prime=x+x^\prime+\alpha(y+y^\prime)+\beta(z+z^\prime)+\gamma(t+t^\prime),
\label{ch1a}
\end{equation}\index{sum!polar fourcomplex}
\begin{eqnarray}
\lefteqn{uu^\prime=
xx^\prime+yt^\prime+zz^\prime+ty^\prime+
\alpha(xy^\prime+yx^\prime+zt^\prime+tz^\prime)+
\beta(xz^\prime+yy^\prime+zx^\prime+tt^\prime)\nonumber}\\
&&+\gamma(xt^\prime+yz^\prime+zy^\prime+tx^\prime).
\label{ch1b}
\end{eqnarray}\index{product!polar fourcomplex}

If $u,u^\prime,u^{\prime\prime}$ are polar fourcomplex numbers, the multiplication 
is associative
\begin{equation}
(uu^\prime)u^{\prime\prime}=u(u^\prime u^{\prime\prime})
\label{ch2}
\end{equation}
and commutative
\begin{equation}
u u^\prime=u^\prime u ,
\label{ch3}
\end{equation}
as can be checked through direct calculation.
The polar fourcomplex zero is $0+\alpha\cdot 0+\beta\cdot 0+\gamma\cdot 0,$ 
denoted simply 0, 
and the polar fourcomplex unity is $1+\alpha\cdot 0+\beta\cdot 0+\gamma\cdot 0,$ 
denoted simply 1.

The inverse of the polar fourcomplex number 
$u=x+\alpha y+\beta z+\gamma t$ is a polar fourcomplex number
$u^\prime=x^\prime+\alpha y^\prime+\beta z^\prime+\gamma t^\prime$
having the property that
\begin{equation}
uu^\prime=1 .
\label{ch4}
\end{equation}
Written on components, the condition, Eq. (\ref{ch4}), is
\begin{equation}
\begin{array}{c}
xx^\prime+ty^\prime+zz^\prime+yt^\prime=1,\\
yx^\prime+xy^\prime+tz^\prime+zt^\prime=0,\\
zx^\prime+yy^\prime+xz^\prime+tt^\prime=0,\\
tx^\prime+zy^\prime+yz^\prime+xt^\prime=0.
\end{array}
\label{ch5}
\end{equation}
The system (\ref{ch5}) has the solution
\begin{equation}
x^\prime=\frac{x(x^2-z^2)+z(y^2+t^2)-2xyt}
{\nu} ,
\label{ch6a}
\end{equation}
\begin{equation}\index{inverse!polar fourcomplex}
y^\prime=\frac{-y(x^2+z^2)+t(y^2-t^2)+2xzt}
{\nu} ,
\label{ch6c}
\end{equation}
\begin{equation}
z^\prime=
\frac{-z(x^2-z^2)+x(y^2+t^2)-2yzt}
{\nu} ,
\label{ch6b}
\end{equation}
\begin{equation}
t^\prime=\frac{-t(x^2+z^2)-y(y^2-t^2)+2xyz}
{\nu} ,
\label{ch6d}
\end{equation}
provided that $\nu\not=0, $ where
\begin{equation}
\nu=x^4+z^4-y^4-t^4-2x^2z^2+2y^2t^2-4x^2yt-4yz^2t+4xy^2z+4xzt^2 .
\label{ch6e}
\end{equation}\index{inverse, determinant!polar fourcomplex}
The quantity $\nu$ can be written as
\begin{equation}
\nu=v_+ v_-\mu_+^2 ,
\label{ch7}
\end{equation}
where
\begin{equation}
v_+=x+y+z+t,\:v_-=x-y+z-t,
\label{ch8a}
\end{equation}
and
\begin{equation}
\mu_+^2=(x-z)^2+(y-t)^2 .
\label{ch8b}
\end{equation}
Then a polar fourcomplex number $q=x+\alpha y+\beta z+\gamma t$ has an inverse,
unless 
\begin{equation}
v_+=0 ,\:\:{\rm or}\:\: v_-=0,
\:\:{\rm or}\:\:\mu_+=0 . 
\label{ch9}
\end{equation}
The condition $v_+=0$ represents the 3-dimensional hyperplane $x+y+z+t=0$, the
condition $v_-=0$ represents the 3-dimensional hyperplane
$x-y+z-t=0$, and the condition $\mu_+=0$ represents the 2-dimensional
hyperplane $x=z, y=t$.
For arbitrary values of the variables $x,y,z,t$, the quantity $\nu$ can be
positive or negative. If $\nu\geq 0$, the quantity $\rho=\nu^{1/4}$
will be called amplitude of the polar fourcomplex number
$x+\alpha y+\beta z +\gamma t$.
Because of conditions (\ref{ch9}), these hyperplanes
will be called nodal hyperplanes. \index{nodal hyperplanes!polar fourcomplex}

It can be shown that if $uu^\prime=0$ then either $u=0$, or $u^\prime=0$, or 
the polar fourcomplex numbers $u, u^\prime$ belong to different members of the pairs
of orthogonal hypersurfaces listed further,
\begin{equation}
x+y+z+t=0 \:\:{\rm and} \:\: x^\prime=y^\prime=z^\prime=t^\prime,
\label{ch10a}
\end{equation}
\begin{equation}
x-y+z-t=0 \:\:{\rm and}\:\: x^\prime=-y^\prime=z^\prime=-t^\prime.
\label{ch10b}
\end{equation}
Divisors of zero also exist if the polar fourcomplex numbers $u,u^\prime$ belong to
different members of the pair of two-dimensional hypersurfaces,
\begin{equation}
x-z=0,\:y-t=0\:\: {\rm and}\:\:  x^\prime+z^\prime=0 , \: y^\prime+t^\prime=0.
\label{ch11}
\end{equation}

\subsection{Geometric representation of polar fourcomplex numbers}

The polar fourcomplex number $x+\alpha y+\beta z+\gamma t$ can be represented by 
the point $A$ of coordinates $(x,y,z,t)$. 
If $O$ is the origin of the four-dimensional space $x,y,z,t,$ the distance 
from $A$ to the origin $O$ can be taken as
\begin{equation}
d^2=x^2+y^2+z^2+t^2 .
\label{ch12}
\end{equation}\index{distance!polar fourcomplex}
The distance $d$ will be called modulus of the polar fourcomplex number $x+\alpha
y+\beta z +\gamma t$, $d=|u|$. 

If $u=x+\alpha y+\beta z +\gamma t, u_1=x_1+\alpha y_1+\beta z_1 +\gamma t_1,
u_2=x_2+\alpha y_2+\beta z_2 +\gamma t_2$, and $u=u_1u_2$, and if
\begin{equation}
s_{j+}=x_j+y_j+z_j+t_j, \: 
s_{j-}=x_j-y_j+z_j-t_j,  
\label{ch13}
\end{equation}
for $j=1,2$, it can be shown that
\begin{equation}
v_+=s_{1+}s_{2+} ,\:\:
v_-=s_{1-}s_{2-}. \:\:
\label{ch14}
\end{equation}\index{transformation of variables!polar fourcomplex}
The relations (\ref{ch14}) are a consequence of the identities
\begin{eqnarray}
\lefteqn{(x_1x_2+z_1z_2+t_1y_2+y_1t_2)+(x_1y_2+y_1x_2+z_1t_2+t_1z_2)
\nonumber}\\
&&+(x_1z_2+z_1x_2+y_1y_2+t_1t_2)+(x_1t_2+t_1x_2+z_1y_2+y_1z_2)\nonumber\\
&&=(x_1+y_1+z_1+t_1)(x_2+y_2+z_2+t_2),
\label{ch15}
\end{eqnarray}
\begin{eqnarray}
\lefteqn{(x_1x_2+z_1z_2+t_1y_2+y_1t_2)-(x_1y_2+y_1x_2+z_1t_2+t_1z_2)
\nonumber}\\
&&+(x_1z_2+z_1x_2+y_1y_2+t_1t_2)-(x_1t_2+t_1x_2+z_1y_2+y_1z_2)\nonumber\\
&&=(x_1+z_1-y_1-t_1)(x_2+z_2-y_2-t_2).
\label{ch16}
\end{eqnarray}

The differences 
\begin{equation}
v_1=x-z,\: \tilde v_1=y-t
\label{ch16a}
\end{equation}
can be written with the
aid of the radius $\mu_+$, Eq. (\ref{ch8b}), and of the azimuthal angle $\phi$, where
$0\leq\phi<2\pi$, as
\begin{equation}
v_1=\mu_+\cos\phi,\:\:\tilde v_1=\mu_+\sin\phi .
\label{ch17}
\end{equation}\index{azimuthal angle!polar fourcomplex}
The variables $v_+, v_-, v_1, \tilde v_1$ will be called canonical 
polar fourcomplex variables.\index{canonical variables!polar fourcomplex}

The distance $d$, Eq. (\ref{ch12}), can be written as
\begin{equation}
d^2
=\frac{1}{4}v_+^2+\frac{1}{4}v_-^2
+\frac{1}{2}\mu_+^2.  
\label{ch24d}
\end{equation}\index{distance, canonical variables!polar fourcomplex}
It can be shown that if $u_1=x_1+\alpha y_1+\beta z_1+\gamma t_1, 
u_2=x_2+\alpha y_2+\beta z_2+\gamma t_2$ are polar fourcomplex
numbers of polar radii and angles $\rho_{1-}, \phi_1$ and
respectively $\rho_{2-}, \phi_2$, then the polar radius $\rho$ and
the angle $\phi$ of the product polar fourcomplex number $u_1u_2$
are 
\begin{equation}
\mu_+=\rho_{1-}\rho_{2-}, 
\label{ch17a}
\end{equation}\index{transformation of variables!polar fourcomplex}
\begin{equation}
\phi=\phi_1+\phi_2.
\label{ch17b}
\end{equation}
The relation (\ref{ch17a}) is a consequence of the identity
\begin{eqnarray}
\lefteqn{\left[(x_1x_2+z_1z_2+y_1t_2+t_1y_2)-(x_1z_2+z_1x_2+y_1y_2+t_1t_2)\right]^2\nonumber}\\
&&+\left[(x_1y_2+y_1x_2+z_1t_2+t_1z_2)-(x_1t_2+t_1x_2+z_1y_2+y_1z_2)\right]^2\nonumber\\
&&=\left[(x_1-z_1)^2+(y_1-t_1)^2\right]\left[(x_2-z_2)^2+(y_2-t_2)^2\right],
\label{ch18}
\end{eqnarray}
and the relation (\ref{ch17b}) is a consequence of the identities
\begin{eqnarray}
\lefteqn{(x_1x_2+z_1z_2+y_1t_2+t_1y_2)-(x_1z_2+z_1x_2+y_1y_2+t_1t_2)\nonumber}\\
&&=(x_1-z_1)(x_2-z_2)-(y_1-t_1)(y_2-t_2) ,
\label{ch19a}
\end{eqnarray}
\begin{eqnarray}
\lefteqn{(x_1y_2+y_1x_2+z_1t_2+t_1z_2)-(x_1t_2+t_1x_2+z_1y_2+y_1z_2)\nonumber}\\
&&=(y_1-t_1)(x_2-z_2)+(x_1-z_1)(y_2-t_2) .
\label{ch19b}
\end{eqnarray}

A consequence of Eqs. (\ref{ch14}) and (\ref{ch17a}) is that if $u=u_1u_2$,
and $\nu_j=s_j s_j^{\prime\prime} \rho_{j-}$ , where $j=1,2$, then
\begin{equation}
\nu=\nu_1\nu_2 .
\label{ch20}
\end{equation}

The angles $\theta_+, \theta_-$ between the line $OA$ and the $v_+$ and
respectively $v_-$ axes are
\begin{equation}
\tan\theta_+=\frac{\sqrt{2}\mu_+}{v_+},
\tan\theta_-=\frac{\sqrt{2}\mu_+}{v_-},
\label{ch20b}
\end{equation}\index{polar angles!polar fourcomplex}
where $0\leq\theta_+\leq\pi,\;0\leq\theta_-\leq\pi $.
The variable $\mu_+$ can be expressed with the aid of Eq. (\ref{ch24d}) as
\begin{equation}
\mu_+^2=2d^2\left(1+\frac{1}{\tan^2\theta_+}
+\frac{1}{\tan^2\theta_-}\right)^{-1}.
\label{ch20c}
\end{equation}
The coordinates $x,y,z,t$ can then be expressed in terms of the distance
$d$, of the polar angles $\theta_+,\theta_-$ and of the azimuthal angle $\phi$
as 
\begin{equation}
x=\frac{\mu_+(\tan\theta_++\tan\theta_-)}
{2\sqrt{2}\tan\theta_+\tan\theta_-}+\frac{1}{2}\mu_+\cos\phi,
\label{ch20i}
\end{equation}
\begin{equation}
y=\frac{\mu_+(-\tan\theta_++\tan\theta_-)}
{2\sqrt{2}\tan\theta_+\tan\theta_-}+\frac{1}{2}\mu_+\sin\phi,
\label{ch20k}
\end{equation}
\begin{equation}
z=\frac{\mu_+(\tan\theta_++\tan\theta_-)}
{2\sqrt{2}\tan\theta_+\tan\theta_-}-\frac{1}{2}\mu_+\cos\phi,
\label{ch20j}
\end{equation}
\begin{equation}
t=\frac{\mu_+(-\tan\theta_++\tan\theta_-)}
{2\sqrt{2}\tan\theta_+\tan\theta_-}-\frac{1}{2}\mu_+\sin\phi.
\label{ch20l}
\end{equation}
If $u=u_1u_2$, then Eqs. (\ref{ch14}) and (\ref{ch17a}) imply that
\begin{equation}
\tan\theta_+=\frac{1}{\sqrt{2}}\tan\theta_{1+}\tan\theta_{2+},\;
\tan\theta_-=\frac{1}{\sqrt{2}}\tan\theta_{1-}\tan\theta_{2-},
\label{ch20d}
\end{equation}
where 
\begin{equation}
\tan\theta_{j+}=\frac{\sqrt{2}\rho_{j+}}{s_j},
\tan\theta_{j-}=\frac{\sqrt{2}\rho_{j-}}{s_j^{\prime\prime}}.
\label{ch20e}
\end{equation}

An alternative choice of the angular variables is 
\begin{equation}
\mu_+=\sqrt{2}d\cos\theta,\;
v_+=2d\sin\theta\cos\lambda,\;v_-=2d\sin\theta\sin\lambda ,
\label{ch20f}
\end{equation}
where $0\leq\theta\leq\pi/2, 0\leq\lambda<2\pi$.
If $u=u_1u_2$, then
\begin{equation}
\tan\lambda=\tan\lambda_1\tan\lambda_2,\;
d\cos\theta=\sqrt{2}d_1d_2\cos\theta_1\cos\theta_2,
\label{ch20g}
\end{equation}
where 
\begin{equation}
\rho_{j-}=\sqrt{2}d_j\cos\theta_j,\;
s_j=2d_j\sin\theta_j\cos\lambda_j,\;
s_j^{\prime\prime}=2d_j\sin\theta_j\sin\lambda_j, 
\label{ch20h}
\end{equation}
for $j=1,2$.
The coordinates $x,y,z,t$ can then be expressed in terms of the distance
$d$, of the polar angles $\theta_,\lambda$ and of the azimuthal angle $\phi$
as
\begin{equation}
x=\frac{d}{\sqrt{2}}\sin\theta\sin(\lambda+\pi/4)+\frac{d}{\sqrt{2}}\cos\theta\cos\phi,
\label{ch20m}
\end{equation}
\begin{equation}
y=\frac{d}{\sqrt{2}}\sin\theta\cos(\lambda+\pi/4)+\frac{d}{\sqrt{2}}\cos\theta\sin\phi,
\label{ch20o}
\end{equation}
\begin{equation}
z=\frac{d}{\sqrt{2}}\sin\theta\sin(\lambda+\pi/4)-\frac{d}{\sqrt{2}}\cos\theta\cos\phi,
\label{ch20n}
\end{equation}
\begin{equation}
t=\frac{d}{\sqrt{2}}\sin\theta\cos(\lambda+\pi/4)-\frac{d}{\sqrt{2}}\cos\theta\sin\phi.
\label{ch20p}
\end{equation}

The polar fourcomplex numbers
\begin{equation}
e_+=\frac{1+\alpha+\beta+\gamma}{4},\:
e_-=\frac{1-\alpha+\beta-\gamma}{4},
\label{ch21}
\end{equation}\index{canonical base!polar fourcomplex}
have the property that
\begin{equation}
e_+^2=e_+, \: 
e_-^2=e_-, \:e_+e_-=0.
\label{ch22}
\end{equation}
The polar fourcomplex numbers 
\begin{equation}
e_1=\frac{1-\beta}{2}, \tilde e_1=\frac{\alpha-\gamma}{2}
\label{ch22b}
\end{equation}
have the property that
\begin{equation}
e_1^2=e_1,\:\:
\tilde e_1^2=-e_1,\:\:
e_1\tilde e_1=\tilde e_1 .
\label{ch23}
\end{equation}
The polar fourcomplex numbers $ e_+, e_-$ are orthogonal to
$e_1, \tilde e_1$, 
\begin{equation}
e_+\:e_1=0,\:\:e_+\:\tilde e_1=0,\:\:
e_-\:e_1=0,
\:\:e_-\:\tilde e_1=0.\:\: 
\label{ch23b}
\end{equation}
The polar fourcomplex number $q=x+\alpha y+\beta z+\gamma t$ can then be
written as 
\begin{equation}
x+\alpha y+\beta z+\gamma t
=v_+e_++v_-e_-+v_1e_1
+\tilde v_1\tilde e_1.
\label{ch24}
\end{equation}\index{canonical form!polar fourcomplex}
The ensemble $e_+, e_-, e_1, \tilde e_1$ will be called the canonical
polar fourcomplex base, and Eq. (\ref{ch24}) gives the canonical form of the
polar fourcomplex number.
Thus, the product of the polar fourcomplex numbers $u, u^\prime$ can be
expressed as 
\begin{equation}
uu^\prime=v_+v_+^\prime e_++v_-v_-^\prime e_-
+(v_1v_1^\prime-\tilde v_1\tilde v_1^\prime)e_1+
(v_1\tilde v_1^\prime+v_1^\prime\tilde v_1)\tilde e_1,
\label{ch24b}
\end{equation}
where $v_+^\prime=x^\prime+y^\prime+z^\prime+t^\prime, 
v_-^\prime=x^\prime-y^\prime+z^\prime-t^\prime,
v_1^\prime=x^\prime-y^\prime, \tilde v_1^\prime=z^\prime-t^\prime$.
The moduli of the bases used in Eq. (\ref{ch24}) are
\begin{equation}
|e_+|=\frac{1}{2},\;|e_-|=\frac{1}{2},\;
\left|e_1\right|=\frac{1}{\sqrt{2}},\;
\left|\tilde e_1\right|=\frac{1}{\sqrt{2}}.
\label{ch24c}
\end{equation}

The fact that the amplitude of the product is equal to the product of the 
amplitudes, as written in Eq. (\ref{ch20}), can 
be demonstrated also by using a representation of the multiplication of the 
polar fourcomplex numbers by matrices, in which the polar fourcomplex number $u=x+\alpha
y+\beta z+\gamma t$ is represented by the matrix
\begin{equation}
A=\left(\begin{array}{cccc}
x &y &z &t\\
t&x &y &z\\
z&t&x &y\\
y&z&t&x 
\end{array}\right) .
\label{ch25a}
\end{equation}\index{matrix representation!polar fourcomplex}
The product $u=x+\alpha y+\beta z+\gamma t$ of the polar fourcomplex numbers
$u_1=x_1+\alpha y_1+\beta z_1+\gamma t_1, u_2=x_2+\alpha y_2+\beta z_2+\gamma
t_2$, can be represented by the matrix multiplication 
\begin{equation}
A=A_1A_2.
\label{ch25b}
\end{equation}
It can be checked that the determinant ${\rm det}(A)$ of the matrix $A$ is
\begin{equation}
{\rm det}A = \nu .
\label{ch25c}
\end{equation}
The identity (\ref{ch20}) is then a consequence of the fact the determinant 
of the product of matrices is equal to the product of the determinants 
of the factor matrices.

\subsection{The polar fourdimensional cosexponential functions}

The exponential function of a hypercomplex variable $u$ and the addition
theorem for the exponential function have been written in Eqs. 
~(1.35)-(1.36).
If $u=x+\alpha y+\beta z+\gamma t$, then $\exp u$ can be calculated as 
$\exp u=\exp x \cdot \exp (\alpha y) \cdot \exp (\beta z)\cdot \exp (\gamma t)$. 
According to Eq. (\ref{ch1}),  
\begin{eqnarray}
\begin{array}{l}
\alpha^{4m}=1, \alpha^{4m+1}=\alpha, \alpha^{4m+2}=\beta, \alpha^{4m+3}=\gamma,  \\
\beta^{2m}=1, \beta^{2m+1}=\beta, \\
\gamma^{4m}=1, \gamma^{4m+1}=\gamma, \gamma^{4m+2}=\beta, \gamma^{4m+3}=\alpha, \\
\end{array}
\label{ch28}
\end{eqnarray}\index{complex units, powers of!polar fourcomplex}
where $m$ is a natural number,
so that $\exp (\alpha y), \: \exp(\beta z)$ and $\exp(\gamma z)$ can be written
as 
\begin{equation}
\exp (\beta z) = \cosh z +\beta \sinh z ,
\label{ch29}
\end{equation}\index{exponential, expressions!polar fourcomplex}
and
\begin{equation}
\exp (\alpha y) = g_{40}(y)+\alpha g_{41}(y)+\beta g_{42}(y) +\gamma g_{43}(y) ,  
\label{ch30a}
\end{equation}
\begin{equation}
\exp (\gamma t) = g_{40}(t)+\gamma g_{41}(t)+\beta g_{42}(t) +\alpha g_{43}(t) ,
\label{ch30b}
\end{equation}
where the four-dimensional cosexponential functions $g_{40}, g_{41}, g_{42}, g_{43}$ are
defined by the series 
\begin{equation}
g_{40}(x)=1+x^4/4!+x^8/8!+\cdots ,
\label{ch30c}
\end{equation}\index{cosexponential function, polar fourcomplex!definitions}
\begin{equation}
g_{41}(x)=x+x^5/5!+x^9/9!+\cdots,
\label{ch30d}
\end{equation}
\begin{equation}
g_{42}(x)=x^2/2!+x^6/6!+x^{10}/10!+\cdots ,
\label{ch30e}
\end{equation}
\begin{equation}
g_{43}(x)=x^3/3!+x^7/7!+x^{11}/11!+\cdots .
\label{ch30f}
\end{equation}
The functions $g_{40}, g_{42}$ are even, the functions $g_{41}, g_{43}$ are
odd, 
\begin{equation}
g_{40}(-u)=g_{40}(u),\:g_{42}(-u)=g_{42}(u),\:g_{41}(-u)=-g_{41}(u),\:g_{43}(-u)=-g_{43}(u).
\label{ch30feo}
\end{equation}\index{cosexponential functions, polar fourcomplex!parity}
It can be seen from Eqs. (\ref{ch30c})-(\ref{ch30f}) that
\begin{equation}
g_{40}+g_{41}+g_{42}+g_{43}=e^x , \:g_{40}-g_{41}+g_{42}-g_{43}=e^{-x} ,
\label{ch30fep}
\end{equation}
and
\begin{equation}
g_{40}-g_{42}=\cos x , \: g_{41}-g_{43}=\sin x ,
\label{ch30feq}
\end{equation}
so that
\begin{equation}
(g_{40}+g_{41}+g_{42}+g_{43})(g_{40}-g_{41}+g_{42}-g_{43})
\left[(g_{40}-g_{42})^2+(g_{41}-g_{43})^2\right]=1 ,
\label{ch30fer}
\end{equation}
which can be also written as
\begin{equation}
g_{40}^4-g_{41}^4+g_{42}^4-g_{43}^4-2(g_{40}^2g_{42}^2-g_{41}^2g_{43}^2)
-4(g_{40}^2g_{41}g_{43}+g_{42}^2g_{41}g_{43}-g_{41}^2g_{40}g_{42}-g_{43}^2g_{40}g_{42})=1. 
\label{ch30fes}
\end{equation}
The combination of terms in Eq. (\ref{ch30fes}) is similar to that in Eq.
(\ref{ch6e}). 

Addition theorems for the four-dimensional cosexponential functions can be
obtained from the relation $\exp \alpha(x+y)=\exp\alpha x\cdot\exp\alpha y $, by
substituting the expression of the exponentials as given in Eq. (\ref{ch30a}),
\begin{equation}
g_{40}(x+y)=g_{40}(x)g_{40}(y)+g_{41}(x)g_{43}(y)+g_{42}(x)g_{42}(y)+g_{43}(x)g_{41}(y) ,
\label{ch30g}
\end{equation}\index{cosexponential functions, polar fourcomplex!addition theorems}
\begin{equation}
g_{41}(x+y)=g_{40}(x)g_{41}(y)+g_{41}(x)g_{40}(y)+g_{42}(x)g_{43}(y)+g_{43}(x)g_{42}(y) ,
\label{ch30h}
\end{equation}
\begin{equation}
g_{42}(x+y)=g_{40}(x)g_{42}(y)+g_{41}(x)g_{41}(y)+g_{42}(x)g_{40}(y)+g_{43}(x)g_{43}(y) ,
\label{ch30i}
\end{equation}
\begin{equation}
g_{43}(x+y)=g_{40}(x)g_{43}(y)+g_{41}(x)g_{42}(y)+g_{42}(x)g_{41}(y)+g_{43}(x)g_{40}(y) .
\label{ch30j}
\end{equation}
For $x=y$ the relations (\ref{ch30g})-(\ref{ch30j}) take the form
\begin{equation}
g_{40}(2x)=g_{40}^2(x)+g_{42}^2(x)+2g_{41}(x)g_{43}(x) ,
\label{ch30gg}
\end{equation}
\begin{equation}
g_{41}(2x)=2g_{40}(x)g_{41}(x)+2g_{42}(x)g_{43}(x) ,
\label{ch30hh}
\end{equation}
\begin{equation}
g_{42}(2x)=g_{41}^2(x)+g_{43}^2(x)+2g_{40}(x)g_{42}(x) ,
\label{ch30ii}
\end{equation}
\begin{equation}
g_{43}(2x)=2g_{40}(x)g_{43}(x)+2g_{41}(x)g_{42}(x) .
\label{ch30jj}
\end{equation}
For $x=-y$ the relations (\ref{ch30g})-(\ref{ch30j}) and (\ref{ch30feo}) yield
\begin{equation}
g_{40}^2(x)+g_{42}^2(x)-2g_{41}(x)g_{43}(x)=1 ,
\label{ch30eog}
\end{equation}
\begin{equation}
g_{41}^2(x)+g_{43}^2(x)-2g_{40}(x)g_{42}(x)=0 .
\label{ch30eoi}
\end{equation}
From Eqs. (\ref{ch29})-(\ref{ch30b}) it can be shown that, for $m$ integer,
\begin{equation}
(\cosh z +\beta \sinh z)^m=\cosh mz +\beta \sinh mz ,
\label{ch2929}
\end{equation}
and
\begin{equation}
[g_{40}(y)+\alpha g_{41}(y)+\beta g_{42}(y) +\gamma g_{43}(y)]^m= 
g_{40}(my)+\alpha g_{41}(my)+\beta g_{42}(my) +\gamma g_{43}(my),  
\label{ch30a30a}
\end{equation}
\begin{equation}
[g_{40}(t)+\gamma g_{41}(t)+\beta g_{42}(t) +\alpha g_{43}(t)]^m=
g_{40}(mt)+\gamma g_{41}(mt)+\beta g_{42}(mt) +\alpha g_{43}(mt) .
\label{ch30b30b}
\end{equation}
Since
\begin{equation}
(\alpha+\gamma)^{2m}=2^{2m-1}(1+\beta),\: (\alpha+\gamma)^{2m+1}=2^{2m}(\alpha+\gamma) ,
\label{ch30k}
\end{equation}
it can be shown from the definition of the exponential function, Eq.
(1.35) that
\begin{equation}
\exp(\alpha+\gamma)x=e_1+\frac{1+\beta}{2}\cosh 2x+\frac{\alpha+\gamma}{2}
\sinh 2x .
\label{ch30l}
\end{equation}
Substituting in the relation $\exp(\alpha+\gamma)x=\exp\alpha x\exp\gamma x$ the
expression of the exponentials from Eqs. (\ref{ch30a}), (\ref{ch30b}) and
(\ref{ch30l}) yields  
\begin{equation}
g_{40}^2+g_{41}^2+g_{42}^2+g_{43}^2=\frac{1+\cosh 2x}{2} ,
\label{ch30m}
\end{equation}
\begin{equation}
g_{40} g_{42} +g_{41} g_{43}=\frac{-1+\cosh 2x}{4} ,
\label{ch30mm}
\end{equation}
\begin{equation}
g_{40}g_{41}+g_{40}g_{43}+g_{41}g_{42}+g_{42}g_{43}=\frac{1}{2}\sinh 2x ,
\label{ch30n}
\end{equation}
where $g_{40}, g_{41}, g_{42}, g_{43}$ are functions of $x$.

Similarly, since
\begin{equation}
(\alpha-\gamma)^{2m}=(-1)^m 2^{2m-1}(1-\beta),\:
(\alpha-\gamma)^{2m+1}=(-1)^m2^{2m}(\alpha-\gamma) , 
\label{ch30k30k}
\end{equation}
it can be shown from the definition of the exponential function, Eq. (1.35)
that
\begin{equation}
\exp(\alpha-\gamma)x=\frac{1+\beta}{2}+e_1\cos 2x+\tilde e_1
\sin 2x .
\label{ch30l30l}
\end{equation}
Substituting in the relation $\exp(\alpha-\gamma)x=\exp\alpha x\exp(-\gamma x)$ the
expression of the exponentials from Eqs. (\ref{ch30a}), (\ref{ch30b}) and
(\ref{ch30l30l}) yields  
\begin{equation}
g_{40}^2-g_{41}^2+g_{42}^2-g_{43}^2=\frac{1+\cos 2x}{2} ,
\label{ch30m30m}
\end{equation}
\begin{equation}
g_{40} g_{42} -g_{41} g_{43}=\frac{1-\cos 2x}{2} ,
\label{ch30mm30mm}
\end{equation}
\begin{equation}
g_{40}g_{41}-g_{40}g_{43}-g_{41}g_{42}+g_{42}g_{43}=\frac{1}{2}\sin 2x ,
\label{ch30n30n}
\end{equation}
where $g_{40}, g_{41}, g_{42}, g_{43}$ are functions of $x$.

The expressions of the four-dimensional cosexponential
functions are
\begin{equation}
g_{40}(x)=\frac{1}{2}(\cosh x + \cos x) ,
\label{ch30c30c}
\end{equation}\index{cosexponential functions, polar fourcomplex!expressions}
\begin{equation}
g_{41}(x)=\frac{1}{2}(\sinh x+\sin x),
\label{ch30d30d}
\end{equation}
\begin{equation}
g_{42}(x)=\frac{1}{2}(\cosh x - \cos x) ,
\label{ch30e30e}
\end{equation}
\begin{equation}
g_{43}(x)=\frac{1}{2}(\sinh x-\sin x).
\label{ch30f30f}
\end{equation}
The graphs of these four-dimensional cosexponential functions are shown in 
Fig. \ref{fig13}.

\begin{figure}
\begin{center}
\epsfig{file=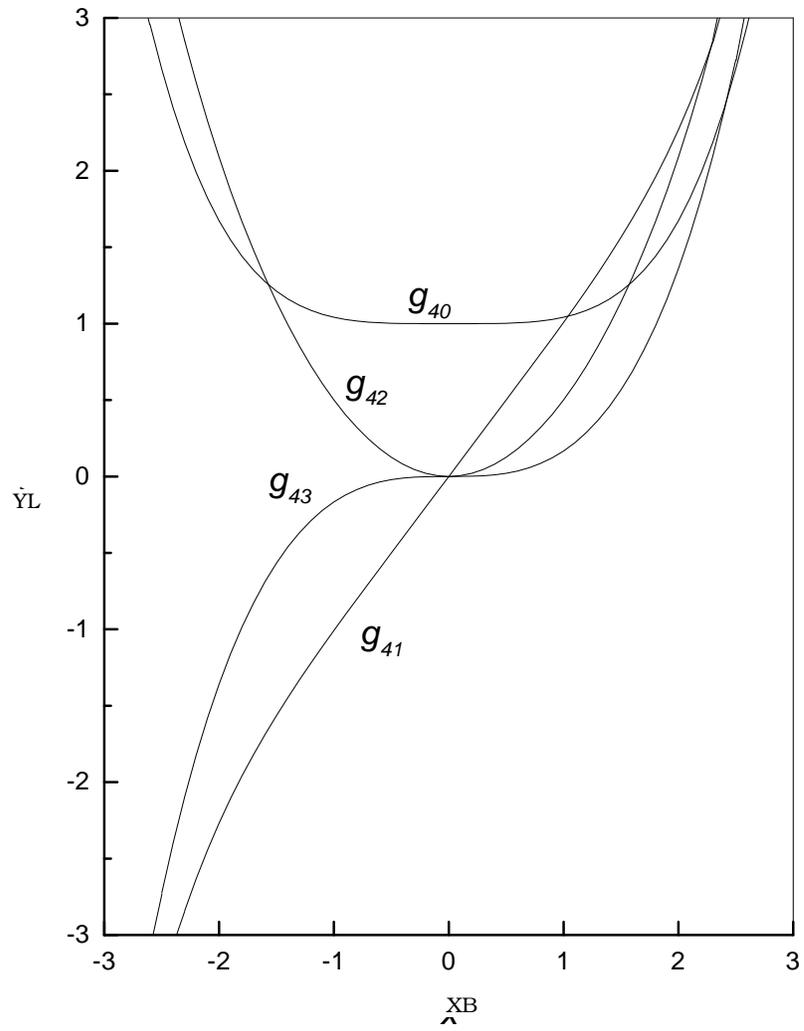,width=12cm}
\caption{The polar fourdimensional cosexponential functions $g_{40}, g_{41}, g_{42}, g_{43}$.}
\label{fig13}
\end{center}
\end{figure}

It can be checked that the cosexponential functions are solutions of the
fourth-order differential equation
\begin{equation}
\frac{{\rm d}^4\zeta}{{\rm d}u^4}=\zeta ,
\label{ch30x}
\end{equation}\index{cosexponential functions, polar fourcomplex!differential equations}
whose solutions are of the form $\zeta(u)=Ag_{40}(u)+Bg_{41}(u)+Cg_{42(u)}+Dg_{43}(u).$
It can also be checked that the derivatives of the cosexponential functions are
related by
\begin{equation}
\frac{dg_{40}}{dw}=g_{43}, \:
\frac{dg_{41}}{dw}=g_{40}, \:
\frac{dg_{42}}{dw}=g_{41} ,
\frac{dg_{43}}{dw}=g_{42} .
\label{ch30x30x}
\end{equation}


\subsection{The exponential and trigonometric forms 
of a polar fourcomplex number}

The polar fourcomplex numbers $u=x+\alpha y+\beta z+\gamma t$ for which
$v_+=x+y+z+t>0, \:  v_-=x-y+z-t>0$ can be written in the form 
\begin{equation}
x+\alpha y+\beta z+\gamma t=e^{x_1+\alpha y_1+\beta z_1+\gamma t_1} .
\label{ch31-34}
\end{equation}

The expressions of $x_1, y_1, z_1, t_1$ as functions of 
$x, y, z, t$ can be obtained by
developing $e^{\alpha y_1}, e^{\beta z_1}$ and $e^{\gamma t_1}$ with the aid of
Eqs. (\ref{ch29})-(\ref{ch30b}), by multiplying these expressions and separating
the hypercomplex components, and then substituting the expressions of the
four-dimensional cosexponential functions, Eqs. (\ref{ch30c30c})-(\ref{ch30f30f}), 
\begin{equation}
x+y+z+t=e^{x_1+y_1+z_1+t_1},
\label{ch35}
\end{equation}
\begin{equation}
x-z=e^{x_1-z_1}\cos(y_1-t_1),
\label{ch36}
\end{equation}
\begin{equation}
x-y+z-t=e^{x_1-y_1+z_1-t_1},
\label{ch37}
\end{equation}
\begin{equation}
y-t=e^{x_1-z_1}\sin(y_1-t_1).
\label{ch38}
\end{equation}
It can be shown from Eqs. (\ref{ch35})-(\ref{ch38}) that
\begin{equation}
x_1=\frac{1}{2} \ln(\mu_+\mu_-) , \:
y_1=\frac{1}{2}(\phi+\omega) ,\:
z_1=\frac{1}{2}\ln\frac{\mu_-}{\mu_+},\:
t_1=\frac{1}{2}(\phi-\omega),
\label{ch39}
\end{equation}
where
\begin{equation}
\mu_-^2=(x+z)^2-(y+t)^2=v_+v_-,\: \:v_+>0, \:\:v_->0.
\label{ch39a}
\end{equation}
The quantitites $\phi$ and $\omega$ are determined by
\begin{equation}
\cos\phi=(x-z)/\mu_+,\:\sin\phi=(y-t)/\mu_+ 
\label{ch39b}
\end{equation}
and
\begin{equation}
\cosh\omega=(x+z)/\mu_-,\:
\sinh\omega=(y+t)/\mu_- .
\label{ch39c}
\end{equation}
The explicit form of $\omega$ is
\begin{equation}
\omega=\frac{1}{2}\ln\frac{x+y+z+t}{x-y+z-t}.
\label{ch39cc}
\end{equation}
If $u=u_1u_2$, and $\mu_{j-}^2=(x_j+z_j)^2-(y_j+t_j)^2, j=1,2$, it can be checked
with the aid of Eqs. (\ref{ch14}) that 
\begin{equation}
\mu_-=\mu_{1-}\mu_{2-}. 
\label{ch39d}
\end{equation}
Moreover, if
$\cosh\omega_j=(x_j+z_j)/\mu_{j-},\:\sinh\omega_j=(y_j+t_j)/\mu_{j-}, j=1,2$,
it can be checked that 
\begin{equation}
\omega=\omega_1+\omega_2.
\label{ch39e}
\end{equation}
The relation (\ref{ch39e}) is a consequence of the identities
\begin{eqnarray}
\lefteqn{(x_1x_2+z_1z_2+y_1t_2+t_1y_2)+(x_1z_2+z_1x_2+y_1y_2+t_1t_2)\nonumber}\\
&&=(x_1+z_1)(x_2+z_2)+(y_1+t_1)(y_2+t_2) ,
\label{ch19c}
\end{eqnarray}
\begin{eqnarray}
\lefteqn{(x_1y_2+y_1x_2+z_1t_2+t_1z_2)+(x_1t_2+t_1x_2+z_1y_2+y_1z_2)\nonumber}\\
&&=(y_1+t_1)(x_2+z_2)+(x_1+z_1)(y_2+t_2) .
\label{ch19d}
\end{eqnarray}

According to Eq. (\ref{ch39b}), $\phi$ is a cyclic variable, $0\leq\phi<2\pi$.
As it has been assumed that $v_+>0, v_->0$, it follows that $x+z>0$
and $x+z>|y+t|$. The range of the variable $\omega$ is $-\infty<\omega<\infty$.
The exponential form of the polar fourcomplex number $u$ is then
\begin{equation}
u=\rho\exp\left[\frac{1}{2}\beta \ln\frac{\mu_-}{\mu_+}
+\frac{1}{2}\alpha (\omega+\phi) 
+\frac{1}{2}\gamma (\omega-\phi)\right] ,
\label{ch40}
\end{equation}\index{exponential form!polar fourcomplex}
where
\begin{equation}
\rho=(\mu_+\mu_-)^{1/2}.
\label{ch40a}
\end{equation}
If $u=u_1u_2$, and $\rho_j=(\mu_{j+}\mu_{j-})^{1/2}, j=1,2$, then 
from Eqs. (\ref{ch17a}) and (\ref{ch39d}) it results that
\begin{equation}
\rho=\rho_1\rho_2.
\label{ch40bb}
\end{equation}\index{transformation of variables!polar fourcomplex}

It can be checked with the aid of Eq. (\ref{ch23}) that
\begin{equation}
\exp\left(\tilde e_1\phi\right)
=\frac{1+\beta}{2}+e_1\cos\phi+\tilde e_1\sin\phi,
\label{ch40x}
\end{equation}
which shows that $e^{(\alpha-\gamma)\phi/2}$ is a periodic function of $\phi$,
with period $2\pi$.
The modulus has the property that
\begin{equation}
\left|u\exp\left(\tilde e_1\phi\right)\right|=|u|.
\label{ch40bc}
\end{equation}

By introducing in Eq. (\ref{ch40}) the polar angles $\theta_+, \theta_-$
defined in Eqs. (\ref{ch20b}), the exponential form of the fourcomplex number
$u$ becomes
\begin{equation}
u=\rho\exp\left[\frac{1}{4}(\alpha+\beta+\gamma) \ln\frac{\sqrt{2}}{\tan\theta_+}
-\frac{1}{4}(\alpha-\beta+\gamma) \ln\frac{\sqrt{2}}{\tan\theta_-}
+\tilde e_1\phi\right] ,
\label{ch40c}
\end{equation}
where $0<\theta_+<\pi/2, 0<\theta_-<\pi/2$.
The relation between the amplitude $\rho$, Eq. (\ref{ch40a}), and the
distance $d$, Eq. (\ref{ch12}), is, according to Eqs. (\ref{ch20b}) and
(\ref{ch20c}), 
\begin{equation}
\rho=\frac{2^{3/4}d}{\left(\tan\theta_+\tan\theta_-\right)^{1/4}}
\left(1+\frac{1}{\tan^2\theta_+}+\frac{1}{\tan^2\theta_-}\right)^{-1/2}.
\label{ch40d}
\end{equation}
Using the properties of the vectors $e_+, e_-$ written in Eq.
(\ref{ch22}), the first part of the exponential, Eq. (\ref{ch40c}) can be
developed as
\begin{eqnarray}
\lefteqn{\exp\left[\frac{1}{4}(\alpha+\beta+\gamma) \ln\frac{\sqrt{2}}{\tan\theta_+}
-\frac{1}{4}(\alpha-\beta+\gamma) \ln\frac{\sqrt{2}}{\tan\theta_-}\right] \nonumber}\\
&&=\left(\frac{1}{2}\tan\theta_+\tan\theta_-\right)^{1/4}
\left(e_1+e_+\frac{\sqrt{2}}{\tan\theta_+}
+e_-\frac{\sqrt{2}}{\tan\theta_-}\right).
\label{ch40e}
\end{eqnarray}
The fourcomplex number $u$, Eq. (\ref{ch40c}), can then be written as
\begin{equation}
u=d\sqrt{2}
\left(1+\frac{1}{\tan^2\theta_+}+\frac{1}{\tan^2\theta_-}\right)^{-1/2}
\left(e_1+e_+\frac{\sqrt{2}}{\tan\theta_+}
+e_-\frac{\sqrt{2}}{\tan\theta_-}\right)
\exp\left(\tilde e_1\phi\right),
\label{ch40f}
\end{equation}\index{trigonometric form!polar fourcomplex}
which is the trigonometric form of the fourcomplex number $u$.

The polar angles $\theta_+, \theta_-$, Eq. (\ref{ch20b}), can be expressed in terms
of the variables $\theta, \lambda$, Eq. (\ref{ch20f}), as
\begin{equation}
\tan\theta_+=\frac{1}{\tan\theta\cos\lambda},\;
\tan\theta_-=\frac{1}{\tan\theta\sin\lambda},
\label{ch40g}
\end{equation}
so that
\begin{equation}
1+\frac{1}{\tan^2\theta_+}+\frac{1}{\tan^2\theta_-}=\frac{1}{\cos^2\theta}.
\label{ch40h}
\end{equation}
The exponential form of the fourcomplex number $u$, written in terms of the
amplitude $\rho$ and of the angles $\theta, \lambda, \phi$ is
\begin{equation}
u=\rho\exp\left[\frac{1}{4}(\alpha+\beta+\gamma)
\ln(\sqrt{2}\tan\theta\cos\lambda) 
-\frac{1}{4}(\alpha-\beta+\gamma) \ln(\sqrt{2}\tan\theta\sin\lambda)
+\tilde e_1\phi
\right],
\label{ch40i}
\end{equation}
where $0\leq\lambda<\pi/2$.
The trigonometric form of the fourcomplex number $u$, written in terms of the
amplitude $\rho$ and of the angles $\theta, \lambda, \phi$ is
\begin{equation}
u=d\sqrt{2}
\left(e_1\cos\theta+e_+\sqrt{2}\sin\theta\cos\lambda
+e_-\sqrt{2}\sin\theta\sin\lambda\right)
\exp\left(\tilde e_1\phi\right).
\label{ch40j}
\end{equation}

If $u=u_1u_2$, it can be shown with the aid of the trigonometric form, Eq.
(\ref{ch40f}), that the modulus of the product as a function of the polar
angles is 
\begin{eqnarray}
\lefteqn{d^2=4d_1^2d_2^2\left(\frac{1}{2}
+\frac{1}{\tan^2\theta_{1+}\tan^2\theta_{2+}}
+\frac{1}{\tan^2\theta_{1-}\tan^2\theta_{2-}}\right)\nonumber}\\
&&\left(1+\frac{1}{\tan^2\theta_{1+}}+\frac{1}{\tan^2\theta_{1-}}\right)^{-1}
\left(1+\frac{1}{\tan^2\theta_{2+}}+\frac{1}{\tan^2\theta_{2-}}\right)^{-1}.
\label{ch40k}
\end{eqnarray}\index{transformation of variables!polar fourcomplex}
The modulus $d$ of the product $u_1u_2$ can be expressed alternatively in terms
of the angles $\theta, \lambda, \phi$ with the aid of the
trigonometric form, Eq. (\ref{ch40j}), as
\begin{eqnarray}
\lefteqn{d^2=4d_1^2d_2^2\left(\frac{1}{2}\cos^2\theta_1\cos^2\theta_2
+\sin^2\theta_1\sin^2\theta_2\cos^2\lambda_1\cos^2\lambda_2
+\sin^2\theta_1\sin^2\theta_2\sin^2\lambda_1\sin^2\lambda_2\right).\nonumber}\\
&&
\label{ch40l}
\end{eqnarray}

\subsection{Elementary functions of polar fourcomplex variables}

The logarithm $u_1$ of the polar fourcomplex number $u$, $u_1=\ln u$, can be defined
as the solution of the equation
\begin{equation}
u=e^{u_1} ,
\label{ch41}
\end{equation}
written explicitly previously in Eq. (\ref{ch31-34}), for $u_1$ as a function
of $u$. From Eq. (\ref{ch40}) it results that 
\begin{equation}
\ln u=\ln\rho+
\frac{1}{2}\beta \ln\frac{\mu_-}{\mu_+}
+\frac{1}{2}\alpha (\omega+\phi) 
+\frac{1}{2}\gamma (\omega-\phi).
\label{ch42}
\end{equation}
If the fourcomplex number $u$ is written in terms of the
amplitude $\rho$ and of the angles $\theta_+, \theta_-, \phi$, the logarithm
is 
\begin{equation}
\ln u=\ln \rho+
\frac{1}{4}(\alpha+\beta+\gamma) \ln\frac{\sqrt{2}}{\tan\theta_+}
-\frac{1}{4}(\alpha-\beta+\gamma) \ln\frac{\sqrt{2}}{\tan\theta_-}
+\tilde e_1\phi,
\label{ch42b}
\end{equation}\index{logarithm!polar fourcomplex}
where $0<\theta_+<\pi/2, 0<\theta_+<\pi/2$.
If the fourcomplex number $u$ is written in terms of the
amplitude $\rho$ and of the angles $\theta, \lambda, \phi$, the logarithm is
\begin{equation}
\ln u=\ln \rho+
\frac{1}{4}(\alpha+\beta+\gamma)
\ln(\sqrt{2}\tan\theta\cos\lambda) 
-\frac{1}{4}(\alpha-\beta+\gamma) \ln(\sqrt{2}\tan\theta\sin\lambda)
+\tilde e_1\phi,
\label{ch42c}
\end{equation}
where $0<\theta<\pi/2, 0\leq\lambda<\pi/2$.
The logarithm is multivalued because of the term proportional to $\phi$.
It can be inferred from Eq. (\ref{ch42b}) that
\begin{equation}
\ln(u_1u_2)=\ln u_1+\ln u_2 ,
\label{ch43}
\end{equation}
up to multiples of $\pi(\alpha-\gamma)$.
If the expressions of $\rho, \mu_+, \mu_-$ and $\omega$ in terms of
$x,y,z,t$ are introduced in
Eq. (\ref{ch42}), the logarithm of the polar fourcomplex number becomes
\begin{equation}
\ln u=\frac{1+\alpha+\beta+\gamma}{4}\ln(x+y+z+t)
+\frac{1-\alpha+\beta-\gamma}{4}\ln(x-y+z-t)
+e_1\ln\mu_+
+\tilde e_1\phi .
\label{ch44}
\end{equation}

The power function $u^m$ can be defined for $v_+>0, 
v_->0$ and real values of $m$ as
\begin{equation}
u^m=e^{m\ln u} .
\label{ch44b}
\end{equation}\index{power function!polar fourcomplex}
The power function is multivalued unless $m$ is an integer. 
For integer $m$, it can be inferred from Eqs. (\ref{ch42}) and (\ref{ch44b}) that
\begin{equation}
(u_1u_2)^m=u_1^m\:u_2^m .
\label{ch45}
\end{equation}

Using the expression (\ref{ch44}) for $\ln u$ and the relations
(\ref{ch22})-(\ref{ch23b}) it can be shown that 
\begin{eqnarray}
u^m=e_+v_+^m+e_-v_-^m
+\mu_+^m\left(e_1\cos m\phi
+\tilde e_1\sin m\phi\right).
\label{ch46}
\end{eqnarray}\index{power function!polar fourcomplex}
For integer $m$, the relation (\ref{ch46}) is valid for any $x,y,z,t$. 
For natural $m$ this relation can be written as
\begin{eqnarray}
u^m=e_+v_+^m+e_-v_-^m
+\left[e_1(x-z)
+\tilde e_1(y-t)\right]^m,
\label{ch46b}
\end{eqnarray}
as can be shown with the aid of the relation
\begin{eqnarray}
e_1\cos m\phi+\tilde e_1\sin m\phi
=\left(e_1\cos \phi+\tilde e_1\sin \phi\right)^m,
\label{ch46c}
\end{eqnarray}
valid for natural $m$.
For $m=-1$ the relation (\ref{ch46}) becomes
\begin{eqnarray}
\lefteqn{\frac{1}{x+\alpha y+\beta z+\gamma t}=
\frac{1}{4}\left(\frac{1+\alpha+\beta+\gamma}{x+y+z+t}
+\frac{1-\alpha+\beta-\gamma}{x-y+z-t}\right)\nonumber}\\
&&+\frac{1}{2}\frac{(1-\beta)(x-z)-(\alpha-\gamma)(y-t)}{(x-z)^2+(y-t)^2}.
\label{ch47}
\end{eqnarray}
If $m=2$, it can be checked that the right-hand side of
Eq. (\ref{ch46b}) is equal to $(x+\alpha y+\beta z+\gamma t)^2=
x^2+z^2+2yt+2\alpha(xy+zt)+\beta(y^2+t^2+2xz)+2\gamma(xt+yz)$.


The trigonometric functions of the hypercomplex variable
$u$ and the addition theorems for these functions have been written in Eqs.
~(1.57)-(1.60). 
The cosine and sine functions of the hypercomplex variables $\alpha y, 
\beta z$ and $ \gamma t$ can be expressed as
\begin{equation}
\cos\alpha y=g_{40}-\beta g_{42}, \: \sin\alpha y=\alpha g_{41}-\gamma g_{43}, 
\label{ch52}
\end{equation}\index{trigonometric functions, expressions!polar fourcomplex}
\begin{equation}
\cos\beta y=\cos y, \: \sin\beta y=\beta\sin y, 
\label{ch51}
\end{equation}
\begin{equation}
\cos\gamma y=g_{40}-\beta g_{42}, \: \sin\gamma y=\gamma g_{41}-\alpha g_{43} .
\label{ch53}
\end{equation}
The cosine and sine functions of a polar fourcomplex number $x+\alpha y+\beta
z+\gamma t$ can then be
expressed in terms of elementary functions with the aid of the addition
theorems Eqs. (1.59), (1.60) and of the expressions in  Eqs. 
(\ref{ch52})-(\ref{ch53}). 

The hyperbolic functions of the hypercomplex variable
$u$ and the addition theorems for these functions have been written in Eqs.
~(1.62)-(1.65). 
The hyperbolic cosine and sine functions of the hypercomplex variables $\alpha y, 
\beta z$ and $ \gamma t$ can be expressed as
\begin{equation}
\cosh\alpha y=g_{40}+\beta g_{42}, \: \sinh\alpha y=\alpha g_{41}+\gamma g_{43}, 
\label{ch59}
\end{equation}\index{hyperbolic functions, expressions!polar fourcomplex}
\begin{equation}
\cosh\beta y=\cosh y, \: \sinh\beta y=\beta\sinh y, 
\label{ch58}
\end{equation}
\begin{equation}
\cosh\gamma y=g_{40}+\beta g_{42}, \: \sinh\gamma y=\gamma g_{41}+\alpha g_{43} .
\label{ch60-63}
\end{equation}
The hyperbolic cosine and sine
functions of a polar fourcomplex number $x+\alpha y+\beta
z+\gamma t$ can then be
expressed in terms of elementary functions with the aid of the addition
theorems Eqs. (1.64), (1.65) and of the expressions in  Eqs. 
(\ref{ch59})-(\ref{ch60-63}). 

\subsection{Power series of polar fourcomplex variables}

A polar fourcomplex series is an infinite sum of the form
\begin{equation}
a_0+a_1+a_2+\cdots+a_l+\cdots , 
\label{ch64}
\end{equation}\index{series!polar fourcomplex}
where the coefficients $a_l$ are polar fourcomplex numbers. The convergence of 
the series (\ref{ch64}) can be defined in terms of the convergence of its 4 real
components. The convergence of a polar fourcomplex series can however be studied
using polar fourcomplex variables. The main criterion for absolute convergence 
remains the comparison theorem, but this requires a number of inequalities
which will be discussed further.

The modulus of a polar fourcomplex number $u=x+\alpha y+\beta z+\gamma t$ can be
defined as 
\begin{equation}
|u|=(x^2+y^2+z^2+t^2)^{1/2} ,
\label{ch65}
\end{equation}\index{modulus, definition!polar fourcomplex}
so that according to Eq. (\ref{ch12}) $d=|u|$. Since $|x|\leq |u|, |y|\leq |u|,
|z|\leq |u|, |t|\leq |u|$, a property of 
absolute convergence established via a comparison theorem based on the modulus
of the series (\ref{ch64}) will ensure the absolute convergence of each real
component of that series.

The modulus of the sum $u_1+u_2$ of the polar fourcomplex numbers $u_1, u_2$ fulfils
the inequality
\begin{equation}
||u_1|-|u_2||\leq |u_1+u_2|\leq |u_1|+|u_2| .
\label{ch66}
\end{equation}\index{modulus, inequalities!polar fourcomplex}
For the product the relation is 
\begin{equation}
|u_1u_2|\leq 2|u_1||u_2| ,
\label{ch67}
\end{equation}
which replaces the relation of equality extant for regular complex numbers.
The equality in Eq. (\ref{ch67}) takes place for 
$x_1=y_1=z_1=t_1, x_2=y_2=z_2=t_2$, or
$x_1=-y_1=z_1=-t_1, x_2=-y_2=z_2=-t_2$.
In particular
\begin{equation}
|u^2|\leq 2(x^2+y^2+z^2+t^2) .
\label{ch68}
\end{equation}
The inequality in Eq. (\ref{ch67}) implies that
\begin{equation}
|u^l|\leq 2^{l-1}|u|^l .
\label{ch69}
\end{equation}
From Eqs. (\ref{ch67}) and (\ref{ch69}) it results that
\begin{equation}
|au^l|\leq 2^l |a| |u|^l .
\label{ch70}
\end{equation}

A power series of the polar fourcomplex variable $u$ is a series of the form
\begin{equation}
a_0+a_1 u + a_2 u^2+\cdots +a_l u^l+\cdots .
\label{ch71}
\end{equation}\index{power series!polar fourcomplex}
Since
\begin{equation}
\left|\sum_{l=0}^\infty a_l u^l\right| \leq  \sum_{l=0}^\infty
2^l|a_l| |u|^l ,
\label{ch72}
\end{equation}
a sufficient condition for the absolute convergence of this series is that
\begin{equation}
\lim_{l\rightarrow \infty}\frac{2|a_{l+1}||u|}{|a_l|}<1 .
\label{ch73}
\end{equation}
Thus the series is absolutely convergent for 
\begin{equation}
|u|<c_0,
\label{ch74}
\end{equation}\index{convergence of power series!polar fourcomplex}
where 
\begin{equation}
c_0=\lim_{l\rightarrow\infty} \frac{|a_l|}{2|a_{l+1}|} .
\label{ch75}
\end{equation}

The convergence of the series (\ref{ch71}) can be also studied with the aid of
the formula (\ref{ch46b}) which, for integer values of $l$, is valid for any $x,
y, z, t$. If $a_l=a_{lx}+\alpha a_{ly}+\beta a_{lz}+\gamma a_{lt}$, and
\begin{eqnarray}
\begin{array}{l}
A_{l+}=a_{lx}+a_{ly}+a_{lz}+a_{lt}, \\
A_{l-}=a_{lx}-a_{ly}+a_{lz}-a_{lt},\\
A_{l1}=a_{lx}-a_{lz},\\
\tilde A_{l1}=a_{ly}-a_{lt},
\end{array}
\label{ch76}
\end{eqnarray}
it can be shown with the aid of relations (\ref{ch22})-(\ref{ch23b}) 
and (\ref{ch46b}) that the expression of the series (\ref{ch71}) is
\begin{equation}
\sum_{l=0}^\infty \left[A_{l+}  v_+^l e_++
A_{l-}v_-e_-+
\left(e_1A_{l1}+\tilde e_1\tilde A_{l1}\right)
\left(e_1v_1+\tilde e_1\tilde v_1\right)^l\right],
\label{ch77-78}
\end{equation}
where the quantities $v_+,  v_-$
have been defined in Eq. (\ref{ch8a}), and the quantities $v_1,\tilde v_1$ have
been defined in Eq. (\ref{ch16a}).

The sufficient conditions for the absolute convergence of the series 
in Eq. (\ref{ch77-78}) are that
\begin{equation}
\lim_{l\rightarrow \infty}\frac{|A_{l+1,+}||v_+|}{|A_{l+}|}<1,
\lim_{l\rightarrow \infty}\frac{|A_{l+1,-}||v_-|}{|A_{l-}|}<1,
\lim_{l\rightarrow \infty}\frac{A_{l+1}\mu_+}{A_l}<1,
\label{ch79}
\end{equation}
where the real and positive quantity $A_{l-}>0$ is given by
\begin{equation}
A_l^2=A_{l1}^2+\tilde A_{l1}^2.
\label{ch79a}
\end{equation}

Thus the series in Eq. (\ref{ch77-78}) is absolutely convergent for 
\begin{equation}
|x+y+z+t|<c_+,\:
|x-y+z-t|<c_-,\:
\mu_+<c_1
\label{ch80}
\end{equation}\index{convergence, region of!polar fourcomplex}
where 
\begin{equation}
c_+=\lim_{l\rightarrow\infty} \frac{|A_{l+}|}{|A_{l+1,+}|} ,\:
c_-=\lim_{l\rightarrow\infty} \frac{|A_{l-}|}{|A_{l+1,-}|} ,\:
c_1=\lim_{l\rightarrow\infty} \frac{A_{l-}}{A_{l+1,-}} .
\label{ch81-87}
\end{equation}
The relations (\ref{ch80}) show that the region of convergence of the series
(\ref{ch77-78}) is a four-dimensional cylinder.
It can be shown that $c_0=(1/2)\;{\rm min}(c_+,c_-,c_1)$, 
where ${\rm min}$ designates the smallest of
the numbers in the argument of this function. Using the expression of $|u|$ in
Eq. (\ref{ch24d}), it can be seen that the spherical region of
convergence defined in Eqs. (\ref{ch74}), (\ref{ch75}) is included in the
cylindrical region of convergence defined in Eqs. (\ref{ch81-87}).

\subsection{Analytic functions of polar fourcomplex variables}

The fourcomplex function $f(u)$ of the fourcomplex variable $u$ has
been expressed in Eq. (\ref{g16}) in terms of 
the real functions $P(x,y,z,t),Q(x,y,z,t),R(x,y,z,t), S(x,y,z,t)$ of real
variables $x,y,z,t$.
The
relations between the partial derivatives of the functions $P, Q, R, S$ are
obtained by setting succesively in   
Eq. (\ref{g17}) $\Delta x\rightarrow 0, \Delta y=\Delta z=\Delta t=0$;
then $ \Delta y\rightarrow 0, \Delta x=\Delta z=\Delta t=0;$  
then $  \Delta z\rightarrow 0,\Delta x=\Delta y=\Delta t=0$; and finally
$ \Delta t\rightarrow 0,\Delta x=\Delta y=\Delta z=0 $. 
The relations are 
\begin{equation}
\frac{\partial P}{\partial x} = \frac{\partial Q}{\partial y} =
\frac{\partial R}{\partial z} = \frac{\partial S}{\partial t},
\label{ch95}
\end{equation}\index{relations between partial derivatives!polar fourcomplex}
\begin{equation}
\frac{\partial Q}{\partial x} = \frac{\partial R}{\partial y} =
\frac{\partial S}{\partial z} = \frac{\partial P}{\partial t},
\label{ch97}
\end{equation}
\begin{equation}
\frac{\partial R}{\partial x} = \frac{\partial S}{\partial y} =
\frac{\partial P}{\partial z} = \frac{\partial Q}{\partial t},
\label{ch96}
\end{equation}
\begin{equation}
\frac{\partial S}{\partial x} =\frac{\partial P}{\partial y} =
\frac{\partial Q}{\partial z} =\frac{\partial R}{\partial t}.
\label{ch98}
\end{equation}


The relations (\ref{ch95})-(\ref{ch98}) are analogous to the Riemann relations
for the real and imaginary components of a complex function. It can be shown
from Eqs. (\ref{ch95})-(\ref{ch98}) that the component $P$ is a solution
of the equations 
\begin{equation}
\frac{\partial^2 P}{\partial x^2}-\frac{\partial^2 P}{\partial z^2}=0,
\:\: 
\frac{\partial^2 P}{\partial y^2}-\frac{\partial^2 P}{\partial t^2}=0,
\:\:
\label{ch99}
\end{equation}\index{relations between second-order derivatives!polar fourcomplex}
and the components $Q, R, S$ are solutions of similar equations.
As can be seen from Eqs. (\ref{ch99})-(\ref{ch99}), the components $P, Q, R, S$ of
an analytic function of polar fourcomplex variable are harmonic 
with respect to the pairs of variables $x,y$ and $ z,t$.
The component $P$ is also a solution of the mixed-derivative
equations 
\begin{equation}
\frac{\partial^2 P}{\partial x^2}=\frac{\partial^2 P}{\partial y\partial t},
\:\: 
\frac{\partial^2 P}{\partial y^2}=\frac{\partial^2 P}{\partial x\partial z},
\:\:
\frac{\partial^2 P}{\partial z^2}=\frac{\partial^2 P}{\partial y\partial t},
\:\:
\frac{\partial^2 P}{\partial t^2}=\frac{\partial^2 P}{\partial x\partial z},
\:\:
\label{ch106b}
\end{equation}
and the components $Q, R, S$ are solutions of similar equations.
The component $P$ is also a solution of the mixed-derivative
equations 
\begin{equation}
\frac{\partial^2 P}{\partial x\partial y}=\frac{\partial^2 P}{\partial
z\partial t} ,
\:\: 
\frac{\partial^2 P}{\partial x\partial t}=\frac{\partial^2 P}{\partial
y\partial z} ,
\label{ch107}
\end{equation}
and the components $Q, R, S$ are solutions of similar equations.

\subsection{Integrals of functions of polar fourcomplex variables}

The singularities of polar fourcomplex functions arise from terms of the form
$1/(u-u_0)^m$, with $m>0$. Functions containing such terms are singular not
only at $u=u_0$, but also at all points of the two-dimensional hyperplanes
passing through $u_0$ and which are parallel to the nodal hyperplanes. 

The integral of a polar fourcomplex function between two points $A, B$ along a path
situated in a region free of singularities is independent of path, which means
that the integral of an analytic function along a loop situated in a region
free from singularities is zero,
\begin{equation}
\oint_\Gamma f(u) du = 0,
\label{ch111}
\end{equation}
where it is supposed that a surface $\Sigma$ spanning 
the closed loop $\Gamma$ is not intersected by any of
the hyperplanes associated with the
singularities of the function $f(u)$. Using the expression, Eq. (\ref{g16})
for $f(u)$ and the fact that $du=dx+\alpha  dy+\beta dz+\gamma dt$, the
explicit form of the integral in Eq. (\ref{ch111}) is
\begin{eqnarray}
\lefteqn{\oint _\Gamma f(u) du = \oint_\Gamma
[(Pdx+Sdy+Rdz+Qdt)+\alpha(Qdx+Pdy+Sdz+Rdt)\nonumber}\\
&&+\beta(Rdx+Qdy+Pdz+Sdt)+\gamma(Sdx+Rdy+Qdz+Pdt)] .
\label{ch112}
\end{eqnarray}\index{integrals, path!polar fourcomplex}
If the functions $P, Q, R, S$ are regular on a surface $\Sigma$
spanning the loop $\Gamma$,
the integral along the loop $\Gamma$ can be transformed with the aid of the
theorem of Stokes in an integral over the surface $\Sigma$ of terms of the form
$\partial P/\partial y - \partial S/\partial x,\:\:
\partial P/\partial z -  \partial R/\partial x, \:\:
\partial P/\partial t - \partial Q/\partial x, \:\:
\partial R/\partial y -  \partial S/\partial z, \:\:
\partial S/\partial t - \partial Q/\partial y, \:\:
\partial R/\partial t - \partial Q/\partial z$
and of similar terms arising
from the $\alpha, \beta$ and $\gamma$ components, 
which are equal to zero by Eqs. (\ref{ch95})-(\ref{ch98}), and this proves Eq.
(\ref{ch111}). 

The integral of the function $(u-u_0)^m$ on a closed loop $\Gamma$ is equal to
zero for $m$ a positive or negative integer not equal to -1,
\begin{equation}
\oint_\Gamma (u-u_0)^m du = 0, \:\: m \:\:{\rm integer},\: m\not=-1 .
\label{ch112b}
\end{equation}
This is due to the fact that $\int (u-u_0)^m du=(u-u_0)^{m+1}/(m+1), $ and to
the fact that the function $(u-u_0)^{m+1}$ is singlevalued for $m$ an integer.

The integral $\oint_\Gamma du/(u-u_0)$ can be calculated using the exponential
form (\ref{ch40c}),
\begin{eqnarray}
u-u_0=
\rho\exp\left[\frac{1}{4}(\alpha+\beta+\gamma) \ln\frac{\sqrt{2}}{\tan\theta_+}
-\frac{1}{4}(\alpha-\beta+\gamma) \ln\frac{\sqrt{2}}{\tan\theta_-}
+\tilde e_1\phi\right] ,
\label{ch113}
\end{eqnarray}
so that 
\begin{equation}
\frac{du}{u-u_0}=\frac{d\rho}{\rho}
+\frac{1}{4}(\alpha+\beta+\gamma) d\ln\frac{\sqrt{2}}{\tan\theta_+}
-\frac{1}{4}(\alpha-\beta+\gamma) d\ln\frac{\sqrt{2}}{\tan\theta_-}
+\tilde e_1d\phi.
\label{ch114}
\end{equation}
Since $\rho$, $\tan\theta_+$ and $\tan\theta_-$ are singlevalued variables,
it follows that 
$\oint_\Gamma d\rho/\rho =0, \oint_\Gamma d\ln\sqrt{2}/\tan\theta_+=0$,
and $\oint_\Gamma d\ln\sqrt{2}/\tan\theta_+=0$. On the other hand,
$\phi$ is a cyclic variables, so that it may give a contribution to
the integral around the closed loop $\Gamma$.
The result of the integrations will be given in the rotated system of
coordinates 
\begin{equation}
\xi=\frac{1}{\sqrt{2}}(x-z),\:
\upsilon=\frac{1}{\sqrt{2}}(y-t),\:
\tau=\frac{1}{2}(x+y+z+t),\:
\upsilon=\frac{1}{2}(x-y+z-t) .
\label{ch115}
\end{equation}
Thus, if $C_\parallel$ is a circle of radius $r$
parallel to the $\xi O\upsilon$ plane, and the
projection of the center of this circle on the $\xi O\upsilon$ plane
coincides with the projection of the point $u_0$ on this plane, the points
of the circle $C_\parallel$ are described according to Eqs.
(\ref{ch16a})-(\ref{ch17}) by the equations
\begin{equation}
\xi=\xi_0+r \cos\phi , \:
\upsilon=\upsilon_0+r \sin\phi , \:
\tau=\tau_0,\:\zeta=\zeta_0,
\label{ch115a}
\end{equation}
where $u_0=x_0+\alpha y_0+\beta 
z_0+\gamma t_0$, and $\xi_0, \upsilon_0, \tau_0, \zeta_0$ are calculated
from $x_0, y_0, z_0, t_0$ according to Eqs. (\ref{ch115}).

Then
\begin{equation}
\oint_{C_\parallel}\frac{du}{u-u_0}=2\pi\tilde e_1. 
\label{ch116}
\end{equation}
The expression of $\oint_\Gamma du/(u-u_0)$ can be written 
with the aid of the functional int($M,C$) defined in Eq. (\ref{118}) as
\begin{equation}
\oint_\Gamma\frac{du}{u-u_0}=
2\pi\tilde e_1 \;{\rm int}(u_{0\xi\upsilon},\Gamma_{\xi\upsilon}) ,
\label{ch119}
\end{equation}\index{poles and residues!polar fourcomplex}
where $u_{0\xi\upsilon}$ and $\Gamma_{\xi\upsilon}$ are respectively the
projections of the point $u_0$ and of 
the loop $\Gamma$ on the plane $\xi \upsilon$.

If $f(u)$ is an analytic polar fourcomplex function which can be expanded in a
series as written in Eq. (1.89), and the expansion holds on the curve
$\Gamma$ and on a surface spanning $\Gamma$, then from Eqs. (\ref{ch112b}) and
(\ref{ch119}) it follows that
\begin{equation}
\oint_\Gamma \frac{f(u)du}{u-u_0}=
2\pi\tilde e_1 \;{\rm int}(u_{0\xi\upsilon},\Gamma_{\xi\upsilon})f(u_0) ,
\label{ch120}
\end{equation}
where $\Gamma_{\xi\upsilon}$ is the projection of 
the curve $\Gamma$ on the plane $\xi \upsilon$,
as shown in Fig. \ref{fig14}.
Substituting in the right-hand side of 
Eq. (\ref{ch120}) the expression of $f(u)$ in terms of the real 
components $P, Q, R, S$, Eq. (\ref{g16}), yields
\begin{equation}
\oint_\Gamma \frac{f(u)du}{u-u_0}=\pi 
\left[(\beta-1)(Q-S)+(\alpha-\gamma)(P-R)\right] 
{\rm int}\left(u_{0\xi\upsilon},\Gamma_{\xi\upsilon}\right),
\label{ch121}
\end{equation}
where $P, Q, R, S$ are the values of the components of $f$ at $u=u_0$.
If the integral is written as
\begin{equation}
\oint_\Gamma \frac{f(u)du}{u-u_0}=
I+\alpha I_\alpha+\beta I_\beta+\gamma I_\gamma,
\label{ch121b}
\end{equation}
it results from Eq. (\ref{ch121}) that
\begin{equation}
I+ I_\alpha+ I_\beta+ I_\gamma=0.
\label{ch121c}
\end{equation}

\begin{figure}
\begin{center}
\epsfig{file=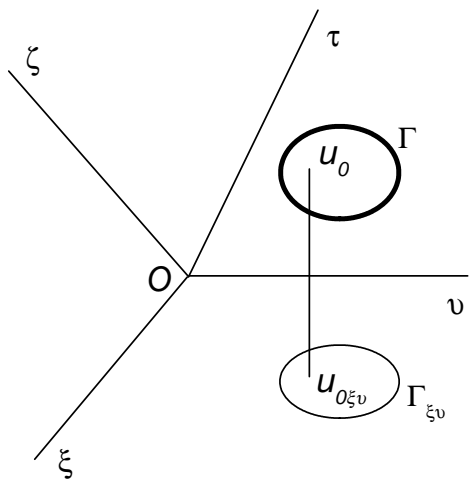,width=12cm}
\caption{Integration path $\Gamma$ and the pole $u_0$, and their projections
$\Gamma_{\xi\upsilon}$ and $u_{0\xi\upsilon}$ on the plane $\xi \upsilon$. }
\label{fig14}
\end{center}
\end{figure}

If $f(u)$ can be expanded as written in Eq. (1.89) on 
$\Gamma$ and on a surface spanning $\Gamma$, then from Eqs. (\ref{ch112b}) and
(\ref{ch119}) it also results that
\begin{equation}
\oint_\Gamma \frac{f(u)du}{(u-u_0)^{m+1}}=
2\frac{\pi}{m!}\tilde e_1 \;{\rm int}(u_{0\xi\upsilon},\Gamma_{\xi\upsilon}) 
\;f^{(m)}(u_0) ,
\label{ch122}
\end{equation}
where the fact has been used  that the derivative $f^{(m)}(u_0)$ of order $m$
of $f(u)$ at $u=u_0$ is related to the expansion coefficient in Eq. (1.89)
according to Eq. (1.93).

If a function $f(u)$ is expanded in positive and negative powers of $u-u_j$,
where $u_j$ are polar fourcomplex constants, $j$ being an index, the integral of $f$
on a closed loop $\Gamma$ is determined by the terms in the expansion of $f$
which are of the form $a_j/(u-u_j)$,
\begin{equation}
f(u)=\cdots+\sum_j\frac{a_j}{u-u_j}+\cdots
\label{ch123}
\end{equation}
Then the integral of $f$ on a closed loop $\Gamma$ is
\begin{equation}
\oint_\Gamma f(u) du = 
2\pi\tilde e_1 \sum_j{\rm int}(u_{j\xi\upsilon},\Gamma_{\xi\upsilon})a_j .
\label{ch124}
\end{equation}

\newpage
\setlength{\oddsidemargin}{0.4cm}      

\subsection{Factorization of polar fourcomplex polynomials}

A polynomial of degree $m$ of the polar fourcomplex variable 
$u=x+\alpha y+\beta z+\gamma t$ has the form
\begin{equation}
P_m(u)=u^m+a_1 u^{m-1}+\cdots+a_{m-1} u +a_m ,
\label{ch125}
\end{equation}
where the constants are in general polar fourcomplex numbers.

If $a_m=a_{mx}+\alpha a_{my}+\beta a_{mz}+\gamma a_{mt}$, and with the
notations of Eqs. (\ref{ch8a}) and (\ref{ch76}) applied for $0, 1, \cdots, m$ , the
polynomial $P_m(u)$ can be written as 
\begin{eqnarray}
\lefteqn{P_m= \left[v_+^m 
+A_1 v_+^{m-1}+\cdots+A_{m-1} v_++ A_m \right] e_+
+\left[v_-^m 
+A_1^{\prime\prime} v_-^{m-1} +\cdots+A_{m-1}^{\prime\prime} 
v_-+ A_m^{\prime\prime} \right]e_-\nonumber}\\
 & &+\left[
\left(e_1v_1+\tilde e_1\tilde v_1\right)^{m}
+\sum_{l=1}^m\left(e_1A_{l1}+\tilde e_1\tilde A_{l1}\right)
\left(e_1v_1+\tilde e_1\tilde v_1\right)^{m-l}\right],\nonumber\\
&&
\label{ch126}
\end{eqnarray}\index{polynomial, canonical variables!polar fourcomplex}
where the constants $A_{l+}, A_{l-}, A_{l1}, \tilde A_{l1}$ 
are real numbers.
The polynomial of degree $m$ in $(e_1v_1+\tilde e_1\tilde v_1)$
can always be written as a product of linear factors of the form
$[e_1(v_1-v_{1l})+\tilde e_1 (\tilde v_1-\tilde v_{1l})]$, where the
constants $v_{1l}, \tilde v_{1l}$ are real.
The two polynomials of degree $m$ with real coefficients in Eq. (\ref{ch126})
which are multiplied by $e_+$ and $e_-$ can be written as a product
of linear or quadratic factors with real coefficients, or as a product of
linear factors which, if imaginary, appear always in complex conjugate pairs.
Using the latter form for the simplicity of notations, the polynomial $P_m$
can be written as
\begin{equation}
P_m=\prod_{l=1}^m (v_+-s_{l+})e_+
+\prod_{l=1}^m (v_--s_{l-})e_-
+\prod_{l=1}^m
\left[e_1(v_1-v_{1l})+
\tilde e_1(\tilde v_1-\tilde v_{1l})\right],
\label{ch127}
\end{equation}
where the quantities $s_{l+}$ appear always in complex conjugate pairs, and the
same is true for the quantities $s_{l-}$.
Due to the properties in Eqs. (\ref{ch22})-(\ref{ch23b}),
the polynomial $P_m(u)$ can be 
\newpage
\setlength{\oddsidemargin}{0.9cm}      
written as a product of factors of
the form  
\begin{equation}
P_m(u)=\prod_{l=1}^m \left[(v_+-s_{l+})e_+
+\left(v_--s_{l-}\right)e_-
+\left(e_1(v_1-v_{1l})+
\tilde e_1(\tilde v_1-\tilde v_{1l})\right)
\right].
\label{ch128}
\end{equation}\index{polynomial, factorization!polar fourcomplex}
These relations can be written with the aid of Eq. (\ref{ch24}) as
\begin{eqnarray}
P_m(u)=\prod_{p=1}^m (u-u_p) ,
\label{ch128c}
\end{eqnarray}
where
\begin{eqnarray}
u_p=s_{p+} e_+
+s_{p-}e_-
+e_1v_{1p}+
\tilde e_1\tilde v_{1p}, p=1,...,m.
\label{ch128d}
\end{eqnarray}
The roots $s_{p+},
s_{p-}, v_{1p}e_1+\tilde v_{1p}
\tilde e_1$ of the corresponding polynomials in Eq.
(\ref{ch127}) may be ordered arbitrarily.
This means that Eq. (\ref{ch128d}) gives sets of $m$ roots
$u_1,...,u_m$ of the polynomial $P_m(u)$, 
corresponding to the various ways in which the roots
$s_{p+}, s_{p-},v_{1p} e_1+\tilde v_{1p}
\tilde e_1$ are ordered according to $p$ in each
group. Thus, while the hypercomplex components in Eq. (\ref{ch126}) taken
separately have unique factorizations, the polynomial $P_m(u)$ can be written
in many different ways as a product of linear factors. 
The result of the polar fourcomplex integration, Eq. (\ref{ch124}), is however
unique.

If $P(u)=u^2-1$, 
the factorization in Eq. (\ref{ch128c}) is $u^2-1=(u-u_1)(u-u_2)$, where 
$u_1=\pm e_+\pm e_-\pm e_1, u_2=-u_1$, so that
there are 4 distinct factorizations of $u^2-1$,
\begin{eqnarray}
\begin{array}{l}
u^2-1=(u-1)(u+1),\\
u^2-1=(u-\beta)(u+\beta),\\
u^2-1=\left(u-\frac{1+\alpha-\beta+\gamma}{2}\right)
\left(u+\frac{1+\alpha-\beta+\gamma}{2}\right),\\
u^2-1=\left(u-\frac{-1+\alpha+\beta+\gamma}{2}\right)
\left(u+\frac{-1+\alpha+\beta+\gamma}{2}\right).
\end{array}
\label{ch129}
\end{eqnarray}

It can be checked that 
$\left\{\pm e_+\pm e_-\pm e_1\right\}^2= 
e_++e_-+e_1=1$.

\subsection{Representation of polar fourcomplex numbers 
by irreducible matrices}

If $T$ is the unitary matrix,
\begin{equation}
T =\left(
\begin{array}{cccc}
\frac{1}{{2}}&\frac{1}{{2}}    &\frac{1}{{2}}    &\frac{1}{{2}}    \\
\frac{1}{{2}}  &-\frac{1}{{2}} &\frac{1}{{2}}    &-\frac{1}{{2}}   \\
\frac{1}{\sqrt{2}}& 0          & -\frac{1}{\sqrt{2}}  & 0               \\
0                 & \frac{1}{\sqrt{2}} & 0       & -\frac{1}{\sqrt{2}}  \\
\end{array}
\right),
\label{ch129x}
\end{equation}
it can be shown 
that the matrix $T U T^{-1}$ has the form 
\begin{equation}
T U T^{-1}=\left(
\begin{array}{ccc}
x+y+z+t &    0    &  0    \\
0       & x-y+z-t &  0    \\
0       &    0    &  V_1  \\
\end{array}
\right),
\label{ch129y}
\end{equation}\index{representation by irreducible matrices!polar fourcomplex}
where $U$ is the matrix in Eq. (\ref{ch23}) used to represent the polar fourcomplex
number $u$. In Eq. (\ref{ch129y}), $V_1$ is the matrix
\begin{equation}
V_1=\left(
\begin{array}{cc}
x-z    &   y-t   \\
-y+t   &   x-z   \\
\end{array}\right).
\label{ch130x}
\end{equation}
The relations between the variables $x+y+z+t, x-y+z-t, x-z, y-t$
for the multiplication
of polar fourcomplex numbers have been written in Eqs. (\ref{ch15}),(\ref{ch16}),
(\ref{ch19a}), (\ref{ch19b}). The matrix 
$T U T^{-1}$ provides an irreducible representation
\cite{4} of the polar fourcomplex number $u$ in terms of matrices with real
coefficients. 

\chapter{Complex Numbers in 5 Dimensions}

A system of complex numbers in 5 dimensions is described in this chapter,
for which the multiplication is associative and commutative, which  
have exponential and trigonometric forms, and for which
the concepts of analytic 5-complex 
function,  contour integration and residue can be defined.
The 5-complex numbers introduced in this chapter have 
the form $u=x_0+h_1x_1+h_2x_2+h_3x_3+h_4x_4$, the variables 
$x_0,x_1,x_2,x_3,x_4$ being real numbers. 
If the 5-complex number $u$ is
represented by the point $A$ of coordinates $x_0,x_1,x_2,x_3,x_4$, 
the position of the point $A$ can be described
by the modulus $d=(x_0^2+x_1^2+x_2^2+x_3^2+x_4^2)^{1/2}$, 
by 2 azimuthal angles $\phi_1, \phi_2$, by 1 planar angle $\psi_1$,
and by 1 polar angle $\theta_+$. 

The exponential function of a 5-complex number can be expanded in terms of
the polar 5-dimensional cosexponential functions
$g_{5k}(y)$, $k=0,1,2,3,4$. The expressions of these functions are obtained
from the properties of the exponential function of a 5-complex variable.
Addition theorems and other relations
are obtained for the polar 5-dimensional cosexponential functions.
Exponential and trigonometric forms are given for the 5-complex numbers.
Expressions are obtained for the elementary functions of 5-complex variable.
The functions $f(u)$ of 5-complex variable which are defined by power series
have derivatives independent of the direction of approach to the point under
consideration. 
If the 5-complex function $f(u)$ 
of the 5-complex variable $u$ is written in terms of 
the real functions $P_k(x_0,x_1,x_2,x_3,x_4), k=0,1,2,3,4$, then
relations of equality  
exist between partial derivatives of the functions $P_k$. 
The integral $\int_A^B f(u) du$ of a 5-complex
function between two points $A,B$ is independent of the path connecting $A,B$,
in regions where $f$ is regular.
The fact that the exponential form
of the 5-complex numbers depends on the cyclic variables $\phi_1,\phi_2$
leads to the 
concept of pole and residue for integrals on closed paths,
and if $f(u)$ is an analytic 5-complex function, then $\oint_\Gamma
f(u)du/(u-u_0)$ is expressed in this chapter in terms of the 5-complex residue
$f(u_0)$. The polynomials of
5-complex variables can be written as products of linear or quadratic
factors. 

The 5-complex numbers described in this chapter are a particular case for 
$n=5$ of the polar complex numbers in $n$ dimensions discussed in Sec. 6.1.

\section{Operations with polar complex numbers in 5 dimensions}

A polar hypercomplex number $u$ in 5 dimensions 
is represented as  
\begin{equation}
u=x_0+h_1x_1+h_2x_2+h_3x_3+h_4x_4. 
\label{5-1a}
\end{equation}
The multiplication rules for the bases 
$h_1, h_2, h_3, h_4$ are 
\begin{eqnarray}
\lefteqn{h_1^2=h_2,\;h_2^2=h_4,\; h_3^2=h_1,\; h_4^2=h_3,\nonumber}\\
&& h_1h_2=h_3, \;h_1h_3=h_4,\;
h_1h_4=1,\; h_2h_3=1,\; h_2h_4=h_1, h_3h_4=h_2.
\label{5-1}
\end{eqnarray}\index{complex units!polar 5-complex}
The significance of the composition laws in Eq.
(\ref{5-1}) can be understood by representing the bases $h_j, h_k$ by points on a
circle at the angles $\alpha_j=2\pi j/5,\alpha_k=2\pi k/5$, as shown in Fig. \ref{fig15},
and the product $h_j h_k$ by the point of the circle at the angle 
$2\pi (j+k)/5$. If $2\pi\leq 2\pi (j+k)/5<4\pi$, the point represents the basis
$h_l$ of angle $\alpha_l=2\pi(j+k)/5-2\pi$.

\begin{figure}
\begin{center}
\epsfig{file=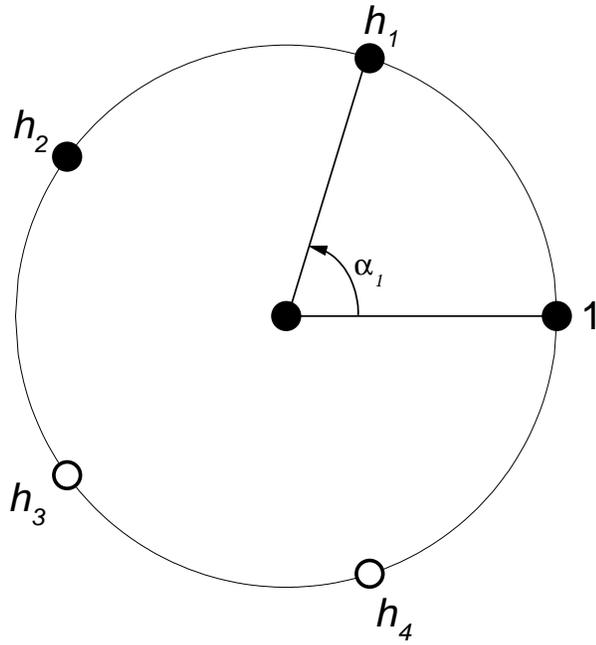,width=12cm}
\caption{Representation of the polar 
hypercomplex bases $1,h_1,h_2,h_3,h_4$
by points on a circle at the angles $\alpha_k=2\pi k/5$.
The product $h_j h_k$ will be represented by the point of the circle at the
angle $2\pi (j+k)/5$, $i,k=0,1,...,4$, where $h_0=1$. If $2\pi\leq 
2\pi (j+k)/5\leq 4\pi$, the point represents the basis
$h_l$ of angle $\alpha_l=2\pi(j+k)/5-2\pi$.}
\label{fig15}
\end{center}
\end{figure}

The sum of the 5-complex numbers $u$ and
$u^\prime$ is
\begin{equation}
u+u^\prime=x_0+x^\prime_0+h_1(x_1+x^\prime_1)
+h_2(x_2+x^\prime_2)+h_3(x_3+x^\prime_3)+h_4(x_4+x^\prime_4).
\label{5-2}
\end{equation}\index{sum!polar 5-complex}
The product of the numbers $u, u^\prime$ is then
\begin{equation}
\begin{array}{l}
uu^\prime=x_0 x_0^\prime +x_1x_4^\prime+x_2 x_3^\prime+x_3x_2^\prime
+x_4x_1^\prime\\
+h_1(x_0 x_1^\prime+x_1x_0^\prime+x_2x_4^\prime+x_3x_3^\prime
+x_4 x_2^\prime) \\
+h_2(x_0 x_2^\prime+x_1x_1^\prime+x_2x_0^\prime+x_3x_4^\prime
+x_4 x_3^\prime) \\
+h_3(x_0 x_3^\prime+x_1x_2^\prime+x_2x_1^\prime+x_3x_0^\prime
+x_4 x_4^\prime) \\
+h_4(x_0 x_4^\prime+x_1x_3^\prime+x_2x_2^\prime+x_3x_1^\prime
+x_4 x_0^\prime). \\
\end{array}
\label{5-3}
\end{equation}\index{product!polar 5-complex}

The relation between the variables
$v_+, v_1, \tilde v_1, v_2, \tilde v_2$ and $x_0, x_1,
x_2, x_3, x_4$ can be written with the aid of the parameters $p=(\sqrt{5}-1)/4,
q=\sqrt{(5+\sqrt{5})/8}$ as
\begin{equation}
\left(
\begin{array}{c}
v_+\\
v_1\\
\tilde v_1\\
v_2\\
\tilde v_2\\
\end{array}\right)
=\left(
\begin{array}{ccccc}
1&1&1&1&1\\
1&p&2p^2-1&2p^2-1&p\\
0&q&2pq&-2pq&-q\\
1&2p^2-1&p&p&2p^2-1\\
0&2pq&-q&q&-2pq\\
\end{array}
\right)
\left(
\begin{array}{c}
x_0\\
x_1\\
x_2\\
x_3\\
x_4
\end{array}
\right).
\label{5-9e}
\end{equation}\index{canonical variables!polar 5-complex}
The other variables are $v_3=v_2, \tilde v_3=-\tilde v_2, v_4=v_1, 
\tilde v_4=-\tilde v_1$. 
The variables $v_+, v_1, \tilde v_1, v_2, \tilde v_2$ will be called canonical
5-complex variables.

\section{Geometric representation of polar complex numbers in 5
dimensions}

The 5-complex number $x_0+h_1x_1+h_2x_2+h_3x_3+h_4x_4$
can be represented by 
the point $A$ of coordinates $(x_0,x_1,x_2, x_3, x_4)$. 
If $O$ is the origin of the 5-dimensional space,  the
distance from the origin $O$ to the point $A$ of coordinates
$(x_0,x_1,x_2, x_3, x_4)$ has the expression
\begin{equation}
d^2=x_0^2+x_1^2+x_2^2+x_3^2+x_4^2.
\label{5-10}
\end{equation}\index{distance!polar 5-complex}
The quantity $d$ will be called modulus of the 5-complex number 
$u$. The modulus of a 5-complex number $u$ will be designated by $d=|u|$.
The modulus has the property that
\begin{equation}
|u^\prime u^{\prime\prime}|\leq \sqrt{5}|u^\prime||u^{\prime\prime}| .
\label{5-79}
\end{equation}\index{modulus, inequalities!polar 5-complex}

The exponential and trigonometric forms of the 5-complex number $u$ can be
obtained conveniently in a rotated system of axes defined by the transformation
\begin{equation}
\left(
\begin{array}{c}
\xi_+\\
\xi_1\\
\eta_1\\
\xi_2\\
\eta_2\\
\end{array}\right)
=\sqrt{\frac{2}{5}}\left(
\begin{array}{ccccc}
\frac{1}{\sqrt{2}}&\frac{1}{\sqrt{2}}&\frac{1}{\sqrt{2}}&\frac{1}{\sqrt{2}}&\frac{1}{\sqrt{2}}\\
1&p&2p^2-1&2p^2-1&p\\
0&q&2pq&-2pq&-q\\
1&2p^2-1&p&p&2p^2-1\\
0&2pq&-q&q&-2pq\\
\end{array}
\right)
\left(
\begin{array}{c}
x_0\\
x_1\\
x_2\\
x_3\\
x_4
\end{array}
\right).
\label{5-12}
\end{equation}
The lines of the matrices in Eq. (\ref{5-12}) gives the components
of the 5 basis vectors of the new system of axes. These vectors have unit
length and are orthogonal to each other.
The relations between the two sets of variables are
\begin{equation}
v_+= \sqrt{5}\xi_+ ,  
v_k= \sqrt{\frac{5}{2}}\xi_k , \tilde v_k= \sqrt{\frac{5}{2}}\eta_k, k=1,2 .
\label{5-12b}
\end{equation}

The radius $\rho_k$ and the azimuthal angle $\phi_k$ in the plane of the axes
$v_k,\tilde v_k$  are
\begin{eqnarray}
\rho_k^2=v_k^2+\tilde v_k^2, \:\cos\phi_k=v_k/\rho_k,
\:\sin\phi_k=\tilde v_k/\rho_k, 
\label{5-19a}
\end{eqnarray}
$0\leq \phi_k<2\pi,\; k=1,2$,
so that there are 2 azimuthal angles.\index{azimuthal angles!polar 5-complex}
The planar angle $\psi_1$ is
\begin{equation}
\tan\psi_1=\rho_1/\rho_2, 
\label{5-19b}
\end{equation}
where $0\leq\psi_1\leq\pi/2$.\index{planar angle!polar 5-complex}
There is a polar angle $\theta_+$, 
\begin{equation}
\tan\theta_+=\frac{\sqrt{2}\rho_1}{v_+}, 
\label{5-19c}
\end{equation}\index{polar angle!polar 5-complex}
where $0\leq\theta_+\leq\pi$.
It can be checked that
\begin{equation}
\frac{1}{5}v_+^2
+\frac{2}{5}(\rho_1^2+\rho_2^2)=d^2 .
\label{5-18}
\end{equation}
The amplitude of a 5-complex number $u$ is
\begin{equation}
\rho=\left(v_+\rho_1^2 \rho_2^2\right)^{1/5}.
\label{5-50bb}
\end{equation}\index{amplitude!polar 5-complex}

If $u=u^\prime u^{\prime\prime}$, the parameters of the hypercomplex numbers
are related by
\begin{equation}
v_+=v_+^\prime v_+^{\prime\prime},  
\label{5-21a}
\end{equation}\index{transformation of variables!polar 5-complex}
\begin{equation}
\rho_k=\rho_k^\prime\rho_k^{\prime\prime}, 
\label{5-21b}
\end{equation}
\begin{equation}
\tan\theta_+=\frac{1}{\sqrt{2}}\tan\theta_+^\prime \tan\theta_+^{\prime\prime},
\label{5-21c}
\end{equation}
\begin{equation}
\tan\psi_1=\tan\psi_1^\prime \tan\psi_1^{\prime\prime},  
\label{5-21d}
\end{equation}
\begin{equation}
\phi_k=\phi_k^\prime+\phi_k^{\prime\prime},
\label{5-21e}
\end{equation}
\begin{equation}
v_k=v_k^\prime v_k^{\prime\prime}-\tilde v_k^\prime \tilde v_k^{\prime\prime},\;
\tilde v_k=v_k^\prime \tilde v_k^{\prime\prime}+\tilde v_k^\prime v_k^{\prime\prime},
\label{5-22}
\end{equation}
\begin{equation}
\rho=\rho^\prime\rho^{\prime\prime} ,
\label{5-24}
\end{equation}
where $k=1,2$.

The 5-complex number
$u=x_0+h_1x_1+h_2x_2+h_3x_3+h_4x_4$ can be  represented by the matrix
\begin{equation}
U=\left(
\begin{array}{ccccc}
x_0     &   x_1     &   x_2   &   x_3  &  x_4\\
x_4     &   x_0     &   x_1   &   x_2  &  x_3\\
x_3     &   x_4     &   x_0   &   x_1  &  x_2\\
x_2     &   x_3     &   x_4   &   x_0  &  x_1\\
x_1     &   x_2     &   x_3   &   x_4  &  x_0\\
\end{array}
\right).
\label{5-24b}
\end{equation}\index{matrix representation!polar 5-complex}
The product $u=u^\prime u^{\prime\prime}$ is
represented by the matrix multiplication $U=U^\prime U^{\prime\prime}$.

\section{The polar 5-dimensional cosexponential functions}

The polar cosexponential functions in 5 dimensions are
\begin{equation}
g_{5k}(y)=\sum_{p=0}^\infty y^{k+5p}/(k+5p)!, 
\label{5-29}
\end{equation}\index{cosexponential function, polar 5-complex!definition}
for $k=0,...,4$.
The polar cosexponential functions $g_{5k}$ do not have a definite
parity. 
It can be checked that
\begin{equation}
\sum_{k=0}^{4}g_{5k}(y)=e^y.
\label{5-29a}
\end{equation}
The exponential of the quantity $h_k y, k=1,...,4$ can be written as
\begin{equation}
\begin{array}{l}
e^{h_1 y}=g_{50}(y)+h_1g_{51}(y)+h_2g_{52}(y)+h_3g_{53}(y)+h_4g_{54}(y),\\
e^{h_2 y}=g_{50}(y)+h_1g_{53}(y)+h_2g_{51}(y)+h_3g_{54}(y)+h_4g_{52}(y),\\
e^{h_3 y}=g_{50}(y)+h_1g_{52}(y)+h_2g_{54}(y)+h_3g_{51}(y)+h_4g_{53}(y),\\
e^{h_4 y}=g_{50}(y)+h_1g_{54}(y)+h_2g_{53}(y)+h_3g_{52}(y)+h_4g_{51}(y).
\end{array}
\label{5-28b}
\end{equation}\index{exponential, expressions!polar 5-complex}

The polar cosexponential functions in 5 dimensions can be obtained by
calculating $e^{(h_1+h_4)y}$ and $e^{(h_1-h_4)y}$ and then by
nultiplying the resulting expression.
The series expansions for $e^{(h_1+h_4)y}$ and $e^{(h_1-h_4)y}$ are
\begin{equation}
e^{(h_1+h_4)y}=\sum_{m=0}^\infty \frac{1}{m!}(h_1+h_4)^m y^m ,
\label{5-g5-1}
\end{equation}\index{cosexponential functions, polar 5-complex!expressions}
\begin{equation}
e^{(h_1-h_4)y}=\sum_{m=0}^\infty \frac{1}{m!}(h_1-h_4)^m y^m .
\label{5-g5-2}
\end{equation}
The powers of $h_1+h_4$ have the form
\begin{equation}
(h_1+h_4)^m=A_m (h_1+h_4)+B_m (h_2+h_3)+C_m.
\label{5-g5-3}
\end{equation}
The recurrence relations for $A_m, B_m, C_m$ are
\begin{equation}
A_{m+1}=B_m+C_m, B_{m+1}=A_m+B_m, C_{m+1}=2A_m,
\label{5-g5-4}
\end{equation}
and $A_1=1, B_1=0, C_1=0, A_2=0, B_2=1, C_2=2, A_3=3, B_3=1, C_3=0$.
The expressions of the coefficients are
\begin{equation}
A_m=\frac{2^m}{5}+\frac{2-3a}{5}a^{m-3}+(-1)^{m-3}\frac{5+3a}{5}(1+a)^{m-3},
m\geq 3 ,
\label{5-g5-5}
\end{equation}
\begin{equation}
B_m=\frac{2^m}{5}+\frac{a-1}{5}a^{m-3}-(-1)^{m-3}\frac{a+2}{5}(1+a)^{m-3},
m\geq 3,
\label{5-g5-6}
\end{equation}
\begin{equation}
C_m=\frac{2^m}{5}+\frac{4-6a}{5}a^{m-4}+(-1)^{m-4}\frac{10+6a}{5}(1+a)^{m-4} ,
m\geq 4,
\label{5-g5-7}
\end{equation}
where $a$ is a solution of the equation 
$a^2+a-1=0$.
Substituting the expressions of $A_m, B_m, C_m$ from Eqs. (\ref{5-g5-5})-(\ref{5-g5-7})
in Eq. (\ref{5-g5-1}) and grouping the terms yields
\begin{eqnarray}
e^{(h_1+h_4)y}
\lefteqn{=\frac{1}{5}e^{2y}+\frac{2}{5}e^{ay}+\frac{2}{5}e^{-(1+a)y}
+(h_1+h_4)
\left[\frac{1}{5}e^{2y}+\frac{a}{5}e^{ay}-\frac{a+1}{5}e^{-(1+a)y}\right]\nonumber}\\
&&+(h_2+h_3)
\left[\frac{1}{5}e^{2y}-\frac{a+1}{5}e^{ay}+\frac{a}{5}e^{-(1+a)y}\right].
\label{5-g5-8}
\end{eqnarray}

The odd powers of $h_1-h_4$ have the form
\begin{equation}
(h_1-h_4)^{2m+1}=D_m (h_1-h_4)+E_m (h_2-h_3).
\label{5-g5-9}
\end{equation}
The recurrence relations for $D_m, E_m$ are
\begin{equation}
D_{m+1}=-3D_m-E_m, E_{m+1}=-D_m-2E_m,
\label{5-g5-10}
\end{equation}
and $D_1=-3, E_1=-1, D_2=10, E_2=5$.
The expressions of the coefficients are
\begin{equation}
D_m=(b+1)b^{m-1}+(-1)^{m-2}(b+4)(5+b)^{m-1}, m\geq 1,
\label{5-g5-11}
\end{equation}
\begin{equation}
E_m=-\frac{b+1}{b+2}b^{m-1}+\frac{(-1)^{m-2}}{b+2}(5+b)^{m-1}, 
m\geq 1,
\label{5-g5-12}
\end{equation}
where $b$ is a solution of the equation 
$b^2+5b+5=0$.
The even powers of $h_1-h_4$ have the form
\begin{equation}
(h_1-h_4)^{2m}=F_m (h_1+h_4)+G_m (h_2+h_3)+H_m.
\label{5-g5-13}
\end{equation}
The recurrence relations for $F_m, G_m, H_m$ are
\begin{equation}
F_{m+1}=-F_m+G_m, G_{m+1}=F_m-2G_m+H_m, H_{m+1}=2(G_m-H_m),
\label{5-g5-14}
\end{equation}
and $F_1=0, G_1=1, H_1=-2, F_2=1, G_2=-4, H_2=6$.
The expressions of the coefficients are
\begin{equation}
F_m=-\frac{1}{5(b+2)}b^m+(-1)^{m-1}\frac{b+1}{5(b+2)}(5+b)^m, m\geq 1,
\label{5-g5-15}
\end{equation}
\begin{equation}
G_m=\frac{4b+5}{5(b+2)}b^{m-1}+(-1)^{m-1}\frac{1}{5(b+2)}(5+b)^m,
m\geq 1,
\label{5-g5-16}
\end{equation}
\begin{equation}
H_m=-\frac{6b+10}{5(b+2)}b^{m-1}+(-1)^m\frac{2}{5}(5+b)^m,
m\geq 1,
\label{5-g5-17}
\end{equation}
where $b$ is a solution of the equation 
$b^2+5b+5=0$.

Substituting the expressions of $D_m, E_m, F_m, G_m, H_m$ from Eqs. (\ref{5-g5-11})-(\ref{5-g5-12})
and (\ref{5-g5-15})-(\ref{5-g5-17})in Eq. (\ref{5-g5-2}) and grouping the terms
yields 
\begin{eqnarray}
\lefteqn{e^{(h_1-h_4)y}
=\frac{1}{5}+\frac{2}{5}\cos(\sqrt{-b}y)+\frac{2}{5}\cos(\sqrt{5+b}y)\nonumber}\\
&&+(h_1+h_4)
\left[\frac{1}{5}-\frac{b+3}{5}\cos(\sqrt{-b}y)+\frac{b+2}{5}\cos(\sqrt{5+b}y)\right]\nonumber\\
&&+(h_2+h_3)
\left[\frac{1}{5}+\frac{b+2}{5}\cos(\sqrt{-b}y)-\frac{b+3}{5}\cos(\sqrt{5+b}y)\right]\nonumber\\
&&+(h_1-h_4)
\left[\frac{\sqrt{-b}}{5}\sin(\sqrt{-b}y)+\frac{1}{\sqrt{-5b}}\sin(\sqrt{5+b}y)\right]\nonumber\\
&&+(h_2-h_3)
\left[-\frac{2b+5}{5\sqrt{-b}}\sin(\sqrt{-b}y)+\frac{b+2}{\sqrt{-5b}}\sin(\sqrt{5+b}y)\right].
\label{5-g5-18}
\end{eqnarray}

On the other hand, $e^{2h_1 y}$ can be written with the aid of the
5-dimensional polar cosexponential functions as
\begin{equation}
e^{2h_1 y}=g_{50}(2y)+h_1 g_{51}(2y)+h_2 g_{52}(2y)+h_3
g_{53}(2y)+h_4 g_{54}(2y).
\label{5-g5-19}
\end{equation}
The multiplication of the expressions of $e^{(h_1+h_4)y}$ and
$e^{(h_1-h_4)y}$ in Eqs. (\ref{5-g5-8}) and (\ref{5-g5-18}) and the
separation of the real components yields the expressions of the 5-dimensional
cosexponential functions, for $a=(\sqrt{5}-1)/2, b=-(5+\sqrt{5})/2$, as
\begin{eqnarray}
g_{50}(2y)=\frac{1}{5}e^{2y}
+\frac{2}{5}e^{ay}\cos(\sqrt{-b}y)
+\frac{2}{5}e^{-(1+a)y}\cos(\sqrt{5+b}y),
\label{5-g5-20}
\end{eqnarray}\index{cosexponential functions, polar 5-complex!expressions}
\begin{eqnarray}
\lefteqn{g_{51}(2y)=\frac{1}{5}e^{2y}
+\frac{1}{5}e^{ay}\left[\frac{-1+\sqrt{5}}{2}\cos(\sqrt{-b}y)+\frac{5+\sqrt{5}}{2\sqrt{-b}}\sin(\sqrt{-b}y)\right]\nonumber}\\
&&+\frac{1}{5}e^{-(1+a)y}\left[-\frac{1+\sqrt{5}}{2}\cos(\sqrt{5+b}y)+\sqrt{\frac{5}{-b}}\sin(\sqrt{5+b}y)\right],
\label{5-g5-21}
\end{eqnarray}
\begin{eqnarray}
\lefteqn{g_{52}(2y)=\frac{1}{5}e^{2y}
+\frac{1}{5}e^{ay}\left[-\frac{1+\sqrt{5}}{2}\cos(\sqrt{-b}y)+\sqrt{\frac{5}{-b}}\sin(\sqrt{-b}y)\right]\nonumber}\\
&&+\frac{1}{5}e^{-(1+a)y}\left[\frac{-1+\sqrt{5}}{2}\cos(\sqrt{5+b}y)-\frac{5+\sqrt{5}}{2\sqrt{-b}}\sin(\sqrt{5+b}y)\right],
\label{5-g5-22}
\end{eqnarray}
\begin{eqnarray}
\lefteqn{g_{53}(2y)=\frac{1}{5}e^{2y}
+\frac{1}{5}e^{ay}\left[-\frac{1+\sqrt{5}}{2}\cos(\sqrt{-b}y)-\sqrt{\frac{5}{-b}}\sin(\sqrt{-b}y)\right]\nonumber}\\
&&+\frac{1}{5}e^{-(1+a)y}\left[\frac{-1+\sqrt{5}}{2}\cos(\sqrt{5+b}y)+\frac{5+\sqrt{5}}{2\sqrt{-b}}\sin(\sqrt{5+b}y)\right],
\label{5-g5-23}
\end{eqnarray}
\begin{eqnarray}
\lefteqn{g_{54}(2y)=\frac{1}{5}e^{2y}
+\frac{1}{5}e^{ay}\left[\frac{-1+\sqrt{5}}{2}\cos(\sqrt{-b}y)-\frac{5+\sqrt{5}}{2\sqrt{-b}}\sin(\sqrt{-b}y)\right]\nonumber}\\
&&+\frac{1}{5}e^{-(1+a)y}\left[-\frac{1+\sqrt{5}}{2}\cos(\sqrt{5+b}y)-\sqrt{\frac{5}{-b}}\sin(\sqrt{5+b}y)\right].
\label{5-g5-24}
\end{eqnarray}
The polar 5-dimensional cosexponential functions can be written as
\begin{equation}
g_{5k}(y)=\frac{1}{5}\sum_{l=0}^{4}\exp\left[y\cos\left(\frac{2\pi l}{5}\right)
\right]
\cos\left[y\sin\left(\frac{2\pi l}{5}\right)-\frac{2\pi kl}{5}\right], k=0,...,4.
\label{5-30}
\end{equation}
The graphs of the polar 5-dimensional cosexponential functions are shown in
Fig. \ref{fig16}.

\begin{figure}
\begin{center}
\epsfig{file=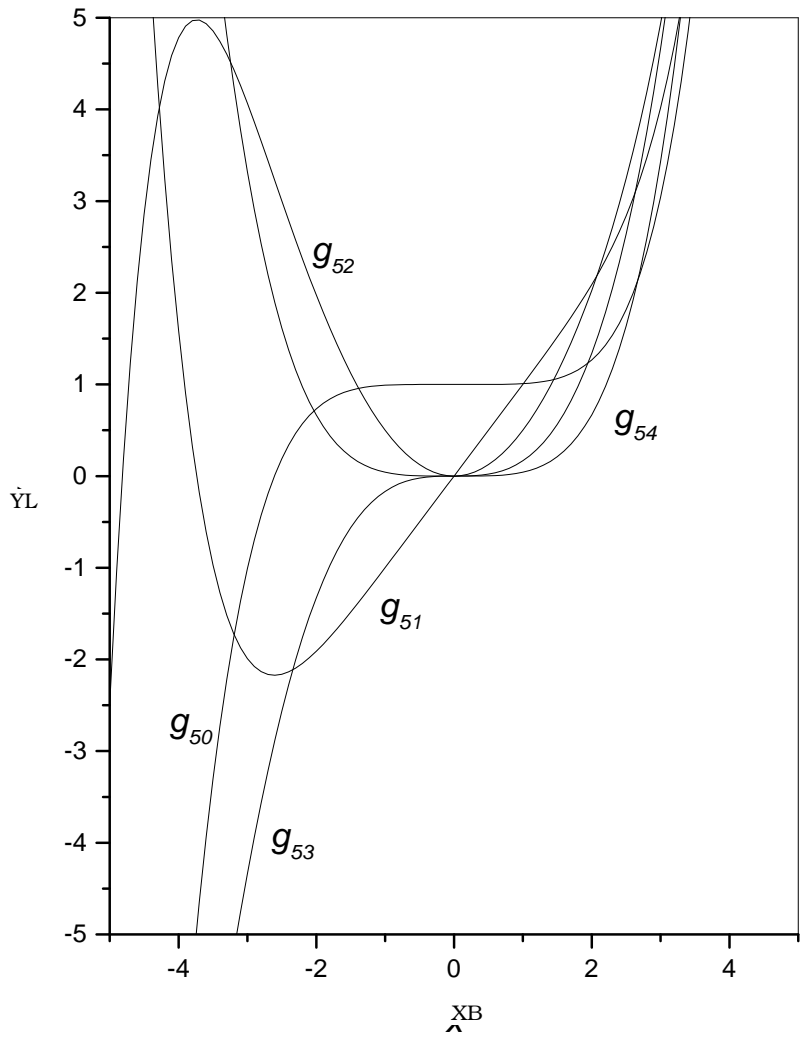,width=12cm}
\caption{Polar cosexponential functions $g_{50}, g_{51},g_{52}, g_{53},g_{54}$.}
\label{fig16}
\end{center}
\end{figure}

It can be checked that
\begin{equation}
\sum_{k=0}^{4}g_k^2(y)=\frac{1}{5}e^{2y}+\frac{2}{5}e^{(\sqrt{5}-1)y/2}
+\frac{2}{5}e^{-(\sqrt{5}+1)y/2}.
\label{5-34a}
\end{equation}

The addition theorems for the polar 5-dimensional cosexponential functions are
\begin{eqnarray}
\lefteqn{\begin{array}{c}
g_{50}(y+z)=g_{50}(y)g_{50}(z)
+g_{51}(y)g_{54}(z)+g_{52}(y)g_{53}(z)+g_{53}(y)g_{52}(z)+g_{54}(y)g_{51}(z) ,\\
g_{51}(y+z)=g_{50}(y)g_{51}(z)+g_{51}(y)g_{50}(z)
+g_{52}(y)g_{54}(z)+g_{53}(y)g_{53}(z)+g_{54}(y)g_{52}(z) ,\\
g_{52}(y+z)=g_{50}(y)g_{52}(z)+g_{51}(y)g_{51}(z)+g_{52}(y)g_{50}(z)
+g_{53}(y)g_{54}(z)+g_{54}(y)g_{53}(z) ,\\
g_{53}(y+z)=g_{50}(y)g_{53}(z)+g_{51}(y)g_{52}(z)+g_{52}(y)g_{51}(z)+g_{53}(y)g_{50}(z)
+g_{54}(y)g_{54}(z) ,\\
g_{54}(y+z)=g_{50}(y)g_{54}(z)
+g_{51}(y)g_{53}(z)+g_{52}(y)g_{52}(z)+g_{53}(y)g_{51}(z)+g_{54}(y)g_{50}(z) .\\
\end{array}\nonumber}\\
&&
\label{5-35a}
\end{eqnarray}\index{cosexponential functions, polar 5-complex!addition theorems}
It can be shown that
\begin{eqnarray}
\begin{array}{l}
\{g_{50}(y)+h_1g_{51}(y)+h_2g_{52}(y)+h_3g_{53}(y)+h_4g_{54}(y)\}^l\\
\hspace*{0.5cm}=g_{50}(ly)+h_1g_{51}(ly)+h_2g_{52}(ly)+h_3g_{53}(ly)+h_4g_{54}(ly),\\
\{g_{50}(y)+h_1g_{53}(y)+h_2g_{51}(y)+h_3g_{54}(y)+h_4g_{52}(y)\}^l\\
\hspace*{0.5cm}=g_{50}(ly)+h_1g_{53}(ly)+h_2g_{51}(ly)+h_3g_{54}(ly)+h_4g_{52}(ly),\\
\{g_{50}(y)+h_1g_{52}(y)+h_2g_{54}(y)+h_3g_{51}(y)+h_4g_{53}(y)\}^l\\
\hspace*{0.5cm}=g_{50}(ly)+h_1g_{52}(ly)+h_2g_{54}(ly)+h_3g_{51}(ly)+h_4g_{53}(ly),\\
\{g_{50}(y)+h_1g_{54}(y)+h_2g_{53}(y)+h_3g_{52}(y)+h_4g_{51}(y)\}^l\\
\hspace*{0.5cm}=g_{50}(ly)+h_1g_{54}(ly)+h_2g_{53}(ly)+h_3g_{52}(ly)+h_4g_{51}(ly).
\end{array}
\label{5-37b}
\end{eqnarray}

The derivatives of the polar cosexponential functions
are related by
\begin{equation}
\frac{dg_{50}}{du}=g_{54}, \:
\frac{dg_{51}}{du}=g_{50}, \:
\frac{dg_{52}}{du}=g_{51}, \:
\frac{dg_{53}}{du}=g_{52} ,\:
\frac{dg_{54}}{du}=g_{53}. \:
\label{5-45}
\end{equation}\index{cosexponential functions, polar 5-complex!differential equations}

\section{Exponential and trigonometric forms of polar 5-complex numbers}

The exponential and trigonometric forms of 5-complex
numbers can be expressed with the aid of the hypercomplex bases 
\begin{equation}
\left(
\begin{array}{c}
e_+\\
e_1\\
\tilde e_1\\
e_2\\
\tilde e_2\\
\end{array}\right)
=
\frac{2}{5}\left(
\begin{array}{ccccc}
\frac{1}{2}&\frac{1}{2}&\frac{1}{2}&\frac{1}{2}&\frac{1}{2}\\
1&p&2p^2-1&2p^2-1&p\\
0&q&2pq&-2pq&-q\\
1&2p^2-1&p&p&2p^2-1\\
0&2pq&-q&q&-2pq\\
\end{array}
\right)
\left(
\begin{array}{c}
1\\
h_1\\
h_2\\
h_3\\ 
h_4
\end{array}
\right).
\label{5-e12}
\end{equation}\index{canonical base!polar 5-complex}

The multiplication relations for these bases are
\begin{eqnarray}
\lefteqn{e_+^2=e_+,\;  e_+e_k=0,\; e_+\tilde e_k=0,\; \nonumber}\\ 
&&e_k^2=e_k,\; \tilde e_k^2=-e_k,\; e_k \tilde e_k=\tilde e_k ,\; e_ke_l=0,\;
e_k\tilde e_l=0,\; \tilde e_k\tilde e_l=0,\; k,l=1,2, \;k\not=l.\nonumber\\
&&
\label{5-e12b}
\end{eqnarray}
The bases have the property that
\begin{equation}
e_+ + e_1 +e_2=1.
\label{5-48b}
\end{equation}
The moduli of the new bases are
\begin{equation}
|e_+|=\frac{1}{\sqrt{5}},\; 
|e_k|=\sqrt{\frac{2}{5}},\; |\tilde e_k|=\sqrt{\frac{2}{5}}, 
\label{5-e12c}
\end{equation}
for $k=1,2$.

It can be checked that
\begin{eqnarray}
x_0+h_1x_1+h_2x_2+h_3x_3+h_4x_4=
e_+v_+ 
+e_1 v_1+\tilde e_1 \tilde v_1+e_2 v_2+\tilde e_2 \tilde v_2.
\label{5-e13b}
\end{eqnarray}\index{canonical form!polar 5-complex}
The ensemble $e_+, e_1, \tilde e_1, e_2, \tilde e_2$ will be called the
canonical 5-complex base, and Eq. (\ref{5-e13b}) gives the canonical form of the
5-complex number.

The exponential form of the 5-complex number $u$ is
\begin{eqnarray}
\lefteqn{u=\rho\exp\left\{\frac{1}{5}(h_1+h_2+h_3+h_4)\ln\frac{\sqrt{2}}{\tan\theta_+}\right.\nonumber}\\
&&\left.+\left[\frac{\sqrt{5}+1}{10}(h_1+h_4)-\frac{\sqrt{5}-1}{10}(h_2+h_3)\right]
\ln\tan\psi_1+\tilde e_1\phi_1+\tilde e_2\phi_2
\right\},
\label{5-50b}
\end{eqnarray}\index{exponential form!polar 5-complex}
for $0<\theta_+<\pi/2$.

The trigonometric form of the 5-complex number $u$ is
\begin{eqnarray}
\lefteqn{u=d\left(\frac{5}{2}\right)^{1/2}
\left(\frac{1}{\tan^2\theta_+}+1
+\frac{1}{\tan^2\psi_1}\right)^{-1/2}
\left(\frac{e_+\sqrt{2}}{\tan\theta_+}
+e_1+\frac{e_2}{\tan\psi_1}\right)
\exp\left(\tilde e_1\phi_1+\tilde e_2\phi_2\right).\nonumber}\\
&&
\label{5-52b}
\end{eqnarray}\index{trigonometric form!polar 5-complex}

The modulus $d$ and the amplitude $\rho$ are related by
\begin{eqnarray}
d=\rho \frac{2^{2/5}}{\sqrt{5}}
\left(\tan\theta_+
\tan^2\psi_1\right)^{1/5}
\left(\frac{1}{\tan^2\theta_+}+1
+\frac{1}{\tan^2\psi_1}\right)^{1/2}.
\label{5-53b}
\end{eqnarray}

\section{Elementary functions of a polar 5-complex variable}

The logarithm and power function exist for $v_+>0$, which means that
$0<\theta_+<\pi/2$, and are given by 
\begin{eqnarray}
\lefteqn{\ln u=\ln\rho+
\frac{1}{5}(h_1+h_2+h_3+h_4)\ln\frac{\sqrt{2}}{\tan\theta_+}\nonumber}\\
&&+\left[\frac{\sqrt{5}+1}{10}(h_1+h_4)-\frac{\sqrt{5}-1}{10}(h_2+h_3)\right]
\ln\tan\psi_1+\tilde e_1\phi_1+\tilde e_2\phi_2,
\label{5-56b}
\end{eqnarray}\index{logarithm!polar 5-complex}
\begin{equation}
u^m=e_+ v_+^m +
\rho_1^m(e_1\cos m\phi_1+\tilde e_1\sin m\phi_1)
+\rho_2^m(e_2\cos m\phi_2+\tilde e_2\sin m\phi_2).
\label{5-59b}
\end{equation}\index{power function!polar 5-complex}

The exponential of the 5-complex variable $u$ is
\begin{eqnarray}
e^u=e_+e^{v_+}  
+e^{v_1}\left(e_1 \cos \tilde v_1+\tilde e_1 \sin\tilde v_1\right)
+e^{v_2}\left(e_2 \cos \tilde v_2+\tilde e_2 \sin\tilde v_2\right).
\label{5-73b}
\end{eqnarray}\index{exponential, expressions!polar 5-complex}
The trigonometric functions of the
5-complex variable $u$ are
\begin{equation}
\cos u=e_+\cos v_+  
+\sum_{k=1}^{2}\left(e_k \cos v_k\cosh \tilde v_k
-\tilde e_k \sin v_k\sinh\tilde v_k\right),
\label{5-74c}
\end{equation}\index{trigonometric functions, expressions!polar 5-complex}
\begin{equation}
\sin u=e_+\sin v_+  
+\sum_{k=1}^{2}\left(e_k \sin v_k\cosh \tilde v_k
+\tilde e_k \cos v_k\sinh\tilde v_k\right).
\label{5-74d}
\end{equation}

The hyperbolic functions of the
5-complex variable $u$ are
\begin{equation}
\cosh u=e_+\cosh v_+  
+\sum_{k=1}^{2}\left(e_k \cosh v_k\cos \tilde v_k
+\tilde e_k \sinh v_k\sin\tilde v_k\right),
\label{5-75c}
\end{equation}\index{hyperbolic functions, expressions!polar 5-complex}
\begin{equation}
\sinh u=e_+\sinh v_+  
+\sum_{k=1}^{2}\left(e_k \sinh v_k\cos \tilde v_k
+\tilde e_k \cosh v_k\sin\tilde v_k\right).
\label{5-75d}
\end{equation}

\section{Power series of 5-complex numbers}

A power series of the 5-complex variable $u$ is a series of the form
\begin{equation}
a_0+a_1 u + a_2 u^2+\cdots +a_l u^l+\cdots .
\label{5-83}
\end{equation}\index{power series!polar 5-complex}
Since
\begin{equation}
|au^l|\leq 5^{l/2} |a| |u|^l ,
\label{5-82}
\end{equation}\index{modulus, inequalities!polar 5-complex}
the series is absolutely convergent for 
\begin{equation}
|u|<c,
\label{5-86}
\end{equation}\index{convergence of power series!polar 5-complex}
where 
\begin{equation}
c=\lim_{l\rightarrow\infty} \frac{|a_l|}{\sqrt{5}|a_{l+1}|} .
\label{5-87}
\end{equation}

If $a_l=\sum_{p=0}^{4}h_p a_{lp}$, where $h_0=1$, and
\begin{equation}
A_{l+}=\sum_{p=0}^{4}a_{lp},
\label{5-88a}
\end{equation}
\begin{equation}
A_{lk}=\sum_{p=0}^{4}a_{lp}\cos\left(\frac{2\pi kp}{5}\right),
\label{5-88b}
\end{equation}
\begin{equation}
\tilde A_{lk}=\sum_{p=0}^{4}a_{lp}\sin\left(\frac{2\pi kp}{5}\right),
\label{5-88c}
\end{equation}
for $k=1,2$,
the series (\ref{5-83}) can be written as
\begin{equation}
\sum_{l=0}^\infty \left[
e_+A_{l+}v_+^l+\sum_{k=1}^{2}
(e_k A_{lk}+\tilde e_k\tilde A_{lk})(e_k v_k+\tilde e_k\tilde v_k)^l 
\right].
\label{5-89b}
\end{equation}
The series in Eq. (\ref{5-83}) is absolutely convergent for 
\begin{equation}
|v_+|<c_+,\:
\rho_k<c_k, k=1,2,
\label{5-90}
\end{equation}\index{convergence, region of!polar 5-complex}
where 
\begin{equation}
c_+=\lim_{l\rightarrow\infty} \frac{|A_{l+}|}{|A_{l+1,+}|} ,\:
c_k=\lim_{l\rightarrow\infty} \frac
{\left(A_{lk}^2+\tilde A_{lk}^2\right)^{1/2}}
{\left(A_{l+1,k}^2+\tilde A_{l+1,k}^2\right)^{1/2}} .
\label{5-91}
\end{equation}

\section{Analytic functions of a polar 5-compex variable}

If $f(u)=\sum_{k=0}^{4}h_kP_k(x_0,x_1,x_2,x_3,x_{4})$,\index{functions, real
components!polar 5-complex} then
\begin{equation}
\frac{\partial P_0}{\partial x_0} 
=\frac{\partial P_1}{\partial x_1} 
=\frac{\partial P_2}{\partial x_2} 
=\frac{\partial P_3}{\partial x_3}
=\frac{\partial P_4}{\partial x_4}, 
\label{5-h95a}
\end{equation}\index{relations between partial derivatives!polar 5-complex}
\begin{equation}
\frac{\partial P_1}{\partial x_0} 
=\frac{\partial P_2}{\partial x_1} 
=\frac{\partial P_3}{\partial x_2} 
=\frac{\partial P_4}{\partial x_3}
=\frac{\partial P_0}{\partial x_4}, 
\label{5-h95b}
\end{equation}
\begin{equation}
\frac{\partial P_2}{\partial x_0} 
=\frac{\partial P_3}{\partial x_1} 
=\frac{\partial P_4}{\partial x_2} 
=\frac{\partial P_0}{\partial x_3}
=\frac{\partial P_1}{\partial x_4}, 
\label{5-h95c}
\end{equation}
\begin{equation}
\frac{\partial P_3}{\partial x_0} 
=\frac{\partial P_4}{\partial x_1} 
=\frac{\partial P_0}{\partial x_2} 
=\frac{\partial P_1}{\partial x_3}
=\frac{\partial P_2}{\partial x_4}, 
\label{5-h95d}
\end{equation}
\begin{equation}
\frac{\partial P_4}{\partial x_0} 
=\frac{\partial P_0}{\partial x_1} 
=\frac{\partial P_1}{\partial x_2} 
=\frac{\partial P_2}{\partial x_3}
=\frac{\partial P_3}{\partial x_4}, 
\label{5-h95e}
\end{equation}
and\index{relations between second-order derivatives!polar 5-complex}
\begin{eqnarray}
\lefteqn{\frac{\partial^2 P_k}{\partial x_0\partial x_l}
=\frac{\partial^2 P_k}{\partial x_1\partial x_{l-1}}
=\cdots=
\frac{\partial^2 P_k}{\partial x_{[l/2]}\partial x_{l-[l/2]}}}\nonumber\\
&&=\frac{\partial^2 P_k}{\partial x_{l+1}\partial x_{4}}
=\frac{\partial^2 P_k}{\partial x_{l+2}\partial x_{3}}
=\cdots
=\frac{\partial^2 P_k}{\partial x_{l+1+[(3-l)/2]}
\partial x_{4-[(3-l)/2]}} ,
\label{5-96}
\end{eqnarray}
for $k,l=0,...,4$.
In Eq. (\ref{5-96}), $[a]$ denotes the integer part of $a$,
defined as $[a]\leq a<[a]+1$. 
In this chapter, brackets larger than the regular brackets
$[\;]$ do not have the meaning of integer part.

\section{Integrals of polar 5-complex functions}

If $f(u)$ is an analytic 5-complex function,
then\index{integrals, path!polar 5-complex}\index{poles and residues!polar 5-complex}
\begin{equation}
\oint_\Gamma \frac{f(u)du}{u-u_0}=
2\pi f(u_0)\left\{\tilde e_1 
\;{\rm int}(u_{0\xi_1\eta_1},\Gamma_{\xi_1\eta_1})+
\tilde e_2 
\;{\rm int}(u_{0\xi_2\eta_2},\Gamma_{\xi_2\eta_2})\right\} ,
\label{5-120}
\end{equation}
where
\begin{equation}
{\rm int}(M,C)=\left\{
\begin{array}{l}
1 \;\:{\rm if} \;\:M \;\:{\rm is \;\:an \;\:interior \;\:point \;\:of} \;\:C ,\\ 
0 \;\:{\rm if} \;\:M \;\:{\rm is \;\:exterior \;\:to}\:\; C ,\\
\end{array}\right.,
\label{5-118}
\end{equation}
and $u_{0\xi_k\eta_k}$, $\Gamma_{\xi_k\eta_k}$ are respectively the
projections of the pole $u_0$ and of 
the loop $\Gamma$ on the plane defined by the axes $\xi_k$ and $\eta_k$,
$k=1,2$.

\section{Factorization of polar 5-complex polynomials}

A polynomial of degree $m$ of the 5-complex variable $u$ has the form
\begin{equation}
P_m(u)=u^m+a_1 u^{m-1}+\cdots+a_{m-1} u +a_m ,
\label{5-125}
\end{equation}
where $a_l$, for $l=1,...,m$, are 5-complex constants.
If $a_l=\sum_{p=0}^{4}h_p a_{lp}$, and with the
notations of Eqs. (\ref{5-88a})-(\ref{5-88c}) applied for $l= 1, \cdots, m$, the
polynomial $P_m(u)$ can be written as \index{polynomial, canonical variables!polar 5-complex}
\begin{eqnarray}
\lefteqn{P_m= 
e_+\left(v_+^m +\sum_{l=1}^{m}A_{l+}v_+^{m-l} \right) \nonumber}\\
&&+\sum_{k=1}^{2}
\left[(e_k v_k+\tilde e_k\tilde v_k)^m+
\sum_{l=1}^m(e_k A_{lk}+\tilde e_k\tilde A_{lk})
(e_k v_k+\tilde e_k\tilde v_k)^{m-l} 
\right].
\label{5-126b}
\end{eqnarray}

The polynomial $P_m(u)$ can be written, as 
\begin{eqnarray}
P_m(u)=\prod_{p=1}^m (u-u_p) ,
\label{5-128c}
\end{eqnarray}\index{polynomial, factorization!polar 5-complex}
where
\begin{eqnarray}
u_p=e_+ v_{p+}
+\left(e_1 v_{1p}+\tilde e_1\tilde v_{1p}\right)
+\left(e_2 v_{2p}+\tilde e_2\tilde v_{2p}\right), p=1,...,m.
\label{5-128e}
\end{eqnarray}
The quantities $v_{p+}$,   
$e_k v_{kp}+\tilde e_k\tilde v_{kp}$,
$p=1,...,m, k=1,2$,
are the roots of the corresponding polynomial in Eq. (\ref{5-126b}). The roots
$v_{p+}$ appear in complex-conjugate pairs, and  
$v_{kp}, \tilde v_{kp}$ are real numbers.
Since all these roots may be ordered arbitrarily, the polynomial $P_m(u)$ can be
written in many different ways as a product of linear factors. 

If $P(u)=u^2-1$, the degree is $m=2$, the coefficients of the polynomial are
$a_1=0, a_2=-1$, the coefficients defined in Eqs. (\ref{5-88a})-(\ref{5-88c})
are $A_{2+}=-1, A_{21}=-1, \tilde A_{21}=0,
A_{22}=-1, \tilde A_{22}=0$. The expression of $P(u)$, Eq. (\ref{5-126b}), is  
$v_+^2-e_++(e_1v_1+\tilde e_1\tilde v_1)^2-e_1+
(e_2v_2+\tilde e_2\tilde v_2)^2-e_2 $. 
The factorization of $P(u)$, Eq. (\ref{5-128c}), is
$P(u)=(u-u_1)(u-u_2)$, where the roots are
$u_1=\pm e_+\pm e_1\pm  e_2, u_2=-u_1$. If $e_+, e_1, e_2$ 
are expressed with the aid of Eq. (\ref{5-e12}) in terms of $h_1, h_2, h_3,
h_4$, the factorizations of $P(u)$ are obtained as
\begin{eqnarray}
\lefteqn{\begin{array}{l}
u^2-1=(u+1)(u-1),\\
u^2-1=\left[u+\frac{1}{5}+\frac{\sqrt{5}+1}{5}(h_1+h_4)
-\frac{\sqrt{5}-1}{5}(h_2+h_3)\right]
\left[u-\frac{1}{5}-\frac{\sqrt{5}+1}{5}(h_1+h_4)
+\frac{\sqrt{5}-1}{5}(h_2+h_3)\right],\\
u^2-1=\left[u+\frac{1}{5}-\frac{\sqrt{5}-1}{5}(h_1+h_4)
+\frac{\sqrt{5}+1}{5}(h_2+h_3)\right]
\left[u-\frac{1}{5}+\frac{\sqrt{5}-1}{5}(h_1+h_4)
-\frac{\sqrt{5}+1}{5}(h_2+h_3)\right],\\
u^2-1=\left[u+\frac{3}{5}-\frac{2}{5}(h_1+h_2+h_3+h_4)\right]
\left[u-\frac{3}{5}+\frac{2}{5}(h_1+h_2+h_3+h_4)\right].
\end{array}\nonumber}\\
&&
\end{eqnarray}
It can be checked that 
$(\pm e_+\pm e_1\pm e_2)^2=e_++e_1+e_2=1$.

\section{Representation of polar 5-complex numbers by irreducible matrices}

If the unitary matrix which can be obtained from the expression, 
Eq. (\ref{5-12}), of the variables $\xi_+, \xi_1, \eta_1, \xi_k, \eta_k$ in terms
of $x_0, x_1, x_2, x_3, x_4$ is called $T$,
the irreducible representation \cite{4} of the hypercomplex number $u$ is
\begin{equation}
T U T^{-1}=\left(
\begin{array}{ccc}
v_+     &     0     &     0      \\
0       &     V_1   &     0      \\
0       &     0     &     V_2    \\
\end{array}
\right),
\label{5-129b}
\end{equation}\index{representation by irreducible matrices!polar 5-complex}
where $U$ is the matrix in Eq. (\ref{5-24b}),
and $V_k$ are the matrices
\begin{equation}
V_k=\left(
\begin{array}{cc}
v_k           &     \tilde v_k   \\
-\tilde v_k   &     v_k          \\
\end{array}\right),\;\; k=1,2.
\label{5-130}
\end{equation}

\chapter{Complex Numbers in 6 Dimensions}

Two distinct systems of commutative complex numbers in 6 dimensions 
having the form $u=x_0+h_1x_1+h_2x_2+h_3x_3+h_4x_4+h_5x_5$ are
described in this chapter, 
for which the multiplication is associative and
commutative,
where the variables $x_0, x_1, x_2, x_3, x_4, x_5$ are
real numbers.  
The first type of 6-complex numbers described in this article is
characterized by the presence
of two polar axes, so that
these numbers will be called polar 6-complex numbers. 
The other type of 6-complex numbers described in this paper 
will be called planar n-complex numbers. 

The polar 6-complex numbers introduced in this chapter can be
specified by the modulus $d$, the amplitude $\rho$, and the polar angles $\theta_+,
\theta_-$, the planar angle $\psi_1$, and the azimuthal angles $\phi_1, \phi_2$. The
planar 6-complex numbers introduced in this paper can be specified by the
modulus $d$, the amplitude $\rho$, the planar angles $\psi_1, \psi_2$, and the
azimuthal angles $\phi_1, \phi_2, \phi_3$.  Exponential and trigonometric forms
are given for the 6-complex numbers.  The 6-complex functions defined by series
of powers are analytic, and the partial derivatives of the components of the
6-complex functions are closely related.  The integrals of polar 6-complex
functions are independent of path in regions where the functions are regular.
The fact that the exponential form of ther 6-complex numbers depends on cyclic
variables leads to the concept of pole and residue for integrals on closed
paths. The polynomials of polar 6-complex variables can be written as products
of linear or quadratic factors, the polynomials of planar 6-complex variables
can always be written as products of linear factors, although the factorization
is not unique.

The polar 6-complex numbers described in this paper are a particular case for 
$n=6$ of the polar hypercomplex numbers in $n$ dimensions discussed in Sec. 6.1,
and the planar 6-complex numbers described in this section are a particular case for 
$n=6$ of the planar hypercomplex numbers in $n$ dimensions discussed in Sec. 6.2.

\section{Polar complex numbers in 6 dimensions}

\subsection{Operations with polar complex numbers in $6$ dimensions}

The polar hypercomplex number $u$ in 6 dimensions 
is represented as  
\begin{equation}
u=x_0+h_1x_1+h_2x_2+h_3x_3+h_4x_4+h_5x_5. 
\label{6ch1a}
\end{equation}
The multiplication rules for the bases 
$h_1, h_2, h_3, h_4, h_5 $ are 
\begin{eqnarray}
\lefteqn{h_1^2=h_2,\;h_2^2=h_4,\;h_3^2=1,\;h_4^2=h_2,\;h_5^2=h_4,\;
h_1h_2=h_3,\;h_1h_3=h_4,\;h_1h_4=h_5,\nonumber}\\
&&\;h_1h_5=1,\;h_2h_3=h_5,\;h_2h_4=1,\;h_2h_5=h_1,\;h_3h_4=h_1,\;h_3h_5=h_2,\;h_4h_5=h_3.
\label{6ch1}
\end{eqnarray}\index{complex units!polar 6-complex}
The significance of the composition laws in Eq.
(\ref{6ch1}) can be understood by representing the bases $h_j, h_k$ by points on a
circle at the angles $\alpha_j=\pi j/3,\alpha_k=\pi k/3$, as shown in Fig. \ref{fig17},
and the product $h_j h_k$ by the point of the circle at the angle 
$\pi (j+k)/3$. If $2\pi\leq \pi (j+k)/3<4\pi$, the point represents the basis
$h_l$ of angle $\alpha_l=\pi(j+k)/3-2\pi$.

\begin{figure}
\begin{center}
\epsfig{file=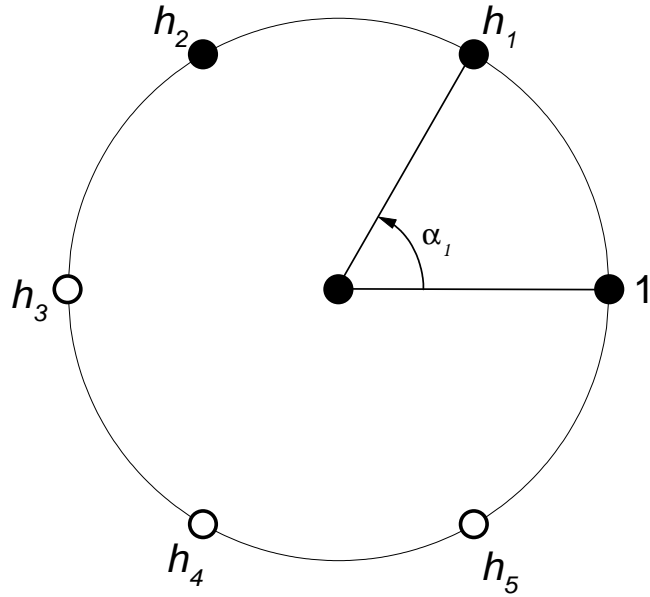,width=12cm}
\caption{Representation of the polar 
hypercomplex bases $1,h_1,h_2,h_3,h_4,h_5$
by points on a circle at the angles $\alpha_k=2\pi k/6$.
The product $h_j h_k$ will be represented by the point of the circle at the
angle $2\pi (j+k)/6$, $i,k=0,1,...,5$, where $h_0=1$. If $2\pi\leq 
2\pi (j+k)/6\leq 4\pi$, the point represents the basis
$h_l$ of angle $\alpha_l=2\pi(j+k)/6-2\pi$.}
\label{fig17}
\end{center}
\end{figure}

The sum of the 6-complex numbers $u$ and $u^\prime$ is
\begin{equation}
u+u^\prime=x_0+x^\prime_0+h_1(x_1+x^\prime_1)+h_1(x_2+x^\prime_2)
+h_3(x_3+x^\prime_3)+h_4(x_4+x^\prime_4)+h_5(x_5+x^\prime_5).
\label{6ch2}
\end{equation}\index{sum!polar 6-complex}
The product of the numbers $u, u^\prime$ is
\begin{equation}
\begin{array}{l}
uu^\prime=x_0 x_0^\prime +x_1x_5^\prime+x_2 x_4^\prime+x_3x_3^\prime
+x_4x_2^\prime+x_5 x_1^\prime\\
+h_1(x_0 x_1^\prime+x_1x_0^\prime+x_2x_5^\prime+x_3x_4^\prime
+x_4 x_3^\prime+x_5 x_2^\prime) \\
+h_2(x_0 x_2^\prime+x_1x_1^\prime+x_2x_0^\prime+x_3x_5^\prime
+x_4 x_4^\prime+x_5 x_3^\prime) \\
+h_3(x_0 x_3^\prime+x_1x_2^\prime+x_2x_1^\prime+x_3x_0^\prime
+x_4 x_5^\prime+x_5 x_4^\prime) \\
+h_4(x_0 x_4^\prime+x_1x_3^\prime+x_2x_2^\prime+x_3x_1^\prime
+x_4 x_0^\prime+x_5 x_5^\prime) \\
+h_5(x_0 x_5^\prime+x_1x_4^\prime+x_2x_3^\prime+x_3x_2^\prime
+x_4 x_1^\prime+x_5 x_0^\prime).
\end{array}
\label{6ch3}
\end{equation}\index{product!polar 6-complex}

The relation between the variables $v_+,v_-,v_1,\tilde v_1, v_2, \tilde v_2$
and $x_0,x_1,x_2,x_3,x_4,x_5$ are
\begin{equation}
\left(
\begin{array}{c}
v_+\\
v_-\\
v_1\\
\tilde v_1\\
v_2\\
\tilde v_2\\
\end{array}\right)
=\left(
\begin{array}{cccccc}
1&1&1&1&1&1\\
1&-1&1&-1&1&-1\\
1&\frac{1}{2} &-\frac{1}{2} &-1 &-\frac{1}{2} &\frac{1}{2} \\
0& \frac{\sqrt{3}}{2}&\frac{\sqrt{3}}{2} &0 &-\frac{\sqrt{3}}{2} &-\frac{\sqrt{3}}{2} \\
1& -\frac{1}{2}&-\frac{1}{2} &1 &-\frac{1}{2} &-\frac{1}{2} \\
0&\frac{\sqrt{3}}{2} &-\frac{\sqrt{3}}{2} &0 &\frac{\sqrt{3}}{2} &-\frac{\sqrt{3}}{2} \\
\end{array}
\right)
\left(
\begin{array}{c}
x_0\\
x_1\\
x_2\\
x_3\\
x_4\\
x_5
\end{array}
\right).
\label{6ch9e}
\end{equation}\index{canonical variables!polar 6-complex}
The other variables are $v_4=v_2, \tilde v_4=-\tilde v_2,
v_5=v_1, \tilde v_5=-\tilde v_1$. 
The variables $v_+, v_-, v_1, \tilde v_1, v_2, \tilde v_2$ will be called
canonical polar 6-complex variables.

\subsection{Geometric representation of polar complex numbers in 6 dimensions}

The 6-complex number $u=x_0+h_1x_1+h_2x_2+h_3x_3+h_4x_4+h_5x_5$
is represented by 
the point $A$ of coordinates $(x_0,x_1,x_2,x_3,x_4,x_5)$. 
The distance from the origin $O$ of the 6-dimensional space to the point $A$
has the expression 
\begin{equation}
d^2=x_0^2+x_1^2+x_2^2+x_3^2+x_4^2+x_5^2.
\label{6ch10}
\end{equation}\index{distance!polar 6-complex}
The distance $d$ is called modulus of the 6-complex number 
$u$, and is designated by $d=|u|$.
The modulus has the property that
\begin{equation}
|u^\prime u^{\prime\prime}|\leq \sqrt{6}|u^\prime||u^{\prime\prime}| .
\label{6ch79}
\end{equation}

The exponential and trigonometric forms of the 6-complex number $u$ can be
obtained conveniently in a rotated system of axes defined by a transformation
which has the form
\begin{equation}
\left(
\begin{array}{c}
\xi_+\\
\xi_-\\
\xi_1\\
\tilde \xi_1\\
\xi_2\\
\tilde \xi_2\\
\end{array}\right)
=\left(
\begin{array}{cccccc}
\frac{1}{\sqrt{6}}&\frac{1}{\sqrt{6}}&\frac{1}{\sqrt{6}}&\frac{1}{\sqrt{6}}&\frac{1}{\sqrt{6}}&\frac{1}{\sqrt{6}}\\
\frac{1}{\sqrt{6}}&-\frac{1}{\sqrt{6}}&\frac{1}{\sqrt{6}}&-\frac{1}{\sqrt{6}}&\frac{1}{\sqrt{6}}&-\frac{1}{\sqrt{6}}\\
\frac{\sqrt{3}}{3}&\frac{\sqrt{3}}{6} &-\frac{\sqrt{3}}{6} &-\frac{\sqrt{3}}{3} &-\frac{\sqrt{3}}{6} &\frac{\sqrt{3}}{6} \\
0& \frac{1}{2}&\frac{1}{2} &0 &-\frac{1}{2} &-\frac{1}{2} \\
\frac{\sqrt{3}}{3}& -\frac{\sqrt{3}}{6}&-\frac{\sqrt{3}}{6} &\frac{\sqrt{3}}{3} &-\frac{\sqrt{3}}{6} &-\frac{\sqrt{3}}{6} \\
0&\frac{1}{2} &-\frac{1}{2} &0 &\frac{1}{2} &-\frac{1}{2} \\
\end{array}
\right)
\left(
\begin{array}{c}
x_0\\
x_1\\
x_2\\
x_3\\
x_4\\
x_5
\end{array}
\right).
\label{6ch9ee}
\end{equation}
The lines of the matrices in Eq. (\ref{6ch9ee})  gives the components
of the 6 basis vectors of the new system of axes. These vectors have unit
length and are orthogonal to each other.
The relations between the two sets of variables are
\begin{equation}
v_+= \sqrt{6}\xi_+ ,  v_-=  \sqrt{6}\xi_-,  
v_k= \sqrt{3}\xi_k , \tilde v_k= \sqrt{3}\eta_k, k=1,2 .
\label{6ch12b}
\end{equation}

The radius $\rho_k$ and the azimuthal angle $\phi_k$ in the plane of the axes
$v_k,\tilde v_k$ are
\begin{equation}
\rho_k^2=v_k^2+\tilde v_k^2, \:\cos\phi_k=v_k/\rho_k,
\:\sin\phi_k=\tilde v_k/\rho_k, 0\leq \phi_k<2\pi ,\;k=1,2,
\label{6ch19a}
\end{equation}
so that there are 2 azimuthal angles.\index{azimuthal angles!polar 6-complex}
The planar angle $\psi_1$ is
\begin{equation}
\tan\psi_1=\rho_1/\rho_2, 0\leq\psi_1\leq\pi/2.
\label{6ch19b}
\end{equation}\index{planar angle!polar 6-complex}
There is a polar angle $\theta_+$, 
\begin{equation}
\tan\theta_+=\frac{\sqrt{2}\rho_1}{v_+}, 0\leq\theta_+\leq\pi , 
\label{6ch19c}
\end{equation}
and there is also a polar angle $\theta_-$,
\begin{equation}
\tan\theta_-=\frac{\sqrt{2}\rho_1}{v_-}, 0\leq\theta_-\leq\pi .
\label{6ch19d}
\end{equation}\index{polar angles!polar 6-complex}
The amplitude of a 6-complex number $u$ is
\begin{equation}
\rho=\left(v_+v_-\rho_1^2 \rho_2^2\right)^{1/6}.
\label{6ch50aa}
\end{equation}\index{amplitude!polar 6-complex}
It can be checked that
\begin{equation}
d^2=\frac{1}{6}v_+^2+\frac{1}{6}v_-^2
+\frac{1}{3}(\rho_1^2+\rho_2^2) .
\label{6ch17}
\end{equation}\index{modulus, canonical variables!polar 6-complex}

If $u=u^\prime u^{\prime\prime}$, the parameters 
of the hypercomplex numbers are related by
\begin{equation}
v_+=v_+^\prime v_+^{\prime\prime},  
\label{6ch21a}
\end{equation}\index{transformation of variables!polar 6-complex}
\begin{equation}
\tan\theta_+=\frac{1}{\sqrt{2}}\tan\theta_+^\prime \tan\theta_+^{\prime\prime},
\label{6ch21c}
\end{equation}
\begin{equation}
v_-=v_-^\prime v_-^{\prime\prime},  
\label{6ch21f}
\end{equation}
\begin{equation}
\tan\theta_-=\frac{1}{\sqrt{2}}\tan\theta_-^\prime \tan\theta_-^{\prime\prime},
\label{6ch21g}
\end{equation}
\begin{equation}
\tan\psi_1=\tan\psi_1^\prime \tan\psi_1^{\prime\prime},
\label{6ch21d}
\end{equation}
\begin{equation}
\rho_k=\rho_k^\prime\rho_k^{\prime\prime}, 
\label{6ch21b}
\end{equation}
\begin{equation}
\phi_k=\phi_k^\prime+\phi_k^{\prime\prime},  
\label{6ch21e}
\end{equation}
\begin{equation}
v_k=v_k^\prime v_k^{\prime\prime}-\tilde v_k^\prime \tilde v_k^{\prime\prime},\;
\tilde v_k=v_k^\prime \tilde v_k^{\prime\prime}+\tilde v_k^\prime v_k^{\prime\prime},
\label{6ch22}
\end{equation}
\begin{equation}
\rho=\rho^\prime\rho^{\prime\prime} ,
\label{6ch24}
\end{equation}
where $k=1,2$.

The 6-complex number
$u=x_0+h_1x_1+h_2x_2+h_3x_3+h_4x_4+h_5x_5$ can be  represented by the matrix
\begin{equation}
U=\left(
\begin{array}{cccccc}
x_0     &   x_1     &   x_2   &   x_3  &  x_4  & x_5\\
x_5     &   x_0     &   x_1   &   x_2  &  x_3  & x_4\\
x_4     &   x_5     &   x_0   &   x_1  &  x_2  & x_3\\
x_3     &   x_4     &   x_5   &   x_0  &  x_1  & x_2\\
x_2     &   x_3     &   x_4   &   x_5  &  x_0  & x_1\\
x_1     &   x_2     &   x_3   &   x_4  &  x_5  & x_0\\
\end{array}
\right).
\label{6ch24b}
\end{equation}\index{matrix representation!polar 6-complex}
The product $u=u^\prime u^{\prime\prime}$ is
represented by the matrix multiplication $U=U^\prime U^{\prime\prime}$.

\subsection{The polar 6-dimensional cosexponential functions}

The polar cosexponential functions in 6 dimensions are
\begin{equation}
g_{6k}(y)=\sum_{p=0}^\infty y^{k+6p}/(k+6p)!, 
\label{6ch29}
\end{equation}\index{cosexponential functions, polar 6-complex!definitions}
for $ k=0,...,5$.
The polar cosexponential functions $g_{6k}$ of even index $k$ are
even functions, $g_{6,2p}(-y)=g_{6,2p}(y)$, 
and the polar cosexponential functions of odd index $k$
are odd functions, $g_{6,2p+1}(-y)=-g_{6,2p+1}(y)$, $p=0,1,2$. 
\index{cosexponential functions, polar 6-complex!parity}

It can be checked that
\begin{equation}
\sum_{k=0}^{5}g_{6k}(y)=e^y,
\label{6ch29a}
\end{equation}
\begin{equation}
\sum_{k=0}^{5}(-1)^k g_{6k}(y)=e^{-y}.
\label{6ch29b}
\end{equation}
\newpage
\setlength{\evensidemargin}{-1cm}       
The exponential function of the quantity $h_k y$ is
\begin{equation}
\begin{array}{l}
e^{h_1 y}=g_{60}(y)+h_1g_{61}(y)+h_2g_{62}(y)+h_3g_{63}(y)+h_4g_{64}(y)+h_5g_{65}(y),\\
e^{h_2 y}=g_{60}(y)+g_{63}(y)+h_2\{g_{61}(y)+g_{64}(y)\}+h_4\{g_{62}(y)+g_{65}(y)\},\\
e^{h_3 y}=g_{60}(y)+g_{62}(y)+g_{64}(y)+h_3\{g_{61}(y)+g_{63}(y)+g_{65}(y)\},\\
e^{h_4 y}=g_{60}(y)+g_{63}(y)+h_2\{g_{62}(y)+g_{65}(y)\}+h_4\{g_{61}(y)+g_{64}(y)\},\\
e^{h_5 y}=g_{60}(y)+h_1g_{65}(y)+h_2g_{64}(y)+h_3g_{63}(y)+h_4g_{62}(y)+h_5g_{61}(y).\\
\end{array}
\label{6ch28b}
\end{equation}\index{exponential, expressions!polar 6-complex}
The relations for $h_2$ and $h_4$ can be written equivalently as 
$e^{h_2 y}=g_{30}+h_2g_{31}+h_4 g_{32}, e^{h_4 y}=g_{30}+h_2g_{32}+h_4 g_{31}$,
and the relation for $h_3$ can be written as $e^{h_3 y}=g_{20}+h_3g_{21}$,
which is the same as $e^{h_3 y}=\cosh y+h_3\sinh y$.

The expressions of the polar 6-dimensional cosexponential functions are
\begin{equation}
\begin{array}{l}
g_{60}(y)=\frac{1}{3}\cosh y +\frac{2}{3}\cosh\frac{y}{2}\cos\frac{\sqrt{3}}{2}y,\\
g_{61}(y)=\frac{1}{3}\sinh y +\frac{1}{3}\sinh\frac{y}{2}\cos\frac{\sqrt{3}}{2}y
+\frac{\sqrt{3}}{3}\cosh\frac{y}{2}\sin\frac{\sqrt{3}}{2}y,\\
g_{62}(y)=\frac{1}{3}\cosh y -\frac{1}{3}\cosh\frac{y}{2}\cos\frac{\sqrt{3}}{2}y
+\frac{\sqrt{3}}{3}\sinh\frac{y}{2}\sin\frac{\sqrt{3}}{2}y,\\
g_{63}(y)=\frac{1}{3}\sinh y -\frac{2}{3}\sinh\frac{y}{2}\cos\frac{\sqrt{3}}{2}y,\\
g_{64}(y)=\frac{1}{3}\cosh y -\frac{1}{3}\cosh\frac{y}{2}\cos\frac{\sqrt{3}}{2}y
-\frac{\sqrt{3}}{3}\sinh\frac{y}{2}\sin\frac{\sqrt{3}}{2}y,\\
g_{65}(y)=\frac{1}{3}\sinh y +\frac{1}{3}\sinh\frac{y}{2}\cos\frac{\sqrt{3}}{2}y
-\frac{\sqrt{3}}{3}\cosh\frac{y}{2}\sin\frac{\sqrt{3}}{2}y.
\end{array}
\label{6ch30}
\end{equation}\index{cosexponential functions, polar 6-complex!expressions}
The cosexponential functions (\ref{6ch30}) can be written as
\begin{equation}
g_{6k}(y)=\frac{1}{6}\sum_{l=0}^{5}\exp\left[y\cos\left(\frac{2\pi l}{6}\right)
\right]
\cos\left[y\sin\left(\frac{2\pi l}{6}\right)-\frac{2\pi kl}{6}\right], 
\label{6ch30x}
\end{equation}
for $k=0,...,5$.
The graphs of the polar 6-dimensional cosexponential functions are shown in 
Fig. \ref{fig18}.

\begin{figure}
\begin{center}
\epsfig{file=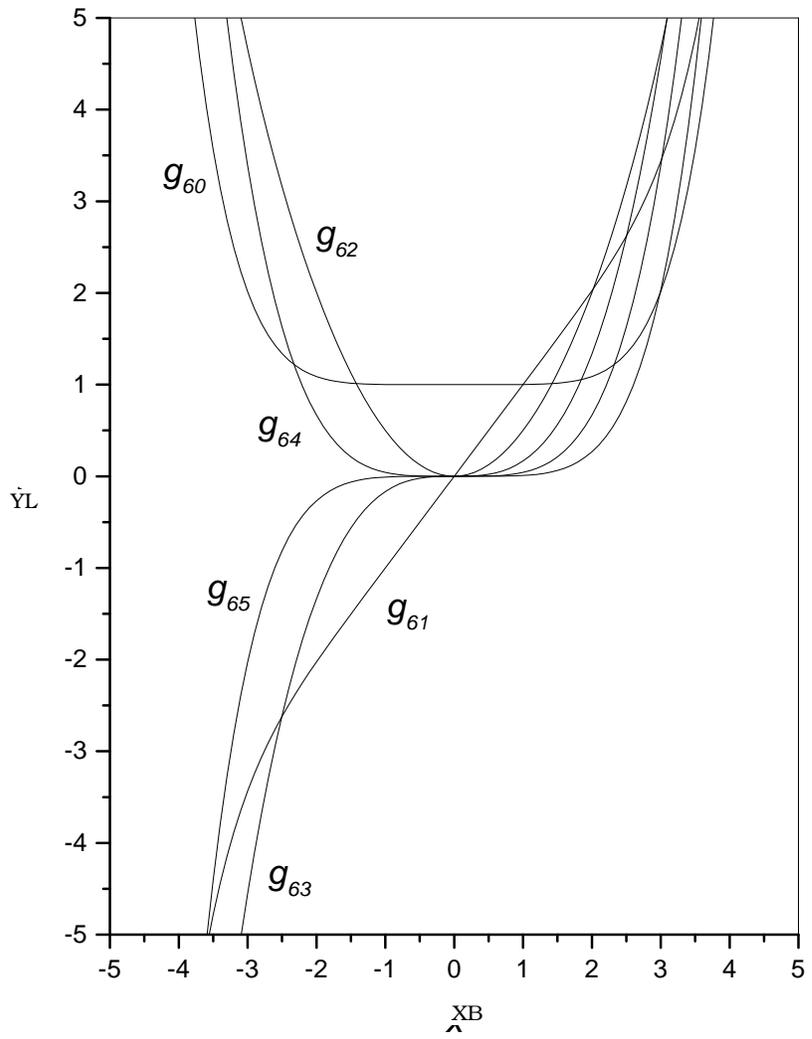,width=12cm}
\caption{Polar cosexponential functions $g_{60}, g_{61},g_{62}, g_{63},g_{64}, g_{65}$.}
\label{fig18}
\end{center}
\end{figure}

It can be checked that
\begin{equation}
\sum_{k=0}^{5}g_{6k}^2(y)=\frac{1}{3}\cosh 2y +\frac{2}{3}\cosh y.
\label{6ch34a}
\end{equation}

The addition theorems for the polar 6-dimensional cosexponential functions are
\begin{eqnarray}
\lefteqn{\begin{array}{l}
g_{60}(y+z)=g_{60}(y)g_{60}(z)+g_{61}(y)g_{65}(z)+g_{62}(y)g_{64}(z)
+g_{63}(y)g_{63}(z)+g_{64}(y)g_{62}(z)+g_{65}(y)g_{61}(z) ,\\
g_{61}(y+z)=g_{60}(y)g_{61}(z)+g_{61}(y)g_{60}(z)+g_{62}(y)g_{65}(z)
+g_{63}(y)g_{64}(z)+g_{64}(y)g_{63}(z)+g_{65}(y)g_{62}(z) ,\\
g_{62}(y+z)=g_{60}(y)g_{62}(z)+g_{61}(y)g_{61}(z)+g_{62}(y)g_{60}(z)
+g_{63}(y)g_{65}(z)+g_{64}(y)g_{64}(z)+g_{65}(y)g_{63}(z) ,\\
g_{63}(y+z)=g_{60}(y)g_{63}(z)+g_{61}(y)g_{62}(z)+g_{62}(y)g_{61}(z)
+g_{63}(y)g_{60}(z)+g_{64}(y)g_{65}(z)+g_{65}(y)g_{64}(z) ,\\
g_{64}(y+z)=g_{60}(y)g_{64}(z)+g_{61}(y)g_{63}(z)+g_{62}(y)g_{62}(z)
+g_{63}(y)g_{61}(z)+g_{64}(y)g_{60}(z)+g_{65}(y)g_{65}(z) ,\\
g_{65}(y+z)=g_{60}(y)g_{65}(z)+g_{61}(y)g_{64}(z)+g_{62}(y)g_{63}(z)
+g_{63}(y)g_{62}(z)+g_{64}(y)g_{61}(z)+g_{65}(y)g_{60}(z) .
\end{array}\nonumber}\\
&&
\label{6ch35a}
\end{eqnarray}
\index{cosexponential functions, polar 6-complex!addition theorems}
\newpage
\setlength{\evensidemargin}{-0.4cm}       
It can be shown that
\begin{eqnarray}
\begin{array}{l}
\{g_{60}(y)+h_1g_{61}(y)+h_2g_{62}(y)+h_3g_{63}(y)+h_4g_{64}(y)+h_5g_{65}(y)\}^l\\
\hspace{0.5cm}=g_{60}(ly)+h_1g_{61}(ly)+h_2g_{62}(ly)+h_3g_{63}(ly)+h_4g_{64}(ly)+h_5g_{65}(ly),\\
\{g_{60}(y)+g_{63}(y)+h_2\{g_{61}(y)+g_{64}(y)\}+h_4\{g_{62}(y)+g_{65}(y)\}\}^l\\
\hspace{0.5cm}=g_{60}(ly)+g_{63}(ly)+h_2\{g_{61}(ly)+g_{64}(ly)\}+h_4\{g_{62}(ly)+g_{65}(ly)\},\\
\{g_{60}(y)+g_{62}(y)+g_{64}(y)+h_3\{g_{61}(y)+g_{63}(y)+g_{65}(y)\}\}^l\\
\hspace{0.5cm}=g_{60}(ly)+g_{62}(ly)+g_{64}(ly)+h_3\{g_{61}(ly)+g_{63}(ly)+g_{65}(ly)\},\\
\{g_{60}(y)+g_{63}(y)+h_2\{g_{62}(y)+g_{65}(y)\}+h_4\{g_{61}(y)+g_{64}(y)\}\}^l\\
\hspace{0.5cm}=g_{60}(ly)+g_{63}(ly)+h_2\{g_{62}(ly)+g_{65}(ly)\}+h_4\{g_{61}(ly)+g_{64}(ly)\},\\
\{g_{60}(y)+h_1g_{65}(y)+h_2g_{64}(y)+h_3g_{63}(y)+h_4g_{62}(y)+h_5g_{61}(y)\}^l\\
\hspace{0.5cm}=g_{60}(ly)+h_1g_{65}(ly)+h_2g_{64}(ly)+h_3g_{63}(ly)+h_4g_{62}(ly)+h_5g_{61}(ly).\\
\end{array}
\label{6ch37b}
\end{eqnarray}

The derivatives of the polar cosexponential functions
are related by
\begin{equation}
\frac{dg_{60}}{du}=g_{65}, \:
\frac{dg_{61}}{du}=g_{60}, \:
\frac{dg_{62}}{du}=g_{61}, \:
\frac{dg_{63}}{du}=g_{62} ,\:
\frac{dg_{64}}{du}=g_{63}, \:
\frac{dg_{65}}{du}=g_{64} .
\label{6ch45}
\end{equation}\index{cosexponential functions, polar 6-complex!differential
equations}

\subsection{Exponential and trigonometric forms of polar 6-complex numbers}

The exponential and trigonometric forms of polar 6-complex
numbers can be expressed with the aid of the hypercomplex bases 
\begin{equation}
\left(
\begin{array}{c}
e_+\\
e_-\\
e_1\\
\tilde e_1\\
e_2\\
\tilde e_2\\
\end{array}\right)
=\left(
\begin{array}{cccccc}
\frac{1}{6}&\frac{1}{6}&\frac{1}{6}&\frac{1}{6}&\frac{1}{6}&\frac{1}{6}\\
\frac{1}{6}&-\frac{1}{6}&\frac{1}{6}&-\frac{1}{6}&\frac{1}{6}&-\frac{1}{6}\\
\frac{1}{3}&\frac{1}{6} &-\frac{1}{6} &-\frac{1}{3} &-\frac{1}{6} &\frac{1}{6} \\
0& \frac{\sqrt{3}}{6}&\frac{\sqrt{3}}{6} &0 &-\frac{\sqrt{3}}{6} &-\frac{\sqrt{3}}{6} \\
\frac{1}{3}& -\frac{1}{6}&-\frac{1}{6} &\frac{1}{3} &-\frac{1}{6} &-\frac{1}{6} \\
0&\frac{\sqrt{3}}{6} &-\frac{\sqrt{3}}{6} &0 &\frac{\sqrt{3}}{6} &-\frac{\sqrt{3}}{6} \\
\end{array}
\right)
\left(
\begin{array}{c}
1\\
h_1\\
h_2\\
h_3\\
h_4\\
h_5
\end{array}
\right).
\label{6che11}
\end{equation}
\index{canonical base!polar 6-complex}

The multiplication relations for these bases are
\begin{eqnarray}
\lefteqn{e_+^2=e_+,\; e_-^2=e_-,\; e_+e_-=0,\; e_+e_k=0,\; e_+\tilde e_k=0,\;
e_-e_k=0,\; 
e_-\tilde e_k=0,\;e_k^2=e_k,\nonumber}\\ 
&& \tilde e_k^2=-e_k,\; e_k \tilde e_k=\tilde e_k ,\; e_ke_l=0,\;
e_k\tilde e_l=0,\; \tilde e_k\tilde e_l=0,\; k,l=1,2, k\not=l.\; 
\label{6che12a}
\end{eqnarray}
The bases have the property that
\begin{equation}
e_+ +e_- + e_1+e_2 =1.
\label{6ch47b}
\end{equation}
The moduli of the new bases are
\begin{equation}
|e_+|=\frac{1}{\sqrt{6}},\; |e_-|=\frac{1}{\sqrt{6}},\; 
|e_k|=\frac{1}{\sqrt{3}},\; |\tilde e_k|=\frac{1}{\sqrt{3}}, k=1,2.
\label{6che12c}
\end{equation}

It can be shown that
\begin{eqnarray}
\lefteqn{x_0+h_1x_1+h_2x_2+h_3x_3+h_4x_4+h_5x_5\nonumber}\\
&&=e_+v_+ + e_-v_- +e_1 v_1+\tilde e_1 \tilde v_1
+e_2 v_2+\tilde e_2 \tilde v_2.
\label{6che13a}
\end{eqnarray}\index{canonical form!polar 6-complex}
The ensemble $e_+, e_-, e_1, \tilde e_1, e_2, \tilde e_2$ will be called the
canonical polar 6-complex base, and Eq. (\ref{6che13a}) gives the canonical
form of the polar 6-complex number.

The exponential form of the 6-complex number $u$ is 
\begin{eqnarray}\lefteqn{
u=\rho\exp\left\{\frac{1}{6}(h_1+h_2+h_3+h_4+h_5)\ln\frac{\sqrt{2}}{\tan\theta_+}
-\frac{1}{6}(h_1-h_2+h_3-h_4+h_5)\ln\frac{\sqrt{2}}{\tan\theta_-}\right.\nonumber}\\
&&
\left.+\frac{1}{6}(h_1+h_2-2h_3+h_4+h_5)\ln\tan\psi_1
+\tilde e_1\phi_1+\tilde e_2\phi_2 \right\},
\label{6ch50a}
\end{eqnarray}\index{exponential form!polar 6-complex}
for $0<\theta_+<\pi/2, 0<\theta_-<\pi/2$.

The
trigonometric form of the 6-complex number $u$ is
\begin{eqnarray}
\lefteqn{u=d
\sqrt{3}
\left(\frac{1}{\tan^2\theta_+}+\frac{1}{\tan^2\theta_-}+1
+\frac{1}{\tan^2\psi_1}\right)^{-1/2}\nonumber}\\
&&\left(\frac{e_+\sqrt{2}}{\tan\theta_+}+\frac{e_-\sqrt{2}}{\tan\theta_-}
+e_1+\frac{e_2}{\tan\psi_1}\right)
\exp\left(\tilde e_1\phi_1+\tilde e_2\phi_2\right).
\label{6ch52a}
\end{eqnarray}\index{trigonometric form!polar 6-complex}

The modulus $d$ and the amplitude $\rho$ are related by
\begin{eqnarray}
\lefteqn{d=\rho \frac{2^{1/3}}{\sqrt{6}}
\left(\tan\theta_+\tan\theta_-
\tan^2\psi_1\right)^{1/6}\nonumber}\\
&&\left(\frac{1}{\tan^2\theta_+}+\frac{1}{\tan^2\theta_-}+1
+\frac{1}{\tan^2\psi_1}\right)^{1/2}.
\label{6ch53a}
\end{eqnarray}

\subsection{Elementary functions of a polar 6-complex variable}

The logarithm and power functions of the 6-complex number $u$ exist for
$v_+>0, v_->0$, which means that $0<\theta_+<\pi/2, 0<\theta_-<\pi/2$,
and are given by
\begin{eqnarray}\lefteqn{
\ln u=\ln \rho+
\frac{1}{6}(h_1+h_2+h_3+h_4+h_5)\ln\frac{\sqrt{2}}{\tan\theta_+}
-\frac{1}{6}(h_1-h_2+h_3-h_4+h_5)\ln\frac{\sqrt{2}}{\tan\theta_-}
\nonumber}\\
&&
+\frac{1}{6}(h_1+h_2-2h_3+h_4+h_5)\ln\tan\psi_1
+\tilde e_1\phi_1+\tilde e_2\phi_2 ,
\label{6ch56a}
\end{eqnarray}\index{logarithm!polar 6-complex}
\begin{equation}
u^m=e_+ v_+^m+e_- v_-^m +
\rho_1^m(e_1\cos m\phi_1+\tilde e_1\sin m\phi_1)
+\rho_2^m(e_2\cos m\phi_2+\tilde e_2\sin m\phi_2).
\label{6ch59a}
\end{equation}\index{power function!polar 6-complex}

The exponential of the 6-complex variable $u$ is 
\begin{eqnarray}
e^u=e_+e^{v_+} + e_-e^{v_-} 
+e^{v_1}\left(e_1 \cos \tilde v_1+\tilde e_1 \sin\tilde v_1\right)
+e^{v_2}\left(e_2 \cos \tilde v_2+\tilde e_2 \sin\tilde v_2\right).
\label{6ch73a}
\end{eqnarray}\index{exponential, expressions!polar 6-complex}
The trigonometric functions of the
6-complex variable $u$ are 
\begin{equation}
\cos u=e_+\cos v_+ + e_-\cos v_- 
+\sum_{k=1}^{2}\left(e_k \cos v_k\cosh \tilde v_k
-\tilde e_k \sin v_k\sinh\tilde v_k\right),
\label{6ch74a}
\end{equation}\index{trigonometric functions, expressions!polar 6-complex}
\begin{equation}
\sin u=e_+\sin v_+ + e_-\sin v_- 
+\sum_{k=1}^{2}\left(e_k \sin v_k\cosh \tilde v_k
+\tilde e_k \cos v_k\sinh\tilde v_k\right).
\label{6ch74b}
\end{equation}

The hyperbolic functions of the
6-complex variable $u$ are
\begin{equation}
\cosh u=e_+\cosh v_+ + e_-\cosh v_- 
+\sum_{k=1}^{2}\left(e_k \cosh v_k\cos \tilde v_k
+\tilde e_k \sinh v_k\sin\tilde v_k\right),
\label{6ch75a}
\end{equation}\index{hyperbolic functions, expressions!polar 6-complex}
\begin{equation}
\sinh u=e_+\sinh v_+ + e_-\sinh v_- 
+\sum_{k=1}^{2}\left(e_k \sinh v_k\cos \tilde v_k
+\tilde e_k \cosh v_k\sin\tilde v_k\right).
\label{6ch75b}
\end{equation}

\subsection{Power series of polar 6-complex numbers}

A power series of the 6-complex variable $u$ is a series of the form
\begin{equation}
a_0+a_1 u + a_2 u^2+\cdots +a_l u^l+\cdots .
\label{6ch83}
\end{equation}\index{power series!polar 6-complex}
Since
\begin{equation}
|au^l|\leq 6^{l/2} |a| |u|^l ,
\label{6ch82}
\end{equation}\index{modulus, inequalities!polar 6-complex}
the series is absolutely convergent for 
\begin{equation}
|u|<c,
\label{6ch86}
\end{equation}\index{convergence of power series!polar 6-complex}
where 
\begin{equation}
c=\lim_{l\rightarrow\infty} \frac{|a_l|}{\sqrt{6}|a_{l+1}|} .
\label{6ch87}
\end{equation}

If $a_l=\sum_{p=0}^{5}h_p a_{lp}$, where $h_0=1$, and
\begin{equation}
A_{l+}=\sum_{p=0}^{5}a_{lp},
\label{6ch88a}
\end{equation}
\begin{equation}
A_{l-}=\sum_{p=0}^{5}(-1)^p a_{lp},
\label{6ch88d}
\end{equation}
\begin{equation}
A_{lk}=\sum_{p=0}^{5}a_{lp}\cos\frac{\pi kp}{3},
\label{6ch88b}
\end{equation}
\begin{equation}
\tilde A_{lk}=\sum_{p=0}^{5}a_{lp}\sin\frac{\pi kp}{3},
\label{6ch88c}
\end{equation}
for $k=1,2$,   
the series (\ref{6ch83}) can be written as
\begin{equation}
\sum_{l=0}^\infty \left[
e_+A_{l+}v_+^l+e_-A_{l-}v_-^l+\sum_{k=1}^{2}
(e_k A_{lk}+\tilde e_k\tilde A_{lk})(e_k v_k+\tilde e_k\tilde v_k)^l 
\right].
\label{6ch89a}
\end{equation}

The series in Eq. (\ref{6ch83}) is absolutely convergent for 
\begin{equation}
|v_+|<c_+,\:
|v_-|<c_-,\:
\rho_k<c_k, k=1,2,
\label{6ch90}
\end{equation}\index{convergence, region of!polar 6-complex}
where 
\begin{equation}
c_+=\lim_{l\rightarrow\infty} \frac{|A_{l+}|}{|A_{l+1,+}|} ,\:
c_-=\lim_{l\rightarrow\infty} \frac{|A_{l-}|}{|A_{l+1,-}|} ,\:
c_k=\lim_{l\rightarrow\infty} \frac
{\left(A_{lk}^2+\tilde A_{lk}^2\right)^{1/2}}
{\left(A_{l+1,k}^2+\tilde A_{l+1,k}^2\right)^{1/2}}, \;k=1,2.
\label{6ch91}
\end{equation}

\subsection{Analytic functions of a polar 6-compex variable}

If 
$f(u)=\sum_{k=0}^{5}h_kP_k(x_0,x_1,x_2,x_3,x_4,x_{5})$,
then\index{relations between partial derivatives!polar 6-complex}
\begin{equation}
\frac{\partial P_0}{\partial x_0} 
=\frac{\partial P_1}{\partial x_1} 
=\frac{\partial P_2}{\partial x_2} 
=\frac{\partial P_3}{\partial x_3}
=\frac{\partial P_4}{\partial x_4}
=\frac{\partial P_5}{\partial x_5}, 
\label{6chh95a}
\end{equation}
\begin{equation}
\frac{\partial P_1}{\partial x_0} 
=\frac{\partial P_2}{\partial x_1} 
=\frac{\partial P_3}{\partial x_2} 
=\frac{\partial P_4}{\partial x_3}
=\frac{\partial P_5}{\partial x_4}
=\frac{\partial P_0}{\partial x_5}, 
\label{6chh95b}
\end{equation}
\begin{equation}
\frac{\partial P_2}{\partial x_0} 
=\frac{\partial P_3}{\partial x_1} 
=\frac{\partial P_4}{\partial x_2} 
=\frac{\partial P_5}{\partial x_3}
=\frac{\partial P_0}{\partial x_4}
=\frac{\partial P_1}{\partial x_5}, 
\label{6chh95c}
\end{equation}
\begin{equation}
\frac{\partial P_3}{\partial x_0} 
=\frac{\partial P_4}{\partial x_1} 
=\frac{\partial P_5}{\partial x_2} 
=\frac{\partial P_0}{\partial x_3}
=\frac{\partial P_1}{\partial x_4}
=\frac{\partial P_2}{\partial x_5}, 
\label{6chh95d}
\end{equation}
\begin{equation}
\frac{\partial P_4}{\partial x_0} 
=\frac{\partial P_5}{\partial x_1} 
=\frac{\partial P_0}{\partial x_2} 
=\frac{\partial P_1}{\partial x_3}
=\frac{\partial P_2}{\partial x_4}
=\frac{\partial P_3}{\partial x_5}, 
\label{6chh95e}
\end{equation}
\begin{equation}
\frac{\partial P_5}{\partial x_0} 
=\frac{\partial P_0}{\partial x_1} 
=\frac{\partial P_1}{\partial x_2} 
=\frac{\partial P_2}{\partial x_3}
=\frac{\partial P_3}{\partial x_4}
=\frac{\partial P_4}{\partial x_5}, 
\label{6chh95f}
\end{equation}
and\index{relations between second-order derivatives!polar 6-complex}
\begin{eqnarray}
\lefteqn{\frac{\partial^2 P_k}{\partial x_0\partial x_l}
=\frac{\partial^2 P_k}{\partial x_1\partial x_{l-1}}
=\cdots=
\frac{\partial^2 P_k}{\partial x_{[l/2]}\partial x_{l-[l/2]}}}\nonumber\\
&&=\frac{\partial^2 P_k}{\partial x_{l+1}\partial x_{5}}
=\frac{\partial^2 P_k}{\partial x_{l+2}\partial x_{4}}
=\cdots
=\frac{\partial^2 P_k}{\partial x_{l+1+[(4-l)/2]}
\partial x_{5-[(4-l)/2]}} ,
\label{6ch96}
\end{eqnarray}
for $k,l=0,...,5$.
In Eq. (\ref{6ch96}), $[a]$ denotes the integer part of $a$,
defined as $[a]\leq a<[a]+1$.
In this work, brackets larger than the regular brackets
$[\;]$ do not have the meaning of integer part.

\subsection{Integrals of polar 6-complex functions}

If $f(u)$ is an analytic 6-complex function,
then\index{integrals, path!polar 6-complex}\index{poles and residues!polar 6-complex}
\begin{equation}
\oint_\Gamma \frac{f(u)du}{u-u_0}=
2\pi f(u_0)\left[\tilde e_1 
\;{\rm int}(u_{0\xi_1\eta_1},\Gamma_{\xi_1\eta_1})
+\tilde e_2 
\;{\rm int}(u_{0\xi_2\eta_2},\Gamma_{\xi_2\eta_2})
 \right],
\label{6ch120}
\end{equation}
where
\begin{equation}
{\rm int}(M,C)=\left\{
\begin{array}{l}
1 \;\:{\rm if} \;\:M \;\:{\rm is \;\:an \;\:interior \;\:point \;\:of} \;\:C ,\\ 
0 \;\:{\rm if} \;\:M \;\:{\rm is \;\:exterior \;\:to}\:\; C ,\\
\end{array}\right.,
\label{6ch118}
\end{equation}
and $u_{0\xi_k\eta_k}$ and $\Gamma_{\xi_k\eta_k}$ are respectively the
projections of the pole $u_0$ and of 
the loop $\Gamma$ on the plane defined by the axes $\xi_k$ and $\eta_k$, $k=1,2$.

\subsection{Factorization of polar 6-complex polynomials}

A polynomial of degree $m$ of the polar 6-complex variable $u$ has the form
\begin{equation}
P_m(u)=u^m+a_1 u^{m-1}+\cdots+a_{m-1} u +a_m ,
\label{6ch125}
\end{equation}
where $a_l$, for $l=1,...,m$, are 6-complex constants.
If $a_l=\sum_{p=0}^{5}h_p a_{lp}$, and with the
notations of Eqs. (\ref{6ch88a})-(\ref{6ch88c}) applied for $l= 1, \cdots, m$, the
polynomial $P_m(u)$ can be written as 
\begin{eqnarray}
\lefteqn{P_m= 
e_+\left(v_+^m +\sum_{l=1}^{m}A_{l+}v_+^{m-l} \right)
+e_-\left(v_-^m +\sum_{l=1}^{m}A_{l-}v_-^{m-l} \right) \nonumber}\\
&&+\sum_{k=1}^{2}
\left[(e_k v_k+\tilde e_k\tilde v_k)^m+
\sum_{l=1}^m(e_k A_{lk}+\tilde e_k\tilde A_{lk})
(e_k v_k+\tilde e_k\tilde v_k)^{m-l} 
\right],
\label{6ch126a}
\end{eqnarray}\index{polynomials, canonical variables!polar 6-complex}
where the constants $A_{l+}, A_{l-}, A_{lk}, \tilde A_{lk}$ 
are real numbers.

The polynomial $P_m(u)$ can be written, as 
\begin{eqnarray}
P_m(u)=\prod_{p=1}^m (u-u_p) ,
\label{6ch128c}
\end{eqnarray}\index{polynomials, factorization!polar 6-complex}
where
\begin{eqnarray}
u_p=e_+ v_{p+}+e_-v_{p-}
+\left(e_1 v_{1p}+\tilde e_1\tilde v_{1p}\right)
+\left(e_2 v_{2p}+\tilde e_2\tilde v_{2p}\right), p=1,...,m.
\label{6ch128d}
\end{eqnarray}
The quantities $v_{p+}$,  $v_{p-}$,  
$e_k v_{kp}+\tilde e_k\tilde v_{kp}$, $p=1,...,m, k=1,2$,
are the roots of the corresponding polynomial in Eq. (\ref{6ch126a}).
The roots $v_{p+}$,  $v_{p-}$ appear in complex-conjugate pairs, and 
$v_{kp}, \tilde v_{kp}$ are real numbers.
Since all these roots may be ordered arbitrarily, the polynomial $P_m(u)$ can
be written in many different ways as a product of linear factors. 

If $P(u)=u^2-1$, the degree is $m=2$, the coefficients of the polynomial are
$a_1=0, a_2=-1$, the coefficients defined in Eqs. (\ref{6ch88a})-(\ref{6ch88c})
are $A_{2+}=-1, A_{2-}=-1, A_{21}=-1, \tilde A_{21}=0,
A_{22}=-1, \tilde A_{22}=0$. The expression of $P(u)$, Eq. (\ref{6ch126a}), is
$v_+^2-e_++v_-^2-e_- +(e_1v_1+\tilde e_1\tilde v_1)^2-e_1+
(e_2v_2+\tilde e_2\tilde v_2)^2-e_2 $. 
The factorization of $P(u)$, Eq. (\ref{6ch128c}), is
$P(u)=(u-u_1)(u-u_2)$, where the roots are
$u_1=\pm e_+\pm e_-\pm e_1\pm  e_2, u_2=-u_1$. If $e_+, e_-, e_1, e_2$ 
are expressed with the aid of Eq. (\ref{6che11}) in terms of $h_1, h_2, h_3,
h_4, h_5$, the factorizations of $P(u)$ are obtained as
\begin{eqnarray}
\lefteqn{\begin{array}{l}
u^2-1=(u+1)(u-1),\\
u^2-1=\left[u+\frac{1}{3}(1+h_1+h_2-2h_3+h_4+h_5)\right]
\left[u-\frac{1}{3}(1+h_1+h_2-2h_3+h_4+h_5)\right],\\
u^2-1=\left[u+\frac{1}{3}(1-h_1+h_2+2h_3+h_4-h_5)\right]
\left[u-\frac{1}{3}(1-h_1+h_2+2h_3+h_4-h_5)\right],\\
u^2-1=\left[u+\frac{1}{3}(2+h_1-h_2+h_3-h_4+h_5)\right]
\left[u-\frac{1}{3}(2+h_1-h_2+h_3-h_4+h_5)\right],\\
u^2-1=\left[u+\frac{1}{3}(-1+2h_2+2h_4)\right]
\left[u-\frac{1}{3}(-1+2h_2+2h_4)\right],\\
u^2-1=\left[u+\frac{1}{3}(2h_1-h_3+2h_5)\right]
\left[u-\frac{1}{3}(2h_1-h_3+2h_5)\right],\\
u^2-1=(u+h_3)(u-h_3),\\
u^2-1=\left[u+\frac{1}{3}(-2+h_1+h_2+h_3+h_4+h_5)\right]
\left[u-\frac{1}{3}(-2+h_1+h_2+h_3+h_4+h_5)\right].
\end{array}\nonumber}\\
&&
\label{ch-factorization}
\end{eqnarray}
It can be checked that 
$(\pm e_+\pm e_-\pm e_1\pm e_2)^2=e_++e_-+e_1+e_2=1$.

\subsection{Representation of polar 6-complex numbers by irreducible matrices}

If the unitary matrix which appears in the expression, 
Eq. (\ref{6ch9ee}), of the variables $\xi_+, \xi_-$, 
$\xi_1, \eta_1$, $\xi_k, \eta_k$ in
terms of $x_0, x_1, x_2, x_3, x_4, x_5$ is called $T$,
the irreducible representation of the hypercomplex number $u$ is
\begin{equation}
T U T^{-1}=\left(
\begin{array}{ccccc}
v_+     &     0     &     0   &    0   \\
0       &     v_-   &     0   &    0   \\
0       &     0     &     V_1 &    0   \\
0       &     0     &     0   &    V_2\\
\end{array}
\right),
\label{6ch129a}
\end{equation}\index{representation by irreducible matrices!polar 6-complex}
where $U$ is the matrix in Eq. (\ref{6ch24b}),
and $V_k$ are the matrices
\begin{equation}
V_k=\left(
\begin{array}{cc}
v_k           &     \tilde v_k   \\
-\tilde v_k   &     v_k          \\
\end{array}\right),\;\; k=1,2.
\label{6ch130}
\end{equation}

\section{Planar complex numbers in 6 dimensions}

\subsection{Operations with planar complex numbers in 6 dimensions}

The planar hypercomplex number $u$ in 6 dimensions 
is represented as  
\begin{equation}
u=x_0+h_1x_1+h_2x_2+h_3x_3+h_4x_4+h_5x_5. 
\label{6c1a}
\end{equation}
The multiplication rules for the bases 
$h_1, h_2, h_3, h_4, h_5 $ are 
\begin{eqnarray}
\lefteqn{h_1^2=h_2,\;h_2^2=h_4,\;h_3^2=1,\;h_4^2=-h_2,\;h_5^2=-h_4,\;
h_1h_2=h_3,\;h_1h_3=h_4,\;h_1h_4=h_5,\nonumber}\\
&&h_1h_5=-1,\;h_2h_3=h_5,\;h_2h_4=-1,\;h_2h_5=-h_1,\;h_3h_4=-h_1,\;h_3h_5=-h_2,\;\nonumber\\
&&h_4h_5=-h_3.
\label{6c1}
\end{eqnarray}\index{complex units!planar 6-complex}
The significance of the composition laws in Eq.
(\ref{6c1}) can be understood by representing the bases 
$1,h_1,h_2,h_3,h_4,h_5$
by points on a circle at the angles $\alpha_k=\pi k/6$, as shown in Fig. \ref{fig19}.
The product $h_j h_k$ will be represented by the point of the circle at the
angle $\pi (j+k)/12$, $j,k=0,1,...,5$. If $\pi\leq\pi (j+k)/12\leq 2\pi$, the
point is opposite to the basis $h_l$ of angle $\alpha_l=\pi (j+k)/6-\pi$. \\

\begin{figure}
\begin{center}
\epsfig{file=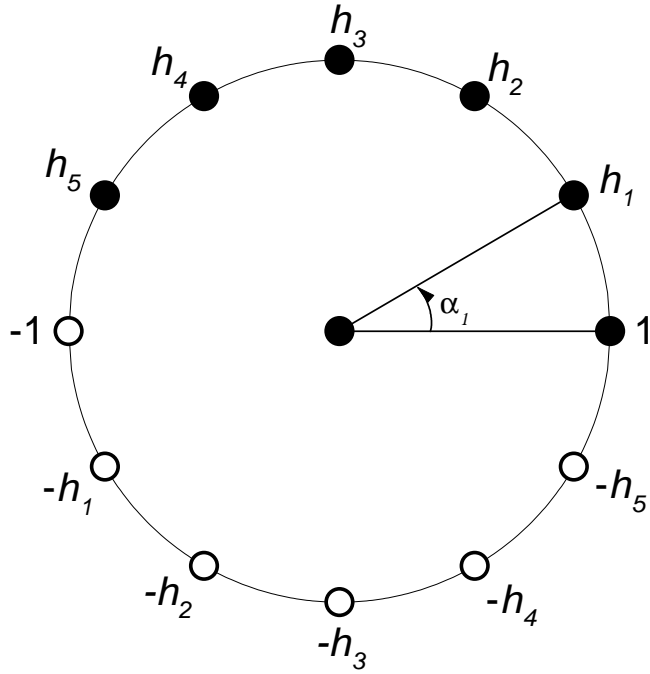,width=12cm}
\caption{Representation of the planar 
hypercomplex bases $1,h_1,h_2,h_3,h_4,h_5$
by points on a circle at the angles $\alpha_k=\pi k/6$.
The product $h_j h_k$ will be represented by the point of the circle at the
angle $\pi (j+k)/12$, $i,k=0,1,...,5$. If $\pi\leq\pi (j+k)/12\leq 2\pi$, the
point is opposite to the basis $h_l$ of angle $\alpha_l=\pi (j+k)/6-\pi$. }
\label{fig19}
\end{center}
\end{figure}

The sum of the 6-complex numbers $u$ and $u^\prime$ is
\begin{equation}
u+u^\prime=x_0+x^\prime_0+h_1(x_1+x^\prime_1)+h_1(x_2+x^\prime_2)
+h_3(x_3+x^\prime_3)+h_4(x_4+x^\prime_4)+h_5(x_5+x^\prime_5).
\label{6c2}
\end{equation}\index{sum!planar 6-complex}
The product of the numbers $u, u^\prime$ is
\begin{equation}
\begin{array}{l}
uu^\prime=x_0 x_0^\prime -x_1x_5^\prime-x_2 x_4^\prime-x_3x_3^\prime
-x_4x_2^\prime-x_5 x_1^\prime\\
+h_1(x_0 x_1^\prime+x_1x_0^\prime-x_2x_5^\prime-x_3x_4^\prime
-x_4 x_3^\prime-x_5 x_2^\prime) \\
+h_2(x_0 x_2^\prime+x_1x_1^\prime+x_2x_0^\prime-x_3x_5^\prime
-x_4 x_4^\prime-x_5 x_3^\prime) \\
+h_3(x_0 x_3^\prime+x_1x_2^\prime+x_2x_1^\prime+x_3x_0^\prime
-x_4 x_5^\prime-x_5 x_4^\prime) \\
+h_4(x_0 x_4^\prime+x_1x_3^\prime+x_2x_2^\prime+x_3x_1^\prime
+x_4 x_0^\prime-x_5 x_5^\prime) \\
+h_5(x_0 x_5^\prime+x_1x_4^\prime+x_2x_3^\prime+x_3x_2^\prime
+x_4 x_1^\prime+x_5 x_0^\prime).
\end{array}
\label{6c3}
\end{equation}\index{product!planar 6-complex}

The relation between the variables $v_1,\tilde v_1, v_2, \tilde v_2,
v_3,\tilde v_3$ and $x_0,x_1,x_2,x_3,x_4,x_5$ are
\begin{equation}
\left(
\begin{array}{c}
v_1\\
\tilde v_1\\
v_2\\
\tilde v_2\\
v_3\\
\tilde v_3\\
\end{array}\right)
=\left(
\begin{array}{cccccc}
1&\frac{\sqrt{3}}{2}&\frac{1}{2}&0&-\frac{1}{2}&-\frac{\sqrt{3}}{2}\\
0&\frac{1}{2}&\frac{\sqrt{3}}{2}&1&\frac{\sqrt{3}}{2}&\frac{1}{2}\\
1&0 &-1 &0 &1 &0 \\
0&1&0 &-1 &0 &1 \\
1&-\frac{\sqrt{3}}{2}&\frac{1}{2}&0&-\frac{1}{2}&\frac{\sqrt{3}}{2}\\
0&\frac{1}{2}&-\frac{\sqrt{3}}{2}&1&-\frac{\sqrt{3}}{2}&\frac{1}{2}\\
\end{array}
\right)
\left(
\begin{array}{c}
x_0\\
x_1\\
x_2\\
x_3\\
x_4\\
x_5
\end{array}
\right).
\label{6c9e}
\end{equation}\index{canonical variables!planar 6-complex}
The other variables are $v_4=v_3, \tilde v_4=-\tilde v_3,
v_5=v_2, \tilde v_5=-\tilde v_2, v_6=v_1, \tilde v_6=-\tilde v_1$. 
The variables $v_1, \tilde v_1, v_2, \tilde v_2, v_3, \tilde v_3$ will be
called canonical planar 6-complex variables.

\subsection{Geometric representation of planar complex numbers in 6 dimensions}

The 6-complex number $u=x_0+h_1x_1+h_2x_2+h_3x_3+h_4x_4+h_5x_5$
is represented by 
the point $A$ of coordinates $(x_0,x_1,x_2,x_3,x_4,x_5)$. 
The distance from the origin $O$ of the 6-dimensional space to the point $A$
has the expression 
\begin{equation}
d^2=x_0^2+x_1^2+x_2^2+x_3^2+x_4^2+x_5^2,
\label{6c10}
\end{equation}\index{distance!planar 6-complex}
is called modulus of the 6-complex number 
$u$, and is designated by $d=|u|$.
The modulus has the property that
\begin{equation}
|u^\prime u^{\prime\prime}|\leq \sqrt{3}|u^\prime||u^{\prime\prime}| .
\label{6c79}
\end{equation}

The exponential and trigonometric forms of the 6-complex number $u$ can be
obtained conveniently in a rotated system of axes defined by a transformation
which has the form
\begin{equation}
\left(
\begin{array}{c}
\xi_1\\
\tilde \xi_1\\
\xi_2\\
\tilde \xi_2\\
\xi_3\\
\tilde \xi_3\\
\end{array}\right)
=\left(
\begin{array}{cccccc}
\frac{1}{\sqrt{3}}&\frac{1}{2}&\frac{1}{2\sqrt{3}}&0&-\frac{1}{2\sqrt{3}}&-\frac{1}{2}\\
0&\frac{1}{2\sqrt{3}}&\frac{1}{2}&\frac{1}{\sqrt{3}}&\frac{1}{2}&\frac{1}{2\sqrt{3}}\\
\frac{1}{\sqrt{3}}&0 &-\frac{1}{\sqrt{3}} &0 &\frac{1}{\sqrt{3}} &0 \\
0&\frac{1}{\sqrt{3}}&0 &-\frac{1}{\sqrt{3}} &0 &\frac{1}{\sqrt{3}} \\
\frac{1}{\sqrt{3}}&-\frac{1}{2}&\frac{1}{2\sqrt{3}}&0&-\frac{1}{2\sqrt{3}}&\frac{1}{2}\\
0&\frac{1}{2\sqrt{3}}&-\frac{1}{2}&\frac{1}{\sqrt{3}}&-\frac{1}{2}&\frac{1}{2\sqrt{3}}\\
\end{array}
\right)
\left(
\begin{array}{c}
x_0\\
x_1\\
x_2\\
x_3\\
x_4\\
x_5
\end{array}
\right).
\label{6c9ee}
\end{equation}

The lines of the matrices in Eq. (\ref{6c9ee}) give the components
of the 6 vectors of the new basis system of axes. These vectors have unit
length and are orthogonal to each other.
The relations between the two sets of variables are
\begin{equation}
v_k= \sqrt{3}\xi_k , \tilde v_k= \sqrt{3}\eta_k, 
\label{6c12b}
\end{equation}
for $k=1,2,3$.
 
The radius $\rho_k$ and the azimuthal angle $\phi_k$ in the plane of the axes
$v_k,\tilde v_k$ are
\begin{equation}
\rho_k^2=v_k^2+\tilde v_k^2, \:\cos\phi_k=v_k/\rho_k,
\:\sin\phi_k=\tilde v_k/\rho_k, 
\label{6c19a}
\end{equation}
where $0\leq \phi_k<2\pi ,  \;k=1,2,3$,
so that there are 3 azimuthal angles.\index{azimuthal angles!planar 6-complex}
The planar angles $\psi_{k-1}$ are
\begin{equation}
\tan\psi_1=\rho_1/\rho_2,  \;\tan\psi_2=\rho_1/\rho_3, 
\label{6c19b}
\end{equation}\index{planar angles!planar 6-complex}
where $0\leq\psi_1\leq\pi/2,\;0\leq\psi_2\leq\pi/2$, 
so that there are 2 planar angles.
The amplitude of an 6-complex number $u$ is
\begin{equation}
\rho=\left(\rho_1\rho_2\rho_3\right)^{1/3}.
\label{6c50aa}
\end{equation}\index{amplitude!planar 6-complex}
It can be checked that
\begin{equation}
d^2=\frac{1}{3}(\rho_1^2+\rho_2^2+\rho_3^2).
\label{6c13}
\end{equation}\index{modulus, canonical variables!planar 6-complex}

If $u=u^\prime u^{\prime\prime}$, the parameters of the hypercomplex numbers
are related by\index{transformation of variables!planar 6-complex}
\begin{equation}
\rho_k=\rho_k^\prime\rho_k^{\prime\prime}, 
\label{6c21b}
\end{equation}
\begin{equation}
\tan\psi_k=\tan\psi_k^\prime \tan\psi_k^{\prime\prime},  
\label{6c21d}
\end{equation}
\begin{equation}
\phi_k=\phi_k^\prime+\phi_k^{\prime\prime},
\label{6c21e}
\end{equation}
\begin{equation}
v_k=v_k^\prime v_k^{\prime\prime}-\tilde v_k^\prime \tilde v_k^{\prime\prime},\;
\tilde v_k=v_k^\prime \tilde v_k^{\prime\prime}+\tilde v_k^\prime v_k^{\prime\prime},
\label{6c22}
\end{equation}
\begin{equation}
\rho=\rho^\prime\rho^{\prime\prime} ,
\label{6c24}
\end{equation}
where $k=1,2,3$.

The 6-complex planar number
$u=x_0+h_1x_1+h_2x_2+h_3x_3+h_4x_4+h_5x_5$ can be  represented by the matrix
\begin{equation}
U=\left(
\begin{array}{cccccc}
x_0      &    x_1     &    x_2   &    x_3  &  x_4   & x_5\\
-x_5     &    x_0     &    x_1   &    x_2  &  x_3   & x_4\\
-x_4     &   -x_5     &    x_0   &    x_1  &  x_2   & x_3\\
-x_3     &   -x_4     &   -x_5   &    x_0  &  x_1   & x_2\\
-x_2     &   -x_3     &   -x_4   &   -x_5  &  x_0   & x_1\\
-x_1     &   -x_2     &   -x_3   &   -x_4  &  -x_5  & x_0\\
\end{array}
\right).
\label{6c24b}
\end{equation}\index{matrix representation!planar 6-complex}
The product $u=u^\prime u^{\prime\prime}$ is
represented by the matrix multiplication $U=U^\prime U^{\prime\prime}$.

\subsection{The planar 6-dimensional cosexponential functions}

The planar cosexponential functions in 6 dimensions are
\begin{equation}
f_{6k}(y)=\sum_{p=0}^\infty (-1)^p \frac{y^{k+6p}}{(k+6p)!}, 
\label{6c29}
\end{equation}\index{cosexponential functions, planar 6-complex!definition}
for $k=0,...,5$.
The planar cosexponential functions  of even index $k$ are
even functions, $f_{6,2l}(-y)=f_{6,2l}(y)$,  
and the planar cosexponential functions of odd index 
are odd functions, $f_{6,2l+1}(-y)=-f_{6,2l+1}(y)$, $l=0,1,2$. 
\index{cosexponential functions, planar 6-complex!parity}
The exponential function of the quantity $h_k y$ is
\begin{equation}
\begin{array}{l}
e^{h_1 y}=f_{60}(y)+h_1f_{61}(y)+h_2f_{62}(y)+h_3f_{63}(y)+h_4f_{64}(y)+h_5f_{65}(y),\\
e^{h_2 y}=g_{60}(y)-g_{63}(y)+h_2\{g_{61}(y)-g_{64}(y)\}+h_4\{g_{62}(y)-g_{65}(y)\},\\
e^{h_3 y}=f_{60}(y)-f_{62}(y)+f_{64}(y)+h_3\{f_{61}(y)-f_{63}(y)+f_{65}(y)\},\\
e^{h_4 y}=g_{60}(y)+g_{63}(y)-h_2\{g_{62}(y)+g_{65}(y)\}+h_4\{g_{61}(y)+g_{64}(y)\},\\
e^{h_5 y}=f_{60}(y)+h_1f_{65}(y)-h_2f_{64}(y)+h_3f_{63}(y)-h_4f_{62}(y)+h_5f_{61}(y).\\
\end{array}
\label{6c28b}
\end{equation}\index{exponential, expressions!planar 6-complex}
The relations for $h_2$ and $h_4$ can be written equivalently as 
$e^{h_2 y}=f_{30}+h_2f_{31}+h_4 f_{32}, e^{h_4 y}=g_{30}-h_2f_{32}+h_4 g_{31}$,
and the relation for $h_3$ can be written as $e^{h_3 y}=f_{20}+h_3f_{21}$,
which is the same as $e^{h_3 y}=\cos y+h_3\sin y$. 

The planar 6-dimensional cosexponential functions $f_{6k}(y)$ are related to the
polar 6-dimensional cosexponential function $g_{6k}(y)$ 
by the relations 
\begin{equation}
f_{6k}(y)=e^{-i\pi k/6}g_{6k}\left(e^{i\pi/6}y\right), 
\label{6c30a}
\end{equation}
for $k=0,...,5$.
The planar 6-dimensional cosexponential functions $f_{6k}(y)$ are related to
the polar 6-dimensional cosexponential function $g_{6k}(y)$ also
by the relations 
\begin{equation}
f_{6k}(y)=e^{-i\pi k/2}g_{6k}(iy), 
\label{6c30ax}
\end{equation}
for $k=0,...,5$.
The expressions of the planar 6-dimensional cosexponential functions are
\begin{equation}
\begin{array}{l}
f_{60}(y)=\frac{1}{3}\cos y+\frac{2}{3}\cosh\frac{\sqrt{3}}{2}y\cos\frac{y}{2},\\
f_{61}(y)=\frac{1}{3}\sin y+\frac{\sqrt{3}}{3}\sinh\frac{\sqrt{3}}{2}y\cos\frac{y}{2}
+\frac{1}{3}\cosh\frac{\sqrt{3}}{2}y\sin\frac{y}{2},\\
f_{62}(y)=-\frac{1}{3}\cos y+\frac{1}{3}\cosh\frac{\sqrt{3}}{2}y\cos\frac{y}{2}
+\frac{\sqrt{3}}{3}\sinh\frac{\sqrt{3}}{2}y\sin\frac{y}{2},\\
f_{63}(y)=-\frac{1}{3}\sin y+\frac{2}{3}\cosh\frac{\sqrt{3}}{2}y\sin\frac{y}{2},\\
f_{64}(y)=\frac{1}{3}\cos y-\frac{1}{3}\cosh\frac{\sqrt{3}}{2}y\cos\frac{y}{2}
+\frac{\sqrt{3}}{3}\sinh\frac{\sqrt{3}}{2}y\sin\frac{y}{2},\\
f_{65}(y)=\frac{1}{3}\sin y-\frac{\sqrt{3}}{3}\sinh\frac{\sqrt{3}}{2}y\cos\frac{y}{2}
+\frac{1}{3}\cosh\frac{\sqrt{3}}{2}y\sin\frac{y}{2}.\\
\end{array}
\label{6c30x}
\end{equation}\index{cosexponential functions, planar 6-complex!expressions}
The planar 6-dimensional cosexponential functions can be written as
\begin{equation}
f_{6k}(y)=\frac{1}{6}\sum_{l=1}^{6}\exp\left[y\cos\left(\frac{\pi (2l-1)}{6}\right)
\right]
\cos\left[y\sin\left(\frac{\pi (2l-1)}{6}\right)-\frac{\pi (2l-1)k}{6}\right], 
\label{6c30}
\end{equation}
for $k=0,...,5$.
The graphs of the planar 6-dimensional cosexponential functions are shown in
Fig. \ref{fig20}.

\begin{figure}
\begin{center}
\epsfig{file=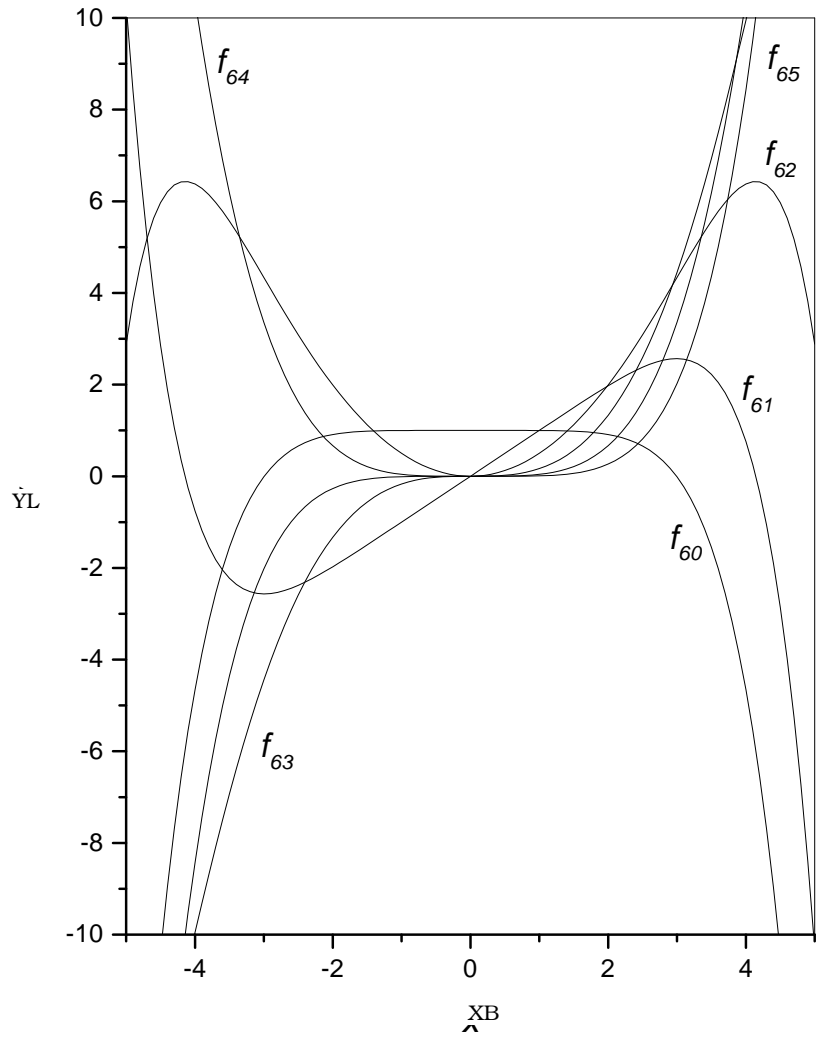,width=12cm}
\caption{Planar cosexponential functions 
$f_{60}, f_{61},f_{62}, f_{63},f_{64}, f_{65}$.}
\label{fig20}
\end{center}
\end{figure}

It can be checked that
\begin{equation}
\sum_{k=0}^{5}f_{6k}^2(y)=\frac{1}{3}+\frac{2}{3}\cosh\sqrt{3}y.
\label{6c34a}
\end{equation}

\newpage
\setlength{\oddsidemargin}{-1.1cm}
The addition theorems for the planar 6-dimensional cosexponential functions are
\begin{eqnarray}
\lefteqn{\begin{array}{l}
g_{60}(y+z)=g_{60}(y)g_{60}(z)-g_{61}(y)g_{65}(z)-g_{62}(y)g_{64}(z)
-g_{63}(y)g_{63}(z)-g_{64}(y)g_{62}(z)-g_{65}(y)g_{61}(z) ,\\
g_{61}(y+z)=g_{60}(y)g_{61}(z)+g_{61}(y)g_{60}(z)-g_{62}(y)g_{65}(z)
-g_{63}(y)g_{64}(z)-g_{64}(y)g_{63}(z)-g_{65}(y)g_{62}(z) ,\\
g_{62}(y+z)=g_{60}(y)g_{62}(z)+g_{61}(y)g_{61}(z)+g_{62}(y)g_{60}(z)
-g_{63}(y)g_{65}(z)-g_{64}(y)g_{64}(z)-g_{65}(y)g_{63}(z) ,\\
g_{63}(y+z)=g_{60}(y)g_{63}(z)+g_{61}(y)g_{62}(z)+g_{62}(y)g_{61}(z)
+g_{63}(y)g_{60}(z)-g_{64}(y)g_{65}(z)-g_{65}(y)g_{64}(z) ,\\
g_{64}(y+z)=g_{60}(y)g_{64}(z)+g_{61}(y)g_{63}(z)+g_{62}(y)g_{62}(z)
+g_{63}(y)g_{61}(z)+g_{64}(y)g_{60}(z)-g_{65}(y)g_{65}(z) ,\\
g_{65}(y+z)=g_{60}(y)g_{65}(z)+g_{61}(y)g_{64}(z)+g_{62}(y)g_{63}(z)
+g_{63}(y)g_{62}(z)+g_{64}(y)g_{61}(z)+g_{65}(y)g_{60}(z) .
\end{array}\nonumber}\\
&&
\label{6c35a}
\end{eqnarray}\index{cosexponential functions, planar 6-complex!addition theorems}
It can be shown that
\begin{equation}
\begin{array}{l}
\{f_{60}(y)+h_1f_{61}(y)+h_2f_{62}(y)+h_3f_{63}(y)+h_4f_{64}(y)+h_5f_{65}(y)\}^l\\
\hspace*{0.5cm}=f_{60}(ly)+h_1f_{61}(ly)+h_2f_{62}(ly)+h_3f_{63}(ly)+h_4f_{64}(ly)+h_5f_{65}(ly),\\
\{g_{60}(y)-g_{63}(y)+h_2\{g_{61}(y)-g_{64}(y)\}+h_4\{g_{62}(y)-g_{65}(y)\}\}^l\\
\hspace*{0.5cm}=g_{60}(ly)-g_{63}(ly)+h_2\{g_{61}(ly)-g_{64}(ly)\}+h_4\{g_{62}(ly)-g_{65}(ly)\},\\
\{f_{60}(y)-f_{62}(y)+f_{64}(y)+h_3\{f_{61}(y)-f_{63}(y)+f_{65}(y)\}\}^l\\
\hspace*{0.5cm}=f_{60}(ly)-f_{62}(ly)+f_{64}(ly)+h_3\{f_{61}(ly)-f_{63}(ly)+f_{65}(ly)\},\\
\{g_{60}(y)+g_{63}(y)-h_2\{g_{62}(y)+g_{65}(y)\}+h_4\{g_{61}(y)+g_{64}(y)\}\}^l\\
\hspace*{0.5cm}=g_{60}(ly)+g_{63}(ly)-h_2\{g_{62}(ly)+g_{65}(ly)\}+h_4\{g_{61}(ly)+g_{64}(ly)\},\\
\{f_{60}(y)+h_1f_{65}(y)-h_2f_{64}(y)+h_3f_{63}(y)-h_4f_{62}(y)+h_5f_{61}(y)\}^l\\
\hspace*{0.5cm}=f_{60}(ly)+h_1f_{65}(ly)-h_2f_{64}(ly)+h_3f_{63}(ly)-h_4f_{62}(ly)+h_5f_{61}(ly).\\
\end{array}
\label{6c37b}
\end{equation}

The derivatives of the planar cosexponential functions
are related by
\begin{equation}
\frac{df_{60}}{du}=-f_{65}, \:
\frac{df_{61}}{du}=f_{60}, \:
\frac{df_{62}}{du}=f_{61}, \:
\frac{df_{63}}{du}=f_{62}, \:
\frac{df_{64}}{du}=f_{63}, \:
\frac{df_{65}}{du}=f_{64}.
\label{6c45}
\end{equation}\index{cosexponential functions, planar 6-complex!differential equations}

\subsection{Exponential and trigonometric forms of planar 6-complex numbers}

The exponential and trigonometric forms of planar 6-complex
numbers can be expressed with the aid of the hypercomplex bases 
\begin{equation}
\left(
\begin{array}{c}
e_1\\
\tilde e_1\\
e_2\\
\tilde e_2\\
e_3\\
\tilde e_3\\
\end{array}\right)
=\left(
\begin{array}{cccccc}
\frac{1}{3}&\frac{\sqrt{3}}{6}&\frac{1}{6}&0&-\frac{1}{6}&-\frac{\sqrt{3}}{6}\\
0&\frac{1}{6}&\frac{\sqrt{3}}{6}&\frac{1}{3}&\frac{\sqrt{3}}{6}&\frac{1}{6}\\
\frac{1}{3}&0 &-\frac{1}{3} &0 &\frac{1}{3} &0 \\
0&\frac{1}{3}&0 &-\frac{1}{3} &0 &\frac{1}{3} \\
\frac{1}{3}&-\frac{\sqrt{3}}{6}&\frac{1}{6}&0&-\frac{1}{6}&\frac{\sqrt{3}}{6}\\
0&\frac{1}{6}&-\frac{\sqrt{3}}{6}&\frac{1}{3}&-\frac{\sqrt{3}}{6}&\frac{1}{6}\\
\end{array}
\right)
\left(
\begin{array}{c}
1\\
h_1\\
h_2\\
h_3\\
h_4\\
h_5
\end{array}
\right).
\label{6ce11}
\end{equation}\index{canonical base!planar 6-complex}

The multiplication relations for the bases $e_k, \tilde e_k$ are
\begin{eqnarray}
e_k^2=e_k, \tilde e_k^2=-e_k, e_k \tilde e_k=\tilde e_k , e_ke_l=0, e_k\tilde
e_l=0, \tilde e_k\tilde e_l=0, \;k,l=1,2,3, \;k\not=l.
\label{6ce12a}
\end{eqnarray}
The moduli of the bases $e_k, \tilde e_k$ are
\begin{equation}
|e_k|=\sqrt{\frac{1}{3}}, |\tilde e_k|=\sqrt{\frac{1}{3}}, 
\label{6ce12c}
\end{equation}
\newpage
\setlength{\oddsidemargin}{0.9cm}
for $k=1,2,3$.
It can be shown that
\begin{eqnarray}
x_0+h_1x_1+h_2x_2+h_3x_3+h_4x_4+h_5x_5
= \sum_{k=1}^3 (e_k v_k+\tilde e_k \tilde v_k).
\label{6ce13a}
\end{eqnarray}\index{canonical form!planar 6-complex}
The ensemble $e_1, \tilde e_1, e_2, \tilde e_2, e_3, \tilde e_3$ will be called the
canonical planar 6-complex base, and Eq. (\ref{6ce13a}) gives the canonical
form of the planar 6-complex number.

The exponential form of the 6-complex number $u$ is
\begin{eqnarray}
\lefteqn{u=\rho\exp\left\{\frac{1}{3}(h_2-h_4)\ln\tan\psi_1
+\frac{1}{6}(\sqrt{3}h_1-h_2+h_4-\sqrt{3}h_5)\ln\tan\psi_2\right.\nonumber}\\
&&\left.+\tilde e_1\phi_1+\tilde e_2\phi_2+\tilde e_3\phi_3\right\}.
\label{6c50a}
\end{eqnarray}\index{exponential form!planar 6-complex}
The
trigonometric form of the 6-complex number $u$ is
\begin{eqnarray}
\lefteqn{u=d
\sqrt{3}
\left(1+\frac{1}{\tan^2\psi_1}+\frac{1}{\tan^2\psi_2}\right)^{-1/2}\nonumber}\\
&&\left(e_1+\frac{e_2}{\tan\psi_1}+\frac{e_3}{\tan\psi_2}\right)
\exp\left(\tilde e_1\phi_1+\tilde e_2\phi_2+\tilde e_3\phi_3\right).
\label{6c52a}
\end{eqnarray}\index{trigonometric form!planar 6-complex}
The modulus $d$ and the amplitude $\rho$ are related by
\begin{eqnarray}
d=\rho \frac{2^{1/3}}{\sqrt6}
\left(\tan\psi_1\tan\psi_2\right)^{1/3}
\left(1+\frac{1}{\tan^2\psi_1}+\frac{1}{\tan^2\psi_2}\right)^{1/2}.
\label{6c53a}
\end{eqnarray}

\subsection{Elementary functions of a planar 6-complex variable}

The logarithm and power functions of the 6-complex number $u$ exist for all
$x_0,...,x_5$ and are
\begin{eqnarray}
\lefteqn{\ln u=\ln \rho+
\frac{1}{3}(h_2-h_4)\ln\tan\psi_1
+\frac{1}{6}(\sqrt{3}h_1-h_2+h_4-\sqrt{3}h_5)\ln\tan\psi_2\nonumber}\\
&&+\tilde e_1\phi_1+\tilde e_2\phi_2+\tilde e_3\phi_3 ,
\label{6c56a}
\end{eqnarray}\index{logarithm!planar 6-complex}
\begin{equation}
u^m=\sum_{k=1}^{3}
\rho_k^m(e_k\cos m\phi_k+\tilde e_k\sin m\phi_k).
\label{6c59a}
\end{equation}\index{power function!planar 6-complex}

The exponential of the 6-complex variable $u$ is
\begin{eqnarray}
e^u= 
\sum_{k=1}^{3}e^{v_k}\left(e_k \cos \tilde v_k+\tilde e_k \sin\tilde
v_k\right).
\label{6c73a}
\end{eqnarray}\index{exponential, expressions!planar 6-complex}

The trigonometric functions of the
6-complex variable $u$ are
\begin{equation}
\cos u=\sum_{k=1}^{3}\left(e_k \cos v_k\cosh \tilde v_k
-\tilde e_k \sin v_k\sinh\tilde v_k\right),
\label{6c74a}
\end{equation}\index{trigonometric functions, expressions!planar 6-complex}
\begin{equation}
\sin u= 
\sum_{k=1}^{3}\left(e_k \sin v_k\cosh \tilde v_k
+\tilde e_k \cos v_k\sinh\tilde v_k\right).
\label{6c74b}
\end{equation}
The hyperbolic functions of the
6-complex variable $u$ are
\begin{equation}
\cosh u=
\sum_{k=1}^{3}\left(e_k \cosh v_k\cos \tilde v_k
+\tilde e_k \sinh v_k\sin\tilde v_k\right),
\label{6c75a}
\end{equation}\index{hyperbolic functions, expressions!planar 6-complex}
\begin{equation}
\sinh u=
\sum_{k=1}^{3}\left(e_k \sinh v_k\cos \tilde v_k
+\tilde e_k \cosh v_k\sin\tilde v_k\right).
\label{6c75b}
\end{equation}

\subsection{Power series of planar 6-complex numbers}

A power series of the 6-complex variable $u$ is a series of the form
\begin{equation}
a_0+a_1 u + a_2 u^2+\cdots +a_l u^l+\cdots .
\label{6c83}
\end{equation}\index{power series!planar 6-complex}
Since
\begin{equation}
|au^l|\leq 3^{l/2} |a| |u|^l ,
\label{6c82}
\end{equation}\index{modulus, inequalities!planar 6-complex}
the series is absolutely convergent for 
\begin{equation}
|u|<c,
\label{6c86}
\end{equation}\index{convergence of power series!planar 6-complex}
where 
\begin{equation}
c=\lim_{l\rightarrow\infty} \frac{|a_l|}{\sqrt{3}|a_{l+1}|} .
\label{6c87}
\end{equation}

If $a_l=\sum_{p=0}^{5} h_p a_{lp}$, and
\begin{equation}
A_{lk}=\sum_{p=0}^{5} a_{lp}\cos\frac{\pi (2k-1)p}{6},
\label{6c88b}
\end{equation}
\begin{equation}
\tilde A_{lk}=\sum_{p=0}^{5} a_{lp}\sin\frac{\pi (2k-1)p}{6},
\label{6c88c}
\end{equation}
where $k=1,2,3$, the series (\ref{6c83}) can be written as
\begin{equation}
\sum_{l=0}^\infty \left[
\sum_{k=1}^{3}
(e_k A_{lk}+\tilde e_k\tilde A_{lk})(e_k v_k+\tilde e_k\tilde v_k)^l 
\right].
\label{6c89a}
\end{equation}
The series is absolutely convergent for   
\begin{equation}
\rho_k<c_k, k=1,2,3,
\label{6c90}
\end{equation}\index{convergence, region of!planar 6-complex}
where 
\begin{equation}
c_k=\lim_{l\rightarrow\infty} \frac
{\left[A_{lk}^2+\tilde A_{lk}^2\right]^{1/2}}
{\left[A_{l+1,k}^2+\tilde A_{l+1,k}^2\right]^{1/2}} .
\label{6c91}
\end{equation}

\subsection{Analytic functions of a planar 6-complex variable}

If 
$f(u)=\sum_{k=0}^{5}h_kP_k(x_0,...,x_5)$,\index{functions, real components!planar 6-complex}
then\index{relations between partial derivatives!planar 6-complex}
\begin{equation}
\frac{\partial P_0}{\partial x_0} 
=\frac{\partial P_1}{\partial x_1} 
=\frac{\partial P_2}{\partial x_2} 
=\frac{\partial P_3}{\partial x_3}
=\frac{\partial P_4}{\partial x_4}
=\frac{\partial P_5}{\partial x_5}, 
\label{6ch95a}
\end{equation}
\begin{equation}
\frac{\partial P_1}{\partial x_0} 
=\frac{\partial P_2}{\partial x_1} 
=\frac{\partial P_3}{\partial x_2} 
=\frac{\partial P_4}{\partial x_3}
=\frac{\partial P_5}{\partial x_4}
=-\frac{\partial P_0}{\partial x_5}, 
\label{6ch95b}
\end{equation}
\begin{equation}
\frac{\partial P_2}{\partial x_0} 
=\frac{\partial P_3}{\partial x_1} 
=\frac{\partial P_4}{\partial x_2} 
=\frac{\partial P_5}{\partial x_3}
=-\frac{\partial P_0}{\partial x_4}
=-\frac{\partial P_1}{\partial x_5}, 
\label{6ch95c}
\end{equation}
\begin{equation}
\frac{\partial P_3}{\partial x_0} 
=\frac{\partial P_4}{\partial x_1} 
=\frac{\partial P_5}{\partial x_2} 
=-\frac{\partial P_0}{\partial x_3}
=-\frac{\partial P_1}{\partial x_4}
=-\frac{\partial P_2}{\partial x_5}, 
\label{6ch95d}
\end{equation}
\begin{equation}
\frac{\partial P_4}{\partial x_0} 
=\frac{\partial P_5}{\partial x_1} 
=-\frac{\partial P_0}{\partial x_2} 
=-\frac{\partial P_1}{\partial x_3}
=-\frac{\partial P_2}{\partial x_4}
=-\frac{\partial P_3}{\partial x_5}, 
\label{6ch95e}
\end{equation}
\begin{equation}
\frac{\partial P_5}{\partial x_0} 
=-\frac{\partial P_0}{\partial x_1} 
=-\frac{\partial P_1}{\partial x_2} 
=-\frac{\partial P_2}{\partial x_3}
=-\frac{\partial P_3}{\partial x_4}
=-\frac{\partial P_4}{\partial x_5}, 
\label{6ch95f}
\end{equation}
and\index{relations between second-order derivatives!planar 6-complex}
\begin{eqnarray}
\lefteqn{\frac{\partial^2 P_k}{\partial x_0\partial x_l}
=\frac{\partial^2 P_k}{\partial x_1\partial x_{l-1}}
=\cdots=
\frac{\partial^2 P_k}{\partial x_{[l/2]}\partial x_{l-[l/2]}}}\nonumber\\
&&=-\frac{\partial^2 P_k}{\partial x_{l+1}\partial x_5}
=-\frac{\partial^2 P_k}{\partial x_{l+2}\partial x_4}
=\cdots
=-\frac{\partial^2 P_k}{\partial x_{l+1+[(4-l)/2]}
\partial x_{5-[(4-l)/2]}} .
\label{6c96}
\end{eqnarray}

\subsection{Integrals of planar 6-complex functions}

If $f(u)$ is an analytic 6-complex function,
then\index{integrals, path!planar 6-complex}\index{poles and residues!planar 6-complex}
\begin{equation}
\oint_\Gamma \frac{f(u)du}{u-u_0}=
2\pi f(u_0)\left\{
\tilde e_1 \;{\rm int}(u_{0\xi_1\eta_1},\Gamma_{\xi_1\eta_1})
+\tilde e_2 \;{\rm int}(u_{0\xi_2\eta_2},\Gamma_{\xi_2\eta_2})
+\tilde e_3 \;{\rm int}(u_{0\xi_3\eta_3},\Gamma_{\xi_3\eta_3})\right\},
\label{6c120}
\end{equation}
where $u_{0\xi_k\eta_k}$ and $\Gamma_{\xi_k\eta_k}$ are respectively the
projections of the point $u_0$ and of 
the loop $\Gamma$ on the plane defined by the axes $\xi_k$ and $\eta_k$,
$k=1,2,3$. 

\subsection{Factorization of planar 6-complex polynomials}

A polynomial of degree $m$ of the planar 6-complex variable $u$ has the form
\begin{equation}
P_m(u)=u^m+a_1 u^{m-1}+\cdots+a_{m-1} u +a_m ,
\label{6c125}
\end{equation}
where $a_l$, for $l=1,...,m$, are 6-complex constants.
If $a_l=\sum_{p=0}^{5}h_p a_{lp}$, and with the
notations of Eqs. (\ref{6c88b})-(\ref{6c88c}) applied for $l= 1, \cdots, m$, the
polynomial $P_m(u)$ can be written as 
\begin{eqnarray}
P_m= 
\sum_{k=1}^{3}
\left[(e_k v_k+\tilde e_k\tilde v_k)^m+
\sum_{l=1}^m(e_k A_{lk}+\tilde e_k\tilde A_{lk})
(e_k v_k+\tilde e_k\tilde v_k)^{m-l} 
\right],
\label{6c126a}
\end{eqnarray}\index{polynomial, canonical variables!planar 6-complex}
\newpage
\setlength{\oddsidemargin}{0cm}
where the constants $A_{lk}, \tilde A_{lk}$ are real numbers.

The polynomial $P_m(u)$ can be written as a product of factors 
\begin{eqnarray}
P_m(u)=\prod_{p=1}^m (u-u_p) ,
\label{6c128c}
\end{eqnarray}\index{polynomial, factorization!planar 6-complex}
where
\begin{eqnarray}
u_p=\sum_{k=1}^{3}\left(e_k v_{kp}+\tilde e_k\tilde v_{kp}\right), 
\label{6c128d}
\end{eqnarray}
for $p=1,...,m$.
The quantities  
$e_k v_{kp}+\tilde e_k\tilde v_{kp}$, $p=1,...,m, k=1,2,3$,
are the roots of the corresponding polynomial in Eq. (\ref{6c126a}) and are
real numbers.
Since these roots may be ordered arbitrarily, the polynomial $P_m(u)$ can be
written in many different ways as a product of linear factors. 

If $P(u)=u^2+1$, the degree is $m=2$, the coefficients of the polynomial are
$a_1=0, a_2=1$, the coefficients defined in Eqs. (\ref{6c88b})-(\ref{6c88c})
are $A_{21}=1, \tilde A_{21}=0,
A_{22}=1, \tilde A_{22}=0, A_{23}=1, \tilde A_{23}=0
$. The expression, Eq. (\ref{6c126a}), is
P(u)=$(e_1v_1+\tilde e_1\tilde v_1)^2+e_1+
(e_2v_2+\tilde e_2\tilde v_2)^2+e_2+(e_3v_3+\tilde e_3\tilde v_3)^2+e_3 $. 
The factorization of $P(u)$, Eq. (\ref{6c128c}), is
$P(u)=(u-u_1)(u-u_2)$, where the roots are
$u_1=\pm \tilde e_1 \pm \tilde e_2\pm \tilde e_3, u_2=-u_1$. 
If $\tilde e_1, \tilde e_2, \tilde e_3$ 
are expressed with the aid of Eq. (\ref{6ce11}) in terms of $h_1, h_2, h_3,
h_4, h_5$, the factorizations of $P(u)$ are obtained as
\begin{eqnarray}
\lefteqn{\begin{array}{l}
u^2+1=\left[u+\frac{1}{3}(2h_1+h_3+2h_5)\right]
\left[u-\frac{1}{3}(2h_1+h_3+2h_5)\right],\\
u^2+1=\left[u+\frac{1}{3}(h_1+\sqrt{3}h_2-h_3+\sqrt{3}h_4+h_5)\right]
\left[u-\frac{1}{3}(h_1+\sqrt{3}h_2-h_3+\sqrt{3}h_4+h_5)\right],\\
u^2+1=(u+h_3)(u-h_3),\\
u^2+1=\left[u+\frac{1}{3}(-h_1+\sqrt{3}h_2+h_3+\sqrt{3}h_4-h_5)\right]
\left[u-\frac{1}{3}(-h_1+\sqrt{3}h_2+h_3+\sqrt{3}h_4-h_5)\right].
\end{array}\nonumber}\\
&&
\label{c-factorization}
\end{eqnarray}
It can be checked that 
$(\pm \tilde e_1\pm \tilde e_2+\pm \tilde e_3)^2=-e_1-e_2-e_3=-1$.

\subsection{Representation of planar 6-complex numbers by irreducible matrices}

If the unitary matrix written in Eq. (\ref{6c9ee}) is called $T$,
the matric $T U T^{-1}$ provides an irreducible representation 
\cite{4} of the planar
hypercomplex number $u$, 
\begin{equation}
T U T^{-1}=\left(
\begin{array}{ccc}
V_1      &     0   &    0   \\
0        &     V_2 &    0   \\
0        &     0   &    V_3 \\
\end{array}
\right),
\label{6c129}
\end{equation}\index{representation by irreducible matrices!planar 6-complex}
where $U$ is the matrix in Eq. (\ref{6c24b}) used to represent the 6-complex
number $u$, and the matrices $V_k$ are
\begin{equation}
V_k=\left(
\begin{array}{cc}
v_k           &     \tilde v_k   \\
-\tilde v_k   &     v_k          \\
\end{array}\right),
\label{6c130}
\end{equation}
for $ k=1,2,3$.
\newpage
\setlength{\oddsidemargin}{0.9cm}

\chapter{Commutative Complex Numbers in $n$ Dimensions}

Two systems of complex numbers in $n$ dimensions are described in 
this chapter,
for which the multiplication is associative and commutative, which 
can be written in exponential and trigonometric forms, and for which
the concepts of analytic n-complex 
function,  contour integration and residue can be defined.
The n-complex numbers introduced in this chapter have 
the form $u=x_0+h_1x_1+h_2x_2+\cdots+h_{n-1}x_{n-1}$, the variables 
$x_0,...,x_{n-1}$ being real numbers. 

The multiplication rules for the complex
units $h_1,...,h_{n-1}$ discussed in Sec. 6.1 
are $h_j h_k =h_{j+k}$ if $0\leq j+k\leq n-1$, and $h_jh_k=h_{j+k-n}$ if
$n\leq j+k\leq 2n-2$.
The product of two n-complex numbers is equal to zero if both numbers are
equal to zero, or if the numbers belong to certain n-dimensional hyperplanes
described further in this chapter. 
If the n-complex number $u=x_0+h_1x_1+h_2x_2+\cdots+h_{n-1}x_{n-1}$ is
represented by the point $A$ of coordinates $x_0,x_1,...,x_{n-1}$, 
the position of the point $A$ can be described, in an even number 
of dimensions, by the modulus $d=(x_0^2+x_1^2+\cdots+x_{n-1}^2)^{1/2}$, 
by $n/2-1$ azimuthal angles $\phi_k$, by $n/2-2$
planar angles $\psi_{k-1}$, and by 2 polar angles $\theta_+,\theta_-$. In an
odd number of dimensions, the position of the point $A$ is described by $d$, by
$(n-1)/2$ azimuthal angles $\phi_k$, by $(n-3)/2$ planar angles $\psi_{k-1}$,
and by 1 polar angle $\theta_+$. 
An amplitude $\rho$ can be defined for
even $n$ as $\rho^n=v_+v_-\rho_1^2\cdots\rho_{n/2-1}^2$,
and for odd $n$ as $\rho^n=v_+\rho_1^2\cdots\rho_{(n-1)/2}^2$,
where $v_+=x_0+\cdots+x_{n-1}, v_-=x_0-x_1+\cdots+x_{n-2}-x_{n-1}$,
and $\rho_k$ are radii in orthogonal two-dimensional planes defined further in
Sec. 6.1. The amplitude $\rho$, the variables $v_+, v_-$, the radii
$\rho_k$, the variables $(1/\sqrt{2})\tan\theta_+,(1/\sqrt{2})\tan\theta_-,
\tan\psi_{k-1}$ are multiplicative, and the azimuthal angles $\phi_k$ are
additive upon the multiplication of n-complex numbers.
Because of the role of the axis $v_+$ and, in
an even number of dimensions, of the axis $v_-$, in the description of the
position of the point $A$ with the aid of the polar angle $\theta_+$ and, in an
even number of dimensions, of the polar angle $\theta_-$, the hypercomplex
numbers studied in Sec. 6.1 will be called polar n-complex number, to
distinguish them from the planar n-complex numbers, which exist in an even
number of dimensions. 

The exponential function of a polar n-complex number can be expanded in terms of
the polar n-dimensional cosexponential functions
$g_{nk}(y)=\sum_{p=0}^\infty y^{k+pn}/(k+pn)!$, $k=0,1,...,n-1$.
It is shown that
$g_{nk}(y) = \frac{1}{n}\sum_{l=0}^{n-1}$$\exp\left\{y\cos\left(2\pi l/n\right)
\right\} $$\cos\left\{y\sin\left(2\pi l/n\right)-2\pi kl/n\right\}$, 
$k=0,1,...,n-1$. Addition theorems and other relations
are obtained for the polar n-dimensional cosexponential functions.

The exponential form of a polar n-complex number, which in an even number of
dimensions $n$ can be defined for
$x_0+\cdots+x_{n-1}>0, x_0-x_1+\cdots+x_{n-2}-x_{n-1}>0$, is\\
$u=\rho \exp\left\{\sum_{p=1}^{n-1}h_p\left[
(1/n)\ln\sqrt{2}/\tan\theta_+
+((-1)^p/n)\ln\sqrt{2}/\tan\theta_-\right.\right.$\\
$\left.\left.-(2/n)\sum_{k=2}^{n/2-1}
\cos\left(2\pi kp/n\right)\ln\tan\psi_{k-1}
\right]\right\}$
 $\exp\left(\sum_{k=1}^{n/2-1}\tilde e_k\phi_k\right)$,
where $\tilde e_k=(2/n)\sum_{p=1}^{n-1}h_p\sin(2\pi pk/n)$. 
In an odd number of dimensions $n$, the exponential form exists for
$x_0+\cdots+x_{n-1}>0$, and is
$u=\rho \exp\left\{\sum_{p=1}^{n-1}h_p\left[
(1/n)\ln\sqrt{2}/\tan\theta_+
-(2/n)\sum_{k=2}^{(n-1)/2}
\cos\left(2\pi kp/n\right)\ln\tan\psi_{k-1}
\right]\right\}$\\
 $\exp\left(\sum_{k=1}^{(n-1)/2}\tilde e_k\phi_k\right)$. 
A trigonometric form also exists for an n-complex number $u$,
when $u$ is written as the product of the modulus $d$, of a factor depending on
the polar and planar angles $\theta_+, \theta_-, \psi_{k-1}$ and of an
exponential factor depending on the azimuthal angles $\phi_k$.

Expressions are given for the elementary functions of polar n-complex variable.
The functions $f(u)$ of n-complex variable which are defined by power series
have derivatives independent of the direction of approach to the point under
consideration. If the n-complex function $f(u)$ 
of the n-complex variable $u$ is written in terms of 
the real functions $P_k(x_0,...,x_{n-1}), k=0,...,n-1$, then
relations of equality  
exist between partial derivatives of the functions $P_k$. 
The integral $\int_A^B f(u) du$ of an n-complex
function between two points $A,B$ is independent of the path connecting $A,B$,
in regions where $f$ is regular.
If $f(u)$ is an analytic n-complex function, then 
$\oint_\Gamma f(u)du/(u-u_0)$
$=2\pi f(u_0)\sum_{k=1}^{[(n-1)/2]}\tilde e_k$ 
$\;{\rm int}(u_{0\xi_k\eta_k},\Gamma_{\xi_k\eta_k})$,
where the functional ${\rm int}$ takes the values 0 or 1 depending on the
relation between $u_{0\xi_k\eta_k}$ and $\Gamma_{\xi_k\eta_k}$, which are
respectively the projections of the point $u_0$ and of 
the loop $\Gamma$ on the plane defined by the orthogonal axes $\xi_k$ and
$\eta_k$, as expained further in this work.

A polar n-complex polynomial can be written as a 
product of linear or quadratic factors, although the factorization may not be
unique. 

Particular cases for $n=2,3,4,5,6$  of the polar
n-complex numbers described in Sec. 6.1 have been studied in previous chapters.


The multiplication rules for the complex
units $h_1,...,h_{n-1}$ discussed in Sec. 6.2 are 
$h_j h_k =h_{j+k}$ if $0\leq j+k\leq n-1$, and $h_jh_k=-h_{j+k-n}$ if
$n\leq j+k\leq 2n-2$, where $h_0=1$.
The product of two n-complex numbers is equal to zero if both numbers are
equal to zero, or if the numbers belong to certain n-dimensional hyperplanes
described further in Sec. 6.2. 
If the n-complex number $u=x_0+h_1x_1+h_2x_2+\cdots+h_{n-1}x_{n-1}$ is
represented by the point $A$ of coordinates $x_0,x_1,...,x_{n-1}$, 
the position of the point $A$ can be described, in an even number 
of dimensions, by the modulus $d=(x_0^2+x_1^2+\cdots+x_{n-1}^2)^{1/2}$, 
by $n/2$ azimuthal angles $\phi_k$ and by $n/2-1$
planar angles $\psi_{k-1}$. 
An amplitude $\rho$ can be defined 
as $\rho^n=\rho_1^2\cdots\rho_{n/2}^2$,
where $\rho_k$ are radii in orthogonal two-dimensional planes defined further
in Sec. 6.2. The amplitude $\rho$, the radii
$\rho_k$ and the variables $\tan\psi_{k-1}$ are multiplicative, and the
azimuthal angles $\phi_k$ are 
additive upon the multiplication of n-complex numbers.
Because the description of the position of the point $A$ requires, in addition
to the azimuthal angles $\phi_k$, only the planar angles $\psi_{k-1}$, 
the hypercomplex numbers studied in Sec. 6.2 will be called planar
n-complex number, to distinguish them from the polar n-complex numbers, which
in an even number of dimensions required two 
polar angles, and in an odd number of dimensions required one polar angle.

The exponential function of an n-complex number can be expanded in terms of
the planar n-dimensional cosexponential functions
$f_{nk}(y)=\sum_{p=0}^\infty (-1)^p y^{k+pn}/(k+pn)!, k=0,1,...,n-1$.
It is shown that
$f_{nk}(y)=(1/n)\sum_{l=1}^{n}$$\exp\left\{y\cos\left(\pi (2l-1)/n\right)\right\}$\\
$\cos\left\{y\sin\left(\pi (2l-1)/n\right)-\pi k(2l-1)/n\right\}, $
$k=0,1,...,n-1.$
Addition theorems and other relations
are obtained for the planar n-dimensional cosexponential functions.

The exponential form of a planar n-complex number, which can be defined for
all $x_0$, ..., $x_{n-1}$, is
$u=\rho \exp\left\{\sum_{p=1}^{n-1}h_p\left[
-(2/n)\sum_{k=2}^{n/2}
\cos\left(\pi (2k-1)p/n\right)\ln\tan\psi_{k-1}
\right]\right\}$
$\exp\left(\sum_{k=1}^{n/2}\tilde e_k\phi_k\right)$,
where $\tilde e_k=(2/n)\sum_{p=1}^{n-1}h_p\sin(\pi (2k-1)p/n)$. 
A trigonometric form also exists for an n-complex number,
$u=d\left(n/2\right)^{1/2}$
$\left(1+1/\tan^2\psi_1+1/\tan^2\psi_2+\cdots
+1/\tan^2\psi_{n/2-1}\right)^{-1/2}$\\
$\left(e_1+\sum_{k=2}^{n/2}e_k/\tan\psi_{k-1}\right)$
$\exp\left(\sum_{k=1}^{n/2}\tilde e_k\phi_k\right)$.

Expressions are given for the elementary functions of planar n-complex variable.
The functions $f(u)$ of planar n-complex variable which are defined by power series
have derivatives independent of the direction of approach to the point under
consideration. If the n-complex function $f(u)$ 
of the n-complex variable $u$ is written in terms of 
the real functions $P_k(x_0,...,x_{n-1}), k=0,...,n-1$, then
relations of equality  
exist between partial derivatives of the functions $P_k$. 
The integral $\int_A^B f(u) du$ of an n-complex
function between two points $A,B$ is independent of the path connecting $A,B$,
in regions where $f$ is regular.
If $f(u)$ is an analytic n-complex function, then 
$\oint_\Gamma f(u)du/(u-u_0)$
$=2\pi f(u_0)\sum_{k=1}^{n/2}\tilde e_k$ 
$\;{\rm int}(u_{0\xi_k\eta_k},\Gamma_{\xi_k\eta_k})$,
where the functional ${\rm int}$ takes the values 0 or 1 depending on the
relation between $u_{0\xi_k\eta_k}$ and $\Gamma_{\xi_k\eta_k}$, which are
respectively the projections of the point $u_0$ and of 
the loop $\Gamma$ on the plane defined by the orthogonal axes $\xi_k$ and
$\eta_k$, as expained further in Sec. 6.2.

A planar n-complex polynomial can always
be written as a product of linear factors, although the factorization may not
be unique. 

For $n=2$, the n-complex numbers discussed in this paper become the usual
2-dimensional complex numbers $x+iy$.
Particular cases for $n=4$ and $n=6$ of the planar
n-complex numbers described in this paper have been studied in previous chapters.

\section{Polar complex numbers in $n$ dimensions}

\subsection{Operations with polar n-complex numbers}

A complex number in $n$ dimensions is determined by its $n$ components
$(x_0,x_1,...,x_{n-1})$. The polar n-complex numbers and
their operations discussed in this section can be represented 
by  writing the n-complex number $(x_0,x_1,...,x_{n-1})$ as  
$u=x_0+h_1x_1+h_2x_2+\cdots+h_{n-1}x_{n-1}$, where 
$h_1, h_2, \cdots, h_{n-1}$ are bases for which the multiplication rules are 
\begin{equation}
h_j h_k =h_l ,\:l=j+k-n[(j+k)/n],
\label{npolar-1}
\end{equation}\index{complex units!polar n-complex}
for $ j,k,l=0,1,..., n-1$.
In Eq. (\ref{npolar-1}), $[(j+k)/n]$ denotes the integer part of $(j+k)/n$, 
the integer part being defined as $[a]\leq a<[a]+1$, so that
$0\leq j+k-n[(j+k)/n]\leq n-1$. 
In this chapter, brackets larger than the regular brackets
$[\;]$ do not have the meaning of integer part.
The significance of the composition laws in Eq.
(\ref{npolar-1}) can be understood by representing the bases $h_j, h_k$ by points on a
circle at the angles $\alpha_j=2\pi j/n,\alpha_k=2\pi k/n$, as shown in Fig. \ref{fig21},
and the product $h_j h_k$ by the point of the circle at the angle 
$2\pi (j+k)/n$. If $2\pi\leq 2\pi (j+k)/n<4\pi$, the point represents the basis
$h_l$ of angle $\alpha_l=2\pi(j+k-n)/n$.

\begin{figure}
\begin{center}
\epsfig{file=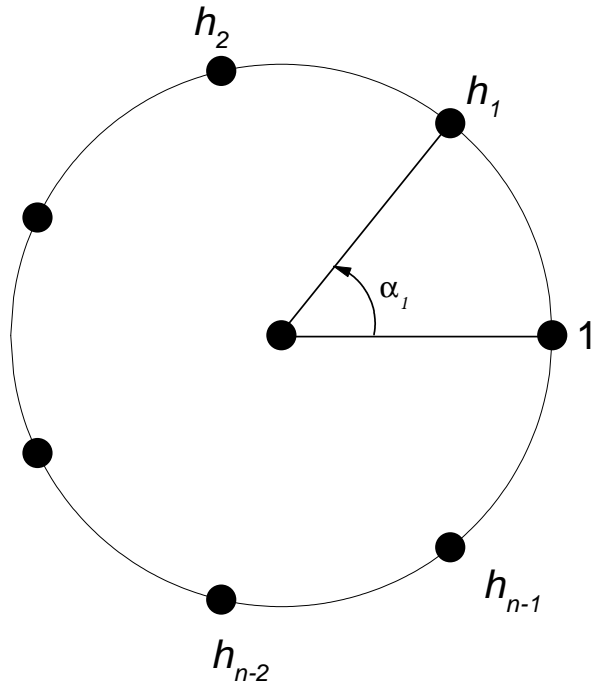,width=12cm}
\caption{Representation of the hypercomplex bases $1, h_1,...,h_{n-1}$
by points on a circle at the angles $\alpha_k=2\pi k/n$.
The product $h_j h_k$ will be represented by the point of the circle at the
angle $2\pi (j+k)/n$, $i,k=0,1,...,n-1$. If $2\pi\leq 
2\pi (j+k)/n\leq 4\pi$, the point represents the basis
$h_l$ of angle $\alpha_l=2\pi(j+k)/n-2\pi$.}
\label{fig21}
\end{center}
\end{figure}

Two n-complex numbers 
$u=x_0+h_1x_1+h_2x_2+\cdots+h_{n-1}x_{n-1}$,
$u^\prime=x^\prime_0+h_1x^\prime_1+h_2x^\prime_2+\cdots+h_{n-1}x^\prime_{n-1}$ 
are equal if and only if $x_i=x^\prime_i, i=0,1,...,n-1$.
The sum of the n-complex numbers $u$
and
$u^\prime$ 
is
\begin{equation}
u+u^\prime=x_0+x^\prime_0+h_1(x_1+x^\prime_1)+\cdots
+h_{n-1}(x_{n-1} +x^\prime_{n-1}) .
\label{npolar-2}
\end{equation}\index{sum!polar n-complex}
The product of the numbers $u, u^\prime$ is 
\begin{equation}
\begin{array}{l}
uu^\prime=x_0 x_0^\prime +x_1x_{n-1}^\prime+x_2 x_{n-2}^\prime+x_3x_{n-3}^\prime
+\cdots+x_{n-1}x_1^\prime\\
+h_1(x_0 x_1^\prime+x_1x_0^\prime+x_2x_{n-1}^\prime+x_3x_{n-2}^\prime
+\cdots+x_{n-1} x_2^\prime) \\
+h_2(x_0 x_2^\prime+x_1x_1^\prime+x_2x_0^\prime+x_3x_{n-1}^\prime
+\cdots+x_{n-1} x_3^\prime) \\
\vdots\\
+h_{n-1}(x_0 x_{n-1}^\prime+x_1x_{n-2}^\prime+x_2x_{n-3}^\prime+x_3x_{n-4}^\prime
+\cdots+x_{n-1} x_0^\prime).
\end{array}
\label{npolar-3}
\end{equation}\index{product!polar n-complex}
The product $uu^\prime$ can be written as
\begin{equation}
uu^\prime=\sum_{k=0}^{n-1}h_k\sum_{l=0}^{n-1}x_l x^\prime_{k-l+n[(n-k-1+l)/n]}.
\label{npolar-3a}
\end{equation}
If $u,u^\prime,u^{\prime\prime}$ are n-complex numbers, the multiplication 
is associative
\begin{equation}
(uu^\prime)u^{\prime\prime}=u(u^\prime u^{\prime\prime})
\label{npolar-3b}
\end{equation}
and commutative
\begin{equation}
u u^\prime=u^\prime u ,
\label{npolar-3c}
\end{equation}
because the product of the bases, defined in Eq. (\ref{npolar-1}), is associative and
commutative. The fact that the multiplication is commutative can be seen also
directly from Eq. (\ref{npolar-3}).
The n-complex zero is $0+h_1\cdot 0+\cdots+h_{n-1}\cdot 0,$ 
denoted simply 0, 
and the n-complex unity is $1+h_1\cdot 0+\cdots+h_{n-1}\cdot 0,$ 
denoted simply 1.

The inverse of the n-complex number $u=x_0+h_1x_1+h_2x_2+\cdots+h_{n-1}x_{n-1}$
is the n-complex number
$u^\prime
=x^\prime_0+h_1x^\prime_1+h_2x^\prime_2+\cdots+h_{n-1}x^\prime_{n-1}$ 
having the property that
\begin{equation}
uu^\prime=1 .
\label{npolar-4}
\end{equation}\index{inverse!polar n-complex}
Written on components, the condition, Eq. (\ref{npolar-4}), is
\begin{equation}
\begin{array}{l}
x_0 x_0^\prime +x_1x_{n-1}^\prime+x_2 x_{n-2}^\prime+x_3x_{n-3}^\prime
+\cdots+x_{n-1}x_1^\prime=1,\\
x_0 x_1^\prime+x_1x_0^\prime+x_2x_{n-1}^\prime+x_3x_{n-2}^\prime
+\cdots+x_{n-1} x_2^\prime=0,\\
x_0 x_2^\prime+x_1x_1^\prime+x_2x_0^\prime+x_3x_{n-1}^\prime
+\cdots+x_{n-1} x_3^\prime=0, \\
\vdots\\
x_0 x_{n-1}^\prime+x_1x_{n-2}^\prime+x_2x_{n-3}^\prime+x_3x_{n-4}^\prime
+\cdots+x_{n-1} x_0^\prime=0.
\end{array}
\label{npolar-5}
\end{equation}
The system (\ref{npolar-5}) has a solution provided that the determinant of the
system, 
\begin{equation}
\nu={\rm det}(A), 
\label{npolar-5b}
\end{equation}\index{inverse, determinant!polar n-complex}
is not equal to zero, $\nu\not=0$, where
\begin{equation}
A=\left(
\begin{array}{ccccc}
x_0     &   x_{n-1} &   x_{n-2}   & \cdots  &x_1\\
x_1     &   x_0     &   x_{n-1}   & \cdots  &x_2\\
x_2     &   x_1     &   x_0       & \cdots  &x_3\\
\vdots  &  \vdots   &  \vdots & \cdots  &\vdots \\
x_{n-1} &   x_{n-2} &   x_{n-3}   & \cdots  &x_0\\
\end{array}
\right).
\label{npolar-6}
\end{equation}
If $\nu>0$, the quantity 
\begin{equation}
\rho=\nu^{1/n}
\label{npolar-6b}
\end{equation}\index{amplitude!polar n-complex}
will be called amplitude of the
n-complex number $u=x_0+h_1x_1+h_2x_2+\cdots+h_{n-1}x_{n-1}$.
The quantity $\nu$ can be written as a product of linear factors
\begin{equation}
\nu=\prod_{k=0}^{n-1}
\left(x_0+\epsilon_k x_1+\epsilon_k^2 x_2+\cdots+\epsilon^{n-1}_k x_{n-1}\right),
\label{npolar-7}
\end{equation}
where $\epsilon_k=e^{2\pi ik/n}$, $i$ being the imaginary
unit. The factors appearing in Eq. (\ref{npolar-7}) are of the form 
\begin{equation}
x_0+\epsilon_k x_1+\epsilon_k^2 x_2+\cdots+\epsilon^{n-1}_k x_{n-1}=v_k+i\tilde v_k,
\label{npolar-8}
\end{equation}
where
\begin{equation}
v_k=\sum_{p=0}^{n-1}x_p\cos\frac{2\pi kp}{n},
\label{npolar-9a}
\end{equation}
\begin{equation}
\tilde v_k=\sum_{p=0}^{n-1}x_p\sin\frac{2\pi kp}{n},
\label{npolar-9b}
\end{equation}\index{canonical variables!polar n-complex}
for $k=1,2,...,n-1$ and, if $n$ is even, $k\not=n/2$.
For $k=0$ the factor in Eq. (\ref{npolar-8}) is 
\begin{equation}
v_+=x_0+x_1+\cdots+x_{n-1},
\label{npolar-9bb}
\end{equation}
and if $n$ is even, for $k=n/2$ the factor in Eq. (\ref{npolar-8}) is
\begin{equation}
v_-=x_0-x_1+\cdots+x_{n-2}-x_{n-1}.
\label{npolar-9bbb}
\end{equation}
It can be seen that $v_k=v_{n-k}, \tilde v_k=-\tilde v_{n-k}$,
$k=1,...,[(n-1)/2]$. 
The variables $v_+, v_-, v_k, \tilde v_k, k=1,...,[(n-1)/2] $ 
will be called canonical polar n-complex variables. 
Therefore,  the factors appear in Eq. (\ref{npolar-7}) in complex-conjugate pairs of
the form $v_k+i\tilde v_k$ and $v_{n-k}+i\tilde v_{n-k}=v_k-i\tilde v_k$, where
$k=1,...,[(n-1)/2]$, so that the product $\nu$ is a real quantity. 
If $n$ is an even number, the quantity $\nu$ is
\begin{equation}
\nu=v_+v_-\prod_{k=1}^{n/2-1}(v_k^2+\tilde v_k^2), 
\label{npolar-9c}
\end{equation}\index{inverse, determinant!polar n-complex}
and if $n$ is an odd number, $\nu$ is
\begin{equation}
\nu=v_+\prod_{k=0}^{(n-1)/2}(v_k^2+\tilde v_k^2).
\label{npolar-9d}
\end{equation}
Thus, in an even number of dimensions $n$, an n-complex number has an inverse
unless it lies on one of the nodal hypersurfaces $x_0+x_1+\cdots+x_{n-1}=0$, or
$x_0-x_1+\cdots+x_{n-2}-x_{n-1}=0$, or $v_1=0, \tilde v_1=0$, ...,
or $v_{n/2-1}=0, \tilde v_{n/2-1}=0$. 
In an odd number of dimensions $n$, an n-complex number has an inverse
unless it lies on one of the nodal hypersurfaces $x_0+x_1+\cdots+x_{n-1}=0$,
or $v_1=0, \tilde v_1=0$, ..., or $v_{(n-1)/2}=0, \tilde v_{(n-1)/2}=0$. 
\index{nodal hypersurfaces!polar n-complex}

\subsection{Geometric representation of polar n-complex numbers}

The n-complex number $x_0+h_1x_1+h_2x_2+\cdots+h_{n-1}x_{n-1}$
can be represented by 
the point $A$ of coordinates $(x_0,x_1,...,x_{n-1})$. 
If $O$ is the origin of the n-dimensional space,  the
distance from the origin $O$ to the point $A$ of coordinates
$(x_0,x_1,...,x_{n-1})$ has the expression
\begin{equation}
d^2=x_0^2+x_1^2+\cdots+x_{n-1}^2 .
\label{npolar-10}
\end{equation}\index{distance!polar n-complex}
The quantity $d$ will be called modulus of the n-complex number 
$u=x_0+h_1x_1+h_2x_2+\cdots+h_{n-1}x_{n-1}$. The modulus of an n-complex number
$u$ will be designated by $d=|u|$.\index{modulus, definition!polar n-complex}

The exponential and trigonometric forms of the n-complex number $u$ can be
obtained conveniently in a rotated system of axes defined by a transformation
which, for even $n$, has the form
\begin{equation}
\left(
\begin{array}{c}
\xi_+\\
\xi_-\\
\vdots\\
\xi_k\\
\eta_k\\
\vdots
\end{array}\right)
=\left(
\begin{array}{ccccc}
\frac{1}{\sqrt{n}}&\frac{1}{\sqrt{n}}&\cdots&\frac{1}{\sqrt{n}}&\frac{1}{\sqrt{n}}\\
\frac{1}{\sqrt{n}}&-\frac{1}{\sqrt{n}}&\cdots&\frac{1}{\sqrt{n}}&-\frac{1}{\sqrt{n}}\\
\vdots&\vdots& &\vdots&\vdots\\
\sqrt{\frac{2}{n}}&\sqrt{\frac{2}{n}}\cos\frac{2\pi k}{n}&\cdots&\sqrt{\frac{2}{n}}\cos\frac{2\pi (n-2)k}{n}&\sqrt{\frac{2}{n}}\cos\frac{2\pi (n-1)k}{n}\\
0&\sqrt{\frac{2}{n}}\sin\frac{2\pi k}{n}&\cdots&\sqrt{\frac{2}{n}}\sin\frac{2\pi (n-2)k}{n}&\sqrt{\frac{2}{n}}\sin\frac{2\pi (n-1)k}{n}\\
\vdots&\vdots&&\vdots&\vdots
\end{array}
\right)
\left(
\begin{array}{c}
x_0\\x_1\\
\vdots\\ 
\vdots\\
\vdots\\
x_{n-1}
\end{array}
\right),
\label{npolar-11}
\end{equation}
where $k=1, 2, ... , n/2-1$.
For odd $n$ the rotation of the axes is described by the relations
\begin{equation}
\left(
\begin{array}{c}
\xi_+\\
\xi_1\\
\eta_1\\
\vdots\\
\xi_k\\
\eta_k\\
\vdots
\end{array}\right)
=\left(
\begin{array}{cccc}
\frac{1}{\sqrt{n}}&\frac{1}{\sqrt{n}}&\cdots&\frac{1}{\sqrt{n}}\\
\sqrt{\frac{2}{n}}&\sqrt{\frac{2}{n}}\cos\frac{2\pi }{n}&\cdots&\sqrt{\frac{2}{n}}\cos\frac{2\pi (n-1)}{n}\\
0&\sqrt{\frac{2}{n}}\sin\frac{2\pi }{n}&\cdots&\sqrt{\frac{2}{n}}\sin\frac{2\pi (n-1)}{n}\\
\vdots&\vdots& &\vdots\\
\sqrt{\frac{2}{n}}&\sqrt{\frac{2}{n}}\cos\frac{2\pi k}{n}&\cdots&\sqrt{\frac{2}{n}}\cos\frac{2\pi (n-1)k}{n}\\
0&\sqrt{\frac{2}{n}}\sin\frac{2\pi k}{n}&\cdots&\sqrt{\frac{2}{n}}\sin\frac{2\pi (n-1)k}{n}\\
\vdots&\vdots&&\vdots
\end{array}
\right)
\left(
\begin{array}{c}
x_0\\
x_1\\
x_2\\
\vdots\\ 
\vdots\\
\vdots\\
x_{n-1}
\end{array}
\right),
\label{npolar-12}
\end{equation}
where $k=0,1,...,(n-1)/2$.
The lines of the matrices in Eqs. (\ref{npolar-11}) or (\ref{npolar-12}) give the components
of the $n$ basis vectors of the new system of axes. These vectors have unit
length and are orthogonal to each other.
By comparing Eqs. (\ref{npolar-9a})-(\ref{npolar-9bbb}) and (\ref{npolar-11})-(\ref{npolar-12}) it can be
seen that
\begin{equation}
v_+= \sqrt{n}\xi_+ ,  v_-=  \sqrt{n}\xi_-,  
v_k= \sqrt{\frac{n}{2}}\xi_k , \tilde v_k= \sqrt{\frac{n}{2}}\eta_k ,
\label{npolar-12b}
\end{equation}
i.e. the two sets of variables differ only by scale factors.

The sum of the squares of the variables $v_k,\tilde v_k$ is, for even $n$,
\begin{equation}
\sum_{k=1}^{n/2-1}(v_k^2+\tilde v_k^2)=\frac{n-2}{2}(x_0^2+\cdots+x_{n-1}^2)
-2(x_0x_2+\cdots+x_{n-4}x_{n-2}+x_1x_3+\cdots+x_{n-3}x_{n-1}) ,
\label{npolar-13}
\end{equation}
and for odd $n$ the sum is
\begin{equation}
\sum_{k=1}^{(n-1)/2}(v_k^2+\tilde v_k^2)=\frac{n-1}{2}(x_0^2+\cdots+x_{n-1}^2)
-(x_0x_1+\cdots+x_{n-2}x_{n-1}) .
\label{npolar-14}
\end{equation}
The relation (\ref{npolar-13}) has been obtained with the aid of the identity, valid
for even $n$,
\begin{equation}
\sum_{k=1}^{n/2-1}\cos\frac{2\pi pk}{n}=\left\{
\begin{array}{l}
-1, \:\:{\rm for \:\;even}\:\:p,\\
0, \:\:{\rm for \:\;odd}\:\:p.
\end{array}
\right.
\label{npolar-15}
\end{equation}
The relation (\ref{npolar-14}) has been obtained with the aid of the identity, valid
for odd values of $n$,
\begin{equation}
\sum_{k=1}^{(n-1)/2}\cos\frac{2\pi pk}{n}=-\frac{1}{2}.
\label{npolar-16}
\end{equation}
From Eq. (\ref{npolar-13}) it results that, for even $n$,
\begin{equation}
d^2=\frac{1}{n}v_+^2+\frac{1}{n}v_-^2
+\frac{2}{n}\sum_{k=1}^{n/2-1}\rho_k^2,
\label{npolar-17}
\end{equation}\index{modulus, canonical variables!polar n-complex}
and from Eq. (\ref{npolar-14}) it results that, for odd $n$,
\begin{equation}
d^2=\frac{1}{n}v_+^2
+\frac{2}{n}\sum_{k=1}^{(n-1)/2}\rho_k^2.
\label{npolar-18}
\end{equation}
The relations (\ref{npolar-17}) and (\ref{npolar-18}) show that the square of the distance
$d$, Eq. (\ref{npolar-10}), is the sum of the squares of the projections
$v_+/\sqrt{n}, \rho_k\sqrt{2/n}$ and, for even $n$, of the square
of $v_-/\sqrt{n}$. This is consistent with the fact that the transformation in
Eqs. (\ref{npolar-11}) or (\ref{npolar-12}) is unitary.

The position of the point $A$ of coordinates $(x_0,x_1,...,x_{n-1})$ can be
also described with the aid of the distance $d$, Eq. (\ref{npolar-10}), and of $n-1$
angles defined further. Thus,
in the plane of the axes $v_k,\tilde v_k$, the radius $\rho_k$ and the azimuthal angle
$\phi_k$ can be introduced by the relations 
\begin{equation}
\rho_k^2=v_k^2+\tilde v_k^2, \:\cos\phi_k=v_k/\rho_k,\:\sin\phi_k=\tilde v_k/\rho_k, 0\leq
\phi_k<2\pi ,
\label{npolar-19a}
\end{equation}\index{azimuthal angles!polar n-complex}
so that there are $[(n-1)/2]$ azimuthal angles.
If the projection of the point $A$ on the plane of the axes $v_k,\tilde v_k$ is $A_k$,
and the projection of the point $A$ on the 4-dimensional space defined by the
axes $v_1, \tilde v_1, v_k,\tilde v_k$ is $A_{1k}$, the angle $\psi_{k-1}$ between the line
$OA_{1k}$ and the 2-dimensional plane defined by the axes $v_k,\tilde v_k$ is 
\begin{equation}
\tan\psi_{k-1}=\rho_1/\rho_k, 
\label{npolar-19b}
\end{equation}\index{planar angles!polar n-complex}
where $0\leq\psi_k\leq\pi/2, k=2,...,[(n-1)/2]$,
so that there are $[(n-3)/2]$ planar angles.
Moreover, there is a polar angle $\theta_+$, which can be
defined as the angle between the line $OA_{1+}$ and the axis $v_+$,
where $A_{1+}$ is the projection of the point $A$ on the 3-dimensional space
generated by the axes $v_1, \tilde v_1, v_+$,
\begin{equation}
\tan\theta_+=\frac{\sqrt{2}\rho_1}{v_+}, 
\label{npolar-19c}
\end{equation}\index{polar angles!polar n-complex}
where $0\leq\theta_+\leq\pi$ , 
and in an even number of dimensions $n$ there is also a polar angle $\theta_-$,
which can be defined as the angle between the line $OA_{1-}$ and the axis
$v_-$, 
where $A_{1-}$ is the projection of the point $A$ on the 3-dimensional space
generated by the axes $v_1, \tilde v_1, v_-$,
\begin{equation}
\tan\theta_-=\frac{\sqrt{2}\rho_1}{v_-}, 
\label{npolar-19d}
\end{equation}
where $0\leq\theta_-\leq\pi$ .
In Eqs. (\ref{npolar-19c}) and (\ref{npolar-19d}), the factor $\sqrt{2}$ appears from the
ratio of the normalization factors in Eq. (\ref{npolar-12b}).
Thus, the position of the point $A$ is described, in an even number 
of dimensions, by the distance $d$, by $n/2-1$ azimuthal angles, by $n/2-2$
planar angles, and by 2 polar angles. In an odd number of dimensions, the
position of the point $A$ is described by $(n-1)/2$ azimuthal angles, by
$(n-3)/2$ planar angles, and by 1 polar angle. These angles are shown in Fig.
22.

\begin{figure}
\begin{center}
\epsfig{file=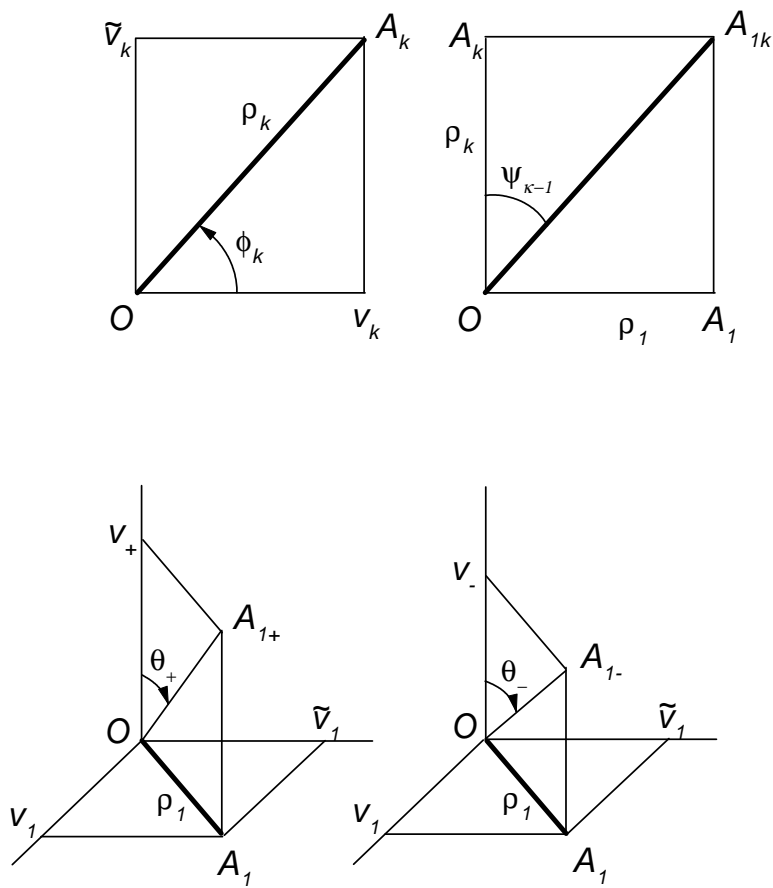,width=12cm}
\caption{Radial distance $\rho_k$  and
azimuthal angle $\phi_k$ in the plane of the axes $v_k,\tilde v_k$, and planar
angle $\psi_{k-1}$ between the line $OA_{1k}$ and the 2-dimensional plane
defined by the axes $v_k,\tilde v_k$. $A_k$ is the projection of the point $A$
on the plane of the axes $v_k,\tilde v_k$, and $A_{1k}$ is the projection of
the point $A$ on the 4-dimensional space defined 
by the axes $v_1, \tilde v_1, v_k,\tilde v_k$. 
The polar angle $\theta_+$ is the angle between the line $OA_{1+}$ and the axis
$v_+$, where $A_{1+}$ is the projection of the point $A$ on the 3-dimensional
space generated by the axes $v_1, \tilde v_1, v_+$.
In an even number of dimensions $n$ there is also a polar angle $\theta_-$,
which is the angle between the line $OA_{1-}$ and the axis $v_-$, 
where $A_{1-}$ is the projection of the point $A$ on the 3-dimensional space
generated by the axes $v_1, \tilde v_1, v_-$.}
\label{fig22}
\end{center}
\end{figure}

The variables $\rho_k$ can be expressed in terms of $d$ and the planar angles
$\psi_k$ as
\begin{equation}
\rho_k=\frac{\rho_1}{\tan\psi_{k-1}}, 
\label{npolar-20a}
\end{equation}
for $k=2,...,[(n-1)/2]$, where, for even $n$,
\begin{eqnarray}
\lefteqn{\rho_1^2=\frac{nd^2}{2}
\left(\frac{1}{\tan^2\theta_+}+\frac{1}{\tan^2\theta_-}+1
+\frac{1}{\tan^2\psi_1}+\frac{1}{\tan^2\psi_2}+\cdots
+\frac{1}{\tan^2\psi_{n/2-2}}\right)^{-1},\nonumber}\\
&&
\label{npolar-20b}
\end{eqnarray}
and for odd $n$
\begin{eqnarray}
\rho_1^2=\frac{nd^2}{2}
\left(\frac{1}{\tan^2\theta_+}+1
+\frac{1}{\tan^2\psi_1}+\frac{1}{\tan^2\psi_2}+\cdots
+\frac{1}{\tan^2\psi_{(n-3)/2}}\right)^{-1}.
\label{npolar-20c}
\end{eqnarray}

If
$u^\prime=x_0^\prime+h_1x_1^\prime+h_2x_2^\prime+\cdots+h_{n-1}x_{n-1}^\prime, 
u^{\prime\prime}=x^{\prime\prime}_0+h_1x^{\prime\prime}_1
+h_2x^{\prime\prime}_2+\cdots+h_{n-1}x^{\prime\prime}_{n-1}$ 
are n-complex numbers of parameters $v_+^\prime,  v_-^\prime,
\rho_k^\prime, \theta_+^\prime,  \theta_-^\prime,
\psi_k^\prime,\phi_k^\prime$ and respectively  
$v_+^{\prime\prime}$,  $v_-^{\prime\prime}$,
$\rho_k^{\prime\prime}$, $\theta_+^{\prime\prime}$, $\theta_-^\prime$,
$\psi_k^{\prime\prime}$, $\phi_k^{\prime\prime}$, then the parameters
$v_+$, $v_-$, $\rho_k$, $\theta_+$,  $\theta_-$, $\psi_k$, $\phi_k$ of the product n-complex
number $u=u^\prime u^{\prime\prime}$ are given by 
\begin{equation}
v_+=v_+^\prime v_+^{\prime\prime},  
\label{npolar-21a}
\end{equation}\index{transformation of variables!polar n-complex}
\begin{equation}
\rho_k=\rho_k^\prime\rho_k^{\prime\prime}, 
\label{npolar-21b}
\end{equation}
for $k=1,..., [(n-1)/2]$,
\begin{equation}
\tan\theta_+=\frac{1}{\sqrt{2}}\tan\theta_+^\prime \tan\theta_+^{\prime\prime},
\label{npolar-21c}
\end{equation}
\begin{equation}
\tan\psi_k=\tan\psi_k^\prime \tan\psi_k^{\prime\prime},
\label{npolar-21d}
\end{equation}
for $k=1,...,[(n-3)/2]$,  
\begin{equation}
\phi_k=\phi_k^\prime+\phi_k^{\prime\prime}, 
\label{npolar-21e}
\end{equation}
for $k=1,...,[(n-1)/2]$, 
and, if $n$ is even,
\begin{equation}
v_-=v_-^\prime v_-^{\prime\prime},  
\label{npolar-21f}
\end{equation}
\begin{equation}
\tan\theta_-=\frac{1}{\sqrt{2}}\tan\theta_-^\prime \tan\theta_-^{\prime\prime}.
\label{npolar-21g}
\end{equation}
The Eqs. (\ref{npolar-21a}) and (\ref{npolar-21f}) can be checked directly, and  
Eqs. (\ref{npolar-21b})-(\ref{npolar-21e}) and (\ref{npolar-21g}) are a consequence of the relations
\begin{equation}
v_k=v_k^\prime v_k^{\prime\prime}-\tilde v_k^\prime \tilde v_k^{\prime\prime},\;
\tilde v_k=v_k^\prime \tilde v_k^{\prime\prime}+\tilde v_k^\prime v_k^{\prime\prime},
\label{npolar-22}
\end{equation}
and of the corresponding relations of definition. Then the product $\nu$ in
Eqs. (\ref{npolar-9c}) and (\ref{npolar-9d}) has the property that
\begin{equation}
\nu=\nu^\prime\nu^{\prime\prime} 
\label{npolar-23}
\end{equation}
and, if $\nu^\prime>0, \nu^{\prime\prime}>0$, the amplitude $\rho$ defined in
Eq. (\ref{npolar-6b}) has the property that
\begin{equation}
\rho=\rho^\prime\rho^{\prime\prime} .
\label{npolar-24}
\end{equation}

The fact that the amplitude of the product is equal to the product of the 
amplitudes, as written in Eq. (\ref{npolar-24}), can 
be demonstrated also by using a representation of the 
n-complex numbers by matrices, in which the n-complex number 
$u=x_0+h_1x_1+h_2x_2+\cdots+h_{n-1}x_{n-1}$ is represented by the matrix
\begin{equation}
U=\left(
\begin{array}{ccccc}
x_0     &   x_1     &   x_2   & \cdots  &x_{n-1}\\
x_{n-1} &   x_0     &   x_1   & \cdots  &x_{n-2}\\
x_{n-2} &   x_{n-1} &   x_0   & \cdots  &x_{n-3}\\
\vdots  &  \vdots   &  \vdots & \cdots  &\vdots \\
x_1     &   x_2     &   x_3   & \cdots  &x_0\\
\end{array}
\right).
\label{npolar-24b}
\end{equation}\index{matrix representation!polar n-complex}
The product $u=u^\prime u^{\prime\prime}$ is
represented by the matrix multiplication $U=U^\prime U^{\prime\prime}$.
The relation (\ref{npolar-23}) is then a consequence of the fact the determinant 
of the product of matrices is equal to the product of the determinants 
of the factor matrices. The use of the representation of the n-complex
numbers with matrices provides an alternative demonstration of the fact
that the product of n-complex numbers is associative, as
stated in Eq. (\ref{npolar-3b}).

According to Eqs. (\ref{npolar-21a}, (\ref{npolar-21b}), (\ref{npolar-21f}), (\ref{npolar-17}) and
(\ref{npolar-18}), the modulus of the product $uu^\prime$ is, for even $n$, 
\begin{equation}
|uu^\prime|^2=
\frac{1}{n}(v_+v_+^\prime)^2+\frac{1}{n}(v_-v_-^\prime)^2
+\frac{2}{n}\sum_{k=1}^{n/2-1}(\rho_k\rho_k^\prime)^2 ,
\label{npolar-25a}
\end{equation}
and for odd $n$
\begin{equation}
|uu^\prime|^2=
\frac{1}{n}(v_+v_+^\prime)^2
+\frac{2}{n}\sum_{k=1}^{(n-1)/2}(\rho_k\rho_k^\prime)^2 .
\label{npolar-25b}
\end{equation}
Thus, if the product of two n-complex numbers is zero, $uu^\prime=0$, then
$v_+v_+^\prime=0, \rho_k\rho_k^\prime=0, k=1,...,[(n-1)/2]$ and, if $n$ is
even, $v_-v_-^\prime=0$. This means that either $u=0$, or $u^\prime=0$, or $u,
u^\prime$ belong to orthogonal hypersurfaces in such a way that the
afore-mentioned products of components should be equal to zero.
\index{divisors of zero!polar n-complex}

\subsection{The polar n-dimensional cosexponential functions}

The exponential function of a hypercomplex variable $u$ and the addition
theorem for the exponential function have been written in Eqs. 
~(1.35)-(1.36).
If $u=x_0+h_1x_1+h_2x_2+\cdots+h_{n-1}x_{n-1}$,
then $\exp u$ can be calculated as 
$\exp u=\exp x_0 \cdot \exp (h_1x_1) \cdots \exp (h_{n-1}x_{n-1})$.

It can be seen with the aid of the representation in Fig. \ref{fig21} that 
\begin{equation}
h_k^{n+p}=h_k^p, \:p\:\:{\rm integer},
\label{npolar-28}
\end{equation}\index{complex units, powers of!polar n-complex}
for $ k=1,...,n-1$.
Then $e^{h_k y}$ can be written as
\begin{equation}
e^{h_k y}=\sum_{p=0}^{n-1}h_{kp-n[kp/n]}g_{np}(y),
\label{npolar-28b}
\end{equation}\index{exponential, expressions!polar n-complex}
where the expression of the functions $g_{nk}$, which will be
called polar cosexponential functions in $n$ dimensions, is
\begin{equation}
g_{nk}(y)=\sum_{p=0}^\infty y^{k+pn}/(k+pn)!, 
\label{npolar-29}
\end{equation}\index{cosexponential function, polar n-complex!definition}
for $ k=0,1,...,n-1$.

If $n$ is even, the polar cosexponential functions of even index $k$ are
even functions, $g_{n,2p}(-y)=g_{n,2p}(y)$, $p=0,1,...,n/2-1$,
and the polar cosexponential functions of odd index 
are odd functions, $g_{n,2p+1}(-y)=-g_{n,2p+1}(y)$, $p=0,1,...,n/2-1$. 
\index{cosexponential functions, polar n-complex!parity}
For odd
values of $n$, the polar cosexponential functions do not have a definite
parity. It can be checked that
\begin{equation}
\sum_{k=0}^{n-1}g_{nk}(y)=e^y
\label{npolar-29a}
\end{equation}
and, for even $n$, 
\begin{equation}
\sum_{k=0}^{n-1}(-1)^k g_{nk}(y)=e^{-y}.
\label{npolar-29b}
\end{equation}

The expression of the polar n-dimensional cosexponential functions is
\begin{equation}
g_{nk}(y)=\frac{1}{n}\sum_{l=0}^{n-1}\exp\left[y\cos\left(\frac{2\pi l}{n}\right)
\right]
\cos\left[y\sin\left(\frac{2\pi l}{n}\right)-\frac{2\pi kl}{n}\right], 
\label{npolar-30}
\end{equation}\index{cosexponential functions, polar n-complex!expressions}
for $k=0,1,...,n-1$.
In order to check that the function in Eq. (\ref{npolar-30}) has the series expansion
written in Eq. (\ref{npolar-29}), the right-hand side of Eq. (\ref{npolar-30}) will be
written as
\begin{equation}
g_{nk}(y)=\frac{1}{n}\sum_{l=0}^{n-1}{\rm Re}\left\{
\exp\left[\left(\cos\frac{2\pi l}{n}+i\sin\frac{2\pi l}{n}\right)y
-i\frac{2\pi kl}{n}\right]\right\}, 
\label{npolar-31}
\end{equation}
for $k=0,1,...,n-1$,
where ${\rm Re} (a+ib)=a$, with $a$ and $b$ real numbers. The part of the
exponential depending on $y$ can be expanded in a series,
\begin{equation}
g_{nk}(y)=\frac{1}{n}\sum_{p=0}^\infty\sum_{l=0}^{n-1}{\rm Re}
\left\{\frac{1}{p!}
\exp\left[i\frac{2\pi l}{n}(p-k)\right]y^p\right\}, 
\label{npolar-32}
\end{equation}
for $k=0,1,...,n-1$.
The expression of $g_{nk}(y)$ becomes
\begin{equation}
g_{nk}(y)=\frac{1}{n}\sum_{p=0}^\infty\sum_{l=0}^{n-1}
\left\{\frac{1}{p!}
\cos\left[\frac{2\pi l}{n}(p-k)\right]y^p\right\}, 
\label{npolar-33}
\end{equation}
for $k=0,1,...,n-1$ and, since
\begin{equation}
\frac{1}{n}\sum_{l=0}^{n-1}\cos\frac{2\pi l}{n}(p-k)
=\left\{
\begin{array}{l}
1, \:\:{\rm if}\:\: p-k\:\: {\rm is \:\;a\:\;multiple\:\;of\:\;}n,\\
0, \:\:{\rm otherwise},
\end{array}
\right.
\label{npolar-34}
\end{equation}
this yields indeed the expansion in Eq. (\ref{npolar-29}).

It can be shown from Eq. (\ref{npolar-30}) that
\begin{equation}
\sum_{k=0}^{n-1}g_{nk}^2(y)=\frac{1}{n}\sum_{l=0}^{n-1}\exp\left[2y\cos\left(
\frac{2\pi l}{n}\right)\right].
\label{npolar-34a}
\end{equation}
It can be seen that the right-hand side of Eq. (\ref{npolar-34a}) does not contain
oscillatory terms. If $n$ is a multiple of 4, it can be shown by replacing $y$
by $iy$ in Eq. (\ref{npolar-34a}) that
\begin{equation}
\sum_{k=0}^{n-1}(-1)^kg_{nk}^2(y)=\frac{2}{n}\left\{1+\cos 2y
+\sum_{l=1}^{n/4-1}\cos\left[2y\cos\left(\frac{2\pi l}{n}\right)\right]\right\},
\label{npolar-34b}
\end{equation}
which does not contain exponential terms.

Addition theorems for the polar n-dimensional cosexponential functions can be
obtained from the relation $\exp h_1(y+z)=\exp h_1 y \cdot\exp h_1 z $, by
substituting the expression of the exponentials as given in Eq. (\ref{npolar-28b})
for $k=1$, $e^{h_1 y}=g_{n0}(y)+h_1g_{n1}(y)+\cdots+h_{n-1} g_{n,n-1}(y)$,
\begin{eqnarray}
\lefteqn{g_{nk}(y+z)=g_{n0}(y)g_{nk}(z)+g_{n1}(y)g_{n,k-1}(z)+\cdots+g_{nk}(y)g_{n0}(z)\nonumber}\\
&&+g_{n,k+1}(y)g_{n,n-1}(z)+g_{n,k+2}(y)g_{n,n-2}(z)+\cdots+g_{n,n-1}(y)g_{n,k+1}(z) ,
\label{npolar-35a}
\end{eqnarray}\index{cosexponential functions, polar n-complex!addition theorems}
where $k=0,1,...,n-1$.
For $y=z$ the relations (\ref{npolar-35a}) take the form
\begin{eqnarray}
\lefteqn{g_{nk}(2y)=g_{n0}(y)g_{nk}(y)+g_{n1}(y)g_{n,k-1}(y)+\cdots+g_{nk}(y)g_{n0}(y)\nonumber}\\
&&+g_{n,k+1}(y)g_{n,n-1}(y)+g_{n,k+2}(y)g_{n,n-2}(y)+\cdots+g_{n,n-1}(y)g_{n,k+1}(y) ,
\label{npolar-35b}
\end{eqnarray}
where $k=0,1,...,n-1$.
For $y=-z$ the relations (\ref{npolar-35a}) and (\ref{npolar-29}) yield
\begin{equation}
g_{n0}(y)g_{n0}(-y)
+g_{n1}(y)g_{n,n-1}(-y)+g_{n2}(y)g_{n,n-2}(-y)+\cdots+g_{n,n-1}(y)g_{n1}(-y)=1 ,
\label{npolar-36a}
\end{equation}
\begin{eqnarray}
\lefteqn{g_{n0}(y)g_{nk}(-y)+g_{n1}(y)g_{n,k-1}(-y)+\cdots+g_{nk}(y)g_{n0}(-y)\nonumber}\\
&&+g_{n,k+1}(y)g_{n,n-1}(-y)+g_{n,k+2}(y)g_{n,n-2}(-y)+\cdots+g_{n,n-1}(y)g_{n,k+1}(-y)=0 ,
\nonumber\\
&&
\label{npolar-36b}
\end{eqnarray}
for $k=1,...,n-1$.

From Eq. (\ref{npolar-28b}) it can be shown, for natural numbers $l$, that
\begin{equation}
\left(\sum_{p=0}^{n-1}h_{kp-n[kp/n]}g_{np}(y)\right)^l
=\sum_{p=0}^{n-1}h_{kp-n[kp/n]}g_{np}(ly), 
\label{npolar-37}
\end{equation}
where $k=0,1,...,n-1$.
For $k=1$ the relation (\ref{npolar-37}) is
\begin{equation}
\left\{g_{n0}(y)+h_1g_{n1}(y)+\cdots+h_{n-1}g_{n,n-1}(y)\right\}^l
=g_{n0}(ly)+h_1g_{n1}(ly)+\cdots+h_{n,n-1}g_{n,n-1}(ly).
\label{npolar-37b}
\end{equation}

If
\begin{equation}
a_k=\sum_{p=0}^{n-1}g_{np}(y)\cos\left(\frac{2\pi kp}{n}\right), 
\label{npolar-38a}
\end{equation}
for $k=0,1,...,n-1$, and
\begin{equation}
b_k=\sum_{p=0}^{n-1}g_{np}(y)\sin\left(\frac{2\pi kp}{n}\right), 
\label{npolar-38b}
\end{equation}
for $k=1,...,n-1$, 
where $g_{nk}(y)$ are the polar cosexponential functions in Eq. (\ref{npolar-30}), it can
be shown that 
\begin{equation}
a_k=\exp\left[y\cos\left(\frac{2\pi k}{n}\right)\right]
\cos\left[y\sin\left(\frac{2\pi k}{n}\right)\right], 
\label{npolar-39a}
\end{equation}
where $k=0,1,...,n-1$,
\begin{equation}
b_k=\exp\left[y\cos\left(\frac{2\pi k}{n}\right)\right]
\sin\left[y\sin\left(\frac{2\pi k}{n}\right)\right], 
\label{npolar-39b}
\end{equation}
where $k=1,...,n-1$.
If
\begin{equation}
G_k^2=a_k^2+b_k^2, 
\label{npolar-40}
\end{equation}
for $k=1,...,n-1$,
then from Eqs. (\ref{npolar-39a}) and (\ref{npolar-39b}) it results that
\begin{equation}
G_k^2=\exp\left[2y\cos\left(\frac{2\pi k}{n}\right)\right], 
\label{npolar-41}
\end{equation}
where $k=1,...,n-1$. If
\begin{equation}
G_+=g_{n0}+g_{n1}+\cdots+g_{n,n-1} ,
\label{npolar-42a}
\end{equation}
from Eq. (\ref{npolar-38a}) it results that $G_+=a_0$, so that $G_+=e^y$, and, in an
even number of dimensions $n$, if 
\begin{equation}
G_-=g_{n0}-g_{n1}+\cdots+g_{n,n-2}-g_{n,n-1} ,
\label{npolar-42b}
\end{equation}
from Eq. (\ref{npolar-38a}) it results that $G_-=a_{n/2}$, so that $G_{n/2}=e^{-y}$.
Then with the aid of Eq. (\ref{npolar-15}) applied for $p=1$ it can be shown that the
polar n-dimensional cosexponential functions have the property that, 
for even $n$,
\begin{equation}
G_+G_-\prod_{k=1}^{n/2-1}G_k^2=1,
\label{npolar-43a}
\end{equation}
and in an odd number of dimensions, with the aid of Eq. (\ref{npolar-16}) it can be
shown that  
\begin{equation}
G_+\prod_{k=1}^{(n-1)/2}G_k^2=1.
\label{npolar-43b}
\end{equation}

The polar n-dimensional cosexponential functions are solutions of the
$n^{\rm th}$-order differential equation
\begin{equation}
\frac{d^n\zeta}{du^n}=\zeta ,
\label{npolar-44}
\end{equation}\index{cosexponential functions, polar n-complex!differential equations}
whose solutions are of the form
$\zeta(u)=A_0g_{n0}(u)+A_1g_{n1}(u)+\cdots+A_{n-1}g_{n,n-1}(u).$ 
It can be checked that the derivatives of the polar cosexponential functions
are related by
\begin{equation}
\frac{dg_{n0}}{du}=g_{n,n-1}, \:
\frac{dg_{n1}}{du}=g_{n0}, \:...,
\frac{dg_{n,n-2}}{du}=g_{n,n-3} ,
\frac{dg_{n,n-1}}{du}=g_{n,n-2} .
\label{npolar-45}
\end{equation}

\subsection{Exponential and trigonometric forms of polar n-complex numbers}

In order to obtain the exponential and trigonometric forms of n-complex
numbers, a canonical base 
$e_+,e_-,e_1,\tilde e_1,...,e_{n/2-1},\tilde e_{n/2-1}$ 
for the polar n-complex numbers will be introduced for
even $n$ by the relations 
\begin{equation}
\left(
\begin{array}{c}
e_+\\
e_-\\
\vdots\\
e_k\\
\tilde e_k\\
\vdots
\end{array}\right)
=\left(
\begin{array}{ccccc}
\frac{1}{n}&\frac{1}{n}&\cdots&\frac{1}{n}&\frac{1}{n}\\
\frac{1}{n}&-\frac{1}{n}&\cdots&\frac{1}{n}&-\frac{1}{n}\\
\vdots&\vdots& &\vdots&\vdots\\
\frac{2}{n}&\frac{2}{n}\cos\frac{2\pi k}{n}&\cdots&\frac{2}{n}\cos\frac{2\pi (n-2)k}{n}&\frac{2}{n}\cos\frac{2\pi (n-1)k}{n}\\
0&\frac{2}{n}\sin\frac{2\pi k}{n}&\cdots&\frac{2}{n}\sin\frac{2\pi (n-2)k}{n}&\frac{2}{n}\sin\frac{2\pi (n-1)k}{n}\\
\vdots&\vdots&&\vdots&\vdots
\end{array}
\right)
\left(
\begin{array}{c}
1\\h_1\\
\vdots\\ 
\vdots\\
\vdots\\
h_{n-1}
\end{array}
\right),
\label{npolar-e11}
\end{equation}\index{canonical base!polar n-complex}
where $k=1, 2, ... , n/2-1$.
For odd $n$, the canonical base 
$e_+, e_1,\tilde e_1,...e_{(n-1)/2},\tilde e_{(n-1)/2}$ 
for the polar n-complex numbers will be introduced by the relations
\begin{equation}
\left(
\begin{array}{c}
e_+\\
e_1\\
\tilde e_1\\
\vdots\\
e_k\\
\tilde e_k\\
\vdots
\end{array}\right)
=\left(
\begin{array}{ccccc}
\frac{1}{n}&\frac{1}{n}&\cdots&\frac{1}{n}\\
\frac{2}{n}&\frac{2}{n}\cos\frac{2\pi }{n}&\cdots&\frac{2}{n}\cos\frac{2\pi (n-1)}{n}\\
0&\frac{2}{n}\sin\frac{2\pi }{n}&\cdots&\frac{2}{n}\sin\frac{2\pi (n-1)}{n}\\
\vdots&\vdots& &\vdots\\
\frac{2}{n}&\frac{2}{n}\cos\frac{2\pi k}{n}&\cdots&\frac{2}{n}\cos\frac{2\pi (n-1)k}{n}\\
0&\frac{2}{n}\sin\frac{2\pi k}{n}&\cdots&\frac{2}{n}\sin\frac{2\pi (n-1)k}{n}\\
\vdots&\vdots&&\vdots
\end{array}
\right)
\left(
\begin{array}{c}
1\\
h_1\\
h_2\\
\vdots\\ 
\vdots\\
\vdots\\
h_{n-1}
\end{array}
\right),
\label{npolar-e12}
\end{equation}
where $k=0,1,...,(n-1)/2$.

The multiplication relations for the new bases are, for even $n$,
\begin{eqnarray}
\lefteqn{e_+^2=e_+,\; e_-^2=e_-,\; e_+e_-=0,\; e_+e_k=0,\; e_+\tilde e_k=0,\;
e_-e_k=0,\; 
e_-\tilde e_k=0,\nonumber}\\ 
&&e_k^2=e_k,\; \tilde e_k^2=-e_k,\; e_k \tilde e_k=\tilde e_k ,\; e_ke_l=0,\;
e_k\tilde e_l=0,\; \tilde e_k\tilde e_l=0,\; k\not=l,\; 
\label{npolar-e12a}
\end{eqnarray}
where $k,l=1,...,n/2-1$.
For odd $n$ the multiplication relations are
\begin{eqnarray}
\lefteqn{e_+^2=e_+,\;  e_+e_k=0,\; e_+\tilde e_k=0,\; \nonumber}\\ 
&&e_k^2=e_k,\; \tilde e_k^2=-e_k,\; e_k \tilde e_k=\tilde e_k ,\; e_ke_l=0,\;
e_k\tilde e_l=0,\; \tilde e_k\tilde e_l=0,\; k\not=l,\; 
\label{npolar-e12b}
\end{eqnarray}
where $k,l=1,...,(n-1)/2$.
The moduli of the new bases are
\begin{equation}
|e_+|=\frac{1}{\sqrt{n}},\; |e_-|=\frac{1}{\sqrt{n}},\; 
|e_k|=\sqrt{\frac{2}{n}},\; |\tilde e_k|=\sqrt{\frac{2}{n}}.
\label{npolar-e12c}
\end{equation}
It can be shown that, for even $n$,
\begin{eqnarray}
x_0+h_1x_1+\cdots+h_{n-1}x_{n-1}=e_+v_+ + e_-v_- 
+\sum_{k=1}^{n/2-1}(e_k v_k+\tilde e_k \tilde v_k),
\label{npolar-e13a}
\end{eqnarray}\index{canonical form!polar n-complex}
and for odd $n$
\begin{eqnarray}
x_0+h_1x_1+\cdots+h_{n-1}x_{n-1}=e_+v_+ 
+\sum_{k=1}^{(n-1)/2}(e_k v_k+\tilde e_k \tilde v_k).
\label{npolar-e13b}
\end{eqnarray}
The relations (\ref{npolar-e13a}),(\ref{npolar-e13b}) give the canonical form
of a polar n-complex number.

Using the properties of the bases in Eqs. (\ref{npolar-e12a}) and (\ref{npolar-e12b}) it can
be shown that
\begin{equation}
\exp(\tilde e_k\phi_k)=1-e_k+e_k\cos\phi_k+\tilde e_k\sin\phi_k ,
\label{npolar-46a}
\end{equation}
\begin{equation}
\exp(e_k\ln\rho_k)=1-e_k+e_k\rho_k ,
\label{npolar-46b}
\end{equation}
\begin{equation}
\exp(e_+\ln v_+)=1-e_++e_+v_+
\label{npolar-46c}
\end{equation}
and, for even $n$,
\begin{equation}
\exp(e_-\ln v_-)=1-e_- +e_-v_- .
\label{npolar-46d}
\end{equation}
In Eq. (\ref{npolar-46c}), $\ln v_+$ exists as a real function 
provided that $v_+=x_0+x_1+\cdots+x_{n-1}>0$, which means that
$0<\theta_+<\pi/2$, 
and for even $n$, $\ln v_-$
exists in Eq. (\ref{npolar-46d}) as a real function provided that
$v_-=x_0-x_1+\cdots+x_{n-2}-x_{n-1}>0$, which means that $0<\theta_-<\pi/2$.
By multiplying the relations (\ref{npolar-46a})-(\ref{npolar-46d}) it results, for even $n$,
that 
\begin{equation}
\exp\left[e_+\ln v_++e_-\ln v_-+\sum_{k=1}^{n/2-1}
(e_k\ln \rho_k+\tilde e_k\phi_k)\right] 
=e_+ v_+ +e_- v_- +\sum_{k=1}^{n/2-1}
(e_k v_k+\tilde e_k \tilde v_k),
\label{npolar-47a}
\end{equation}
where the fact has ben used that 
\begin{equation}
e_+ +e_- +\sum_{k=1}^{n/2-1} e_k =1, 
\label{npolar-47b}
\end{equation}
the latter relation being a consequence of Eqs. (\ref{npolar-e11}) and (\ref{npolar-15}).
Similarly, by multiplying the relations (\ref{npolar-46a})-(\ref{npolar-46c}) it results, for
odd $n$, that
\begin{equation}
\exp\left[e_+\ln v_+ +\sum_{k=1}^{(n-1)/2}
(e_k\ln \rho_k+\tilde e_k\phi_k)\right] 
=e_+ v_+  +\sum_{k=1}^{(n-1)/2}
(e_k v_k+\tilde e_k \tilde v_k),
\label{npolar-48a}
\end{equation}
where the fact has ben used that 
\begin{equation}
e_+ +\sum_{k=1}^{(n-1)/2} e_k =1, 
\label{npolar-48b}
\end{equation}
the latter relation being a consequence of Eqs. (\ref{npolar-e12}) and (\ref{npolar-16}).

By comparing Eqs. (\ref{npolar-e13a}) and (\ref{npolar-47a}), it can be seen that, for even
$n$, 
\begin{eqnarray}
\lefteqn{x_0+h_1x_1+\cdots+h_{n-1}x_{n-1}=
\exp\left[e_+\ln v_++e_-\ln v_-+\sum_{k=1}^{n/2-1}
(e_k\ln \rho_k+\tilde e_k\phi_k)\right] ,\nonumber}\\
&&
\label{npolar-49a}
\end{eqnarray}
and by comparing Eqs. (\ref{npolar-e13b}) and (\ref{npolar-48a}), it can be seen that, for
odd $n$,
\begin{eqnarray}
x_0+h_1x_1+\cdots+h_{n-1}x_{n-1}
=\exp\left[e_+\ln v_+ +\sum_{k=1}^{(n-1)/2}
(e_k\ln \rho_k+\tilde e_k\phi_k)\right] .
\label{npolar-49b}
\end{eqnarray}
Using the expression of the bases in Eqs. (\ref{npolar-e11}) and (\ref{npolar-e12}) yields,
for even values of $n$, the exponential form of the n-complex number
$u=x_0+h_1x_1+\cdots+h_{n-1}x_{n-1}$ as
\begin{eqnarray}\lefteqn{
u=\rho\exp\left\{\sum_{p=1}^{n-1}h_p\left[
\frac{1}{n}\ln\frac{\sqrt{2}}{\tan\theta_+}
+\frac{(-1)^p}{n}\ln\frac{\sqrt{2}}{\tan\theta_-}\right.\right.\nonumber}\\
&&
\left.\left.-\frac{2}{n}\sum_{k=2}^{n/2-1}
\cos\left(\frac{2\pi kp}{n}\right)\ln\tan\psi_{k-1}
\right]
+\sum_{k=1}^{n/2-1}\tilde e_k\phi_k 
\right\},
\label{npolar-50a}
\end{eqnarray}\index{exponential form!polar n-complex}
where $\rho$ is the amplitude defined in Eq. (\ref{npolar-6b}), which for even $n$ has
according to Eq. (\ref{npolar-9c}) the expression
\begin{equation}
\rho=\left(v_+v_-\rho_1^2\cdots \rho_{n/2-1}^2\right)^{1/n}.
\label{npolar-50aa}
\end{equation}
For odd values of $n$, the exponential form of the n-complex number $u$ is
\begin{eqnarray}
\lefteqn{
u=\rho\exp\left\{\sum_{p=1}^{n-1}h_p\left[
\frac{1}{n}\ln\frac{\sqrt{2}}{\tan\theta_+}
\right.\right.
\left.\left.-\frac{2}{n}\sum_{k=2}^{(n-1)/2}
\cos\left(\frac{2\pi kp}{n}\right)\ln\tan\psi_{k-1}
\right]
+\sum_{k=1}^{(n-1)/2}\tilde e_k\phi_k
\right\},\nonumber}\\
&&
\label{npolar-50b}
\end{eqnarray}
where for odd $n$, $\rho$ has according to Eq. (\ref{npolar-9d}) the expression
\begin{equation}
\rho=\left(v_+\rho_1^2\cdots \rho_{(n-1)/2}^2\right)^{1/n}.
\label{npolar-50bb}
\end{equation}

It can be checked with the aid of Eq. (\ref{npolar-46a}) that
the n-complex number $u$ can also be written, for even $n$, as
\begin{eqnarray}
x_0+h_1x_1+\cdots+h_{n-1}x_{n-1}
=\left(e_+ v_++e_- v_-+\sum_{k=1}^{n/2-1}e_k \rho_k\right)
\exp\left(\sum_{k=1}^{n/2-1}\tilde e_k\phi_k\right),
\label{npolar-51a}
\end{eqnarray}
and for odd $n$, as
\begin{eqnarray}
x_0+h_1x_1+\cdots+h_{n-1}x_{n-1}
=\left(e_+ v_++\sum_{k=1}^{(n-1)/2}e_k \rho_k\right)
\exp\left(\sum_{k=1}^{(n-1)/2}\tilde e_k\phi_k\right).
\label{npolar-51b}
\end{eqnarray}
Writing in Eqs. (\ref{npolar-51a}) and (\ref{npolar-51b}) the radius $\rho_1$, Eqs.
(\ref{npolar-20b}) and (\ref{npolar-20c}), as a factor and expressing the 
variables in terms of the polar and planar
angles with the aid of Eqs. (\ref{npolar-19b})-(\ref{npolar-19d}) yields the
trigonometric form of the n-complex number $u$, for even $n$, as
\begin{eqnarray}
\lefteqn{u=d
\left(\frac{n}{2}\right)^{1/2}
\left(\frac{1}{\tan^2\theta_+}+\frac{1}{\tan^2\theta_-}+1
+\frac{1}{\tan^2\psi_1}+\frac{1}{\tan^2\psi_2}+\cdots
+\frac{1}{\tan^2\psi_{n/2-2}}\right)^{-1/2}\nonumber}\\
&&\left(\frac{e_+\sqrt{2}}{\tan\theta_+}+\frac{e_-\sqrt{2}}{\tan\theta_-}
+e_1+\sum_{k=2}^{n/2-1}\frac{e_k}{\tan\psi_{k-1}}\right)
\exp\left(\sum_{k=1}^{n/2-1}\tilde e_k\phi_k\right),
\label{npolar-52a}
\end{eqnarray}\index{trigonometric form!polar n-complex}
and for odd $n$ as
\begin{eqnarray}
\lefteqn{u=d
\left(\frac{n}{2}\right)^{1/2}
\left(\frac{1}{\tan^2\theta_+}+1
+\frac{1}{\tan^2\psi_1}+\frac{1}{\tan^2\psi_2}+\cdots
+\frac{1}{\tan^2\psi_{(n-3)/2}}\right)^{-1/2}\nonumber}\\
&&\left(\frac{e_+\sqrt{2}}{\tan\theta_+}
+e_1+\sum_{k=2}^{(n-1)/2}\frac{e_k}{\tan\psi_{k-1}}\right)
\exp\left(\sum_{k=1}^{(n-1)/2}\tilde e_k\phi_k\right).
\label{npolar-52b}
\end{eqnarray}
In Eqs. (\ref{npolar-52a}) and {\ref{npolar-52b}),  the n-complex
number $u$, written in trigonometric form, is the
product of the modulus $d$, of a part depending on the polar and planar
angles $\theta_+, \theta_-,\psi_1,...,\psi_{[(n-3)/2]},$ and of a factor depending
on the azimuthal angles $\phi_1,...,\phi_{[(n-1)/2]}$. 
Although the modulus of a product of n-complex numbers is not equal in
general to the product of the moduli of the factors,
it can be checked that the modulus of the factor in Eq. (\ref{npolar-52a}) is
\begin{eqnarray}
\lefteqn{
\left|\frac{e_+\sqrt{2}}{\tan\theta_+}+\frac{e_-\sqrt{2}}{\tan\theta_-}
+e_1+\sum_{k=2}^{n/2-1}\frac{e_k}{\tan\psi_{k-1}}\right|\nonumber}\\
&&=\left(\frac{2}{n}\right)^{1/2}
\left(\frac{1}{\tan^2\theta_+}+\frac{1}{\tan^2\theta_-}+1
+\frac{1}{\tan^2\psi_1}+\frac{1}{\tan^2\psi_2}+\cdots
+\frac{1}{\tan^2\psi_{n/2-2}}\right)^{1/2},\nonumber\\
&&
\label{npolar-52c}
\end{eqnarray}
and the modulus of the factor in Eq. (\ref{npolar-52b}) is
\begin{eqnarray}
\lefteqn{
\left|\frac{e_+\sqrt{2}}{\tan\theta_+}
+e_1+\sum_{k=2}^{(n-1)/2}\frac{e_k}{\tan\psi_{k-1}}\right|\nonumber}\\
&&=\left(\frac{2}{n}\right)^{1/2}
\left(\frac{1}{\tan^2\theta_+}+1
+\frac{1}{\tan^2\psi_1}+\frac{1}{\tan^2\psi_2}+\cdots
+\frac{1}{\tan^2\psi_{(n-3)/2}}\right)^{1/2}.
\label{npolar-52d}
\end{eqnarray}
Moreover, it can be checked that
\begin{eqnarray}
\left|\exp\left[\sum_{k=1}^{[(n-1)/2]}\tilde e_k\phi_k\right]\right|=1.
\label{npolar-52e}
\end{eqnarray}

The modulus $d$ in Eqs. (\ref{npolar-52a}) and (\ref{npolar-52b}) can be expressed in terms
of the amplitude $\rho$, for even $n$, as
\begin{eqnarray}
\lefteqn{d=\rho \frac{2^{(n-2)/2n}}{\sqrt{n}}
\left(\tan\theta_+\tan\theta_-
\tan^2\psi_1\cdots\tan^2\psi_{n/2-2}\right)^{1/n}\nonumber}\\
&&\left(\frac{1}{\tan^2\theta_+}+\frac{1}{\tan^2\theta_-}+1
+\frac{1}{\tan^2\psi_1}+\frac{1}{\tan^2\psi_2}+\cdots
+\frac{1}{\tan^2\psi_{n/2-2}}\right)^{1/2},
\label{npolar-53a}
\end{eqnarray}
and for odd $n$ as
\begin{eqnarray}
\lefteqn{d=\rho \frac{2^{(n-1)/2n}}{\sqrt{n}}
\left(\tan\theta_+
\tan^2\psi_1\cdots\tan^2\psi_{(n-3)/2}\right)^{1/n}\nonumber}\\
&&\left(\frac{1}{\tan^2\theta_+}+1
+\frac{1}{\tan^2\psi_1}+\frac{1}{\tan^2\psi_2}+\cdots
+\frac{1}{\tan^2\psi_{(n-3)/2}}\right)^{1/2}.
\label{npolar-53b}
\end{eqnarray}

\subsection{Elementary functions of a polar n-complex variable}

The logarithm $u_1$ of the n-complex number $u$, $u_1=\ln u$, can be defined
as the solution of the equation
\begin{equation}
u=e^{u_1} .
\label{npolar-54}
\end{equation}
For even $n$ the relation (\ref{npolar-47a}) shows that $\ln u$ exists 
as an n-complex function with real components if $v_+=x_0+x_1+\cdots+x_{n-1}>0$
and $v_-=x_0-x_1+\cdots+x_{n-2}-x_{n-1}>0$, which means that $0<\theta_+<\pi/2,
0<\theta_-<\pi/2$.
For odd $n$ the relation (\ref{npolar-48a}) shows that $\ln u$ exists 
as an n-complex function with real components if
$v_+=x_0+x_1+\cdots+x_{n-1}>0$, which means that $0<\theta_+<\pi/2$. 
The expression of the logarithm, obtained from Eqs. (\ref{npolar-49a}) and
(\ref{npolar-49b}), is, for even $n$,
\begin{equation}
\ln u=e_+\ln v_++e_-\ln v_-+\sum_{k=1}^{n/2-1}
(e_k\ln \rho_k+\tilde e_k\phi_k),
\label{npolar-55a}
\end{equation}\index{logarithm!polar n-complex}
and for odd $n$ the expression is
\begin{equation}
\ln u=e_+\ln v_+ +\sum_{k=1}^{(n-1)/2}
(e_k\ln \rho_k+\tilde e_k\phi_k).
\label{npolar-55b}
\end{equation}
An expression of the logarithm depending on the amplitude $\rho$ can be
obtained from the exponential forms in Eqs. (\ref{npolar-50a}) and (\ref{npolar-50b}), for
even $n$, as
\begin{eqnarray}\lefteqn{
\ln u=\ln \rho+\sum_{p=1}^{n-1}h_p\left[
\frac{1}{n}\ln\frac{\sqrt{2}}{\tan\theta_+}
+\frac{(-1)^p}{n}\ln\frac{\sqrt{2}}{\tan\theta_-}\right.
\left.-\frac{2}{n}\sum_{k=2}^{n/2-1}
\cos\left(\frac{2\pi kp}{n}\right)\ln\tan\psi_{k-1}
\right]\nonumber}\\
&&+\sum_{k=1}^{n/2-1}\tilde e_k\phi_k,
\label{npolar-56a}
\end{eqnarray}\index{logarithm!polar n-complex}
and for odd $n$ as
\begin{eqnarray}
\lefteqn{\ln u=\ln\rho+\sum_{p=1}^{n-1}h_p\left[
\frac{1}{n}\ln\frac{\sqrt{2}}{\tan\theta_+}
-\frac{2}{n}\sum_{k=2}^{(n-1)/2}
\cos\left(\frac{2\pi kp}{n}\right)\ln\tan\psi_{k-1}
\right]
+\sum_{k=1}^{(n-1)/2}\tilde e_k\phi_k.\nonumber}\\
&&
\label{npolar-56b}
\end{eqnarray}

The function $\ln u$ is multivalued because of the presence of the terms 
$\tilde e_k\phi_k$.
It can be inferred from Eqs. (\ref{npolar-21a})-(\ref{npolar-21g}) and (\ref{npolar-24}) that
\begin{equation}
\ln(uu^\prime)=\ln u+\ln u^\prime ,
\label{npolar-57}
\end{equation}
up to integer multiples of $2\pi\tilde e_k, k=1,...,[(n-1)/2]$.

The power function $u^m$ can be defined for real values of $m$ as
\begin{equation}
u^m=e^{m\ln u} .
\label{npolar-58}
\end{equation}
Using the expression of $\ln u$ in Eqs. (\ref{npolar-55a}) and (\ref{npolar-55b}) yields, for
even values of $n$,
\begin{equation}
u^m=e_+ v_+^m+e_- v_-^m +\sum_{k=1}^{n/2-1}
\rho_k^m(e_k\cos m\phi_k+\tilde e_k\sin m\phi_k),
\label{npolar-59a}
\end{equation}\index{power function!polar n-complex}
and for odd values of $n$
\begin{equation}
u^m=e_+ v_+^m +\sum_{k=1}^{(n-1)/2}
\rho_k^m(e_k\cos m\phi_k+\tilde e_k\sin m\phi_k).
\label{npolar-59b}
\end{equation}
For integer values of $m$, the relations (\ref{npolar-59a}) and (\ref{npolar-59b}) are valid 
for any $x_0,...,x_{n-1}$.
The power function is multivalued unless $m$ is an integer. 
For integer $m$, it can be inferred from Eq. (\ref{npolar-57}) that
\begin{equation}
(uu^\prime)^m=u^m\:u^{\prime m} .
\label{npolar-59}
\end{equation}

The trigonometric functions of the hypercomplex variable
$u$ and the addition theorems for these functions have been written in Eqs.
~(1.57)-(1.60). 
In order to obtain expressions for the trigonometric functions of n-complex
variables, these will be expressed with the aid of the imaginary unit $i$ as
\begin{equation}
\cos u=\frac{1}{2}(e^{iu}+e^{-iu}),\:\sin u=\frac{1}{2i}(e^{iu}-e^{-iu}).
\label{npolar-64}
\end{equation}
The imaginary unit $i$ is used for the convenience of notations, and it does
not appear in the final results.
The validity of Eq. (\ref{npolar-64}) can be checked by comparing the series for the
two sides of the relations.
Since the expression of the exponential function $e^{h_k y}$ in terms of the
units $1, h_1, ... h_{n-1}$ given in Eq. (\ref{npolar-28b}) depends on the polar
cosexponential functions $g_{np}(y)$, the expression of the trigonometric
functions will depend on the functions
$g_{p+}^{(c)}(y)=(1/2)[g_{np}(iy)+g_{np}(-iy)]$ and 
$g_{p-}^{(c)}(y)=(1/2i)[g_{np}(iy)-g_{np}(-iy)]$, 
\begin{equation}
\cos(h_k y)=\sum_{p=0}^{n-1}h_{kp-n[kp/n]}g_{p+}^{(c)}(y),
\label{npolar-66a}
\end{equation}\index{trigonometric functions, expressions!polar n-complex}
\begin{equation}
\sin(h_k y)=\sum_{p=0}^{n-1}h_{kp-n[kp/n]}g_{p-}^{(c)}(y),
\label{npolar-66b}
\end{equation}
where
\begin{eqnarray}
\lefteqn{g_{p+}^{(c)}(y)=\frac{1}{n}\sum_{l=0}^{n-1}\left\{
\cos\left[y\cos\left(\frac{2\pi l}{n}\right)\right]
\cosh\left[y\sin\left(\frac{2\pi l}{n}\right)\right]
\cos\left(\frac{2\pi lp}{n}\right)\right.\nonumber}\\
&&\left.-\sin\left[y\cos\left(\frac{2\pi l}{n}\right)\right]
\sinh\left[y\sin\left(\frac{2\pi l}{n}\right)\right]
\sin\left(\frac{2\pi lp}{n}\right)
\right\},
\label{npolar-65a}
\end{eqnarray}
\begin{eqnarray}
\lefteqn{g_{p-}^{(c)}(y)=\frac{1}{n}\sum_{l=0}^{n-1}\left\{
\sin\left[y\cos\left(\frac{2\pi l}{n}\right)\right]
\cosh\left[y\sin\left(\frac{2\pi l}{n}\right)\right]
\cos\left(\frac{2\pi lp}{n}\right)\right.\nonumber}\\
&&\left.+\cos\left[y\cos\left(\frac{2\pi l}{n}\right)\right]
\sinh\left[y\sin\left(\frac{2\pi l}{n}\right)\right]
\sin\left(\frac{2\pi lp}{n}\right)
\right\}.
\label{npolar-65b}
\end{eqnarray}

The hyperbolic functions of the hypercomplex variable
$u$ and the addition theorems for these functions have been written in Eqs.
~(1.62)-(1.65). 
In order to obtain expressions for the hyperbolic functions of n-complex
variables, these will be expressed as
\begin{equation}
\cosh u=\frac{1}{2}(e^{u}+e^{-u}),\:\sinh u=\frac{1}{2}(e^{u}-e^{-u}).
\label{npolar-70}
\end{equation}
The validity of Eq. (\ref{npolar-70}) can be checked by comparing the series for the
two sides of the relations.
Since the expression of the exponential function $e^{h_k y}$ in terms of the
units $1, h_1, ... h_{n-1}$ given in Eq. (\ref{npolar-28b}) depends on the polar
cosexponential functions $g_{np}(y)$, the expression of the hyperbolic
functions will depend on the even part $g_{p+}(y)=(1/2)[g_{np}(y)+g_{np}(-y)]$ and on
the odd part $g_{p-}(y)=(1/2)[g_{np}(y)-g_{np}(-y)]$ of $g_{np}$, 
\begin{equation}
\cosh(h_k y)=\sum_{p=0}^{n-1}h_{kp-n[kp/n]}g_{p+}(y),
\label{npolar-71a}
\end{equation}\index{hyperbolic functions, expressions!polar n-complex}
\begin{equation}
\sinh(h_k y)=\sum_{p=0}^{n-1}h_{kp-n[kp/n]}g_{p-}(y),
\label{npolar-71b}
\end{equation}
where
\begin{eqnarray}
\lefteqn{g_{p+}(y)=\frac{1}{n}\sum_{l=0}^{n-1}\left\{
\cosh\left[y\cos\left(\frac{2\pi l}{n}\right)\right]
\cos\left[y\sin\left(\frac{2\pi l}{n}\right)\right]
\cos\left(\frac{2\pi lp}{n}\right)\right.\nonumber}\\
&&\left.+\sinh\left[y\cos\left(\frac{2\pi l}{n}\right)\right]
\sin\left[y\sin\left(\frac{2\pi l}{n}\right)\right]
\sin\left(\frac{2\pi lp}{n}\right)
\right\},
\label{npolar-72a}
\end{eqnarray}
\begin{eqnarray}
\lefteqn{g_{p-}(y)=\frac{1}{n}\sum_{l=0}^{n-1}\left\{
\sinh\left[y\cos\left(\frac{2\pi l}{n}\right)\right]
\cos\left[y\sin\left(\frac{2\pi l}{n}\right)\right]
\cos\left(\frac{2\pi lp}{n}\right)\right.\nonumber}\\
&&\left.+\cosh\left[y\cos\left(\frac{2\pi l}{n}\right)\right]
\sin\left[y\sin\left(\frac{2\pi l}{n}\right)\right]
\sin\left(\frac{2\pi lp}{n}\right)
\right\}.
\label{npolar-72b}
\end{eqnarray}

The exponential, trigonometric and hyperbolic functions can also be expressed
with the aid of the bases introduced in Eqs. (\ref{npolar-e11}) and (\ref{npolar-e12}).
Using the expression of the n-complex number in Eq. (\ref{npolar-e13a}), for even $n$,
yields for the exponential of the n-complex variable $u$
\begin{eqnarray}
e^u=e_+e^{v_+} + e_-e^{v_-} 
+\sum_{k=1}^{n/2-1}e^{v_k}\left(e_k \cos \tilde v_k+\tilde e_k \sin\tilde
v_k\right).
\label{npolar-73a}
\end{eqnarray}\index{exponential, expressions!polar n-complex}
For odd $n$, the expression of the n-complex variable in Eq. (\ref{npolar-e13b})
yileds for the exponential 
\begin{eqnarray}
e^u=e_+e^{v_+}  
+\sum_{k=1}^{(n-1)/2}e^{v_k}\left(e_k \cos \tilde v_k+\tilde e_k \sin\tilde
v_k\right).
\label{npolar-73b}
\end{eqnarray}

The trigonometric functions can be obtained from Eqs. (\ref{npolar-73a}) and
(\ref{npolar-73b} with the aid of Eqs. (\ref{npolar-64}). The trigonometric functions of the
n-complex variable $u$ are, for even $n$,
\begin{equation}
\cos u=e_+\cos v_+ + e_-\cos v_- 
+\sum_{k=1}^{n/2-1}\left(e_k \cos v_k\cosh \tilde v_k
-\tilde e_k \sin v_k\sinh\tilde v_k\right),
\label{npolar-74a}
\end{equation}\index{trigonometric functions, expressions!polar n-complex}
\begin{equation}
\sin u=e_+\sin v_+ + e_-\sin v_- 
+\sum_{k=1}^{n/2-1}\left(e_k \sin v_k\cosh \tilde v_k
+\tilde e_k \cos v_k\sinh\tilde v_k\right),
\label{npolar-74b}
\end{equation}
and for odd $n$ the trigonometric functions are
\begin{equation}
\cos u=e_+\cos v_+  
+\sum_{k=1}^{(n-1)/2}\left(e_k \cos v_k\cosh \tilde v_k
-\tilde e_k \sin v_k\sinh\tilde v_k\right),
\label{npolar-74c}
\end{equation}
\begin{equation}
\sin u=e_+\sin v_+  
+\sum_{k=1}^{(n-1)/2}\left(e_k \sin v_k\cosh \tilde v_k
+\tilde e_k \cos v_k\sinh\tilde v_k\right).
\label{npolar-74d}
\end{equation}

The hyperbolic functions can be obtained from Eqs. (\ref{npolar-73a}) and
(\ref{npolar-73b} with the aid of Eqs. (\ref{npolar-70}). The hyperbolic functions of the
n-complex variable $u$ are, for even $n$,
\begin{equation}
\cosh u=e_+\cosh v_+ + e_-\cosh v_- 
+\sum_{k=1}^{n/2-1}\left(e_k \cosh v_k\cos \tilde v_k
+\tilde e_k \sinh v_k\sin\tilde v_k\right),
\label{npolar-75a}
\end{equation}\index{hyperbolic functions, expressions!polar n-complex}
\begin{equation}
\sinh u=e_+\sinh v_+ + e_-\sinh v_- 
+\sum_{k=1}^{n/2-1}\left(e_k \sinh v_k\cos \tilde v_k
+\tilde e_k \cosh v_k\sin\tilde v_k\right),
\label{npolar-75b}
\end{equation}
and for odd $n$ the hyperbolic functions are
\begin{equation}
\cosh u=e_+\cosh v_+  
+\sum_{k=1}^{(n-1)/2}\left(e_k \cosh v_k\cos \tilde v_k
+\tilde e_k \sinh v_k\sin\tilde v_k\right),
\label{npolar-75c}
\end{equation}
\begin{equation}
\sinh u=e_+\sinh v_+  
+\sum_{k=1}^{(n-1)/2}\left(e_k \sinh v_k\cos \tilde v_k
+\tilde e_k \cosh v_k\sin\tilde v_k\right).
\label{npolar-75d}
\end{equation}

\subsection{Power series of polar n-complex numbers}

An n-complex series is an infinite sum of the form
\begin{equation}
a_0+a_1+a_2+\cdots+a_n+\cdots , 
\label{npolar-76}
\end{equation}\index{series!polar n-complex}
where the coefficients $a_n$ are n-complex numbers. The convergence of 
the series (\ref{npolar-76}) can be defined in terms of the convergence of its $n$
real components. The convergence of a n-complex series can also be studied
using n-complex variables. The main criterion for absolute convergence 
remains the comparison theorem, but this requires a number of inequalities
which will be discussed further.

The modulus $d=|u|$ of an n-complex number $u$ has been defined in Eq.
(\ref{npolar-10}). Since $|x_0|\leq |u|, |x_1|\leq |u|,..., |x_{n-1}|\leq |u|$, a
property of absolute convergence established via a comparison theorem based on
the modulus of the series (\ref{npolar-76}) will ensure the absolute convergence of
each real component of that series.

The modulus of the sum $u_1+u_2$ of the n-complex numbers $u_1, u_2$ fulfils
the inequality
\begin{equation}
||u^\prime|-|u^{\prime\prime}||\leq |u^\prime+u^{\prime\prime}|\leq
|u^\prime|+|u^{\prime\prime}| . 
\label{npolar-78}
\end{equation}\index{modulus, inequalities!polar n-complex}
For the product, the relation is 
\begin{equation}
|u^\prime u^{\prime\prime}|\leq \sqrt{n}|u^\prime||u^{\prime\prime}| ,
\label{npolar-79}
\end{equation}
as can be shown from Eqs. (\ref{npolar-17}) and (\ref{npolar-18}). The relation (\ref{npolar-79})
replaces the relation of equality extant between 2-dimensional regular complex
numbers. The equality in Eq. (\ref{npolar-79}) takes place for
$\rho_1\rho_1^\prime=0,..., \rho_{[(n-1)/2]}\rho_{[(n-1)/2]}^\prime=0$ and, for
even $n$, for $v_+v_-^\prime=0$, $v_-v_+^\prime=0$.  

From Eq. (\ref{npolar-79}) it results, for $u=u^\prime$, that
\begin{equation}
|u^2|\leq \sqrt{n} |u|^2 .
\label{npolar-80}
\end{equation}
The relation in Eq. (\ref{npolar-80}) becomes an equality for 
$\rho_1=0,...,\rho_{[(n-1)/2]}=0$ and, for even $n$, $v_+=0$ or $v_-=0$.
The inequality in Eq. (\ref{npolar-79}) implies that
\begin{equation}
|u^l|\leq n^{(l-1)/2}|u|^l ,
\label{npolar-81}
\end{equation}
where $l$ is a natural number.
From Eqs. (\ref{npolar-79}) and (\ref{npolar-81}) it results that
\begin{equation}
|au^l|\leq n^{l/2} |a| |u|^l .
\label{npolar-82}
\end{equation}

A power series of the n-complex variable $u$ is a series of the form
\begin{equation}
a_0+a_1 u + a_2 u^2+\cdots +a_l u^l+\cdots .
\label{npolar-83}
\end{equation}\index{power series!polar n-complex}
Since
\begin{equation}
\left|\sum_{l=0}^\infty a_l u^l\right| \leq  \sum_{l=0}^\infty
n^{l/2} |a_l| |u|^l ,
\label{npolar-84}
\end{equation}
a sufficient condition for the absolute convergence of this series is that
\begin{equation}
\lim_{l\rightarrow \infty}\frac{\sqrt{n}|a_{l+1}||u|}{|a_l|}<1 .
\label{npolar-85}
\end{equation}
Thus the series is absolutely convergent for 
\begin{equation}
|u|<c,
\label{npolar-86}
\end{equation}\index{convergence of power series!polar n-complex}
where 
\begin{equation}
c=\lim_{l\rightarrow\infty} \frac{|a_l|}{\sqrt{n}|a_{l+1}|} .
\label{npolar-87}
\end{equation}

The convergence of the series (\ref{npolar-83}) can be also studied with the aid of
the formulas (\ref{npolar-59a}), (\ref{npolar-59b}) which for integer values of $m$ are
valid for any values of $x_0,...,x_{n-1}$, as mentioned previously.
If $a_l=\sum_{p=0}^{n-1}h_p a_{lp}$, and
\begin{equation}
A_{l+}=\sum_{p=0}^{n-1}a_{lp},
\label{npolar-88a}
\end{equation}
\begin{equation}
A_{lk}=\sum_{p=0}^{n-1}a_{lp}\cos\frac{2\pi kp}{n},
\label{npolar-88b}
\end{equation}
\begin{equation}
\tilde A_{lk}=\sum_{p=0}^{n-1}a_{lp}\sin\frac{2\pi kp}{n},
\label{npolar-88c}
\end{equation}
for $k=1,...,[(n-1)/2]$, and for even $n$ 
\begin{equation}
A_{l-}=\sum_{p=0}^{n-1}(-1)^p a_{lp},
\label{npolar-88d}
\end{equation}
the series (\ref{npolar-83}) can be written, for even $n$, as
\begin{equation}
\sum_{l=0}^\infty \left[
e_+A_{l+}v_+^l+e_-A_{l-}v_-^l+\sum_{k=1}^{n/2-1}
(e_k A_{lk}+\tilde e_k\tilde A_{lk})(e_k v_k+\tilde e_k\tilde v_k)^l 
\right],
\label{npolar-89a}
\end{equation}
and for odd $n$ as
\begin{equation}
\sum_{l=0}^\infty \left[
e_+A_{l+}v_+^l+\sum_{k=1}^{(n-1)/2}
(e_k A_{lk}+\tilde e_k\tilde A_{lk})(e_k v_k+\tilde e_k\tilde v_k)^l 
\right].
\label{npolar-89b}
\end{equation}

The series in Eq. (\ref{npolar-83}) is absolutely convergent for 
\begin{equation}
|v_+|<c_+,\:
|v_-|<c_-,\:
\rho_k<c_k, 
\label{npolar-90}
\end{equation}\index{convergence, region of!polar n-complex}
for $k=1,..., [(n-1)/2]$, where 
\begin{equation}
c_+=\lim_{l\rightarrow\infty} \frac{|A_{l+}|}{|A_{l+1,+}|} ,\:
c_-=\lim_{l\rightarrow\infty} \frac{|A_{l-}|}{|A_{l+1,-}|} ,\:
c_k=\lim_{l\rightarrow\infty} \frac
{\left(A_{lk}^2+\tilde A_{lk}^2\right)^{1/2}}
{\left(A_{l+1,k}^2+\tilde A_{l+1,k}^2\right)^{1/2}} .
\label{npolar-91}
\end{equation}
The relations (\ref{npolar-90}) show that the region of convergence of the series
(\ref{npolar-83}) is an n-dimensional cylinder.

It can be shown that, for even $n$, $c=(1/\sqrt{n})\;{\rm
min}(c_+,c_-,c_1,...,c_{n/2-1})$, and for odd $n$ $c=(1/\sqrt{n})\;{\rm
min}(c_+,c_1,...,c_{(n-1)/2})$, where ${\rm min}$ designates the smallest of
the numbers in the argument of this function. Using the expression of $|u|$ in
Eqs. 
(\ref{npolar-17}) or (\ref{npolar-18}), it can be seen that the spherical region of
convergence defined in Eqs. (\ref{npolar-86}), (\ref{npolar-87}) is a subset of the
cylindrical region of convergence defined in Eqs. (\ref{npolar-90}) and (\ref{npolar-91}).

\subsection{Analytic functions of polar n-complex variables}

The analytic functions of the hypercomplex variable $u$ and the series 
expansion of functions have been discussed in Eqs. ~(1.85)-(1.93).
If the n-complex function $f(u)$
of the n-complex variable $u$ is written in terms of 
the real functions $P_k(x_0,...,x_{n-1}), k=0,1,...,n-1$ of the real
variables $x_0,x_1,...,x_{n-1}$ as 
\begin{equation}
f(u)=\sum_{k=0}^{n-1}h_kP_k(x_0,...,x_{n-1}),
\label{npolar-h93}
\end{equation}\index{functions, real components!polar n-complex}
then relations of equality 
exist between the partial derivatives of the functions $P_k$. 
The derivative of the function $f$ can be written as
\begin{eqnarray}
\lim_{\Delta u\rightarrow 0}\frac{1}{\Delta u} 
\sum_{k=0}^{n-1}\left(h_k\sum_{l=0}^{n-1}
\frac{\partial P_k}{\partial x_l}\Delta x_l\right),
\label{npolar-h94}
\end{eqnarray}\index{derivative, independence of direction!polar n-complex}
where
\begin{equation}
\Delta u=\sum_{k=0}^{n-1}h_l\Delta x_l.
\label{npolar-h94a}
\end{equation}
The
relations between the partials derivatives of the functions $P_k$ are
obtained by setting successively in   
Eq. (\ref{npolar-h94}) $\Delta u=h_l\Delta x_l$, for $l=0,1,...,n-1$, and equating the
resulting expressions. 
The relations are \index{relations between partial derivatives!polar n-complex}
\begin{equation}
\frac{\partial P_k}{\partial x_0} = \frac{\partial P_{k+1}}{\partial x_1} 
=\cdots=\frac{\partial P_{n-1}}{\partial x_{n-k-1}} 
= \frac{\partial P_0}{\partial x_{n-k}}=\cdots
=\frac{\partial P_{k-1}}{\partial x_{n-1}}, 
\label{npolar-h95}
\end{equation}
for $k=0,1,...,n-1$.
The relations (\ref{npolar-h95}) are analogous to the Riemann relations
for the real and imaginary components of a complex function. 
It can be shown from Eqs. (\ref{npolar-h95}) that the components $P_k$ fulfil the
second-order equations\index{relations between second-order derivatives!polar n-complex}
\begin{eqnarray}
\lefteqn{\frac{\partial^2 P_k}{\partial x_0\partial x_l}
=\frac{\partial^2 P_k}{\partial x_1\partial x_{l-1}}
=\cdots=
\frac{\partial^2 P_k}{\partial x_{[l/2]}\partial x_{l-[l/2]}}}\nonumber\\
&&=\frac{\partial^2 P_k}{\partial x_{l+1}\partial x_{n-1}}
=\frac{\partial^2 P_k}{\partial x_{l+2}\partial x_{n-2}}
=\cdots
=\frac{\partial^2 P_k}{\partial x_{l+1+[(n-l-2)/2]}
\partial x_{n-1-[(n-l-2)/2]}} ,
\label{npolar-96}
\end{eqnarray}
for $k,l=0,1,...,n-1$.

\subsection{Integrals of polar n-complex functions}

The singularities of n-complex functions arise from terms of the form
$1/(u-u_0)^n$, with $n>0$. Functions containing such terms are singular not
only at $u=u_0$, but also at all points of the hypersurfaces
passing through the pole $u_0$ and which are parallel to the nodal hypersurfaces. 

The integral of an n-complex function between two points $A, B$ along a path
situated in a region free of singularities is independent of path, which means
that the integral of an analytic function along a loop situated in a region
free of singularities is zero,
\begin{equation}
\oint_\Gamma f(u) du = 0,
\label{npolar-111}
\end{equation}
where it is supposed that a surface $\Sigma$ spanning 
the closed loop $\Gamma$ is not intersected by any of
the hypersurfaces associated with the
singularities of the function $f(u)$. Using the expression, Eq. (\ref{npolar-h93}),
for $f(u)$ and the fact that 
\begin{eqnarray}
du=\sum_{k=0}^{n-1}h_k dx_k, 
\label{npolar-111a}
\end{eqnarray}
the explicit form of the integral in Eq. (\ref{npolar-111}) is
\begin{eqnarray}
\oint _\Gamma f(u) du = \oint_\Gamma
\sum_{k=0}^{n-1}h_k\sum_{l=0}^{n-1}P_l dx_{k-l+n[(n-k-1+l)/n]}.
\label{npolar-112}
\end{eqnarray}\index{integrals, path!polar n-complex}

If the functions $P_k$ are regular on a surface $\Sigma$
spanning the loop $\Gamma$,
the integral along the loop $\Gamma$ can be transformed in an integral over the
surface $\Sigma$ of terms of the form
$\partial P_l/\partial x_{k-m+n[(n-k+m-1)/n]} 
-  \partial P_m/\partial x_{k-l+n[(n-k+l-1)/n]}$.
These terms are equal to zero by Eqs. (\ref{npolar-h95}), and this
proves Eq. (\ref{npolar-111}). 

The integral of the function $(u-u_0)^m$ on a closed loop $\Gamma$ is equal to
zero for $m$ a positive or negative integer not equal to -1,
\begin{equation}
\oint_\Gamma (u-u_0)^m du = 0, \:\: m \:\:{\rm integer},\: m\not=-1 .
\label{npolar-112b}
\end{equation}
This is due to the fact that $\int (u-u_0)^m du=(u-u_0)^{m+1}/(m+1), $ and to
the fact that the function $(u-u_0)^{m+1}$ is singlevalued for $m$ an integer.

The integral $\oint_\Gamma du/(u-u_0)$ can be calculated using the exponential
form, Eqs. (\ref{npolar-50a}) and (\ref{npolar-50b}), for the difference $u-u_0$, which for
even $n$ is 
\begin{eqnarray}\lefteqn{
u-u_0=\rho\exp\left\{\sum_{p=1}^{n-1}h_p\left[
\frac{1}{n}\ln\frac{\sqrt{2}}{\tan\theta_+}
+\frac{(-1)^p}{n}\ln\frac{\sqrt{2}}{\tan\theta_-}\right.\right.\nonumber}\\
&&
\left.\left.-\frac{2}{n}\sum_{k=2}^{n/2-1}
\cos\left(\frac{2\pi kp}{n}\right)\ln\tan\psi_{k-1}
\right]
+\sum_{k=1}^{n/2-1}\tilde e_k\phi_k\right\},
\label{npolar-113a}
\end{eqnarray}
and for odd $n$ is
\begin{eqnarray}
\lefteqn{
u-u_0=\rho\exp\left\{\sum_{p=1}^{n-1}h_p\left[
\frac{1}{n}\ln\frac{\sqrt{2}}{\tan\theta_+}
\right.\right.
\left.\left.-\frac{2}{n}\sum_{k=2}^{(n-1)/2}
\cos\left(\frac{2\pi kp}{n}\right)\ln\tan\psi_{k-1}
\right]
+\sum_{k=1}^{(n-1)/2}\tilde e_k\phi_k\right\}.\nonumber}\\
&&
\label{npolar-113b}
\end{eqnarray}
Thus for even $n$ the quantity $du/(u-u_0)$ is
\begin{eqnarray}\lefteqn{
\frac{du}{u-u_0}=
\frac{d\rho}{\rho}
+\sum_{p=1}^{n-1}h_p\left[
\frac{1}{n}d\ln\frac{\sqrt{2}}{\tan\theta_+}
+\frac{(-1)^p}{n}d\ln\frac{\sqrt{2}}{\tan\theta_-}\right.\nonumber}\\
&&
\left.-\frac{2}{n}\sum_{k=2}^{n/2-1}
\cos\left(\frac{2\pi kp}{n}\right)d\ln\tan\psi_{k-1}
\right]
+\sum_{k=1}^{n/2-1}\tilde e_kd\phi_k,
\label{npolar-114a}
\end{eqnarray}
and for odd $n$
\begin{eqnarray}
\lefteqn{
\frac{du}{u-u_0}=\frac{d\rho}{\rho}
+\sum_{p=1}^{n-1}h_p\left[
\frac{1}{n}d\ln\frac{\sqrt{2}}{\tan\theta_+}
\right.
\left.-\frac{2}{n}\sum_{k=2}^{(n-1)/2}
\cos\left(\frac{2\pi kp}{n}\right)d\ln\tan\psi_{k-1}
\right]
+\sum_{k=1}^{(n-1)/2}\tilde e_kd\phi_k
.\nonumber}\\
&&
\label{npolar-114b}
\end{eqnarray}
Since $\rho, \ln(\sqrt{2}/\tan\theta_+),\ln(\sqrt{2}/\tan\theta_-),
\ln(\tan\psi_{k-1})$ are singlevalued variables, it follows that
$\oint_\Gamma d\rho/\rho =0, 
\oint_\Gamma d(\ln\sqrt{2}/\tan\theta_+)=0,
\oint_\Gamma d(\ln\sqrt{2}/\tan\theta_-)=0,
\oint_\Gamma d(\ln\tan\psi_{k-1})=0$.
On the other hand since,
$\phi_k$ are cyclic variables, they may give contributions to
the integral around the closed loop $\Gamma$.

The expression of $\oint_\Gamma du/(u-u_0)$ can be written 
with the aid of a functional which will be called int($M,C$), defined for a
point $M$ and a closed curve $C$ in a two-dimensional plane, such that 
\begin{equation}
{\rm int}(M,C)=\left\{
\begin{array}{l}
1 \;\:{\rm if} \;\:M \;\:{\rm is \;\:an \;\:interior \;\:point \;\:of} \;\:C ,\\ 
0 \;\:{\rm if} \;\:M \;\:{\rm is \;\:exterior \;\:to}\:\; C .\\
\end{array}\right.
\label{npolar-118}
\end{equation}
With this notation the result of the integration on a closed path $\Gamma$
can be written as 
\begin{equation}
\oint_\Gamma\frac{du}{u-u_0}=
\sum_{k=1}^{[(n-1)/2]}2\pi\tilde e_k 
\;{\rm int}(u_{0\xi_k\eta_k},\Gamma_{\xi_k\eta_k}) ,
\label{npolar-119}
\end{equation}\index{poles and residues!polar n-complex}
where $u_{0\xi_k\eta_k}$ and $\Gamma_{\xi_k\eta_k}$ are respectively the
projections of the point $u_0$ and of 
the loop $\Gamma$ on the plane defined by the axes $\xi_k$ and $\eta_k$,
as shown in Fig. \ref{fig23}.

\begin{figure}
\begin{center}
\epsfig{file=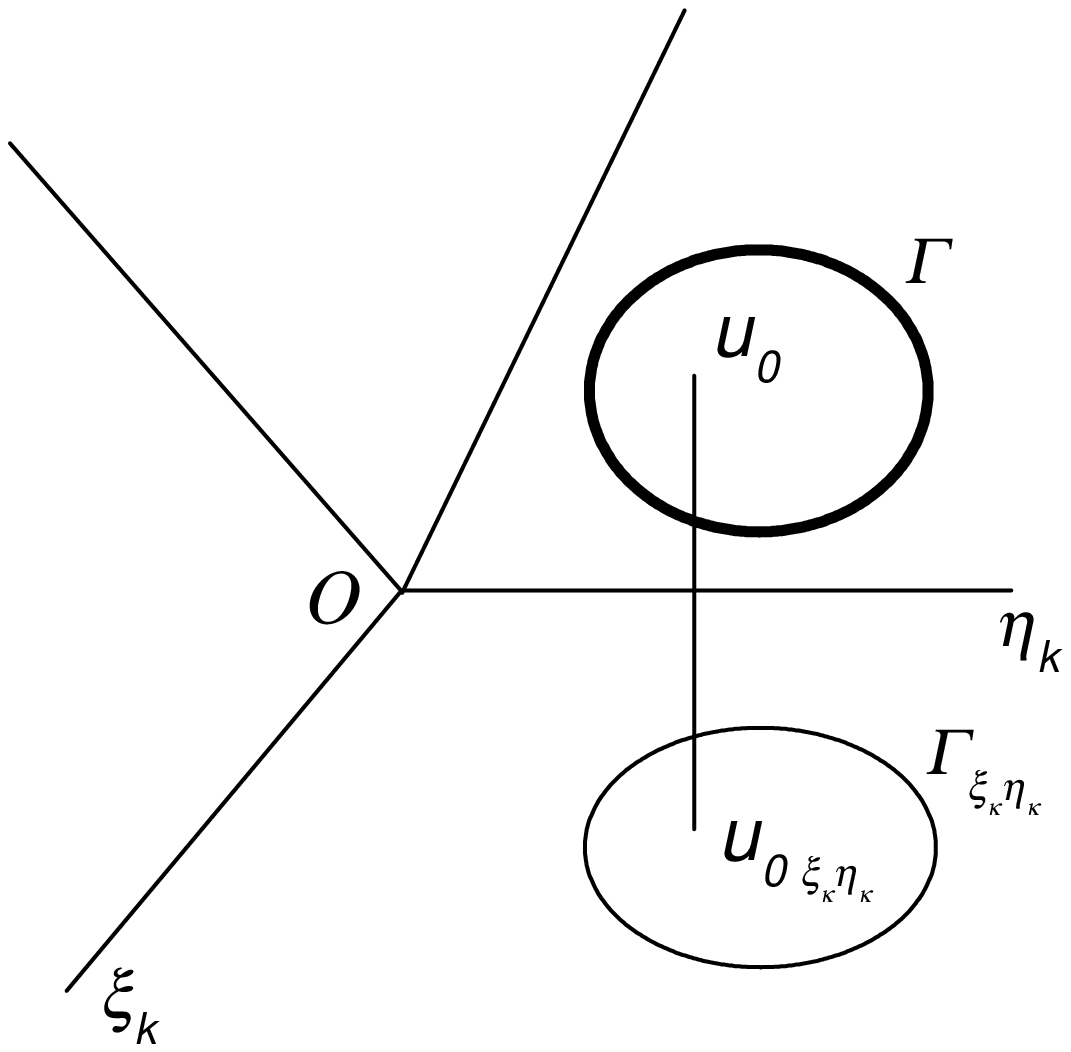,width=12cm}
\caption{Integration path $\Gamma$ and pole $u_0$, and their projections
$\Gamma_{\xi_k\eta_k}$ and $u_{0\xi_k\eta_k}$ on the plane $\xi_k \eta_k$. }
\label{fig23}
\end{center}
\end{figure}

If $f(u)$ is an analytic n-complex function which can be expanded in a
series as written in Eq. (1.89), and the expansion holds on the curve
$\Gamma$ and on a surface spanning $\Gamma$, then from Eqs. (\ref{npolar-112b}) and
(\ref{npolar-119}) it follows that
\begin{equation}
\oint_\Gamma \frac{f(u)du}{u-u_0}=
2\pi f(u_0)\sum_{k=1}^{[(n-1)/2]}\tilde e_k 
\;{\rm int}(u_{0\xi_k\eta_k},\Gamma_{\xi_k\eta_k}) .
\label{npolar-120}
\end{equation}

Substituting in the right-hand side of 
Eq. (\ref{npolar-120}) the expression of $f(u)$ in terms of the real 
components $P_k$, Eq. (\ref{npolar-h93}), yields
\begin{equation}
\oint_\Gamma \frac{f(u)du}{u-u_0}=
\frac{2}{n}\sum_{k=1}^{[(n-1)/2]}\sum_{l,m=0}^{n-1}
h_l\sin\left[\frac{2\pi(l-m)k}{n}\right] P_m(u_0)
\;{\rm int}(u_{0\xi_k\eta_k},\Gamma_{\xi_k\eta_k}) .
\label{npolar-121}
\end{equation}
It the integral in Eq. (\ref{npolar-121}) is written as 
\begin{equation}
\oint_\Gamma \frac{f(u)du}{u-u_0}=\sum_{l=0}^{n-1}h_l I_l,
\label{npolar-122a}
\end{equation}
it can be checked that
\begin{equation}
\sum_{l=0}^{n-1} I_l=0.
\label{npolar-122b}
\end{equation}

If $f(u)$ can be expanded as written in Eq. (1.89) on 
$\Gamma$ and on a surface spanning $\Gamma$, then from Eqs. (\ref{npolar-112b}) and
(\ref{npolar-119}) it also results that
\begin{equation}
\oint_\Gamma \frac{f(u)du}{(u-u_0)^{n+1}}=
\frac{2\pi}{n!}f^{(n)}(u_0)\sum_{k=1}^{[(n-1)/2]}\tilde e_k 
\;{\rm int}(u_{0\xi_k\eta_k},\Gamma_{\xi_k\eta_k}) ,
\label{npolar-122}
\end{equation}
where the fact has been used  that the derivative $f^{(n)}(u_0)$ is related to
the expansion coefficient in Eq. (1.89) according to Eq. (1.93).

If a function $f(u)$ is expanded in positive and negative powers of $u-u_l$,
where $u_l$ are n-complex constants, $l$ being an index, the integral of $f$
on a closed loop $\Gamma$ is determined by the terms in the expansion of $f$
which are of the form $r_l/(u-u_l)$,
\begin{equation}
f(u)=\cdots+\sum_l\frac{r_l}{u-u_l}+\cdots.
\label{npolar-123}
\end{equation}
Then the integral of $f$ on a closed loop $\Gamma$ is
\begin{equation}
\oint_\Gamma f(u) du = 
2\pi \sum_l\sum_{k=1}^{[(n-1)/2]}\tilde e_k 
\;{\rm int}(u_{l\xi_k\eta_k},\Gamma_{\xi_k\eta_k})r_l .
\label{npolar-124}
\end{equation}

\subsection{Factorization of polar n-complex polynomials}

A polynomial of degree $m$ of the n-complex variable $u$ has the form
\begin{equation}
P_m(u)=u^m+a_1 u^{m-1}+\cdots+a_{m-1} u +a_m ,
\label{npolar-125}
\end{equation}
where $a_l$, for $l=1,...,m$, are in general n-complex constants.
If $a_l=\sum_{p=0}^{n-1}h_p a_{lp}$, and with the
notations of Eqs. (\ref{npolar-88a})-(\ref{npolar-88d}) applied for $l= 1, \cdots, m$, the
polynomial $P_m(u)$ can be written, for even $n$, as 
\begin{eqnarray}
\lefteqn{P_m= 
e_+\left(v_+^m +\sum_{l=1}^{m}A_{l+}v_+^{m-l} \right)
+e_-\left(v_-^m +\sum_{l=1}^{m}A_{l-}v_-^{m-l} \right) \nonumber}\\
&&+\sum_{k=1}^{n/2-1}
\left[(e_k v_k+\tilde e_k\tilde v_k)^m+
\sum_{l=1}^m(e_k A_{lk}+\tilde e_k\tilde A_{lk})
(e_k v_k+\tilde e_k\tilde v_k)^{m-l} 
\right],
\label{npolar-126a}
\end{eqnarray}\index{polynomial, canonical variables!polar n-complex}
where the constants $A_{l+}, A_{l-}, A_{lk}, \tilde A_{lk}$ 
are real numbers.
For odd $n$ the expression of the polynomial is
\begin{eqnarray}
\lefteqn{P_m= 
e_+\left(v_+^m +\sum_{l=1}^{m}A_{l+}v_+^{m-l} \right) \nonumber}\\
&&+\sum_{k=1}^{(n-1)/2}
\left[(e_k v_k+\tilde e_k\tilde v_k)^m+
\sum_{l=1}^m(e_k A_{lk}+\tilde e_k\tilde A_{lk})
(e_k v_k+\tilde e_k\tilde v_k)^{m-l} 
\right].
\label{npolar-126b}
\end{eqnarray}

The polynomials of degree $m$ in $e_k v_k+\tilde e_k\tilde v_k$ 
in Eqs. (\ref{npolar-126a}) and (\ref{npolar-126b})
can always be written as a product of linear factors of the form
$e_k (v_k-v_{kp})+\tilde e_k(\tilde v_k-\tilde v_{kp})$, where the
constants $v_{kp}, \tilde v_{kp}$ are real.
The polynomials of degree $m$ with real coefficients in Eqs. (\ref{npolar-126a}) 
and (\ref{npolar-126b})
which are multiplied by $e_+$ and $e_-$ can be written as a product
of linear or quadratic factors with real coefficients, or as a product of
linear factors which, if imaginary, appear always in complex conjugate pairs.
Using the latter form for the simplicity of notations, the polynomial $P_m$
can be written, for even $n$, as
\begin{equation}
P_m=e_+\prod_{p=1}^m (v_+ -v_{p+})+e_-\prod_{p=1}^m (v_- -v_{p-})
+\sum_{k=1}^{n/2-1}\prod_{p=1}^m
\left\{e_k (v_k-v_{kp})+\tilde e_k(\tilde v_k-\tilde v_{kp})\right\},
\label{npolar-127a}
\end{equation}
where the quantities $v_{p+}$ appear always in complex conjugate
pairs, and the quantities $\tilde v_{p-}$ appear always in complex
conjugate pairs.
For odd $n$ the polynomial can be written as
\begin{equation}
P_m=e_+\prod_{p=1}^m (v_+ -v_{p+})
+\sum_{k=1}^{(n-1)/2}\prod_{p=1}^m
\left\{e_k (v_k-v_{kp})+\tilde e_k(\tilde v_k-\tilde v_{kp})\right\},
\label{npolar-127b}
\end{equation}
where the quantities $v_{p+}$ appear always in complex conjugate
pairs.
Due to the relations  (\ref{npolar-e12a}),(\ref{npolar-e12b}),
the polynomial $P_m(u)$ can be written, for even $n$, as a product of factors
of the form  
\begin{eqnarray}
\lefteqn{P_m(u)=\prod_{p=1}^m \left\{e_+(v_+ -v_{p+})+e_-(v_- -v_{p-})
+\sum_{k=1}^{n/2-1}\left\{e_k (v_k-v_{kp})+\tilde e_k(\tilde v_k-
\tilde v_{kp})\right\}\right\}.\nonumber}\\
&&
\label{npolar-128a}
\end{eqnarray}\index{polynomial, factorization!polar n-complex}
For odd $n$, the polynomial $P_m(u)$ can be written as the product
\begin{eqnarray}
P_m(u)=\prod_{p=1}^m \left\{e_+(v_+ -v_{p+})
+\sum_{k=1}^{(n-1)/2}\left\{e_k (v_k-v_{kp})
+\tilde e_k(\tilde v_k-\tilde v_{kp})\right\}\right\}.
\label{npolar-128b}
\end{eqnarray}
These relations can be written with the aid of Eqs. (\ref{npolar-e13a}) and
(\ref{npolar-e13b}) as
\begin{eqnarray}
P_m(u)=\prod_{p=1}^m (u-u_p) ,
\label{npolar-128c}
\end{eqnarray}
where, for even $n$,
\begin{eqnarray}
u_p=e_+ v_{p+}+e_-v_{p-}
+\sum_{k=1}^{n/2-1}\left(e_k v_{kp}+\tilde e_k\tilde v_{kp}\right), 
\label{npolar-128d}
\end{eqnarray}
and for odd $n$
\begin{eqnarray}
u_p=e_+ v_{p+}
+\sum_{k=1}^{(n-1)/2}\left(e_k v_{kp}+\tilde e_k\tilde v_{kp}\right), 
\label{npolar-128e}
\end{eqnarray}
for $p=1,...,m$.
The roots $v_{p+}$, the roots $v_{p-}$ and, for a given $k$, the roots 
$e_k v_{k1}+\tilde e_k\tilde v_{k1}, ...,  e_k v_{km}+\tilde e_k\tilde v_{km}$
defined in Eqs. (\ref{npolar-127a}) or (\ref{npolar-127b}) may be ordered arbitrarily.
This means that Eqs. (\ref{npolar-128d}) or (\ref{npolar-128e}) give sets of $m$ roots
$u_1,...,u_m$ of the polynomial $P_m(u)$, 
corresponding to the various ways in which the roots $v_{p+}, v_{p-}$,
$e_k v_{kp}+\tilde e_k\tilde v_{kp}$ are ordered according to $p$ in each
group. Thus, while the n-complex components in Eq. (\ref{npolar-126b}) taken
separately have unique factorizations, the polynomial $P_m(u)$ can be written
in many different ways as a product of linear factors. 

If $P(u)=u^2-1$, the degree is $m=2$, the coefficients of the polynomial are
$a_1=0, a_2=-1$, the n-complex components of $a_2$ are $a_{20}=-1, a_{21}=0,
... ,a_{2,n-1}=0$, the components $A_{2+}, A_{2-}, A_{2k}, \tilde A_{2k}$
calculated according 
to Eqs. (\ref{npolar-88a})-(\ref{npolar-88d}) are $A_{2+}=-1, A_{2-}=-1, A_{2k}=-1, \tilde
A_{2k}=0, k=1,...,[(n-1)/2]$. The expression of $P(u)$ for even $n$, Eq.
(\ref{npolar-126a}), is  
$e_+(v_+^2-1)+e_-(v_-^2-1)+\sum_{k=1}^{n/2-1}\{(e_k v_k+\tilde e_k\tilde
v_k)^2-e_k$\}, and Eq. (\ref{npolar-127a}) has the form 
$u^2-1=e_+(v_++1)(v_+-1)+e_-(v_-+1)(v_--1)+
\sum_{k=1}^{n/2-1}\left\{e_k (v_k+1)+\tilde e_k\tilde v_k\right\}
\left\{e_k (v_k-1)+\tilde e_k\tilde v_k\right\}$.
For odd $n$, the expression of $P(u)$, Eq. (\ref{npolar-126b}), is  
$e_+(v_+^2-1)+\sum_{k=1}^{(n-1)/2}\{(e_k v_k+\tilde e_k\tilde v_k)^2-e_k\}$,
and Eq. (\ref{npolar-127b}) has the form 
$u^2-1=e_+(v_++1)(v_+-1)+
\sum_{k=1}^{(n-1)/2}\left\{e_k (v_k+1)+\tilde e_k\tilde v_k\right\}
\left\{e_k (v_k-1)+\tilde e_k\tilde v_k\right\}$.
The factorization in Eq. (\ref{npolar-128c}) is $u^2-1=(u-u_1)(u-u_2)$, where for even
$n$, $u_1=\pm e_+\pm e_-\pm e_1\pm e_2\pm\cdots
\pm e_{n/2-1}, u_2=-u_1$, so that
there are $2^{n/2}$ independent sets of roots $u_1,u_2$
of $u^2-1$. 
It can be checked that 
$(\pm e_+\pm e_-\pm e_1\pm e_2\pm\cdots\pm e_{n/2-1})^2= 
e_++e_-+e_1+e_2+\cdots+e_{n/2-1}=1$.
For odd $n$, $u_1=\pm e_+\pm e_1\pm e_2\pm\cdots
\pm e_{(n-1)/2}, u_2=-u_1$, so that there are $2^{(n-1)/2}$ independent
sets of roots $u_1,u_2$ of $u^2-1$.
It can be checked that 
$(\pm e_+\pm e_1\pm e_2\pm\cdots\pm e_{(n-1)/2})^2= 
e_++e_1+e_2+\cdots+e_{(n-1)/2}=1$.

\subsection{Representation of polar n-complex numbers by irreducible matrices}

If the unitary matrix written in Eq. (\ref{npolar-11}), for even $n$, is called $T_e$,
and the unitary matrix written in Eq. (\ref{npolar-12}), for odd $n$, is called $T_o$,
it can be shown that, for even $n$, the matrix $T_e U T_e^{-1}$ has the form 
\begin{equation}
T_e U T_e^{-1}=\left(
\begin{array}{ccccc}
v_+     &     0     &     0   & \cdots  &   0   \\
0       &     v_-   &     0   & \cdots  &   0   \\
0       &     0     &     V_1 & \cdots  &   0   \\
\vdots  &  \vdots   &  \vdots & \cdots  &\vdots \\
0       &     0     &     0   & \cdots  &   V_{n/2-1}\\
\end{array}
\right)
\label{npolar-129a}
\end{equation}\index{representation by irreducible matrices!polar n-complex}
and, for odd $n$, the matrix $T_o U T_o^{-1}$ has the form 
\begin{equation}
T_o U T_o^{-1}=\left(
\begin{array}{ccccc}
v_+     &     0     &     0   & \cdots  &   0   \\
0       &     V_1   &     0   & \cdots  &   0   \\
0       &     0     &     V_2 & \cdots  &   0   \\
\vdots  &  \vdots   &  \vdots & \cdots  &\vdots \\
0       &     0     &     0   & \cdots  &   V_{(n-1)/2}\\
\end{array}
\right),
\label{npolar-129b}
\end{equation}
where $U$ is the matrix in Eq. (\ref{npolar-24b}) used to represent the n-complex
number $u$. In Eqs. (\ref{npolar-129a}) and (\ref{npolar-129b}), $V_k$ are
the matrices
\begin{equation}
V_k=\left(
\begin{array}{cc}
v_k           &     \tilde v_k   \\
-\tilde v_k   &     v_k          \\
\end{array}\right),
\label{npolar-130}
\end{equation}
for $ k=1,...,[(n-1)/2]$, where $v_k, \tilde v_k$ are the variables introduced in Eqs. (\ref{npolar-9a}) and
(\ref{npolar-9b}), and the symbols 0 denote, according to the case, the real number
zero, or one of the matrices
\begin{equation}
\left(
\begin{array}{c}
0  \\
0  \\
\end{array}\right)
\;\; {\rm or}\;\;
\left(
\begin{array}{cc}
0   &  0   \\
0   &  0   \\
\end{array}\right).
\label{npolar-131}
\end{equation}
The relations between the variables $v_k, \tilde v_k$ for the multiplication of
n-complex numbers have been written in Eq. (\ref{npolar-22}). The matrices 
$T_e U T_e^{-1}$ and $T_o U T_o^{-1}$ provide an irreducible representation
\cite{4} of the n-complex numbers $u$ in terms of matrices with real
coefficients. 

\section{Planar complex numbers in even $n$ dimensions}

\subsection{Operations with planar n-complex numbers}

A hypercomplex number in $n$ dimensions is determined by its $n$ components
$(x_0,x_1,...,x_{n-1})$. The planar n-complex numbers and
their operations discussed in this section can be represented 
by  writing the n-complex number $(x_0,x_1,...,x_{n-1})$ as  
$u=x_0+h_1x_1+h_2x_2+\cdots+h_{n-1}x_{n-1}$, where 
$h_1, h_2, \cdots, h_{n-1}$ are bases for which the multiplication rules are 
\begin{equation}
h_j h_k =(-1)^{[(j+k)/n]}h_l ,\:l=j+k-n[(i+k)/n],
\label{nplanar-1}
\end{equation}\index{complex units!planar n-complex}
for $ j,k,l=0,1,..., n-1$,
where $h_0=1$.
In Eq. (\ref{nplanar-1}), $[(j+k)/n]$ denotes the integer part of $(j+k)/n$, 
the integer part being defined as $[a]\leq a<[a]+1$, so that
$0\leq j+k-n[(j+k)/n]\leq n-1$. 
As already mentioned, brackets larger than the regular brackets
$[\;]$ do not have the meaning of integer part.
The significance of the composition laws in Eq.
(\ref{nplanar-1}) can be understood by representing the bases $h_j, h_k$ by points on a
circle at the angles $\alpha_j=\pi j/n,\alpha_k=\pi k/n$, as shown in Fig. \ref{fig24},
and the product $h_j h_k$ by the point of the circle at the angle 
$\pi (j+k)/n$. If $\pi\leq \pi (j+k)/n<2\pi$, the point is opposite to the
basis $h_l$ of angle $\alpha_l=\pi (j+k)/n-\pi$.

\begin{figure}
\begin{center}
\epsfig{file=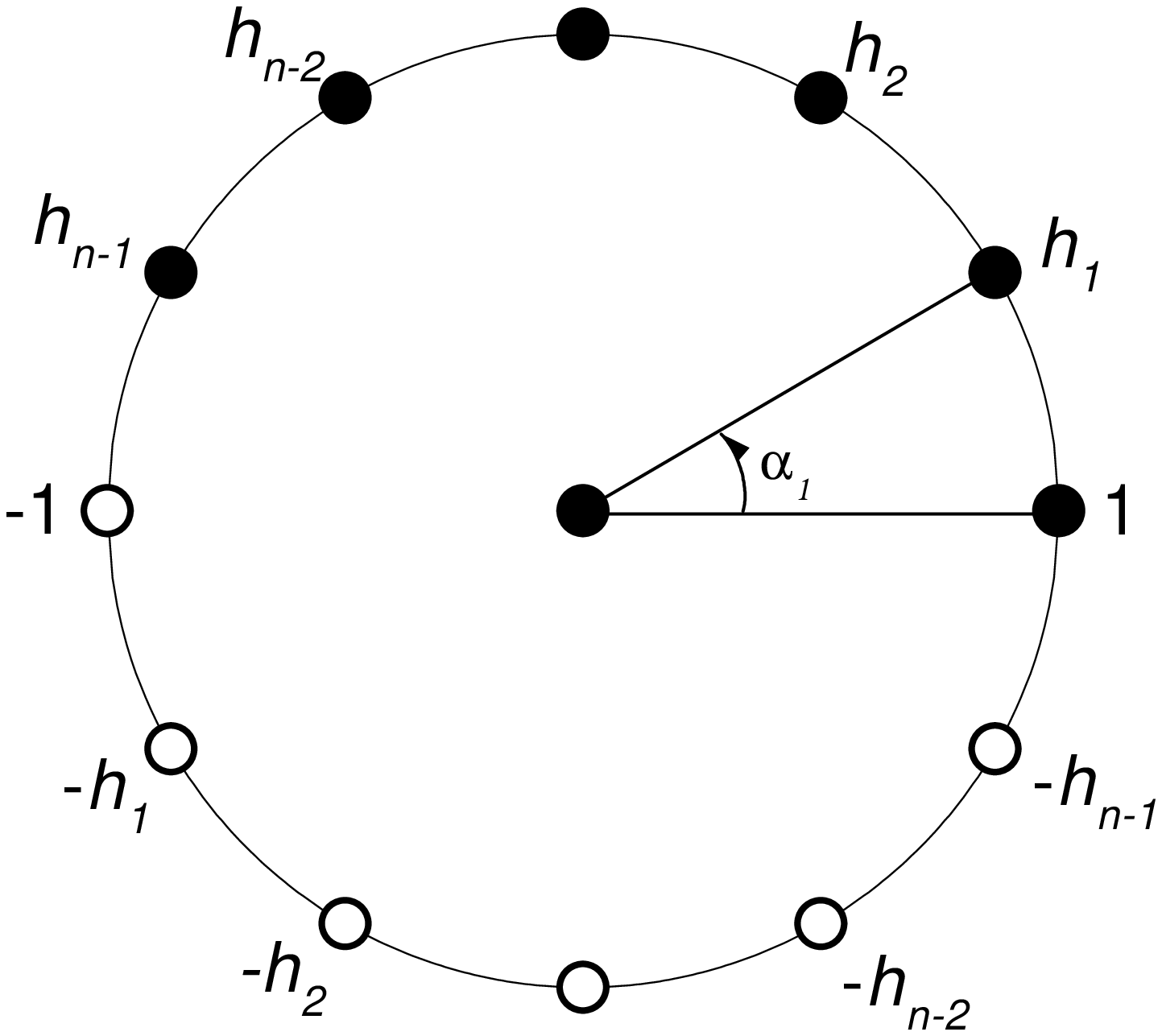,width=12cm}
\caption{Representation of the hypercomplex bases $1, h_1,...,h_{n-1}$
by points on a circle at the angles $\alpha_k=\pi k/n$.
The product $h_j h_k$ will be represented by the point of the circle at the
angle $\pi (j+k)/2n$, $j,k=0,1,...,n-1$. If $\pi\leq\pi (j+k)/2n\leq 2\pi$, the
point is opposite to the basis $h_l$ of angle $\alpha_l=\pi (j+k)/n-\pi$. }
\label{fig24}
\end{center}
\end{figure}

In an odd number of dimensions $n$, a transformation of coordinates according to
\begin{equation}
x_{2l}=x^\prime_l, x_{2m-1}=-x^\prime_{(n-1)/2+m}, 
\label{nplanar-1a}
\end{equation}
and of the bases according to
\begin{equation}
h_{2l}=h^\prime_l, h_{2m-1}=-h^\prime_{(n-1)/2+m}, 
\label{nplanar-1b}
\end{equation}
where $l=0,...,(n-1)/2, \; m=1,...,(n-1)/2$,
leaves the expression of an n-complex number unchanged,
\begin{equation}
\sum_{k=0}^{n-1}h_k x_k=\sum_{k=0}^{n-1}h^\prime_k x^\prime_k,
\label{nplanar-1c}
\end{equation}
and the products of the bases $h^\prime_k$ are
\begin{equation}
h^\prime_j h^\prime_k =h^\prime_l ,\:l=j+k-n[(j+k)/n],
\label{nplanar-1d}
\end{equation}
for $j,k,l=0,1,..., n-1$. 
Thus, the n-complex numbers with the rules (\ref{nplanar-1}) are equivalent in an
odd number of dimensions to the polar n-complex numbers described in the
previous chapter.
Therefore, in this section it will be supposed that $n$ is an even number, unless
otherwise stated.

Two n-complex numbers 
$u=x_0+h_1x_1+h_2x_2+\cdots+h_{n-1}x_{n-1}$,
$u^\prime=x^\prime_0+h_1x^\prime_1+h_2x^\prime_2+\cdots+h_{n-1}x^\prime_{n-1}$ 
are equal if and only if $x_j=x^\prime_j, j=0,1,...,n-1$.
The sum of the n-complex numbers $u$
and
$u^\prime$ 
is
\begin{equation}
u+u^\prime=x_0+x^\prime_0+h_1(x_1+x^\prime_1)+\cdots
+h_{n-1}(x_{n-1} +x^\prime_{n-1}) .
\label{nplanar-2}
\end{equation}\index{sum!planar n-complex}
The product of the numbers $u, u^\prime$ is
\begin{equation}
\begin{array}{l}
uu^\prime=x_0 x_0^\prime -x_1x_{n-1}^\prime-x_2 x_{n-2}^\prime-x_3x_{n-3}^\prime
-\cdots-x_{n-1}x_1^\prime\\
+h_1(x_0 x_1^\prime+x_1x_0^\prime-x_2x_{n-1}^\prime-x_3x_{n-2}^\prime
-\cdots-x_{n-1} x_2^\prime) \\
+h_2(x_0 x_2^\prime+x_1x_1^\prime+x_2x_0^\prime-x_3x_{n-1}^\prime
-\cdots-x_{n-1} x_3^\prime) \\
\vdots\\
+h_{n-1}(x_0 x_{n-1}^\prime+x_1x_{n-2}^\prime+x_2x_{n-3}^\prime+x_3x_{n-4}^\prime
+\cdots+x_{n-1} x_0^\prime).
\end{array}
\label{nplanar-3}
\end{equation}\index{product!planar n-complex}
The product $uu^\prime$ can be written as
\begin{equation}
uu^\prime=\sum_{k=0}^{n-1}h_k\sum_{l=0}^{n-1}(-1)^{[(n-k-1+l)/n]}x_l x^\prime_{k-l+n[(n-k-1+l)/n]}.
\label{nplanar-3a}
\end{equation}
If $u,u^\prime,u^{\prime\prime}$ are n-complex numbers, the multiplication 
is associative
\begin{equation}
(uu^\prime)u^{\prime\prime}=u(u^\prime u^{\prime\prime})
\label{nplanar-3b}
\end{equation}
and commutative
\begin{equation}
u u^\prime=u^\prime u ,
\label{nplanar-3c}
\end{equation}
because the product of the bases, defined in Eq. (\ref{nplanar-1}), is associative and
commutative. The fact that the multiplication is commutative can be seen also
directly from Eq. (\ref{nplanar-3}). 
The n-complex zero is $0+h_1\cdot 0+\cdots+h_{n-1}\cdot 0,$ 
denoted simply 0, 
and the n-complex unity is $1+h_1 \cdot 0+\cdots+h_{n-1}\cdot 0,$ 
denoted simply 1.

The inverse of the n-complex number $u=x_0+h_1x_1+h_2x_2+\cdots+h_{n-1}x_{n-1}$
is the n-complex number
$u^\prime
=x^\prime_0+h_1x^\prime_1+h_2x^\prime_2+\cdots+h_{n-1}x^\prime_{n-1}$ 
having the property that
\begin{equation}
uu^\prime=1 .
\label{nplanar-4}
\end{equation}\index{inverse!planar n-complex}
Written on components, the condition, Eq. (\ref{nplanar-4}), is
\begin{equation}
\begin{array}{l}
x_0 x_0^\prime -x_1x_{n-1}^\prime-x_2 x_{n-2}^\prime-x_3x_{n-3}^\prime
-\cdots-x_{n-1}x_1^\prime=1,\\
x_0 x_1^\prime+x_1x_0^\prime-x_2x_{n-1}^\prime-x_3x_{n-2}^\prime
-\cdots-x_{n-1} x_2^\prime=0, \\
x_0 x_2^\prime+x_1x_1^\prime+x_2x_0^\prime-x_3x_{n-1}^\prime
-\cdots-x_{n-1} x_3^\prime=0, \\
\vdots\\
x_0 x_{n-1}^\prime+x_1x_{n-2}^\prime+x_2x_{n-3}^\prime+x_3x_{n-4}^\prime
+\cdots+x_{n-1} x_0^\prime=0.
\end{array}
\label{nplanar-5}
\end{equation}
The system (\ref{nplanar-5}) has a solution provided that the determinant of the
system, 
\begin{equation}
\nu={\rm det}(A), 
\label{nplanar-5b}
\end{equation}\index{inverse, determinant!planar n-complex}
is not equal to zero, $\nu\not=0$, where
\begin{equation}
A=\left(
\begin{array}{ccccc}
x_0     &  -x_{n-1} &   -x_{n-2}  & \cdots  &-x_1\\
x_1     &   x_0     &   -x_{n-1}  & \cdots  &-x_2\\
x_2     &   x_1     &   x_0       & \cdots  &-x_3\\
\vdots  &  \vdots   &  \vdots     & \cdots  &\vdots \\
x_{n-1} &   x_{n-2} &   x_{n-3}   & \cdots  &x_0\\
\end{array}
\right).
\label{nplanar-6}
\end{equation}
It will be shown that $\nu>0$, and the quantity 
\begin{equation}
\rho=\nu^{1/n}
\label{nplanar-6b}
\end{equation}\index{amplitude!planar n-complex}
will be called amplitude of the
n-complex number $u=x_0+h_1x_1+h_2x_2+\cdots+h_{n-1}x_{n-1}$.
The quantity $\nu$ can be written as a product of linear factors
\begin{equation}
\nu=\prod_{k=1}^{n}
\left(x_0+\epsilon_k x_1+\epsilon_k^2 x_2+\cdots+\epsilon^{n-1}_k x_{n-1}\right),
\label{nplanar-7}
\end{equation}\index{inverse, determinant!planar n-complex}
where $\epsilon_k=e^{i\pi (2k-1)/n}$, $k=1,...,n$, and $i$ being the imaginary
unit. The factors appearing in Eq. (\ref{nplanar-7}) are of the form 
\begin{equation}
x_0+\epsilon_k x_1+\epsilon_k^2 x_2+\cdots+\epsilon^{n-1}_k x_{n-1}=v_k+i\tilde v_k,
\label{nplanar-8}
\end{equation}
where
\begin{equation}
v_k=\sum_{p=0}^{n-1}x_p\cos\frac{\pi (2k-1)p}{n},
\label{nplanar-9a}
\end{equation}
\begin{equation}
\tilde v_k=\sum_{p=0}^{n-1}x_p\sin\frac{\pi (2k-1)p}{n},
\label{nplanar-9b}
\end{equation}
for $k=1,...,n$.\index{canonical variables!planar n-complex}
The variables $v_k, \tilde v_k, k=1,...,n/2$ will be called canonical polar
n-complex variables. 
It can be seen that $v_k=v_{n-k+1}, \tilde v_k=-\tilde v_{n-k+1}$, for
$k=1,...,n/2$.
Therefore,  the factors appear in Eq. (\ref{nplanar-7}) in complex-conjugate pairs of
the form $v_k+i\tilde v_k$ and $v_{n-k+1}+i\tilde v_{n-k+1}=v_k-i\tilde v_k$,
where $k=1,...n/2$, so that the determinant $\nu$ is a real and positive
quantity, $\nu>0$, 
\begin{equation}
\nu=\prod_{k=1}^{n/2}\rho_k^2,
\label{nplanar-9c}
\end{equation}\index{inverse, determinant!planar n-complex}
where
\begin{equation}
\rho_k^2=v_k^2+\tilde v_k^2 .
\label{nplanar-9d}
\end{equation}
Thus, an n-complex number has an inverse
unless it lies on one of the nodal hypersurfaces 
$\rho_1=0$, or $\rho_2=0$, or ... or $\rho_{n/2}=0$.
\index{nodal hypersurfaces!planar n-complex}

\newpage
\setlength{\oddsidemargin}{-1cm}
\subsection{Geometric representation of planar n-complex numbers}

The n-complex number $x_0+h_1x_1+h_2x_2+\cdots+h_{n-1}x_{n-1}$
can be represented by 
the point $A$ of coordinates $(x_0,x_1,...,x_{n-1})$. 
If $O$ is the origin of the n-dimensional space,  the
distance from the origin $O$ to the point $A$ of coordinates
$(x_0,x_1,...,x_{n-1})$ has the expression
\begin{equation}
d^2=x_0^2+x_1^2+\cdots+x_{n-1}^2 .
\label{nplanar-10}
\end{equation}\index{distance!planar n-complex}
The quantity $d$ will be called modulus of the n-complex number 
$u=x_0+h_1x_1+h_2x_2+\cdots+h_{n-1}x_{n-1}$. 
\index{modulus, definition!planar n-complex}
The modulus of an n-complex number
$u$ will be designated by $d=|u|$.

The exponential and trigonometric forms of the n-complex number $u$ can be
obtained conveniently in a rotated system of axes defined by a transformation
which has the form
\begin{equation}
\left(
\begin{array}{c}
\vdots\\
\xi_k\\
\eta_k\\
\vdots
\end{array}\right)
=\left(
\begin{array}{ccccc}
\vdots&\vdots& &\vdots&\vdots\\
\sqrt{\frac{2}{n}}&\sqrt{\frac{2}{n}}\cos\frac{\pi (2k-1)}{n}&\cdots&\sqrt{\frac{2}{n}}\cos\frac{\pi (2k-1)(n-2)}{n}&\sqrt{\frac{2}{n}}\cos\frac{\pi (2k-1)(n-1)}{n}\\
0&\sqrt{\frac{2}{n}}\sin\frac{\pi (2k-1)}{n}&\cdots&\sqrt{\frac{2}{n}}\sin\frac{\pi (2k-1)(n-2)}{n}&\sqrt{\frac{2}{n}}\sin\frac{\pi (2k-1)(n-1)}{n}\\
\vdots&\vdots&&\vdots&\vdots
\end{array}
\right)
\left(
\begin{array}{c}
x_0\\
\vdots\\ 
\vdots\\
x_{n-1}
\end{array}
\right),
\label{nplanar-11}
\end{equation}
where $k=1, 2, ... , n/2$.
The lines of the matrices in Eq. (\ref{nplanar-11}) give the components
of the $n$ vectors of the new basis system of axes. These vectors have unit
length and are orthogonal to each other.
By comparing Eqs. (\ref{nplanar-9a})-(\ref{nplanar-9b}) and (\ref{nplanar-11}) it can be
seen that
\begin{equation}
v_k= \sqrt{\frac{n}{2}}\xi_k , \tilde v_k= \sqrt{\frac{n}{2}}\eta_k ,
\label{nplanar-12b}
\end{equation}
i.e. the two sets of variables differ only by a scale factor.

The sum of the squares of the variables $v_k,\tilde v_k$ is
\begin{equation}
\sum_{k=1}^{n/2}(v_k^2+\tilde v_k^2)=\frac{n}{2}d^2.
\label{nplanar-13}
\end{equation}
The relation (\ref{nplanar-13}) has been obtained with the aid of the relation
\begin{equation}
\sum_{k=1}^{n/2}\cos\frac{\pi (2k-1)p}{n}=0, 
\label{nplanar-15}
\end{equation}
for $p=1,...,n-1$.
From Eq. (\ref{nplanar-13}) it results that
\begin{equation}
d^2=\frac{2}{n}\sum_{k=1}^{n/2}\rho_k^2 .
\label{nplanar-17}
\end{equation}\index{modulus, canonical variables!planar n-complex}
The relation (\ref{nplanar-17}) shows that the square of the distance
$d$, Eq. (\ref{nplanar-10}), is equal to the sum of the squares of the projections
$\rho_k\sqrt{2/n}$. This is consistent with the fact that the transformation in
Eq. (\ref{nplanar-11}) is unitary.

\newpage
\setlength{\oddsidemargin}{0.9cm}

The position of the point $A$ of coordinates $(x_0,x_1,...,x_{n-1})$ can be
also described with the aid of the distance $d$, Eq. (\ref{nplanar-10}), and of $n-1$
angles defined further. Thus, in the plane of the axes $v_k,\tilde v_k$, the
azimuthal angle $\phi_k$ can be introduced by the relations 
\begin{equation}
\cos\phi_k=v_k/\rho_k,\:\sin\phi_k=\tilde v_k/\rho_k, 
\label{nplanar-19a}
\end{equation}\index{azimuthal angles!planar n-complex}
where $0\leq \phi_k<2\pi, \;k=1,...,n/2$, so that there are $n/2$ azimuthal angles.
The radial distance $\rho_k$ in the plane of the axes $v_k,\tilde v_k$ has been
defined in Eq. (\ref{nplanar-9d}).
If the projection of the point $A$ on the plane of the axes $v_k,\tilde v_k$ is
$A_k$, and the projection of the point $A$ on the 4-dimensional space defined
by the axes $v_1, \tilde v_1, v_k,\tilde v_k$ is $A_{1k}$, the angle
$\psi_{k-1}$ between the line $OA_{1k}$ and the 2-dimensional plane defined by
the axes $v_k,\tilde v_k$ is  
\begin{equation}
\tan\psi_{k-1}=\rho_1/\rho_k, 
\label{nplanar-19b}
\end{equation}\index{planar angles!planar n-complex}
where $0\leq\psi_k\leq\pi/2, k=2,...,n/2$,
so that there are $n/2-1$ planar angles. Thus, the
position of the point $A$ is described by the distance $d$, by $n/2$
azimuthal angles and by $n/2-1$ 
planar angles. These angles are shown in Fig. \ref{fig25}.

\begin{figure}
\begin{center}
\epsfig{file=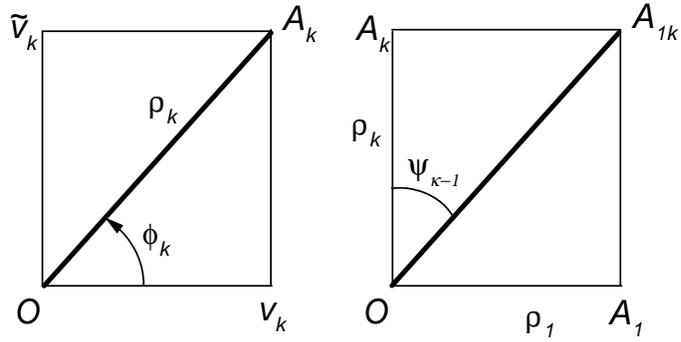,width=12cm}
\caption{Radial distance $\rho_k$  and
azimuthal angle $\phi_k$ in the plane of the axes $v_k,\tilde v_k$, and planar
angle $\psi_{k-1}$ between the line $OA_{1k}$ and the 2-dimensional plane
defined by the axes $v_k,\tilde v_k$. $A_k$ is the projection of the point $A$
on the plane of the axes $v_k,\tilde v_k$, and $A_{1k}$ is the projection of
the point $A$ on the 4-dimensional space defined 
by the axes $v_1, \tilde v_1, v_k,\tilde v_k$. }
\label{fig25}
\end{center}
\end{figure}

The variables $\rho_k$ can be expressed in terms of $d$ and the planar angles
$\psi_k$ as
\begin{equation}
\rho_k=\frac{\rho_1}{\tan\psi_{k-1}}, 
\label{nplanar-20a}
\end{equation}
for $k=2,...,n/2$, where
\begin{eqnarray}
\rho_1^2=\frac{nd^2}{2}
\left(1
+\frac{1}{\tan^2\psi_1}+\frac{1}{\tan^2\psi_2}+\cdots
+\frac{1}{\tan^2\psi_{n/2-1}}\right)^{-1}.
\label{nplanar-20b}
\end{eqnarray}

If
$u^\prime=x_0^\prime+h_1x_1^\prime+h_2x_2^\prime+\cdots+h_{n-1}x_{n-1}^\prime, 
u^{\prime\prime}=x^{\prime\prime}_0+h_1x^{\prime\prime}_1
+h_2x^{\prime\prime}_2+\cdots+h_{n-1}x^{\prime\prime}_{n-1}$ 
are n-complex numbers of parameters
$\rho_k^\prime,\psi_k^\prime,\phi_k^\prime$ and respectively  
$\rho_k^{\prime\prime}, \psi_k^{\prime\prime},
\phi_k^{\prime\prime}$, then the parameters
$v_+,\rho_k,\psi_k,\phi_k$ of the product n-complex number
$u=u^\prime u^{\prime\prime}$ are given by
\begin{equation}
\rho_k=\rho_k^\prime\rho_k^{\prime\prime}, 
\label{nplanar-21b}
\end{equation}\index{transformation of variables!planar n-complex}
for $k=1,..., n/2$,
\begin{equation}
\tan\psi_k=\tan\psi_k^\prime \tan\psi_k^{\prime\prime}, 
\label{nplanar-21d}
\end{equation}
for $k=1,...,n/2-1$, 
\begin{equation}
\phi_k=\phi_k^\prime+\phi_k^{\prime\prime}, 
\label{nplanar-21e}
\end{equation}
for $k=1,..., n/2$.
The Eqs. (\ref{nplanar-21b})-(\ref{nplanar-21e}) are a consequence of the relations
\begin{equation}
v_k=v_k^\prime v_k^{\prime\prime}-\tilde v_k^\prime \tilde v_k^{\prime\prime},\;
\tilde v_k=v_k^\prime \tilde v_k^{\prime\prime}+\tilde v_k^\prime v_k^{\prime\prime},
\label{nplanar-22}
\end{equation}
and of the corresponding relations of definition. Then the product $\nu$ in
Eq. (\ref{nplanar-9c}) has the property that
\begin{equation}
\nu=\nu^\prime\nu^{\prime\prime} ,
\label{nplanar-23}
\end{equation}
and the amplitude $\rho$ defined in Eq. (\ref{nplanar-6b}) has the property that
\begin{equation}
\rho=\rho^\prime\rho^{\prime\prime} .
\label{nplanar-24}
\end{equation}

The fact that the amplitude of the product is equal to the product of the 
amplitudes, as written in Eq. (\ref{nplanar-24}), can 
be demonstrated also by using a representation of the 
n-complex numbers by matrices, in which the n-complex number 
$u=x_0+h_1x_1+h_2x_2+\cdots+h_{n-1}x_{n-1}$ is represented by the matrix
\begin{equation}
U=\left(
\begin{array}{ccccc}
x_0      &   x_1      &   x_2    & \cdots  &x_{n-1}\\
-x_{n-1} &   x_0      &   x_1    & \cdots  &x_{n-2}\\
-x_{n-2} &  - x_{n-1} &   x_0    & \cdots  &x_{n-3}\\
\vdots   &  \vdots    &  \vdots  & \cdots  &\vdots \\
-x_1     &   -x_2     &  - x_3   & \cdots  &x_0\\
\end{array}
\right).
\label{nplanar-24a}
\end{equation}\index{matrix representation!planar n-complex}
The product $u=u^\prime u^{\prime\prime}$ is be
represented by the matrix multiplication $U=U^\prime U^{\prime\prime}$.
The relation (\ref{nplanar-23}) is then a consequence of the fact the determinant 
of the product of matrices is equal to the product of the determinants 
of the factor matrices. 
The use of the representation  of the n-complex numbers with matrices
provides an alternative demonstration of the fact that the product of
n-complex numbers is associative, as stated in 
Eq. (\ref{nplanar-3b}).

According to Eqs. (\ref{nplanar-13} and (\ref{nplanar-9d}), the modulus of the product
$uu^\prime$ is given by  
\begin{equation}
|uu^\prime|^2=
\frac{2}{n}\sum_{k=1}^{n/2}(\rho_k\rho_k^\prime)^2 .
\label{nplanar-25a}
\end{equation}
Thus, if the product of two n-complex numbers is zero, $uu^\prime=0$, then
$\rho_k\rho_k^\prime=0, k=1,...,n/2$. This means that either $u=0$, or
$u^\prime=0$, or $u, u^\prime$ belong to orthogonal hypersurfaces in such a way
that the afore-mentioned products of components should be equal to zero.
\index{divisors of zero!planar n-complex}

\subsection{The planar n-dimensional cosexponential functions}

The exponential function of a hypercomplex variable $u$ and the addition
theorem for the exponential function have been written in Eqs. 
~(1.35)-(1.36).
It can be seen with the aid of the representation in Fig. \ref{fig24} that 
\begin{equation}
h_k^{n+p}=(-1)^k h_k^p, \:p\:\:{\rm integer},
\label{nplanar-28}
\end{equation}\index{complex units, powers of!planar n-complex}
where $k=1,...,n-1$. For $k$ even, $e^{h_k y}$ can be written as
\begin{equation}
e^{h_k y}=\sum_{p=0}^{n-1}(-1)^{[kp/n]}h_{kp-n[kp/n]}g_{np}(y),
\label{nplanar-28bx}
\end{equation}
where $h_0=1$, and where $g_{np}$ are the polar n-dimensional cosexponential
functions.  For odd $k$,  $e^{h_k y}$ is
\begin{equation}
e^{h_k y}=\sum_{p=0}^{n-1}(-1)^{[kp/n]}h_{kp-n[kp/n]}f_{np}(y),
\label{nplanar-28b}
\end{equation}\index{exponential, expressions!planar n-complex}
where the functions $f_{nk}$, which will be
called planar cosexponential functions in $n$ dimensions, are
\begin{equation}
f_{nk}(y)=\sum_{p=0}^\infty (-1)^p \frac{y^{k+pn}}{(k+pn)!}, 
\label{nplanar-29}
\end{equation}\index{cosexponential functions, planar n-complex!definition}
for $ k=0,1,...,n-1$.

The planar cosexponential functions  of even index $k$ are
even functions, $f_{n,2l}(-y)=f_{n,2l}(y)$,  
and the planar cosexponential functions of odd index 
are odd functions, $f_{n,2l+1}(-y)=-f_{n,2l+1}(y)$, $l=0,...,n/2-1$ . 
\index{cosexponential functions, planar n-complex!parity}

The planar n-dimensional cosexponential function $f_{nk}(y)$ is related to the
polar n-dimensional cosexponential function $g_{nk}(y)$ discussed in the
previous chapter by
the relation 
\begin{equation}
f_{nk}(y)=e^{-i\pi k/n}g_{nk}\left(e^{i\pi/n}y\right), 
\label{nplanar-30a}
\end{equation}
for $k=0,...,n-1$.
The expression of the planar n-dimensional cosexponential functions is then
\begin{equation}
f_{nk}(y)=\frac{1}{n}\sum_{l=1}^{n}\exp\left[y\cos\left(\frac{\pi (2l-1)}{n}\right)
\right]
\cos\left[y\sin\left(\frac{\pi (2l-1)}{n}\right)-\frac{\pi (2l-1)k}{n}\right], 
\label{nplanar-30}
\end{equation}\index{cosexponential functions, planar n-complex!expressions}
for $k=0,1,...,n-1$.
The planar cosexponential function defined in Eq. (\ref{nplanar-29}) has the expression
given in Eq. (\ref{nplanar-30}) for any natural value of $n$, this result not being
restricted to even values of $n$.
In order to check that the function in Eq. (\ref{nplanar-30}) has the series expansion
written in Eq. (\ref{nplanar-29}), the right-hand side of Eq. (\ref{nplanar-30}) will be
written as
\begin{equation}
f_{nk}(y)=\frac{1}{n}\sum_{l=1}^{n}{\rm Re}\left\{
\exp\left[\left(\cos\frac{\pi (2l-1)}{n}+i\sin\frac{\pi (2l-1)}{n}\right)y
-i\frac{\pi k(2l-1)}{n}\right]\right\}, 
\label{nplanar-31}
\end{equation}
for $k=0,1,...,n-1$, where ${\rm Re} (a+ib)=a$, with $a$ and $b$ real numbers. The part of the
exponential depending on $y$ can be expanded in a series,
\begin{equation}
f_{nk}(y)=\frac{1}{n}\sum_{p=0}^\infty\sum_{l=1}^{n}{\rm Re}
\left\{\frac{1}{p!}
\exp\left[i\frac{\pi (2l-1)}{n}(p-k)\right]y^p\right\}, 
\label{nplanar-32}
\end{equation}
for $k=0,1,...,n-1$.
The expression of $f_{nk}(y)$ becomes
\begin{equation}
f_{nk}(y)=\frac{1}{n}\sum_{p=0}^\infty\sum_{l=1}^{n}
\left\{\frac{1}{p!}
\cos\left[\frac{\pi (2l-1)}{n}(p-k)\right]y^p\right\}, 
\label{nplanar-33}
\end{equation}
where $k=0,1,...,n-1$ and, since
\begin{equation}
\frac{1}{n}\sum_{l=1}^{n}\cos\left[\frac{\pi (2l-1)}{n}(p-k)\right]
=\left\{
\begin{array}{l}
1, \:\:{\rm if}\:\: p-k\:\: {\rm is \:\;an\:\;even\:\;multiple\:\;of\:\;}n,\\
-1, \:\:{\rm if}\:\: p-k\:\: {\rm is \:\;an\:\;odd\:\;multiple\:\;of\:\;}n,\\
0, \:\:{\rm otherwise},
\end{array}
\right.
\label{nplanar-34}
\end{equation}
this yields indeed the expansion in Eq. (\ref{nplanar-29}).

It can be shown from Eq. (\ref{nplanar-30}) that
\begin{equation}
\sum_{k=0}^{n-1}f_{nk}^2(y)=\frac{1}{n}\sum_{l=1}^{n}\exp\left[2y\cos\left(
\frac{\pi (2l-1)}{n}\right)\right].
\label{nplanar-34a}
\end{equation}
It can be seen that the right-hand side of Eq. (\ref{nplanar-34a}) does not contain
oscillatory terms. If $n$ is a multiple of 4, it can be shown by replacing $y$
by $iy$ in Eq. (\ref{nplanar-34a}) that
\begin{equation}
\sum_{k=0}^{n-1}(-1)^kf_{nk}^2(y)=\frac{4}{n}
\sum_{l=1}^{n/4}\cos\left[2y\cos\left(\frac{\pi (2l-1)}{n}\right)\right],
\label{nplanar-34b}
\end{equation}
which does not contain exponential terms.

For odd $n$, the planar n-dimensional cosexponential function $f_{nk}(y)$ is
related to the n-dimensional cosexponential function $g_{nk}(y)$ discussed in
the previous chapter also by the relation 
\begin{equation}
f_{nk}(y)=(-1)^kg_{nk}(-y),
\label{nplanar-34c}
\end{equation}
as can be seen by comparing the series for the two classes of functions.
For values of the form $n=4p+2$, $p=0,1,2,...$, the planar n-dimensional
cosexponential function $f_{nk}(y)$ is related to the n-dimensional
cosexponential function $g_{nk}(y)$ by the relation  
\begin{equation}
f_{nk}(y)=e^{-i\pi k/2}g_{nk}(iy).
\label{nplanar-34cx}
\end{equation}

Addition theorems for the planar n-dimensional cosexponential functions can be
obtained from the relation $\exp h_1(y+z)=\exp h_1 y \cdot\exp h_1 z $, by
substituting the expression of the exponentials as given in Eq. (\ref{nplanar-28b})
for $k=1$, $e^{h_1 y}=f_{n0}(y)+h_1f_{n1}(y)+\cdots+h_{n-1} f_{n,n-1}(y)$,
\begin{eqnarray}
\lefteqn{f_{nk}(y+z)=f_{n0}(y)f_{nk}(z)+f_{n1}(y)f_{n,k-1}(z)+\cdots+f_{nk}(y)f_{n0}(z)\nonumber}\\
&&-f_{n,k+1}(y)f_{n,n-1}(z)-f_{n,k+2}(y)f_{n,n-2}(z)-\cdots-f_{n,n-1}(y)f_{n,k+1}(z) ,
\nonumber\\
&&
\label{nplanar-35a}
\end{eqnarray}
where $k=0,1,...,n-1$.\index{cosexponential functions, planar
n-complex!addition theorems}
For $y=z$ the relations (\ref{nplanar-35a}) take the form
\begin{eqnarray}
\lefteqn{f_{nk}(2y)=f_{n0}(y)f_{nk}(y)+f_{n1}(y)f_{n,k-1}(y)+\cdots+f_{nk}(y)f_{n0}(y)\nonumber}\\
&&-f_{n,k+1}(y)f_{n,n-1}(y)-f_{n,k+2}(y)f_{n,n-2}(y)-\cdots-f_{n,n-1}(y)f_{n,k+1}(y) ,
\nonumber\\
&&
\label{nplanar-35b}
\end{eqnarray}
where $k=0,1,...,n-1$.
For $y=-z$ the relations (\ref{nplanar-35a}) and (\ref{nplanar-29}) yield
\begin{equation}
f_{n0}(y)f_{n0}(-y)
-f_{n1}(y)f_{n,n-1}(-y)-f_{n2}(y)f_{n,n-2}(-y)-\cdots-f_{n,n-1}(y)f_{n1}(-y)=1 ,
\label{nplanar-36a}
\end{equation}
\begin{eqnarray}
\lefteqn{f_{n0}(y)f_{nk}(-y)+f_{n1}(y)f_{n,k-1}(-y)+\cdots+f_{nk}(y)f_{n0}(-y)\nonumber}\\
&&-f_{n,k+1}(y)f_{n,n-1}(-y)-f_{n,k+2}(y)f_{n,n-2}(-y)-\cdots-f_{n,n-1}(y)f_{n,k+1}(-y)=0 ,
\nonumber\\
&&
\label{nplanar-36b}
\end{eqnarray}
where $k=1,...,n-1$.

From Eq. (\ref{nplanar-28bx}) it can be shown, for even $k$ and 
natural numbers $l$, that
\begin{equation}
\left(\sum_{p=0}^{n-1}(-1)^{[kp/n]}h_{kp-n[kp/n]}g_{np}(y)\right)^l
=\sum_{p=0}^{n-1}(-1)^{[kp/n]}h_{kp-n[kp/n]}g_{np}(ly), 
\label{nplanar-37x}
\end{equation}
where $k=0,1,...,n-1$.
For odd $k$ and natural numbers $l$, Eq. (\ref{nplanar-28b}) implies
\begin{equation}
\left(\sum_{p=0}^{n-1}(-1)^{[kp/n]}h_{kp-n[kp/n]}f_{np}(y)\right)^l
=\sum_{p=0}^{n-1}(-1)^{[kp/n]}h_{kp-n[kp/n]}f_{np}(ly), 
\label{nplanar-37}
\end{equation}
where $k=0,1,...,n-1$.
For $k=1$ the relation (\ref{nplanar-37}) is
\begin{eqnarray}
\lefteqn{\left\{f_{n0}(y)+h_1f_{n1}(y)+\cdots+h_{n-1}f_{n,n-1}(y)\right\}^l
=f_{n0}(ly)+h_1f_{n1}(ly)+\cdots+h_{n-1}f_{n,n-1}(ly).\nonumber}\\
&&
\label{nplanar-37b}
\end{eqnarray}

If
\begin{equation}
a_k=\sum_{p=0}^{n-1}f_{np}(y)\cos\left(\frac{\pi (2k-1)p}{n}\right), 
\label{nplanar-38a}
\end{equation}
and
\begin{equation}
b_k=\sum_{p=0}^{n-1}f_{np}(y)\sin\left(\frac{\pi (2k-1)p}{n}\right), 
\label{nplanar-38b}
\end{equation}
for $k=1,...,n$, where $f_{np}(y)$ are the planar cosexponential functions in Eq. (\ref{nplanar-30}), it
can be shown that 
\begin{equation}
a_k=\exp\left[y\cos\left(\frac{\pi (2k-1)}{n}\right)\right]
\cos\left[y\sin\left(\frac{\pi (2k-1)}{n}\right)\right], 
\label{nplanar-39a}
\end{equation}
\begin{equation}
b_k=\exp\left[y\cos\left(\frac{\pi (2k-1)}{n}\right)\right]
\sin\left[y\sin\left(\frac{\pi (2k-1)}{n}\right)\right], 
\label{nplanar-39b}
\end{equation}
for $k=1,...,n$. If
\begin{equation}
G_k^2=a_k^2+b_k^2, 
\label{nplanar-40}
\end{equation}
from Eqs. (\ref{nplanar-39a}) and (\ref{nplanar-39b}) it results that
\begin{equation}
G_k^2=\exp\left[2y\cos\left(\frac{\pi (2k-1)}{n}\right)\right], 
\label{nplanar-41}
\end{equation}
for $k=1,...,n$.
Then the planar n-dimensional cosexponential functions have the property that
\begin{equation}
\prod_{p=1}^{n/2}G_p^2=1.
\label{nplanar-43a}
\end{equation}

The planar n-dimensional cosexponential functions are solutions of the
$n^{\rm th}$-order differential equation
\begin{equation}
\frac{d^n\zeta}{du^n}=-\zeta ,
\label{nplanar-44}
\end{equation}\index{cosexponential functions, planar n-complex!differential equations}
whose solutions are of the form
$\zeta(u)=A_0f_{n0}(u)+A_1f_{n1}(u)+\cdots+A_{n-1}f_{n,n-1}(u).$ 
It can be checked that the derivatives of the planar cosexponential functions
are related by
\begin{equation}
\frac{df_{n0}}{du}=-f_{n,n-1}, \:
\frac{df_{n1}}{du}=f_{n0}, \:...,
\frac{df_{n,n-2}}{du}=f_{n,n-3} ,
\frac{df_{n,n-1}}{du}=f_{n,n-2} .
\label{nplanar-45}
\end{equation}

\subsection{Exponential and trigonometric forms of planar n-complex numbers}

In order to obtain the exponential and trigonometric forms of n-complex
numbers, a canonical base
$e_1,\tilde e_1,...,e_{n/2},\tilde e_{n/2}$ for the planar n-complex numbers
will be introduced by the relations 
\begin{equation}
\left(
\begin{array}{c}
\vdots\\
e_k\\
\tilde e_k\\
\vdots
\end{array}\right)
=\left(
\begin{array}{ccccc}
\vdots&\vdots& &\vdots&\vdots\\
\frac{2}{n}&\frac{2}{n}\cos\frac{\pi (2k-1)}{n}&\cdots&\frac{2}{n}\cos\frac{\pi (2k-1)(n-2)}{n}&\frac{2}{n}\cos\frac{\pi (2k-1)(n-1)}{n}\\
0&\frac{2}{n}\sin\frac{\pi (2k-1)}{n}&\cdots&\frac{2}{n}\sin\frac{\pi (2k-1)(n-2)}{n}&\frac{2}{n}\sin\frac{\pi (2k-1)(n-1)}{n}\\
\vdots&\vdots&&\vdots&\vdots
\end{array}
\right)
\left(
\begin{array}{c}
1\\h_1\\
\vdots\\
h_{n-1}
\end{array}
\right),
\label{nplanar-e11}
\end{equation}\index{canonical base!planar n-complex}
where $k=1, 2, ... , n/2$.

The multiplication relations for the bases $e_k, \tilde e_k$ are
\begin{eqnarray}
e_k^2=e_k, \tilde e_k^2=-e_k, e_k \tilde e_k=\tilde e_k , e_ke_l=0, e_k\tilde e_l=0, \tilde e_k\tilde e_l=0, k\not=l,
\label{nplanar-e12a}
\end{eqnarray}
where $k,l=1,...,n/2$.
The moduli of the bases $e_k, \tilde e_k$ are
\begin{equation}
|e_k|=\sqrt{\frac{2}{n}}, |\tilde e_k|=\sqrt{\frac{2}{n}}.
\label{nplanar-e12c}
\end{equation}
It can be shown that
\begin{eqnarray}
x_0+h_1x_1+\cdots+h_{n-1}x_{n-1}= 
\sum_{k=1}^{n/2}(e_k v_k+\tilde e_k \tilde v_k).
\label{nplanar-e13a}
\end{eqnarray}\index{canonical form!planar n-complex}
The relation (\ref{nplanar-e13a} gives the canonical form of a planar n-complex
number.

Using the properties of the bases in Eqs. (\ref{nplanar-e11}) it can
be shown that
\begin{equation}
\exp(\tilde e_k\phi_k)=1-e_k+e_k\cos\phi_k+\tilde e_k\sin\phi_k ,
\label{nplanar-46a}
\end{equation}
\begin{equation}
\exp(e_k\ln\rho_k)=1-e_k+e_k\rho_k ,
\label{nplanar-46b}
\end{equation}
By multiplying the relations (\ref{nplanar-46a}), (\ref{nplanar-46b}) it results that 
\begin{equation}
\exp\left[\sum_{k=1}^{n/2}
(e_k\ln \rho_k+\tilde e_k\phi_k)\right] 
=\sum_{k=1}^{n/2}(e_k v_k+\tilde e_k \tilde v_k),
\label{nplanar-47a}
\end{equation}
where the fact has ben used that 
\begin{equation}
\sum_{k=1}^{n/2} e_k =1, 
\label{nplanar-47b}
\end{equation}
the latter relation being a consequence of Eqs. (\ref{nplanar-e11}) and (\ref{nplanar-15}).

By comparing Eqs. (\ref{nplanar-e13a}) and (\ref{nplanar-47a}), it can be seen that
\begin{eqnarray}
x_0+h_1x_1+\cdots+h_{n-1}x_{n-1}=
\exp\left[\sum_{k=1}^{n/2}
(e_k\ln \rho_k+\tilde e_k\phi_k)\right] .
\label{nplanar-49a}
\end{eqnarray}
Using the expression of the bases in Eq. (\ref{nplanar-e11}) yields
the exponential form of the n-complex number
$u=x_0+h_1x_1+\cdots+h_{n-1}x_{n-1}$ as
\begin{eqnarray}
u=\rho\exp\left\{\sum_{p=1}^{n-1}h_p\left[
-\frac{2}{n}\sum_{k=2}^{n/2}
\cos\left(\frac{\pi (2k-1)p}{n}\right)\ln\tan\psi_{k-1}
\right]
+\sum_{k=1}^{n/2}\tilde e_k\phi_k
\right\},
\label{nplanar-50a}
\end{eqnarray}\index{exponential form!planar n-complex}
where $\rho$ is the amplitude defined in Eq. (\ref{nplanar-6b}), and has according to
Eq. (\ref{nplanar-9c}) the expression
\begin{equation}
\rho=\left(\rho_1^2\cdots\rho_{n/2}^2\right)^{1/n}.
\label{nplanar-50aa}
\end{equation}

It can be checked with the aid of Eq. (\ref{nplanar-46a}) that
the n-complex number $u$ can also be written as
\begin{eqnarray}
x_0+h_1x_1+\cdots+h_{n-1}x_{n-1}
=\left(\sum_{k=1}^{n/2}e_k \rho_k\right)
\exp\left(\sum_{k=1}^{n/2}\tilde e_k\phi_k\right).
\label{nplanar-51a}
\end{eqnarray}
Writing in Eq. (\ref{nplanar-51a}) the radius $\rho_1$, Eq.
(\ref{nplanar-20b}), as a factor and expressing the 
variables in terms of the planar
angles with the aid of Eq. (\ref{nplanar-19b}) yields the
trigonometric form of the n-complex number $u$ as
\begin{eqnarray}
\lefteqn{u=d
\left(\frac{n}{2}\right)^{1/2}
\left(1+\frac{1}{\tan^2\psi_1}+\frac{1}{\tan^2\psi_2}+\cdots
+\frac{1}{\tan^2\psi_{n/2-1}}\right)^{-1/2}\nonumber}\\
&&\left(e_1+\sum_{k=2}^{n/2}\frac{e_k}{\tan\psi_{k-1}}\right)
\exp\left(\sum_{k=1}^{n/2}\tilde e_k\phi_k\right).
\label{nplanar-52a}
\end{eqnarray}\index{trigonometric form!planar n-complex}
In Eq. (\ref{nplanar-52a}),  the n-complex
number $u$, written in trigonometric form, is the
product of the modulus $d$, of a part depending on the planar
angles $\psi_1,...,\psi_{n/2-1}$, and of a factor depending
on the azimuthal angles $\phi_1,...,\phi_{n/2}$. 
Although the modulus of a product of n-complex numbers is not equal in
general to the product of the moduli of the factors,
it can be checked that the modulus of the factors in Eq. (\ref{nplanar-52a}) are
\begin{eqnarray}
\lefteqn{
\left|e_1+\sum_{k=2}^{n/2}\frac{e_k}{\tan\psi_{k-1}}\right|\nonumber}\\
&&=\left(\frac{2}{n}\right)^{1/2}
\left(1+\frac{1}{\tan^2\psi_1}+\frac{1}{\tan^2\psi_2}+\cdots
+\frac{1}{\tan^2\psi_{n/2-1}}\right)^{1/2},
\label{nplanar-52c}
\end{eqnarray}
and
\begin{eqnarray}
\left|\exp\left(\sum_{k=1}^{n/2}\tilde e_k\phi_k\right)\right|=1.
\label{nplanar-52e}
\end{eqnarray}

The modulus $d$ in Eqs. (\ref{nplanar-52a}) can be expressed in terms
of the amplitude $\rho$ as
\begin{eqnarray}
\lefteqn{d=\rho \frac{2^{(n-2)/2n}}{\sqrt{n}}
\left(\tan\psi_1\cdots\tan\psi_{n/2-1}\right)^{2/n}\nonumber}\\
&&\left(1+\frac{1}{\tan^2\psi_1}+\frac{1}{\tan^2\psi_2}+\cdots
+\frac{1}{\tan^2\psi_{n/2-1}}\right)^{1/2}.
\label{nplanar-53a}
\end{eqnarray}

\subsection{Elementary functions of a planar n-complex variable}

The logarithm $u_1$ of the n-complex number $u$, $u_1=\ln u$, can be defined
as the solution of the equation
\begin{equation}
u=e^{u_1} .
\label{nplanar-54}
\end{equation}
The relation (\ref{nplanar-47a}) shows that $\ln u$ exists  
as an n-complex function with real components for all values of
$x_0,...,x_{n-1}$ for which $\rho\not=0$.
The expression of the logarithm, obtained from Eq. (\ref{nplanar-49a}) is
\begin{equation}
\ln u=\sum_{k=1}^{n/2}
(e_k\ln \rho_k+\tilde e_k\phi_k).
\label{nplanar-55a}
\end{equation}\index{logarithm!planar n-complex}
An expression of the logarithm depending on the amplitude $\rho$ can be
obtained from the exponential forms in Eq. (\ref{nplanar-50a}) as
\begin{eqnarray}
\ln u=\ln \rho+\sum_{p=1}^{n-1}h_p\left[
-\frac{2}{n}\sum_{k=2}^{n/2}
\cos\left(\frac{\pi (2k-1)p}{n}\right)\ln\tan\psi_{k-1}
\right]+\sum_{k=1}^{n/2}\tilde e_k\phi_k.
\label{nplanar-56a}
\end{eqnarray}

The function $\ln u$ is multivalued because of the presence of the terms 
$\tilde e_k\phi_k$.
It can be inferred from Eqs. (\ref{nplanar-21b})-(\ref{nplanar-21e}) and (\ref{nplanar-24}) that
\begin{equation}
\ln(uu^\prime)=\ln u+\ln u^\prime ,
\label{nplanar-57}
\end{equation}
up to integer multiples of $2\pi\tilde e_k, k=1,...,n/2$.

The power function $u^m$ can be defined for real values of $m$ as
\begin{equation}
u^m=e^{m\ln u} .
\label{nplanar-58}
\end{equation}
Using the expression of $\ln u$ in Eq. (\ref{nplanar-55a}) yields
\begin{equation}
u^m=\sum_{k=1}^{n/2}
\rho_k^m(e_k\cos m\phi_k+\tilde e_k\sin m\phi_k).
\label{nplanar-59a}
\end{equation}\index{power function!planar n-complex}
The power function is multivalued unless $m$ is an integer. 
For integer $m$, it can be inferred from Eq. (\ref{nplanar-57}) that
\begin{equation}
(uu^\prime)^m=u^m\:u^{\prime m} .
\label{nplanar-59}
\end{equation}

The trigonometric functions of the hypercomplex variable
$u$ and the addition theorems for these functions have been written in Eqs.
~(1.57)-(1.60). 
In order to obtain expressions for the trigonometric functions of n-complex
variables, these will be expressed with the aid of the imaginary unit $i$ as
\begin{equation}
\cos u=\frac{1}{2}(e^{iu}+e^{-iu}),\:\sin u=\frac{1}{2i}(e^{iu}-e^{-iu}).
\label{nplanar-64}
\end{equation}
The imaginary unit $i$ is used for the convenience of notations, and it does
not appear in the final results.
The validity of Eq. (\ref{nplanar-64}) can be checked by comparing the series for the
two sides of the relations.
Since the expression of the exponential function $e^{h_k y}$ in terms of the
units $1, h_1, ... h_{n-1}$ given in Eq. (\ref{nplanar-28b}) depends on the planar
cosexponential functions $f_{np}(y)$, the expression of the trigonometric
functions will depend on the functions
$f_{p+}^{(c)}(y)=(1/2)[f_{np}(iy)+f_{np}(-iy)]$ and 
$f_{p-}^{(c)}(y)=(1/2i)[f_{np}(iy)-f_{np}(-iy)]$, 
\begin{equation}
\cos(h_k y)=\sum_{p=0}^{n-1}(-1)^{[kp/n]}h_{kp-n[kp/n]}f_{p+}^{(c)}(y),
\label{nplanar-66a}
\end{equation}\index{trigonometric functions, expressions!planar n-complex}
\begin{equation}
\sin(h_k y)=\sum_{p=0}^{n-1}(-1)^{[kp/n]}h_{kp-n[kp/n]}f_{p-}^{(c)}(y),
\label{nplanar-66b}
\end{equation}
where
\begin{eqnarray}
\lefteqn{f_{p+}^{(c)}(y)=\frac{1}{n}\sum_{l=1}^{n}\left\{
\cos\left[y\cos\left(\frac{\pi (2l-1)}{n}\right)\right]
\cosh\left[y\sin\left(\frac{\pi (2l-1)}{n}\right)\right]
\cos\left(\frac{\pi (2l-1)p}{n}\right)\right.\nonumber}\\
&&\left.-\sin\left[y\cos\left(\frac{\pi (2l-1)}{n}\right)\right]
\sinh\left[y\sin\left(\frac{\pi (2l-1)}{n}\right)\right]
\sin\left(\frac{\pi (2l-1)p}{n}\right)
\right\},
\label{nplanar-65a}
\end{eqnarray}
\begin{eqnarray}
\lefteqn{f_{p-}^{(c)}(y)=\frac{1}{n}\sum_{l=1}^{n}\left\{
\sin\left[y\cos\left(\frac{\pi (2l-1)}{n}\right)\right]
\cosh\left[y\sin\left(\frac{\pi (2l-1)}{n}\right)\right]
\cos\left(\frac{\pi (2l-1)p}{n}\right)\right.\nonumber}\\
&&\left.+\cos\left[y\cos\left(\frac{\pi (2l-1)}{n}\right)\right]
\sinh\left[y\sin\left(\frac{\pi (2l-1)}{n}\right)\right]
\sin\left(\frac{\pi (2l-1)p}{n}\right)
\right\}.
\label{nplanar-65b}
\end{eqnarray}

The hyperbolic functions of the hypercomplex variable
$u$ and the addition theorems for these functions have been written in Eqs.
~(1.62)-(1.65). 
In order to obtain expressions for the hyperbolic functions of n-complex
variables, these will be expressed as
\begin{equation}
\cosh u=\frac{1}{2}(e^{u}+e^{-u}),\:\sinh u=\frac{1}{2}(e^{u}-e^{-u}).
\label{nplanar-70}
\end{equation}
The validity of Eq. (\ref{nplanar-70}) can be checked by comparing the series for the
two sides of the relations.
Since the expression of the exponential function $e^{h_k y}$ in terms of the
units $1, h_1, ... h_{n-1}$ given in Eq. (\ref{nplanar-28b}) depends on the planar
cosexponential functions $f_{np}(y)$, the expression of the hyperbolic
functions will depend on the even part $f_{p+}(y)=(1/2)[f_{np}(y)+f_{np}(-y)]$ and on
the odd part $f_{p-}(y)=(1/2)[f_{np}(y)-f_{np}(-y)]$ of $f_{np}$, 
\begin{equation}
\cosh(h_k y)=\sum_{p=0}^{n-1}(-1)^{[kp/n]}h_{kp-n[kp/n]}f_{p+}(y),
\label{nplanar-71a}
\end{equation}\index{hyperbolic functions, expressions!planar n-complex}
\begin{equation}
\sinh(h_k y)=\sum_{p=0}^{n-1}(-1)^{[kp/n]}h_{kp-n[kp/n]}f_{p-}(y),
\label{nplanar-71b}
\end{equation}
where
\begin{eqnarray}
\lefteqn{f_{p+}(y)=\frac{1}{n}\sum_{(2l-1)=1}^{n}\left\{
\cosh\left[y\cos\left(\frac{\pi (2l-1)}{n}\right)\right]
\cos\left[y\sin\left(\frac{\pi (2l-1)}{n}\right)\right]
\cos\left(\frac{\pi (2l-1)p}{n}\right)\right.\nonumber}\\
&&\left.+\sinh\left[y\cos\left(\frac{\pi (2l-1)}{n}\right)\right]
\sin\left[y\sin\left(\frac{\pi (2l-1)}{n}\right)\right]
\sin\left(\frac{\pi (2l-1)p}{n}\right)
\right\},
\label{nplanar-72a}
\end{eqnarray}
\begin{eqnarray}
\lefteqn{f_{p-}(y)=\frac{1}{n}\sum_{l=1}^{n}\left\{
\sinh\left[y\cos\left(\frac{\pi (2l-1)}{n}\right)\right]
\cos\left[y\sin\left(\frac{\pi (2l-1)}{n}\right)\right]
\cos\left(\frac{\pi (2l-1)p}{n}\right)\right.\nonumber}\\
&&\left.+\cosh\left[y\cos\left(\frac{\pi (2l-1)}{n}\right)\right]
\sin\left[y\sin\left(\frac{\pi (2l-1)}{n}\right)\right]
\sin\left(\frac{\pi (2l-1)p}{n}\right)
\right\}.
\label{nplanar-72b}
\end{eqnarray}

The exponential, trigonometric and hyperbolic functions can also be expressed
with the aid of the bases introduced in Eq. (\ref{nplanar-e11}).
Using the expression of the n-complex number in Eq. (\ref{nplanar-e13a})
yields for the exponential of the n-complex variable $u$
\begin{eqnarray}
e^u= 
\sum_{k=1}^{n/2}e^{v_k}\left(e_k \cos \tilde v_k+\tilde e_k \sin\tilde
v_k\right).
\label{nplanar-73a}
\end{eqnarray}

The trigonometric functions can be obtained from Eq. (\ref{nplanar-73a}) 
with the aid of Eqs. (\ref{nplanar-64}). The trigonometric functions of the
n-complex variable $u$ are
\begin{equation}
\cos u=\sum_{k=1}^{n/2}\left(e_k \cos v_k\cosh \tilde v_k
-\tilde e_k \sin v_k\sinh\tilde v_k\right),
\label{nplanar-74a}
\end{equation}\index{trigonometric functions, expressions!planar n-complex}
\begin{equation}
\sin u= 
\sum_{k=1}^{n/2}\left(e_k \sin v_k\cosh \tilde v_k
+\tilde e_k \cos v_k\sinh\tilde v_k\right).
\label{nplanar-74b}
\end{equation}

The hyperbolic functions can be obtained from Eq. (\ref{nplanar-73a}) 
with the aid of Eqs. (\ref{nplanar-70}). The hyperbolic functions of the
n-complex variable $u$ are
\begin{equation}
\cosh u=
\sum_{k=1}^{n/2}\left(e_k \cosh v_k\cos \tilde v_k
+\tilde e_k \sinh v_k\sin\tilde v_k\right),
\label{nplanar-75a}
\end{equation}\index{hyperbolic functions, expressions!planar n-complex}
\begin{equation}
\sinh u=
\sum_{k=1}^{n/2}\left(e_k \sinh v_k\cos \tilde v_k
+\tilde e_k \cosh v_k\sin\tilde v_k\right).
\label{nplanar-75b}
\end{equation}

\subsection{Power series of planar n-complex numbers}

An n-complex series is an infinite sum of the form
\begin{equation}
a_0+a_1+a_2+\cdots+a_n+\cdots , 
\label{nplanar-76}
\end{equation}\index{series!planar n-complex}
where the coefficients $a_n$ are n-complex numbers. The convergence of 
the series (\ref{nplanar-76}) can be defined in terms of the convergence of its $n$
real components. The convergence of a n-complex series can also be studied
using n-complex variables. The main criterion for absolute convergence 
remains the comparison theorem, but this requires a number of inequalities
which will be discussed further.

The modulus $d=|u|$ of an n-complex number $u$ has been defined in Eq.
(\ref{nplanar-10}). Since $|x_0|\leq |u|, |x_1|\leq |u|,..., |x_{n-1}|\leq |u|$, a
property of absolute convergence established via a comparison theorem based on
the modulus of the series (\ref{nplanar-76}) will ensure the absolute convergence of
each real component of that series.

The modulus of the sum $u_1+u_2$ of the n-complex numbers $u_1, u_2$ fulfils
the inequality
\begin{equation}
||u^\prime|-|u^{\prime\prime}||\leq |u^\prime+u^{\prime\prime}|\leq
|u^\prime|+|u^{\prime\prime}| . 
\label{nplanar-78}
\end{equation}\index{modulus, inequalities!planar n-complex}
For the product, the relation is 
\begin{equation}
|u^\prime u^{\prime\prime}|\leq \sqrt{\frac{n}{2}}|u^\prime||u^{\prime\prime}| ,
\label{nplanar-79}
\end{equation}
as can be shown from Eq. (\ref{nplanar-17}). The relation (\ref{nplanar-79})
replaces the relation of equality extant between 2-dimensional regular complex
numbers. 

For $u=u^\prime$ Eq. (\ref{nplanar-79}) becomes
\begin{equation}
|u^2|\leq \sqrt{\frac{n}{2}} |u|^2 ,
\label{nplanar-80}
\end{equation}
and in general
\begin{equation}
|u^l|\leq \left(\frac{n}{2}\right)^{(l-1)/2}|u|^l ,
\label{nplanar-81}
\end{equation}
where $l$ is a natural number.
From Eqs. (\ref{nplanar-79}) and (\ref{nplanar-81}) it results that
\begin{equation}
|au^l|\leq \left(\frac{n}{2}\right)^{l/2} |a| |u|^l .
\label{nplanar-82}
\end{equation}

A power series of the n-complex variable $u$ is a series of the form
\begin{equation}
a_0+a_1 u + a_2 u^2+\cdots +a_l u^l+\cdots .
\label{nplanar-83}
\end{equation}\index{power series!planar n-complex}
Since
\begin{equation}
\left|\sum_{l=0}^\infty a_l u^l\right| \leq  \sum_{l=0}^\infty
(n/2)^{l/2} |a_l| |u|^l ,
\label{nplanar-84}
\end{equation}
a sufficient condition for the absolute convergence of this series is that
\begin{equation}
\lim_{l\rightarrow \infty}\frac{\sqrt{n/2}|a_{l+1}||u|}{|a_l|}<1 .
\label{nplanar-85}
\end{equation}
Thus the series is absolutely convergent for 
\begin{equation}
|u|<c,
\label{nplanar-86}
\end{equation}\index{convergence of power series!planar n-complex}
where 
\begin{equation}
c=\lim_{l\rightarrow\infty} \frac{|a_l|}{\sqrt{n/2}|a_{l+1}|} .
\label{nplanar-87}
\end{equation}

The convergence of the series (\ref{nplanar-83}) can be also studied with the aid of
the formula (\ref{nplanar-59a}) which is
valid for any values of $x_0,...,x_{n-1}$, as mentioned previously.
If $a_l=\sum_{p=0}^{n-1} h_p a_{lp}$, and
\begin{equation}
A_{lk}=\sum_{p=0}^{n-1} a_{lp}\cos\frac{\pi (2k-1)p}{n},
\label{nplanar-88b}
\end{equation}
\begin{equation}
\tilde A_{lk}=\sum_{p=0}^{n-1} a_{lp}\sin\frac{\pi (2k-1)p}{n},
\label{nplanar-88c}
\end{equation}
where $k=1,...,n/2$, the series (\ref{nplanar-83}) can be written as
\begin{equation}
\sum_{l=0}^\infty \left[
\sum_{k=1}^{n/2}
(e_k A_{lk}+\tilde e_k\tilde A_{lk})(e_k v_k+\tilde e_k\tilde v_k)^l 
\right].
\label{nplanar-89a}
\end{equation}
The series in Eq. (\ref{nplanar-89a}) can be regarded as the sum of the $n/2$ series
obtained from each value of $k$, so that the series in Eq. (\ref{nplanar-83}) is
absolutely convergent for   
\begin{equation}
\rho_k<c_k, 
\label{nplanar-90}
\end{equation}\index{convergence, region of!planar n-complex}
for $k=1,..., n/2$, where 
\begin{equation}
c_k=\lim_{l\rightarrow\infty} \frac
{\left[A_{lk}^2+\tilde A_{lk}^2\right]^{1/2}}
{\left[A_{l+1,k}^2+\tilde A_{l+1,k}^2\right]^{1/2}} .
\label{nplanar-91}
\end{equation}
The relations (\ref{nplanar-90}) show that the region of convergence of the series
(\ref{nplanar-83}) is an n-dimensional cylinder.

It can be shown that $c=\sqrt{2/n}\;{\rm
min}(c_+,c_-,c_1,...,c_{n/2-1})$, where ${\rm min}$ designates the smallest of
the numbers in the argument of this function. Using the expression of $|u|$ in
Eq. (\ref{nplanar-17}),  it can be seen that the spherical region of
convergence defined in Eqs. (\ref{nplanar-86}), (\ref{nplanar-87}) is a subset of the
cylindrical region of convergence defined in Eqs. (\ref{nplanar-90}) and (\ref{nplanar-91}).

\subsection{Analytic functions of planar n-complex variables}

The analytic functions of the hypercomplex variable $u$ and the series 
expansion of functions have been discussed in Eqs. ~(1.85)-(1.93).
If the n-complex function $f(u)$
of the n-complex variable $u$ is written in terms of 
the real functions $P_k(x_0,...,x_{n-1}), k=0,1,...,n-1$ of the real
variables $x_0,x_1,...,x_{n-1}$ as 
\begin{equation}
f(u)=\sum_{k=0}^{n-1}h_kP_k(x_0,...,x_{n-1}),
\label{nplanar-h93}
\end{equation}\index{functions, real components!planar n-complex}
where $h_0=1$, then relations of equality 
exist between the partial derivatives of the functions $P_k$. 
The derivative of the function $f$ can be written as
\begin{eqnarray}
\lim_{\Delta u\rightarrow 0}\frac{1}{\Delta u} 
\sum_{k=0}^{n-1}\left(h_k\sum_{l=0}^{n-1}
\frac{\partial P_k}{\partial x_l}\Delta x_l\right),
\label{nplanar-h94}
\end{eqnarray}\index{derivative, independence of direction!planar n-complex}
where
\begin{equation}
\Delta u=\sum_{k=0}^{n-1}h_l\Delta x_l.
\label{nplanar-h94a}
\end{equation}
The relations between the partials derivatives of the functions $P_k$ are
obtained by setting successively in   
Eq. (\ref{nplanar-h94}) $\Delta u=h_l\Delta x_l$, for $l=0,1,...,n-1$, and equating the
resulting expressions. 
The relations are \index{relations between partial derivatives!planar n-complex}
\begin{equation}
\frac{\partial P_k}{\partial x_0} = \frac{\partial P_{k+1}}{\partial x_1} 
=\cdots=\frac{\partial P_{n-1}}{\partial x_{n-k-1}} 
= -\frac{\partial P_0}{\partial x_{n-k}}=\cdots
=-\frac{\partial P_{k-1}}{\partial x_{n-1}}, 
\label{nplanar-h95}
\end{equation}
for $k=0,1,...,n-1$.
The relations (\ref{nplanar-h95}) are analogous to the Riemann relations
for the real and imaginary components of a complex function. 
It can be shown from Eqs. (\ref{nplanar-h95}) that the components $P_k$ fulfil the
second-order equations\index{relations between second-order derivatives!planar n-complex}
\begin{eqnarray}
\lefteqn{\frac{\partial^2 P_k}{\partial x_0\partial x_l}
=\frac{\partial^2 P_k}{\partial x_1\partial x_{l-1}}
=\cdots=
\frac{\partial^2 P_k}{\partial x_{[l/2]}\partial x_{l-[l/2]}}}\nonumber\\
&&=-\frac{\partial^2 P_k}{\partial x_{l+1}\partial x_{n-1}}
=-\frac{\partial^2 P_k}{\partial x_{l+2}\partial x_{n-2}}
=\cdots
=-\frac{\partial^2 P_k}{\partial x_{l+1+[(n-l-2)/2]}
\partial x_{n-1-[(n-l-2)/2]}} ,
\label{nplanar-96}
\end{eqnarray}
for $k,l=0,1,...,n-1$.

\subsection{Integrals of planar n-complex functions}

The singularities of n-complex functions arise from terms of the form
$1/(u-u_0)^m$, with $m>0$. Functions containing such terms are singular not
only at $u=u_0$, but also at all points of the hypersurfaces
passing through the pole $u_0$ and which are parallel to the nodal hypersurfaces. 

The integral of an n-complex function between two points $A, B$ along a path
situated in a region free of singularities is independent of path, which means
that the integral of an analytic function along a loop situated in a region
free of singularities is zero,
\begin{equation}
\oint_\Gamma f(u) du = 0,
\label{nplanar-111}
\end{equation}
where it is supposed that a surface $\Sigma$ spanning 
the closed loop $\Gamma$ is not intersected by any of
the hypersurfaces associated with the
singularities of the function $f(u)$. Using the expression, Eq. (\ref{nplanar-h93}),
for $f(u)$ and the fact that 
\begin{eqnarray}
du=\sum_{k=0}^{n-1}h_k dx_k, 
\label{nplanar-111a}
\end{eqnarray}
the explicit form of the integral in Eq. (\ref{nplanar-111}) is
\begin{eqnarray}
\oint _\Gamma f(u) du = \oint_\Gamma
\sum_{k=0}^{n-1}h_k\sum_{l=0}^{n-1}(-1)^{[(n-k-1+l)/n]}
P_l dx_{k-l+n[(n-k-1+l)/n]}.
\label{nplanar-112}
\end{eqnarray}\index{integrals, path!planar n-complex}

If the functions $P_k$ are regular on a surface $\Sigma$
spanning the loop $\Gamma$,
the integral along the loop $\Gamma$ can be transformed in an integral over the
surface $\Sigma$ of terms of the form
$\partial P_l/\partial x_{k-m+n[(n-k+m-1)/n]} 
- (-1)^s \partial P_m/\partial x_{k-l+n[(n-k+l-1)/n]}$, where 
$s=[(n-k+m-1)/n]-[(n-k+l-1)/n]$.
These terms are equal to zero by Eqs. (\ref{nplanar-h95}), and this
proves Eq. (\ref{nplanar-111}). 

The integral of the function $(u-u_0)^m$ on a closed loop $\Gamma$ is equal to
zero for $m$ a positive or negative integer not equal to -1,
\begin{equation}
\oint_\Gamma (u-u_0)^m du = 0, \:\: m \:\:{\rm integer},\: m\not=-1 .
\label{nplanar-112b}
\end{equation}
This is due to the fact that $\int (u-u_0)^m du=(u-u_0)^{m+1}/(m+1), $ and to
the fact that the function $(u-u_0)^{m+1}$ is singlevalued for $n$ an integer.

The integral $\oint_\Gamma du/(u-u_0)$ can be calculated using the exponential
form, Eq. (\ref{nplanar-50a}), for the difference $u-u_0$, 
\begin{eqnarray}
u-u_0=\rho\exp\left\{\sum_{p=1}^{n-1}h_p\left[
-\frac{2}{n}\sum_{k=2}^{n/2}
\cos\left(\frac{2\pi kp}{n}\right)\ln\tan\psi_{k-1}
\right]
+\sum_{k=1}^{n/2}\tilde e_k\phi_k\right\}.
\label{nplanar-113a}
\end{eqnarray}
Thus the quantity $du/(u-u_0)$ is
\begin{eqnarray}
\frac{du}{u-u_0}=
\frac{d\rho}{\rho}+\sum_{p=1}^{n-1}h_p\left[
-\frac{2}{n}\sum_{k=2}^{n/2}
\cos\left(\frac{2\pi kp}{n}\right)d\ln\tan\psi_{k-1}\right]
+\sum_{k=1}^{n/2}\tilde e_kd\phi_k.
\label{nplanar-114a}
\end{eqnarray}
Since $\rho$ and $\ln(\tan\psi_{k-1})$ are singlevalued variables, it follows
that 
$\oint_\Gamma d\rho/\rho =0$, and 
$\oint_\Gamma d(\ln\tan\psi_{k-1})=0$.
On the other hand, since $\phi_k$ are cyclic variables, they may give
contributions to the integral around the closed loop $\Gamma$.

The expression of $\oint_\Gamma du/(u-u_0)$ can be written 
with the aid of a functional which will be called int($M,C$), defined for a
point $M$ and a closed curve $C$ in a two-dimensional plane, such that 
\begin{equation}
{\rm int}(M,C)=\left\{
\begin{array}{l}
1 \;\:{\rm if} \;\:M \;\:{\rm is \;\:an \;\:interior \;\:point \;\:of} \;\:C ,\\ 
0 \;\:{\rm if} \;\:M \;\:{\rm is \;\:exterior \;\:to}\:\; C .\\
\end{array}\right.
\label{nplanar-118}
\end{equation}
With this notation the result of the integration on a closed path $\Gamma$
can be written as 
\begin{equation}
\oint_\Gamma\frac{du}{u-u_0}=
\sum_{k=1}^{n/2}2\pi\tilde e_k 
\;{\rm int}(u_{0\xi_k\eta_k},\Gamma_{\xi_k\eta_k}) ,
\label{nplanar-119}
\end{equation}\index{poles and residues!planar n-complex}
where $u_{0\xi_k\eta_k}$ and $\Gamma_{\xi_k\eta_k}$ are respectively the
projections of the point $u_0$ and of 
the loop $\Gamma$ on the plane defined by the axes $\xi_k$ and $\eta_k$,
as shown in Fig. \ref{fig26}.

\begin{figure}
\begin{center}
\epsfig{file=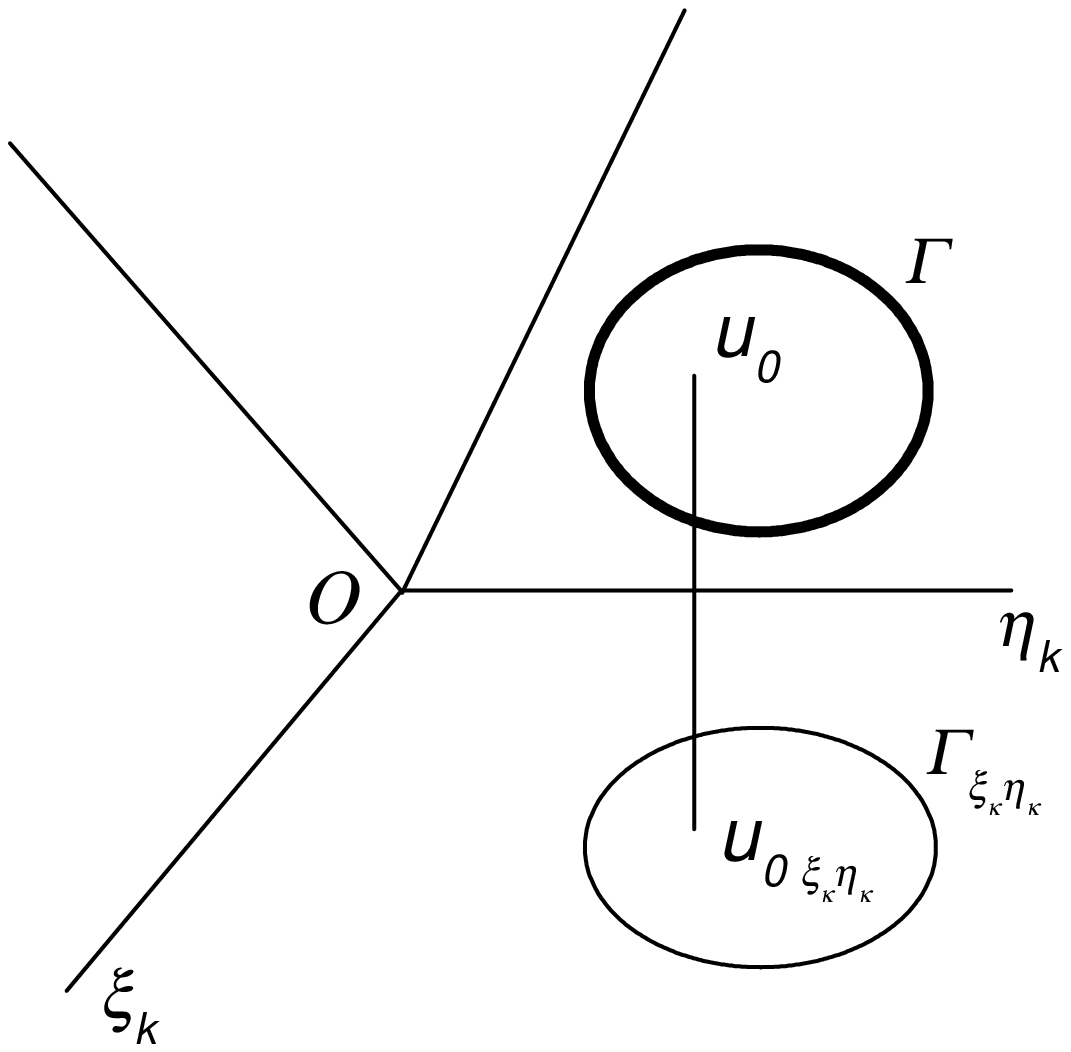,width=12cm}
\caption{Integration path $\Gamma$ and pole $u_0$, and their projections
$\Gamma_{\xi_k\eta_k}$ and $u_{0\xi_k\eta_k}$ on the plane $\xi_k \eta_k$. }
\label{fig26}
\end{center}
\end{figure}

If $f(u)$ is an analytic n-complex function which can be expanded in a
series as written in Eq. (1.89), and the expansion holds on the curve
$\Gamma$ and on a surface spanning $\Gamma$, then from Eqs. (\ref{nplanar-112b}) and
(\ref{nplanar-119}) it follows that
\begin{equation}
\oint_\Gamma \frac{f(u)du}{u-u_0}=
2\pi f(u_0)\sum_{k=1}^{n/2}\tilde e_k 
\;{\rm int}(u_{0\xi_k\eta_k},\Gamma_{\xi_k\eta_k}) .
\label{nplanar-120}
\end{equation}

Substituting in the right-hand side of 
Eq. (\ref{nplanar-120}) the expression of $f(u)$ in terms of the real 
components $P_k$, Eq. (\ref{nplanar-h93}), yields
\begin{eqnarray}
\lefteqn{\oint_\Gamma \frac{f(u)du}{u-u_0}
=\frac{2}{n}\sum_{k=1}^{n/2}\sum_{l=0}^{n-1}h_l
\nonumber}\\
&&\sum_{p=1}^{n-1}
(-1)^{[(l-p)/n]}\sin\left[\frac{\pi (2k-1)p}{n}\right] 
P_{n-p+l-n[(n-p+l)/n]}(u_0)
\:{\rm int}(u_{0\xi_k\eta_k},\Gamma_{\xi_k\eta_k}) .
\label{nplanar-121}
\end{eqnarray}
It the integral in Eq. (\ref{nplanar-121}) is written as 
\begin{equation}
\oint_\Gamma \frac{f(u)du}{u-u_0}=\sum_{l=0}^{n-1}h_l I_l,
\label{nplanar-122a}
\end{equation}
it can be checked that
\begin{equation}
\sum_{l=0}^{n-1} I_l=0.
\label{nplanar-122b}
\end{equation}

If $f(u)$ can be expanded as written in Eq. (1.89) on 
$\Gamma$ and on a surface spanning $\Gamma$, then from Eqs. (\ref{nplanar-112b}) and
(\ref{nplanar-119}) it also results that
\begin{equation}
\oint_\Gamma \frac{f(u)du}{(u-u_0)^{n+1}}=
\frac{2\pi}{n!}f^{(n)}(u_0)\sum_{k=1}^{[(n-1)/2]}\tilde e_k 
\;{\rm int}(u_{0\xi_k\eta_k},\Gamma_{\xi_k\eta_k}) ,
\label{nplanar-122}
\end{equation}
where the fact has been used  that the derivative $f^{(n)}(u_0)$ is related to
the expansion coefficient in Eq. (1.89) according to Eq. (1.93).

If a function $f(u)$ is expanded in positive and negative powers of $u-u_l$,
where $u_l$ are n-complex constants, $l$ being an index, the integral of $f$
on a closed loop $\Gamma$ is determined by the terms in the expansion of $f$
which are of the form $r_l/(u-u_l)$,
\begin{equation}
f(u)=\cdots+\sum_l\frac{r_l}{u-u_l}+\cdots.
\label{nplanar-123}
\end{equation}
Then the integral of $f$ on a closed loop $\Gamma$ is
\begin{equation}
\oint_\Gamma f(u) du = 
2\pi \sum_l\sum_{k=1}^{n/2}\tilde e_k 
\;{\rm int}(u_{l\xi_k\eta_k},\Gamma_{\xi_k\eta_k})r_l .
\label{nplanar-124}
\end{equation}

\subsection{Factorization of planar n-complex polynomials}

A polynomial of degree $m$ of the n-complex variable $u$ has the form
\begin{equation}
P_m(u)=u^m+a_1 u^{m-1}+\cdots+a_{m-1} u +a_m ,
\label{nplanar-125}
\end{equation}
where $a_l$, for $l=1,...,m$, are in general n-complex constants.
If $a_l=\sum_{p=0}^{n-1} h_p a_{lp}$, and with the
notations of Eqs. (\ref{nplanar-88b}), (\ref{nplanar-88c}) applied for $l= 1, \cdots, m$, the
polynomial $P_m(u)$ can be written as
\begin{eqnarray}
P_m= \sum_{k=1}^{n/2}\left\{(e_k v_k+\tilde e_k\tilde v_k)^m+
\sum_{l=1}^m(e_k A_{lk}+\tilde e_k\tilde A_{lk})
(e_k v_k+\tilde e_k\tilde v_k)^{m-l} 
\right\},
\label{nplanar-126a}
\end{eqnarray}\index{polynomial, canonical variables!planar n-complex}
where the constants $A_{lk}, \tilde A_{lk}$ are real numbers.

The polynomials of degree $m$ in $e_k v_k+\tilde e_k\tilde v_k$ 
in Eq. (\ref{nplanar-126a}) 
can always be written as a product of linear factors of the form
$e_k (v_k-v_{kp})+\tilde e_k(\tilde v_k-\tilde v_{kp})$, where the
constants $v_{kp}, \tilde v_{kp}$ are real,
\begin{eqnarray}
\lefteqn{(e_k v_k+\tilde e_k\tilde v_k)^m+
\sum_{l=1}^m(e_k A_{lk}+\tilde e_k\tilde A_{lk})
(e_k v_k+\tilde e_k\tilde v_k)^{m-l} \nonumber}\\
&&=\prod_{p=1}^{m}\left\{e_k (v_k-v_{kp})
+\tilde e_k(\tilde v_k-\tilde v_{kp})\right\}.
\label{nplanar-126b}
\end{eqnarray}
Then the polynomial $P_m$ can be written as
\begin{equation}
P_m=\sum_{k=1}^{n/2}\prod_{p=1}^m
\left\{e_k (v_k-v_{kp})+\tilde e_k(\tilde v_k-\tilde v_{kp})\right\}.
\label{nplanar-127a}
\end{equation}
Due to the relations  (\ref{nplanar-e12a}),
the polynomial $P_m(u)$ can be written as a product of factors of
the form  
\begin{eqnarray}
P_m(u)=\prod_{p=1}^m \left\{\sum_{k=1}^{n/2}
\left\{e_k (v_k-v_{kp})
+\tilde e_k(\tilde v_k-\tilde v_{kp})\right\}\right\}.
\label{nplanar-128}
\end{eqnarray}
This relation can be written with the aid of Eq. (\ref{nplanar-e13a}) as
\begin{eqnarray}
P_m(u)=\prod_{p=1}^m (u-u_p) ,
\label{nplanar-129a}
\end{eqnarray}\index{polynomial, factorization!planar n-complex}
where
\begin{eqnarray}
u_p=\sum_{k=1}^{n/2}\left(e_k v_{kp}
+\tilde e_k\tilde v_{kp}\right), 
\label{nplanar-129b}
\end{eqnarray}
for $p=1,...,m$.
For a given $k$, the roots  
$e_k v_{k1}+\tilde e_k\tilde v_{k1}, ...,  e_k v_{km}+\tilde e_k\tilde v_{km}$
defined in Eq. (\ref{nplanar-126b}) may be ordered arbitrarily. This means that Eq.
(\ref{nplanar-129b}) gives sets of $m$ roots 
$u_1,...,u_m$ of the polynomial $P_m(u)$, 
corresponding to the various ways in which the roots $e_k v_{kp}+\tilde
e_k\tilde v_{kp}$ are ordered according to $p$ for each value of $k$. 
Thus, while the n-complex components in Eq. (\ref{nplanar-126b}) taken separately
have  
unique factorizations, the polynomial $P_m(u)$ can be written in many different
ways as a product of linear factors. 

If $P(u)=u^2+1$, the degree is $m=2$, the coefficients of the polynomial are
$a_1=0, a_2=1$, the n-complex components of $a_2$ are $a_{20}=1, a_{21}=0,
... ,a_{2,n-1}=0$, the components $A_{2k}, \tilde A_{2k}$ calculated according
to Eqs. (\ref{nplanar-88b}), (\ref{nplanar-88c}) are $A_{2k}=1, \tilde A_{2k}=0, k=1,...,n/2$.
The left-hand side of Eq. (\ref{nplanar-126b}) has the form 
$(e_k v_k+\tilde e_k\tilde v_k)^2+e_k$, and since $e_k=-\tilde e_k^2$, the
right-hand side of Eq. (\ref{nplanar-126b}) is 
$\left\{e_k v_k+\tilde e_k(\tilde v_k+1)\right\}
\left\{e_k v_k+\tilde e_k(\tilde v_k-1)\right\}$, so that
$v_{kp}=0, \tilde v_{kp}=\pm 1, k=1,...,n/2, p=1,2$.
Then Eq. (\ref{nplanar-127a}) has the form $u^2+1=\sum_{k=1}^{n/2}
\left\{e_k v_k+\tilde e_k(\tilde v_k+1)\right\}
\left\{e_k v_k+\tilde e_k(\tilde v_k-1)\right\}$.
The factorization in Eq. (\ref{nplanar-129a}) is $u^2+1=(u-u_1)(u-u_2)$, where 
$u_1=\pm \tilde e_1\pm\tilde e_2\pm\cdots\pm \tilde e_{n/2}, u_2=-u_1$, so that
there are $2^{n/2-1}$ independent sets of roots $u_1,u_2$
of $u^2+1$. It can be checked
that $(\pm \tilde e_1\pm\tilde e_2\pm\cdots\pm \tilde e_{n/2})^2=
-e_1-e_2-\cdots-e_{n/2}=-1$.


\subsection{Representation of planar n-complex numbers by irreducible matrices}

If the unitary matrix written in Eq. (\ref{nplanar-11}) is called $T$,
it can be shown that the matrix $T U T^{-1}$ has the form 
\begin{equation}
T U T^{-1}=\left(
\begin{array}{cccc}
V_1      &     0   & \cdots  &   0   \\
0        &     V_2 & \cdots  &   0   \\
\vdots   &  \vdots & \cdots  &\vdots \\
0        &     0   & \cdots  &   V_{n/2}\\
\end{array}
\right),
\label{nplanar-129}
\end{equation}\index{representation by irreducible matrices!planar n-complex}
where $U$ is the matrix in Eq. (\ref{nplanar-24a}) used to represent the n-complex
number $u$. In Eq. (\ref{nplanar-129}), the matrices $V_k$ are
the matrices
\begin{equation}
V_k=\left(
\begin{array}{cc}
v_k           &     \tilde v_k   \\
-\tilde v_k   &     v_k          \\
\end{array}\right),
\label{nplanar-130}
\end{equation}
for $ k=1,...,n/2$, where $v_k, \tilde v_k$ are the variables introduced in Eqs. (\ref{nplanar-9a}) and
(\ref{nplanar-9b}), and the symbols 0 denote
the matrix
\begin{equation}
\left(
\begin{array}{cc}
0   &  0   \\
0   &  0   \\
\end{array}\right).
\label{nplanar-131}
\end{equation}
The relations between the variables $v_k, \tilde v_k$ for the multiplication of
n-complex numbers have been written in Eq. (\ref{nplanar-22}). The matrix
$T U T^{-1}$ provides an irreducible representation
\cite{4} of the n-complex number $u$ in terms of matrices with real
coefficients. For $n=2$, Eqs. (\ref{nplanar-9a}) and (\ref{nplanar-9b}) give $v_1=x_0, 
\tilde v_1=x_1$, and Eq. (\ref{nplanar-e11}) gives $e_1=1, \tilde e_1=h_1$, 
where according to Eq. (\ref{nplanar-1}) $h_1^2=-1$, so that the matrix $V_1$, Eq.
(\ref{nplanar-130}), is 
\begin{equation}
v_1=\left(
\begin{array}{cc}
x_0    &  x_1   \\
-x_1   &  x_0   \\
\end{array}\right),
\label{nplanar-132}
\end{equation}
which shows that, for $n=2$, the hypercomplex numbers $x_0+h_1 x_1$ are
identical to the usual 2-dimensional complex numbers $x+iy$.


\begin{theindex}

  \item amplitude
    \subitem circular fourcomplex, 56
    \subitem hyperbolic fourcomplex, 78
    \subitem planar 6-complex, 170
    \subitem planar fourcomplex, 91
    \subitem planar n-complex, 214
    \subitem polar 5-complex, 142
    \subitem polar 6-complex, 158
    \subitem polar n-complex, 184
    \subitem tricomplex, 28
    \subitem twocomplex, 9
  \item analytic function, hypercomplex
    \subitem definition, 17
    \subitem expansion in series, 17
    \subitem power series, 17
  \item azimuthal angle
    \subitem polar fourcomplex, 116
  \item azimuthal angle, tricomplex, 25
  \item azimuthal angles
    \subitem circular fourcomplex, 57
    \subitem planar 6-complex, 170
    \subitem planar fourcomplex, 92
    \subitem planar n-complex, 216
    \subitem polar 5-complex, 142
    \subitem polar 6-complex, 158
    \subitem polar n-complex, 187

  \indexspace

  \item canonical base
    \subitem circular fourcomplex, 67
    \subitem hyperbolic fourcomplex, 80
    \subitem planar 6-complex, 173
    \subitem planar fourcomplex, 105
    \subitem planar n-complex, 222
    \subitem polar 5-complex, 148
    \subitem polar 6-complex, 162
    \subitem polar fourcomplex, 118
    \subitem polar n-complex, 194
    \subitem tricomplex, 42
    \subitem twocomplex, 10
  \item canonical form
    \subitem circular fourcomplex, 67
    \subitem hyperbolic fourcomplex, 80
    \subitem planar 6-complex, 174
    \subitem planar fourcomplex, 105
    \subitem planar n-complex, 222
    \subitem polar 5-complex, 148
    \subitem polar 6-complex, 163
    \subitem polar fourcomplex, 118
    \subitem polar n-complex, 195
    \subitem twocomplex, 11
  \item canonical variables
    \subitem circular fourcomplex, 57
    \subitem hyperbolic fourcomplex, 78
    \subitem planar 6-complex, 169
    \subitem planar fourcomplex, 92
    \subitem planar n-complex, 214
    \subitem polar 5-complex, 140
    \subitem polar 6-complex, 156
    \subitem polar fourcomplex, 116
    \subitem polar n-complex, 184
    \subitem tricomplex, 42
    \subitem twocomplex, 9
  \item complex units
    \subitem circular fourcomplex, 55
    \subitem hyperbolic fourcomplex, 76
    \subitem planar 6-complex, 167
    \subitem planar fourcomplex, 90
    \subitem planar n-complex, 211
    \subitem polar 5-complex, 140
    \subitem polar 6-complex, 156
    \subitem polar fourcomplex, 113
    \subitem polar n-complex, 181
    \subitem tricomplex, 22
    \subitem twocomplex, 8
  \item complex units, powers of
    \subitem circular fourcomplex, 61
    \subitem hyperbolic fourcomplex, 81
    \subitem planar fourcomplex, 95
    \subitem planar n-complex, 218
    \subitem polar fourcomplex, 119
    \subitem polar n-complex, 190
    \subitem tricomplex, 30
    \subitem twocomplex, 12
  \item convergence of power series
    \subitem circular fourcomplex, 67
    \subitem hyperbolic fourcomplex, 85
    \subitem planar 6-complex, 175
    \subitem planar fourcomplex, 105
    \subitem planar n-complex, 228
    \subitem polar 5-complex, 150
    \subitem polar 6-complex, 164
    \subitem polar fourcomplex, 129
    \subitem polar n-complex, 203
    \subitem tricomplex, 42
    \subitem twocomplex, 16
  \item convergence, region of
    \subitem circular fourcomplex, 68
    \subitem hyperbolic fourcomplex, 86
    \subitem planar 6-complex, 175
    \subitem planar fourcomplex, 106
    \subitem planar n-complex, 228
    \subitem polar 5-complex, 150
    \subitem polar 6-complex, 165
    \subitem polar fourcomplex, 130
    \subitem polar n-complex, 203
    \subitem tricomplex, 44
    \subitem twocomplex, 16
  \item cosexponential function, polar 5-complex
    \subitem definition, 143
  \item cosexponential function, polar fourcomplex
    \subitem definitions, 119
  \item cosexponential function, polar n-complex
    \subitem definition, 190
  \item cosexponential functions, planar 6-complex
    \subitem addition theorems, 173
    \subitem definition, 170
    \subitem differential equations, 173
    \subitem expressions, 171
    \subitem parity, 171
  \item cosexponential functions, planar fourcomplex
    \subitem addition theorems, 96
    \subitem definitions, 96
    \subitem differential equations, 98
    \subitem expressions, 98
    \subitem parity, 96
  \item cosexponential functions, planar n-complex
    \subitem addition theorems, 220
    \subitem definition, 219
    \subitem differential equations, 222
    \subitem expressions, 219
    \subitem parity, 219
  \item cosexponential functions, polar 5-complex
    \subitem addition theorems, 146
    \subitem differential equations, 148
    \subitem expressions, 144, 145
  \item cosexponential functions, polar 6-complex
    \subitem addition theorems, 160
    \subitem definitions, 159
    \subitem differential equations, 162
    \subitem expressions, 160
    \subitem parity, 159
  \item cosexponential functions, polar fourcomplex
    \subitem addition theorems, 120
    \subitem differential equations, 122
    \subitem expressions, 122
    \subitem parity, 120
  \item cosexponential functions, polar n-complex
    \subitem addition theorems, 192
    \subitem differential equations, 193
    \subitem expressions, 191
    \subitem parity, 190
  \item cosexponential functions, tricomplex
    \subitem addition theorems, 32
    \subitem definitions, 32
    \subitem differential equations, 34
    \subitem expressions, 34

  \indexspace

  \item derivative, hypercomplex
    \subitem definition, 17
    \subitem of power function, 17
  \item derivative, independence of direction
    \subitem fourcomplex, 69
    \subitem planar n-complex, 229
    \subitem polar n-complex, 204
    \subitem tricomplex, 46
    \subitem twocomplex, 18
  \item distance
    \subitem circular fourcomplex, 57
    \subitem hyperbolic fourcomplex, 79
    \subitem planar 6-complex, 169
    \subitem planar fourcomplex, 92
    \subitem planar n-complex, 215
    \subitem polar 5-complex, 142
    \subitem polar 6-complex, 156
    \subitem polar fourcomplex, 115
    \subitem polar n-complex, 185
    \subitem tricomplex, 26
    \subitem twocomplex, 9
  \item distance, canonical variables
    \subitem polar fourcomplex, 116
  \item divisors of zero
    \subitem hyperbolic fourcomplex, 78
    \subitem planar n-complex, 218
    \subitem polar n-complex, 190
    \subitem tricomplex, 23
    \subitem twocomplex, 9

  \indexspace

  \item exponential form
    \subitem circular fourcomplex, 63
    \subitem hyperbolic fourcomplex, 82
    \subitem planar 6-complex, 174
    \subitem planar fourcomplex, 101
    \subitem planar n-complex, 223
    \subitem polar 5-complex, 148
    \subitem polar 6-complex, 163
    \subitem polar fourcomplex, 124
    \subitem polar n-complex, 196
    \subitem tricomplex, 37
    \subitem twocomplex, 12
  \item exponential function, hypercomplex
    \subitem addition theorem, 11
    \subitem definition, 11
  \item exponential, expression
    \subitem twocomplex, 12
  \item exponential, expressions
    \subitem circular fourcomplex, 61
    \subitem hyperbolic fourcomplex, 81
    \subitem planar 6-complex, 171, 174
    \subitem planar fourcomplex, 95
    \subitem planar n-complex, 219
    \subitem polar 5-complex, 143, 149
    \subitem polar 6-complex, 160, 163
    \subitem polar fourcomplex, 119
    \subitem polar n-complex, 190, 200
    \subitem tricomplex, 30

  \indexspace

  \item functions, real components
    \subitem circular fourcomplex, 68
    \subitem planar 6-complex, 176
    \subitem planar n-complex, 229
    \subitem polar 5-complex, 151
    \subitem polar n-complex, 204
    \subitem tricomplex, 44
    \subitem twocomplex, 18

  \indexspace

  \item hyperbolic functions, expressions
    \subitem circular fourcomplex, 65
    \subitem hyperbolic fourcomplex, 84
    \subitem planar 6-complex, 175
    \subitem planar fourcomplex, 103
    \subitem planar n-complex, 226
    \subitem polar 5-complex, 149
    \subitem polar 6-complex, 164
    \subitem polar fourcomplex, 128
    \subitem polar n-complex, 200, 201
    \subitem tricomplex, 40
    \subitem twocomplex, 14
  \item hyperbolic functions, hypercomplex
    \subitem addition theorems, 14
    \subitem definitions, 14

  \indexspace

  \item integrals, path
    \subitem circular fourcomplex, 70
    \subitem hyperbolic fourcomplex, 88
    \subitem planar 6-complex, 176
    \subitem planar fourcomplex, 108
    \subitem planar n-complex, 230
    \subitem polar 5-complex, 151
    \subitem polar 6-complex, 165
    \subitem polar fourcomplex, 132
    \subitem polar n-complex, 205
    \subitem tricomplex, 47
    \subitem twocomplex, 19
  \item inverse
    \subitem circular fourcomplex, 56
    \subitem hyperbolic fourcomplex, 77
    \subitem planar fourcomplex, 91
    \subitem planar n-complex, 213
    \subitem polar fourcomplex, 114
    \subitem polar n-complex, 182
    \subitem tricomplex, 23
    \subitem twocomplex, 9
  \item inverse, determinant
    \subitem hyperbolic fourcomplex, 78
    \subitem planar fourcomplex, 91
    \subitem planar n-complex, 214
    \subitem polar fourcomplex, 114
    \subitem polar n-complex, 184
    \subitem tricomplex, 23, 30
    \subitem twocomplex, 9

  \indexspace

  \item logarithm
    \subitem circular fourcomplex, 65
    \subitem hyperbolic fourcomplex, 83
    \subitem planar 6-complex, 174
    \subitem planar fourcomplex, 103
    \subitem planar n-complex, 224
    \subitem polar 5-complex, 149
    \subitem polar 6-complex, 163
    \subitem polar fourcomplex, 126
    \subitem polar n-complex, 198
    \subitem tricomplex, 38
    \subitem twocomplex, 13

  \indexspace

  \item matrix representation
    \subitem circular fourcomplex, 60
    \subitem hyperbolic fourcomplex, 81
    \subitem planar 6-complex, 170
    \subitem planar fourcomplex, 95
    \subitem planar n-complex, 218
    \subitem polar 5-complex, 143
    \subitem polar 6-complex, 159
    \subitem polar fourcomplex, 119
    \subitem polar n-complex, 189
    \subitem tricomplex, 28
    \subitem twocomplex, 11
  \item modulus
    \subitem hyperbolic fourcomplex, 79
  \item modulus, canonical variables
    \subitem circular fourcomplex, 67
    \subitem hyperbolic fourcomplex, 80
    \subitem planar 6-complex, 170
    \subitem planar fourcomplex, 105
    \subitem planar n-complex, 215
    \subitem polar 6-complex, 158
    \subitem polar n-complex, 186
    \subitem tricomplex, 44
    \subitem twocomplex, 10
  \item modulus, definition
    \subitem circular fourcomplex, 66
    \subitem hyperbolic fourcomplex, 85
    \subitem planar fourcomplex, 104
    \subitem planar n-complex, 215
    \subitem polar fourcomplex, 128
    \subitem polar n-complex, 185
    \subitem tricomplex, 41
    \subitem twocomplex, 15
  \item modulus, inequalities
    \subitem circular fourcomplex, 66
    \subitem hyperbolic fourcomplex, 85
    \subitem planar 6-complex, 175
    \subitem planar fourcomplex, 104
    \subitem planar n-complex, 227
    \subitem polar 5-complex, 142, 150
    \subitem polar 6-complex, 164
    \subitem polar fourcomplex, 129
    \subitem polar n-complex, 202
    \subitem tricomplex, 41
    \subitem twocomplex, 15

  \indexspace

  \item nodal hyperplanes
    \subitem circular fourcomplex, 57
    \subitem planar fourcomplex, 92
    \subitem polar fourcomplex, 115
  \item nodal hypersurfaces
    \subitem planar n-complex, 214
    \subitem polar n-complex, 185
  \item nodal line, tricomplex, 23
  \item nodal lines
    \subitem twocomplex, 9
  \item nodal plane, tricomplex, 23

  \indexspace

  \item planar angle
    \subitem circular fourcomplex, 57
    \subitem planar fourcomplex, 92
    \subitem polar 5-complex, 142
    \subitem polar 6-complex, 158
  \item planar angles
    \subitem planar 6-complex, 170
    \subitem planar n-complex, 216
    \subitem polar n-complex, 187
  \item polar angle
    \subitem polar 5-complex, 142
  \item polar angle, tricomplex, 26
  \item polar angles
    \subitem polar 6-complex, 158
    \subitem polar fourcomplex, 117
    \subitem polar n-complex, 187
  \item poles and residues
    \subitem circular fourcomplex, 71
    \subitem planar 6-complex, 176
    \subitem planar fourcomplex, 109
    \subitem planar n-complex, 232
    \subitem polar 5-complex, 151
    \subitem polar 6-complex, 165
    \subitem polar fourcomplex, 133
    \subitem polar n-complex, 206
    \subitem tricomplex, 47
  \item polynomial, canonical variables
    \subitem circular fourcomplex, 75
    \subitem hyperbolic fourcomplex, 88
    \subitem planar 6-complex, 176
    \subitem planar fourcomplex, 111
    \subitem planar n-complex, 233
    \subitem polar 5-complex, 152
    \subitem polar fourcomplex, 135
    \subitem polar n-complex, 208
    \subitem tricomplex, 51
    \subitem twocomplex, 19
  \item polynomial, factorization
    \subitem circular fourcomplex, 75
    \subitem hyperbolic fourcomplex, 89
    \subitem planar 6-complex, 177
    \subitem planar fourcomplex, 112
    \subitem planar n-complex, 233
    \subitem polar 5-complex, 152
    \subitem polar fourcomplex, 136
    \subitem polar n-complex, 209
    \subitem tricomplex, 51
    \subitem twocomplex, 20
  \item polynomials, canonical variables
    \subitem polar 6-complex, 166
  \item polynomials, factorization
    \subitem polar 6-complex, 166
  \item power function
    \subitem circular fourcomplex, 65
    \subitem hyperbolic fourcomplex, 83, 84
    \subitem planar 6-complex, 174
    \subitem planar fourcomplex, 103
    \subitem planar n-complex, 224
    \subitem polar 5-complex, 149
    \subitem polar 6-complex, 163
    \subitem polar fourcomplex, 127
    \subitem polar n-complex, 199
    \subitem tricomplex, 38
    \subitem twocomplex, 13
  \item power series
    \subitem circular fourcomplex, 67
    \subitem hyperbolic fourcomplex, 85
    \subitem planar 6-complex, 175
    \subitem planar fourcomplex, 105
    \subitem planar n-complex, 227
    \subitem polar 5-complex, 150
    \subitem polar 6-complex, 164
    \subitem polar fourcomplex, 129
    \subitem polar n-complex, 202
    \subitem tricomplex, 42
    \subitem twocomplex, 15
  \item product
    \subitem circular fourcomplex, 55, 56
    \subitem hyperbolic fourcomplex, 76
    \subitem planar 6-complex, 169
    \subitem planar fourcomplex, 90
    \subitem planar n-complex, 213
    \subitem polar 5-complex, 140
    \subitem polar 6-complex, 156
    \subitem polar fourcomplex, 113
    \subitem polar n-complex, 182
    \subitem tricomplex, 22, 23
    \subitem twocomplex, 8

  \indexspace

  \item relations between partial derivatives
    \subitem circular fourcomplex, 69
    \subitem hyperbolic fourcomplex, 87
    \subitem planar 6-complex, 176
    \subitem planar fourcomplex, 106
    \subitem planar n-complex, 229
    \subitem polar 5-complex, 151
    \subitem polar 6-complex, 165
    \subitem polar fourcomplex, 130
    \subitem polar n-complex, 204
    \subitem tricomplex, 46
  \item relations between second-order derivatives
    \subitem circular fourcomplex, 69
    \subitem hyperbolic fourcomplex, 87
    \subitem planar 6-complex, 176
    \subitem planar fourcomplex, 107
    \subitem planar n-complex, 229
    \subitem polar 5-complex, 151
    \subitem polar 6-complex, 165
    \subitem polar fourcomplex, 131
    \subitem polar n-complex, 204
    \subitem tricomplex, 46
    \subitem twocomplex, 18
  \item representation by irreducible matrices, 52
    \subitem circular fourcomplex, 76
    \subitem hyperbolic fourcomplex, 90
    \subitem planar 6-complex, 177
    \subitem planar fourcomplex, 113
    \subitem planar n-complex, 234
    \subitem polar 5-complex, 153
    \subitem polar 6-complex, 167
    \subitem polar fourcomplex, 136
    \subitem polar n-complex, 210
    \subitem twocomplex, 20
  \item roots
    \subitem circular fourcomplex, 73
    \subitem planar fourcomplex, 110
  \item rules for derivation and integration, 18

  \indexspace

  \item series
    \subitem circular fourcomplex, 66
    \subitem hyperbolic fourcomplex, 84
    \subitem planar fourcomplex, 104
    \subitem planar n-complex, 227
    \subitem polar fourcomplex, 128
    \subitem polar n-complex, 202
    \subitem tricomplex, 40
    \subitem twocomplex, 15
  \item sum
    \subitem circular fourcomplex, 55
    \subitem hyperbolic fourcomplex, 76
    \subitem planar 6-complex, 167
    \subitem planar fourcomplex, 90
    \subitem planar n-complex, 213
    \subitem polar 5-complex, 140
    \subitem polar 6-complex, 156
    \subitem polar fourcomplex, 113
    \subitem polar n-complex, 182
    \subitem tricomplex, 22
    \subitem twocomplex, 8

  \indexspace

  \item transformation of variables
    \subitem circular fourcomplex, 58, 67
    \subitem hyperbolic fourcomplex, 79
    \subitem planar 6-complex, 170
    \subitem planar fourcomplex, 93, 102, 106
    \subitem planar n-complex, 216
    \subitem polar 5-complex, 142
    \subitem polar 6-complex, 158
    \subitem polar fourcomplex, 115, 116, 125, 126
    \subitem polar n-complex, 189
    \subitem tricomplex, 26, 28, 44
    \subitem twocomplex, 10
  \item trigonometric form
    \subitem circular fourcomplex, 64
    \subitem expression, 12
    \subitem hyperbolic fourcomplex, 82
    \subitem planar 6-complex, 174
    \subitem planar fourcomplex, 102
    \subitem planar n-complex, 223
    \subitem polar 5-complex, 149
    \subitem polar 6-complex, 163
    \subitem polar fourcomplex, 125
    \subitem polar n-complex, 197
    \subitem tricomplex, 37
  \item trigonometric functions, expressions
    \subitem circular fourcomplex, 65
    \subitem hyperbolic fourcomplex, 84
    \subitem planar 6-complex, 175
    \subitem planar fourcomplex, 103
    \subitem planar n-complex, 225, 226
    \subitem polar 5-complex, 149
    \subitem polar 6-complex, 164
    \subitem polar fourcomplex, 128
    \subitem polar n-complex, 199, 201
    \subitem tricomplex, 38
    \subitem twocomplex, 14
  \item trigonometric functions, hypercomplex
    \subitem addition theorems, 14
    \subitem definitions, 13
  \item trisector line, 25

\end{theindex}


\newpage


\begin{thebibliography}{9}
\bibitem{1} G. Birkhoff and S. MacLane, {\it Modern Algebra} (Macmillan, New York,
Third Edition 1965), p. 222.
\bibitem{2a} B. L. van der Waerden, {\it Modern Algebra} (F. Ungar, New York, 
Third Edition 1950),
vol. II, p. 133.
\bibitem{2} O. Taussky, Algebra, in {\it Handbook of Physics}, edited by E. U.
Condon and H. Odishaw (McGraw-Hill, New York, Second Edition 1958), p. I-22.
\bibitem{2b} 
D. Kaledin, arXiv:alg-geom/9612016;
K. Scheicher, R. F. Tichy, and K. W. Tomantschger, Anzeiger Abt.
II 134, 3 (1997); 
S. De Leo and P. Rotelli, arXiv:funct-an/9701004, 9703002;
M. Verbitsky, arXiv:alg-geom/9703016;
S. De Leo, arXiv:physics/9703033;
J. D. E. Grant and I. A. B. Strachan, arXiv:solv-int/9808019;
D. M. J. Calderbank and P. Tod, arXiv:math.DG/9911121;
L. Ornea and P. Piccinni, arXiv:math.DG/0001066.
\bibitem{2c} S. Olariu, arXiv:math.OA/0007180, math.CV/0008119-0008125.
\bibitem{3} E. T. Whittaker and G. N. Watson {\it A Course of Modern
Analysis}, (Cambridge University Press, Fourth Edition 1958), p. 83.
\bibitem{4} E. Wigner, {\it Group Theory} (Academic Press, New York, 1959), p.
73.

\end{thebibliography}
\end{document}